\documentclass[english,10pt]{amsbook}

\usepackage{amsmath, amssymb, amsthm}
\usepackage[all]{xy}
\usepackage[activeacute]{babel}

\newtheorem{thm}{Theorem}[section]
\newtheorem{lem}[thm]{Lemma}
\newtheorem{prop}[thm]{Proposition}
\newtheorem{cor}[thm]{Corollary}

\newtheorem{defi}[thm]{Definition}
\newtheorem{exam}[thm]{Example}
\newtheorem{rmk}[thm]{Remark}
\newtheorem{rmks}[thm]{Remark}
\newtheorem{note}[thm]{Note}

\newcommand{\modelcat}{\mathcal A}
\newcommand{\modelcattwo}{\mathcal B}
\newcommand{\modelcatthree}{\mathcal D}
\newcommand{\simpsets}{\mathbf{SSets}}						
\newcommand{\prescellcplx}{i:A\rightarrow B, \; A=A_{0}\rightarrow A_{1} \rightarrow \cdots \rightarrow A_{\beta}
						\rightarrow \cdots (\beta <\lambda),\; \{T^{\beta},e^{\beta},h^{\beta}\}_{\beta <\lambda}}
\newcommand{\prescellcplxtwo}{f:X\rightarrow Y, \; X=X_{0}\rightarrow X_{1} \rightarrow \cdots \rightarrow X_{\beta}
						\rightarrow \cdots (\beta <\lambda),\; \{T^{\beta},e^{\beta},h^{\beta}\}_{\beta <\lambda}}
\newcommand{\presabscellcplx}{i:\emptyset\rightarrow X, \; \emptyset=X_{0}\rightarrow X_{1} \rightarrow \cdots \rightarrow X_{\beta}
						\rightarrow \cdots (\beta <\lambda),\; \{T^{\beta},e^{\beta},h^{\beta}\}_{\beta <\lambda}}
\newcommand{\colimit}{\mathrm{colim}}
\newcommand{\lambdaseq}{A_{0}\rightarrow A_{1}\rightarrow \cdots \rightarrow A_{\beta}\rightarrow \cdots
		\; \; (\beta <\lambda)}
\newcommand{\lambdaseqtwo}{X_{0}\rightarrow X_{1}\rightarrow \cdots \rightarrow X_{\beta}\rightarrow \cdots
		\; \; (\beta <\lambda)}
\newcommand{\smoothS}{Sm|_{S}}							
\newcommand{\nissite}{(Sm|_{S})_{{Nis}}}
\newcommand{\simppre}{\Delta ^{op}Pre\nissite}
\newcommand{\motivic}{\mathcal{M}}
\newcommand{\pointedmotivic}{\mathcal{M}_{\ast}}
\newcommand{\simpprepointed}{\Delta ^{op}Pre_{\ast}\nissite}
\newcommand{\Hom}{\mathrm{Hom}}
\newcommand{\hocat}{\mathrm{Ho}\modelcat}
\newcommand{\hocattwo}{\mathrm{Ho}\modelcattwo}
\newcommand{\hocatleft}{\mathrm{Ho}L_{\mathcal V}\modelcat}
\newcommand{\hocatright}{\mathrm{Ho}R_{\mathcal V}\modelcat}
\newcommand{\inthompre}{\mathbf{Hom}_{Pre}}
\newcommand{\inthomprepointed}{\mathbf{Hom}_{\pointedmotivic}}
\newcommand{\inthomspectrapresheaf}{\mathbf{hom}_{r}}
\newcommand{\affineS}{\mathbb A ^{1}_{S}}
\newcommand{\tee}{S^{1}\wedge \mathbb G _{m}}
\newcommand{\gm}{\mathbb G_{m}}
\newcommand{\multgroup}{\mathbb G _{m}}
\newcommand{\Tspectra}{\mathrm{Spt}_{T}\nissite}
\newcommand{\motivicTspectra}{\mathrm{Spt}_{T}\pointedmotivic}
\newcommand{\Tsuspension}{\Sigma ^{\infty}_{T}}
\newcommand{\Tsuspfunctor}{\Sigma _{T}}
\newcommand{\fakeTsuspfunctor}{\Sigma _{T}^{\ell}}
\newcommand{\Tloops}{\Omega _{T}}
\newcommand{\fakeTloops}{\Omega _{T}^{\ell}}
\newcommand{\stablehomotopy}{\mathcal{SH}(S)}
\newcommand{\stablehomotopyeff}{\mathcal{SH}^{eff}(S)}
\newcommand{\stablehomotopyeffq}{\Sigma _{T}^{q}\mathcal{SH}^{eff}(S)}
\newcommand{\stablehomotopyeffqplusone}{\Sigma _{T}^{q+1}\mathcal{SH}^{eff}(S)}
\newcommand{\stablehomotopyeffqminusone}{\Sigma _{T}^{q-1}\mathcal{SH}^{eff}(S)}
\newcommand{\spherespectrum}{F_{0}(S^{0})}
\newcommand{\generatorNRS}{F_{n}(S^{r}\wedge \gm ^{s}\wedge U_{+})}

\newcommand{\generatorNNRSS}{F_{n+1}(S^{r}\wedge \gm ^{s+1}\wedge U_{+})}

\newcommand{\compactgenerators}{C=\bigcup _{n,r,s\geq 0}\; \; \bigcup _{U\in (\smoothS)}\generatorNRS}
\newcommand{\effcompactgenerators}{C_{eff}=\bigcup _{n,r,s\geq 0;s-n\geq 0}\; \; \bigcup _{U\in (\smoothS)}\generatorNRS}
\newcommand{\qeffcompactgenerators}{C_{eff}^{q}=\bigcup _{n,r,s\geq 0;s-n\geq q}\; \; \bigcup _{U\in (\smoothS)}\generatorNRS}
\newcommand{\zeroconnectedTspectra}{R_{C^{0}_{eff}}\motivicTspectra}
\newcommand{\qconnectedTspectra}{R_{C^{q}_{eff}}\motivicTspectra}
\newcommand{\qconnectedstablehomotopy}{R_{C^{q}_{eff}}\stablehomotopy}
\newcommand{\qplusoneconnectedstablehomotopy}{R_{C^{q+1}_{eff}}\stablehomotopy}

\newcommand{\diskNRS}{F_{n}(D^{r}\wedge \gm ^{s}\wedge U_{+})}
\newcommand{\diskNRRS}{F_{n}(D^{r+1}\wedge \gm ^{s}\wedge U_{+})}
\newcommand{\killNRS}{\iota ^{U}_{n,r,s}:\generatorNRS \rightarrow \diskNRRS}
\newcommand{\killNRSmap}{\iota ^{U}_{n,r,s}}
\newcommand{\qleftlocmaps}{L(<q)}
\newcommand{\qplusoneleftlocmaps}{L(<q+1)}

\newcommand{\weightqTspectra}{L_{<q}\motivicTspectra}
\newcommand{\weightqplusoneTspectra}{L_{<q+1}\motivicTspectra}
\newcommand{\weightqstablehomotopy}{L_{<q}\stablehomotopy}
\newcommand{\weightqplusonestablehomotopy}{L_{<q+1}\stablehomotopy}

\newcommand{\qslicegenerators}{S(q)}
\newcommand{\qsliceTspectra}{S^{q}\motivicTspectra}
\newcommand{\qslicestablehomotopy}{S^{q}\stablehomotopy}
\newcommand{\nsymmgroup}{\Sigma _{n}}
\newcommand{\symmTspectra}{\mathrm{Spt}_{T}^{\Sigma}\nissite}
\newcommand{\motivicsymmTspectra}{\mathrm{Spt}_{T}^{\Sigma}\pointedmotivic}
\newcommand{\inthomsymmTspectrapresheaf}{\mathbf{hom}_{r}^{\Sigma}}
\newcommand{\symmsequences}{(\pointedmotivic )^{\Sigma}}
\newcommand{\inthomsymmTspectra}{\mathbf{Hom}_{Spt_{T}^{\Sigma}}}
\newcommand{\symmstablehomotopy}{\mathcal{SH}^{\Sigma}(S)}
\newcommand{\ksymmstablehomotopy}{\mathcal{SH}^{\Sigma}(k)}
\newcommand{\symmspherespectrum}{\mathbf{1}}
\newcommand{\symmgeneratorNRS}{F_{n}^{\Sigma}(S^{r}\wedge \gm ^{s}\wedge U_{+})}
\newcommand{\symmgeneratorKLM}{F_{j}^{\Sigma}(S^{k}\wedge \gm ^{l}\wedge V_{+})}
\newcommand{\symmgeneratorKLMNRS}{F_{j+n}^{\Sigma}(S^{k+r}\wedge \gm ^{l+s}\wedge U\times_{S}V_{+})}
\newcommand{\symmcompactgenerators}{C^{\Sigma}=\bigcup _{n,r,s\geq 0}\; \; \bigcup _{U\in (\smoothS)}\symmgeneratorNRS}
\newcommand{\qeffsymmcompactgenerators}{C_{eff}^{q,\Sigma}=\bigcup _{n,r,s\geq 0;s-n\geq q}\; \; \bigcup _{U\in (\smoothS)}\symmgeneratorNRS}
\newcommand{\zeroconnectedsymmTspectra}{R_{C^{0}_{eff}}\motivicsymmTspectra}
\newcommand{\qconnectedsymmTspectra}{R_{C^{q}_{eff}}\motivicsymmTspectra}
\newcommand{\qplusoneconnectedsymmTspectra}{R_{C^{q+1}_{eff}}\motivicsymmTspectra}
\newcommand{\qminusoneconnectedsymmTspectra}{R_{C^{q-1}_{eff}}\motivicsymmTspectra}
\newcommand{\pconnectedsymmTspectra}{R_{C^{p}_{eff}}\motivicsymmTspectra}
\newcommand{\pplusoneconnectedsymmTspectra}{R_{C^{p+1}_{eff}}\motivicsymmTspectra}
\newcommand{\pplusoneplusqconnectedsymmTspectra}{R_{C^{p+q+1}_{eff}}\motivicsymmTspectra}
\newcommand{\pplusqconnectedsymmTspectra}{R_{C^{p+q}_{eff}}\motivicsymmTspectra}
\newcommand{\zeroconnectedsymmstablehomotopy}{R_{C^{0}_{eff}}\symmstablehomotopy}
\newcommand{\pconnectedsymmstablehomotopy}{R_{C^{p}_{eff}}\symmstablehomotopy}
\newcommand{\qconnectedsymmstablehomotopy}{R_{C^{q}_{eff}}\symmstablehomotopy}
\newcommand{\qplusoneconnectedsymmstablehomotopy}{R_{C^{q+1}_{eff}}\symmstablehomotopy}
\newcommand{\pplusqconnectedsymmstablehomotopy}{R_{C^{p+q}_{eff}}\symmstablehomotopy}
\newcommand{\qleftsymmlocmaps}{L^{\Sigma}(<q)}
\newcommand{\qplusoneleftsymmlocmaps}{L^{\Sigma}(<q+1)}
\newcommand{\pplusqleftsymmlocmaps}{L^{\Sigma}(<p+q)}

\newcommand{\weightqsymmTspectra}{L_{<q}\motivicsymmTspectra}

\newcommand{\weightonesymmTspectra}{L_{<1}\motivicsymmTspectra}
\newcommand{\weightrsymmTspectra}{L_{<r}\motivicsymmTspectra}
\newcommand{\weightqplusonesymmTspectra}{L_{<q+1}\motivicsymmTspectra}
\newcommand{\weightpplusonesymmTspectra}{L_{<p+1}\motivicsymmTspectra}
\newcommand{\weightpsymmTspectra}{L_{<p}\motivicsymmTspectra}
\newcommand{\weightpplusqsymmTspectra}{L_{<p+q}\motivicsymmTspectra}
\newcommand{\weightpplusqplusonesymmTspectra}{L_{<p+q+1}\motivicsymmTspectra}
\newcommand{\weightqsymmstablehomotopy}{L_{<q}\symmstablehomotopy}
\newcommand{\weightqplusonesymmstablehomotopy}{L_{<q+1}\symmstablehomotopy}
\newcommand{\symmdiskNRRS}{F_{n}^{\Sigma}(D^{r+1}\wedge \gm ^{s}\wedge U_{+})}
\newcommand{\symmkillNRS}{V(\iota ^{U}_{n,r,s}):\symmgeneratorNRS \rightarrow \symmdiskNRRS}
\newcommand{\symmqslicegenerators}{S^{\Sigma}(q)}
\newcommand{\zeroslicesymmTspectra}{S^{0}\motivicsymmTspectra}
\newcommand{\qslicesymmTspectra}{S^{q}\motivicsymmTspectra}
\newcommand{\qminusoneslicesymmTspectra}{S^{q-1}\motivicsymmTspectra}
\newcommand{\pslicesymmTspectra}{S^{p}\motivicsymmTspectra}
\newcommand{\pplusqslicesymmTspectra}{S^{p+q}\motivicsymmTspectra}
\newcommand{\pslicesymmstablehomotopy}{S^{p}\symmstablehomotopy}
\newcommand{\qslicesymmstablehomotopy}{S^{q}\symmstablehomotopy}
\newcommand{\pplusqslicesymmstablehomotopy}{S^{p+q}\symmstablehomotopy}
\newcommand{\zeroslicesymmstablehomotopy}{S^{0}\symmstablehomotopy}
\newcommand{\modules}{\mathrm{mod}}
\newcommand{\inthomAmodpresheaf}{\mathbf{hom}_{r}^{A\text{-}\modules}}
\newcommand{\inthomAAmodpresheaf}{\mathbf{hom}_{r}^{A'\text{-}\modules}}
\newcommand{\inthomAmod}{\mathbf{Hom}_{A\text{-}\modules}}
\newcommand{\motivicAmod}{A\text{-}\modules (\pointedmotivic)}
\newcommand{\motivicAAmod}{A'\text{-}\modules (\pointedmotivic)}

\newcommand{\algebras}{\mathrm{alg}}

\newcommand{\Amodstablehomotopy}{\mathcal{SH}({A\text{-}\modules})}
\newcommand{\Amodcompactgenerators}{C^{m}=\bigcup _{n,r,s\geq 0}\; \; \bigcup _{U\in (\smoothS)}A\wedge \symmgeneratorNRS}
\newcommand{\qconnectedAmod}{R_{C^{q}_{eff}}\motivicAmod}
\newcommand{\pconnectedAmod}{R_{C^{p}_{eff}}\motivicAmod}
\newcommand{\pplusqconnectedAmod}{R_{C^{p+q}_{eff}}\motivicAmod}
\newcommand{\qconnectedAAmod}{R_{C^{q}_{eff}}\motivicAAmod}

\newcommand{\qconnectedAmodstablehomotopy}{R_{C^{q}_{eff}}\Amodstablehomotopy}
\newcommand{\Amodqeffsymmcompactgenerators}{\bigcup _{n,r,s\geq 0;s-n\geq q}\; \; \bigcup _{U\in (\smoothS)}A\wedge \symmgeneratorNRS}
\newcommand{\qplusoneconnectedAmod}{R_{C^{q+1}_{eff}}\motivicAmod}

\newcommand{\qminusoneconnectedAmod}{R_{C^{q-1}_{eff}}\motivicAmod}
\newcommand{\qleftAmodlocmaps}{L^{m}(<q)}
\newcommand{\qplusoneleftAmodlocmaps}{L^{m}(<q+1)}
\newcommand{\qleftAAmodlocmaps}{L^{m'}(<q)}
\newcommand{\AmodkillNRS}{id\wedge V(\iota ^{U}_{n,r,s}):A\wedge \symmgeneratorNRS \rightarrow \\ A\wedge \symmdiskNRRS}
\newcommand{\weightqAmod}{L_{<q}\motivicAmod}
\newcommand{\weightqAAmod}{L_{<q}\motivicAAmod}
\newcommand{\weightqplusoneAmod}{L_{<q+1}\motivicAmod}
\newcommand{\weightpplusoneAmod}{L_{<p+1}\motivicAmod}
\newcommand{\weightqplusoneAAmod}{L_{<q+1}\motivicAAmod}

\newcommand{\weightoneAmod}{L_{<1}\motivicAmod}
\newcommand{\weightqAmodstablehomotopy}{L_{<q}\Amodstablehomotopy}
\newcommand{\weightqplusoneAmodstablehomotopy}{L_{<q+1}\Amodstablehomotopy}
\newcommand{\Amodqslicegenerators}{S^{m}(q)}
\newcommand{\qsliceAmod}{S^{q}\motivicAmod}
\newcommand{\psliceAmod}{S^{p}\motivicAmod}
\newcommand{\pplusqsliceAmod}{S^{p+q}\motivicAmod}
\newcommand{\zerosliceAmod}{S^{0}\motivicAmod}
\newcommand{\qsliceAAmod}{S^{q}\motivicAAmod}

\newcommand{\qsliceAmodstablehomotopy}{S^{q}\Amodstablehomotopy}
\newcommand{\weightqAmodinduced}{\widetilde{L_{<q}}\motivicAmod}

\newcommand{\motivicAalg}{A\text{-}\algebras (\pointedmotivic)}

\newcommand{\stableHZmodules}{\mathcal{SH}({H\mathbb Z \text{-}\modules})}

\numberwithin{section}{chapter}

\begin{document}

\frontmatter
\title{Multiplicative Properties of the Slice Filtration}
\author{Pablo Pelaez}
\maketitle
\chapter*{Introduction}

	Let $S$ be a Noetherian separated scheme of finite Krull dimension, and let
	$\stablehomotopy$ denote the motivic stable homotopy category of Morel
	and Voevodsky.  In order to get a motivic version of the Postnikov tower, Voevodsky \cite{MR1977582} constructs
	a filtered family of triangulated subcategories of $\stablehomotopy$:
		\begin{equation}
				\label{equation.introd.slicefiltration}
			\cdots \subseteq \stablehomotopyeffqplusone \subseteq 
			\stablehomotopyeffq \subseteq \stablehomotopyeffqminusone \subseteq \cdots
		\end{equation}
	The work of Neeman \cite{MR1308405}, \cite{MR1812507}, shows that the inclusion:
		$$\xymatrix{i_{q}:\stablehomotopyeffq \ar@{^{(}->}[r]& \stablehomotopy}
		$$
	has a right adjoint $r_{q}$, and that the following functors are exact:
		$$\xymatrix{f_{q},s_{q}:\stablehomotopy \ar[r]& \stablehomotopy}
		$$
	where $f_{q}=i_{q}r_{q}$, and for every $T$-spectrum $X$, $s_{q}(X)$ fits in
	the following distinguished triangle:
		$$\xymatrix{f_{q+1}X\ar[r]& f_{q}X \ar[r]& s_{q}X \ar[r]& \Sigma _{T}^{1,0}f_{q+1}X}
		$$
	We say that $f_{q}(X)$ is the ($q-1$)-connective cover of $X$ and that $s_{q}(X)$ is the
	$q$-slice of $X$.
	
	Let $\motivicsymmTspectra$ be Jardine's model category of symmetric $T$-spectra 
	\cite{MR1787949} and $\symmstablehomotopy$ its associated homotopy category.  
	The symmetrization functor induces an equivalence of categories between $\stablehomotopy$ and $\symmstablehomotopy$
	(see \cite[theorem 4.31]{MR1787949}), we will denote
	by $f_{q}^{\Sigma}$, $s_{q}^{\Sigma}$ the functors on $\symmstablehomotopy$
	that correspond to $f_{q}$, $s_{q}$.
	
	We consider 
	a cofibrant ring spectrum $A$ with unit in 
	$\motivicsymmTspectra$, such that the natural map $f_{0}^{\Sigma}(A)\rightarrow A$
	is an isomorphism in $\symmstablehomotopy$.  Let 
	$Y, Y', Z, Z'$ be arbitrary symmetric
	$T$-spectra, and $p, p', q, q' \in \mathbb Z$.
	Our main results (see theorems \ref{thm.3.5.colocal===>covers-slices-inherit-module-structures},
	\ref{thm.3.5.change-of-coefficients}, \ref{thm.3.5.q-slices=====HZ-modules}, \ref{thm.3.5.q-slices=====big-motives},
	\ref{thm.3.5.external-pairings}, \ref{thm.3.5.qslices-modules-over-zeroslice}, 
	and \ref{thm.3.5.pairings-Atiyah-Hirzebruch-ss}) are the following:
		\begin{enumerate}
			\item	For every $A$-module $M$ in $\motivicsymmTspectra$, its $q$-slice
						$s_{q}^{\Sigma}(M)$ 
						is again an $A$-module in $\motivicsymmTspectra$. 
						(see lemma 
						\ref{lemma.3.6.connecting-functors-with-model-structures}(\ref{lemma.3.6.connecting-functors-with-model-structures.a}) 
						and theorem \ref{thm.3.5.colocal===>covers-slices-inherit-module-structures}).  This means
						that the $q$-slice of every $A$-module inherits an $A$-module structure 
						which is defined not just up to homotopy,
						but in a very strict sense. 
			\item	If the unit map $u:\symmspherespectrum \rightarrow A$
						becomes an isomorphism after applying the zero slice functor $s_{0}^{\Sigma}$
						in $\symmstablehomotopy$,
						then $s_{q}^{\Sigma}(Y)$ has a natural structure of
						$A$-module in $\motivicsymmTspectra$
						(see lemma
						\ref{lemma.3.6.connecting-functors-with-model-structures}(\ref{lemma.3.6.connecting-functors-with-model-structures.d}) 
						and theorem \ref{thm.3.5.change-of-coefficients}), i.e. the $q$-slice of every symmetric $T$-spectrum
						is enriched with an $A$-module structure which is defined not just up to homotopy, but in
						a very strict sense.
			\item	As a consequence, if the base scheme $S$ is a perfect field,
						we are able to prove a conjecture of M. Levine (see \cite[corollary 11.1.3]{MR2365658}),
						which says that for every symmetric $T$-spectrum $X$, 
						its $q$-slice $s_{q}^{\Sigma}(X)$ is naturally equipped with a module structure in $\motivicsymmTspectra$
						over the motivic Eilenberg-MacLane spectrum $H\mathbb Z$ (see theorem \ref{thm.3.5.q-slices=====HZ-modules}).
						Restricting the field even further to the case of characteristic zero, we get that
						all the slices $s_{q}^{\Sigma}(X)$ are big motives in the sense of Voevodsky
						(see theorem \ref{thm.3.5.q-slices=====big-motives}).
			\item	The smash product of symmetric $T$-spectra
						induces the following natural external pairing in $\symmstablehomotopy$
						(i.e. up to homotopy):
							$$\xymatrix{s_{p}^{\Sigma}(Y)\wedge s_{q}^{\Sigma}(Z)\ar[r]^-{\cup ^{s}_{p,q}}& s_{p+q}^{\Sigma}(Y\wedge Z)}
							$$
			\item	As a consequence, if $B$ is a ring in $\symmstablehomotopy$ (i.e. up to homotopy)
						and $N$ is a $B$-module
						in $\symmstablehomotopy$ (i.e. also up to homotopy) then:
							\begin{enumerate} 
								\item	The zero slice of $B$, $s_{0}^{\Sigma}(B)$ is a ring spectrum in $\symmstablehomotopy$
											(i.e. up to homotopy).
								\item	The $q$-slice of $N$,
											$s_{q}^{\Sigma}(N)$ is a module in $\symmstablehomotopy$ (i.e up to homotopy)
											over $s_{0}^{\Sigma}(B)$.
								\item	The direct sum of all the slices of $B$, $s^{\Sigma}(B)=\oplus _{n\in \mathbb Z}s_{n}^{\Sigma}(B)$
											is a graded ring spectrum in $\symmstablehomotopy$
											(i.e. up to homotopy).
								\item	The direct sum of all the slices of $N$, $s^{\Sigma}(N)=\oplus _{n\in \mathbb Z}s_{n}^{\Sigma}(N)$
											is a graded module in $\symmstablehomotopy$ (i.e. up to homotopy)
											over $s^{\Sigma}(B)$.
							\end{enumerate}
			\item	The smash product of symmetric $T$-spectra induces natural external pairings in the motivic
						Atiyah-Hirzebruch spectral sequence generated by the slice filtration 
						(see definition \ref{def.3.5.Atiyah-Hirzebruch-spectral-sequence}):
							$$\xymatrix{E_{r}^{p,q}(Y;Z)\otimes E_{r}^{p',q'}(Y';Z') \ar[r]& E_{r}^{p+p',q+q'}(Y\wedge Y';Z\wedge Z')}
							$$
		\end{enumerate}
			
	To prove the results mentioned above, we need to carry out a very detailed analysis of the
	multiplicative properties (with respect to the smash product of spectra) of the filtration
	(\ref{equation.introd.slicefiltration}) considered above.  It turns out that the natural framework to do this, is provided by
	Jardine's category of motivic symmetric $T$-spectra \cite{MR1787949}.
	Our approach consists basically of three steps:
		\begin{enumerate}
			\item First, we lift Voevodsky's slice filtration to the usual category of $T$-spectra equipped
						with the Morel-Voevodsky motivic stable model structure (see section \ref{section.3.2.modstrslicefilt}).
			\item	Then, using the Quillen equivalence given by the symmetrization and forgetful functors
						\cite{MR1787949}, we are able to
						promote the previous lifting to the category of symmetric $T$-spectra (see section \ref{section.3.3.symmmodstrslicefilt}).
			\item	Finally, we describe the multiplicative properties of the slice filtration using
						the symmetric monoidal structure given by the smash product of symmetric $T$-spectra
						(see sections \ref{section.3.4.multiplicativeproperties-slicefiltration}
						and \ref{section.3.5furthermultiplicativeproperties}).
		\end{enumerate}	
	
	We use Hirschhorn's approach to localization of model categories for	
	the construction of the lifting of the slice filtration to the model category setting.
	In order to apply Hirschhorn's techniques, it is necessary to show that the Morel-Voevodsky motivic
	stable model structure is cellular; for this we rely on Hovey's general approach to spectra
	\cite{MR1860878} and on an unpublished result of Hirschhorn (see theorem \ref{thm2.2.injmod=>cellular}).
	For the description of the multiplicative properties of the slice filtration in the model category
	setting, we use Hovey's results on symmetric monoidal model categories \cite{MR1650134}*{chapter 4}.
	
	We now give an outline of this thesis.  In chapter \ref{chapter1}, we just recall some standard results
	about Quillen model categories.  The reader who is familiar with the terminology of model categories
	may skip this chapter.
	
	In chapter \ref{chapter-mot-hom-theory}, we review the definitions of the Morel-Voevodsky stable model structure
	for simplicial presheaves and Jardine's stable model structure for symmetric $T$-spectra.  We also show that
	these two model structures are cellular, therefore it is possible to apply Hirschhorn's technology to construct
	Bousfield localizations.
	In section \ref{section.2.8Modules-Algebras}
	we recall the construction of the model structures for the categories
	of $A$-modules and $A$-algebras,
	where $A$ denotes a cofibrant ring spectrum 
	with unit in Jardine's motivic symmetric stable model category.
	We verify that the category of $A$-modules equipped with this model structure also
	satisfies Hirschhorn's cellularity condition.
	The reader who is familiar with these model structures may either skip this chapter
	or simply look at
	sections \ref{section.cellularity-injective}, 
	\ref{section-cellularity-motivic-stable}, \ref{section-cellularity-motivic-symmetric-stable} and 
	\ref{section.2.8Modules-Algebras} where we prove that the cellularity condition holds.
	
	Finally in chapter \ref{chap-slice-filtration}, we carry out the program sketched above.
	In section \ref{section-slice-filtration}, we review Voevodsky's construction for the slice
	filtration in the setting of simplicial presheaves.  In section \ref{section.3.2.modstrslicefilt},
	we apply Hirschhorn's localizations techniques to the Morel-Voevodsky stable model structure 
	in order to construct
	three families of model structures, namely
	$\qconnectedTspectra$, $\weightqTspectra$ and $\qsliceTspectra$ ($q\in \mathbb Z$).  The first
	family, $\qconnectedTspectra$ is constructed by a right Bousfield localization
	with respect to the Morel-Voevodsky stable model structure
	(see theorem \ref{thm.3.2.connective-model-structure}), and it provides a lifting
	of Voevodsky's slice filtration to the model category level 
	(see theorem \ref{thm.3.2.liftingslicefilt-modelcats}).  Moreover, this family has the property
	that the cofibrant replacement functor $C_{q}$ provides an alternative description
	for the functor $f_{q}$ (($q-1$)-connective cover) defined above (see theorem \ref{thm.3.2.Rq-models-f<q}).
	In order to get a lifting for the slice functors $s_{q}$ to the model category level,
	we need to introduce the model structures $\weightqTspectra$ and $\qsliceTspectra$.
	The model category $\weightqTspectra$ is defined as a left Bousfield localization with respect to the Morel-Voevodsky
	stable model structure
	(see theorem \ref{thm.3.2.Lqmodelstructures}); its main property is that its fibrant replacement functor $W_{q}$ gives
	an alternative description for the cone of the natural map $f_{q}X\rightarrow X$
	(see theorems \ref{thm.3.1.motivictower} and \ref{thm.3.2.Lq-models-s<q}).  On the other hand, the model structure
	$\qsliceTspectra$ is constructed using right Bousfield localization with respect to the model category
	$\weightqplusoneTspectra$ (see theorem \ref{thm.3.2.Sq-modelcats}), and it gives the desired lifting 
	for the slice functor $s_{q}$ to the
	model category level (see theorem \ref{thm.3.2.Sq-models-sq}).

	In section \ref{section.3.3.symmmodstrslicefilt}, we promote the model structures defined above
	(section \ref{section.3.2.modstrslicefilt}) to the setting of symmetric $T$-spectra.
	In this case, Hirschhorn's localization technology applied to Jardine's stable model structure
	for symmetric $T$-spectra allows us to introduce three families of model structures
	which we denote by $\qconnectedsymmTspectra$,
	$\weightqsymmTspectra$ and $\qslicesymmTspectra$; where the underlying category
	is given by symmetric $T$-spectra
	(see theorems \ref{thm.3.3.symmetricconnective-model-structure}, 
	\ref{thm.3.3.Lqsymmetricmodelstructures} and \ref{thm.3.3.symmetricSq-modelcats}).  Using the
	Quillen equivalence \cite{MR1787949}  given by the symmetrization and the forgetful functors, 
	we are then able to show
	that these new families of model structures are also Quillen equivalent to the ones introduced
	in section \ref{section.3.2.modstrslicefilt}
	(see theorems \ref{thm.3.3.symmetrization-qconnected-Quillen-equiv}, \ref{thm.3.3.symmetrization-Quillen-equivalence-weight<q} 
	and \ref{thm.3.3.symmetrization-qslice-Quillen-equiv}).  
	Therefore, these model structures give liftings for the functors $f_{q}$ and $s_{q}$
	to the model category level
	(see corollary \ref{cor.3.3.symmetric-fq,s<q,sq}, and theorems \ref{thm.3.3.symmRq-models-symmf<q}(\ref{thm.3.3.symmRq-models-symmf<q.c}),
	\ref{thm.3.3.symmSq-models-symm-sq}(\ref{thm.3.3.symmSq-models-symm-sq.c})), 
	with the great technical advantage that the underlying categories
	are now symmetric monoidal.  Hence, we have a natural framework for the study of the multiplicative properties
	of Voevodsky's slice filtration.
	
	In section \ref{section.3.4.multiplicativeproperties-slicefiltration}, we show that
	the smash product of symmetric $T$-spectra
		$$\xymatrix{\pconnectedsymmTspectra \times \qconnectedsymmTspectra \ar[r]^-{-\wedge -}& \pplusqconnectedsymmTspectra}						
		$$
		
		$$\xymatrix{\pslicesymmTspectra \times \qslicesymmTspectra \ar[r]^-{-\wedge -}& \pplusqslicesymmTspectra}						
		$$
	is in both cases a Quillen bifunctor in the sense of Hovey 
	(see theorems \ref{thm.3.4.smash-Quillenbifunctor-RpxRq----->Rqplusp} 
	and \ref{thm.3.4.smash-Quillenbifunctor-SpxSq----->Sqplusp}).  
	
	In section \ref{section.3.5furthermultiplicativeproperties}, we will promote 
	(using the free $A$-module functor $A\wedge -$) the model structures
	constructed in section \ref{section.3.3.symmmodstrslicefilt} to the category of $A$-modules,
	where $A$ is a cofibrant ring spectrum with unit in $\motivicsymmTspectra$.  We will denote these 
	new model structures by $\qconnectedAmod$, $\weightqAmod$ and $\qsliceAmod$.  These
	new model structures will be used as an analogue of the slice filtration for
	the motivic stable homotopy category of $A$-modules, as well as a tool to describe
	the behavior of the slice functors $s_{q}^{\Sigma}$ when they are restricted to the category
	of $A$-modules.
	We will see that if one imposes some natural additonal conditions on the ring spectrum $A$, then
	the free $A$-module functor $A\wedge -$ induces a strict compatibility
	between the slice filtration in the categories of symmetric $T$-spectra and $A$-modules
	(see theorems \ref{thm.3.5.inheritingmodulestructures1}, \ref{thm.3.5.inheritmodelstructures2},
	\ref{thm.3.5.inheritmodelstructures-slices} and \ref{thm.3.5.sliceinvarianceofcoefficients}).
	
	In section \ref{section.3.5.applications}, we will rely in all the previous results
	to show that
	if we have a cofibrant ring spectrum $A$ with unit in $\motivicsymmTspectra$
	which also satisfies some additional hypothesis, then for every $q\in \mathbb Z$
	and for every $A$-module $M$ in $\motivicsymmTspectra$
	(see theorem \ref{thm.3.5.colocal===>covers-slices-inherit-module-structures}):
		\begin{enumerate}
			\item $f_{q}^{\Sigma}(M)$ is again an $A$-module in $\motivicsymmTspectra$ (not just up to homotopy,
						but in a very strict sense).
			\item $s_{q}^{\Sigma}(M)$ is again an $A$-module in $\motivicsymmTspectra$ (not just up to homotopy,
						but in a very strict sense).
		\end{enumerate}
	Furthermore, if the unit map $u:\symmspherespectrum \rightarrow A$ fulfills some mild conditions,
	then the free $A$-module functor $A\wedge -$ induces for every symmetric
	$T$-spectrum $X$ (see theorem \ref{thm.3.5.change-of-coefficients}), 
	a natural structure of $A$-module in $\motivicsymmTspectra$ (i.e. not just up to homotopy, but in a very strict sense) 
	on its $q$-slice $s_{q}^{\Sigma}(X)$.
	
	Finally, we will be able to prove a conjecture of M. Levine (see \cite[corollary 11.1.3]{MR2365658}),
	which says that if the base scheme $S$ is a perfect field, then for every $q\in \mathbb Z$ and for every symmetric
	$T$-spectrum $X$, its $q$-slice $s_{q}^{\Sigma}(X)$ is naturally equipped with a module structure in $\motivicsymmTspectra$
	over the motivic Eilenberg-MacLane spectrum $H\mathbb Z$ (see theorem \ref{thm.3.5.q-slices=====HZ-modules}).
	If we restrict the field even further, considering a field of characteristic zero, then as a consequence we will prove that
	all the slices $s_{q}^{\Sigma}(X)$ are big motives in the sense of Voevodsky
	(see theorem \ref{thm.3.5.q-slices=====big-motives}).
	
	We will also show that for every $p,q\in \mathbb Z$, the smash product of symmetric $T$-spectra induces
	up to homotopy natural pairings 
	(see theorem \ref{thm.3.5.external-pairings}):
		$$\xymatrix{f_{p}^{\Sigma}(X)\wedge f_{q}^{\Sigma}(Y)\ar[r]^-{\cup ^{c}_{p,q}}& f_{p+q}^{\Sigma}(X\wedge Y)}
		$$
							
		$$\xymatrix{s_{p}^{\Sigma}(X)\wedge s_{q}^{\Sigma}(Y)\ar[r]^-{\cup ^{s}_{p,q}}& s_{p+q}^{\Sigma}(X\wedge Y)}
		$$
	As a consequence, if $A$ is a ring spectrum in $\symmstablehomotopy$ (i.e. up to homotopy)
	and $M$ is an $A$-module in $\symmstablehomotopy$, then
	(see theorem \ref{thm.3.5.qslices-modules-over-zeroslice}):
		\begin{enumerate}
			\item  The ($-1$)-connective cover of $A$, $f_{0}^{\Sigma}(A)$ is a ring spectrum (up to homotopy)
							in $\symmstablehomotopy$.
			\item	 For every $q\in \mathbb Z$, the ($q-1$)-connective cover of $M$, $f_{q}^{\Sigma}(M)$ is a module 
							(up to homotopy) over $f_{0}^{\Sigma}(A)$.
			\item	 The direct sum of all the connective covers of $A$, $f^{\Sigma}(A)=\oplus _{n\in \mathbb Z}f_{n}^{\Sigma}(A)$
						 is a graded ring (up to homotopy) in $\symmstablehomotopy$.
			\item	 The direct sum of all the connective covers of $M$, $f^{\Sigma}(M)=\oplus _{n\in \mathbb Z}f_{n}^{\Sigma}(M)$
						 is a graded module (up to homotopy) over $f^{\Sigma}(A)$.
			\item  The zero slice of $A$, $s_{0}^{\Sigma}(A)$ is a ring spectrum (up to homotopy) in $\symmstablehomotopy$.
			\item	 For every $q\in \mathbb Z$, the $q$-slice of $M$, $s_{q}^{\Sigma}(M)$ is a module 
							(up to homotopy) over $s_{0}^{\Sigma}(A)$.
			\item  The direct sum of all the slices of $A$, $s^{\Sigma}(A)=\oplus _{n\in \mathbb Z}s_{n}^{\Sigma}(A)$
					  	is a graded ring (up to homotopy) in $\symmstablehomotopy$.
			\item	The direct sum of all the slices of $M$, $s^{\Sigma}(M)=\oplus _{n\in \mathbb Z}s_{n}^{\Sigma}(M)$
						is a graded module (up to homotopy) over $s^{\Sigma}(A)$.
		\end{enumerate}
	We will also see that the smash product of symmetric $T$-spectra induces 
	(via the external pairings $\cup ^{c}$ and $\cup ^{s}$)
	natural external pairings in the motivic Atiyah-Hirzebruch spectral sequence
	(see definition \ref{def.3.5.Atiyah-Hirzebruch-spectral-sequence} and 
	theorem \ref{thm.3.5.pairings-Atiyah-Hirzebruch-ss}):
		$$\xymatrix@R=0.5pt{E_{r}^{p,q}(Y;X)\otimes E_{r}^{p',q'}(Y';X') \ar[r]& E_{r}^{p+p',q+q'}(Y\wedge Y';X\wedge X')\\
									     (\alpha, \beta) \ar@{|->}[r]& \alpha \smile \beta }
		$$
		
\section*{Acknowledgements}

	I would like to thank my advisor Marc Levine for his generosity in sharing
	his infinite knowledge, and for his unconditional support.  I would also like to
	thank the Mathematics Department at the University of Essen, especially
	professors H\'{e}l\`{e}ne Esnault and Eckart Viehweg
	for their hospitality and their financial support during the year in which 
	this thesis was prepared.  
	
	I would like to thank the Mathematics Department
	at Northeastern University for the financial support 
	during my graduate studies.  Finally, I would like to thank
	my friends Shih-Wei Yang and Shouxin Dai for all their help and encouragement.
\tableofcontents

\mainmatter
\begin{chapter}{Preliminaries}
				\label{chapter1}
				
	All the results in this chapter are classical (see \cite{MR0223432}, \cite{MR1650134}, \cite{MR1944041},
	\cite{MR1711612}, \cite{MR1650938})
	and are included here just to fix notation and to make this note
	self contained.
	
\begin{section}{Model Categories}
			\label{section1}
			
	Model categories were first introduced by Quillen in \cite{MR0223432}, his
	original definition has been slightly modified along the years, we will use the
	definition introduced in \cite{MR2102294}.
	
\begin{defi}
			\label{def1.1.1}
	A model category $\modelcat$ is a category equipped with three classes
	of maps $(\mathcal W ,\mathcal C ,\mathcal F )$ called weak equivalences,
	cofibrations and fibrations, such that the following axioms hold:
	
	\begin{description}
		\item[MC1]	$\modelcat$ is closed under small limits and colimits.
		\item[MC2]	If $f,g$ are two composable maps in $\modelcat$ and two out
								of $f, g, g\circ f$ are weak equivalences then so is the
								third one.
		\item[MC3]	The classes of weak equivalences, cofibrations and fibrations
								are closed under retracts.
		\item[MC4]	Suppose we have a solid commutative diagram:
								$$\xymatrix{A \ar[r] \ar[d]_{i}& X\ar[d]^{p}\\
														B \ar[r] \ar@{-->}[ur]& Y}$$
								where $i$ is a cofibration, $p$ is a fibration, and either
								$i$ or $p$ is a weak equivalence, then the dotted arrow
								making the diagram commutative exists.
		\item[MC5]	Given any arrow $f:A\rightarrow B$ in $\modelcat$, there exist
								two functorial factorizations, $f=p\circ i$ and $f=q\circ j$,
								where $p$ is a fibration and a weak equivalence, $i$ is a
								cofibration, $q$ is a fibration and $j$ is a cofibration and a weak
								equivalence.
	\end{description}
\end{defi}

	A map $j:A\rightarrow B$ will be called a \emph{trivial cofibration}
	(respectively \emph{trivial fibration})
	if it is both a cofibration and a weak equivalence (respectively
	a fibration and a weak equivalence).

	If a given category $\modelcat$ has a model structure, then
	we get immediately the following consequences:
	
\begin{rmks}
		\label{rmk1.1.triv-conseq}
	\begin{enumerate}
		\item	The limit axiom $\mathbf{MC1}$ implies that
					there is an \emph{initial} and a \emph{final}
					object in $\modelcat$, which we will denote by
					$\emptyset$ and $*$ respectively.  We say that
					the category $\modelcat$ is \emph{pointed}
					if the canonical map $\emptyset \rightarrow *$
					is an isomorphism.
		\item	The axioms for a model category are self dual,
					therefore the opposite category $\modelcat ^{op}$
					has also a model structure, where a map
					$i:A\rightarrow B$ in $\modelcat ^{op}$
					is a weak equivalence, cofibration or fibration
					if its dual $i:B\rightarrow A$ is a weak
					equivalence, fibration or cofibration
					in $\modelcat$.  
					This implies in particular
					that any result we prove about model categories
					will have a dual version.
		\item	Let $X$ be an object in $\modelcat$.  Then the category 
					$(\modelcat \downarrow X)$ of objects in $\modelcat$ over $X$
					has also a model structure,
					where the weak equivalences, cofibrations and fibrations
					are maps which become weak equivalences, cofibrations and fibrations
					after applying the forgetful functor
					$(\modelcat \downarrow X) \rightarrow \modelcat$.
		\item	Similarly the category $(X\downarrow \modelcat)$ has a model structure
					induced from the one in $\modelcat$.
					We will denote by $\modelcat _{*}$ the category $(*\downarrow \modelcat)$
					of objects under the final object of $\modelcat$.
		\item	Let $A,X$ be two objects in $\modelcat$, then the
					category $(A\downarrow \modelcat \downarrow X)$ 
					of objects which are simultaneously under $A$ and over $X$ has also
					a model structure induced from the one in $\modelcat$.
	\end{enumerate}
\end{rmks}

	Let $X$ be an object in $\modelcat$.  We say that $X$ is
	\emph{cofibrant} if the natural map $\emptyset \rightarrow X$
	is a cofibration.  Similarly, we say that $X$ is \emph{fibrant}
	if the natural map $X\rightarrow \ast$ is a fibration.

	Consider two objects $A,X$ in $\modelcat$.  We say
	that $A$ is a \emph{cofibrant replacement} for $X$, if
	$A$ is cofibrant and there is a map $A\rightarrow X$
	which is a weak equivalence in $\modelcat$.  Dually, we say that
	$X$ is a \emph{fibrant replacement} for $A$, if $X$
	is fibrant and there is a map $A\rightarrow X$ which
	is a weak equivalence in $\modelcat$.

	Let $i:A\rightarrow B$, $p:X\rightarrow Y$ be two maps in $\modelcat$.
	We say that $i$ has the \emph{left lifting property with respect to}
	$p$ (or that $p$ has the \emph{right lifting property with respect to} $i$)
	if for every solid commutative diagram:
	$$\xymatrix{A \ar[r] \ar[d]_-{i}& X \ar[d]^-{p}\\
							B \ar[r] \ar@{-->}[ur]^-{s}& Y}
	$$
	the dotted arrow making the diagram commutative exists.

	The following are two elementary but extremely useful results about
	model categories.

\begin{prop}[Retract Argument, \cite{MR1650134}]
		\label{prop1.1.retract}
		Let $\modelcat$ be a model category and $f=p\circ i$ a factorization
		of $f$ such that $f$ has the left lifting property with respect to $p$
		(respectively $f$ has the right lifting property with respect to $i$).
		Then $f$ is a retract of $i$ (respectively $f$ is a retract of $p$).
\end{prop}
\begin{proof}
	By duality it is enough to show the case where $f$ has the left lifting property
	with respect to $p$.
	
	Consider the following solid commutative diagram:
	$$\xymatrix{A \ar[r]^{i} \ar[d]_{f}& X\ar[d]^{p}\\
							B \ar[r]_{id} \ar@{-->}[ur]^{j}& B}$$
	By hypothesis the dotted arrow $j$ making the diagram commutative exists.
	But then the following commutative diagram shows that $f$ is
	a retract of $i$.
	$$\xymatrix{A \ar[r]^{id} \ar[d]_{f}& A \ar[r]^{id} \ar[d]^{i}& A \ar[d]^{f}\\
							B \ar[r]_{j}& X \ar[r]_{p}& B}$$
\end{proof}

\begin{lem}[Ken Brown's lemma, \cite{MR1650134}]
		\label{lem1.1.KenBrown}
	Let $F:\modelcat \rightarrow \mathcal D$ be a functor, where $\modelcat$ is a model category.
	Assume that there exists a class $\mathcal V$ of maps in $\mathcal D$ which has the two out of three
	property, and that $F(i)\in \mathcal V$ for all trivial cofibrations $i:A\rightarrow B$ between
	cofibrant objects $A$ and $B$ in $\modelcat$.
	Then $F(g)\in \mathcal V$ for all weak equivalences $g:A\rightarrow B$ between cofibrant
	objects $A$ and $B$ in $\modelcat$.
\end{lem}
\begin{proof}
	Consider the following commutative diagram:
	$$\xymatrix{\emptyset \ar[r] \ar[d]& A \ar[d]_{i_{A}} \ar@/^2pc/[dddrr]^{g}& & \\
							B \ar[r]^-{i_{B}} \ar@/_2.5pc/[ddrrr]_{id}& A\coprod B \ar[dr]^{i} \ar@/_1pc/[ddrr]_{(g,id)}& &\\
							  &            & C\ar[dr]^{p}&\\
							  &&&B}$$
	where we have a factorization of $(g,id)=p\circ i$, with $i$ a cofibration and $p$  a trivial fibration. 
	
	Since $A$ and $B$ are cofibrant, it follows that $i_{A}$ and $i_{B}$ are cofibrations.  This implies
	that $i\circ i_{A}$, $i\circ i_{B}$ are both cofibrations, and hence $C$ is a cofibrant
	object in $\modelcat$.
	
	By the two
	out of three property in $\modelcat$, $i\circ i_{A}$  and $i\circ i_{B}$ are  weak equivalences, 
	since $g$, $p$ and $id_{B}$
	are weak equivalences.  Therefore $i\circ i_{A}$ and $i\circ i_{B}$ are
	both trivial cofibrations.
	It follows that $F(i\circ i_{A})$ and $F(i\circ i_{B})$ are both in
	$\mathcal V$.  But then $F(p)\circ F(i\circ i_{B})=F(p\circ i \circ i_{B})=F(id)=id$,
	and since $\mathcal V$ has the two out of three property, we have that
	$F(p)$ is in $\mathcal V$.
	Then the two out of three property for $\mathcal V$ implies that
	$F(g)=F(p)\circ F(i \circ i_{A})$ is also in $\mathcal V$.
\end{proof}

	By duality we get immediately the following lemma:
	
\begin{lem}
		\label{lemma1.1.KenBrown2}
	Let $F:\modelcat \rightarrow \mathcal D$ be a functor, where
	$\modelcat$ is a model category.  Assume that there exists a class
	$\mathcal V$ of maps in $\mathcal D$ which has the two out of three
	property, and that $F(p)\in \mathcal V$ for all trivial fibrations
	$p:X\rightarrow Y$ between fibrant objects $X$, $Y$ in $\modelcat$.
	Then $F(g)\in \mathcal V$ for all weak equivalences $g:X\rightarrow Y$
	between fibrant objects $X$ and $Y$ 
	in $\modelcat$.  \flushright{$\square$}
\end{lem}

	The retract argument has the following consequences, which give
	nice characterizations for the cofibrations and trivial cofibrations
	(respectively fibrations and trivial fibrations)
	in terms of a left lifting property (respectively right lifting property).

\begin{cor}
		\label{cor.1.1.class-cof-fibs}
	The class of cofibrations (respectively trivial cofibrations)
	in a model category $\modelcat$ is equal to the class of
	maps having the left lifting property with respect to
	any trivial fibration in $\modelcat$ (respectively
	any fibration in $\modelcat$).
	The class of fibrations (respectively trivial fibrations)
	in a model category $\modelcat$ is equal to the class
	of maps having the right lifting property with
	respect to any trivial cofibration in $\modelcat$
	(respectively any cofibration in $\modelcat$).
\end{cor}
\begin{proof}
	By duality it is enough to prove the case of
	cofibrations and trivial cofibrations.
	Suppose that $i:A\rightarrow B$ is a cofibration in $\modelcat$,
	then the lifting axiom $\mathbf{MC4}$ implies that $i$
	has the left lifting property with respect to any trivial
	fibration in $\modelcat$.
	Conversely, if $i:A\rightarrow B$ has the left lifting property
	with respect to any trivial fibration in $\modelcat$,
	then the factorization axiom $\mathbf{MC5}$ implies that
	$i=ql$ where $l$ is a cofibration in $\modelcat$
	and $q$ is a trivial fibration in $\modelcat$.
	Since $i$ has the left lifting property with respect to $q$,
	the retract argument (see proposition \ref{prop1.1.retract})
	implies that $i$ is a retract of $l$.
	Therefore, the retract axiom $\mathbf{MC3}$ implies
	that $i$ is also a cofibration.
	The case for trivial cofibrations is similar.
\end{proof}

\begin{cor}
		\label{cor1.1.stabilility-cof-fib}
	Any isomorphism in a model category $\modelcat$ is
	a cofibration, a fibration, and a weak equivalence.
	The class of cofibrations and the class of trivial cofibrations in $\modelcat$
	are closed under retracts and pushouts.
	The class of fibrations and the class of trivial fibrations in $\modelcat$
	are closed under retracts and pullbacks.
\end{cor}
\begin{proof}
	Follows immediately from the lifting property characterization (corollary \ref{cor.1.1.class-cof-fibs})
	for cofibrations, trivial cofibrations, fibrations and trivial fibrations.
\end{proof}

\begin{rmk}
	Let $\modelcat$ be a model category.
	Given any object $X$ in $\modelcat$,
	we can apply the factorization axiom $\mathbf{MC5}$
	to the natural map $\emptyset \rightarrow X$ to get
	a cofibrant replacement for $X$:
		$$\xymatrix{\emptyset \ar[r] & QX \ar[r]^-{Q^{X}}& X}
		$$
	where $QX$ is cofibrant and $Q^{X}$ is a trivial fibration.
	We also get fibrant replacements for $X$ when we factor the
	natural map $X \rightarrow \ast$:
		$$\xymatrix{X \ar[r]^-{R^{X}}& RX \ar[r]& \ast}
		$$
	where $RX$ is fibrant and $R^{X}$ is a trivial cofibration.
	The factorization
	axiom $\mathbf{MC5}$ implies also that these two
	constructions are functorial.
\end{rmk}

\begin{defi}
		\label{def.1.1.Quillenfunctors}
	Let $\modelcat$, $\modelcattwo$ be two model categories.
	A functor $F:\modelcat \rightarrow \modelcattwo$ is called
	a \emph{left Quillen functor} if it has a right adjoint
	$G:\modelcattwo \rightarrow \modelcat$, and satisfies
	the following conditions:
	\begin{enumerate}
		\item	If $i$ is a cofibration in $\modelcat$,
					then $F(i)$ is also a cofibration in $\modelcattwo$.
		\item	If $j$ is a trivial cofibration in $\modelcat$,
					then $F(j)$ is also a trivial cofibration
					in $\modelcattwo$.
	\end{enumerate}
	The right adjoint $G$ is called a
	\emph{right Quillen functor},
	and the adjunction 
		$$\xymatrix{(F,G,\varphi):\modelcat \ar[r]& \modelcattwo}$$
	is called a \emph{Quillen adjunction}.
\end{defi}

\begin{defi}
		\label{def.1.1.Quillen-equiv}
	Let $(F,G,\varphi):\modelcat \rightarrow \modelcattwo$ be a Quillen adjunction.
	We say that $F$ is a \emph{left Quillen equivalence} if
	for every cofibrant object $X$ in $\modelcat$ and every
	fibrant object $Y$ in $\modelcattwo$ the following condition holds:
	\begin{itemize}
		\item	A map $f:X\rightarrow GY$ is a weak equivalence in $\modelcat$
					if and  only if its adjoint $f^{\sharp}:FX\rightarrow Y$
					is a weak equivalence in $\modelcattwo$.
	\end{itemize}
	In this case $G$ will be called a \emph{right Quillen equivalence}, and
	$(F,G,\varphi)$ a \emph{Quillen equivalence}.
\end{defi}

\begin{defi}
		\label{def.1.1.cylinderobject}
	Let $\modelcat$ be a model category, and let $X$ be an object of $\modelcat$.
	We say that $\tilde{X}$ is a \emph{cylinder object} for $X$, if we have a
	factorization of the fold map
	$$\xymatrix{X\coprod X \ar[r]^-{\nabla} \ar[d]_-{i} & X\\
							\tilde{X} \ar[ur]_-{s}&}
	$$
	where $i$ is a cofibration and $s$ is a weak equivalence.
\end{defi}

\begin{defi}
		\label{def.1.1.pathobject}
	Let $\modelcat$ be a model category, and let $X$ be an object of $\modelcat$.
	We say that $\hat{X}$ is a \emph{path object} for $X$, if we have a
	factorization of the diagonal map
	$$\xymatrix{X \ar[r]^-{\Delta} \ar[d]_-{r} & X\times X\\
							\hat{X} \ar[ur]_-{p}&}
	$$
	where $p$ is a fibration and $r$ is a weak equivalence.
\end{defi}

\begin{defi}
		\label{def.1.1.lefthomotopy}
	Let $\modelcat$ be a model category and consider two
	maps $f,g:X\rightarrow Y$.  We say that
	$f$ is \emph{left homotopic} to $g$ ($f\stackrel{l}{\sim}g$) if there
	exists a cylinder object $\tilde{X}$ for $X$, together with the following factorization:
	$$\xymatrix{X\coprod X \ar[r]^-{(f,g)} \ar[d]_-{i}& Y\\
							\tilde{X} \ar[ur]_-{H}&}
	$$
	The map $H$ is called a \emph{left homotopy} from $f$
	to $g$.
\end{defi}

\begin{defi}
		\label{def.1.1.righthomotopy}
	Let $\modelcat$ be a model category and consider two
	maps $f,g:X\rightarrow Y$.  We say that
	$f$ is \emph{right homotopic} to $g$ ($f\stackrel{r}{\sim}$g) if there
	exists a path object $\hat{Y}$ for $Y$,	together with the following factorization:
		$$\xymatrix{X \ar[r]^-{(f,g)} \ar[d]_-{H}& Y\times Y\\
							\hat{Y} \ar[ur]_-{p}&}
		$$
	The map $H$ is called a \emph{right homotopy} from $f$
	to $g$.
\end{defi}

\begin{defi}
		\label{def.1.1.homotopymaps}
	Let $\modelcat$ be a model category and consider
	two maps $f,g:A\rightarrow B$.
	We say that $f$ is \emph{homotopic} to $g$ ($f\sim g$) if
	$f$ and $g$ are both left and right homotopic.
\end{defi}

\begin{defi}[cf. \cite{MR0223432}]
		\label{def.1.2.loc-of-cats}
	Let $\modelcat$ be an arbitrary category and
	$\mathcal W$ a class of maps in $\modelcat$.
	The \emph{localization} of
	$\modelcat$ with respect to $\mathcal W$
	will be a category $\mathcal W ^{-1}\modelcat$
	together with a functor 
		$$\xymatrix{\gamma :\modelcat \ar[r]& \mathcal W ^{-1}\modelcat}$$
	having the following universal property:
	for every $w\in \mathcal W$, $\gamma (w)$ is an isomorphism, and
	given any functor $F:\modelcat \rightarrow \mathcal D$ such that
	$F(w)$ is an isomorphism for every $w\in \mathcal W$, there is a unique
	functor $\theta :\mathcal W ^{-1}\modelcat \rightarrow \mathcal D$,
	such that $\theta \circ \gamma =F$, i.e. the following
	diagram commutes:
		$$\xymatrix{\modelcat \ar[r]^-{F} \ar[d]_-{\gamma}& \mathcal D\\
							\mathcal W ^{-1}\modelcat \ar@{-->}[ur]_-{\theta}&}
		$$
\end{defi}

\begin{thm}[Quillen]
		\label{thm.1.1.HomotCat-exist}
	Let $\modelcat$ be a model category.
	Then there exists a category $\hocat$, which is the localization
	of $\modelcat$ with respect to the class $\mathcal W$ of weak equivalences,
	and is called the \emph{homotopy category}
	of $\modelcat$.  $\hocat$ is
	defined as follows:
	\begin{enumerate}
		\item	The objects of $\hocat$ are
					just the objects in $\modelcat$.
		\item	The set of maps in $\hocat$ between two
					objects $X,Y$ is given by
					the set of homotopy classes between
					cofibrant-fibrant replacements for
					$X$ and $Y$:
					$$\Hom_{\hocat}(X,Y)=\pi_{\modelcat}(RQX,RQY)$$
					and the composition law is induced by the composition
					in $\modelcat$.
	\end{enumerate}
	Let $\hocat _{c}$, $\hocat _{f}$, $\hocat _{cf}$ be the full
	subcategories of $\hocat$ generated by the cofibrant, fibrant
	and cofibrant-fibrant objects of $\modelcat$ respectively.
	In the following diagram, all the functors are equivalences of
	categories:
	$$\xymatrix{&& \hocat_{c} \ar@<-1ex>@{_{(}->}[drr]_-{\sim} \ar@<1ex>[lld]^-{\mathrm{Ho}R}&& \\
							 \hocat_{cf} 
							\ar@<1ex>@{^{(}->}[urr]^-{\sim} \ar@<-1ex>@{_{(}->}[drr]_-{\sim} &&&& 
							\hocat \ar@<-1ex>[llu]_-{\mathrm{Ho}Q} \ar@<1ex>[dll]^-{\mathrm{Ho}R}\\
							&& \hocat_{f} \ar@<1ex>@{^{(}->}[urr]^-{\sim} \ar@<-1ex>[ull]_-{\mathrm{Ho}Q}&&}
	$$
	where the adjoints to the equivalences given above are constructed
	taking cofibrant, fibrant and cofibrant-fibrant replacements.
\end{thm}
\begin{proof}
	We refer the reader to \cite[I.1 theorem 1]{MR0223432}.
\end{proof}

\begin{thm}[Quillen]
		\label{thm.1.1.Quillen-adjunction=>descends}
	Let $(F,G,\varphi):\modelcat \rightarrow \modelcattwo$ be a Quillen
	adjunction.  Then the adjunction $(F,G,\varphi)$ descends to the homotopy categories, i.e.
	we get an adjuntion:
		$$\xymatrix{(QF,RG,\varphi):\hocat \ar[r]& \hocattwo}$$
	Furthermore, if $(F,G,\varphi)$ is a Quillen equivalence, then
	$(QF,RG,\varphi)$ is an equivalence of categories.
\end{thm}
\begin{proof}
	We refer the reader to \cite[I.4 theorem 3]{MR0223432}.
\end{proof}

\end{section}

\begin{section}{Cofibrantly Generated Model Categories}

	In this section we recall the definition of a cofibrantly generated model category.
	
	In order to get the functorial factorizations required in axiom $\mathbf{MC5}$,
	we need to introduce ordinals, cardinals, and regular cardinals.
	For a definition of these, see
	\cite[chapter 10]{MR1944041}.
	It will be convenient in some situations to consider an ordinal $\lambda$ as a small category, 
	with objects equal to the elements of $\lambda$, and a unique map from $a$ to $b$ if
	$a\leq b$.
	
\begin{defi}
		\label{def.1.1.lambda-sequence}
	Let $\mathcal C$ be a category that is closed under small colimits, and let
	$\mathcal V$ be a class of maps in $\mathcal C$.
	If $\lambda$ is an ordinal, then a \emph{$\lambda$-sequence} in
	$\mathcal C$ is a functor $A:\lambda \rightarrow \mathcal C$, i.e.
	a diagram
	$$A_{0}\rightarrow A_{1}\rightarrow \cdots \rightarrow A_{\beta}\rightarrow \cdots
		\; (\beta <\lambda)
	$$
	such that for every limit ordinal $\gamma <\lambda$ the induced map
	$$\colimit _{\beta <\gamma}A_{\beta}\rightarrow A_{\gamma}
	$$
	is an isomorphism.
	
	The \emph{composition} of the $\lambda$-sequence is the map
	$A_{0}\rightarrow \colimit _{\beta <\lambda}A_{\beta}$.
	
	If $A_{\beta}\rightarrow A_{\beta +1}$ is in $\mathcal V$ for
	any $\beta <\lambda$, we say that the $\lambda$-sequence is a
	\emph{$\lambda$-sequence of maps in} $\mathcal V$, and the transfinite composition
	$A_{0}\rightarrow \colimit _{\beta <\lambda}A_{\beta}$ is called
	a \emph{transfinite composition} of maps in $\mathcal V$.
\end{defi}

\begin{prop}
		\label{prop1.1.1.cofs-wrt-transcomp}
	Let $\modelcat$ be a model category, then the cofibrations
	and trivial cofibrations in $\modelcat$ are both
	closed under transfinite composition.
\end{prop}
\begin{proof}
	The cofibrations and trivial cofibrations in $\modelcat$
	are characterized by a left lifting property.
	But the universal property of the colimit
	clearly preserves this lifting property under transfinite
	composition.
\end{proof}

\begin{defi}
		\label{def.1.1.1.small-objects}
	Let $\mathcal C$ be a category closed under small colimits, and let
	$\mathcal V$ be a class of maps in $\mathcal C$.
	\begin{enumerate}
		\item	If $\kappa$ is a cardinal, then an object $D$ in $\mathcal C$
					is \emph{$\kappa$-small relative to} $\mathcal V$, if
					for every regular cardinal $\lambda \geq \kappa$ and every
					$\lambda$-sequence
					$$\lambdaseq
					$$
					of maps in $\mathcal V$, we have a bijection of sets:
					$$\colimit _{\beta <\lambda}\Hom _{\mathcal C}(D,A_{\beta})
						\rightarrow \Hom _{\mathcal C}(D,\colimit _{\beta <\lambda}A_{\beta})
					$$
		\item	An object $D$ in $\mathcal C$ is \emph{small relative to}
					$\mathcal V$ if it is $\kappa$-small relative to $\mathcal V$ for some
					cardinal $\kappa$, and it is \emph{small} if it is small
					relative to the class of all maps in $\mathcal C$.
	\end{enumerate}
\end{defi}

\begin{defi}
		\label{def.1.1.1.I-inj--I-cof}
	Let $\mathcal C$ be a category, and let $I$ be a \emph{set}
	of maps in $\mathcal C$.
	\begin{enumerate}
		\item	We define $I$-inj as the class of maps in $\mathcal C$
					that have the right lifting property with respect
					to every map in $I$.
		\item	We define $I$-cof as the class of maps in $\mathcal C$
					that have the left lifting property with respect
					to every map in $I$-inj.
	\end{enumerate}
\end{defi}

\begin{defi}
		\label{def1.2.1.rel-cell-complex}
	Let $\mathcal C$ be a category closed under small colimits, and let $I$
	be a \emph{set} of maps in $\mathcal C$, then
	\begin{enumerate}
		\item	The \emph{relative $I$-cell complexes} are the maps that can be
					constructed as a transfinite composition of pushouts of
					elements of $I$.
		\item	An object $A$ of $\mathcal C$ is an \emph{$I$-cell complex}, if the map
					$\emptyset \rightarrow A$ is a relative $I$-cell complex.
		\item	\label{def1.2.1.rel-cell-complex.a}	A map is an \emph{inclusion of $I$-cell complexes} if it is a relative
					$I$-cell complex whose domain is an $I$-cell complex.
	\end{enumerate}
\end{defi}

	We will denote the class of relative $I$-cell complexes as
	\emph{$I$-cells}.

\begin{rmk}
		\label{remark.1.2.1.Icells=>I-cofs}
	Since the left lifting property is preserved under pushouts
	and transfinite compositions we have that
	$I\text{-cells}\subseteq I\text{-cof}$.
\end{rmk}
	
\begin{thm}[Quillen's small object argument]
		\label{thm.1.1.1.Quillen-small-object-arg}
	Let $\mathcal C$ be a category closed under small colimits, and let $I$ be
	a \emph{set} of maps in $\mathcal C$.
	Assume that the domains of all the maps in $I$ are small
	with respect to  $I$-cells.  Then for every map $f:X\rightarrow Y$
	in $\mathcal C$, there is a functorial factorization
	$$\xymatrix{X \ar[r]^-{i}& E^{f}_{I} \ar[r]^-{p}& Y}
	$$
	where $i$ is in $I$-cells, and
	$p$ is in $I$-inj.
\end{thm}
\begin{proof}
	We refer the reader to \cite{MR0223432}, \cite{MR1944041}, or \cite{MR1650134}.
\end{proof}
	
\begin{defi}
		\label{def1.2.cofgen}
	A model category $\modelcat$ is \emph{cofibrantly generated}
	if there exist \emph{sets} $I$ and $J$ of maps in $\modelcat$,
	such that:
	
	\begin{enumerate}
		\item	The domains of all the maps in $I$ are small
					with respect to the $I$-cells.
		\item	The domains of all the maps in $J$ are small
					with respect to the $J$-cells.
		\item	The class $\mathcal F  \cap \mathcal W$ of trivial fibrations in $\modelcat$ is 
					equal to $I$-inj.
		\item	The class $\mathcal F$ of fibrations in $\modelcat$ is equal
					to $J$-inj.
	\end{enumerate}
	In this situation, $I$ will be called the set of \emph{generating cofibrations},
	and $J$ will be called the set of \emph{generating
	trivial cofibrations}.
\end{defi}

	To work with spectra, we need to start with a pointed model category.
	The following result will allow us to go from an unpointed cofibrantly
	generated model category to a pointed one.

\begin{thm}[Hirschhorn]
		\label{thm1.3.pointedcat-cofgen}
	Let $\modelcat$ be a cofibrantly generated model category
	with set of generating cofibrations $I$ and 
	set of generating trivial cofibrations $J$.
	Then the associated pointed model category $\modelcat _{*}$
	(see remark \ref{rmk1.1.triv-conseq}) is also a cofibrantly
	generated model category, with set of generating
	cofibrations $F(I)=I_{+}$ and set of generating trivial
	cofibrations $F(J)=J_{+}$, where $F$ is the functor
	$F:\modelcat \rightarrow \modelcat _{*}$
	defined on objects $A$ in $\modelcat$ as
	the pushout in the commutative diagram:
	$$\xymatrix{\emptyset \ar[r] \ar[d]& A \ar[d]\\
							\ast \ar[r]& F(A)=A_{+}}
	$$
	and on maps $i:A\rightarrow B$ in $\modelcat$
	as:
	$$\xymatrix{F(A)=A\coprod \ast \ar[rr]^-{F(i)\coprod id} && B\coprod \ast=F(B)}
	$$
\end{thm}
\begin{proof}
	We refer the reader to \cite[theorem 2.7]{Over-Und}.
\end{proof}

\end{section}
\begin{section}{Cellular Model Categories}
		\label{subsection-cell-mod-cats}
		
	In this section we review Hirschhorn's cellularity, which is the main property
	that a  model category has to satisfy if we want to construct Bousfield localizations.
	 
\begin{defi}
		\label{def.1.2.1.presented-rel-cell-complex}
	Let $\mathcal C$ be a category closed under small colimits, and let
	$I$ be a \emph{set} of maps in $\mathcal C$.  If $i:A\rightarrow B$ is
	a relative $I$-cell complex, then a \emph{presentation} of $i$ is a pair
	consisting of a $\lambda$-sequence
	$$\lambdaseq
	$$
	for some ordinal $\lambda$, and a sequence of ordered triples
	$$\{(T^{\beta},e^{\beta},h^{\beta})\}
	$$
	such that:
	\begin{enumerate}
		\item	The composition of the $\lambda$-sequence is isomorphic to $i$
		\item	For every $\beta <\lambda$
					\begin{enumerate}
						\item	$T^{\beta}$ is a set.
						\item	$e^{\beta}$ is a function $e^{\beta}:T^{\beta}\rightarrow I$.
						\item	If $i\in T^{\beta}$ and $e^{\beta}_{i}$ is the element $C_{i}\rightarrow D_{i}$
									of $I$, then $h^{\beta}_{i}$ is a map
									$h^{\beta}_{i}:C_{i}\rightarrow A_{\beta}$, such that there is a pushout diagram
									\begin{equation*}
										\xymatrix{\displaystyle{\coprod _{i\in T^{\beta}}C_{i}} \ar[rr]^-{\coprod e^{\beta}_{i}} \ar[d]_-{\coprod h^{\beta}_{i}}
															&& \displaystyle{\coprod _{i\in T^{\beta}}D_{i}} \ar[d] \\
															A_{\beta} \ar[rr] && A_{\beta +1}}
									\end{equation*}
					\end{enumerate}
	\end{enumerate}
\end{defi}

\begin{defi}
		\label{def1.2.1.size-cell-complex}
	Let $\mathcal C$ be a category closed under small colimits, and let $I$ be a \emph{set}
	of maps in $\mathcal C$.  If
	$$\prescellcplx
	$$
	is a presented relative $I$-cell complex, then
	\begin{enumerate}
		\item	The \emph{presentation ordinal} of $i$ is $\lambda$.
		\item	The \emph{set of cells} of $i$ is $\coprod _{\beta <\lambda}T^{\beta}$.
		\item	The \emph{size} of $i$ is the cardinal of the set of cells of $i$.
		\item	If $e$ is a cell of $i$, the \emph{presentation ordinal} of $e$ is the
					ordinal $\beta$ such that $e\in T^{\beta}$.
		\item	If $\beta <\lambda$, then the \emph{$\beta$-skeleton} of $i$ is $A_{\beta}$.
	\end{enumerate}
\end{defi}

	The next remark follows directly from the previous definitions.
	
\begin{rmk}
		\label{rmk1.2.1.class-pres-cell-complex}
	If $\mathcal C$ is a category closed under small colimits, and $I$ is a \emph{set}
	of maps in $\mathcal C$, then a presented relative $I$-cell complex is entirely
	determined by its presentation ordinal $\lambda$, and its sequence of
	triples $\{(T^{\beta},e^{\beta},h^{\beta})\}_{\beta <\lambda}$.
\end{rmk}

\begin{defi}
		\label{def.1.2.1.subcomplex}
	Let $\mathcal C$ be a category closed under small colimits, and $I$ a \emph{set}
	of maps in $\mathcal C$.  If
	$$\prescellcplx
	$$
	is a presented relative $I$-cell complex, then a \emph{subcomplex} of $i$
	consists of a presented relative $I$-cell complex
	$$\tilde{\imath}:A\rightarrow \tilde{B}, \; A=\tilde{A}_{0}\rightarrow \tilde{A}_{1} \rightarrow \cdots \rightarrow \tilde{A}_{\beta}
		\rightarrow \cdots (\beta <\lambda),\; \{\tilde{T}^{\beta},\tilde{e}^{\beta},\tilde{h}^{\beta}\}_{\beta <\lambda}
	$$
	such that
	\begin{enumerate}
		\item	For every $\beta <\lambda$, $\tilde{T}^{\beta}\subseteq T^{\beta}$ and
					$\tilde{e}^{\beta}$ is the restriction of $e^{\beta}$ to $\tilde{T}^{\beta}$.
		\item	There is a map of $\lambda$-sequences
					$$\xymatrix{A \ar[r]^-{id} \ar[d]_-{id} & \tilde{A}_{0} \ar[r] \ar[d]& \tilde{A}_{1} \ar[r] \ar[d]
											& \tilde{A}_{2} \ar[r] \ar[d] & \cdots \\
											A \ar[r]_-{id}& A_{0} \ar[r] & A_{1} \ar[r]& A_{2} \ar[r]& \cdots}
					$$
					such that, for every $\beta <\lambda$ and every $i\in \tilde{T}^{\beta}$,
					the map $\tilde{h}^{\beta}_{i}:C_{i}\rightarrow \tilde{A}_{\beta}$
					is a factorization of the map $h^{\beta}_{i}:C_{i}\rightarrow A_{\beta}$
					through the map $\tilde{A}_{\beta}\rightarrow A_{\beta}$.
	\end{enumerate}
\end{defi}

\begin{prop}
		\label{prop1.2.1.cell-cplxs-gen-monos}
	Let $\mathcal C$ be a category closed under small colimits, and $I$
	a \emph{set} of maps in $\mathcal C$ such that the relative $I$-cell complexes
	are monomorphisms, then a subcomplex of a presented
	relative $I$-cell complex is entirely determined by its set
	of cells $\{\tilde{T}^{\beta}\}_{\beta <\lambda}$.
\end{prop}
\begin{proof}
	The definition of a subcomplex implies that the maps
	$\tilde{A}_{\beta}\rightarrow A_{\beta}$ are all
	inclusions of subcomplexes (see definition \ref{def1.2.1.rel-cell-complex}(\ref{def1.2.1.rel-cell-complex.a})).  
	Since inclusions of subcomplexes are monomorphisms,
	there is at most one possible factorization $\tilde{h}^{\beta}_{i}$ of each
	$h^{\beta}_{i}$ through $\tilde{A}_{\beta}\rightarrow A_{\beta}$
\end{proof}

\begin{prop}
		\label{prop1.2.1.constr-subcplxs}
	Let $\mathcal C$ be a category closed under small colimits, and let $I$ be a \emph{set}
	of maps in $\mathcal C$ such that the relative $I$-cell complexes are monomorphisms.
	If
	$$\prescellcplx
	$$
	is a presented relative $I$-cell complex, then an arbitrary subcomplex of $i$
	can be constructed by the following inductive procedure:
	\begin{enumerate}
		\item	Choose an arbitrary subset $\tilde{T}^{0}$ of $T^{0}$.
		\item	If $\beta <\lambda$ and we have defined $\{\tilde{T}^{\gamma}\}_{\gamma <\beta}$,
					then we have determined the object $\tilde{A}_{\beta}$ and the map
					$\tilde{A}_{\beta}\rightarrow A_{\beta}$.  Consider the set
					$$\{i\in T^{\beta}|h^{\beta}_{i}:C_{i}\rightarrow A_{\beta}
					\text{ factors through }\tilde{A}_{\beta}\rightarrow A_{\beta}\}
					$$
					Choose an arbitrary subset $\tilde{T}^{\beta}$ of this set.  For every $i\in \tilde{T}^{\beta}$
					there is a unique map $\tilde{h}^{\beta}_{i}:C_{i}\rightarrow \tilde{A}_{\beta}$
					that makes the diagram
					$$\xymatrix{C_{i} \ar[d]_-{\tilde{h}^{\beta}_{i}} \ar[dr]^-{h^{\beta}_{i}}&\\
											\tilde{A}_{\beta} \ar[r]& A_{\beta}}
					$$
					commute.  Let $\tilde{A}_{\beta +1}$ be defined by the pushout diagram
					$$\xymatrix{\displaystyle{\coprod _{i\in \tilde{T}^{\beta}}}C_{i} \ar[r] \ar[d]_-{\coprod \tilde{h}^{\beta}_{i}}& 
											\displaystyle{\coprod _{i\in \tilde{T}^{\beta}}}D_{i} \ar[d] \\
											\tilde{A}_{\beta} \ar[r]& \tilde{A}_{\beta +1}}
					$$
	\end{enumerate}
\end{prop}
\begin{proof}
	Follows immediately from the definitions and proposition \ref{prop1.2.1.cell-cplxs-gen-monos}.
\end{proof}

\begin{cor}
		\label{cor.1.1.2.existence-union-cellcplxs}
	Let $\mathcal C$ be a category closed under small colimits, and let
	$I$ be a \emph{set} of maps in $\mathcal C$ such that the
	relative $I$-cell complexes are monomorphisms.
	Consider an arbitrary
	$$\prescellcplx
	$$
	presented cell complex. 
	Assume that $S$ is a \emph{set} and
	take an arbitrary family $\{A_{s}\}_{s\in S}$ of subcomplexes
	of $i:A\rightarrow B$, then there exists
	a subcomplex $\cup _{s\in S}A_{s}$ which
	represents the union of the given family.	
\end{cor}
\begin{proof}
	Follows immediately from proposition \ref{prop1.2.1.constr-subcplxs}.
\end{proof}

\begin{defi}
		\label{def.1.2.1.relative-compact}
	Let $\mathcal C$ be a category closed under small colimits, and let $I$ be a \emph{set} of maps
	in $\mathcal C$.
	\begin{enumerate}
		\item	If $\gamma$ is a cardinal, then an object $A$ of $\mathcal C$ is
					\emph{$\gamma$-compact relative to} $I$ if, for every presented relative
					$I$-cell complex $i:X\rightarrow Y$, every map
					from $A$ to $Y$ factors through a subcomplex of $i$ of size at most
					$\gamma$.
		\item	An object $A$ of $\mathcal C$ is \emph{compact relative to} $I$ if it
					is $\gamma$-compact relative to $I$ for some cardinal $\gamma$.
	\end{enumerate}
\end{defi}

\begin{defi}
		\label{def.1.1.2.compactness}
	Let $\modelcat$ be a cofibrantly generated model category
	with set of generating cofibrations $I$.
	\begin{enumerate}
		\item	If $\gamma$ is a cardinal, then an object $X$ of $\modelcat$
					is \emph{$\gamma$-compact} if it is $\gamma$-compact relative
					to $I$ (see definition \ref{def.1.2.1.relative-compact}).
		\item	An object $X$ of $\modelcat$ is \emph{compact} if there is a cardinal
					$\gamma$ for which it is $\gamma$-compact.
	\end{enumerate}
\end{defi}		
		
	To complete the definition of a cellular model category, we need
	to introduce the concept of \emph{effective monomorphism}.
	
\begin{defi}
		\label{def.1.1.2.effective-mono}
	Let $\mathcal C$ be a category that is closed under pushouts.
	The map $i:A\rightarrow B$ is an \emph{effective monomorphism}
	if $i$ is the equalizer of the pair of natural inclusions
	$B\rightrightarrows B\coprod _{A}B$.
\end{defi}

\begin{rmk}
		\label{rmk1.1.2.sets-effec=inj}
	In the category of sets, the class of effective monomorphisms
	is just the class of injective maps.
\end{rmk}

\begin{defi}[cf. \cite{MR1944041}]
		\label{def.1.1.2.cell-mod-cats}
	Let $\modelcat$ be a model category.
	We say that $\modelcat$ is \emph{cellular}
	if it satisfies the following conditions:
	\begin{enumerate}
		\item	\label{def.1.1.2.cell-mod-cats.a}$\modelcat$ is cofibrantly generated (see definition \ref{def1.2.cofgen})
					with set of generating cofibrations $I$ and set of generating
					trivial cofibrations $J$.
		\item	\label{def.1.1.2.cell-mod-cats.b}Both the domains and codomains of the maps in $I$ are compact
					(see definition \ref{def.1.1.2.compactness}).
		\item	\label{def.1.1.2.cell-mod-cats.c}The domains of the maps in $J$ are small relative to $I$
					(see definition \ref{def.1.1.1.small-objects}).
		\item	\label{def.1.1.2.cell-mod-cats.d}The cofibrations in $\modelcat$ are effective monomorphisms
					(see definition \ref{def.1.1.2.effective-mono}).
	\end{enumerate}
\end{defi}

	When we have a cellular model category $\modelcat$ with
	set of generating cofibrations $I$, the relative
	$I$-cell complexes will be called relative cell complexes.
	
\begin{thm}[Hirschhorn]
		\label{thm.1.5.pointedmodelcat-cellular}
	Let $\modelcat$ be a cellular model category.
	Then the associated pointed model category
	$\modelcat _{*}$ equipped with the model structure
	considered in theorem \ref{thm1.3.pointedcat-cofgen}
	is also cellular.
\end{thm}
\begin{proof}
	We refer the reader to \cite[theorem 2.8]{Over-Und}.
\end{proof}

\end{section}
\begin{section}{Proper Model Categories}
		\label{subsec-propmodcats}
		
	In this section we just recall the definition of proper model categories.
	
\begin{defi}
		\label{def.leftproper}
	Let $\modelcat$ be a model category.
	We say that $\modelcat$ is \emph{left proper}
	if the class of weak equivalences is closed
	under pushouts along cofibrations, i.e.
	in any pushout diagram
	$$\xymatrix{A \ar[r]^-{h} \ar[d]_{i}& X \ar[d]\\
							B \ar[r]_-{h_{*}}& Y}
	$$
	where $i$ is a cofibration and $h$ is a weak equivalence,
	we then have that $h_{*}$ is also a weak equivalence. 
\end{defi}

\begin{defi}
		\label{def.rightproper}
	Let $\modelcat$ be a model category.
	We say that $\modelcat$ is \emph{right proper}
	if the class of weak equivalences is closed
	under pullbacks along fibrations, i.e.
	in any pullback diagram
	$$\xymatrix{A \ar[r]^-{h^{*}} \ar[d]& X \ar[d]^-{p}\\
							B \ar[r]_-{h}& Y}
	$$
	where $p$ is a fibration and $h$ is a weak equivalence,
	we then have that $h^{*}$ is also a weak equivalence.
\end{defi}

\begin{defi}
		\label{def.proper}
	Let $\modelcat$ be a model category.
	We say that $\modelcat$ is \emph{proper}
	if it is both left and right proper.
\end{defi}

\begin{thm}[Hirschhorn]
		\label{thm1.6.pointedmodelcat-proper}
	Let $\modelcat$ be a left proper, right proper, or proper
	model category.  Then the associated pointed model
	category $\modelcat _{*}$ (see remark \ref{rmk1.1.triv-conseq}) 
	is also left proper, right proper, or proper. 
\end{thm}
\begin{proof}
	We refer the reader to \cite[theorem 2.8]{Over-Und}.
\end{proof}

\end{section}
\begin{section}{Simplicial Sets}

	Let $\Delta$ denote the category of well ordered finite sets, i.e. the category with objets:
			$$\mathbf{n}=\{ 0<1<\cdots < n\}$$
	where $n\geq 0$; and maps the weakly order preserving functions, i.e.:
			$$\Hom _{\Delta}(\mathbf{m},\mathbf{n})=\{f:\mathbf{m}\rightarrow
			\mathbf{n} | i\leq j \Rightarrow f(i)\leq f(j)\}$$
			
	There exists a canonical set of generators for the maps in $\Delta$, called
	\emph{cofaces} ($\delta ^{i}:\mathbf{n}\rightarrow \mathbf{n+1}$), 
	and codegeneracies ($\sigma ^{i}:\mathbf{n+1}\rightarrow \mathbf{n}$), defined as:
			$$\delta ^{i}(j)=
				\begin{cases}
					j, &\text{if } j<i\\
					j+1, &\text{if } j\geq i 
				\end{cases}$$
			$$\sigma ^{i}(j)=
				\begin{cases}
					j, &\text{if } j\leq i\\
					j-1, &\text{if } j>i
				\end{cases}$$
				
	The cofaces and degeneracies satisfy a list of relations called the
	\emph{cosimplicial identities}:	
			\begin{eqnarray}
					\label{cosimplicial-identities}
				\begin{array}{lll}
					\delta ^{j}\delta ^{i} &=\delta ^{i} \delta ^{j-1} &\text{for } i<j\\
					\sigma ^{j}\delta ^{i} &=\delta ^{i}\sigma ^{j-1} &\text{for } i<j\\
					\sigma ^{i}\delta ^{i} &= id &\\
					\sigma ^{i}\delta ^{i+1} &= id &\\
					\sigma ^{j}\delta ^{i} &= \delta ^{i-1}\sigma ^{j} &\text{for } i> j+1\\
					\sigma ^{j}\sigma ^{i} &= \sigma ^{i}\sigma ^{j+1} &\text{for } i<j\\
				\end{array}
			\end{eqnarray}
			
\begin{defi}
	A simplicial set $X$ is a contravariant functor from the category
	$\Delta$ to the category of sets.
\end{defi}

	We will denote the category of simplicial sets by $\mathbf{SSets}$, where the
	maps between to simplicial sets $X$ and $Y$ are just the natural transformations
	$\eta :X\rightarrow Y$.

	It follows from the cosimplicial identities that to specify a simplicial set $X$,
	it is enough to give sets $X_{0}, X_{1}, \ldots , X_{n}, \ldots$;
	where $X_{i}=X(i)$ together with \emph{face} maps $d_{i}:X_{n}\rightarrow X_{n-1}$
	($d_{i}=X(\delta ^{i})$) and \emph{degeneracy} maps $s_{i}:X_{n}\rightarrow X_{n+1}$
	($s_{i}=X(\sigma ^{i})$), satisfying the following relations which are called
	\emph{simplicial identities} (these are just the duals with respect to the 
	cosimplicial identities):
			\begin{eqnarray}
					\label{simplicial-identities}
				\begin{array}{lll}
					d_{i}d_{j} &=d_{j-1}d_{i} &\text{for } i<j\\
					d_{i}s_{j} &=s_{j-1}d_{i} &\text{for } i<j\\
					d_{i}s_{i}&= id &\\
					d_{i+1}s_{i}&= id &\\
					d_{i}s_{j}&= s_{j}d_{i-1} &\text{for } i> j+1\\
					s_{i}s_{j}&= s_{j+1}s_{i} &\text{for } i<j\\
				\end{array}
			\end{eqnarray}
	
	There exist three particular interesting families of simplicial sets:
	$\Delta ^{n}$, $\partial \Delta ^{n}$ and $\bigwedge ^{n}_{k}$; they are
	defined in the following way:	
		\begin{equation}
					\label{simplicial-disk}
			\Delta ^{n}=\Hom _{\Delta}(-,\mathbf{n})
		\end{equation}
	$\partial \Delta ^{n}$ is the subobject of $\Delta ^{n}$ characterized by:	
		\begin{equation}
				\label{simplicial-spheres}
			(\partial \Delta ^{n})_{m}=\{ f:\mathbf{m}\rightarrow \mathbf{n}|f \text{
			is not surjective}\}
		\end{equation}
	and finally $\bigwedge ^{n}_{k}$ is the subobject of $\partial \Delta ^{n}$
	given by:	
		\begin{equation}
				\label{simplicial-horns}
			(\wedge ^{n}_{k})_{m}=\{f:\mathbf{m}\rightarrow \mathbf{n}|
			\{0<1< \cdots < \hat{k}<\cdots <n\} \nsubseteq im(f)\}
		\end{equation}
	where $\{0<1<\cdots <\hat{k}<\cdots <n\}$ denotes the well
	ordered set $\mathbf{n}$ with the $k$ element removed.		
	
	We also have the dual notion of \emph{cosimplicial set}:
	
\begin{defi}
	A cosimplicial set $X$ is a covariant functor from the category
	$\Delta$ to the category of sets.
\end{defi}

	Given any category $\mathcal C$, we can also define \emph{simplical} and
	\emph{cosimplicial} objects in $\mathcal C$, where a simplicial  (respectively cosimplicial) object
	$X$ in $\mathcal C$ is just a contravariant (respectively covariant) functor from $\Delta$
	to $\mathcal C$.
	
	Let $\mathbf{Top}$ be the category of compactly generated 
	Hausdorff topological spaces.
	Consider the following family of objects in $\mathbf{Top}$:
	
	$$|\Delta ^{n}|=\{(t_{0},t_{1},\ldots ,t_{n})|t_{i}\geq 0, \sum t_{i}=1\}\subseteq \mathbb R ^{n+1}$$
	
	We get a cosimplicial object $|\Delta ^{\bullet}|$ in $\mathbf{Top}$ if we define the
	coface and codegeneracy maps for $|\Delta ^{n}|$ as:	
		\begin{eqnarray}
					\label{top-cofaces}
			\begin{array}{l}
				\xymatrix@R=.5ex{\delta ^{i}: |\Delta ^{n}| \ar[r]& |\Delta ^{n+1}|\\
												(t_{0},t_{1},\ldots ,t_{n}) \ar@{|->}[r]&(t_{0},\ldots ,t_{i-1},0,t_{i},\ldots ,t_{n})}
			\end{array}
		\end{eqnarray}
	and
		\begin{eqnarray}
					\label{top-codegeneracies}
			\begin{array}{l}
				\xymatrix@R=.5ex{\sigma^{i}:|\Delta ^{n+1}| \ar[r]& |\Delta ^{n}|\\
												(t_{0},t_{1},\ldots ,t_{n+1}) \ar@{|->}[r]&(t_{0}, \ldots ,t_{i-1},t_{i}+t_{i+1},t_{i+2},\ldots, t_{n+1})}
			\end{array}
		\end{eqnarray}

	Now we are ready to define the \emph{geometric realization}
	functor: 
			$$|-|:\mathbf{SSets}\rightarrow \mathbf{Top}$$
	
	Let $X$ be a simplicial set, then its geometric realization $|X|$ is the
	following topological space:	
		\begin{equation}
				\label{geom-real}
			|X|=\varinjlim _{\Delta ^{n}\downarrow X}|\Delta ^{n}|
		\end{equation}
	where the indexing category to compute the colimit has objects
	given by the simplices over $X$, i.e. maps of simplicial sets
	$\Delta ^{n}\rightarrow X$; and morphisms given by commutative
	triangles:
		$$\xymatrix{\Delta^{n} \ar[rr]^{\theta _{*}} \ar[dr]& & \Delta ^{m}\ar[ld]\\
									&X&}$$
	for $\theta :\mathbf{n}\rightarrow \mathbf{m}$
	
	The geometric realization functor $|-|$ has a right adjoint:
			$$Sing:\mathbf{Top}\rightarrow \mathbf{SSets}$$
	called the \emph{singular functor} and defined in the following way:	
		\begin{equation}
					\label{singular-functor}
			\begin{array}{c}
				\xymatrix@R=.5pt{Sing(T):\Delta ^{op} \ar[r]& Sets\\
									n \ar@{|->}[r] & \Hom _{\mathbf{Top}}(|\Delta^{n}|,T)}
			\end{array}
		\end{equation}
	with faces and degeneracies induced by the cofaces and codegeneracies of the
	cosimplical object $|\Delta ^{\bullet}|$.
	
	We say that a map of simplicial sets $\theta :X \rightarrow Y$ is
	a weak equivalence if its geometric realization
	$|\theta |:|X|\rightarrow |Y|$ is a weak equivalence of
	topological spaces, i.e. $\pi _{i}(|\theta|,*)$ is an isomorphism
	for any $i\geq 0$, and for every choice of base point $*\in |X|$.
	
	With all the previous definitions, we are ready to give a cofibrantly generated
	model category structure on the category of simplicial sets.
	Take $I=\{\partial \Delta ^{n}\hookrightarrow \Delta ^{n}\}$ and
	$J=\{\bigwedge ^{n}_{k}\hookrightarrow \Delta ^{n}\}$.
	
\begin{thm}[Quillen]
		\label{simp-modelcat-structure}
	The category of simplicial sets $\mathbf{SSets}$ has a cofibrantly generated
	model category structure, where the weak equivalences,
	the  set of generating cofibrations $I$ and the set of generating trivial cofibrations
	$J$ are defined as above.
\end{thm}
\begin{proof}
	The proof is probably one of the most difficult ones in abstract homotopy theory.
	We refer the reader to \cite[II.3 theorem 3]{MR0223432}, \cite{MR1711612} or \cite{MR1650134}.
\end{proof}
\end{section}
\begin{section}{Simplicial Model Categories}
		\label{subsection.simp.mod.cats}
	Simplicial model categeries were defined by Quillen in \cite{MR0223432},
	we will follow the approach in \cite[chapter 2]{MR1711612} and \cite[chapter 9]{MR1944041}.
	
\begin{defi}
		\label{def.1.1.4.simplicial-cat}
	Let $\modelcat$ be a category.  We say that $\modelcat$
	is \emph{simplicial} if it satisfies the following axioms:
	\begin{enumerate}
		\item	There exists a functor 
					$$\xymatrix@R=.5ex{\modelcat ^{op}\times \modelcat \ar[r]& \simpsets \\
											X,Y \ar@{|->}[r]& Map(X,Y)}$$
					such that
		\item	The set of $0$-simplices in $Map(X,Y)$ is equal to the set of maps in $\modelcat$
					from $X$ to $Y$, i.e. $Map(X,Y)_{0}=\Hom _{\modelcat}(X,Y)$.
		\item	For every triple $X,Y,Z$ of objects in $\modelcat$,
					there exists a map of simplicial sets called \emph{composition law}
					$$\xymatrix{ \circ _{X,Y,Z}:Map(Y,Z)\times Map(X,Y) \ar[r]& Map(X,Z)}
					$$
					which is compatible with the composition in $\modelcat$.
		\item	There exists a map of simplicial sets $i_{X}:*\rightarrow Map(X,X)$,
					for every object $X\in \modelcat$.
		\item	There exists three commutative diagrams (see \cite[definition 9.1.2]{MR1944041}),
					which give the associativity of the composition law, and right and left unit
					properties for the map $i_{X}$.		
	\end{enumerate}
\end{defi}

\begin{defi}
		\label{def.1.1.4.simp-model-cat}
	Let $\modelcat$ be a model category, we say that
	$\modelcat$ is a \emph{simplicial model category}
	if it is a simplicial category (see definition \ref{def.1.1.4.simplicial-cat})
	and satisfies the following two axioms:
	\begin{description}
		\item	[SM0]\label{def.1.1.4.simp-model-cat.SM0}	
					\begin{enumerate}
						\item	For every $X\in \modelcat$, the functor
									$$\xymatrix@R=.5ex{Map(X,-):\modelcat \ar[r]& \mathbf{SSets}\\
															Y \ar@{|->}[r]& Map(X,Y)}
									$$
									has a left adjoint
									$$\xymatrix@R=.5ex{X\otimes -:\mathbf{SSets} \ar[r]& \modelcat \\
															K \ar@{|->}[r]& X\otimes K}
									$$
									such that the adjuntion is compatible with the simpicial structure
									on $\modelcat$, i.e.
									$Map(X\otimes K,Y)\cong Map(K,Map(X,Y))$, where the simplicial set on the right hand side
									is the one defined in remark \ref{rmk1.8simpsets-simpmodcat}(\ref{rmk1.8simpsets-simpmodcat.a}).
						\item For every $Y\in \modelcat$, the functor
									$$\xymatrix@R=.5ex{Map(-,Y):\modelcat ^{op} \ar[r]& \mathbf{SSets}\\
															X \ar@{|->}[r]& Map(X,Y)}
									$$
									has a left adjoint
									$$\xymatrix@R=.5ex{Y^{-}:\mathbf{SSets} \ar[r]& \modelcat ^{op}\\
															K \ar@{|->}[r]& Y^{K}}
									$$
									such that the adjunction is compatible with the simplicial structure
									on $\modelcat$, i.e.
									$Map(X,Y^{K})\cong Map(K,Map(X,Y))$, where the simplicial set on the right hand side
									is the one defined in remark \ref{rmk1.8simpsets-simpmodcat}(\ref{rmk1.8simpsets-simpmodcat.a}).
					\end{enumerate}
		\item	[SM7]	\label{def.1.1.4.simp-model-cat.SM7}	For any cofibration $i:A\rightarrow B$ in $\modelcat$
								and fibration $p:X\rightarrow Y$ in $\modelcat$, the map
								$$\xymatrix{Map(B,X) \ar[rr]^-{(i^{*}, p_{*})}&& Map(A,X)\times _{Map(A,Y)}Map(B,Y)}
								$$
								is a fibration of simplicial sets, which is trivial if either $i$ or $p$ is a weak equivalence.
	\end{description}
\end{defi}

\begin{rmk}
		\label{rmk1.8simpsets-simpmodcat}
	\begin{enumerate}
		\item \label{rmk1.8simpsets-simpmodcat.a}The category of simplicial sets $\simpsets$ has
					a canonical simplicial model category structure
					where $Map(X,Y)$ is the simplicial set
					having $n$-simplices
					$$Map(X,Y)_{n}=\Hom _{\simpsets}(X\times \Delta ^{n},Y)
					$$
					with faces and degeneracies induced from the cosimplicial
					object $\Delta ^{\bullet}$.
		\item	The associated category of pointed simplicial sets $\simpsets _{*}$
					equipped with the induced model structure from $\simpsets$
					(see remark \ref{rmk1.1.triv-conseq}) has a natural
					simplicial model category structure.
	\end{enumerate}
\end{rmk}

\begin{lem}
		\label{lem1.1.4.simplicial-liftings}
	Let $\modelcat$ be a simplicial model category.
	Suppose that $i:A\rightarrow B$, $p:X\rightarrow Y$ are maps in $\modelcat$
	and $j:L\rightarrow K$ is a map of simplicial sets.  Then the
	following are equivalent:
	\begin{enumerate}
		\item	For every solid commutative diagram of simplicial sets
					$$\xymatrix{L \ar[rr] \ar[d]_-{j}&& Map(B,X) \ar[d]^-{(i^{*}, p_{*})}\\
											K \ar[rr] \ar@{-->}[urr]&& Map(A,X)\times _{Map(A,Y)}Map(B,Y)}
					$$
					the dotted arrow making the diagram commutative exists.
		\item	For every solid commutative diagram in $\modelcat$
					$$\xymatrix{A\otimes K\coprod _{A\otimes L}B\otimes L \ar[rr] \ar[d]_{i\Box j}&& X \ar[d]^-{p} \\
											B\otimes K \ar[rr] \ar@{-->}[urr]&& Y}
					$$
					the dotted arrow making the diagram commutative exists.
		\item	For every solid commutative diagram in $\modelcat$
					$$\xymatrix{A \ar[d]_-{i} \ar[r]& X^{K} \ar[d]^-{(j^{*}, p_{*})}\\
											B \ar[r] \ar@{-->}[ur]& X^{L}\times _{Y^{L}}Y^{K}}
					$$
					the dotted arrow making the diagram commutative exists.
	\end{enumerate}
\end{lem}
\begin{proof}
	Follows directly from the existence of the adjunctions in axiom $\mathbf{SM0}$.
\end{proof}

	The following is a useful criterion to check
	axiom $\mathbf{SM7}$.

\begin{prop}
		\label{prop1.1.4.equiv-cond-SM7}
	Let $\modelcat$ be a model category with a simplicial structure
	(see definition \ref{def.1.1.4.simplicial-cat}), satisfying
	axiom $\mathbf{SM0}$, then the
	following are equivalent:
	\begin{enumerate}
		\item	$\modelcat$ satisfies axiom $\mathbf{SM7}$.
		\item	Suppose that $i:A\rightarrow B$ is a cofibration in $\modelcat$,
					and $j:L\rightarrow K$ is a cofibration of simplicial sets, then the map
					$$\xymatrix{A\otimes K\coprod _{A\otimes L}B\otimes L \ar[r]^-{i\Box j}& B\otimes K}
					$$
					is a cofibration in $\modelcat$, which is trivial if
					either $i$ or $j$ is a weak equivalence.
		\item	Suppose that $p:X\rightarrow Y$ is a fibration in $\modelcat$,
					and $j:L\rightarrow K$ is a cofibration of simplicial sets,
					then the map
					$$\xymatrix{X^{K} \ar[rr]^-{(j^{*}, p_{*})}&& X^{L}\times _{Y^{L}}Y^{K}}
					$$
					is a fibration in $\modelcat$, which is trivial if
					either $p$ or $j$ is a weak equivalence.
	\end{enumerate}
\end{prop}
\begin{proof}
	Follows from lemma \ref{lem1.1.4.simplicial-liftings} and 
	corollary \ref{cor.1.1.class-cof-fibs}.
\end{proof}

	These characterizations of axiom $\mathbf{SM7}$, allow
	to construct ``simplicial'' cylinder (respectively path) objects
	for any cofibrant (respectively fibrant) object $A$ of $\modelcat$.
	
\begin{prop}
		\label{prop.1.1.4.simp-cylinders}
	Let $\modelcat$ be a simplicial model category,
	and let $A$ be a cofibrant object in $\modelcat$.
	Then the following diagram represents a
	cylinder object for $A$
	$$\xymatrix{A\otimes \partial \Delta ^{1}\cong A\coprod A \ar[d]_-{i} \ar[dr]^-{\nabla}& \\
							A\otimes \Delta ^{1} \ar[r]_-{s}& A\otimes *\cong A}
	$$
\end{prop}
\begin{proof}
	Proposition \ref{prop1.1.4.equiv-cond-SM7} implies that
	$i$ is a cofibration.  In the following commutative diagram
	$$\xymatrix{A\otimes *\cong A \ar[d]_-{t} \ar[dr]^-{id}&\\
						A\otimes \Delta ^{1} \ar[r]_-{s}& A\otimes *\cong A}
	$$
	proposition
	\ref{prop1.1.4.equiv-cond-SM7} implies that $t$ is a trivial cofibration,
	so by the two out of three property for weak equivalences we have
	that $s$ is a weak equivalence.  It only remains to show that
	$A\otimes \partial \Delta ^{1}\rightarrow A\otimes *$
	is the fold map $A\coprod A\rightarrow A$, but this follows from
	the next commutative diagram:
	$$\xymatrix{A\otimes \ar[d]_-{id\otimes d_{1}} \ar[drr]^-{id}* &&\\
							 A\otimes \partial \Delta ^{1}\ar[r]^-{i}& A\otimes \Delta ^{1} \ar[r]^-{s}& A\otimes *\\
							A\otimes \ar[u]^-{id \otimes d_{0}} \ar[urr]_{id}* &&}
	$$
\end{proof}

	The dual statement for path objects is the following.
	
\begin{prop}
		\label{prop.1.1.4.simp-pathobjects}
	Let $\modelcat$ be a simplicial model category,
	and let $X$ be a fibrant object in $\modelcat$.
	Then the following diagram represents a
	path object for $X$
	$$\xymatrix{& X^{\Delta ^{1}} \ar[d]^-{p}\\
							X\cong X^{*} \ar[ur]^-{r} \ar[r]_-{\Delta}& X^{\partial \Delta ^{1}}\cong X\times X}
	$$
	\flushright $\square$
\end{prop}
	
	One of the interesting consequences we get when we have a simplicial
	model category $\modelcat$, is that we can compute the
	maps in the homotopy category $\mathbf{Ho}\modelcat$
	simplicially.
	
\begin{prop}
		\label{prop1.1.4.simp-homot-maps}
	Let $X,Y$ be a pair of objects in $\modelcat$,
	where $X$ is cofibrant and $Y$ is fibrant.
	Then $[X,Y]=\pi _{0}Map(X,Y)$,
	where $[X,Y]=\Hom _{\mathbf{Ho}\modelcat}(X,Y)$. 
\end{prop}
\begin{proof}
	Since $X$ is cofibrant and $Y$ is fibrant,
	we have that $[X,Y]$ is just the set
	of homotopy classes of maps between $X$ and $Y$.
	On the other hand, axiom $\mathbf{SM7}$
	implies that $Map(X,Y)$ is a fibrant simplicial set
	(Kan complex), so
	$\pi _{0}Map(X,Y)$ is computed using the simplicial homotopies
	given by $\Delta ^{1}\rightarrow Map(X,Y)$, which by the adjunction
	are in bijection with the homotopies given by $X\otimes \Delta ^{1}\rightarrow Y$.
	But these are just homotopies between $X$ and $Y$,
	since proposition \ref{prop.1.1.4.simp-cylinders}
	implies that $X\otimes \Delta ^{1}$ is a cylinder object for $X$.
\end{proof}

\begin{cor}
		\label{cor1.1.4.simp-homot-maps.2}
	Let $\modelcat$ be a simplicial model category,
	and consider a couple of objects $X,Y$ in $\modelcat$.
	Then $[X,Y]=\pi _{0}Map(RQX,RQY)$.
\end{cor}
\begin{proof}
	By construction $[X,Y]$ is equal to set of
	homotopy classes of maps between $RQX$ and $RQY$.
	But $RQX,RQY$ are both cofibrant and fibrant objects
	in $\modelcat$, so proposition \ref{prop1.1.4.simp-homot-maps}
	implies that this set of homotopy classes of maps is equal to
	$\pi _{0}Map(RQX,RQY)$.
\end{proof}
	
	Another simple but very useful consequence of having a simplicial model
	category $\modelcat$, is that we can also detect
	weak equivalences in $\modelcat$ at the level of simplicial sets.
	
\begin{prop}
		\label{prop.1.1.4.detect-weak-equiv.1}
	Let $\modelcat$ be a simplicial model category,
	and let $h:A\rightarrow B$ be a map between two
	cofibrant (respectively fibrant) objects in $\modelcat$.
	Then $h$ is a weak equivalence if and only if
	for every fibrant (respectively cofibrant) object $X$ in $\modelcat$,
	$h^{*}:Map(B,X)\rightarrow Map(A,X)$ (respectively $h_{*}:Map(X,A)\rightarrow Map(X,B)$)
	is a weak equivalence of simplicial sets.
\end{prop}
\begin{proof}
	By duality, it is enough to consider the case in which $A,B$ are cofibrant
	objects in $\modelcat$.
	Assume that $h$ is a weak equivalence.
	Since weak equivalences of simplicial sets have the two out of three property,
	then by Ken Brown's lemma (see lemma \ref{lem1.1.KenBrown}) we can assume that
	$h$ is a trivial cofibration.  The conclusion then follows from axiom $\mathbf{SM7}$
	which implies that for any fibrant object $X$ in $\modelcat$,
	$h^{*}:Map(B,X)\rightarrow Map(A,X)$ is a trivial fibration of simplicial sets,
	so in particular $h^{*}$ is a weak equivalence.
	
	For the converse, it is enough to show that
	$h^{*}:[B,X]\rightarrow [A,X]$ is a bijection for every fibrant object $X$
	in $\modelcat$.  But since for every fibrant object $X$ in $\modelcat$,
	$h^{*}:Map(B,X)\rightarrow Map(A,X)$ is a weak equivalence of simplicial sets,
	in particular we have that
	$h^{*}:\pi _{0}Map(B,X)\rightarrow \pi _{0}Map(A,X)$ is a bijection, and the
	result follows from proposition \ref{prop1.1.4.simp-homot-maps} since
	$A,B$ are cofibrant in $\modelcat$ and $X$ is fibrant in $\modelcat$.
\end{proof}

\begin{cor}
		\label{cor.1.1.4.detect-weak-equiv.2}
	Let $\modelcat$ be a simplicial model category and consider a couple
	of objects $A,B$ in $\modelcat$, and a map $h:A\rightarrow B$
	between them.
	Then the following conditions are equivalent:
	\begin{enumerate}
		\item	\label{cor.1.1.4.detect-weak-equiv.2.a}	$h$ is a weak equivalence in $\modelcat$.
		\item	\label{cor.1.1.4.detect-weak-equiv.2.b}	For every fibrant object $X$ in $\modelcat$,
					$(Qh)^{*}:Map(QB,X)\rightarrow Map(QA,X)$ is a weak
					equivalence of simplicial sets.
		\item	\label{cor.1.1.4.detect-weak-equiv.2.c}	For every cofibrant object $C$ in $\modelcat$,
					$(Rh)_{*}:Map(C,RA)\rightarrow Map(C,RB)$ is
					a weak equivalence of simplicial sets.
	\end{enumerate}	
\end{cor}
\begin{proof}
	(\ref{cor.1.1.4.detect-weak-equiv.2.a}) $\Leftrightarrow$ (\ref{cor.1.1.4.detect-weak-equiv.2.b}).
	We have that $h$ is a weak equivalence if and only if
	every (or some) cofibrant approximation $Qh:QA\rightarrow QB$ is also a weak equivalence.
	Since $QA,QB$ are cofibrant the result follows from proposition \ref{prop.1.1.4.detect-weak-equiv.1}.
	
	(\ref{cor.1.1.4.detect-weak-equiv.2.a}) $\Leftrightarrow$ (\ref{cor.1.1.4.detect-weak-equiv.2.c}).
	We know that $h$ is a weak equivalence if and only if every (or some)
	fibrant approximation $Rh:RA\rightarrow RB$ is also a weak equivalence.
	But $RA,RB$ are fibrant, so the result follows from
	proposition \ref{prop.1.1.4.detect-weak-equiv.1}.
\end{proof}
\end{section}

\begin{section}{Symmetric Monoidal Model Categories}
		\label{subsection-monmodcats}
		
		Symmetric monoidal model categories were introduced
		by Hovey in \cite[chapter 4]{MR1650134}.
		In this section we just recall some of his definitions
		and results without proof.  This is the language that we will
		use in section \ref{section.3.5.applications} to construct external 
		pairings for the slice filtration.
		
\begin{defi}
		\label{def.leftmods/moncats}
	Let $\mathcal C$ be a monoidal category.
	We say that a category $\mathcal D$ is a 
	\emph{left $\mathcal C$-module} if the following
	conditions are satisfied:
	\begin{enumerate}
		\item	There exists a bifunctor
					$\otimes :\mathcal C \times \mathcal D \rightarrow \mathcal D$
		\item	For every pair of objects $X,Y$ in $\mathcal C$ and every object $A$ in $\mathcal D$
					there exists a natural isomorphism
					$a:(X\otimes Y)\otimes A \rightarrow X\otimes (Y\otimes A)$.
		\item	For every object $A$ in $\mathcal D$ there exists a natural isomorphism
					$l:\mathbf{1}\otimes A\rightarrow A$, where $\mathbf{1}$ denotes the
					unit for the monoidal structure on $\mathcal C$.
		\item	Three coherence diagrams commute (see \cite[definition 4.1.6]{MR1650134}).
	\end{enumerate}	
\end{defi}

	We also have right modules over a given monoidal category.
	
\begin{defi}
		\label{def.adj-two-vars}
	Given three categories $\mathcal{C,D,E}$, we 
	define an \emph{adjunction of two variables}
	as a bifunctor $\otimes:\mathcal C \times \mathcal D \rightarrow \mathcal E$
	together with two extra functors
	$\mathrm{Hom}_{r}:\mathcal D ^{op}\times \mathcal E \rightarrow \mathcal C$
	and $\mathrm{Hom}_{l}:\mathcal C ^{op}\times \mathcal E \rightarrow \mathcal D$,
	such that there exist the following two adjunctions:
	\begin{enumerate}
		\item	$\xymatrix{\mathrm{Hom}_{\mathcal E}(X\otimes Y,Z) \ar[r]^-{\varphi _{r}}& 
											\mathrm{Hom}_{\mathcal C}(X,\mathrm{Hom}_{r}(Y,Z))}$
		\item	$\xymatrix{\mathrm{Hom}_{\mathcal E}(X\otimes Y,Z) \ar[r]^-{\varphi _{l}}& 
											\mathrm{Hom}_{\mathcal D}(Y,\mathrm{Hom}_{l}(X,Z))}$
	\end{enumerate}
\end{defi}

\begin{defi}
		\label{def.closed-mon-cat}
	We say that a category $\mathcal C$
	is \emph{closed monoidal} if it is a monoidal
	category such that the bifunctor $\otimes :\mathcal C \times \mathcal C \rightarrow \mathcal C$
	giving the monoidal structure
	is an adjunction of two variables. 
\end{defi}

\begin{defi}
		\label{def.Quillen-bifunct}
	Given model categories $\modelcat$, $\modelcattwo$, $\modelcatthree$ an adjunction of two
	variables $\otimes :\modelcat \times \modelcattwo \rightarrow \modelcatthree$
	is called a \emph{Quillen adjunction of two variables}, if given a 
	cofibration $i:A\rightarrow B$ in $\modelcat$ and a cofibration
	$j:C\rightarrow D$ in $\modelcattwo$, the induced map
		$$\xymatrix{i\Box j:A\otimes D \coprod _{A\otimes C} B\otimes C \ar[r]& B\otimes D}
		$$
	is a cofibration in $\modelcatthree$ which is trivial if either
	$i$ or $j$ is a weak equivalence.
	In this case, we will refer to the functor $\otimes$ as a
	\emph{Quillen bifunctor}.
\end{defi}

\begin{lem}[Hovey]
		\label{lem.cond-Quillen-bifunc}
	Let $\modelcat$, $\modelcattwo$, $\modelcatthree$ be three model categories and
	let $\otimes :\modelcat \times \modelcattwo \rightarrow \modelcatthree$
	be an adjunction of two variables.  Then the following conditions are equivalent:
	\begin{enumerate}
		\item	\label{lem.cond-Quillen-bifunc.a}$  \otimes$ is a Quillen bifunctor.
		\item	\label{lem.cond-Quillen-bifunc.b}  Given a cofibration $j:C\rightarrow D$ in $\modelcattwo$ and 
					a fibration $p:X\rightarrow Y$ in $\modelcatthree$, the induced
					map
					$$\xymatrix{(j^{*}, p_{*}):\Hom _{r}(D,X) \ar[r]& \Hom _{r}(C,X)\times _{\Hom _{r}(C,Y)}\Hom _{r}(D,Y)}
					$$
					is a fibration in $\modelcat$ which is trivial if either $j$ or $p$ is a weak equivalence.
		\item	\label{lem.cond-Quillen-bifunc.c}  Given a cofibration $i:A\rightarrow B$ in $\modelcat$ and a fibration $p:X\rightarrow Y$
					in $\modelcatthree$, the induced map
					$$\xymatrix{(i^{*}, p_{*}):\Hom _{l}(B,X) \ar[r]& \Hom _{l}(A,X)\times _{\Hom _{l}(A,Y)}\Hom _{l}(B,Y)}
					$$
					is a fibration in $\modelcattwo$ which is trivial if either $i$ or $p$ is a weak equivalence.
	\end{enumerate}
\end{lem}
\begin{proof}
	Follows immediately from the adjunctions that appear in the definition of an
	adjunction of two variables (see definition \ref{def.adj-two-vars}), and the lifting property
	characterization for cofibrations, fibrations, trivial cofibrations and trivial fibrations.
\end{proof}
	
\begin{rmk}[cf. \cite{MR1650134}]
		\label{rmk.Quillenbifunc=>Quillenfuncts}
	Let $\otimes :\modelcat \times \modelcattwo \rightarrow \modelcatthree$ be a Quillen
	bifunctor.  Then if $A$ is a cofibrant object in $\modelcat$,
	the functor $A\otimes -:\modelcattwo \rightarrow \modelcatthree$ is a Quillen
	functor with right adjoint $\mathrm{Hom}_{l}(A,-):\modelcatthree \rightarrow \modelcattwo$.
	Similarly if $B$ is a cofibrant object in $\modelcattwo$, we get a Quillen
	functor $-\otimes B:\modelcat \rightarrow \modelcatthree$ with right adjoint
	$\mathrm{Hom}_{r}(B,-)$.  Finally, if $X$ is a fibrant object in $\modelcatthree$,
	we get a Quillen functor $\mathrm{Hom}_{r}(-,X):\modelcattwo \rightarrow \modelcat ^{op}$
	with right adjoint $\mathrm{Hom}_{l}(-,X):\modelcat ^{op}\rightarrow \modelcattwo$.
\end{rmk}

\begin{defi}
		\label{def.mon-mod-cats}
	A \emph{monoidal model category} $\modelcat$ is a
	closed monoidal category with a model category structure, such
	that the following conditions are satisfied:
	\begin{enumerate}
		\item	\label{def.mon-mod-cats.a}The bifunctor $\otimes :\modelcat \times \modelcat \rightarrow \modelcat$
					giving the monoidal structure is a Quillen bifunctor.
		\item	\label{def.mon-mod-cats.b}Let $q:Q\mathbf{1}\rightarrow \mathbf{1}$ be a cofibrant replacement
					for the unit $\mathbf{1}$.  Then the natural maps
					$q\otimes id:Q\mathbf{1}\otimes A\rightarrow \mathbf{1}\otimes A$,
					$id\otimes q:A\otimes Q\mathbf{1}\rightarrow A\otimes \mathbf{1}$
					are weak equivalences for any cofibrant object $A$ in $\modelcat$.
	\end{enumerate}
\end{defi}

	We have an analogous definition for \emph{symmetric monoidal categories}.
	
\begin{prop}[Quillen]
		\label{prop.ssets-sym-mon-mod-cat}
	The category of simplicial sets $\simpsets$ is
	a symmetric monoidal model category.
\end{prop}
\begin{proof}
	We refer the reader to \cite[II.3 theorem 3]{MR0223432}.
\end{proof}

\begin{prop}[Hovey]
		\label{prop.pointed-mon-cats}
	Let $\modelcat$ be a monoidal model category,
	with unit $\mathbf{1}$ equal to the terminal
	object $*$, and assume that $*$ is cofibrant.
	Then the associated pointed category
	$\modelcat _{*}$ (equipped with the monoidal structure
	described in remark \ref{rmk1.1.triv-conseq})
	is also a monoidal model category,
	which is symmetric if $\modelcat$ is.
\end{prop}
\begin{proof}
	We refer the reader to \cite[proposition 4.2.9]{MR1650134}.
\end{proof}

\begin{cor}
		\label{cor.pointed-simpsets-symmmon}
	The category of pointed simplicial sets
	$\simpsets _{*}$ is a symmetric monoidal model category.
\end{cor}
\begin{proof}
	Follows immediately from propositions \ref{prop.ssets-sym-mon-mod-cat} and \ref{prop.pointed-mon-cats}.
\end{proof}

\begin{defi}
		\label{def.mon-Quillen-func}
	Let $(F,G,\varphi):\modelcat \rightarrow \modelcattwo$ be a Quillen adjunction
	between two monoidal model categories.  We say that $(F,G,\varphi)$ is a
	\emph{monoidal Quillen adjunction} if $F$ is a monoidal functor (see \cite[definition 4.1.2]{MR1650134})
	and the map $F(q_{\mathbf{1}}):F(Q\mathbf{1})\rightarrow F\mathbf{1}$ is a
	weak equivalence.  In this situation
	we say that $F$ is a \emph{left Quillen monoidal functor}.
\end{defi}

\begin{defi}
		\label{def.module-modcats}
	Let $\modelcat$ be a monoidal model category.
	A \emph{$\modelcat$-model category} is a left 
	$\modelcat$-module $\modelcattwo$ equipped with a model category structure
	such that the following conditions hold:
	\begin{enumerate}
		\item	\label{def.module-modcats.a}The action map 
					$-\otimes -:\modelcat \times \modelcattwo \rightarrow \modelcattwo$
					is a Quillen bifunctor.
		\item	\label{def.module-modcats.b}If $q:Q\mathbf{1}\rightarrow \mathbf{1}$ 
					is a cofibrant replacement for $\mathbf{1}$
					in $\modelcat$, then the map
					$q\otimes id:Q\mathbf{1}\otimes A \rightarrow \mathbf{1}\otimes A$ is a weak equivalence
					for every cofibrant object $A$ in $\modelcattwo$.
	\end{enumerate}
\end{defi}

	The simplicial model categories discussed in section \ref{subsection.simp.mod.cats}
	are just $\simpsets$-model categories.

\begin{prop}[Hovey]
		\label{prop.1.10.M-modelcats-descend-pointed}
	Let $\modelcat$ be a monoidal model category
	where the unit $\mathbf{1}$ is equal to the
	terminal object $\ast$.  Assume that $\ast$
	is cofibrant.  If $\modelcattwo$ is
	a $\modelcat$-model category, then the
	associated pointed category $\modelcattwo _{*}$ has a
	natural $\modelcat _{*}$-model category structure.
\end{prop}
\begin{proof}
	We refer the reader to \cite[proposition 4.2.19]{MR1650134}.
\end{proof}

\begin{prop}[Hovey]
		\label{prop.Quillenbifunc-descends-homot}
	Let $\modelcat$, $\modelcattwo$, $\modelcatthree$ be three model categories, and let
	$-\otimes -:\modelcat \times \modelcattwo \rightarrow \modelcatthree$
	be a Quillen bifunctor.  Then the total
	derived functors define an adjunction of two variables
	$\otimes ^{\mathbf L}:\mathrm{Ho}\modelcat \times \mathrm{Ho}\modelcattwo
	\rightarrow \mathrm{Ho}\modelcatthree$, with adjoints given by
	$\mathbf{R}\mathrm{Hom}_{l}:(\mathrm{Ho}\modelcat)^{op}\times \mathrm{Ho}\modelcatthree 
	\rightarrow \modelcattwo$ and
	$\mathbf{R}\mathrm{Hom_{r}:(\mathrm{Ho}\modelcattwo)^{op}\times \mathrm{Ho}\modelcatthree 
	\rightarrow \mathrm{Ho}\modelcat}$.
\end{prop}
\begin{proof}
	We refer the reader to \cite[proposition 4.3.1]{MR1650134}.
\end{proof}

\begin{thm}[Hovey]
		\label{thm.symmmon-descends-homot}
	Let $\modelcat$ be a (symmetric) monoidal model category.
	Then $\mathrm{Ho}\modelcat$ can be given the structure of a
	closed (symmetric) monoidal category.
	The adjunction of two variables $(\otimes ^{\mathbf L},\mathbf{R}\mathrm{Hom}_{l},
	\mathbf{R}\mathrm{Hom}_{r})$ which gives the closed structure on $\mathrm{Ho}\modelcat$
	is the total derived adjunction of $(\otimes, \mathrm{Hom}_{l},\mathrm{Hom}_{r})$
	described in proposition \ref{prop.Quillenbifunc-descends-homot}.
	The associativity and unit isomorphisms (and the commutativity isomorphism
	in case $\modelcat$ is symmetric) on $\mathrm{Ho}\modelcat$ are derived from the corresponding
	isomorphisms of $\modelcat$.
\end{thm}
\begin{proof}
	We refer the reader to \cite[theorem 4.3.2]{MR1650134}.
\end{proof}

\end{section}
\begin{section}{Localization of Model Categories}
		\label{section.loc-mod-cats}

	In this section we recall some of Hirschhorn's constructions \cite[sections 3.1, 3.2]{MR1944041}
	restricted to the case where all the model categories are simplicial.

\begin{defi}
		\label{def-1.leftlocmodcats}
	Let $\modelcat$ be a model category and let $\mathcal V$ be a class of maps
	in $\modelcat$.  A \emph{left localization of $\modelcat$ with respect
	to $\mathcal V$} is a model category $L_{\mathcal V}\modelcat$
	equipped with a left Quillen functor $\lambda :\modelcat \rightarrow
	L_{\mathcal V}\modelcat$ satisfying the following properties:
	\begin{enumerate}
		\item	The total left derived functor $\mathbf{L}\lambda :\hocat \rightarrow
					\hocatleft$ takes the images in $\hocat$ of
					the elements in $\mathcal V$ into isomorphisms in $\hocatleft$.
		\item	If $\modelcattwo$ is a model category and $\tau :\modelcat \rightarrow \modelcattwo$
					is a left Quillen functor such that $\mathbf{L}\tau :\hocat \rightarrow \hocattwo$
					takes the images in $\hocat$ of the elements of $\mathcal V$ into
					isomorphisms in $\hocattwo$, then there exists a unique left 
					Quillen functor $\sigma :L_{\mathcal V}\modelcat \rightarrow \modelcattwo$
					with $\sigma \lambda =\tau$.
	\end{enumerate}
\end{defi}

\begin{defi}
		\label{def-1.rightlocmodcats}
	Let $\modelcat$ be a model category and let $\mathcal V$ be a class of maps
	in $\modelcat$.  A \emph{right localization of $\modelcat$ with respect
	to $\mathcal V$} is a model category $R_{\mathcal V}\modelcat$
	equipped with a right Quillen functor $\rho :\modelcat \rightarrow
	R_{\mathcal V}\modelcat$ satisfying the following properties:
	\begin{enumerate}
		\item	The total right derived functor $\mathbf{R}\rho :\hocat \rightarrow
					\hocatright$ takes the images in $\hocat$ of
					the elements in $\mathcal V$ into isomorphisms in $\hocatright$.
		\item	if $\modelcattwo$ is a model category and $\tau :\modelcat \rightarrow \modelcattwo$
					is a right Quillen functor such that $\mathbf{R}\tau :\hocat \rightarrow \hocattwo$
					takes the images in $\hocat$ of the elements of $\mathcal V$ into
					isomorphisms in $\hocattwo$, then there exists a unique right 
					Quillen functor $\sigma :R_{\mathcal V}\modelcat \rightarrow \modelcattwo$
					with $\sigma \rho =\tau$.
	\end{enumerate}
\end{defi}

	From the universal property, we immediately get the following uniqueness
	statement.
\begin{rmk}
		\label{rmk-1.2-univproplocmodcats}
	Let $\modelcat$ be a model category and $\mathcal V$ a class of maps
	in $\modelcat$.  If a left or right localization of $\modelcat$ with
	respect to $\mathcal V$ exists, then it is unique up to a unique isomorphism.
\end{rmk}

\begin{defi}
		\label{def1.2localob-weq}
	Let $\modelcat$ be a model category and $\mathcal V$ a class of maps in
	$\modelcat$.
	\begin{enumerate}
		\item	An object $A$ of $\modelcat$ is \emph{$\mathcal V$-local} if $A$
					is fibrant and for every map $f:C\rightarrow D$ in $\mathcal V$,
					the induced map of simplicial sets
					$Map(QD,A)\rightarrow Map(QC,A)$ is a weak equivalence.
		\item	A map $f:C\rightarrow D$ in $\modelcat$ is a
					\emph{$\mathcal V$-local equivalence} if for every
					$\mathcal V$-local object $A$, the induced map
					of simplicial sets
					$Map(QD,A)\rightarrow Map(QC,A)$ is a weak equivalence.
	\end{enumerate} 
\end{defi}

\begin{defi}
		\label{def1.2colocalob-weq}
	Let $\modelcat$ be a model category and $\mathcal V$ a class of maps in
	$\modelcat$.
	\begin{enumerate}
		\item	An object $A$ of $\modelcat$ is \emph{$\mathcal V$-colocal} if $A$
					is cofibrant and for every map $f:C\rightarrow D$ in $\mathcal V$,
					the induced map of simplicial sets
					$Map(A,RC)\rightarrow Map(A,RD)$ is a weak equivalence.
		\item	A map $f:C\rightarrow D$ in $\modelcat$ is a
					\emph{$\mathcal V$-colocal equivalence} if for every
					$\mathcal V$-colocal object $A$, the induced map
					of simplicial sets
					$Map(A,RC)\rightarrow Map(A,RD)$ is a weak equivalence.
	\end{enumerate} 
\end{defi}

	The following definition will be necessary for
	the construction of right Bousfield localizations.

\begin{defi}
		\label{def1.2.cellularization}
	Let $\modelcat$ be a model category and let $K$ be
	a \emph{set} of objects in $\modelcat$.
	\begin{enumerate}
		\item	A map $g:X\rightarrow Y$ is a $K$-colocal equivalence
					if for every object $A$ in $K$ the induced map
					of simplicial sets $(Rg)_{*}:Map(QA,RX)\rightarrow Map(QA,RY)$
					is a weak equivalence.
		\item	If $\mathcal V$ is the class of $K$-colocal equivalences,
					then a $\mathcal V$-colocal object will be called $K$-colocal.
	\end{enumerate}
\end{defi}

\begin{prop}[Hirschhorn]
		\label{prop1.2.2outof3locwe}
	Let $\modelcat$ be a model category and let $\mathcal V$ be
	a class of maps in $\modelcat$.
	\begin{enumerate}
		\item	The class of $\mathcal V$-local equivalences
					satisfies the \emph{two out of three} property
					(see $\mathbf{MC2}$ in definition \ref{def1.1.1}).
		\item	The class of $\mathcal V$-colocal equivalences
					satisfies the \emph{two out of three} property.
	\end{enumerate}
\end{prop}
\begin{proof}
	We refer the reader to \cite[proposition 3.2.3]{MR1944041}.
\end{proof}

\end{section}

\begin{section}{Bousfield Localization}
		\label{subsec-Bousloc}

	In this section we review Hirschhorn's construction of Bousfield localizations
	\cite[section 3.3]{MR1944041} in the restricted situation where all the model categories
	are simplicial.  These constructions will be the main technical ingredient in our
	approach to produce a lifting of the slice filtration to the model category
	setting (see chapter \ref{chap-slice-filtration}).
	
\begin{defi}
		\label{def1.2.1.leftBousloc}
	Let $\modelcat$ be a model category and let $\mathcal V$ be a class of maps
	in $\modelcat$.  The \emph{left Bousfield localization} of $\modelcat$ with
	respect to $\mathcal V$ (in case it exists) is a model category structure
	$L_{\mathcal V}\modelcat$ on the underlying category of $\modelcat$ such that
	\begin{enumerate}
		\item	the class of weak equivalences of $L_{\mathcal V}\modelcat$ is defined as the
					class of $\mathcal V$-local equivalences of $\modelcat$ (see definition \ref{def1.2localob-weq}).
		\item	the class of cofibrations of $L_{\mathcal V}\modelcat$ is the same
					as the class of cofibrations of $\modelcat$.
		\item	the class of fibrations of $L_{\mathcal V}\modelcat$ is defined
					as the class of maps that have the right lifting property with respect
					to the maps which are cofibrations and $\mathcal V$-local equivalences.
	\end{enumerate}
\end{defi}

	We will also need the dual notion of right Bousfield localization.
	
\begin{defi}
		\label{def1.2.1.rightBousloc}
	Let $\modelcat$ be a model category and let $\mathcal V$ be a class of maps
	in $\modelcat$.  The \emph{right Bousfield localization} of $\modelcat$ with
	respect to $\mathcal V$ (in case it exists) is a model category structure
	$R_{\mathcal V}\modelcat$ on the underlying category of $\modelcat$ such that
	\begin{enumerate}
		\item	the class of weak equivalences of $R_{\mathcal V}\modelcat$ is defined as the
					class of $\mathcal V$-colocal equivalences of $\modelcat$ (see definition \ref{def1.2colocalob-weq}).
		\item	the class of fibrations of $R_{\mathcal V}\modelcat$ is the same
					as the class of fibrations of $\modelcat$.
		\item	the class of cofibrations of $R_{\mathcal V}\modelcat$ is defined
					as the class of maps that have the left lifting property with respect
					to the maps which are fibrations and $\mathcal V$-colocal equivalences.
	\end{enumerate}
\end{defi}

\begin{prop}[Hirschhorn]
		\label{prop1.2.1.rel-leftBous-orig}
	Let $\modelcat$ be a model category and $\mathcal V$ a class of maps in $\modelcat$.
	Let $L_{\mathcal V}\modelcat$ be the left Bousfield localization of $\modelcat$
	with respect to $\mathcal V$, then
	\begin{enumerate}
		\item	\label{prop1.2.1.rel-leftBous-orig.a}every weak equivalence in $\modelcat$ is a weak equivalence
					in $L_{\mathcal V}\modelcat$.
		\item	\label{prop1.2.1.rel-leftBous-orig.b}the class of trivial fibrations of $L_{\mathcal V}\modelcat$ equals the class
					of trivial fibrations of $\modelcat$.
		\item	every fibration of $L_{\mathcal V}\modelcat$ is a fibration of $\modelcat$.
		\item	\label{prop1.2.1.rel-leftBous-orig.d}every trivial cofibration of $\modelcat$ is a trivial cofibration of
					$L_{\mathcal V}\modelcat$.
	\end{enumerate}
\end{prop}
\begin{proof}
	We refer the reader to proposition 3.3.3 in \cite{MR1944041}.
\end{proof}

	We then get the dual version for right Bousfield localizations.

\begin{prop}[Hirschhorn]
		\label{prop1.2.1.rel-rightBous-orig}
	Let $\modelcat$ be a model category and $\mathcal V$ a class of maps in $\modelcat$.
	Let $R_{\mathcal V}\modelcat$ be the right Bousfield localization of $\modelcat$
	with respect to $\mathcal V$, then
	\begin{enumerate}
		\item	\label{prop1.2.1.rel-rightBous-orig.a}every weak equivalence in $\modelcat$ is a weak equivalence
					in $R_{\mathcal V}\modelcat$.
		\item	\label{prop1.2.1.rel-rightBous-orig.b}the class of trivial cofibrations of $R_{\mathcal V}\modelcat$ equals the class
					of trivial cofibrations of $\modelcat$.
		\item	every cofibration of $R_{\mathcal V}\modelcat$ is a cofibration of $\modelcat$.
		\item	\label{prop1.2.1.rel-rightBous-orig.d}every trivial fibration of $\modelcat$ is a trivial fibration of
					$R_{\mathcal V}\modelcat$.
	\end{enumerate}
\end{prop}
\begin{proof}
	We refer the reader to proposition 3.3.3 in \cite{MR1944041}.
\end{proof}

\begin{prop}[Hirschhorn]
		\label{prop1.2.1.Quillenfunct-Bousloc}
	Let $\modelcat$ be a model category and $\mathcal V$ a class of maps
	in $\modelcat$.
	\begin{enumerate}
		\item	If $L_{\mathcal V}\modelcat$ is the left Bousfield localization of
					$\modelcat$ with respect to $\mathcal V$, then the identity functor
					$id:\modelcat \rightarrow L_{\mathcal V}\modelcat$ is a left Quillen
					functor with right adjoint $id:L_{\mathcal V}\modelcat \rightarrow \modelcat$.
		\item	If $R_{\mathcal V}\modelcat$ is the right Bousfield localization of
					$\modelcat$ with respect to $\mathcal V$, then the identity functor
					$id:R_{\mathcal V}\modelcat \rightarrow \modelcat$ is a left Quillen
					functor with right adjoint $id:\modelcat \rightarrow R_{\mathcal V}\modelcat$.
	\end{enumerate}
\end{prop}
\begin{proof}
	Follows immediately from propositions \ref{prop1.2.1.rel-leftBous-orig} and \ref{prop1.2.1.rel-rightBous-orig}.
\end{proof}

\begin{thm}[Hirschhorn]
		\label{thm1.2.1.Bousloc=loc}
	Let $\modelcat$ be a model category and let $\mathcal V$ be a class of maps
	in $\modelcat$.
	\begin{enumerate}
		\item	If $L_{\mathcal V}\modelcat$ is the left Bousfield localization of
					$\modelcat$ with respect to $\mathcal V$, then the
					identity functor $id:\modelcat \rightarrow L_{\mathcal V}\modelcat$
					is a left localization of $\modelcat$
					with respect to $\mathcal V$ (see definition \ref{def-1.leftlocmodcats}).
		\item	If $R_{\mathcal V}\modelcat$ is the right Bousfield localization of 
					$\modelcat$ with respect to $\mathcal V$ then the identity
					functor $id:\modelcat \rightarrow R_{\mathcal V}\modelcat$
					is a right localization of $\modelcat$ with respect to $\mathcal V$.
	\end{enumerate}
\end{thm}
\begin{proof}
	We refer the reader to \cite[theorem 3.3.19]{MR1944041}.
\end{proof}

\end{section}

\end{chapter}
\begin{chapter}{Motivic Unstable and Stable Homotopy Theory}
		\label{chapter-mot-hom-theory}

	In this chapter we review the construction of the Morel-Voevodsky motivic stable model
	structure and the construction of Jardine's motivic symmetric stable model structure
	(see sections \ref{section.motivic.spectra} and \ref{sec.2.6.symmetricTspectra}).
	We also show that these two model structures satisfy Hirschhorn's cellularity condition
	(see sections \ref{section-cellularity-motivic-stable} and \ref{section-cellularity-motivic-symmetric-stable}).
	Therefore, it is possible to apply Hirschhorn's localization techniques to get Bousfield
	localizations with respect to these two model structures.  Finally, in section \ref{section.2.8Modules-Algebras}
	we recall the construction of the model structures for the categories
	of $A$-modules and $A$-algebras,
	where $A$ denotes a cofibrant ring spectrum in Jardine's motivic symmetric stable model category.
	We will see that the category of $A$-modules equipped with this model structure also
	satisfies Hirschhorn's cellularity condition.

\begin{section}{The Injective Model Structure}
		\label{sec-glob-inj-mod}

	Let $S$ be a Noetherian separated scheme of finite Krull dimension, 
	and consider the category
	$\smoothS$ of smooth schemes
	of finite type over $S$.  $\nissite$ will denote the site
	with underlying 
	category $\smoothS$ equipped with the Nisnevich topology. 	
	We are interested in the category $\simppre$
	of presheaves of simplicial sets
	on  $\smoothS$.  The objects in $\simppre$ can also be described
	as simplicial presheaves on $\smoothS$.  
	The work of
	Jardine (see \cite{MR906403}) shows in particular that
	$\simppre$ has the structure of a proper simplicial cofibrantly generated 
	model category.
	
	We will denote by $\Delta ^{n}_{U}$ the representable simplicial presheaf
	corresponding to the objects $U$ in $\smoothS$ and $\mathbf{n}$ in $\Delta$,
	i.e.
	$$\xymatrix@R=.5ex{\Delta ^{n}_{U}:(\smoothS \times \Delta)^{op} \ar[r]& Sets\\
										(V,m) \ar@{|->}[r]& (\mathrm{Hom}_{\smoothS}(V,U))\times (\Delta ^{n})_{m}}
	$$
	The following functor gives
	a fully faithful embedding of $\smoothS$ into $\simppre$:
	$$\xymatrix@R=.5pt{Y:\smoothS \ar[r]& \simppre \\
										U \ar@{|->}[r]& \Delta ^{0}_{U}}
	$$
	we will abuse notation and write
	$U$ instead of $\Delta ^{0}_{U}$.  Given any simplicial set $K$
	we can consider the associated constant presheaf of simplicial sets  which we also
	denote by $K$, i.e.
	$$\xymatrix@R=.5pt{K:(\smoothS \times \Delta)^{op} \ar[r]& Sets \\
											(U,n) \ar[r]& K_{n}}
	$$

	The category of simplicial presheaves $\simppre$ inherits a natural simplical structure
	from the one on pointed simplicial sets.
	
	Given a simplicial presheaf $X$, the
	tensor objects for the simplicial structure on $\simppre$ are defined as follows:
	$$\xymatrix{X\otimes -:\simpsets \ar[r]& \simppre}
	$$
	where $X\otimes K$ is the following simplicial presheaf:	 
	$$\xymatrix@R=.5ex{X\otimes K: (\smoothS \times \Delta)^{op} \ar[r]& Sets\\
							(U,n) \ar@{|->}[r]& X_{n}(U)\times K_{n}}$$
	The simplicial functor in two variables is:		
	$$\xymatrix{Map(-,-):(\simppre)^{op}\times \simppre \ar[r]& \simpsets}
	$$
	where $Map(X,Y)$ is the simplicial set given by:
	$$\xymatrix@R=.5ex{Map(X,Y):\Delta ^{op} \ar[r]& Sets \\
							\mathbf{n} \ar@{|->}[r]& \mathrm{Hom}_{\simppre}(X\otimes \Delta ^{n},Y)}
	$$
	and finally for any simplicial presheaf $Y$ we have the following functor
	$$\xymatrix{Y^{-}:\simpsets \ar[r]& (\simppre)^{op}}
	$$
	where $Y^{K}$ is the simplicial presheaf given as follows:
	$$\xymatrix@R=.5ex{Y^{K}: (\smoothS \times \Delta)^{op} \ar[r]& Sets\\
							(U,n) \ar@{|->}[r]& \mathrm{Hom}_{\simpsets}(K \times \Delta ^{n},Y(U))}
	$$
					
	Let $t$ be a point in $\nissite$.  Denote by
	$\theta_{t}$ the fibre functor which assigns to every
	simplicial presheaf its stalk at $t$:
	$$\xymatrix@R=.5pt{\theta _{t}:\simppre \ar[r]& \simpsets\\
							X \ar@{|->}[r]& \theta _{t}(X)=X_{t}}	
	$$

	Now we proceed to define the model structure on $\simppre$
	constructed by Jardine.
	A map $f:X\rightarrow Y$ in $\simppre$ is defined to be a weak equivalence, if
	$f$ induces a weak equivalence of simplicial sets in all the stalks on $\nissite$,
	i.e. if for every point $t$ in $\nissite$ the map
	$$\xymatrix{\theta _{t}(X) \ar[rr]^-{\theta _{t}(f)}&& \theta _{t}(Y)}
	$$
	is a weak equivalence of simplicial sets.
	
	The set $I$ of generating cofibrations is given by all the subobjects of
	$\Delta ^{n}_{U}$ for $U$ in $\smoothS$ and $n\geq 0$, i.e.
	$$I=\{Y\hookrightarrow \Delta ^{n}_{U}|U\in (\smoothS), n\geq 0\}
	$$
	it is easy to see that a map
	$i:X\rightarrow Y$ is in $I$-cell
	if and only if it is a monomorphism, i.e.
	$i_{n}(U):X_{n}(U)\rightarrow Y_{n}(U)$ is an injective map of sets,
	for every $U$ in $\smoothS$, $n\geq 0$.	

	Let $\lambda$ be a cardinal, and $X$ a simplicial presheaf on $\smoothS$.
	We say that $X$ is \emph{$\lambda$-bounded} if the cardinal of all the
	simplices of $X$ is bounded by $\lambda$, i.e.
	$|X_{n}(U)|<\lambda$ for every $U$ in $\smoothS$, $n\geq 0$.
	The site $\nissite$ is essentially small, so we can find a cardinal $\kappa$
	such that $\kappa$ is greater than $2^{\alpha}$, where $\alpha$ 
	is the cardinality of 
	the set $Map(\smoothS)$ of maps in $\smoothS$.
	
	We say that a map $j:X\rightarrow Y$ of simplicial presheaves in $\smoothS$
	is a \emph{trivial cofibration},
	if it is both a cofibration and a weak equivalence.
	The set $J$ of generating trivial cofibrations is given
	by all the trivial cofibrations where the codomain
	is bounded by the cardinal $\kappa$ described above,
	i.e.
	$$J=\{j:X\rightarrow Y|j \text{ is a trivial cofibration and }
		Y \text{ is }\kappa \text{-bounded}\}
	$$
	
\begin{thm}[Jardine]
		\label{thm.2.1.1.glob-inj-mod}
	The category $\simppre$ of simplicial presheaves
	on the Nisnevich site $\nissite$, has the
	structure of a proper simplicial cofibrantly generated model category
	where the class $\mathcal W$ of weak equivalences, and the sets $I,J$
	of generating cofibrations and generating trivial cofibrations
	are defined as above.
\end{thm}
\begin{proof}
	We refer the reader to \cite[theorem 2.3]{MR906403}.
\end{proof}

	The model structure defined above will be called the \emph{injective model
	structure} for $\simppre$.  
	
\begin{rmks}
		\label{rmks.1.1.class-cofs-inj-mod}
	The cofibrations for the injective model structure on
	$\simppre$ have the following properties:	
	\begin{enumerate}
		\item	\label{rmks.1.1.class-cofs-inj-mod.a}The class of cofibrations coincides with the class
					of relative $I$-cell complexes, therefore a map is a cofibration if and only
					if it is a monomorphism.
		\item	\label{rmks.1.1.class-cofs-inj-mod.b}If a map $i:A\rightarrow B$ in $\simppre$ is a cofibration
					then for every point $t$ in $\nissite$ the associated
					map $\theta _{t}(i):\theta _{t}(A)\rightarrow
					\theta _{t}(B)$ is a cofibration of simplicial sets.
		\item	\label{rmks.1.1.class-cofs-inj-mod.c}Every object $A$ in $\simppre$ is an
					$I$-cell complex, therefore every object in  $\simppre$ is cofibrant.
	\end{enumerate}
\end{rmks}

	The category $\simppre$ of simplicial presheaves on the smooth Nisnevich site $\nissite$
	also has a closed symmetric monoidal structure which is compatible
	with the injective model structure, i.e. $\simppre$
	equipped with the injective structure is a symmetric
	monoidal model category in the sense of Hovey (see definition \ref{def.mon-mod-cats}).
	
	The closed symmetric monoidal structure is defined as follows:
	$$\xymatrix@R=.5pt{\simppre \times \simppre \ar[r]& 
							\simppre\\
							(X,Y) \ar@{|->}[r]& X\times Y}
	$$
	where $X\times Y$ is the presheaf of simplicial sets defined as follows:
	$$\xymatrix@R=.5pt{X\times Y:(\smoothS \times \Delta)^{op} \ar[r]& Sets\\
							(U,n) \ar@{|->}[r]& X_{n}(U)\times Y_{n}(U)}
	$$
	and the functor that gives the adjunction of two variables is
	the following:
	$$\xymatrix{\inthompre (-,-):(\simppre)^{op}\times \simppre \ar[r]& 
							\simppre}
	$$
	where $\inthompre(X,Y)$ is the simplicial presheaf given by:
	$$\xymatrix@R=.5pt{\inthompre (X,Y):(\smoothS \times \Delta)^{op} \ar[r]& Sets\\
										(U,n) \ar@{|->}[r]& 
										\Hom_{\simppre}(X\times \Delta _{U}^{n},Y)}
	$$

\begin{prop}
		\label{prop.2.1.adjunction-simplicial-internal}
	Let $X,Y,Z$ be simplicial presheaves on $\nissite$.
	\begin{enumerate}
		\item	\label{prop.2.1.adjunction-simplicial-internal.a}
					There is a natural isomorphism of simplicial sets
					$$\xymatrix{Map(X\times Y,Z) \ar[r]^-{\cong}& Map(X,\inthompre(Y,Z))}
					$$
		\item	\label{prop.2.1.adjunction-simplicial-internal.b}
					There is a natural isomorphism of simplicial presheaves	
					on $\nissite$
					$$\xymatrix{\inthompre(X\times Y,Z) \ar[r]^-{\cong}& \inthompre(X,\inthompre(Y,Z))}
					$$
	\end{enumerate}
\end{prop}
\begin{proof}
	(\ref{prop.2.1.adjunction-simplicial-internal.a}).  
	To any $n$-simplex $\alpha$ in $Map(X\times Y,Z)$
		$$\xymatrix{(X\otimes \Delta ^{n})\times Y \ar[r]^-{\cong}& (X\times Y)\otimes \Delta ^{n} \ar[r]^-{\alpha}& Z}$$
	associate
	the $n$-simplex $\alpha _{*}$ in $Map(X,\inthompre(Y,Z))$
		$$\alpha_{*}:X\otimes \Delta ^{n}\rightarrow \inthompre(Y,Z)$$
	coming from the adjunction between $-\times Y$ and $\inthompre(Y,-)$.
	
	(\ref{prop.2.1.adjunction-simplicial-internal.b}).  
	To any $n$-simplex $\alpha$ in $\inthompre(X\times Y,Z)_{n}(U)$
		$$\xymatrix{(X\times \Delta ^{n}_{U})\times Y \ar[r]^{\cong}& (X\times Y)\times \Delta ^{n}_{U} \ar[r]^-{\alpha}& Z}$$ 
	associate the $n$-simplex $\alpha _{*}$ in $\inthompre(X,\inthompre(Y,Z))_{n}(U)$
		$$\alpha_{*}:X\times \Delta ^{n}_{U}\rightarrow \inthompre(Y,Z)$$
	coming from the adjunction
	between $-\times Y$ and $\inthompre(Y,-)$.
\end{proof}
	
\begin{lem}[cf. \cite{MR1787949}]
		\label{lem2.1.monoidalmodcat}
	The category $\simppre$ of simplicial presheaves on the smooth
	Nisnevich site $\nissite$ equipped with the injective model structure
	is a symmetric monoidal model
	category (see definition \ref{def.mon-mod-cats}).
\end{lem}
\begin{proof}
	We need to check that the conditions (\ref{def.mon-mod-cats.a})-(\ref{def.mon-mod-cats.b})
	in definition \ref{def.mon-mod-cats} are satisfied.
	Since every object is cofibrant in $\simppre$, condition (\ref{def.mon-mod-cats.b}) is
	trivially satisfied.  To check condition (\ref{def.mon-mod-cats.a}),  we need to show
	that if we take
	two cofibrations $i:A\rightarrow B$ and $j:C\rightarrow D$ for the injective model structure
	on $\simppre$, then the induced map
	$$i\square j:A\times D\coprod _{A\times C}B\times C \rightarrow B\times D$$ 
	is also a cofibration, which is trivial if either $i$ or $j$ is a weak equivalence.
	To see that $i\square j$ is a cofibration, it is enough to show that $i\square j(U)$
	is a cofibration of simplicial sets for every $U$ in $\smoothS$, but this is
	true since the category of simplicial sets is a symmetric monoidal model category.
	
	Now we show that $i\square j$ is a trivial cofibration if either $i$ or
	$j$ is a weak equivalence.
	The definition of weak equivalences for the injective model structure implies
	that is enough to prove it at the level of the stalks, so let
	$t$ be any point in $\nissite$, and consider the induced map of simplicial sets
	$$\theta _{t}(i\square j):\theta _{t}(A)\times \theta _{t}(D)\coprod
		_{\theta _{t}(A)\times \theta _{t}(C)}\theta _{t}(B)\times \theta _{t}(C)
		\rightarrow \theta _{t}(B)\times \theta _{t}(D)
	$$
	Now
	since the category of simplicial sets is in particular a
	symmetric monoidal model category,
	we have that $\theta _{t}(i\square j)$ is a trivial cofibration if either 
	$i$ or $j$ is a weak equivalence.  
	Since this holds for every point $t$ in $\nissite$, we have that
	$i\square j$ is a cofibration in $\simppre$ which is trivial
	if either $i$ or $j$ is a weak equivalence, hence the result follows.
\end{proof}

\begin{lem}[Morel-Voevodsky, cf. \cite{MR1813224}]
		\label{lem2.1.criterion-weakequivs}
	Let $X,Y$ be two
	fibrant simplicial presheaves in the injective
	model structure, and consider a map
	$f:X\rightarrow Y$.
	The following are equivalent:
	\begin{enumerate}
		\item	$f$ is a weak equivalence in the injective model
					structure for $\simppre$.
		\item	For every $U$ in $\smoothS$ the map
					$$f(U):X(U)\rightarrow Y(U)$$ 
					is a weak equivalence
					of simplicial sets.
	\end{enumerate}
\end{lem}
\begin{proof}
	Assume that $f$ is a weak equivalence in $\simppre$.
	Since $X,Y$ are fibrant and weak equivalences of simplicial sets
	satisfy the two out of three property, by Ken Brown's lemma (see
	lemma \ref{lem1.1.KenBrown})
	we can assume that $f$ is a trivial fibration in $\simppre$.
	Now consider $U$ as an element of $\simppre$.
	Since every object in $\simppre$ is cofibrant,
	axiom $\mathbf{SM7}$ for simplicial model categories implies that:
	$f_{*}:Map(U,X)\rightarrow Map(U,Y)$
	is a trivial fibration of simplicial sets, but this is just equal to
	$f(U):X(U)\rightarrow Y(U)$.
	
	Conversely, suppose now that for every $U$ in $\simppre$,
	$f(U):X(U)\rightarrow Y(U)$ is
	a weak equivalence of simplicial sets. Let $t$ be an arbitrary point
	in $\nissite$.  We know that $t$ is associated to a pro-object
	$\{U_{\alpha}\}$ in $\nissite$.  Therefore
	$\theta _{t}(f):\theta _{t}(X)\rightarrow \theta _{t}(Y)$ is a
	filtered colimit of weak equivalences of simplicial sets, hence
	a weak equivalence of simplicial sets.
	But this implies that $f$ is a weak equivalence in
	$\simppre$.
\end{proof}

\begin{defi}
		\label{def.2.1.Nisnevich-descent}
	Let $X$ be a simplicial presheaf on $\nissite$.
	We say that $X$ satisfies the \emph{B.G. property}
	if any \emph{elementary Cartesian square}
		$$\xymatrix{U\times _{W}V \ar[r] \ar[d]& V \ar[d]^-{p}\\
								U \ar[r]_-{i}& W}
		$$
	of smooth schemes over $S$ with $p$ \'{e}tale, $i$ an open immersion and
	$p^{-1}(W-U)\cong W-U$ (both equipped with the reduced scheme structure) maps to a homotopy Cartesian diagram of
	simplicial sets after applying $X$
		$$\xymatrix{X(W) \ar[r] \ar[d]& X(V) \ar[d] \\
								X(U) \ar[r]& X(U\times _{W}V)}
		$$
\end{defi}

\begin{thm}[Jardine]
		\label{thm.2.1.Nisnevich-descent}
	Let $X$ be a simplicial presheaf on $\nissite$.
	Then $X$ satisfies the \emph{B.G. property}
	if and only if  any fibrant replacement
	$X\rightarrow GX$ in the injective model structure for
	$\simppre$ is a sectionwise weak equivalence of simplicial
	sets, i.e. for any $U$ in $\smoothS$,
		$$\xymatrix{X(U) \ar[rr]^-{g_{X}(U)}&& GX(U)}
		$$
	$g_{X}(U)$ is a weak equivalence of simplicial sets.
\end{thm}
\begin{proof}
	We refer the reader to \cite[theorem 1.3]{MR1787949}.
\end{proof}

\begin{defi}
		\label{def.2.1.basechangefunctor}
	Consider $U \in \smoothS$ with structure map
	$\phi :U\rightarrow S$, i.e. $\phi$ is a smooth
	map of finite type.  Then we have the following
	adjunction (see \cite[proposition I.5.1]{MR0354652}):	
		$$\xymatrix{(\phi^{-1}, \phi_{\ast},\varphi):\simppre \ar[r]& \Delta ^{op}Pre(Sm|_{U})_{Nis}}
		$$
	where $\phi ^{-1}$ and $\phi _{\ast}$ are defined as follows:	
		$$\xymatrix@R=.5pt{\phi ^{-1}:\simppre \ar[r]& \Delta ^{op}Pre(Sm|_{U})_{Nis}\\
								X \ar@{|->}[r]& \phi ^{-1}X}
		$$
	with $\phi ^{-1}X$ defined as the composition of $\phi$ and $X$:		
		$$\xymatrix{(Sm|_{U}\times \Delta)^{op} \ar[rr]^-{\phi \times id} \ar[dr]_-{\phi ^{-1}X} && 
								(\smoothS \times \Delta)^{op} \ar[dl]^-{X}\\
									& Sets &}
		$$
	and the right adjoint $\phi _{\ast}$ is given by:	
		$$\xymatrix@R=.5pt{\phi _{\ast}:\Delta ^{op}Pre(Sm|_{U})_{Nis} \ar[r]& \simppre \\
								X \ar@{|->}[r]& \phi _{\ast}X}
		$$
	where $\phi _{\ast}X$ is the following simplicial presheaf:		
		$$\xymatrix@R=.5pt{\phi _{\ast}X:(\smoothS \times \Delta)^{op} \ar[r]& Sets \\
								(V,n) \ar@{|->}[r]& X_{n}(V\times _{S} U)}
		$$
\end{defi}

\begin{rmk}
		\label{rmk2.1.counit=isomorphism}
	Let $\phi: U\rightarrow S$ be a smooth map of finite type,
	and let $Y$ be an arbitrary simplicial presheaf on $(Sm|_{U})_{Nis}$.
	It follows immediately from the description of the functors
	$\phi ^{-1}$ and $\phi _{\ast}$ that the counit of the adjunction
	$(\phi ^{-1}, \phi _{\ast}, \varphi)$
		$$\xymatrix{\phi ^{-1}\phi _{\ast}Y \ar[r] \ar[r]^-{\cong}& Y}
		$$
	is an isomorphism which can be
	naturally identified with the identity map on $Y$, in particular
	$\phi ^{-1}\phi _{\ast}Y$ is canonically isomorphic to $Y$.
\end{rmk}

\begin{prop}
		\label{prop.2.1.inthom==basechange}
	Let $\phi :U\rightarrow S$ be a smooth map of finite type, and let
	$X$ be an arbitrary simplicial presheaf on $\nissite$.
	Then we have a canonical isomorphism:
		$$\xymatrix{\phi _{\ast}\phi ^{-1}X \ar[r]^-{\cong}& \inthompre (U,X)}
		$$
\end{prop}
\begin{proof}
	To any $n$-simplex $\alpha$ in $(\phi _{\ast}\phi ^{-1}X)_{n}(V)=X_{n}(V\times _{S}U)$
		$$\xymatrix{\Delta ^{n}_{V}\times U\cong \Delta ^{n}_{V\times _{S}U} \ar[r]^-{\alpha}& X}
		$$
	associate the $n$-simplex $\alpha _{\ast}$ in $\inthompre (U,X)_{n}(V)$
		$$\xymatrix{\Delta ^{n}_{V} \ar[r]^-{\alpha _{\ast}}& \inthompre (U,X)}
		$$
	coming from the adjunction between $-\times U$ and $\inthompre (U,-)$.
\end{proof}

\begin{prop}
		\label{prop.2.1.basechange-enriched-adj}
	Let $\phi :U\rightarrow S$ be a smooth map of finite type, let
	$X$ be a simplicial presheaf on $\nissite$ and $Y$ a simplicial
	presheaf on $(Sm|_{U})_{Nis}$.
	Then we have the following enriched adjunctions:
	\begin{enumerate}
		\item	\label{prop.2.1.basechange-enriched-adj.a}There is a natural isomorphism of simplicial sets
						$$\xymatrix{Map(\phi ^{-1}X,Y) \ar[r]^-{\cong}& Map(X,\phi _{\ast}Y)}
						$$
		\item	\label{prop.2.1.basechange-enriched-adj.b}There is a natural isomorphism of simplicial presheaves on $\nissite$
						$$\xymatrix{\inthompre (X,\phi _{\ast}Y) \ar[r]^-{\cong}& \phi _{\ast}(\inthompre (\phi ^{-1}X,Y))}
						$$
		\item	\label{prop.2.1.basechange-enriched-adj.c}There is a natural isomorphism of simplicial presheaves on $(Sm|_{U})_{Nis}$
						$$\xymatrix{\phi ^{-1}(\inthompre (X, \phi _{\ast}Y)) \ar[r]^-{\cong}& \inthompre (\phi ^{-1}X, Y)}
						$$
	\end{enumerate}
\end{prop}
\begin{proof}
	(\ref{prop.2.1.basechange-enriched-adj.a}): To any $n$-simplex $\alpha$ in $Map(\phi ^{-1}X,Y)$
		$$\xymatrix{\phi ^{-1}(X\otimes \Delta ^{n})\cong \phi ^{-1}X\otimes \Delta ^{n} \ar[r]^-{\alpha}& Y}
		$$
	associate the $n$-simplex $\alpha _{*}$ in $Map(X, \phi _{\ast}Y)$
		$$\xymatrix{X\otimes \Delta ^{n} \ar[r]^-{\alpha _{\ast}}& \phi _{\ast}Y}
		$$
	coming from the adjunction between $\phi ^{-1}$ and $\phi _{\ast}$.
	
	(\ref{prop.2.1.basechange-enriched-adj.b}):  To any $n$-simplex $\alpha$ in
	$\inthompre (X,\phi _{\ast}Y)_{n}(V)$ (where $V\in (\smoothS )$)
		$$\xymatrix{X\times \Delta ^{n}_{V} \ar[r]^-{\alpha}& \phi _{\ast}Y}
		$$
	associate the $n$-simplex $\alpha ^{*}$ in $\phi _{\ast}(\inthompre (\phi ^{-1}X,Y))_{n}(V)$
		$$\xymatrix{\phi ^{-1}X \times \Delta ^{n}_{V\times _{S}U}\cong \phi ^{-1}(X\times \Delta ^{n}_{V})
								\ar[r]^-{\alpha ^{\ast}}& Y}
		$$
	coming from the adjunction between $\phi ^{-1}$ and $\phi _{\ast}$.
	
	(\ref{prop.2.1.basechange-enriched-adj.c}): To any $n$-simplex $\alpha$ in
	$\phi ^{-1}(\inthompre (X, \phi _{\ast}Y))_{n}(V)$ (where $V\in (Sm|_{U})$)
		$$\xymatrix{X\times \Delta ^{n}_{V} \ar[r]^{\alpha}& \phi _{\ast}Y}
		$$
	associate the $n$-simplex $\alpha ^{\ast}$ in $\inthompre (\phi ^{-1}X, Y)_{n}(V)$
		$$\xymatrix{\phi ^{-1}X \times \Delta ^{n}_{V}\cong \phi ^{-1}(X\times \Delta ^{n}_{V}) \ar[r]^-{\alpha ^{\ast}}& Y}
		$$
	coming from the adjunction between $\phi ^{-1}$ and $\phi _{\ast}$.
\end{proof}

\begin{defi}[cf. \cite{MR1787949}]
		\label{def.2.1.flasque-simplicial-presheaves}
	Let $X$ be a simplicial presheaf on $\nissite$.
	We say that $X$ is \emph{flasque} if:
	\begin{enumerate}
		\item	$X$ is a presheaf of Kan complexes.
		\item	Every finite collection $V_{i}\hookrightarrow V$,
					$i=1,\ldots ,n$ of subschemes of a scheme $V$ induces
					a Kan fibration
						$$\xymatrix{X(V)\cong Map(V,X) \ar[r]^-{i^{\ast}}& Map(\cup _{i=1}^{n}V_{i},X)}
						$$
	\end{enumerate}
\end{defi}

\begin{rmk}
		\label{rmk.2.1.inj-fib=>flasque}
	\begin{enumerate}
		\item \label{rmk.2.1.inj-fib=>flasque.a}Let $X$ be a simplicial presheaf on $\nissite$
					which is fibrant in the injective model structure
					for $\simppre$.  Then $X$ is flasque and satisfies the B.G. property.
		\item	\label{rmk.2.1.inj-fib=>flasque.b}The class of flasque simplicial presheaves is closed under
					filtered colimits.
		\item	\label{rmk.2.1.inj-fib=>flasque.c}The B.G. property is stable under filtered colimits.
		\item	\label{rmk.2.1.inj-fib=>flasque.d}The functors $\phi ^{-1}$ and $\phi _{\ast}$ preserve
					flasque simplicial presheaves.
		\item	\label{rmk.2.1.inj-fib=>flasque.e}The functors $\phi ^{-1}$ and $\phi _{\ast}$ preserve
					the B.G. property.
	\end{enumerate}
\end{rmk}
\end{section}

\begin{section}{Cellularity of the Injective Model Structure}
		\label{section.cellularity-injective}

	In this section we prove that the injective model structure
	on $\simppre$ is cellular (see definition \ref{def.1.1.2.cell-mod-cats}).	
	This is an unpublished result due to Hirschhorn, which also appears in 
	\cite[theorem 1.4]{MR2197578}.  
	The author would like to thank Jens Hornbostel
	for the discussion related to Hirschhorn's cellularity results.

\begin{lem}
		\label{lem2.2.all-obj-small}
	Let $A$ be a simplicial presheaf  on the smooth Nisnevich site
	$\nissite$.  Then $A$ is small
	(see definition \ref{def.1.1.1.small-objects}).
\end{lem}
\begin{proof}
	Let $\mu$ be the cardinal of the set $S_{A}$
	of simplices of $A$, i.e.
	$$S_{A}=\coprod _{V\in (\smoothS),n\geq 0}A_{n}(V)
	$$
	and let $\kappa$ be the successor cardinal of $\mu$.
	Since $\kappa$ is a succesor cardinal, we have that $\kappa$
	is a regular cardinal
	(see \cite[proposition 10.1.14]{MR1944041}).
	
	We claim that $A$ is $\kappa$ small with respect to the class
	of all maps in $\simppre$.  In effect, consider an arbitrary
	$\lambda$-sequence where $\lambda$ is a regular cardinal
	greater than $\kappa$,
	$$\lambdaseqtwo
	$$
	we need to show that the map
	$\colimit _{\beta <\lambda}\Hom_{\simppre}(A,X_{\beta})\rightarrow
	\Hom(A,X_{\lambda})$ is a bijection.  To check the injectivity,
	we just take sections on every $U\in (\smoothS)$,
	and use the fact
	that every simplicial set is small (see \cite[lemma 3.1.1]{MR1650134}).
	To check the surjectivity, consider an arbitrary map
	$f:A\rightarrow X_{\lambda}$, now the restriction of $f$ to
	every simplex of $A$ ($\Delta ^{n}_{U}\rightarrow A$), factors
	through some $X_{\beta}$ with $\beta <\lambda$.  Since
	$\lambda$ is a regular cardinal and there are fewer than
	$\kappa$ simplices in $A$ (where $\kappa <\lambda$), there
	exists $X_{\alpha}$ with $\alpha <\lambda$ such that the restriction
	of $f$ to every simplex of $A$ factors through $X_{\alpha}$.
	But this implies that $f$ factors through $X_{\alpha}$, and
	therefore we get the surjectivity.
\end{proof}

\begin{lem}
		\label{lemma2.2.injcofs=>eff-monos}
	Consider the category 
	$\simppre$ of simplicial presheaves
	on the smooth Nisnevich site $\nissite$ equipped
	with the injective model structure.
	Then all the cofibrations in $\simppre$
	are effective monomorphisms.
\end{lem}
\begin{proof}
	A map $i:A\rightarrow B$
	is an effective monomorphism if and only if for every $U\in (\smoothS)$, $n\geq 0$
	the induced map $i_{n}(U):A_{n}(U)\rightarrow B_{n}(U)$ is an effective monomorphism
	of sets, this is true since all small limits and colimits are computed termwise.
	Now in the injective model structure for
	$\simppre$ the class of cofibrations coincides
	with the class of monomorphisms.
	But this implies that all the cofibrations
	are effective monomorphisms in $\simppre$,
	since in the category of sets any injective map
	is an effective monomorphism (see remark \ref{rmk1.1.2.sets-effec=inj}).
\end{proof}

	The next proposition is an unpublished result due to Hirschhorn, which also
	appears in \cite[lemma 1.5]{MR2197578}, nevertheless the proof given here
	is slightly different since it also handles the case of a relative $I$-cell
	complex, which is necessary according to Hirschhorn's definition of compactness
	(see definition \ref{def.1.1.2.cell-mod-cats}).

\begin{prop}
		\label{prop2.2.compactness-injmod}
	Let $I$ be the \emph{set} of generating
	cofibrations for the injective model structure
	in the category of simplicial presheaves $\simppre$
	(see theorem \ref{thm.2.1.1.glob-inj-mod}).
	The domains and codomains of the maps in $I$
	are compact relative to $I$.
\end{prop}
\begin{proof}	
	Let $\mu$ be the cardinal of the set $S_{I}$  of simplices
	corresponding to all the domains and codomains of the maps in $I$, i.e.
	$$S_{I}=\coprod _{(i:A\rightarrow B)\in I}\; \; \coprod _{U\in (\smoothS),n\geq 0}A_{n}(U)\sqcup B_{n}(U)
	$$
	and let $\kappa$ be the successor cardinal of $\mu$.
	Since $\kappa$ is a successor cardinal, we have that
	it is a regular cardinal (see \cite[proposition 10.1.14]{MR1944041}).
	
	If $X$ is a presented $I$-cell complex with presentation ordinal
	$\lambda$,
	$$\presabscellcplx
	$$
	we claim that
	every cell $e$ of $X$ is contained in a subcomplex $X_{e}$ of $X$ of size
	less than $\kappa$.  This follows from a transfinite induction
	argument over the presentation ordinal of $e$
	(see definition \ref{def1.2.1.size-cell-complex}).  If the presentation
	ordinal of $e$ is $0$, then the cell $e$ defines a subcomplex of $X$ of
	size $1$, this gets the induction started.  Now assume that
	the result holds for every cell of presentation ordinal less than $\beta<\lambda$, and
	consider an arbitrary cell $e$ of presentation ordinal $\beta$.  The attaching
	map $h_{e}$ of this cell has image contained in the union of fewer than
	$\kappa$ simplices $\{s^{e}\}$ of $X$ (since the domain of $h_{e}$ is also a domain for a map in $I$), 
	now
	each such simplex $s^{e}$ is contained in a cell $e_{s}$
	of presentation ordinal less than $\beta$ and the induction hypothesis implies that
	each such cell $e_{s}$ is contained
	in a subcomplex $X_{s}$ of size less than $\kappa$, thus taking the union
	of all these subcomplexes $X_{s}$ (which is possible by corollary \ref{cor.1.1.2.existence-union-cellcplxs}
	since all the $I$-cells are monomorphisms in this case)
	we get a subcomplex $X_{e}'$ of size less than $\kappa$
	(since $\kappa$ is regular) which contains the image of the attaching map $h_{e}$.
	Therefore if we define $X_{e}$ as the subcomplex obtained from $X_{e}'$
	after attaching the cell $e$ via $h_{e}$, we get a subcomplex of size less than $\kappa$
	containing the given cell $e$.  This proves the claim.
	
	Now if $A$ is a simplicial presheaf on $(\smoothS)$ which is a domain
	or codomain of a map in $I$, we have that
	$A$ has less than $\kappa$ simplices.  Consider a map $j:A\rightarrow X$
	where $X$ is a presented $I$-cell complex,
	$$\presabscellcplx
	$$
	then the image of $j$
	has less than $\kappa$ simplices $\{s_{j}\}$, each such simplex $s_{j}$ is contained in some
	cell $e_{s}$ of $X$ which by the previous argument is contained in a subcomplex $X_{s}$ of $X$
	of size less than $\kappa$.  We take now the union of all these subcomplexes $X_{s}$
	to get a subcomplex $X_{j}$ of $X$ of size less than $\kappa$ (since $\kappa$ is regular)
	which contains the image of $j$.  Therefore $j$ factors through the subcomplex
	$X_{j}$ which has size less than $\kappa$.
	
	Finally, we consider a relative cell complex $f:X\rightarrow Y$,
	$$\prescellcplxtwo
	$$
	Take any map $j:A\rightarrow Y$ where $A$ is a domain or codomain of a map in $I$.
	Since all the inclusions are $I$-cells for the injective model structure,
	we have that $X$ is a  cell complex, 
	$$i:\emptyset\rightarrow X, \; \emptyset=X_{0}\rightarrow X_{1} \rightarrow \cdots \rightarrow X_{\beta}
						\rightarrow \cdots (\beta <\nu),\; \{T^{\beta},e^{\beta},h^{\beta}\}_{\beta <\nu}
	$$
	Combining this presentation of $X$ with the presentation of $f$ we
	get a presentation for $Y$ as a cell complex, where
	$X$ is a subcomplex.
	The previous argument shows that the image of $j$ is contained in a subcomplex
	$W'$ of $Y$ where the size of $W'$ is less than $\kappa$.
	Taking the union of $W'$ and $X$ we get a subcomplex $X_{f}$ of $f$ having the
	same size as $W'$ (as a subcomplex of $f$) which contains the image of $j$.  Therefore $j$ factors through
	$X_{f}$ where $X_{f}$ is a subcomplex of $f$ of size less than $\kappa$, and this
	shows that $A$ is $\kappa$-compact relative to $I$.
\end{proof}

	Finally we are ready to prove Hirschhorn's cellularity theorem.
	
\begin{thm}
		\label{thm2.2.injmod=>cellular}
	The category $\simppre$ of simplicial presheaves
	on the smooth Nisnevich site $\nissite$
	is a cellular model category when it is
	equipped with the injective model structure, the sets
	of generating cofibrations and  generating trivial cofibrations
	are the ones considered in theorem \ref{thm.2.1.1.glob-inj-mod}.
\end{thm}
\begin{proof}
	We have to check that the conditions 
	(\ref{def.1.1.2.cell-mod-cats.a})-(\ref{def.1.1.2.cell-mod-cats.d}) of
	definition \ref{def.1.1.2.cell-mod-cats} hold.
	(\ref{def.1.1.2.cell-mod-cats.a}) follows from theorem \ref{thm.2.1.1.glob-inj-mod}
	which shows that the injective model structure is cofibrantly generated.
	(\ref{def.1.1.2.cell-mod-cats.b}) follows from proposition
	\ref{prop2.2.compactness-injmod} and (\ref{def.1.1.2.cell-mod-cats.c}) follows from
	lemma \ref{lem2.2.all-obj-small} which says that every simplicial presheaf is small.
	Finally (\ref{def.1.1.2.cell-mod-cats.d}) follows from
	lemma \ref{lemma2.2.injcofs=>eff-monos}.
\end{proof}

	Theorem \ref{thm2.2.injmod=>cellular} will be used
	to show that the category $\Tspectra$ of $T$-spectra on $\smoothS$
	equipped with the motivic stable model structure
	is cellular
	(see theorem \ref{thm.2.5.cellularity-motivic-stable-str}).
	This will allow us to apply all the 
	localization technology of Hirschhorn \cite{MR1944041}
	to construct new model structures for 
	$\Tspectra$.

\end{section}

\begin{section}{The Motivic Model Structure}

	Let $\affineS$ be the affine line over $S$.
	Consider the following set of maps
	$$\mathcal V _{M} =\{ \pi_{U}:U\times \affineS \rightarrow U\; |\; U\in (\smoothS)\}
	$$
	In \cite{MR1813224} Morel and Voevodsky show in particular that
	for simplicial sheaves  on $\nissite$ the left Bousfield localization
	for the injective model structure
	with respect to $\mathcal V _{M}$ exists, and furthermore
	they show it is a proper simplicial model structure.
	Their work was extended to the case of simplicial presheaves
	by Jardine in \cite[section 1]{MR1787949}.
	Following Jardine we call this localized model structure
	the 	
	\emph{motivic model structure on $\simppre$}.
	
\begin{thm}[Morel-Voevodsky, Jardine]
		\label{thm2.3.motmodstruc}
	Consider the category  of simplicial presheaves
	on the smooth Nisnevich site $\nissite$
	equipped with the injective model
	structure.  Then  the left Bousfield localization
	(see section \ref{subsec-Bousloc})
	with respect to the set of maps $\mathcal V _{M}$
	defined above exists.  This model
	structure will be called \emph{motivic}, and the category
	$\simppre$ equipped with the motivic model structure will
	be denoted by $\motivic$.
	Furthermore $\motivic$ is a
	proper and simplicial model category.
\end{thm}
\begin{proof}
	We refer the reader to \cite[theorem 1.1]{MR1787949}.
\end{proof}
	
	The following theorem gives
	explicit sets of generating cofibrations
	and trivial cofibrations for $\motivic$;
	and it also shows that with this choice
	of generators, $\motivic$ has the structure
	of a cellular model category.	
	In \cite[corollary 1.6]{MR2197578} it is also proved
	that  $\motivic$
	is cellular.  
		
\begin{thm}
		\label{thm1.3.motmodstru-cellular}
	$\motivic$ is a cellular model category, where the set $I_{M}$ of
	generating cofibrations and the set $J_{M}$
	of generating trivial cofibrations are defined as follows:
	\begin{enumerate}
		\item	\label{thm1.3.motmodstru-cellular.a}$I_{M}=I$ where $I$ is the set of generating cofibrations
					for the injective model structure on $\simppre$ (see theorem \ref{thm.2.1.1.glob-inj-mod}).
		\item	\label{thm1.3.motmodstru-cellular.b}$J_{M}=\{j:A\rightarrow B\}$ such that:
					\begin{enumerate}
						\item	$j$  is a monomorphism.
						\item	$j$ is a $\mathcal V _{M}$-local equivalence.
						\item	The size of  $B$ as an $I$-cell complex
									(see definition \ref{def1.2.1.size-cell-complex}) 
									is less than  $\kappa$,
									where $\kappa$ is the cardinal defined 
									by Hirschhorn in \cite[definition 4.5.3]{MR1944041}.
					\end{enumerate}
	\end{enumerate}
\end{thm}
\begin{proof}
	By theorem \ref{thm2.2.injmod=>cellular} the injective model structure
	on $\simppre$ is cellular.  Therefore
	we can use Hirschhorn's techniques (see section \ref{subsec-Bousloc})
	to construct the left Bousfield localization with respect to the set
	of maps $\mathcal V _{M}$ defined above.  This model structure is identical to
	the motivic model structure of theorem \ref{thm2.3.motmodstruc} since
	both are left Bousfield localizations with respect to the same set of maps.  Now using
	\cite[theorem 4.1.1]{MR1944041} we have that
	the motivic model structure is cellular.  So it only remains to
	show that the sets of generating cofibrations and trivial cofibrations
	are the ones described above. For the set of generating cofibrations
	it is clear.  Theorem  4.1.1 in \cite{MR1944041} implies
	that the generating trivial cofibrations are the maps
	$j:A\rightarrow B$ where $j$ is an inclusion of $I$-cell complexes
	and a $\mathcal V _{M}$-local equivalence, and the size of $B$ is  
	less than $\kappa$.
	The result follows from the fact that 
	in the injective model structure for $\simppre$, $I$-cell is just
	the class of monomorphisms and that every object in $\simppre$ is an
	$I$-cell complex (see remark \ref{rmks.1.1.class-cofs-inj-mod}).
\end{proof}
			
	Following Jardine we say that a simplicial
	presheaf $X$ is 
	\emph{motivic fibrant} if $X$ is
	$\mathcal V _{M}$-local.
	
\begin{prop}
		\label{prop2.4.mot-fib-obj}
	The following conditions are equivalent:
	\begin{enumerate}
		\item	\label{prop2.4.mot-fib-obj.a}$X$ is motivic fibrant.
		\item	\label{prop2.4.mot-fib-obj.b}$X$ is fibrant in the injective structure and
					for every $U$ in $\smoothS$ the map induced by  $U\times \affineS \rightarrow U$ 
					$$\xymatrix{Map(U,X) \ar[r] & Map(U\times \affineS, X)}
					$$
					is a weak equivalence of simplicial sets.
		\item	\label{prop2.4.mot-fib-obj.c}$X$ is fibrant in the injective structure and
					for every $U$ in $\smoothS$ the map induced by $U\times \ast \rightarrow U\times \affineS$
					$$\xymatrix{Map(U\times \affineS,X) \ar[r] & Map(U\times \ast,X)}
					$$
					is a weak equivalence of simplicial sets, where
					$\ast \rightarrow \affineS$ is any rational point for $\affineS$.
		\item	\label{prop2.4.mot-fib-obj.d}$X$ is fibrant in the injective structure and
					for every $U$ in $\smoothS$ the map induced by $U\times \ast \rightarrow U\times \affineS$
					$$\xymatrix{Map(U\times \affineS,X) \ar[r] & Map(U\times \ast,X)}
					$$
					is a trivial fibration of Kan complexes, where
					$\ast \rightarrow \affineS$ is any rational point for $\affineS$.
		\item	\label{prop2.4.mot-fib-obj.e}$X$ is fibrant in the injective structure and
					for every $U$ in $\smoothS$ the map induced by $U\times \ast \rightarrow U\times \affineS$
					$$\xymatrix{\inthompre(U\times \affineS,X) \ar[r]& \inthompre(U\times \ast,X)}
					$$
					is a trivial fibration between fibrant objects in the injective model structure for
					$\simppre$, where $\ast \rightarrow \affineS$ is any rational point for $\affineS$.
	\end{enumerate}
\end{prop}
\begin{proof}
	The claim that (\ref{prop2.4.mot-fib-obj.a}) and (\ref{prop2.4.mot-fib-obj.b}) are equivalent
	follows from the definition of $\mathcal V _{M}$-local and the fact that
	every simplicial presheaf is cofibrant in the injective model structure.
	(\ref{prop2.4.mot-fib-obj.b}) and (\ref{prop2.4.mot-fib-obj.c}) are equivalent
	since the following diagram is commutative
	$$\xymatrix{U\cong U\times \ast \ar[r] \ar[d]_-{id} & U\times \affineS \ar[dl]\\
							U &}
	$$
	and weak equivalences of simplicial sets satisfy the
	two out of three property.
	(\ref{prop2.4.mot-fib-obj.c}) and (\ref{prop2.4.mot-fib-obj.d})
	are equivalent since the injective structure is in particular
	a simplicial model category. 
	
	(\ref{prop2.4.mot-fib-obj.d}) $\Rightarrow$
	(\ref{prop2.4.mot-fib-obj.e}):  Since $\simppre$
	equipped with the injective model structure is a symmetric monoidal
	model category we have that
	$$\xymatrix{\inthompre(U\times \affineS,X) \ar[r]^-{p} & \inthompre(U\times \ast,X)}
	$$
	is a fibration between fibrant objects in the injective structure.
	It only remains to show that $p$ is a weak equivalence in
	the injective model structure.
	Lemma \ref{lem2.1.criterion-weakequivs} implies that it is enough to show
	that
	$$\xymatrix{\inthompre(U\times \affineS,X)(V) \ar[r]^-{p(V)} & \inthompre(U\times \ast,X)(V)}
	$$
	is a weak equivalence of simplicial sets for every $V$ in $(\smoothS)$.
	But for any simplicial presheaf $Z$ we have a natural isomorphism of simplicial
	sets $Z(V)\cong Map(V,Z)$, therefore $p(V)$ is just
	$$\xymatrix{Map(V,\inthompre(U\times \affineS,X)) \ar[r]^-{p(V)} & Map(V,\inthompre(U\times \ast,X))}
	$$
	Now using the enriched adjuntions of \ref{prop.2.1.adjunction-simplicial-internal},
	$p(V)$ becomes
	$$\xymatrix{Map(V\times U\times \affineS,X) \ar[r]^-{p(V)}& Map(V\times U\times \ast,X)}
	$$
	and by hypothesis we know that this map is a weak equivalence of simplicial sets.
	
	(\ref{prop2.4.mot-fib-obj.e}) $\Rightarrow$
	(\ref{prop2.4.mot-fib-obj.d}):  Since the injective model structure
	is simplicial, we have that
	$$\xymatrix{Map(U\times \affineS,X) \ar[r]^-{f} & Map(U\times \ast,X)}
	$$
	is a fibration between Kan complexes.  So it
	only remains to show that $f$ is a weak equivalence of simplicial sets.
	By hypothesis we have that
	$$\xymatrix{\inthompre(U\times \affineS,X) \ar[r]^-{p}& \inthompre(U\times \ast,X)}
	$$
	is a trivial fibration between fibrant objects in the injective model structure.
	Lemma \ref{lem2.1.criterion-weakequivs} implies that if we take global
	sections at $S$:
	$$\xymatrix{\inthompre(U\times \affineS,X)(S) \ar[r]^-{p(S)}& \inthompre(U\times \ast,X)(S)}
	$$
	we get a weak equivalence of simplicial sets.
	But $p(S)$ is natural isomorphic to
	$$\xymatrix{Map(U\times \affineS,X) \ar[r]^-{f} & Map(U\times \ast,X)}
	$$
	so this proves the result.
\end{proof}

\begin{prop}
		\label{prop.2.3.motivic-fibrant-stable-inthoms}
	Let $X$ be a motivic fibrant simplicial presheaf
	on the smooth Nisnevich site $\nissite$.
	Then for any $Y$ in $\simppre$, the simplicial presheaf
	$\inthompre(Y,X)$ is also motivic fibrant.
\end{prop}
\begin{proof}
	Since the injective structure is a
	symmetric monoidal model category (see lemma \ref{lem2.1.monoidalmodcat})
	we have that $\inthompre (Y,X)$
	is a fibrant object for the injective model
	structure.  Proposition \ref{prop2.4.mot-fib-obj}(\ref{prop2.4.mot-fib-obj.e})
	implies that for every $U$ in $\smoothS$, the map
	$$\xymatrix{\inthompre(U\times \affineS,X) \ar[r]^-{p}& \inthompre(U\times \ast,X)}
	$$
	is a trivial fibration between fibrant objects in the injective model structure
	for $\simppre$, and since the injective model structure is simplicial
	we have that
	$$\xymatrix{Map(Y,\inthompre(U\times \affineS,X)) \ar[r]^-{p_{\ast}}& Map(Y,\inthompre(U\times \ast,X))}
	$$
	is a trivial fibration of Kan complexes.
	Now using the enriched adjunctions of proposition
	\ref{prop.2.1.adjunction-simplicial-internal}, $p_{*}$ becomes
	$$\xymatrix{Map(Y\times U\times \affineS,X) \ar[r]^-{p_{\ast}}& Map(Y\times U\times \ast,X)}
	$$
	and finally
	$$\xymatrix{Map(U\times \affineS,\inthompre(Y,X)) \ar[r]^-{p_{*}}& Map(U\times \ast,\inthompre(Y,X))}
	$$
	therefore proposition \ref{prop2.4.mot-fib-obj}(\ref{prop2.4.mot-fib-obj.d})
	implies that $\inthompre(Y,X)$ is motivic fibrant since $p_{*}$ is a trivial
	fibration of Kan complexes for every $U$ in $(\smoothS)$.
\end{proof}

	Since the motivic and the injective model
	structures have the same class of cofibrations and the
	same set of generating cofibrations, it follows that
	the cofibrations for the motivic model structure also
	have the properties described in 
	remark \ref{rmks.1.1.class-cofs-inj-mod}.

\begin{cor}
		\label{cor.2.3.motmodstr-symmmon}
	$\motivic$
	is a symmetric monoidal model category.
\end{cor}
\begin{proof}
	The cofibrations for the motivic and injective
	model structures coincide, therefore it only remains
	to show that if we have two cofibrations
	$i:A\rightarrow B$, $j:C\rightarrow D$ where $j$
	is a motivic weak equivalence, the induced map
	$$\xymatrix{i\square j:A\times D\coprod_{A\times C}B\times C \rightarrow B\times D}
	$$
	is a trivial cofibration in $\motivic$.
	Since every simplicial presheaf is cofibrant in the motivic
	model structure, it is enough to prove the following
	claim:  For any trivial cofibration $j:C\rightarrow D$ in $\motivic$
	and for any simplicial presheaf $A$, the induced map
	$j\times id:C\times A\rightarrow D\times A$ is a trivial cofibration
	in $\motivic$.
	Since the injective model structure for $\simppre$ is a symmetric monoidal
	model category (see lemma \ref{lem2.1.monoidalmodcat})
	we have that $j\times id$ is a cofibration, so it only
	remains to show that it is a weak equivalence in $\motivic$.
	Let $X$ be any motivic fibrant simplicial presheaf, 
	proposition \ref{prop.2.3.motivic-fibrant-stable-inthoms} implies
	that $\inthompre(A,X)$ is also motivic fibrant, therefore since
	$j$ is a weak equivalence in $\motivic$, the map
	$$\xymatrix{Map(D,\inthompre(A,X)) \ar[r]^-{j^{\ast}} & Map(C,\inthompre(A,X))}
	$$
	is a weak equivalence of simplicial sets.  Now
	using the enriched adjunctions of proposition \ref{prop.2.1.adjunction-simplicial-internal},
	$j^{*}$ becomes
	$$\xymatrix{Map(D\times A,X) \ar[r]^-{j^{\ast}} & Map(C\times A,X)}
	$$
	and this implies that $j\times id:C\times A \rightarrow D\times A$
	is a weak equivalence in $\motivic$, hence the result follows.
\end{proof}
	
\begin{rmk}
		\label{rmk.2.3.pointedmonstrpresh}
	Proposition \ref{prop.pointed-mon-cats} implies that
	the associated pointed category $\simpprepointed$ of
	pointed simplicial presheaves is also closed symmetric
	monoidal, we denote by $X\wedge Y$ the functor giving
	the monoidal structure, and by
	$\inthomprepointed (X,Y)$ the adjunction of two variables.
\end{rmk}	
	
\begin{prop}
		\label{prop.2.3.cell-prop-pointedmotcat}
	Let $\pointedmotivic$ denote the pointed
	category associated to $\motivic$ (see remark \ref{rmk1.1.triv-conseq}), i.e.
	the category with pointed simplicial presheaves
	as objects and base point preserving
	maps.		The model structure on
	$\pointedmotivic$ induced from the
	model structure on $\motivic$ is
	cellular, proper, simplicial and symmetric monoidal.
	Furthermore, $\pointedmotivic$ is a
	$\simpsets _{\ast}$-model category
	(see definition \ref{def.module-modcats}).
	The sets $I_{M_{\ast}}$, $J_{M_{\ast}}$ of generating cofibrations and trivial
	cofibrations respectively, are defined as follows:
	\begin{enumerate}
		\item	$$I_{M_{\ast}}=\{ i_{+}:Y_{+}\hookrightarrow (\Delta ^{n}_{U})_{+}\}
					$$
					where $i:Y\hookrightarrow \Delta ^{n}_{U}$ is a generating
					cofibration for $\motivic$
					(see theorem \ref{thm1.3.motmodstru-cellular}(\ref{thm1.3.motmodstru-cellular.a})).
		\item	$$J_{M_{\ast}}=\{ j_{+}:A_{+}\rightarrow B_{+}\}
					$$
					where $j$ is a map in the set $J$ defined in theorem 
					\ref{thm1.3.motmodstru-cellular}(\ref{thm1.3.motmodstru-cellular.b}),
					i.e. $j$ is a generating
					trivial cofibration for $\motivic$.
	\end{enumerate}
\end{prop}
\begin{proof}
	Theorems \ref{thm2.3.motmodstruc},
	\ref{thm1.3.motmodstru-cellular}
	together with corollary \ref{cor.2.3.motmodstr-symmmon} 	
	imply that $\motivic$
	is cellular, proper, simplicial and symmetric monoidal.
	Then theorem
	\ref{thm.1.5.pointedmodelcat-cellular} and
	theorem
	\ref{thm1.6.pointedmodelcat-proper}	
	imply that the associated pointed 
	category $\pointedmotivic$ with the induced  model structure
	is cellular (with the sets of generating cofibrations and trivial
	cofibrations as defined above) and proper.
	Now proposition \ref{prop.1.10.M-modelcats-descend-pointed}
	implies that $\pointedmotivic$ is a $\simpsets _{*}$-model category,
	and this induces a simplicial model structure in $\pointedmotivic$,
	since the natural functor $\simpsets \rightarrow \simpsets _{\ast}$
	which adds a disjoint base point is a left Quillen monoidal functor.
	Finally proposition \ref{prop.pointed-mon-cats}
	implies that $\pointedmotivic$
	is symmetric monoidal.
\end{proof}

\begin{defi}[cf. \cite{MR1787949}]
		\label{def.2.3.motivic-flasque}
	Let $X\in \motivic$ be a simplicial presheaf.
	We say that $X$ is \emph{motivic flasque} if:
	\begin{enumerate}
		\item	$X$ is flasque (see definition \ref{def.2.1.flasque-simplicial-presheaves}).
		\item	For every $U\in \smoothS$ the map
						$$\xymatrix{X(U)\cong Map(U,X) \ar[r]& Map(U\times \affineS,X)\cong X(U\times \affineS)}
						$$
					induced by the projection $U\times \affineS \rightarrow U$
					is a weak equivalence of simplicial sets.
	\end{enumerate}
\end{defi}

\begin{rmk}
		\label{rmk.2.3.stabilitymotivicflasque}
	\begin{enumerate}
		\item	\label{rmk.2.3.stabilitymotivicflasque.a}The class of motivic
					flasque simplicial presheaves is closed under filtered colimits.
		\item	\label{rmk.2.3.stabilitymotivicflasque.b}The functors $\phi ^{-1}$
					and $\phi _{\ast}$ (see definition \ref{def.2.1.basechangefunctor}) preserve motivic flasque simplicial
					presheaves.
		\item \label{rmk.2.3.stabilitymotivicflasque.c}If $X$ is fibrant in the motivic model structure for
					$\simpprepointed$ then $X$ is also motivic flasque.
	\end{enumerate}
\end{rmk}

\begin{defi}[cf. \cite{MR1787949}]
		\label{def.2.3.compactness}
	Let $X\in \pointedmotivic$ be a pointed simplicial presheaf.
	We say that $X$ is \emph{compact} if:
	\begin{enumerate}
		\item	All inductive systems $Z_{1}\rightarrow Z_{2}\rightarrow \cdots$
					of pointed simplicial presheaves induce isomorphisms
						$$\inthomprepointed (X,\varinjlim Z_{i})\cong \varinjlim \inthomprepointed (X,Z_{i})
						$$
		\item	If $Z$ is motivic flasque, then $\inthomprepointed (X,Z)$ is also
					motivic flasque.
		\item	The functor 
						$$\xymatrix{\inthomprepointed (X,-):\pointedmotivic \ar[r]& \pointedmotivic}
						$$
					takes sectionwise weak equivalences of motivic flasque pointed simplicial
					presheaves to sectionwise weak equivalences.
	\end{enumerate}
\end{defi}

\begin{prop}
		\label{prop.2.3.compact=>commutes-colimits.b}
	Let $X\in \pointedmotivic$ be a pointed simplicial presheaf,
	and let 
		$$\xymatrix{Z_{1} \ar[r]& Z_{2} \ar[r]& \cdots}
		$$
	be an inductive system of pointed simplicial presheaves.
	If $X$ is compact in the sense of Jardine (see definition \ref{def.2.3.compactness})
	then:
		$$[X,\varinjlim Z_{i}]\cong \varinjlim [X,Z_{i}]
		$$
	where $[-,-]$ denotes the set of maps in the homotopy category associated
	to $\pointedmotivic$.
\end{prop}
\begin{proof}
	Let $R$ denote a functorial fibrant replacement in $\pointedmotivic$, 
	such that the natural map
	$R_{Y}:Y\rightarrow RY$ is always a trivial cofibration.
	Consider the following
	commutative diagram:
		$$\xymatrix{Z_{1} \ar[r] \ar[d]& Z_{2}\ar[r] \ar[d]& \cdots \ar[r] & 
								\varinjlim Z_{i}\ar[r]^-{j} \ar[d]_-{i}& R(\varinjlim Z_{i}) \ar[d]^-{i_{R}}\\
								RZ_{1} \ar[r] & RZ_{2}\ar[r] & \cdots \ar[r] & 
								\varinjlim RZ_{i}\ar[r]_-{j_{R}} & R(\varinjlim RZ_{i})}
		$$
	Since all the maps $Z_{i}\rightarrow RZ_{i}$ are trivial cofibrations,
	it follows that the induced map $i:\varinjlim Z_{i}\rightarrow \varinjlim RZ_{i}$
	is also a trivial cofibration.  Therefore:
		\begin{equation}
				\label{equation.2.3.compact-colimit1}
			[X,\varinjlim Z_{i}]\cong [X,\varinjlim RZ_{i}]\cong [X,R(\varinjlim RZ_{i})]
		\end{equation}
	
	We have that the pointed simplicial presheaves $RZ_{i}$ are motivic fibrant, then
	remark \ref{rmk.2.1.inj-fib=>flasque}(\ref{rmk.2.1.inj-fib=>flasque.c}) implies that
	$\varinjlim RZ_{i}$ satisfies the B.G. property.
	Therefore using theorem \ref{thm.2.1.Nisnevich-descent}
	we get that the map $j_{R}:\varinjlim RZ_{i}\rightarrow R(\varinjlim RZ_{i})$
	is a sectionwise weak equivalence.  On the other hand
	$\varinjlim RZ_{i}$ and $R(\varinjlim RZ_{i})$ are both motivic flasque 
	(see remark \ref{rmk.2.3.stabilitymotivicflasque}(\ref{rmk.2.3.stabilitymotivicflasque.a})),
	and since $X$ is compact we have that
		$$\xymatrix{\varinjlim \inthomprepointed (X,RZ_{i})\cong \inthomprepointed (X,\varinjlim RZ_{i})
			\ar[r]& \inthomprepointed (X,R(\varinjlim RZ_{i}))}
			$$			
	is a sectionwise weak equivalence of simplicial sets.  Taking global sections at $S$ we get
	the following weak equivalence of simplicial sets:
		$$\varinjlim Map(X,RZ_{i})\rightarrow Map(X,R(\varinjlim Z_{i}))
		$$
	Therefore
		\begin{eqnarray}
				\label{equation.2.3.compact-colimit2}
			[X,R(\varinjlim RZ_{i})]&\cong & \pi _{0}Map(X,R(\varinjlim RZ_{i}))\\ &\cong & 
			\pi _{0}\varinjlim Map(X,RZ_{i})\cong \varinjlim \pi _{0}Map(X,RZ_{i}) \nonumber
		\end{eqnarray}
	On the other hand:
		\begin{equation}
				\label{equation.2.3.compact-colimit3}
			\varinjlim \pi _{0}Map(X,RZ_{i}) \cong \varinjlim [X,RZ_{i}]
		\end{equation}
	Hence equations (\ref{equation.2.3.compact-colimit1}), (\ref{equation.2.3.compact-colimit2})
	and (\ref{equation.2.3.compact-colimit3}) imply that
		$$[X,\varinjlim Z_{i}]\cong [X, R(\varinjlim RZ_{i})]\cong \varinjlim [X,RZ_{i}]
			\cong \varinjlim [X, Z_{i}]
		$$		
	as we wanted.
\end{proof}

\begin{defi}
		\label{def.2.3.loopfunctor}
	Let $A \in \pointedmotivic$ be an arbitrary pointed simplicial presheaf.
	We define the \emph{functor of $A$-loops}
	as follows:
	$$\xymatrix@R=.5ex{\Omega _{A}:\pointedmotivic \ar[r]& \pointedmotivic \\
											X \ar@{|->}[r]& \inthomprepointed (A,X)}
	$$
\end{defi}

\begin{rmk}
		\label{rmk.2.3.adjointforloops}
	\begin{enumerate}
		\item	The functor of $A$-loops $\Omega _{A}$ has
					a left adjoint given by smash product with $A$,
					i.e.
					$$\xymatrix@R=.5pt{-\wedge A:\pointedmotivic \ar[r]& \pointedmotivic \\
														X \ar@{|->}[r]& X\wedge A}
					$$
		\item The adjunction 
					$$\xymatrix{(-\wedge A,\Omega _{A},\varphi):\pointedmotivic \ar[r]& \pointedmotivic}
					$$
					is a Quillen adjunction.
	\end{enumerate}
\end{rmk}


\end{section}
\begin{section}{The Motivic Stable Model Structure}
		\label{section.motivic.spectra}
		
	In \cite{MR1787949} Jardine constructs a stable model structure
	for the category of $T$-spectra on $\smoothS$.
	In order to define this stable model structure,
	he constructs two auxiliary model structures
	called projective and injective.  In this section
	we recall Jardine's definitions for these three model structures
	on the category of $T$-spectra.

	Let $S^{1}$ denote the constant presheaf associated to
	the pointed simplicial set
	$\Delta ^{1}/\partial \Delta ^{1}$, let $S^{n}$ denote $S^{1}\wedge\cdots \wedge S^{1}$
	($n$-factors) and let
	$\multgroup$ denote the multiplicative group
	over the base scheme $S$, i.e. $\multgroup =\affineS -\{0\}$
	pointed by the unit $e$ for the group operation.
	Let $T=\tee$.

\begin{prop}
		\label{prop.2.4.T=compact}
	\begin{enumerate}
		\item	$T=\tee$ is compact in the sense of Jardine (see definition \ref{def.2.3.compactness}).
		\item	Consider $U\in \smoothS$ and $r,s\geq 0$.  Then the pointed simplicial presheaf
					$S^{r}\wedge \gm ^{s}\wedge U_{+}$ is compact in the sense of Jardine, where
					$\gm ^{s}$ denotes $\gm \wedge \cdots \wedge \gm$ ($s$-factors).
	\end{enumerate}
\end{prop}
\begin{proof}
	Follows immediately from \cite[lemma 2.2]{MR1787949}.
\end{proof}
	
\begin{defi}
		\label{defi2.4.T-spectra}
	\begin{enumerate}
		\item	A \emph{$T$-spectrum} $X$ is a collection
					of pointed simplicial presheaves  $(X^{n})_{n\geq 0}$
					on the smooth Nisnevich site $\smoothS$, together
					with \emph{bonding maps}
					$$\xymatrix{T\wedge X^{n} \ar[r]^-{\sigma ^{n}}& X^{n+1}}
					$$
		\item	A map $f:X\rightarrow Y$ of $T$-spectra
					is a collection of maps
					$$\xymatrix{X^{n} \ar[r]^-{f^{n}} & Y^{n}}
					$$
					in $\pointedmotivic$ which are compatible
					with the bonding maps, i.e.
					the following diagram:
					$$\xymatrix{T\wedge X^{n} \ar[rr]^-{id\wedge f^{n}} \ar[d]_-{\sigma^{n}}&& 
											T\wedge Y^{n} \ar[d]^-{\sigma ^{n}}\\
											X^{n+1} \ar[rr]_-{f^{n+1}}&& Y^{n+1}}
					$$
					commutes for all $n\geq 0$.
		\item	With the previous definitions we get a category
					called the category of \emph{$T$-spectra} which
					will be denoted by $\Tspectra$.
	\end{enumerate}
\end{defi}

	The category of $T$-spectra has a natural simplicial
	structure induced from the one
	on pointed simplicial presheaves.
	
	Given a $T$-spectrum $X$,	the tensor objects are defined as follows:
	$$\xymatrix@R=.5pt{X\wedge -:\simpsets \ar[r]& \Tspectra \\
										 K \ar@{|->}[r]& X\wedge K}
	$$
	where $(X\wedge K)^{n}=X^{n}\wedge K_{+}$ and the bonding maps are 
	$$\xymatrix{T\wedge (X^{n}\wedge K_{+}) \ar[r]^-{\cong}& 
							(T\wedge X^{n})\wedge K_{+} \ar[rr]^-{\sigma^{n}\wedge id_{K_{+}}} && X^{n+1}\wedge K_{+}}
	$$
	The simplicial functor in two variables is:
	$$\xymatrix@R=.5pt{Map(-,-):(\Tspectra)^{op}\times \Tspectra \ar[r] & \simpsets \\
										(X,Y) \ar@{|->}[r]& Map(X,Y)}
	$$
	where $Map(X,Y)_{n}=\Hom _{\Tspectra}(X\wedge \Delta ^{n}_{+},Y)$,
	and finally for any $T$-spectrum $Y$ we have the following functor
	$$\xymatrix@R=.5pt{Y^{-}:\simpsets \ar[r]& (\Tspectra )^{op}\\
										K \ar@{|->}[r] & Y^{K}}
	$$
	where $(Y^{K})^{n}=(Y^{n})^{K_{+}}$ with bonding maps
	$$\xymatrix{T\wedge(Y^{n})^{K_{+}} \ar[r]^-{\alpha}& (T\wedge Y^{n})^{K_{+}} \ar[r]^-{(\sigma ^{n})_{*}}& (Y^{n+1})^{K_{+}}}
	$$
	where for $U\in (\smoothS)$, $\alpha(U)$ is adjoint to
	$$\xymatrix{T(U)\wedge (Y^{n}(U))^{K_{+}}\wedge K_{+} \ar[rr]^-{id_{T(U)}\wedge ev_{K_{+}}} && T(U)\wedge Y^{n}(U)}
	$$

\begin{rmk}
		\label{rmk.2.4.simp-simppointed-structures-Tspectra}
	\begin{enumerate}
		\item \label{rmk.2.4.simp-simppointed-structures-Tspectra.a}In fact there exists
					an adjunction of two variables (see definition \ref{def.adj-two-vars}):
						$$\xymatrix{-\wedge -:\Tspectra \times \simpsets _{*} \ar[r]& \Tspectra}
						$$
					which induces
					the simplicial structure for $T$-spectra described above 
					via the
					monoidal functor $\simpsets \rightarrow \simpsets _{\ast}$
					which adds a disjoint base point.
		\item	\label{rmk.2.4.simp-simppointed-structures-Tspectra.b}For any two given
					spectra $X$, $Y$, the simplicial set $Map(X,Y)$ is just
					$Map_{\ast}(X,Y)$ (i.e. the pointed simplicial set coming from the adjunction of two
					variables described above)
					after forgetting its base point 
						$$\xymatrix{\omega _{0}:X\ar[r] & \ast \ar[r] & Y}$$
	\end{enumerate}
\end{rmk}

	We have the following family of shift functors
	between $T$-spectra defined for every $n\in \mathbb Z$
	$$\xymatrix@R=.5pt{s_{n}:\Tspectra \ar[r]& \Tspectra \\
											X \ar@{|->}[r]& X[n]}
	$$
	where $X[n]$ is defined as follows:
	$$(X[n])^{m}=
		\begin{cases}
			\ast & \text{if $m+n<0$.}\\
			X^{m+n} & \text{if  $m+n\geq 0$.}
		\end{cases}
	$$
	with the obvious bonding maps induced by $X$.
	It is clear that  $s_{0}=id$ and
	that for $n\geq 0$, $s_{n}$ is right adjoint to $s_{-n}$, i.e.
	$$\Hom _{\Tspectra}(X,Y[n])\cong \Hom _{\Tspectra}(X[-n],Y)
	$$

	We define the projective model structure as follows.

\begin{defi}
		\label{def.2.4.infinitesuspension+shift}
	Consider the following family of functors 
	from the category of pointed simplicial presheaves
	to the category of $T$-spectra:
	\begin{equation}
				\label{diagram.2.infinite-suspensions}
		\begin{array}{c}
			$$\xymatrix@R=.5pt{F_{n}:\pointedmotivic \ar[r]& \Tspectra\\
							X \ar@{|->}[r]& (\Tsuspension X)[-n]}
			$$
		\end{array}
	\end{equation}
	where $\Tsuspension X$ is defined as follows:
	$$(\Tsuspension X)^{k}=T^{k}\wedge X
	$$
	where the bonding maps are the canonical isomorphisms
	$T\wedge (T^{k}\wedge X)\stackrel{\cong}{\rightarrow}T^{k+1}X$
	and $T^{0}\wedge X$ is just $X$.
\end{defi}
	
	We also have the following evaluation functors from
	the category of $T$-spectra to the category of pointed
	simplicial presheaves:
	$$\xymatrix@R=.5pt{Ev_{n}:\Tspectra \ar[r]& \pointedmotivic \\
										X \ar@{|->}[r]& X^{n}}
	$$
	where $n\geq 0$.
	It is clear that $F_{0}$ is left adjoint to $Ev_{0}$.  
	This implies that for every $n>0$, $F_{n}$ is left adjoint to $Ev_{n}$
	and $F_{-n}$ is left adjoint to $\Tloops ^{n}\circ Ev_{0}$.
			
	We say that a map of $T$-spectra $f:X\rightarrow Y$
	is a \emph{level equivalence} if for every $n\geq 0$,
	$f^{n}:X^{n}\rightarrow Y^{n}$ is a weak equivalence
	in $\pointedmotivic$.
	
	Let $I_{M_{\ast}}$ and $J_{M_{\ast}}$ denote the sets of generating
	cofibrations and trivial cofibrations for $\pointedmotivic$
	(see proposition \ref{prop.2.3.cell-prop-pointedmotcat}).
	
\begin{thm}[Jardine]
		\label{thm2.4.projstabmodstr}
	There exists a cofibrantly generated model structure for the category
	$\Tspectra$ of $T$-spectra with the following choices:
	\begin{enumerate}
		\item The weak equivalences are the level weak equivalences defined above.
		
		\item	The set $I$ of generating cofibrations is 
					$$I=\bigcup _{n\geq 0}F_{n}(I_{M_{\ast}})$$
		
		\item	The set $J$ of generating trivial cofibrations is
					$$J=\bigcup _{n\geq 0}F_{n}(J_{M_{\ast}})$$
	\end{enumerate}
	This model structure will be called the \emph{projective
	model structure} for $T$-spectra.  Furthermore,
	the projective model structure is proper and simplicial.
\end{thm}
\begin{proof}
	We refer the reader to \cite[lemma 2.1]{MR1787949}.
\end{proof}

\begin{rmk}
		\label{rmk.2.4.classprojcofs}
	Let $f:A\rightarrow B$ be a map
	of $T$-spectra.
	\begin{enumerate}
		\item	$f$ is a cofibration in the projective model structure
					if and only if $f^{0}:A^{0}\rightarrow B^{0}$ and
					the induced maps 
					$$\xymatrix{T\wedge B^{n}\coprod _{T\wedge A^{n}}A^{n+1} 
											\ar[rr]^-{(\sigma ^{n},f^{n+1})}&& B^{n+1}}
					$$
					are all cofibrations in $\pointedmotivic$.
		\item	$f$ is a fibration in the projective model structure
					if and only if $f$ is a level motivic fibration, i.e. for every $n\geq 0$,
					$f^{n}:A^{n}\rightarrow B^{n}$ is a fibration
					in $\pointedmotivic$.
	\end{enumerate}
\end{rmk}

\begin{prop}
		\label{prop.2.4.Ev-n=Quillen-Proj-functor}
	Let $n\geq 0$.  Consider $\pointedmotivic$ and $\Tspectra$
	equipped with the projective model structure (see 
	theorem \ref{thm2.4.projstabmodstr}).
	Then the adjunction
		$$\xymatrix{(F_{n},Ev_{n},\varphi ):\pointedmotivic \ar[r]& \Tspectra}
		$$
	is a Quillen adjunction.
\end{prop}
\begin{proof}
	It is enough to show that $Ev_{n}$ is a right Quillen functor.
	Let $p:X\rightarrow Y$ be a fibration in the projective model structure for
	$\Tspectra$, then $p$ is a level motivic fibration.
	In particular, $Ev_{n}(p)=p^{n}:X^{n}\rightarrow Y^{n}$ is a
	fibration in $\pointedmotivic$.
	
	Now let $q:X\rightarrow Y$ be a trivial fibration in the projective model structure
	for $\Tspectra$.  Then $q$ is a level motivic trivial fibration.
	In particular, $Ev_{n}(q)=q^{n}:X^{n}\rightarrow Y^{n}$ is a trivial
	fibration in $\pointedmotivic$. 
\end{proof}

	We now proceed to define the injective model structure
	for the category of $T$-spectra.
	
	We say that a map of $T$-spectra $i:A\rightarrow B$ is
	a \emph{level cofibration} (respectively \emph{level trivial cofibration})
	if for every $n\geq 0$, $i^{n}:A^{n}\rightarrow B^{n}$
	is a cofibration (respectively trivial cofibration) in $\pointedmotivic$.
	Notice that a map $i:A\rightarrow B$ is a level cofibration if and only
	if it is a monomorphism in the category of $T$-spectra.
	
	Let $A$ be an arbitrary $T$-spectra.  We say that $A$ is
	$\lambda$-bounded if for every $n\geq 0$,
	the presheaf of pointed simplicial sets $A^{n}$ is
	$\lambda$-bounded.
	
\begin{thm}[Jardine]
		\label{thm.2.4.injmodstrspt}
	Let $\kappa$ be a regular cardinal larger than
	$2^{\alpha}$ where $\alpha$ is the cardinality
	of the set $Map(\smoothS)$ of maps in $\smoothS$.
	There exists a cofibrantly generated
	model structure for the category $\Tspectra$ of $T$-spectra
	with the following choices:
	\begin{enumerate}
		\item The weak equivalences are the level weak equivalences.
		
		\item	The set $I$ of generating cofibrations is 
					$$I=\{i:A\rightarrow B\}
					$$
					where $i$ satisfies the following conditions:
					\begin{enumerate}
						\item $i$ is a level cofibration.
						\item The codomain $B$ of $i$ is
									$\kappa$-bounded.
					\end{enumerate}
		\item	The set $J$ of generating trivial cofibrations is
					$$J=\{j:A\rightarrow B\}
					$$
					where $j$ satisfies the following conditions:
					\begin{enumerate}
						\item	$j$ is a level trivial cofibration.
						\item	The codomain $B$ of $j$ is 
									$\kappa$-bounded.
					\end{enumerate}
	\end{enumerate}
	This model structure will be called the \emph{injective
	model structure} for $T$-spectra.  Furthermore, the
	injective model structure is proper and simplicial.
\end{thm}
\begin{proof}
	We refer the reader \cite[lemma 2.1]{MR1787949}.
\end{proof}

\begin{rmk}
		\label{rmk2.4.injcofclass}
	\begin{enumerate}
		\item	\label{rmk2.4.injcofclass.a}Let $f:A\rightarrow B$ be a map of $T$-spectra.
					Then $f$ is a cofibration
					in the injective model structure for $\Tspectra$
					if and only if $f$ is a level cofibration.
		\item	\label{rmk2.4.injcofclass.b}The identity functor on $\Tspectra$
					induces a left Quillen functor from the
					projective model structure to the injective model
					structure.
	\end{enumerate}
\end{rmk}

	Proposition \ref{prop.2.3.cell-prop-pointedmotcat} implies 
	in particular that $\pointedmotivic$ is a closed symmetric monoidal category.
	The category of $T$-spectra $\Tspectra$ has the structure of
	a closed $\pointedmotivic$-module, which
	is obtained by extending the symmetric monoidal structure for $\pointedmotivic$
	levelwise.
	
	The bifunctor giving the adjunction of two variables is defined
	as follows:
	$$\xymatrix@R=.5pt{- \wedge -:\Tspectra \times \pointedmotivic \ar[r]& \Tspectra \\
										(X,A) \ar@{|->}[r]& X\wedge A}
	$$
	with $(X\wedge A)^{n}=X^{n}\wedge A$ and bonding maps given by
	$$\xymatrix{T\wedge (X^{n}\wedge A) \ar[r]^-{\cong}& (T\wedge X^{n})\wedge A \ar[rr]^-{\sigma ^{n}\wedge id _{A}}
							&& X^{n+1}\wedge A}
	$$
	
	The adjoints are given by:
	$$\xymatrix@R=.5pt{\Omega _{-}-:\pointedmotivic ^{op}\times \Tspectra \ar[r]& \Tspectra \\
										(A,X) \ar@{|->}[r]& \Omega _{A}X}
	$$
	
	$$\xymatrix@R=.5pt{\inthomspectrapresheaf (-,-):(\Tspectra )^{op}\times \Tspectra \ar[r]& \pointedmotivic \\
										(X,Y) \ar@{|->}[r]& \inthomspectrapresheaf (X,Y)}
	$$
	where $(\Omega _{A}X)^{n}=\Omega _{A}X^{n}$ and the bonding maps 
	$T\wedge(\Omega _{A}X^{n})\rightarrow \Omega _{A}X^{n+1}$ are adjoint to
	$$\xymatrix{T\wedge (\Omega _{A}X^{n})\wedge A \ \ar[rr]^-{id\wedge ev_{A}}
							&& T\wedge X^{n} \ar[r]^{\sigma ^{n}}& X^{n+1}}$$
	and $\inthomspectrapresheaf (X,Y)$ is the following pointed simplicial presheaf
	on $\smoothS$:
	$$\xymatrix@R=.5pt{\inthomspectrapresheaf (X,Y):(\smoothS \times \Delta)^{op} \ar[r] &
										Sets \\
										(U,n) \ar@{|->}[r] & \Hom _{\Tspectra}(X \wedge (\Delta ^{n}_{U})_{+},Y)}
	$$

\begin{prop}
		\label{prop.2.4.Tspectra-pointpresheaf-modcat}
	\begin{enumerate}
		\item	\label{prop.2.4.Tspectra-pointpresheaf-modcat.a} Let
					$\Tspectra$ denote the category of $T$-spectra equipped
					with the projective model structure.  Then $\Tspectra$ is
					a $\pointedmotivic$-model category (see definition \ref{def.module-modcats}).
		\item	\label{prop.2.4.Tspectra-pointpresheaf-modcat.b}  Let
					$\Tspectra$ denote the category of $T$-spectra equipped
					with the injective model structure.  Then $\Tspectra$ is
					a $\pointedmotivic$-model category.
	\end{enumerate}
\end{prop}
\begin{proof}
	In both cases we need to check that conditions (\ref{def.module-modcats.a})
	and (\ref{def.module-modcats.b}) in definition \ref{def.module-modcats} are satisfied.
	Condition (\ref{def.module-modcats.b}) is automatic since
	the unit $\ast \coprod \ast$ is cofibrant in $\pointedmotivic$.
	
	(\ref{prop.2.4.Tspectra-pointpresheaf-modcat.a}):      
	To check condition (\ref{def.module-modcats.a}) in definition \ref{def.module-modcats}
	we use
	lemma \ref{lem.cond-Quillen-bifunc}(\ref{lem.cond-Quillen-bifunc.c}) which
	implies that it is enough to prove the following claim:
	Given a cofibration $i:A\rightarrow B$ in $\pointedmotivic$
	and a fibration $p:X\rightarrow Y$ in $\Tspectra$ then
	$(i^{*},p_{*}):\Omega _{B}X\rightarrow \Omega _{A}X\times _{\Omega _{A}Y}\Omega _{B}Y$
	is a fibration of $T$-spectra (in the projective model structure), which is
	trivial if either $i$ or $p$ is a weak equivalence.  But fibrations
	in the projective model structure are level motivic fibrations, by
	proposition \ref{prop.2.3.cell-prop-pointedmotcat} we have that
	$\pointedmotivic$ is a symmetric monoidal model category, so in particular
	$(i^{*},p_{*})$ is a level motivic fibration which is trivial if
	either $i$ or $p$ is a weak equivalence.  This proves the claim.
	
	(\ref{prop.2.4.Tspectra-pointpresheaf-modcat.b}):    
	We will prove directly
	that we have a Quillen bifunctor, i.e. given a cofibration $i:A\rightarrow B$
	in $\pointedmotivic$ and a level cofibration $j:C\rightarrow D$ of $T$-spectra,
	we will show that $i\square j:D\wedge A \coprod _{C\wedge A}C\wedge B \rightarrow D\wedge B$
	is a level cofibration (i.e a cofibration in the injective model structure) which is trivial if
	either $i$ or $j$ is a weak equivalence.  But cofibrations in the injective
	model structure are level cofibrations, and since $\pointedmotivic$
	is a symmetric monoidal model category,
	we have that $i\square j$ is a level cofibration which is trivial
	if either $i$ or $j$ is a weak equivalence.  This finishes the proof.
\end{proof}
			
	If we fix $A$ in $\pointedmotivic$, we get an adjunction
	$$\xymatrix{(-\wedge A,\Omega _{A},\varphi _{A}):\Tspectra \ar[r]& \Tspectra}
	$$		
	
\begin{prop}
		\label{prop.2.4.smash-left-Quillen-functor-spectra}
	Let $A$ in $\pointedmotivic$ be an arbitrary presheaf
	of pointed simplicial sets on $\smoothS$.
	\begin{enumerate}
		\item	\label{prop.2.4.smash-left-Quillen-functor-spectra.a}The
					adjunction $(-\wedge A,\Omega _{A},\varphi _{A})$ defined above
					is a Quillen adjunction for the projective model structure
					on $\Tspectra$.
		\item	\label{prop.2.4.smash-left-Quillen-functor-spectra.b}The
					adjunction $(-\wedge A,\Omega _{A},\varphi _{A})$ defined above
					is a Quillen adjunction for the injective model structure
					on $\Tspectra$.
	\end{enumerate}
\end{prop}
\begin{proof}
	Since every object $A$ in $\pointedmotivic$ is cofibrant,
	the result follows immediately from proposition \ref{prop.2.4.Tspectra-pointpresheaf-modcat}.
\end{proof}

\begin{prop}
		\label{prop.2.4.simp-adjunc-spectraloops}
	Let $X,Y$ be two arbitrary $T$-spectra and let
	$A$ in $\pointedmotivic$ be an arbitrary presheaf
	of pointed simplicial sets.  Then we have the
	following enriched adjunctions:
	$$\xymatrix{Map(A, \inthomspectrapresheaf (X,Y))\ar[r]^-{\alpha}_-{\cong}&
							Map(X\wedge A,Y) \ar[r]^-{\beta}_-{\cong}& Map(X,\Omega _{A}Y)}
	$$
	
	$$\xymatrix{\inthomprepointed (A, \inthomspectrapresheaf (X,Y)) \ar[r]^-{\delta}_-{\cong}& 
							\inthomspectrapresheaf(X\wedge A,Y) \ar[r]^-{\epsilon}_-{\cong}& \inthomspectrapresheaf (X,\Omega _{A}Y)}
	$$
	where the maps in the first row are isomorphisms of simplicial sets and the maps in the second row
	are isomorphisms in $\pointedmotivic$.
\end{prop}
\begin{proof}
	We consider first the simplicial adjunctions:
	To any $n$-simplex $t$ in $Map(A, \inthomspectrapresheaf (X,Y))$
		$$\xymatrix{A\otimes \Delta ^{n} \ar[r]^-{t}& \inthomspectrapresheaf (X,Y)}
		$$
	associate the following $n$-simplex in $Map(X\wedge A,Y)$:
	$$\xymatrix{ X \wedge A \otimes \Delta ^{n} \ar[r]^-{\alpha (t)}& Y}
	$$
	corresponding to the adjunction between $X\wedge -$ and $\inthomspectrapresheaf (X,-)$.
	
	To any $n$-simplex $t$ in $Map(X\wedge A,Y)$
	$$\xymatrix{\Delta ^{n}\otimes X\wedge A \ar[r]^-{\cong}& X\wedge A\otimes \Delta ^{n} \ar[r]^-{t}& Y}
	$$
	associate the following $n$-simplex in $Map(X,\Omega _{A}Y)$:
	$$\xymatrix{ X\otimes \Delta ^{n} \ar[r]^-{\cong}& \Delta ^{n}\otimes X \ar[r]^-{\beta (t)}& \Omega _{A}Y}
	$$
	corresponding to the adjunction between $-\wedge A$ and $\Omega _{A}$.
	
	We consider now the isomorphisms of simplicial presheaves:
	To any simplex $s$
	in $\inthomprepointed (A, \inthomspectrapresheaf (X,Y))$
	$$\xymatrix{ A\wedge \Delta ^{n}_{U}  \ar[r]^-{s}& \inthomspectrapresheaf (X,Y)}
	$$
	we associate the following simplex in $\inthomspectrapresheaf (X\wedge A, Y)$
	$$\xymatrix{X\wedge A\wedge \Delta ^{n}_{U} \ar[r]^-{\delta (s)}& Y}
	$$
	corresponding to the adjunction between $X\wedge -$ and $\inthomspectrapresheaf (X,-)$.
	
	To any simplex $s$
	in $\inthomspectrapresheaf (X\wedge A,Y)$
	$$\xymatrix{ X\wedge \Delta ^{n}_{U}\wedge A \ar[r]^{\cong} & X\wedge A\wedge \Delta ^{n}_{U} \ar[r]^-{s}& Y}
	$$
	we associate the following simplex in $\inthomspectrapresheaf (X,\Omega _{A}Y)$
	$$\xymatrix{X\wedge \Delta ^{n}_{U} \ar[r]^-{\epsilon (s)}& \Omega _{A}Y}
	$$
	corresponding to the adjunction between $-\wedge A$ and $\Omega _{A}$.
\end{proof}
			
	We now proceed to define the stable model structure
	for the category of $T$-spectra.
	Consider the functor $\Tloops$ of $T$-loops in
	$\Tspectra$.  There is another way to promote
	the $T$-loops functor from the category of pointed
	simplicial presheaves to the category of $T$-spectra.
	
\begin{defi}
		\label{def.2.4.loopfunctor-on-Tspectra}
		We define the functor $\fakeTloops$ as follows:
				$$\xymatrix@R=.5pt{\fakeTloops :\Tspectra \ar[r]& \Tspectra \\
														X \ar@{|->}[r]& \fakeTloops X}
				$$
				where $(\fakeTloops X)^{n}=\Tloops X^{n}$ and the bonding
				maps $T\wedge \Tloops X^{n}\rightarrow \Tloops X^{n+1}$
				are given by the adjoints to
				$$\xymatrix{\Tloops X^{n} \ar[rr]^-{\Tloops (\sigma ^{n}_{*})} && \Tloops(\Tloops X^{n+1})}
				$$
				where $\sigma ^{n}_{*}:X^{n}\rightarrow \Tloops X^{n+1}$
				is adjoint to the bonding map 
				$$\xymatrix{X^{n}\wedge T \ar[r]^-{\cong} &T\wedge X^{n} \ar[r]^-{\sigma ^{n}} & X^{n+1}}
				$$
\end{defi}

	Following Jardine we call
	the functor $\fakeTloops$ the \emph{fake $T$-loops functor}.

\begin{rmk}
		\label{rmk.2.4.adjoints-for-loops}
		The fake $T$-loops functor $\fakeTloops$ has a left adjoint $\fakeTsuspfunctor$
		called the \emph{fake $T$-suspension functor} defined as follows:
		$$\xymatrix@R=.5pt{\fakeTsuspfunctor : \Tspectra \ar[r] & \Tspectra \\
											X \ar@{|->}[r] & \fakeTsuspfunctor X}
		$$
		where $(\fakeTsuspfunctor X)^{n}=T\wedge X^{n}$ and the bonding
		maps are 
		$$id\wedge \sigma ^{n}:T\wedge (T\wedge X^{n})\rightarrow T\wedge X^{n+1}$$
\end{rmk}

	We will denote by $\Tsuspfunctor$ the left adjoint $(-\wedge T)$ to
	$\Tloops$.

	For any $T$-spectrum $X$, the adjoints
	$\sigma ^{n}_{\ast}:X^{n}\rightarrow \Tloops X^{n+1}$ of the bonding maps
	are the levelwise components of
	a map $\sigma _{*}:X\rightarrow \fakeTloops X[1]$.
	Consider the following inductive system of $T$-spectra:
	$$\xymatrix{X \ar[r] \ar[r]^-{\sigma _{\ast}}& 
							\fakeTloops X[1] \ar[rr]^-{\fakeTloops \sigma _{\ast}[1]} && 
							(\fakeTloops)^{2}X[2] \ar[rr]^-{(\fakeTloops)^{2}\sigma _{*}[2]} && \cdots}
	$$
	and denote its colimit by $Q_{T}X$.  
	The functor $Q_{T}$ is called the \emph{stabilization functor}.
	
	Following Jardine, $J$ will denote a fibrant replacement functor
	for the projective model structure and $I$ will denote
	the corresponding fibrant replacement functor for
	the injective model structure on $\Tspectra$.
	The tranfinite composition
	$X\rightarrow Q_{T}X$ will be denoted $\eta _{X}$, and we 
	define $\tilde{\eta}_{X}$ as the composition
	$$\xymatrix{X \ar[r]^-{\eta_{X}} & Q_{T}X \ar[rr]^-{Q_{T}(j_{X})}&& Q_{T}JX}
	$$
	
	We say that a map $f:X\rightarrow Y$ of $T$-spectra is a
	\emph{stable equivalence} if it becomes a level
	equivalence after taking a fibrant replacement and applying
	the stabilization functor, i.e. if
	$Q_{T}J(f):Q_{T}JX\rightarrow Q_{T}JY$ is a level
	equivalence of $T$-spectra.
	 
\begin{rmk}
		\label{rmk.2.4.class-motstableweakequivs}
	Let $f:X\rightarrow Y$ be a map of $T$-spectra.
	\begin{enumerate}
		\item	$f$
					is a stable equivalence if and only if the map
					$$IQ_{T}J(f):IQ_{T}JX\rightarrow IQ_{T}JY$$ 
					is a level
					equivalence of $T$-spectra.
		\item	If $f$ is a level motivic equivalence then
					$f$ is also a stable equivalence.
	\end{enumerate}	
\end{rmk}	 
		
\begin{thm}[Jardine]
		\label{thm.2.4.stableTspectramodelstr}
	Let $\Tspectra$ be the category of $T$-spectra equipped
	with the projective model structure (see theorem \ref{thm2.4.projstabmodstr}).
	Then the left Bousfield localization of $\Tspectra$
	with respect to the class of stable equivalences exists,
	and furthermore it is proper and simplicial.
	This model structure will be called  \emph{motivic stable}, and the category
	of $T$-spectra $\Tspectra$, equipped with the motivic stable model structure
	will be denoted by $\motivicTspectra$.
\end{thm}
\begin{proof}
	We refer the reader to \cite[theorem 2.9]{MR1787949}.
\end{proof}

\begin{prop}
		\label{prop.2.4.Ev-n=Quillen-stable-functor}
	Let $n\geq 0$.  Consider the adjunction
		$$\xymatrix{(F_{n},Ev_{n},\varphi):\pointedmotivic \ar[r]& \motivicTspectra}
		$$
	described in proposition \ref{prop.2.4.Ev-n=Quillen-Proj-functor}.
	Then $(F_{n},Ev_{n},\varphi)$ is a Quillen adjunction. 
\end{prop}
\begin{proof}
	Follows immediately from proposition \ref{prop.2.4.Ev-n=Quillen-Proj-functor}
	and the following fact:
	\begin{itemize}
		\item	The identity functor on $\Tspectra$ is a left Quillen functor
					from the projective model structure to the motivic stable model structure.
	\end{itemize}
\end{proof}
	 
\begin{lem}[Jardine]
		\label{lem.2.4.char-stable-fibs}
	Let $p:X\rightarrow Y$ be a map of $T$-spectra.
	Then $p$ is a fibration in  $\motivicTspectra$ (we then say that $p$ is a
	\emph{stable fibration}) if the following conditions
	are satisfied:
	\begin{enumerate}
		\item	$p$ is a fibration in the projective model structure
					for $\Tspectra$, i.e. $p$ is a level motivic fibration.
		\item	The following diagram is level homotopy Cartesian:
					$$\xymatrix{X \ar[r]^-{\tilde{\eta}_{X}} \ar[d]_-{p}& Q_{T}JX \ar[d]^-{Q_{T}J(p)} \\
											Y \ar[r]_-{\tilde{\eta}_{Y}}& Q_{T}JY}
					$$
	\end{enumerate}
\end{lem}
\begin{proof}
	We refer the reader to \cite[lemma 2.7]{MR1787949}.
\end{proof}

\begin{lem}[Jardine]
		\label{lem.2.4.fibrant-objects-stablemodstr}
	Let $X$ be a $T$-spectrum.
	The following are equivalent:
	\begin{enumerate}
		\item	\label{lem.2.4.fibrant-objects-stablemodstr.a}$X$ is a fibrant 
					object in $\motivicTspectra$ (we then say
					that $X$ is \emph{stably fibrant}).
		\item	\label{lem.2.4.fibrant-objects-stablemodstr.b}$X$ is a fibrant 
					object in the projective model
					structure for $T$ spectra (i.e. $X$ is level motivic fibrant)
					and the adjoints to the bonding maps
					$\sigma ^{n}_{*}:X^{n}\rightarrow \Tloops X^{n+1}$
					are weak equivalences in $\pointedmotivic$.
		\item \label{lem.2.4.fibrant-objects-stablemodstr.c}$X$ is a fibrant 
					object in the projective model
					structure for $T$-spectra and the adjoints to the
					bonding maps are sectionwise weak equivalences
					of simplicial sets, i.e. for any $U\in (\smoothS)$
					the induced map $\sigma ^{n}_{*}(U):X^{n}(U)\rightarrow \Tloops X^{n+1}(U)$
					is a weak equivalence of simplicial sets.
	\end{enumerate}
\end{lem}
\begin{proof}
	We refer the reader to \cite[lemma 2.8]{MR1787949}.
\end{proof}

	We say that a $T$-spectrum $X$ is \emph{stably fibrant injective},
	if $X$ is a fibrant object in both the motivic stable and the injective
	model structures for $\Tspectra$.

\begin{cor}
		\label{cor.2.4.stablyfibrant-injective-model}
	Let $X$ be a $T$-spectrum.  Then
	$IQ_{T}JX$ is a stably fibrant injective
	replacement for $X$, i.e. the natural map
		$$\xymatrix{X \ar[r]^-{r_{X}}& IQ_{T}JX}
		$$
	is a level weak equivalence (in particular a stable weak equivalence)
	and $IQ_{T}JX$ is stably fibrant injective.
\end{cor}
\begin{proof}
	It is clear that $r_{X}$ is a level weak equivalence and
	that $IQ_{T}JX$ is fibrant in the injective model structure
	for $\Tspectra$, so we only need
	to show that $IQ_{T}JX$ is stably fibrant.
	Since the identity functor on $\Tspectra$
	is a left Quillen functor from the projective to the
	injective model structure (see remark 
	\ref{rmk2.4.injcofclass}(\ref{rmk2.4.injcofclass.b})), we have that
	$IQ_{T}JX$ is in particular 
	a fibrant object in the projective
	model structure for $T$ spectra.
	Lemma \ref{lem.2.4.fibrant-objects-stablemodstr}(\ref{lem.2.4.fibrant-objects-stablemodstr.c})
	implies that it is enough to show that the adjoints
	to the bonding maps for $IQ_{T}JX$,
	$\sigma ^{n}_{\ast}:(IQ_{T}JX)^{n}\rightarrow \Tloops(IQ_{T}JX)^{n+1}$
	are all sectionwise weak equivalences of simplicial sets.
	We will prove that using the following commutative diagram, and showing
	that the top row and the vertical maps are all sectionwise weak
	equivalences of simplicial sets:
		$$\xymatrix{(Q_{T}JX)^{n} \ar[r]^-{\tau ^{n}_{\ast}} \ar[d] & \Tloops (Q_{T}JX)^{n+1} \ar[d] \\
								(IQ_{T}JX)^{n} \ar[r]^-{\sigma ^{n}_{\ast}}& \Tloops (IQ_{T}JX)^{n+1}}
		$$
	
	A cofinal argument implies that the adjoints
	of the bonding maps for $Q_{T}JX$: 
		$$\xymatrix{(Q_{T}JX)^{n} \ar[r]^-{\tau ^{n}_{\ast}}& \Tloops (Q_{T}JX)^{n+1}}
		$$
	are isomorphisms, so in particular these maps are sectionwise weak equivalences
	of simplicial sets.
	
	Since the B.G. property (see definition \ref{def.2.1.Nisnevich-descent})
	is preserved under filtered colimits and the fibrant objects for $\pointedmotivic$ 
	in particular satisfy the B.G. property
	(see theorem \ref{thm.2.1.Nisnevich-descent}),
	we have that the pointed simplicial presheaves $(Q_{T}JX)^{n}$
	satisfy the B.G. property.  Therefore theorem \ref{thm.2.1.Nisnevich-descent}
	implies that the maps
	$(Q_{T}JX)^{n}\rightarrow (IQ_{T}JX)^{n}$ are sectionwise weak
	equivalences of simplicial sets.

	Remark \ref{rmk.2.3.stabilitymotivicflasque}(\ref{rmk.2.3.stabilitymotivicflasque.a})
	implies that
	the pointed simplicial presheaves $(Q_{T}JX)^{n}$
	are all motivic flasque,
	and since the simplicial presheaves $(IQ_{T}JX)^{n}$ are fibrant in
	$\pointedmotivic$, we have that $(IQ_{T}JX)^{n}$ are also
	motivic flasque.
	Now since $T$ is compact in the sense of Jardine (see proposition \ref{prop.2.4.T=compact}),
	we have that  the maps $\Tloops (Q_{T}JX)^{n+1}\rightarrow \Tloops (IQ_{T}JX)^{n+1}$
	are sectionwise weak equivalences of simplicial sets.
	This finishes the proof. 
\end{proof}

\begin{cor}
		\label{cor.2.4.loops-preserve-stablyfib-objects}
	Let $A$ in $\pointedmotivic$ be an arbitrary pointed
	simplicial presheaf, and let $X$ be a stably fibrant
	$T$-spectrum.  Then $\Omega _{A}X$ is also
	stably fibrant.
\end{cor}
\begin{proof}
	Using proposition \ref{prop.2.4.smash-left-Quillen-functor-spectra}
	we have that $\Omega _{A}$ is a right Quillen functor
	for the projective model structure on $\Tspectra$, therefore
	in particular $\Omega _{A}X$ is level fibrant.
	Lemma \ref{lem.2.4.fibrant-objects-stablemodstr}(\ref{lem.2.4.fibrant-objects-stablemodstr.b})
	implies that $\sigma ^{n}_{*}:X^{n}\rightarrow \Tloops X^{n+1}$ are motivic
	weak equivalences between motivic fibrant objects.
	$\pointedmotivic$ is a symmetric monoidal model category, then
	Ken Brown's lemma \ref{lemma1.1.KenBrown2}
	implies that $\Omega_{A}(\sigma ^{n}_{*}):\Omega _{A}X\rightarrow \Omega _{A}\Tloops X^{n+1}$
	is a motivic weak equivalence.
	Let $\theta^{n}:\Omega _{A}X^{n}\rightarrow \Tloops \Omega _{A}X^{n+1}$
	be the adjoint to the bonding map $T\wedge \Omega _{A}X^{n}\rightarrow \Omega _{A}X^{n+1}$
	for the spectrum $\Omega _{A}X$, then we have the following commutative diagram:
	$$\xymatrix{\Omega _{A}X^{n} \ar[r]^-{\theta ^{n}} \ar[dr]_-{\Omega _{A}(\sigma _{*}^{n})}
							& \Tloops \Omega _{A}X^{n+1} \ar[d]_-{\cong}^-{t} \\
							& \Omega _{A}\Tloops X^{n+1}}
	$$
	where $t$ is the isomorphism which flips loop factors.
	Then the two out of three property for weak equivalences
	in $\pointedmotivic$ implies that the maps
	$\theta^{n}:\Omega _{A}X^{n}\rightarrow \Tloops \Omega _{A}X^{n+1}$
	are  motivic weak equivalences.  Finally, lemma 
	\ref{lem.2.4.fibrant-objects-stablemodstr}(\ref{lem.2.4.fibrant-objects-stablemodstr.b})
	implies that $\Omega _{A}X$ is stably fibrant as we wanted.
\end{proof}
		
\begin{lem}[Jardine]
		\label{lem.2.4.class-stableequivs1}
	Let $f:A\rightarrow B$ be a map of $T$-spectra.
	The following are equivalent:
	\begin{enumerate}
		\item	\label{lem.2.4.class-stableequivs1.a}$f$ is a weak equivalence 
					in $\motivicTspectra$.
		\item \label{lem.2.4.class-stableequivs1.b}For every stably 
					fibrant injective object $X$,
					$f$ induces a bijection
					$$\xymatrix{f^{\ast}:[B,X]_{Spt} \ar[r] & [A,X]_{Spt}}
					$$
					in the homotopy category
					associated to $\motivicTspectra$.
		\item	\label{lem.2.4.class-stableequivs1.c}For every stably 
					fibrant injective object $X$,
					$f$ induces a bijection
					$$\xymatrix{f^{\ast}:[B,X] \ar[r]& [A,X]}
					$$
					in the projective homotopy category for $\Tspectra$.
		\item	\label{lem.2.4.class-stableequivs1.d}For every stably fibrant 
					injective object $X$,
					$f$ induces a weak equivalence of simplicial sets
					$$\xymatrix{f^{*}:Map(B,X) \ar[r] & Map(A,X)}
					$$
	\end{enumerate}
\end{lem}
\begin{proof}
	We refer the reader to \cite[lemma 2.11 and corollary 2.12]{MR1787949}.
\end{proof}

\begin{prop}
		\label{prop.2.4.loops-Quillen-functor-stablespectra}
	Let $A$ in $\pointedmotivic$ be an arbitrary presheaf
	of pointed simplicial sets.  Then the adjunction
		$$\xymatrix{(-\wedge A,\Omega _{A},\varphi _{A}):\motivicTspectra \ar[r]& \motivicTspectra}$$
	is a Quillen adjunction.
\end{prop}
\begin{proof}
	Since the cofibrations in the stable and projective model structures
	for $T$-spectra coincide, we have that $-\wedge A$ preserves
	stable cofibrations (since $-\wedge A$ is a left Quillen functor
	for the projective model structure).  So it only remains to show
	that if $j:B\rightarrow C$ is a trivial cofibration in $\motivicTspectra$, then
	$j\wedge id:B\wedge A\rightarrow C\wedge A$ is a weak
	equivalence in $\motivicTspectra$.  Let 
	$X$ be an arbitrary stably fibrant injective $T$-spectrum,  corollary
	\ref{cor.2.4.loops-preserve-stablyfib-objects} implies that $\Omega _{A}X$
	is also stably fibrant,  and since $\Omega _{A}$ is a right Quillen
	functor for the injective model structure on $\Tspectra$
	(see proposition \ref{prop.2.4.smash-left-Quillen-functor-spectra}), we have that
	$\Omega _{A}X$ is also fibrant in the injective model structure.
	Thus $\Omega _{A}X$ is stably fibrant injective, then
	lemma \ref{lem.2.4.class-stableequivs1}
	implies that $j^{*}:Map(C,\Omega _{A}X)\rightarrow Map(B,\Omega _{A}X)$
	is a weak equivalence of simplicial sets.  Using the enriched adjunction
	of propositon \ref{prop.2.4.simp-adjunc-spectraloops}, $j^{*}$ becomes
	$(j\wedge id)^{*}:Map(C\wedge A,X)\rightarrow Map(B\wedge A,X)$.
	Finally since $(j\wedge id)^{*}$ is a weak equivalence for every stably fibrant
	injective spectrum $X$, we get that $j\wedge id:C\wedge A\rightarrow B\wedge A$
	is a weak equivalence in $\motivicTspectra$.
\end{proof}

\begin{prop}
		\label{prop.2.4.motTspectra-pointedpresheaf-modelcat}
	$\motivicTspectra$ is a $\pointedmotivic$-model category
	(see definition \ref{def.module-modcats}).
\end{prop}
\begin{proof}
	Condition (\ref{def.module-modcats.b}) in definition \ref{def.module-modcats}
	follows automatically since the unit in $\pointedmotivic$ is cofibrant.
	It only remains to prove that if $i:A\rightarrow B$ is a cofibration in
	$\pointedmotivic$ and $j:C\rightarrow D$ is a cofibration in $\motivicTspectra$
	then $i\square j:D\wedge A\coprod _{C\wedge A}C\wedge B\rightarrow D\wedge B$
	is a cofibration in $\motivicTspectra$ which is trivial if either $i$ or $j$
	is a weak equivalence.  Since the cofibrations in the projective and
	the motivic stable model structure for $\Tspectra$ coincide, and
	proposition \ref{prop.2.4.Tspectra-pointpresheaf-modcat} 
	implies in particular that the category
	of $T$-spectra equipped with the projective model structure is a
	$\pointedmotivic$-model category, we have that
	$i\square j$ is a cofibration in the motivic stable structure.
	It only remains to show that $i\square j$ is a stable weak equivalence when
	either $i$ or $j$ is a weak equivalence.  If $i$ is a weak equivalence
	(i.e. a trivial cofibration) then using proposition \ref{prop.2.4.Tspectra-pointpresheaf-modcat}
	again we have that $i\square j$ is a level weak equivalence,
	therefore $i\square j$ is also a stable equivalence 
	(see remark \ref{rmk.2.4.class-motstableweakequivs}).
	Finally if $j$ is a stable equivalence (i.e. a trivial
	cofibration in the motivic stable structure) then we consider
	the following commutative diagram
	$$\xymatrix{C\wedge A \ar[rr]^-{id_{C}\wedge i} \ar[d]_-{j\wedge id_{A}} && 
							C\wedge B \ar[d]^-{f} \ar@/^2pc/[ddr]^-{j\wedge id_{B}}& \\
							D\wedge A \ar[rr] \ar@/_2pc/[drrr] && D\wedge A\coprod _{C\wedge A}C\wedge B \ar[dr]_-{i\square j}&\\
							&&& D\wedge B}
	$$
	Proposition \ref{prop.2.4.loops-Quillen-functor-stablespectra}
	implies that $j\wedge id_{A}$ and $j\wedge id_{B}$ are
	both trivial cofibrations in $\motivicTspectra$.  Thus $f$ is also a trivial cofibration
	(since it is the pushout of $j\wedge id_{A}$ along $id_{C}\wedge i$),
	and therefore the two out of three property for
	stable weak equivalences implies that
	$i\square j$ is a stable equivalence.  This finishes the proof.
\end{proof}

	In order to prove that the motivic stable model structure on 
	$\Tspectra$ is in fact ``stable'', i.e. that the
	$T$-suspension functor $\Tsuspfunctor$ is indeed a Quillen equivalence,
	Jardine introduces bigraded stable homotopy groups
	which allow to give another criterion to detect motivic
	stable weak equivalences.
	
\begin{defi}
		\label{def.2.4.bigraded-stable-homotopygroups}
	Let $X$ be an arbitrary $T$-spectrum.
	The \emph{weighted stable homotopy groups}
	of $X$ are presheaves of abelian groups
	$\pi_{t,s}X$ (where $t,s\in \mathbb Z$) on $\smoothS$.
	For $U\in (\smoothS)$ the sections $\pi_{t,s}X(U)$
	are defined as the colimit of the inductive system:
		$$\xymatrix{[S^{t+n}\wedge \gm ^{s+n},X^{n}|_{U}] \ar[r] & [S^{t+n+1}\wedge \gm ^{s+n+1},X^{n+1}|_{U}] \ar[r] 
							&  \cdots}
		$$
	where $[-,X^{i}|_{U}]$ denotes the set of maps in the homotopy category associated
	to the motivic model structure on the category $\Delta ^{op}Pre_{*}(Sm|_{U})_{Nis}$
	of pointed simplicial presheaves on
	the smooth Nisnevich site over the base scheme $U$,  
	and the transition maps are given by taking suspension with $T$
	and composing with the bonding maps of $X$.  The index
	$t$ is called the \emph{degree} and the index $s$ is called the
	\emph{weight} of $\pi_{t,s}X$.
\end{defi}

\begin{prop}
		\label{prop.2.4.representability-sections-bigraded-homotopy-groups}
	Consider $t,s\in \mathbb Z$ and $U\in (\smoothS )$.  Then
	the following functor:
		$$\xymatrix@R=.5pt{ \motivicTspectra \ar[r]& \text{Abelian Groups}\\
									X \ar@{|->}[r]& \pi _{t,s}X(U)}
		$$
	is representable in the homotopy category associated to $\motivicTspectra$.  
	To represent it we can choose any
	spectrum of the form (see definition \ref{def.2.4.infinitesuspension+shift})
		$$F_{n}(S^{p}\wedge \gm ^{q}\wedge U_{+})$$
	where $n,p,q\geq 0$, $p-n=t$ and $q-n=s$.
\end{prop}
\begin{proof}
	Since every pointed simplicial presheaf on $\smoothS$ is cofibrant in $\pointedmotivic$,
	proposition \ref{prop.2.4.Ev-n=Quillen-stable-functor} 
	and corollary \ref{cor.2.4.stablyfibrant-injective-model} imply that
		$$[F_{n}(S^{p}\wedge \gm ^{q}\wedge U_{+}),X]_{Spt} \cong 
		  [S^{p}\wedge \gm ^{q}\wedge U_{+}, (IQ_{T}JX)^{n}]$$
	where $[-,-]_{Spt}$ denotes the set of maps between two objects
	in the homotopy category associated to $\motivicTspectra$,
	and $[-,-]$ denotes the set of maps in the homotopy category
	associated to $\pointedmotivic$.
	Since $Q_{T}JX\rightarrow IQ_{T}JX$
	is in particular a level motivic trivial fibration we have
	the natural isomorphism
		$$[S^{p}\wedge \gm ^{q}\wedge U_{+}, (IQ_{T}JX)^{n}]\cong 
			[S^{p}\wedge \gm ^{q}\wedge U_{+}, (Q_{T}JX)^{n}]$$
	
	Now since $S^{p}\wedge \gm ^{q}\wedge U_{+}$ is compact in the sense of Jardine
	(see proposition \ref{prop.2.4.T=compact}), 
	using proposition \ref{prop.2.3.compact=>commutes-colimits.b}
	we have that
		$$[S^{p}\wedge \gm ^{q}\wedge U_{+}, (Q_{T}JX)^{n}]\cong
			\varinjlim _{j\geq 0}[S^{p+j}\wedge \gm ^{q+j}\wedge U_{+}, (JX)^{n+j}]$$
	Since $\pointedmotivic$  is
	in particular a symmetric monoidal model category
	(see proposition \ref{prop.2.3.cell-prop-pointedmotcat}) and $U_{+}\in \simpprepointed$
	is cofibrant, we have that
		$$\varinjlim _{j\geq 0}[S^{p+j}\wedge \gm ^{q+j}\wedge U_{+}, (JX)^{n+j}]
			\cong \varinjlim _{j\geq 0}[S^{p+j}\wedge \gm ^{q+j}, \inthomprepointed (U_{+},(JX)^{n+j})]$$
	Proposition \ref{prop.2.1.inthom==basechange} implies that
		$$\varinjlim _{j\geq 0}[S^{p+j}\wedge \gm ^{q+j}, \inthomprepointed (U_{+},(JX)^{n+j})]\cong
			\varinjlim _{j\geq 0}[S^{p+j}\wedge \gm ^{q+j}, \phi _{\ast}\phi ^{-1}(JX)^{n+j}]$$
	where $\phi :U\rightarrow S$ is the structure map defining $U$ as an object in $\smoothS$.
	
	Now since $\pointedmotivic$  is in particular
	a simplicial model category, and $\phi _{\ast}\phi ^{-1}(JX)^{n+j}\cong \inthomprepointed (U_{+},(JX)^{n+j})$
	is a fibrant object,
	we have that
		$$\pi _{0}(Map(S^{p+j}\wedge \gm ^{q+j}, \phi _{\ast}\phi ^{-1}(JX)^{n+j}))$$
	computes
	$[S^{p+j}\wedge \gm ^{q+j}, \phi _{\ast}\phi ^{-1}(JX)^{n+j}]$.
	The enriched adjunctions of proposition
	\ref{prop.2.1.basechange-enriched-adj} imply that
		\begin{eqnarray*}
			Map(S^{p+j}\wedge \gm ^{q+j}, \phi _{\ast}\phi ^{-1}(JX)^{n+j}) & \cong &
			Map(\phi ^{-1}(S^{p+j}\wedge \gm ^{q+j}), \phi ^{-1}(JX)^{n+j}) \\
			& = & Map(S^{p+j}\wedge \gm ^{q+j}, \phi ^{-1}(JX)^{n+j})
		\end{eqnarray*}
	Let $r_{U}:\phi ^{-1}(JX)^{n+j}\rightarrow R_{U}\phi ^{-1}(JX)^{n+j}$
	be a functorial fibrant replacement for $\phi ^{-1}(JX)^{n+j}$ in
	the category of pointed simplicial presheaves $\Delta ^{op}Pre_{*}(Sm|_{U})_{Nis}$
	on the smooth Nisnevich site over $U$ equipped with the motivic model structure.
	It is clear that $(JX)^{n+j}$ is motivic flasque (see definition \ref{def.2.3.motivic-flasque})
	and satisfies the B.G. property (see definition \ref{def.2.1.Nisnevich-descent})
	on $\simpprepointed$, and since  $\phi ^{-1}$ preserves both properties
	we have that $\phi ^{-1}(JX)^{n+j}$ is motivic flasque and satisfies the B.G.
	property on $\Delta ^{op}Pre_{*}(Sm|_{U})_{Nis}$.  Thus
	$r_{U}$ is a sectionwise weak equivalence,
	and since
	$S^{p+j}\wedge \gm ^{q+j}$ is compact in the sense of Jardine in
	$\Delta ^{op}Pre_{*}(Sm|_{U})_{Nis}$ (see proposition \ref{prop.2.4.T=compact})
	we have that
		\begin{equation}
					\label{global-induced-compatibility}
			\begin{array}{c}
				\xymatrix{\inthomprepointed (S^{p+j}\wedge \gm ^{q+j}, \phi ^{-1}(JX)^{n+j} \ar[d]^-{r_{U}\ast}) \\
									\inthomprepointed (S^{p+j}\wedge \gm ^{q+j}, R_{U}\phi ^{-1}(JX)^{n+j})}
			\end{array}
		\end{equation}
	is also a sectionwise weak equivalence.  Taking global sections at
	$U$ we get a weak equivalence of simplicial sets:
		\begin{equation}
					\label{local-induced-compatibility}
			\begin{array}{c}
				\xymatrix{Map(S^{p+j}\wedge \gm ^{q+j}, \phi ^{-1}(JX)^{n+j}) \ar[d]^-{r_{U}\ast} \\
									Map(S^{p+j}\wedge \gm ^{q+j}, R_{U}\phi ^{-1}(JX)^{n+j})}
			\end{array}
		\end{equation}
	Thus $Map(S^{p+j}\wedge \gm ^{q+j}, \phi _{\ast}\phi ^{-1}(JX)^{n+j})$ and
	$Map(S^{p+j}\wedge \gm ^{q+j}, R_{U}\phi ^{-1}(JX)^{n+j})$ are naturally weakly
	equivalent simplicial sets.  Since $\Delta ^{op}Pre_{*}(Sm|_{U})_{Nis}$ is
	a simplicial model category we have that
		$$\pi _{0}Map(S^{p+j}\wedge \gm ^{q+j}, R_{U}\phi ^{-1}(JX)^{n+j})$$
	computes $[S^{p+j}\wedge \gm ^{q+j}, \phi ^{-1}(JX)^{n+j}]_{U}=[S^{p+j}\wedge \gm ^{q+j}, (JX)^{n+j}|_{U}]$,
	where $[-,-]_{U}$ denotes the set of maps in the homotopy category associated to the
	motivic model structure on $\Delta ^{op}Pre_{*}(Sm|_{U})_{Nis}$.
	Thus $[S^{p+j}\wedge \gm ^{q+j}, \phi _{\ast}\phi ^{-1}(JX)^{n+j}]$ is naturally
	isomorphic to $[S^{p+j}\wedge \gm ^{q+j}, (JX)^{n+j}|_{U}]$.
	This implies that
		\begin{eqnarray}
			[F_{n}(S^{p}\wedge \gm ^{q}\wedge U_{+}),X]_{Spt}  & \cong &
			\varinjlim _{j\geq 0}[S^{p+j}\wedge \gm ^{q+j}, \phi _{\ast}\phi ^{-1}(JX)^{n+j}] \nonumber\\
			& \cong & \varinjlim _{j\geq 0}[S^{p+j}\wedge \gm ^{q+j}, (JX)^{n+j}|_{U}] \label{iso-compatible}\\
			& \cong & \varinjlim _{j\geq 0}[S^{p+j}\wedge \gm ^{q+j}, X^{n+j}|_{U}] \nonumber \\
			& \cong & \pi _{p-n,q-n}X(U)=\pi _{t,s}X(U) \nonumber
		\end{eqnarray}
	Therefore the functors $[F_{n}(S^{p}\wedge \gm ^{q}\wedge U_{+}),-]_{Spt}$ and
	$\pi _{t,s}(-)(U)$ have canonically isomorphic image for every $T$-spectrum $X$.
	To finish the proof we will give an element 
	$\alpha \in \pi _{t,s}(F_{n}(S^{p}\wedge \gm ^{q}\wedge U_{+}))(U)$
	which induces	an isomorphism of functors
		$$\xymatrix{[F_{n}(S^{p}\wedge \gm ^{q}\wedge U_{+}),-]_{Spt} \ar[r]^-{\alpha _{\ast}}& \pi _{t,s}(-)(U)}
		$$
	Consider the identity map
	$id:S^{t+j}\wedge \gm ^{s+j}\wedge U_{+}\rightarrow S^{t+j}\wedge \gm ^{s+j}\wedge U_{+}$.
	Since $-\wedge U_{+}$ and $\inthomprepointed (U_{+},-)$ are adjoint functors, we have
	an associated adjoint $\beta ^{j}$:
		$$\xymatrix{S^{t+j}\wedge \gm ^{s+j} \ar[r]^-{\beta ^{j}}& 
								\inthomprepointed (U_{+},S^{t+j}\wedge \gm ^{s+j}\wedge U_{+})\cong 
								\phi _{\ast}\phi ^{-1}(S^{t+j}\wedge \gm ^{s+j}\wedge U_{+})}
		$$
	Now let $\gamma ^{j}$ be the adjoint to $\beta ^{j}$ corresponding to the adjunction
	between $\phi ^{-1}$ and $\phi _{\ast}$:
		$$\xymatrix{\phi ^{-1}(S^{t+j}\wedge \gm ^{s+j}) \ar[r]^-{\gamma ^{j}}& 
			\phi ^{-1}(S^{t+j}\wedge \gm ^{s+j}\wedge U_{+})}
		$$
	Let $[\gamma ^{j}]\in [\phi ^{-1}(S^{t+j}\wedge \gm ^{s+j}),\phi ^{-1}(S^{t+j}\wedge \gm ^{s+j}\wedge U_{+})]_{U}
	=[S^{t+j}\wedge \gm ^{s+j}, (S^{t+j}\wedge \gm ^{s+j}\wedge U_{+})|_{U}]$
	denote the map induced by $\gamma ^{j}$ in the homotopy category associated to
	$\Delta ^{op}Pre_{*}(Sm|_{U})_{Nis}$ equipped with the motivic model structure.
	It is clear that the maps $[\gamma ^{j}]$ define an element
		$$\alpha \in \varinjlim _{j\geq 0}[S^{t+j}\wedge \gm ^{s+j}, (S^{t+j}\wedge \gm ^{s+j}\wedge U_{+})|_{U}]
		$$
	But 
		\begin{eqnarray*}
			[S^{t+j}\wedge \gm ^{s+j}, (S^{t+j}\wedge \gm ^{s+j}\wedge U_{+})|_{U}] & =&
			[S^{t+j}\wedge \gm ^{s+j}, (S^{p+j-n}\wedge \gm ^{q+j-n}\wedge U_{+})|_{U}] \\
			& = & [S^{t+j}\wedge \gm ^{s+j}, (F_{n}(S^{p}\wedge \gm ^{q}\wedge U_{+}))^{j}|_{U}]
		\end{eqnarray*}
	Thus
		$$\alpha \in \varinjlim _{j\geq 0}[S^{t+j}\wedge \gm ^{s+j}, (F_{n}(S^{p}\wedge \gm ^{q}\wedge U_{+}))^{j}|_{U}]
			=\pi _{t,s}(F_{n}(S^{p}\wedge \gm ^{q}\wedge U_{+}))$$
			
	Finally it is clear that $\alpha$ induces the required isomorphism of functors
		$$\xymatrix{[F_{n}(S^{p}\wedge \gm ^{q}\wedge U_{+}),-]_{Spt} \ar[r]^-{\alpha _{\ast}}& \pi _{t,s}(-)(U)}
		$$
	since by construction $\alpha$ is compatible with the isomorphisms in
	(\ref{iso-compatible}) which are induced by $r_{U}:\phi ^{-1}(JX)^{n+j}\rightarrow R_{U}\phi ^{-1}(JX)^{n+j}$
	via the natural maps $r_{U}\ast$ in the diagrams (\ref{global-induced-compatibility})
	and (\ref{local-induced-compatibility}), where $r_{U}$ denotes
	a functorial fibrant replacement in the category
	$\Delta ^{op}Pre_{*}(Sm|_{U})_{Nis}$ equipped with the motivic model structure. 
\end{proof}

\begin{prop}[Jardine]
		\label{prop2.4.weightedgroups-detectwequis}
	Let $f:X\rightarrow Y$ be a map of $T$-spectra.
	The following are equivalent:
	\begin{enumerate}
		\item $f$ is a weak equivalence in $\motivicTspectra$.
		\item	For every $t,s\in \mathbb Z$, $f$ induces an isomorphism
					$$\xymatrix{\pi_{t,s}(f):\pi_{t,s}X \ar[r]& \pi_{t,s}Y}
					$$
					of presheaves of abelian groups on $\smoothS$.
	\end{enumerate}
\end{prop}
\begin{proof}
	We refer the reader to \cite[lemma 3.7]{MR1787949}.
\end{proof}

\begin{cor}
		\label{cor.2.4.detecting-wequivs-repsobjects}
	Let $f:X\rightarrow Y$ be a map of $T$-spectra.
	The following are equivalent:
	\begin{enumerate}
		\item	$f$ is a weak equivalence in $\motivicTspectra$.
		\item	For every $n,p,q\geq 0$ and every $U\in \smoothS$, $f$ induces an
					isomorphism
						$$\xymatrix{[F_{n}(S^{p}\wedge \gm ^{q}\wedge U _{+}),X]_{Spt} \ar[r]^-{f_{\ast}}&
												[F_{n}(S^{p}\wedge \gm ^{q}\wedge U_{+}),Y]_{Spt}}
						$$
					in the homotopy category associated to $\motivicTspectra$.
	\end{enumerate}
\end{cor}
\begin{proof}
	Follows immediately from propositions \ref{prop2.4.weightedgroups-detectwequis} and 
	\ref{prop.2.4.representability-sections-bigraded-homotopy-groups}.
\end{proof}

\begin{thm}[Jardine]
		\label{thm.2.4.Tloops-Quillenequiv}
	The Quillen adjunction:
	$$\xymatrix{(\Tsuspfunctor,\Tloops,\varphi):\motivicTspectra \ar[r]& \motivicTspectra}
	$$
	is a Quillen equivalence.
\end{thm}
\begin{proof}
	We refer the reader to \cite[theorem 3.11 and corollary 3.17]{MR1787949}.
\end{proof}

\begin{prop}[Jardine]
		\label{prop.2.4.fakesusp-shift-equiv}
	The natural map $\fakeTsuspfunctor X\rightarrow X[1]$
	from the fake suspension functor to the shift functor
	is a weak equivalence in $\motivicTspectra$.  
	Therefore the fake suspension functor
	and the shift functor are naturally equivalent
	in the homotopy category associated to $\motivicTspectra$.
\end{prop}
\begin{proof}
	We refer the reader to \cite[lemma 3.19]{MR1787949}.
\end{proof}

\begin{prop}[Jardine]
		\label{prop.2.4.fakesusp-susp-equiv}
		The fake suspension functor $\fakeTsuspfunctor$ and
		the suspension functor $\Tsuspfunctor$ are naturally
		equivalent in the homotopy category associated
		to $\motivicTspectra$.
\end{prop}
\begin{proof}
	We refer the reader to \cite[lemma 3.20]{MR1787949}.
\end{proof}

\begin{cor}[Jardine]
		\label{cor.2.4.loopfunctors-equiv}
	The $T$-loops functor $\Tloops$, fake $T$-loops functor $\fakeTloops$,
	and shift functor $s_{-1}$ ($s_{-1}X=X[-1]$) are all
	naturally equivalent in the homotopy category
	associated to $\motivicTspectra$.
\end{cor}
\begin{proof}
	Follows immediately from propositions \ref{prop.2.4.fakesusp-shift-equiv} 
	and \ref{prop.2.4.fakesusp-susp-equiv}.
\end{proof}

\begin{prop}
		\label{prop.2.4.compact-respects-colimits-spectra}
	Let $X\in \simpprepointed$ be a pointed simplicial presheaf which
	is compact in the sense of Jardine (see definition \ref{def.2.3.compactness}), and let
	$F_{n}(X)$ be the $T$-spectrum constructed in definition \ref{def.2.4.infinitesuspension+shift}.
	Consider an inductive system of $T$-spectra:
		$$\xymatrix{Z_{0} \ar[r]& Z_{1} \ar[r]& Z_{2} \ar[r]& \cdots &}
		$$
	Then 
		$$[F_{n}(X),\varinjlim Z_{i}]_{Spt}\cong \varinjlim[F_{n}(X),Z_{i}]_{Spt}
		$$
	where $[-,-]_{Spt}$ denotes the set of maps in the homotopy category
	associated to $\motivicTspectra$.
\end{prop}
\begin{proof}
	Since $X$ is cofibrant in $\pointedmotivic$, 
	proposition \ref{prop.2.4.Ev-n=Quillen-stable-functor} and corollary
	\ref{cor.2.4.stablyfibrant-injective-model} imply that
		\begin{eqnarray*}
			[F_{n}(X),\varinjlim Z_{i}]_{Spt} &\cong& 
			[X, (IQ_{T}J\varinjlim Z_{i})^{n}] \\
			&\cong & [X, (Q_{T}J\varinjlim Z_{i})^{n}]
		\end{eqnarray*}
	where $[-,-]$ denotes the set of maps in the homotopy category
	associated to $\pointedmotivic$.
	Since $X$  is compact in the sense
	of Jardine,
	we have that proposition \ref{prop.2.3.compact=>commutes-colimits.b} implies the following:
		\begin{eqnarray*}
			[X, (Q_{T}J\varinjlim Z_{i})^{n}] &\cong &
			\varinjlim _{j\geq 0} [S^{j}\wedge \gm ^{j}\wedge X, (J\varinjlim Z_{i})^{n+j}]\\
			&\cong & \varinjlim _{j\geq 0} [S^{j}\wedge \gm ^{j}\wedge X, (\varinjlim Z_{i})^{n+j}]
		\end{eqnarray*}
	Now lemma 2.2(4) in \cite{MR1787949} implies that
	$S^{j}\wedge \gm ^{j}\wedge X$ are all compact in the sense of Jardine, 
	therefore using proposition 
	\ref{prop.2.3.compact=>commutes-colimits.b} again, we have:
		\begin{eqnarray*}
			\varinjlim _{j\geq 0} [S^{j}\wedge \gm ^{j}\wedge X, (\varinjlim Z_{i})^{n+j}] &\cong &
			\varinjlim _{j\geq 0} \varinjlim _{i\geq 0}[S^{j}\wedge \gm ^{j}\wedge X, (Z_{i})^{n+j}]\\
			&\cong & \varinjlim _{i\geq 0} \varinjlim _{j\geq 0}[S^{j}\wedge \gm ^{j}\wedge X, (Z_{i})^{n+j}]\\
			&\cong & \varinjlim _{i\geq 0} [X, (Q_{T}JZ_{i})^{n}]\\
			&\cong & \varinjlim _{i\geq 0} [F_{n}(X), Z_{i}]_{Spt}
		\end{eqnarray*}
	and this finishes the proof.
\end{proof}

\end{section}
\begin{section}{Cellularity of the Motivic Stable Model Structure}
		\label{section-cellularity-motivic-stable}
		
	In this section we will show
	that $\motivicTspectra$ is a
	cellular model category.  For this we will use 
	the cellularity of 
	$\pointedmotivic$ (see proposition \ref{prop.2.3.cell-prop-pointedmotcat})
	together with some results of Hovey \cite{MR1860878}.
	
	The cellularity for the motivic stable model structure
	is also proved in \cite[corollary 1.6]{MR2197578}.
	However, our proof is different since we use the characterization
	for weak equivalences given in corollary
	\ref{cor.1.1.4.detect-weak-equiv.2}(\ref{cor.1.1.4.detect-weak-equiv.2.b})
	(which holds in any simplicial model category) whereas the argument
	given in \cite[corollary 1.6]{MR2197578} relies  on 
	\cite[theorem 4.12]{MR1860878} which does not
	apply to the model category $\pointedmotivic$ 
	described in proposition \ref{prop.2.3.cell-prop-pointedmotcat}
	(see \cite[p. 83]{MR1860878}).
	
\begin{thm}[Hovey]
		\label{thm.2.5.proj-str-cellular}
	Let $\Tspectra$ be the category of $T$-spectra
	equipped with the projective model structure
	(see theorem \ref{thm2.4.projstabmodstr}).
	Then the category $\Tspectra$ is a cellular
	model category where the sets of generating
	cofibrations and trivial cofibrations
	are the ones described in theorem \ref{thm2.4.projstabmodstr}.
\end{thm}
\begin{proof}
	Proposition \ref{prop.2.3.cell-prop-pointedmotcat} implies that
	the model category $\pointedmotivic$ is in particular cellular
	and left proper.  Therefore we can apply
	theorem A.9 in
	\cite{MR1860878}, which says that the category of $T$-spectra
	equipped with the projective model structure
	is also cellular under our conditions.
\end{proof}

	Theorem \ref{thm2.4.projstabmodstr} together 
	with theorem \ref{thm.2.5.proj-str-cellular}
	imply that the projective model structure on
	$\Tspectra$ is cellular, proper and simplicial.
	Therefore we can apply Hirschhorn's localization technology to it.
	If we are able to find a suitable \emph{set} of maps
	such that the left Bousfield localization with respect to this
	set recovers the motivic stable model structure, then
	an immediate corollary of this will be the cellularity of
	the motivic stable model structure for $\Tspectra$.
	
\begin{defi}[Hovey, cf. \cite{MR1860878}]
		\label{def.2.5.localizingset-motstab}
	Let $I_{M_{*}}=\{Y_{+}\hookrightarrow (\Delta ^{n}_{U})_{+}\}$
	be the set of generating cofibrations for $\pointedmotivic$
	(see propositon \ref{prop.2.3.cell-prop-pointedmotcat}).
	Notice that $Y_{+}$ may be equal to $(\Delta ^{n}_{U})_{+}$.
	We consider the following set of maps
	of $T$-spectra
	$$\xymatrix{S=\{ F_{k+1}(T\wedge Y_{+}) \ar[r]^-{\zeta ^{Y}_{k}} & F_{k}Y_{+} \}}
	$$
	where $\zeta ^{Y}_{k}$ is the adjoint to the identity map
	(in $\simpprepointed$)
	$$id: T\wedge Y_{+}\rightarrow Ev_{k+1}(F_{k}Y_{+})=T\wedge Y_{+}$$
	coming from the adjunction between $F_{k+1}$ and $Ev_{k+1}$
	(see definition \ref{def.2.4.infinitesuspension+shift}).
\end{defi}

\begin{prop}[Hovey]
		\label{prop.2.5.S-local==stably-fibrant}
	Let $X$ be a $T$-spectrum.  The
	following conditions are equivalent:
	\begin{enumerate}
		\item \label{prop.2.5.S-local==stably-fibrant.a}  $X$ is stably fibrant, i.e.
					$X$ is a fibrant object in $\motivicTspectra$.
		\item	\label{prop.2.5.S-local==stably-fibrant.b}	$X$ is $S$-local.
	\end{enumerate}
\end{prop}
\begin{proof}
	Follows from \cite[theorem 3.4]{MR1860878} and lemma
	\ref{lem.2.4.fibrant-objects-stablemodstr}.
\end{proof}

	Now it is very easy to show that the
	motivic stable model structure for $T$-spectra
	is in fact cellular.
	
\begin{thm}
		\label{thm.2.5.cellularity-motivic-stable-str}
	$\motivicTspectra$ is a cellular model category with the
	following sets $I^{T}_{M_{*}}$, $J^{T}_{M_{*}}$ of
	generating cofibrations and trivial cofibrations respectively:
	$$\begin{array}{l}
	I^{T}_{M_{*}}  =\bigcup _{k\geq 0}\{F_{k}(Y_{+}\hookrightarrow (\Delta ^{n}_{U})_{+})\mid
									U\in (\smoothS ),n\geq 0 \} \\
	 \\
	J^{T}_{M_{*}}  =\{j:A\rightarrow B\}
	\end{array}
	$$
	where $j$ satisfies the following conditions:
	\begin{enumerate}
		\item	$j$ is an inclusion of $I^{T}_{M_{*}}$-complexes.
		\item	$j$ is a stable weak equivalence.
		\item	the size of $B$ as an $I^{T}_{M_{*}}$-complex is less than $\kappa$, 
					where $\kappa$ is the regular cardinal described by Hirschhorn in 
					\cite[definition 4.5.3]{MR1944041}.
	\end{enumerate}
\end{thm}
\begin{proof}
	By theorem \ref{thm.2.5.proj-str-cellular}
	we know that $\Tspectra$ is cellular when it is equipped with the projective
	model structure.  Therefore we can apply Hirschhorn's localization techniques
	to construct the left Bousfield localization with respect to the set $S$
	of definition \ref{def.2.5.localizingset-motstab}.  
	We claim that this localization coincides with $\motivicTspectra$.
	In effect, using
	proposition \ref{prop.2.5.S-local==stably-fibrant},
	we have that the fibrant objects in the left Bousfield localization
	with respect to $S$ coincide with the fibrant objects in $\motivicTspectra$.
	Therefore a map $f:X\rightarrow Y$ of $T$-spectra is a weak equivalence
	in the left Bousfield localization with respect to $S$ if and only if
	$Qf^{*}:Map(QY,Z)\rightarrow Map(QX,Z)$ is a weak equivalence of simplicial
	sets for every stably fibrant object $Z$ (here $Q$ denotes the cofibrant
	replacement functor in $\Tspectra$ equipped with the projective model
	structure).  But since $\motivicTspectra$ is a simplicial model category
	and the cofibrations coincide with the projective cofibrations,
	using corollary
	\ref{cor.1.1.4.detect-weak-equiv.2}(\ref{cor.1.1.4.detect-weak-equiv.2.b})
	we get exactly the same characterization for the stable equivalences.
	Hence the weak equivalences in both the motivic stable structure and the left
	Bousfield localization with respect to $S$ coincide.  This implies
	that the motivic stable model structure and the left Bousfield localization with
	respect to $S$ are identical, since the cofibrations in both cases
	are just the cofibrations for the projective model structure on $\Tspectra$.
	
	Therefore using theorem 4.1.1 in \cite{MR1944041}
	we have that $\motivicTspectra$
	is cellular, since it is constructed applying Hirschhorn technology
	with respect to the set $S$.
	
	The claim with respect to the sets of generating cofibrations and trivial
	cofibrations also follows from \cite[theorem 4.1.1]{MR1944041}
	and the fact that $I^{T}_{M_{*}}$ is just the set
	of generating cofibrations for the projective model structure
	on $\Tspectra$.
\end{proof}

	Theorem \ref{thm.2.5.cellularity-motivic-stable-str} will be
	one of the main technical ingredients for the construction
	of new model structures on $\Tspectra$ which lift
	Voevodsky's slice filtration to the model category level.

\end{section}
\begin{section}{The Motivic Symmetric Stable Model Structure}
		\label{sec.2.6.symmetricTspectra}

	One of the technical disadvantages of the category of $T$-spectra
	$\Tspectra$ (see definition \ref{defi2.4.T-spectra}) is that
	it does not inherit a closed symmetric monoidal structure from the category
	of pointed simplicial presheaves $\pointedmotivic$.
	Symmetric spectra were introduced by Hovey, Shipley and Smith in
	\cite{MR1695653} to solve this problem in the context of simplicial sets.
	
	Their construction was lifted to the motivic setting by Jardine
	in \cite{MR1787949}, where he constructs a closed symmetric monoidal category of
	$T$-spectra together with a suitable model structure which is Quillen
	equivalent to the category $\motivicTspectra$ (see theorem \ref{thm.2.4.stableTspectramodelstr}).
	In this section we describe some of his constructions and results that will be necessary
	for our study of the multiplicative properties of the slice filtration.
	
\begin{defi}
		\label{def.2.6.symmetric-groups}
	For $n\geq 0$, let $\nsymmgroup$ denote the symmetric group on $n$ letters where
	$\Sigma _{0}$ is by definition the group with only one element.
	
	The \emph{$(q,p)$-shuffle} $c_{q,p}\in \Sigma _{p+q}$ is given by the following formula:
		$$c_{q,p}(i)=
		\begin{cases}
			i+p & \text{if } \: 1\leq i\leq q.\\
			i-q & \text{if } \: q+1\leq i\leq p+q.
		\end{cases}
		$$
\end{defi}

\begin{defi}[Jardine, cf. \cite{MR1787949}]
		\label{defi2.6.symmetric-T-spectra}
	\begin{enumerate}
		\item	A \emph{symmetric $T$-spectrum} $X$ is a collection
					of pointed simplicial presheaves  $(X^{n})_{n\geq 0}$
					on the smooth Nisnevich site $\smoothS$, together
					with: 
						\begin{enumerate}
							\item	Left actions
											$$\xymatrix{\nsymmgroup \times X^{n} \ar[r]& X^{n}}
											$$

							\item Bonding maps
											$$\xymatrix{T\wedge X^{n} \ar[r]^-{\sigma ^{n}}& X^{n+1}}
											$$
										such that the iterated composition
											$$\xymatrix{T^{r}\wedge X^{n} \ar[r]& X^{n+r}}
											$$
										is $\Sigma _{r}\times \nsymmgroup$-equivariant for $r\geq 1$ and
										$n\geq 0$.
						\end{enumerate}
		\item	A map $f:X\rightarrow Y$ of symmetric $T$-spectra
					is a collection of maps
					$$\xymatrix{X^{n} \ar[r]^-{f^{n}} & Y^{n}}
					$$
					in $\pointedmotivic$ satisfying the following conditions:
						\begin{enumerate}
							\item Compatibility
										with the bonding maps, i.e.
										the following diagram:
											$$\xymatrix{T\wedge X^{n} \ar[rr]^-{id\wedge f^{n}} \ar[d]_-{\sigma^{n}}&& 
																	T\wedge Y^{n} \ar[d]^-{\sigma ^{n}}\\
																	X^{n+1} \ar[rr]_-{f^{n+1}}&& Y^{n+1}}
											$$
										commutes for all $n\geq 0$
							\item	$f^{n}$ is $\nsymmgroup$-equivariant.
						\end{enumerate}
		\item	With the previous definitions we get a category,
					called the category of \emph{symmetric $T$-spectra} which
					will be denoted by $\symmTspectra$.
	\end{enumerate}
\end{defi}

\begin{exam}
		\label{exam.2.6.sphere-T-spectrum}
	Given any pointed simplicial presheaf $X$ in $\pointedmotivic$, the
	$T$-spectrum $F_{0}(X)$ has the structure of a symmetric
	$T$-spectrum; where the left action of $\nsymmgroup$ on
	$F_{0}(X)^{n}=T^{n}\wedge X$ is given by the permutation of
	the $T$ factors.
	
	In particular if we take $X=S^{0}$, we get the \emph{sphere $T$-spectrum};
	which will be denoted by $\symmspherespectrum$.
\end{exam}

	The category of symmetric $T$-spectra has a  simplicial
	structure similar to the one that exists for $T$-spectra,
	which is induced from the one
	on pointed simplicial presheaves.
	
	Given a symmetric $T$-spectrum $X$, the tensor objects are defined as follows:
	$$\xymatrix@R=.5pt{X\wedge -:\simpsets \ar[r]& \symmTspectra \\
										 K \ar@{|->}[r]& X\wedge K}
	$$
	where $(X\wedge K)^{n}=X^{n}\wedge K_{+}$ which has an action of $\nsymmgroup$
	induced by the one in $X^{n}$ and the functor $-\wedge K_{+}$,
	and with bonding maps 
	$$\xymatrix{T\wedge (X^{n}\wedge K_{+}) \ar[r]^-{\cong}& 
							(T\wedge X^{n})\wedge K_{+} \ar[rr]^-{\sigma^{n}\wedge id_{K_{+}}} && X^{n+1}\wedge K_{+}}
	$$
	The simplicial functor in two variables is:
	$$\xymatrix@R=.5pt{Map\: _{\Sigma}(-,-):(\symmTspectra)^{op}\times \symmTspectra \ar[r] & \simpsets \\
										(X,Y) \ar@{|->}[r]& Map\: _{\Sigma}(X,Y)}
	$$
	where $Map\: _{\Sigma}(X,Y)_{n}=\Hom _{\symmTspectra}(X\wedge \Delta ^{n}_{+},Y)$,
	and finally for any symmetric $T$-spectrum $Y$ we have the following functor
	$$\xymatrix@R=.5pt{Y^{-}:\simpsets \ar[r]& (\symmTspectra )^{op}\\
										K \ar@{|->}[r] & Y^{K}}
	$$
	where $(Y^{K})^{n}=(Y^{n})^{K_{+}}$ which has an  action of $\nsymmgroup$ induced by
	the one in $Y^{n}$ and the  $K_{+}$-loops functor,
	and with
	bonding maps
	$$\xymatrix{T\wedge(Y^{n})^{K_{+}} \ar[r]^-{\alpha}& (T\wedge Y^{n})^{K_{+}} \ar[r]^-{(\sigma ^{n})_{*}}& (Y^{n+1})^{K_{+}}}
	$$
	where for $U\in (\smoothS)$, $\alpha(U)$ is adjoint to
	$$\xymatrix{T(U)\wedge (Y^{n}(U))^{K_{+}}\wedge K_{+} \ar[rr]^-{id_{T(U)}\wedge ev_{K_{+}}} && T(U)\wedge Y^{n}(U)}
	$$

	In a similar way, it is possible to promote the action
	of $\pointedmotivic$ on the category of $T$-spectra to the
	category of symmetric $T$-spectra, i.e.
	the category of symmetric $T$-spectra $\symmTspectra$ has the structure of
	a closed $\pointedmotivic$-module, which
	is obtained by extending the symmetric monoidal structure for $\pointedmotivic$
	levelwise.
	
	The bifunctor giving the adjunction of two variables is defined
	as follows:
		$$\xymatrix@R=.5pt{- \wedge -:\symmTspectra \times \pointedmotivic \ar[r]& \symmTspectra \\
											(X,A) \ar@{|->}[r]& X\wedge A}
		$$
	with $(X\wedge A)^{n}=X^{n}\wedge A$ which has an action of $\nsymmgroup$
	induced by the one in $X^{n}$ and the functor $-\wedge A$,
	and with bonding maps 
		$$\xymatrix{T\wedge (X^{n}\wedge A) \ar[r]^-{\cong}& (T\wedge X^{n})\wedge A \ar[rr]^-{\sigma ^{n}\wedge id _{A}}
							&& X^{n+1}\wedge A}
		$$
	
	The adjoints are given by:
	$$\xymatrix@R=.5pt{\Omega _{-}-:\pointedmotivic ^{op}\times \symmTspectra \ar[r]& \symmTspectra \\
										(A,X) \ar@{|->}[r]& \Omega _{A}X}
	$$
	
	$$\xymatrix@R=.5pt{\inthomsymmTspectrapresheaf (-,-):(\symmTspectra )^{op}\times \symmTspectra \ar[r]& \pointedmotivic \\
										(X,Y) \ar@{|->}[r]& \inthomsymmTspectrapresheaf (X,Y)}
	$$
	where $(\Omega _{A}X)^{n}=\Omega _{A}X^{n}$ which has an  action of $\nsymmgroup$ induced by
	the one in $X^{n}$ and the  $A$-loops functor, with
	bonding maps 
	$T\wedge(\Omega _{A}X^{n})\rightarrow \Omega _{A}X^{n+1}$  adjoint to
		$$\xymatrix{T\wedge (\Omega _{A}X^{n})\wedge A \ \ar[rr]^-{id\wedge ev_{A}}
								&& T\wedge X^{n} \ar[r]^{\sigma ^{n}}& X^{n+1}}$$
	and $\inthomsymmTspectrapresheaf (X,Y)$ is the following pointed simplicial presheaf
	on $\smoothS$:
	$$\xymatrix@R=.5pt{\inthomsymmTspectrapresheaf (X,Y):(\smoothS \times \Delta)^{op} \ar[r] &
										Sets \\
										(U,n) \ar@{|->}[r] & \Hom _{\symmTspectra}(X \wedge (\Delta ^{n}_{U})_{+},Y)}
	$$

	The main difference between the categories of $T$-spectra and symmetric $T$-spectra
	is that the latter has a closed symmetric monoidal structure, i.e. it is
	possible to construct the smash product of two symmetric $T$-spectra.
	
\begin{defi}[cf. \cite{MR1787949}]
		\label{def.2.6.symmetric-sequences}
	\begin{enumerate}
		\item	A \emph{symmetric sequence} $X$ is a collection
					of pointed simplicial presheaves  $(X^{n})_{n\geq 0}$
					on the smooth Nisnevich site $\smoothS$, together
					with left actions
						$$\xymatrix{\nsymmgroup \times X^{n} \ar[r]& X^{n}}
						$$
		\item	A map $f:X\rightarrow Y$ of symmetric sequences consists of
					a collection of $\nsymmgroup$-equivariant maps
						$$\xymatrix{f^{n}:X^{n}\ar[r] & Y^{n}}
						$$
					in $\pointedmotivic$.
		\item	With these definitions we get a category, called the
					category of symmetric sequences which will be denoted by
					$\symmsequences$.
	\end{enumerate}
\end{defi}

\begin{defi}
		\label{def.2.6.product-symmetric-sequences}
	Let $X$ and $Y$ be two symmetric sequences.  Then the
	\emph{product} $X\otimes Y$ is given by the following
	symmetric sequence:
		$$(X\otimes Y)^{n}=\bigvee _{p+q=n}\nsymmgroup \otimes _{\Sigma _{p}\times \Sigma _{q}}X^{p}\wedge Y^{q}
		$$
\end{defi}

\begin{rmk}
		\label{rmk.2.6.symmTspectra===Tmodule-symmetric-sequence}
	A symmetric $T$-spectrum $X$ can be identified with a symmetric
	sequence $X$ equipped with a module structure over the sphere spectrum, i.e. with
	a map of symmetric sequences:
		$$\xymatrix{\symmspherespectrum \otimes X \ar[r]^-{\sigma _{X}}& X}
		$$
	satisfying the usual associativity conditions.
\end{rmk}

\begin{defi}[cf. \cite{MR1695653}]
		\label{prop.2.6.Gn-Evn-adjunction}
	For every $n\geq 0$, we have the following adjunction:
		$$\xymatrix{(G_{n},Ev_{n},\varphi):\pointedmotivic \ar[r]& \symmsequences}
		$$
	where $Ev_{n}$ is the $n$-evaluation functor
		$$\xymatrix@R=.5pt{Ev_{n}:\symmsequences \ar[r]& \pointedmotivic \\
											 X \ar@{|->}[r]& X^{n}}
		$$
	and $G_{n}$ is the $n$-\emph{free symmetric sequence functor}:
		$$\xymatrix@R=.5pt{G_{n}:\pointedmotivic \ar[r]& \symmsequences \\
											 X \ar@{|->}[r]& G_{n}(X)}
		$$
	where
		$$G_{n}(X)^{m}=
		\begin{cases}
			\ast & \text{if }\: m\neq n.\\
			\bigvee _{\sigma \in \nsymmgroup}X & \text{if }\: m=n.
		\end{cases}
		$$
\end{defi}

\begin{defi}[cf. \cite{MR1787949}]
		\label{def.2.6.Fnsigma-Evn-adjunction}
	For every $n\geq 0$, we have the following adjunction:
		$$\xymatrix{(F_{n}^{\Sigma},Ev_{n},\varphi):\pointedmotivic \ar[r]& \symmTspectra}
		$$
	where $Ev_{n}$ is the $n$-evaluation functor
		$$\xymatrix@R=.5pt{Ev_{n}:\symmTspectra \ar[r]& \pointedmotivic \\
											 X \ar@{|->}[r]& X^{n}}
		$$
	and $F_{n}^{\Sigma}$ is the $n$-\emph{free symmetric $T$-spectrum functor}:
		$$\xymatrix@R=.5pt{F_{n}^{\Sigma}:\pointedmotivic \ar[r]& \symmTspectra \\
											 X \ar@{|->}[r]& \symmspherespectrum \otimes G_{n}(X)}
		$$
\end{defi}

\begin{defi}[cf. \cite{MR1787949}]
		\label{def.2.6.smashproduct--symmTspectra}
	Let $X$ and $Y$ be two symmetric $T$-spectra.
	Then the \emph{smash product} $X\wedge Y$ is given by
	the colimit of the following diagram
		$$\xymatrix{\symmspherespectrum \otimes X\otimes Y \ar@<1ex>[rr]^-{\sigma _{X}\otimes id} \ar@<-1ex>[rr]&& X\otimes Y}
		$$
	where the bottom arrow is the following composition
		$$\xymatrix{\symmspherespectrum \otimes X\otimes Y \ar[rr]^-{t} && 
			X\otimes \symmspherespectrum \otimes Y \ar[rr]^-{id\otimes \sigma _{Y}}&& X\otimes Y}
		$$
\end{defi}	

\begin{prop}[Jardine]
		\label{prop.2.6.symmTspectra==>closed-symmetric-monoidal-category}
	The category of symmetric $T$-spectra $\symmTspectra$ has a
	closed symmetric monoidal structure where the product is given
	by the smash product described in definition \ref{def.2.6.smashproduct--symmTspectra},
	and the functor  that gives the adjunction of two variables is the following:
		$$\xymatrix@R=.5pt{\inthomsymmTspectra (-,-):(\symmTspectra)^{op}\times \symmTspectra \ar[r]&
								\symmTspectra \\
								(X,Y) \ar@{|->}[r]& \inthomsymmTspectra (X,Y)}
		$$
	where $\inthomsymmTspectra (X,Y)^{n}=\inthomsymmTspectrapresheaf (F_{n}^{\Sigma}(S^{0})\wedge X,Y)$, and the
	adjoints $\sigma ^{n}_{\ast}$ to the bonding maps are given as follows:
	Let $\zeta : F_{n+1}^{\Sigma}(T)\cong F_{n+1}^{\Sigma}(S^{0})\wedge T\rightarrow F_{n}^{\Sigma}(S^{0})$ be the adjoint
	corresponding to the inclusion determined by the identity in $\Sigma _{n+1}$:
		$$\imath _{e}:T \hookrightarrow Ev_{n+1}(F_{n}^{\Sigma}(S^{0}))=
			\Sigma _{n+1}\otimes  _{\Sigma _{1}\times \Sigma _{n}}(T\wedge \bigvee _{\sigma \in \Sigma _{n}}S^{0})
			=\bigvee _{\sigma \in \Sigma _{n+1}} T
		$$
	then $\sigma ^{n}_{\ast}$ is the following map induced by $\zeta \wedge id$:
		$$\xymatrix{\inthomsymmTspectrapresheaf (F_{n}^{\Sigma}(S^{0})\wedge X,Y)\ar[rr]^-{(\zeta \wedge id) ^{\ast}}&&
								\inthomsymmTspectrapresheaf (F_{n+1}^{\Sigma}(S^{0})\wedge T\wedge X,Y)}
		$$
	
	The twist isomorphism $\tau :X\wedge Y\rightarrow Y\wedge X$ is induced levelwise by:
		$$\xymatrix{X^{p}\wedge Y^{q} \ar[r]^-{t} \ar[d]_-{\alpha}& Y^{q}\wedge X^{p} \ar[d]^-{\alpha c_{q,p}}\\ 
								(X\otimes Y)^{p+q} & (Y\otimes X)^{p+q}}
		$$
	Finally, the unit is given by the sphere $T$-spectrum $F_{0}^{\Sigma}(S^{0})=\symmspherespectrum$.
\end{prop}
\begin{proof}
	We refer the reader to \cite[section 4.3]{MR1787949}.
\end{proof}

\begin{prop}
		\label{prop.2.6.simp-adjunc-symmspectraloops}
	Let $X,Y$ be two arbitrary symmetric $T$-spectra and let
	$A$ in $\pointedmotivic$ be an arbitrary pointed simplicial presheaf.  
	Then we have the
	following enriched adjunctions:
		\begin{equation}
					\label{diagram.2.6.adj.a}
				\xymatrix{Map(A, \inthomsymmTspectrapresheaf (X,Y))\ar[r]^-{\alpha}_-{\cong}&
								  Map\: _{\Sigma}(X\wedge A,Y) \ar[r]^-{\beta}_-{\cong}& Map\: _{\Sigma}(X,\Omega _{A}Y)}
		\end{equation}
		
		\begin{equation}
					\label{diagram.2.6.adj.b}
				\xymatrix{\inthomprepointed (A, \inthomsymmTspectrapresheaf (X,Y)) \ar[r]^-{\delta}_-{\cong}& 
									\inthomsymmTspectrapresheaf(X\wedge A,Y) \ar[r]^-{\epsilon}_-{\cong}& \inthomsymmTspectrapresheaf (X,\Omega _{A}Y)}
		\end{equation}
		
		\begin{equation}
					\label{diagram.2.6.adj.c}
				\xymatrix{\inthomsymmTspectra (X\wedge A,Y)\ar[r]^-{\gamma}_-{\cong}& \inthomsymmTspectra (X, \Omega _{A}Y)}
		\end{equation}
	where the maps in (\ref{diagram.2.6.adj.a}) are isomorphisms of simplicial sets, the maps in (\ref{diagram.2.6.adj.b})
	are isomorphisms of simplicial presheaves, and the map in (\ref{diagram.2.6.adj.c}) is an isomorphism of symmetric
	$T$-spectra.
\end{prop}
\begin{proof}
	We consider first the simplicial adjunctions:
	To any $n$-simplex $t$ in $Map(A, \inthomsymmTspectrapresheaf (X,Y))$
		$$\xymatrix{A\wedge \Delta ^{n} \ar[r]^-{t}& \inthomsymmTspectrapresheaf (X,Y)}
		$$
	associate the following $n$-simplex in $Map\: _{\Sigma}(X\wedge A,Y)$:
	$$\xymatrix{ X \wedge A \wedge \Delta ^{n} 
							 \ar[r]^-{\alpha (t)}& Y}
	$$
	corresponding to the adjunction between $X\wedge -$ and $\inthomsymmTspectrapresheaf (X,-)$.
	
	To any $n$-simplex $t$ in $Map\: _{\Sigma}(X\wedge A,Y)$
	$$\xymatrix{\Delta ^{n}\wedge X\wedge A \ar[r]^-{\cong}& X\wedge A\wedge \Delta ^{n} \ar[r]^-{t}& Y}
	$$
	associate the following $n$-simplex in $Map\: _{\Sigma}(X,\Omega _{A}Y)$:
	$$\xymatrix{ X\wedge \Delta ^{n} \ar[r]^-{\cong}& \Delta ^{n}\wedge X \ar[r]^-{\beta (t)}& \Omega _{A}Y}
	$$
	corresponding to the adjunction between $-\wedge A$ and $\Omega _{A}$.
	
	We consider now the isomorphisms of simplicial presheaves:
	To any simplex $s$
	in $\inthomprepointed (A, \inthomsymmTspectrapresheaf (X,Y))$
	$$\xymatrix{ A\wedge \Delta ^{n}_{U}  \ar[r]^-{s}& \inthomsymmTspectrapresheaf (X,Y)}
	$$
	we associate the following simplex in $\inthomsymmTspectrapresheaf (X\wedge A, Y)$
	$$\xymatrix{X\wedge A\wedge \Delta ^{n}_{U}  
							 \ar[r]^-{\delta (s)}& Y}
	$$
	corresponding to the adjunction between $X\wedge -$ and $\inthomsymmTspectrapresheaf (X,-)$.
	
	To any simplex $s$
	in $\inthomsymmTspectrapresheaf (X\wedge A,Y)$
	$$\xymatrix{ X\wedge \Delta ^{n}_{U}\wedge A \ar[r]^{\cong} & X\wedge A\wedge \Delta ^{n}_{U} \ar[r]^-{s}& Y}
	$$
	we associate the following simplex in $\inthomsymmTspectrapresheaf (X,\Omega _{A}Y)$
	$$\xymatrix{X\wedge \Delta ^{n}_{U} \ar[r]^-{\epsilon (s)}& \Omega _{A}Y}
	$$
	corresponding to the adjunction between $-\wedge A$ and $\Omega _{A}$.
	
	Finally, we consider the isomorphism of symmetric $T$-spectra:
	Using the adjunction given by $\epsilon$ in (\ref{diagram.2.6.adj.b}), we get
	for every $n\geq 0$
	the following commutative diagram, where the vertical maps are isomorphisms of
	simplicial presheaves:
		$$\xymatrix{\inthomsymmTspectrapresheaf (F_{n}^{\Sigma}(S^{0})\wedge X\wedge A,Y) 
								\ar[rr]^-{(\alpha \wedge id_{X\wedge A})^{\ast}} \ar[d]^-{\epsilon}_-{\cong}&& 
								\inthomsymmTspectrapresheaf (F_{n+1}^{\Sigma}(T)\wedge X\wedge A,Y) \ar[d]^-{\epsilon}_-{\cong}\\
								\inthomsymmTspectrapresheaf (F_{n}^{\Sigma}(S^{0})\wedge X,\Omega _{A}Y) 
								\ar[rr]_-{(\alpha \wedge id_{X})^{\ast}}&&
								\inthomsymmTspectrapresheaf (F_{n+1}^{\Sigma}(T)\wedge X,\Omega _{A}Y)}
		$$
	By definition (see proposition \ref{prop.2.6.symmTspectra==>closed-symmetric-monoidal-category})
	the diagram above is equal to:
		$$\xymatrix{\inthomsymmTspectra (X\wedge A,Y)^{n} 
								\ar[r]^-{\sigma ^{n}_{\ast}} \ar[d]^-{\epsilon}_-{\cong}& 
								\Tloops \inthomsymmTspectra (X\wedge A,Y)^{n+1} \ar[d]^-{\epsilon}_-{\cong}\\
								\inthomsymmTspectra (X,\Omega _{A}Y)^{n} 
								\ar[r]_-{\sigma ^{n}_{\ast}}& 
								\Tloops \inthomsymmTspectra (X,\Omega _{A}Y)^{n+1}}
		$$
	This induces the isomorphism $\gamma$.
\end{proof}

\begin{prop}
		\label{prop.2.6.enriched-adjunc-symmspectraloops}
	Let $X,Y,Z$ be three arbitrary symmetric $T$-spectra.  
	Then we have the
	following enriched adjunctions:
		\begin{equation}
				\label{diagram.enriched-adjunc-symmspectraloops.a}
			\xymatrix{Map\: _{\Sigma}(X\wedge Y,Z)\ar[r]^-{\lambda}_-{\cong}&
								Map\: _{\Sigma}(X,\inthomsymmTspectra (Y,X))}
		\end{equation}
	
		\begin{equation}
				\label{diagram.enriched-adjunc-symmspectraloops.b}
			\xymatrix{\inthomsymmTspectrapresheaf (X\wedge Y,Z) \ar[r]^-{\kappa}_-{\cong}& 
								\inthomsymmTspectrapresheaf (X,\inthomsymmTspectra (Y,Z))}
		\end{equation}
	
		\begin{equation}
				\label{diagram.enriched-adjunc-symmspectraloops.c}
			\xymatrix{\inthomsymmTspectra (X\wedge Y,Z)\ar[r]^-{\mu}_-{\cong}& \inthomsymmTspectra (X, \inthomsymmTspectra (Y,Z))}
		\end{equation}
	where the map in (\ref{diagram.enriched-adjunc-symmspectraloops.a}) is an isomorphism of simplicial sets, 
	the map in (\ref{diagram.enriched-adjunc-symmspectraloops.b})
	is an isomorphism of simplicial presheaves, and the map in (\ref{diagram.enriched-adjunc-symmspectraloops.c}) 
	is an isomorphism of symmetric $T$-spectra.
\end{prop}
\begin{proof}
	We consider first the simplicial adjunctions:
	To any $n$-simplex $t$ in $Map\: _{\Sigma}(X\wedge Y,Z)$
		$$\xymatrix{\Delta ^{n}\wedge X\wedge Y \ar[r]^-{\cong}&
								X\wedge Y\wedge \Delta ^{n} \ar[r]^-{t}& Z}
		$$
	associate the following $n$-simplex in $Map\: _{\Sigma}(X,\inthomsymmTspectra (Y,Z))$:
	$$\xymatrix{ X \wedge  \Delta ^{n} \ar[r]^-{\cong}&
							\Delta ^{n}\wedge X \ar[r]^-{\lambda (t)}& \inthomsymmTspectra (Y,Z)}
	$$
	corresponding to the adjunction between $-\wedge Y$ and $\inthomsymmTspectra (Y,-)$.
	
	We consider now the isomorphisms of simplicial presheaves:
	To any simplex $s$
	in $\inthomsymmTspectrapresheaf (X\wedge Y, Z)$
	$$\xymatrix{\Delta ^{n}_{U}\wedge X\wedge Y \ar[r]^-{\cong}& 
							X\wedge Y\wedge \Delta ^{n}_{U}  \ar[r]^-{s}& Z}
	$$
	we associate the following simplex in $\inthomsymmTspectrapresheaf (X, \inthomsymmTspectra (Y,Z))$
	$$\xymatrix{X\wedge \Delta ^{n}_{U} \ar[r]^-{\cong}& 
							 \Delta ^{n}_{U}\wedge X \ar[r]^-{\kappa (s)}& \inthomsymmTspectra (Y,Z)}
	$$
	corresponding to the adjunction between $-\wedge Y$ and $\inthomsymmTspectra (Y,-)$.
	
	Finally, we consider the isomorphism of symmetric $T$-spectra:
	Using the adjunction given by $\kappa$ in (\ref{diagram.enriched-adjunc-symmspectraloops.b}), we get
	for every $n\geq 0$
	the following commutative diagram, where the vertical maps are isomorphisms of
	simplicial presheaves:
		$$\xymatrix{\inthomsymmTspectrapresheaf (F_{n}^{\Sigma}(S^{0})\wedge X\wedge Y,Z) 
								\ar[dr]^-{(\alpha \wedge id_{X\wedge Y})^{\ast}} \ar[dd]^-{\kappa}_-{\cong}&\\ 
								& \inthomsymmTspectrapresheaf (F_{n+1}^{\Sigma}(T)\wedge X\wedge Y, Z) \ar[dd]^-{\kappa}_-{\cong}\\
								\inthomsymmTspectrapresheaf (F_{n}^{\Sigma}(S^{0})\wedge X,\inthomsymmTspectra (Y,Z)) 
								\ar[dr]_-{(\alpha \wedge id_{X})^{\ast}}&\\
								& \inthomsymmTspectrapresheaf (F_{n+1}^{\Sigma}(T)\wedge X,\inthomsymmTspectra (Y,Z))}
		$$
	By definition (see proposition \ref{prop.2.6.symmTspectra==>closed-symmetric-monoidal-category})
	the diagram above is equal to:
		$$\xymatrix{\inthomsymmTspectra (X\wedge Y,Z)^{n} 
								\ar[r]^-{\sigma ^{n}_{\ast}} \ar[d]^-{\kappa}_-{\cong}& 
								\Tloops \inthomsymmTspectra (X\wedge Y,Z)^{n+1} \ar[d]^-{\kappa}_-{\cong}\\
								\inthomsymmTspectra (X,\inthomsymmTspectra (Y,Z))^{n} 
								\ar[r]_-{\sigma ^{n}_{\ast}}& 
								\Tloops \inthomsymmTspectra (X, \inthomsymmTspectra (Y,Z))^{n+1}}
		$$
	This induces the isomorphism $\mu$.
\end{proof}

	The following proposition will have remarkable consequences in
	our study of the multiplicative properties for Voevodsky's slice
	filtration.

\begin{prop}[Jardine]
		\label{prop.2.6.free-symmetric-functors--compatiblewithsmash}
	Let $A$, $B$ be two arbitrary pointed simplicial presheaves in $\pointedmotivic$.  
	Then we have an isomorphism:
		$$\xymatrix{F_{n}^{\Sigma}(A)\wedge F_{m}^{\Sigma}(B)\ar[r]^-{\cong}& F_{m+n}^{\Sigma}(A\wedge B)}
		$$
	which is natural in $A$ and $B$.
\end{prop}
\begin{proof}
	We refer the reader to \cite[corollary 4.18]{MR1787949}.
\end{proof}

	For the construction of the motivic stable model structure on the category of
	$T$-spectra, it was necessary to introduce the projective and injective model structures
	(see theorem \ref{thm.2.4.stableTspectramodelstr}).  In \cite{MR1787949}, Jardine
	considers an injective model structure for symmetric $T$-spectra as a preliminary step 
	in the construction of a model
	structure  which turns out to be Quillen equivalent to $\motivicTspectra$.  
	We will also
	need to consider a projective model structure for symmetric $T$-spectra,
	in order to show that
	this stable model structure for symmetric $T$-spectra is cellular.
	
\begin{defi}
		\label{def.2.6.levelcofs-weakequivs}
	Let $f:X\rightarrow Y$ be a map of symmetric $T$-spectra.
	We say that $f$ is a level cofibration (respectively level fibration, level weak equivalence), 
	if for every
	$n\geq 0$, the map $f^{n}:X^{n}\rightarrow Y^{n}$ is
	a cofibration (respectively a fibration, a weak equivalence)
	in $\pointedmotivic$.
\end{defi}	
	
	In proposition \ref{prop.2.3.cell-prop-pointedmotcat}
	we used $I_{M_{\ast}}$ and $J_{M_{\ast}}$ to denote the sets of generating
	cofibrations and trivial cofibrations for $\pointedmotivic$.
	
\begin{thm}[Hovey]
		\label{thm2.6.projstabmodstrsymmTspectra}
	There exists a cofibrantly generated model structure for the category
	$\symmTspectra$ of symmetric $T$-spectra with the following choices:
	\begin{enumerate}
		\item The weak equivalences are the level weak equivalences.
		
		\item	The set $I$ of generating cofibrations is 
					$$I=\bigcup _{n\geq 0}F_{n}^{\Sigma}(I_{M_{\ast}})$$
		
		\item	The set $J$ of generating trivial cofibrations is
					$$J=\bigcup _{n\geq 0}F_{n}^{\Sigma}(J_{M_{\ast}})$$
	\end{enumerate}
	This model structure will be called the \emph{projective
	model structure} for symmetric $T$-spectra.  Furthermore,
	the projective model structure is left proper and simplicial.
\end{thm}
\begin{proof}
	Proposition \ref{prop.2.3.cell-prop-pointedmotcat} implies that
	the model category $\pointedmotivic$ 
	is in particular pointed, proper, simplicial and
	symmetric monoidal.  We also have that every pointed simplicial presheaf
	in $\pointedmotivic$ is cofibrant.
	Then the result follows immediately
	from theorems 8.2 and 8.3 in \cite{MR1860878}.
\end{proof}

\begin{rmk}
		\label{rmk.2.6.classifying-proj-fibrations}
	Let $f:X\rightarrow Y$ be a map of symmetric $T$-spectra.
	\begin{enumerate}
		\item	$f$ is a fibration in $\symmTspectra$ equipped with the projective
					model structure if and only if $f$ is a level fibration.
		\item	$f$ is a trivial fibration in $\symmTspectra$ equipped with the
					projective model structure if and only if $f$ is both a level
					fibration and a level weak equivalence.
	\end{enumerate}
\end{rmk}

	It follows directly from the definition that every symmetric $T$-spectrum
	after forgetting the $\nsymmgroup$-actions becomes a  
	$T$-spectrum in $\Tspectra$.  Therefore we get a functor:
		$$\xymatrix{U:\symmTspectra \ar[r]& \Tspectra}
		$$
	It turns out that this forgetful functor has a left adjoint.
	
\begin{defi}[Jardine, cf. \cite{MR1787949}]
		\label{def.2.6.layerfiltration}
	Let $X$ be an arbitrary $T$-spectrum in $\Tspectra$.
	Then $X$ has a natural filtration $\{ L_{n}X\}_{n\geq 0}$
	called the \emph{layer filtration}, where $L_{n}X$ is
	defined as
		$$X^{0}, X^{1}, \ldots , X^{n}, T\wedge X^{n}, T^{2}\wedge X^{n},\ldots
		$$
	and furthermore
		$$X\cong \varinjlim L_{n}X
		$$
		
	It is also possible to give an inductive definition for the layers $L_{n}X$
	using the following pushout diagrams (see definition \ref{def.2.4.infinitesuspension+shift}):
		$$\xymatrix{F_{n+1}(T\wedge X^{n}) \ar[r] \ar[d]& L_{n}X \ar[d]\\
								F_{n+1}(X^{n+1}) \ar[r]& L_{n+1}X}
		$$
\end{defi}
	
\begin{prop}[Jardine]
		\label{prop.2.6.symmetrizationfunctorV}
	We have the following adjunction
		$$\xymatrix{(V,U,\varphi):\Tspectra \ar[r]& \symmTspectra}
		$$
	The functor $V$ is called the \emph{symmetrization functor} and
	is defined as follows:
	\begin{enumerate}
		\item	For every pointed simplicial presheaf
					$X$ on the smooth Nisnevich site $\nissite$ we have
						$$V(F_{n}(X))=F_{n}^{\Sigma}(X)$$
		\item $V$ is constructed inductively using the layer filtration
					(see definition \ref{def.2.6.layerfiltration}) together with
					the following pushout diagrams (see definition \ref{def.2.6.Fnsigma-Evn-adjunction}):
						$$\xymatrix{F_{n+1}^{\Sigma}(T\wedge X^{n}) \ar[r] \ar[d]& V(L_{n}X)\ar[d]\\
												F_{n+1}^{\Sigma}(X^{n+1}) \ar[r]& V(L_{n+1}X)}
						$$
		\item	Finally, $V(X)=\varinjlim V(L_{n}X)$
	\end{enumerate}
\end{prop}
\begin{proof}
	We refer the reader to \cite[p. 507]{MR1787949}
\end{proof}

\begin{prop}
		\label{prop.2.6.enriched-symmetrization-adjunction}
	The adjunction
		$$\xymatrix{(V,U,\varphi):\Tspectra \ar[r]& \symmTspectra}
		$$
	is enriched in the categories of simplicial sets and pointed
	simplicial presheaves on $\nissite$, i.e. for every $T$-spectrum $X$ and
	for every symmetric $T$-spectrum $Y$ we have the following natural isomorphisms:
		$$\xymatrix{Map\: _{\Sigma}(VX,Y)\ar[r]_-{\cong}^-{\epsilon}& Map(X,UY)}
		$$
		
		$$\xymatrix{\inthomsymmTspectrapresheaf (VX,Y) \ar[r]_-{\cong}^-{\eta}& \inthomspectrapresheaf (X,UY)}
		$$
\end{prop}
\begin{proof}
	We consider first the simplicial isomorphism:
	Given any $n$-simplex $t$ in $Map\: _{\Sigma} (VX,Y)$
		$$\xymatrix{VX\wedge \Delta ^{n}\ar[r]^-{t}& Y}
		$$
	consider the map corresponding to the adjuntion between $-\wedge \Delta ^{n}$ and
	$-^{\Delta ^{n}}$ in $\symmTspectra$
		$$\xymatrix{VX \ar[r]^-{t'}& Y^{\Delta ^{n}}}
		$$
	Now use the adjunction between $V$ and $U$ to get the map:
		$$\xymatrix{X \ar[r]^-{t''}& U(Y^{\Delta ^{n}})=(UY)^{\Delta ^{n}}}
		$$
	and finally use the adjunction between $-\wedge \Delta ^{n}$ and
	$-^{\Delta ^{n}}$ in $\Tspectra$ to get the associated $n$-simplex $\epsilon (t)$
	in $Map(X,UY)$:
		$$\xymatrix{X\wedge \Delta ^{n}\ar[r]^-{\epsilon (t)}& UY}
		$$
		
	We consider now the isomorphism of simplicial presheaves:
	Given any simplex $s$ in $\inthomsymmTspectrapresheaf (VX,Y)$
		$$\xymatrix{VX\wedge \Delta ^{n}_{W}\ar[r]^-{s}& Y}
		$$
	consider the map corresponding to the adjunction between
	$-\wedge \Delta ^{n}_{W}$ and $\Omega _{\Delta ^{n}_{W}}-$ in $\symmTspectra$
		$$\xymatrix{VX \ar[r]^-{s'}& \Omega _{\Delta ^{n}_{W}}Y}
		$$
	Now use the adjunction between $V$ and $U$ to get the map:
		$$\xymatrix{X \ar[r]^-{s''}& U(\Omega _{\Delta ^{n}_{W}}Y)}=\Omega _{\Delta ^{n}_{W}}UY
		$$
	and finally use the adjunction between $-\wedge \Delta ^{n}_{W}$ and $\Omega _{\Delta ^{n}_{W}}$
	in $\Tspectra$ to get the associated simplex $\eta (s)$ in $\inthomspectrapresheaf (X,UY)$:
		$$\xymatrix{X\wedge \Delta ^{n}_{W}\ar[r]^-{\eta (s)}& UY}
		$$
\end{proof}

We say that a map $f:X\rightarrow Y$ of symmetric $T$-spectra is
	an \emph{injective fibration} if it has the right lifting property
	with respect to the class of maps which are both level cofibrations and
	level weak equivalences.
	
\begin{thm}[Jardine]
		\label{thm.2.6.injmodstrsymmTspectra}
	There exists a 
	model structure for the category $\symmTspectra$ of symmetric $T$-spectra
	with the following choices:
	\begin{enumerate}
		\item The weak equivalences are the level weak equivalences.
		
		\item	The cofibrations are the level cofibrations.
		
		\item	The fibrations are the injective fibrations.					
	\end{enumerate}
	This model structure will be called the \emph{injective
	model structure} for symmetric $T$-spectra.  Furthermore, the
	injective model structure is proper, simplicial, and cofibrantly
	generated with the following sets $I$, $J$ of generating
	cofibrations and trivial cofibrations, respectively
	(see theorem \ref{thm.2.4.injmodstrspt}):
		\begin{enumerate}
			\item	The set $I$ of generating cofibrations is 
							$$I=\{V(i):VA\rightarrow VB\}
							$$
						where $i$ satisfies the following conditions:
							\begin{enumerate}
								\item $i$ is a level cofibration in $\Tspectra$.
								\item The codomain $B$ of $i$ is
											$\kappa$-bounded.								
							\end{enumerate}
			\item	The set $J$ of generating trivial cofibrations is
							$$J=\{V(j):V(A)\rightarrow V(B) \}
							$$
						where $j$ satisfies the following conditions:
							\begin{enumerate}
								\item	$j$ is a level trivial cofibration in $\Tspectra$.
								\item	The codomain $B$ of $j$ is 
											$\kappa$-bounded.
							\end{enumerate}
		\end{enumerate}
\end{thm}
\begin{proof}
	We refer the reader to \cite[theorem 4.2]{MR1787949}.
\end{proof}

\begin{rmk}
		\label{rmk.2.6.id-leftQuillen-proj-inj-symmspectra}
	The identity functor on $\symmTspectra$
	induces a left Quillen functor from the
	projective model structure to the injective model
	structure.
\end{rmk}

\begin{defi}
		\label{def.2.6.symmetricstableweakequiv}
	\begin{enumerate}
		\item \label{def.2.6.symmetricstableweakequiv.a}Let $Z$ be a symmetric $T$-spectrum.  We say that $Z$
					is \emph{injective stably fibrant} if $Z$ satisfies the following
					conditions:
					\begin{enumerate}
						\item	$Z$ is fibrant in $\symmTspectra$ equipped with the injective
									model structure.
						\item	$UZ$ is fibrant in $\motivicTspectra$.
					\end{enumerate}
		\item \label{def.2.6.symmetricstableweakequiv.b}Let $f:X\rightarrow Y$ be a map of symmetric $T$-spectra.  We say that
					$f$ is a \emph{stable weak equivalence} if for every injective stably fibrant
					symmetric $T$-spectrum $Z$, the induced map
						$$\xymatrix{Map\: _{\Sigma}(Y,Z)\ar[r]^{f^{\ast}}& Map\: _{\Sigma}(X,Z)}
						$$
					is a weak equivalence of simplicial sets.
		\item	\label{def.2.6.symmetricstableweakequiv.c}Let $f:X\rightarrow Y$ be a map of symmetric $T$-spectra.  We say that
					$f$ is a \emph{stable fibration} if $Uf$ is a fibration in $\motivicTspectra$
					(see theorem 
					\ref{thm.2.4.stableTspectramodelstr}).
	\end{enumerate}
\end{defi}

	In theorem \ref{thm.2.5.cellularity-motivic-stable-str}
	we used $I^{T}_{M_{\ast}}$ and $J^{T}_{M_{\ast}}$ to denote the sets of generating
	cofibrations and trivial cofibrations for $\motivicTspectra$.

\begin{thm}[Jardine]
		\label{thm2.6.stablesymmetricmodelstructure}
	There exists a model structure for the category
	$\symmTspectra$ of symmetric $T$-spectra with the following choices:
	\begin{enumerate}
		\item The weak equivalences are the stable weak equivalences.
		
		\item	The cofibrations are the projective cofibrations
					(see theorem \ref{thm2.6.projstabmodstrsymmTspectra}), i.e.
					they are generated by the
					set  
						$$\bigcup _{n\geq 0}F_{n}^{\Sigma}(I_{M_{\ast}})
						=V(I^{T}_{M_{\ast}})$$
		
		\item	The fibrations are the stable fibrations.
	\end{enumerate}
	This model structure will be called \emph{motivic symmetric stable}, 
	and the category of symmetric $T$-spectra, 
	equippped with the motivic symmetric stable model structure
	will be denoted by $\motivicsymmTspectra$.  Furthermore,
	$\motivicsymmTspectra$ is a proper and simplicial model category.
\end{thm}
\begin{proof}
	We refer the reader to \cite[proposition 4.4 and theorem 4.15]{MR1787949}.
\end{proof}

\begin{rmk}
		\label{rmk.2.6.classifying-stable-trivial-fibrations}
	Let $p:X\rightarrow Y$ be a map of symmetric $T$-spectra.
	Then
	$p$ is a trivial fibration in $\motivicsymmTspectra$  if and only if
					$Up$ is a trivial fibration in $\motivicTspectra$.
\end{rmk}

\begin{prop}
		\label{prop.2.6.motsymmTspectra-pointedpresheaf-modelcat}
	$\motivicsymmTspectra$ is a $\pointedmotivic$-model category
	(see definition \ref{def.module-modcats}).
\end{prop}
\begin{proof}
	Condition (\ref{def.module-modcats.b}) in definition \ref{def.module-modcats}
	follows automatically since the unit in $\pointedmotivic$ is cofibrant.
	It remains to show that
		$$\xymatrix{-\wedge -:\symmTspectra \times \pointedmotivic \ar[r]& \symmTspectra}
		$$
	is a Quillen bifunctor.  By lemma \ref{lem.cond-Quillen-bifunc} it is enough
	to prove the following claim:
	
	Given a cofibration $i:A\rightarrow B$ in $\pointedmotivic$ and a fibration
	$p:X\rightarrow Y$ in $\motivicsymmTspectra$, then the map
		$$\xymatrix{\Omega _{B}X \ar[rr]^-{(i^{\ast},p_{\ast})}&& \Omega _{B}Y\times _{\Omega _{A}Y}\Omega _{A}X}
		$$
	is a fibration in $\motivicsymmTspectra$ which is trivial if either $i$ or $p$ is
	a weak equivalence.
	
	But this follows immediately from the following facts:
		\begin{enumerate}
			\item	A map of symmetric $T$-spectra $f:X\rightarrow Y$ is a fibration
						(respectively a trivial fibration) in $\motivicsymmTspectra$
						if and only if $Uf:UX\rightarrow UY$
						is a fibration (respectively a trivial fibration) in $\motivicTspectra$.
			\item	For every symmetric $T$-spectrum $X$ and for any pointed simplicial
						presheaf $A$ in $\pointedmotivic$, we have that $U(\Omega _{A}X)=\Omega _{A}UX$, where the
						right hand side denotes the action of $\pointedmotivic$ in $\motivicTspectra$.
			\item	$\motivicTspectra$ is a $\pointedmotivic$-model category 
						(see proposition \ref{prop.2.4.motTspectra-pointedpresheaf-modelcat}).
		\end{enumerate}
\end{proof}

\begin{cor}
		\label{cor.2.6.smash-Quillenfunctor}
	For every pointed simplicial presheaf $A\in \pointedmotivic$,
	the adjunction
		$$\xymatrix{(-\wedge A, \Omega _{A}-, \varphi ):\motivicsymmTspectra \ar[r]& \motivicsymmTspectra}
		$$
	is a Quillen adjunction.
\end{cor}
\begin{proof}
	We have that every pointed simplicial presheaf is cofibrant in $\pointedmotivic$.  
	Then the result follows from
	proposition \ref{prop.2.6.motsymmTspectra-pointedpresheaf-modelcat}.
\end{proof}

\begin{thm}[Jardine]
		\label{thm.2.6.Tsusp==QuillenequivonSymmTspectra}
	Let $T=S^{1}\wedge \gm \in \pointedmotivic$.  Then
	the Quillen adjunction:
		$$\xymatrix{(-\wedge T, \Omega _{T} , \varphi):\motivicsymmTspectra
								\ar[r]& \motivicsymmTspectra}
		$$
	is a Quillen equivalence.
\end{thm}
\begin{proof}
	Let $\eta$, $\epsilon$ denote the unit and counit of the
	adjunction $(-\wedge T, \Omega _{T}, \varphi)$.
	By proposition 1.3.13 in \cite{MR1650134}, it suffices to check
	that the following conditions hold:
		\begin{enumerate}
			\item \label{thm.2.6.Tsusp==QuillenequivonSymmTspectra.a}For every cofibrant
						symmetric $T$-spectrum $A$ in $\motivicsymmTspectra$, the following composition
							$$\xymatrix{A\ar[r]^-{\eta _{A}}& \Omega _{T}(T\wedge A)
								\ar[rr]^-{\Omega _{T}(R^{T\wedge A})}&& \Omega _{T}R(T\wedge A)}
							$$
						is a weak equivalence in $\motivicsymmTspectra$, where $R$ denotes a fibrant replacement functor
						in $\motivicsymmTspectra$.
			\item \label{thm.2.6.Tsusp==QuillenequivonSymmTspectra.b}For every fibrant
						symmetric $T$-spectrum $X$ in $\motivicsymmTspectra$, the following composition
							$$\xymatrix{T\wedge Q(\Omega _{T}X)\ar[rr]^-{id\wedge Q^{\Omega _{T}X}}&& 
													T\wedge (\Omega _{T}X)\ar[r]^-{\epsilon _{X}}& X}
							$$
						is a weak equivalence in $\motivicsymmTspectra$, where $Q$ denotes a cofibrant replacement functor
						in $\motivicsymmTspectra$.
		\end{enumerate}
		
	(\ref{thm.2.6.Tsusp==QuillenequivonSymmTspectra.a}): Follows directly from
	corollary 4.26 in \cite{MR1787949}.
	
	(\ref{thm.2.6.Tsusp==QuillenequivonSymmTspectra.b}):  By construction
	the map $Q^{\Omega _{T}X}:Q(\Omega _{T}X)\rightarrow \Omega _{T}X$ is
	a weak equivalence in $\motivicsymmTspectra$.  Therefore by lemma 4.25 in \cite{MR1787949},
	we have that $id\wedge Q^{\Omega _{T}X}$ is also a weak equivalence
	in $\motivicsymmTspectra$.  Then by the two out of three property for weak equivalences,
	it suffices to show that $\epsilon _{X}$ is a weak equivalence in $\motivicsymmTspectra$.
	
	Since $X$ is fibrant in $\motivicsymmTspectra$,
	it follows that $UX$ is fibrant in $\motivicTspectra$.  
	Therefore by lemma \ref{lem.2.4.fibrant-objects-stablemodstr}(\ref{lem.2.4.fibrant-objects-stablemodstr.b}) 
	we have that $UX$ is in particular level fibrant.
	Then by corollary 3.16 in \cite{MR1787949} it follows that the map:
		$$\xymatrix{\epsilon _{UX}:T\wedge (\Omega _{T}UX)\ar[r]& UX}
		$$
	is a weak equivalence in $\motivicTspectra$, but this is just $U(\epsilon _{X})$.
	Hence by proposition 4.8 in \cite{MR1787949}, we have that $\epsilon _{X}$
	is a weak equivalence in $\motivicsymmTspectra$, as we wanted.
\end{proof}

\begin{prop}[Jardine]
		\label{prop.2.6.symmTspectra-monoidalmodelcategory}
	$\motivicsymmTspectra$ is a symmetric monoidal model category (with respect to the smash product of symmetric
	$T$-spectra) in the sense of Hovey
	(see definition \ref{def.mon-mod-cats}).
\end{prop}
\begin{proof}
	We refer the reader to \cite[proposition 4.19]{MR1787949}.
\end{proof}

\begin{cor}
		\label{cor.2.6.smash-cofibrant-leftQuillenfunctor}
	Let $A$ be a cofibrant symmetric $T$-spectrum in $\motivicsymmTspectra$.
	Then the adjunction:
		$$\xymatrix{(-\wedge A, \inthomsymmTspectra (A,-),\varphi):\motivicsymmTspectra \ar[r]& \motivicsymmTspectra}
		$$
	is a Quillen adjunction.
\end{cor}
\begin{proof}
	Follows directly from proposition \ref{prop.2.6.symmTspectra-monoidalmodelcategory}.
\end{proof}

\begin{thm}[Jardine]
		\label{thm.2.6.T-spectra===symmTspectra}
	The adjunction:
		$$\xymatrix{(V,U,\varphi):\motivicTspectra \ar[r]& \motivicsymmTspectra}
		$$
	given by the symmetrization and the forgetful functor is a
	Quillen equivalence.
\end{thm}
\begin{proof}
	We refer the reader to \cite[theorem 4.31]{MR1787949}.
\end{proof}

\end{section}
\begin{section}{Cellularity of the Motivic Symmetric Stable Model Structure}
		\label{section-cellularity-motivic-symmetric-stable}
		
	In this section we will show
	that the model category $\motivicsymmTspectra$ is cellular.
	For this we will use 
	the cellularity of 
	$\pointedmotivic$ (see proposition \ref{prop.2.3.cell-prop-pointedmotcat})
	together with some results of Hovey \cite{MR1860878}.
	
\begin{thm}[Hovey]
		\label{thm.2.7.proj-symmstr-cellular}
	Let $\symmTspectra$ be the category of symmetric $T$-spectra
	equipped with the projective model structure
	(see theorem \ref{thm2.6.projstabmodstrsymmTspectra}).
	Then $\symmTspectra$ is a cellular
	model category where the sets of generating
	cofibrations and trivial cofibrations
	are the ones described in theorem \ref{thm2.6.projstabmodstrsymmTspectra}.
\end{thm}
\begin{proof}
	Proposition \ref{prop.2.3.cell-prop-pointedmotcat} implies that
	$\pointedmotivic$ is in particular a cellular, left proper
	and symmetric monoidal model category. 
	We also have that $T=\tee$ is cofibrant in $\pointedmotivic$.
	Therefore we can apply
	theorem A.9 in
	\cite{MR1860878}, which says that the category of symmetric $T$-spectra
	equipped with the projective model structure
	is also cellular under our conditions.
\end{proof}

	Theorem \ref{thm2.6.projstabmodstrsymmTspectra} together 
	with theorem \ref{thm.2.7.proj-symmstr-cellular}
	imply that the projective model structure on
	$\symmTspectra$ is cellular, left proper and simplicial.
	Therefore we can apply Hirschhorn's localization technology to construct
	left Bousfield localizations.
	If we are able to find a suitable \emph{set} of maps
	such that the left Bousfield localization with respect to this
	set recovers the motivic stable model structure on $\symmTspectra$, then
	an immediate corollary of this will be the cellularity of
	the motivic stable model structure for symmetric $T$-spectra.
	
\begin{defi}[Hovey, cf. \cite{MR1860878}]
		\label{def.2.7.localizingset-symmmotstab}
	Let $I_{M_{*}}=\{Y_{+}\hookrightarrow (\Delta ^{n}_{U})_{+}\}$
	be the set of generating cofibrations for $\pointedmotivic$
  (see propositon \ref{prop.2.3.cell-prop-pointedmotcat}).
	Notice that $Y_{+}$ may be equal to $(\Delta ^{n}_{U})_{+}$.
	We consider the following set of maps
	of symmetric $T$-spectra
	$$\xymatrix{S_{\Sigma}=\{ F_{k+1}^{\Sigma}(T\wedge Y_{+}) \ar[r]^-{\zeta ^{Y}_{\Sigma ,k}} & F_{k}^{\Sigma}(Y_{+}) \}}
	$$
	where $\zeta ^{Y}_{\Sigma ,k}$ is the adjoint 
	corresponding to the inclusion determined by the identity in $\Sigma _{k+1}$
	$$\imath _{e}: T\wedge Y_{+}\hookrightarrow Ev_{k+1}(F_{k}^{\Sigma}(Y_{+}))=\Sigma _{k+1}\otimes _{\Sigma _{1}\times \Sigma _{k}}
				(T\wedge \bigvee _{\sigma \in \Sigma _{k}}Y_{+})=\bigvee _{\sigma \in \Sigma _{k+1}} T\wedge Y_{+}$$
	coming from the adjunction between $F_{k+1}^{\Sigma}$ and $Ev_{k+1}$
	(see definition \ref{def.2.6.Fnsigma-Evn-adjunction})
\end{defi}

\begin{prop}[Hovey]
		\label{prop.2.7.Ssigma-local==stably-fibrant}
	Let $X$ be a symmetric $T$-spectrum.  The
	following conditions are equivalent:
	\begin{enumerate}
		\item \label{prop.2.7.Ssigma-local==stably-fibrant.a}  $X$ is stably fibrant, i.e.
					$X$ is a fibrant object in $\motivicsymmTspectra$.
		\item	\label{prop.2.7.Ssigma-local==stably-fibrant.b}	$X$ is $S_{\Sigma}$-local.
	\end{enumerate}
\end{prop}
\begin{proof}
	Follows from definition 8.6 and theorem 8.8 in \cite{MR1860878}, together with
	definition \ref{def.2.6.symmetricstableweakequiv}(\ref{def.2.6.symmetricstableweakequiv.c}) and lemma
	\ref{lem.2.4.fibrant-objects-stablemodstr}.
\end{proof}

	Now it is very easy to show that the
	motivic symmetric stable model structure for symmetric $T$-spectra
	is in fact cellular.
	
\begin{thm}
		\label{thm.2.7.cellularity-motivicsymm-stable-str}
	$\motivicsymmTspectra$ is a cellular model category with the following
	sets $I_{\Sigma}^{T}$, $J_{\Sigma}^{T}$
	of generating cofibrations and trivial cofibrations respectively:
		\begin{eqnarray*}
			I^{T}_{\Sigma} &=& V(I^{T}_{M_{\ast}}) = \bigcup _{k\geq 0}F_{k}^{\Sigma}(I_{M_{\ast}})
										\\ &=& \bigcup _{k\geq 0}\{F_{k}^{\Sigma}(Y_{+}) \hookrightarrow F_{k}^{\Sigma}((\Delta ^{n}_{U})_{+}) 
										\mid U\in (\smoothS), n\geq 0\}  \\
			J_{\Sigma}^{T} &=& \{ j:A\rightarrow B\}
		\end{eqnarray*}
	where $j$ satisfies the following conditions:
		\begin{enumerate}
			\item $j$ is an inclusion of $I_{\Sigma}^{T}$-complexes.
			\item	$j$ is a stable weak equivalence of symmetric $T$-spectra.
			\item	the size of $B$ as an $I_{\Sigma}^{T}$-complex is less than $\kappa$,
						where $\kappa$ is the regular cardinal described by Hirschhorn in 
						\cite[definition 4.5.3]{MR1944041}.
		\end{enumerate}
\end{thm}
\begin{proof}
	By theorem \ref{thm.2.7.proj-symmstr-cellular}
	we know that $\symmTspectra$ is cellular when it is equipped with the projective
	model structure.  Therefore we can apply Hirschhorn's localization techniques
	to construct the left Bousfield localization with respect to the set $S_{\Sigma}$
	of definition \ref{def.2.7.localizingset-symmmotstab}.  
	We claim that this localization coincides with $\motivicsymmTspectra$.
	In effect, using
	proposition \ref{prop.2.7.Ssigma-local==stably-fibrant},
	we have that the fibrant objects in the left Bousfield localization
	with respect to $S_{\Sigma}$ coincide with the fibrant objects in $\motivicsymmTspectra$.
	Therefore a map $f:X\rightarrow Y$ of symmetric $T$-spectra is a weak equivalence
	in the left Bousfield localization with respect to $S_{\Sigma}$ if and only if
	$Qf^{*}:Map(QY,Z)\rightarrow Map(QX,Z)$ is a weak equivalence of simplicial
	sets for every stably fibrant object $Z$ (here $Q$ denotes the cofibrant
	replacement functor in $\symmTspectra$ equipped with the projective model
	structure).  But since $\motivicsymmTspectra$ is a simplicial model category
	and the cofibrations coincide with the projective cofibrations,
	using corollary
	\ref{cor.1.1.4.detect-weak-equiv.2}(\ref{cor.1.1.4.detect-weak-equiv.2.b})
	we get exactly the same characterization for the stable equivalences.
	Hence the weak equivalences in both the motivic symmetric stable structure and the left
	Bousfield localization with respect to $S_{\Sigma}$ coincide.  This implies
	that the motivic symmetric stable model structure and the left Bousfield localization with
	respect to $S_{\Sigma}$ are identical, since the cofibrations in both cases
	are just the cofibrations for the projective model structure on $\Tspectra$.
	
	Therefore using \cite[theorem 4.1.1]{MR1944041}
	we have that the motivic symmetric stable model structure on $\symmTspectra$
	is cellular, since it is constructed applying Hirschhorn technology
	with respect to the set $S_{\Sigma}$.
	
	The claim with respect to the sets of generating cofibrations and trivial
	cofibrations also follows from  \cite[theorem 4.1.1]{MR1944041}.
\end{proof}

	Theorem \ref{thm.2.7.cellularity-motivicsymm-stable-str} will be
	used for the construction
	of new model structures on $\symmTspectra$ which are adequate
	to study the multiplicative properties of
	Voevodsky's slice filtration.

\end{section}
\begin{section}{Modules and Algebras of Motivic Symmetric Spectra}
		\label{section.2.8Modules-Algebras}

	In this section $A$ will always denote a ring spectrum with unit in $\symmTspectra$, and $A$-$\modules$ will denote the
	category of left (or right) $A$-modules.  In case $A$ is a commutative ring spectrum, we will denote the category of $A$-algebras
	by $A$-$\algebras$.  Our goal is to define the model structures induced by the motivic symmetric
	stable model structure on the categories of $A$-modules and $A$-algebras, and to show
	some of their properties.
	
\begin{prop}
		\label{prop.2.8.Free-module-functor}
	We have the following adjunction between the categories of symmetric $T$-spectra and $A$-modules:
		$$\xymatrix@R=.5pt{(A\wedge -, U,\varphi ):\symmTspectra \ar[r] & A\text{-}\modules}
		$$
	where $U(N)=N$ after forgetting the $A$-module structure, and $A\wedge X$ has a structure of $A$-module
	induced by the ring structure on $A$.
\end{prop}
\begin{proof}
	The unit $\eta$ and counit $\delta$ of the adjunction are defined as follows:
		$$\xymatrix@R=.5pt{\eta _{X}:X \cong \symmspherespectrum \wedge X \ar[r]^-{u_{A}\wedge id} & U(A\wedge X)=A\wedge X \\
											 \delta _{N}:A\wedge U(N)=A\wedge N \ar[r]_-{\mu_{N}} & N}
		$$
	where $u_{A}$ is the unit of $A$ and $\mu_{N}$ is the map inducing the $A$-module structure on $N$.
\end{proof}

	The category of $A$-modules inherits a  simplicial
	structure from the one that exists on symmetric $T$-spectra (see section \ref{sec.2.6.symmetricTspectra}).
	
	Given an $A$-module $M$, the tensor objects are defined as follows:
		$$\xymatrix@R=.5pt{M\wedge -:\simpsets \ar[r]& A\text{-}\modules \\
											 K \ar@{|->}[r]& M\wedge K}
		$$
	where $(M\wedge K)^{n}=M^{n}\wedge K_{+}$, i.e. it coincides with the tensor object
	defined for symmetric $T$-spectra and 
	has a structure of $A$-module induced by the one in $M$.
	
	The simplicial functor in two variables is:
		$$\xymatrix@R=.5pt{Map\: _{A\text{-}\modules}(-,-):(A\text{-}\modules)^{op}\times A\text{-}\modules \ar[r] & \simpsets \\
											(M,N) \ar@{|->}[r]& Map\: _{A\text{-}\modules}(M,N)}
		$$
	where $Map\: _{A\text{-}\modules}(M,N)_{n}=\Hom _{A\text{-}\modules}(M\wedge \Delta ^{n}_{+},N)$,
	and finally for any $A$-module $N$ we have the following functor
		$$\xymatrix@R=.5pt{N^{-}:\simpsets \ar[r]& (A\text{-}\modules)^{op}\\
											K \ar@{|->}[r] & N^{K}}
		$$
	where $(N^{K})^{n}=(N^{n})^{K_{+}}$, i.e. it coincides with the cotensor object
	defined for symmetric $T$-spectra and 
	has a structure of $A$-module	
		$A\wedge(N)^{K_{+}}\rightarrow  N^{K_{+}}$
	adjoint to
		$$\xymatrix{A\wedge (N)^{K_{+}}\wedge K_{+} \ar[rr]^-{id\wedge ev_{K_{+}}} && A\wedge N\ar[r]^-{\mu} & N}
		$$
	where $\mu$ is the map that induces the $A$-module structure on $N$.

	Similarly, it is possible to promote the action
	of $\pointedmotivic$ on the category of symmetric $T$-spectra to the
	category of $A$-modules, i.e.
	the category of $A$-modules $A$-$\modules$ has the structure of
	a closed $\pointedmotivic$-module, which
	is obtained by extending the symmetric monoidal structure for $\pointedmotivic$
	levelwise.
	
	The bifunctor giving the adjunction of two variables is defined
	as follows:
		$$\xymatrix@R=.5pt{- \wedge -:A\text{-}\modules \times \pointedmotivic \ar[r]& A\text{-}\modules \\
											(M,D) \ar@{|->}[r]& M\wedge D}
		$$
	with $(M\wedge D)^{n}=M^{n}\wedge D$, i.e. it coincides with the tensor object
	defined for symmetric $T$-spectra and 
	has a structure of $A$-module induced by the one in $M$.
	
	The adjoints are given by:
		$$\xymatrix@R=.5pt{\Omega _{-}-:\pointedmotivic ^{op}\times A\text{-}\modules \ar[r]& A\text{-}\modules \\
											(D,N) \ar@{|->}[r]& \Omega _{D}N}
		$$
	
		$$\xymatrix@R=.5pt{\inthomAmodpresheaf (-,-):(A\text{-}\modules )^{op}\times A\text{-}\modules \ar[r]& \pointedmotivic \\
											(M,N) \ar@{|->}[r]& \inthomAmodpresheaf (M,N)}
		$$
	where $(\Omega _{D}N)^{n}=\Omega _{D}N^{n}$, i.e. it coincides with the cotensor object
	defined for symmetric $T$-spectra and 
	has a structure of $A$-module
	$A\wedge(\Omega _{D}N)\rightarrow \Omega _{D}N$  adjoint to
		$$\xymatrix{A\wedge (\Omega _{D}N)\wedge D \ \ar[rr]^-{id\wedge ev_{D}}
								&& A\wedge N \ar[r]^-{\mu}& N}$$
	and $\inthomAmodpresheaf (M,N)$ is the following pointed simplicial presheaf
	on $\smoothS$:
		$$\xymatrix@R=.5pt{\inthomAmodpresheaf (M,N):(\smoothS \times \Delta)^{op} \ar[r] &
											Sets \\
											(U,n) \ar@{|->}[r] & \Hom _{A\text{-}\modules}(M \wedge (\Delta ^{n}_{U})_{+},N)}
		$$
		
\begin{prop}
		\label{prop.2.8.enriched-freemodule-adjunction}
	The adjunction (see proposition \ref{prop.2.8.Free-module-functor})
		$$\xymatrix{(A\wedge -,U,\varphi):\symmTspectra \ar[r]& A\text{-}\modules}
		$$
	is enriched in the categories of simplicial sets and pointed
	simplicial presheaves on $\nissite$, i.e. for every symmetric $T$-spectrum $X$ and
	for every $A$-module $N$ we have the following natural isomorphisms:
		$$\xymatrix{Map\: _{A\text{-}\modules}(A\wedge X,N)\ar[r]_-{\cong}^-{\epsilon}& Map_{\Sigma}(X,UN)}
		$$
		
		$$\xymatrix{\inthomAmodpresheaf (A\wedge X,N) \ar[r]_-{\cong}^-{\eta}& \inthomsymmTspectrapresheaf (X,UN)}
		$$
\end{prop}
\begin{proof}
	We consider first the simplicial isomorphism:
	Given any $n$-simplex $t$ in $Map\: _{A\text{-}\modules} (A\wedge X,N)$
		$$\xymatrix{A\wedge X\wedge \Delta ^{n}_{+}\ar[r]^-{t}& N}
		$$
	use the adjunction between $A\wedge-$ and
	$U$ to get the associated $n$-simplex $\epsilon (t)$
	in $Map_{\Sigma}(X,UY)$:
		$$\xymatrix{X\wedge \Delta ^{n}_{+}\ar[r]^-{\epsilon (t)}& UN}
		$$
		
	We consider now the isomorphism of simplicial presheaves:
	Given any simplex $s$ in $\inthomAmodpresheaf (A\wedge X,N)$
		$$\xymatrix{A\wedge X\wedge (\Delta ^{n}_{W})_{+}\ar[r]^-{s}& N}
		$$
	use the adjunction between $A\wedge -$ and $U$
	to get the associated simplex $\eta (s)$ in $\inthomsymmTspectrapresheaf (X,UN)$:
		$$\xymatrix{X\wedge (\Delta ^{n}_{W})_{+}\ar[r]^-{\eta (s)}& UN}
		$$
\end{proof}
		
	If $A$ is a commutative ring spectrum then
	$A$-$\modules$ is a closed symmetric monoidal category, where the monoidal structure is induced
	by the one exisiting on $\symmTspectra$.  Namely,
		$$\xymatrix@R=.5pt{-\wedge _{A}-:A\text{-}\modules \times A\text{-}\modules \ar[r]& A\text{-}\modules \\
											(M,N) \ar@{|->}[r]& M\wedge _{A}N}
		$$
	
		$$\xymatrix@R=.5pt{\inthomAmod (-,-):(A\text{-}\modules )^{op}\times A\text{-}\modules \ar[r]& A\text{-}\modules \\
											(M,N) \ar@{|->}[r]& \inthomAmod (M,N)}
		$$
	where $M\wedge _{A}N$ is defined as the colimit of the following diagram
		$$\xymatrix{A\wedge M\wedge N \ar@<1ex>[rr]^-{\mu _{M}\wedge id} \ar@<-1ex>[rr] && M\wedge N}
		$$
	with the bottom arrow  given by the following composition
		$$\xymatrix{A\wedge M\wedge N \ar[rr]^-{t\wedge id} && M\wedge A\wedge N \ar[rr]^-{id\wedge \mu _{N}} 
								&& M\wedge N}
		$$
	and $\inthomAmod (M,N)$ is defined as the limit of the following diagram
		$$\xymatrix{\inthomsymmTspectra (M,N) \ar@<1ex>[rr]^-{(\mu _{M})^{\ast}} \ar@<-1ex>[rr]_-{(\mu _{N})_{\ast}}
								&& \inthomsymmTspectra (A\wedge M,N)}
		$$

	If $A$ is not commutative, the bifunctor $-\wedge _{A}-$ defines instead
	an adjunction of two variables from the categories of right and left $A$-modules
	to the category of symmetric $T$-spectra:
		$$\xymatrix@R=.5pt{-\wedge _{A}-:A\text{-}\modules _{r}\times A\text{-}\modules _{l}\ar[r]& \symmTspectra \\
											(M,N) \ar@{|->}[r]& M\wedge _{A}N}
		$$
	given a right $A$-module $M$, the right adjoint to 
		$$M\wedge _{A}-:A\text{-}\modules _{l} \rightarrow \symmTspectra$$ 
	is given by
		$$\inthomsymmTspectra (M,-):\symmTspectra \rightarrow A\text{-}\modules _{l}$$
	where $\inthomsymmTspectra (M,Z)$ has a structure of left $A$-module
		$$\mu:A\wedge \inthomsymmTspectra (M,Z)\rightarrow \inthomsymmTspectra (M,Z)$$
	defined as the adjoint of the following composition
		$$\xymatrix{M\wedge A\wedge \inthomsymmTspectra (M,Z) \ar[rr]^-{\mu _{M}\wedge id}&& M\wedge \inthomsymmTspectra (M,Z)
								\ar[r]^-{\epsilon _{Z}}& Z}$$
	where $\mu _{M}$ denotes the map defining the right $A$-module structure for $M$ and $\epsilon$ denotes the counit of the adjunction
	between $M \wedge Z$ and $\inthomsymmTspectra (M,-)$.   The construction of the remaining adjoint is similar.
	
\begin{thm}
		\label{thm.2.8.motivicsymmetricA-modules}
	Let $A$ be a cofibrant ring object in $\motivicsymmTspectra$.  Then
	the adjuntion (see proposition \ref{prop.2.8.Free-module-functor}):
		$$\xymatrix{(A\wedge -, U,\varphi ):\motivicsymmTspectra \ar[r]
									& A\text{-}\modules}
		$$
	induces a model structure for the category
	$A$-$\modules$ of $A$-modules, i.e.
	a map $f$ in $A$-$\modules$ is a fibration or a weak equivalence
	if and only if $U(f)$ is a fibration or a weak equivalence in $\motivicsymmTspectra$
	(see theorem \ref{thm2.6.stablesymmetricmodelstructure}).
		
	This model structure will be called \emph{motivic stable}, and the category
	of $A$-modules equipped with the motivic stable model structure will be denoted
	by $\motivicAmod$.
\end{thm}
\begin{proof}
	We have that $\motivicsymmTspectra$ is a cellular model category (see
	theorem \ref{thm.2.7.cellularity-motivicsymm-stable-str}), 
	i.e. in particular a cofibrantly generated model category,
	and a monoidal model category in the sense of Hovey
	(see proposition \ref{prop.2.6.symmTspectra-monoidalmodelcategory}).
	Therefore, since $A$ is cofibrant the result follows from
	\cite[corollary 2.2]{HoveyPreprint}.
\end{proof}

\begin{lem}
		\label{lemma.2.5.changecoefficients-Quillenfunctor}
	Let $f:A\rightarrow A'$ be a map between cofibrant ring spectra in $\motivicsymmTspectra$,
	which is compatible with the ring structures.
	Then the adjunction:
		$$(A'\wedge _{A}-, U,\varphi ):\motivicAmod \rightarrow \motivicAAmod
		$$
	is a Quillen adjunction.  Furthermore, a map $w:M\rightarrow M'$ in $\motivicAAmod$ is a weak equivalence if
	and only if $Uw$ is a weak equivalence in $\motivicAmod$.
\end{lem}
\begin{proof}
	It is clear that $U:\motivicAAmod \rightarrow \motivicAmod$ is a right Quillen functor,
	since the fibrations (respectively, trivial fibrations) 
	for both model structures are detected in $\motivicsymmTspectra$.
	Finally, the claim related to the weak equivalences follows immediately 
	from theorem \ref{thm.2.8.motivicsymmetricA-modules}.
\end{proof}

\begin{prop}
		\label{prop.2.8.invarianceofcoefficients}
	Let $f:A\rightarrow A'$ be a weak equivalence between cofibrant ring spectra in $\motivicsymmTspectra$,
	which is compatible with the ring structures.
	Then $f$ induces a Quillen equivalence between the motivic stable model structures of $A$ and $A'$ modules:
		$$\xymatrix{(A'\wedge _{A} -, U,\varphi ):\motivicAmod \ar[r] & \motivicAAmod}
		$$
\end{prop}
\begin{proof}
	It follows immediately from theorem 2.4 in \cite{HoveyPreprint} together with the fact that
	the domains of the generating cofibrations for $\motivicsymmTspectra$ are cofibrant
	(see theorem \ref{thm.2.7.cellularity-motivicsymm-stable-str}).
\end{proof}

\begin{prop}
		\label{prop.2.8.enriched-freemodule-adjunction2}
	Let $f:A\rightarrow A'$ be a map between cofibrant ring spectra in $\motivicsymmTspectra$,
	which is compatible with the ring structures.
	Then the adjunction 
		$$\xymatrix{(A'\wedge _{A}-,U,\varphi):\motivicAmod \ar[r]& \motivicAAmod}
		$$
	is enriched in the categories of simplicial sets and pointed
	simplicial presheaves on $\nissite$, i.e. for every $A$-module $M$ and
	for every $A'$-module $N$ we have the following natural isomorphisms:
		$$\xymatrix{Map\: _{A'\text{-}\modules}(A'\wedge _{A}M,N)\ar[r]_-{\cong}^-{\epsilon}& 
								Map_{A\text{-}\modules}(M,UN)}
		$$
		
		$$\xymatrix{\inthomAAmodpresheaf (A'\wedge _{A}M,N) \ar[r]_-{\cong}^-{\eta}& \inthomAmodpresheaf (M,UN)}
		$$
\end{prop}
\begin{proof}
	The proof is exactly the same as the one in propositon \ref{prop.2.8.enriched-freemodule-adjunction}.
\end{proof}

\begin{prop}
		\label{prop.2.8.motivicAmod-cofibrations=>motivic-cofibrations}
	Let $A$ be a cofibrant ring spectrum in $\motivicsymmTspectra$, and let
	$i$ be a cofibration in $\motivicAmod$.  Then $U(i)$
	is also a cofibration
	in $\motivicsymmTspectra$.
\end{prop}
\begin{proof}
	Theorem \ref{thm.2.7.cellularity-motivicsymm-stable-str} implies in particular that
	$\motivicsymmTspectra$ is a cofibrantly generated model category.  Therefore the proposition
	follows directly from \cite[corollary 2.2]{HoveyPreprint}.
\end{proof}

\begin{prop}
			\label{prop.2.8.motivicA-modules=>M-symmSpt-modelcat}
	$\motivicAmod$ is a: 
		\begin{enumerate}
			\item \label{prop.2.8.motivicA-modules=>M-symmSpt-modelcat.a}	proper model category.
			\item \label{prop.2.8.motivicA-modules=>M-symmSpt-modelcat.b}  simplicial model category.
			\item	\label{prop.2.8.motivicA-modules=>M-symmSpt-modelcat.c}  $\pointedmotivic$-model category
																																		(see definition \ref{def.module-modcats}).
		\end{enumerate}	
\end{prop}
\begin{proof}
	(\ref{prop.2.8.motivicA-modules=>M-symmSpt-modelcat.a}):  It follows directly from the fact that
	$\motivicsymmTspectra$ is a proper model category 
	(see theorem \ref{thm2.6.stablesymmetricmodelstructure}), together with 
	theorem \ref{thm.2.8.motivicsymmetricA-modules} and
	proposition \ref{prop.2.8.motivicAmod-cofibrations=>motivic-cofibrations}.
	
	(\ref{prop.2.8.motivicA-modules=>M-symmSpt-modelcat.b}):  Since the cotensor objects $N^{K}$ for the simplicial structure
	are identical in $\motivicAmod$ and $\motivicsymmTspectra$, the results follows from theorem
	\ref{thm.2.8.motivicsymmetricA-modules} and theorem
	\ref{thm2.6.stablesymmetricmodelstructure} which implies in particular that $\motivicsymmTspectra$
	is a simplicial model category.
	
	(\ref{prop.2.8.motivicA-modules=>M-symmSpt-modelcat.c}):  Since the cotensor objects $\Omega _{D}N$
	for the $\pointedmotivic$-action are identical in $\motivicAmod$ and $\motivicsymmTspectra$,
	the results follows from the fact that $\motivicsymmTspectra$ is a $\pointedmotivic$-model category
	(see proposition \ref{prop.2.6.motsymmTspectra-pointedpresheaf-modelcat}) together with
	theorem \ref{thm.2.8.motivicsymmetricA-modules}.
\end{proof}

\begin{thm}
		\label{thm.2.8.motivicAmods=>cellular}
	$\motivicAmod$ is a cellular
	model category with the following sets $I_{A\text{-}\modules}$, $J_{A\text{-}\modules}$
	of generating cofibrations and trivial cofibrations respectively 
	(see theorem \ref{thm.2.7.cellularity-motivicsymm-stable-str}):  
		\begin{eqnarray*}
			I_{A\text{-}\modules} &=& A\wedge I^{T}_{\Sigma}
														\\ &=& \bigcup _{k\geq 0}
														\{ id\wedge i:A\wedge F_{k}^{\Sigma}(Y_{+})\rightarrow
														A\wedge F_{k}^{\Sigma}((\Delta _{n}^{U})_{+}) \mid U\in (\smoothS), n\geq 0\} \\		
			J_{A\text{-}\modules} &=& A\wedge J^{T}_{\Sigma}=\{ id\wedge j:A\wedge X\rightarrow A\wedge Y\}
		\end{eqnarray*}
	where $j:X\rightarrow Y$ satisfies the following conditions:
		\begin{enumerate}
			\item $j$ is an inclusion of $I_{\Sigma}^{T}$-complexes.
			\item	$j$ is a stable weak equivalence of symmetric $T$-spectra.
			\item	the size of $Y$ as an $I_{\Sigma}^{T}$-complex is less than $\kappa$,
						where $\kappa$ is the regular cardinal described by Hirschhorn in 
						\cite[definition 4.5.3]{MR1944041}.
		\end{enumerate}
\end{thm}
\begin{proof}
	We have to check that the conditions 
	(\ref{def.1.1.2.cell-mod-cats.a})-(\ref{def.1.1.2.cell-mod-cats.d}) of
	definition \ref{def.1.1.2.cell-mod-cats} hold.
	
	By construction (see theorem \ref{thm.2.8.motivicsymmetricA-modules})
	it is clear that $I_{A\text{-}\modules}$ and $J_{A\text{-}\modules}$ are generators for the model
	structure on $\motivicAmod$.  This takes care of (\ref{def.1.1.2.cell-mod-cats.a}).
	
	By adjointness, to prove (\ref{def.1.1.2.cell-mod-cats.b}) it suffices to show that
	the domains and codomains of $I^{T}_{\Sigma}$ are compact relative to $I_{A\text{-}\modules}$.
	However, the domains and codomains of  $I^{T}_{\Sigma}$ are cofibrant in $\motivicsymmTspectra$,
	which is in particular a cellular model category (see theorem \ref{thm.2.7.cellularity-motivicsymm-stable-str}).
	Hence \cite[corollary 12.3.4]{MR1944041} implies that the domains and codomains of $I^{T}_{\Sigma}$
	are compact with respect to the class of cofibrations in $\motivicsymmTspectra$.  Finally,
	proposition \ref{prop.2.8.motivicAmod-cofibrations=>motivic-cofibrations} implies that all the maps
	in $I_{A\text{-}\modules}$ are cofibrations in $\motivicsymmTspectra$.  Thus, the domains and codomains
	of $I^{T}_{\Sigma}$ are compact with respect to $I_{A\text{-}\modules}$, as we wanted.
	
	Again by adjointness, to prove (\ref{def.1.1.2.cell-mod-cats.c}) it suffices to show that
	the domains of $J^{T}_{\Sigma}$ are small relative to $I_{A\text{-}\modules}$.  But
	proposition \ref{prop.2.8.motivicAmod-cofibrations=>motivic-cofibrations} implies that all the maps
	in $I_{A\text{-}\modules}$ are cofibrations in $\motivicsymmTspectra$.  Therefore by
	\cite[theorem 12.4.4]{MR1944041} we have that the domains of $J^{T}_{\Sigma}$ are small relative to 
	$I_{A\text{-}\modules}$, since $\motivicsymmTspectra$ is a cellular model category
	(see theorem \ref{thm.2.7.cellularity-motivicsymm-stable-str}).
	
	Finally, proposition \ref{prop.2.8.motivicAmod-cofibrations=>motivic-cofibrations} 
	implies that the cofibrations in $\motivicAmod$ are in particular
	cofibrations in $\motivicsymmTspectra$, which is a cellular model category
	(see theorem \ref{thm.2.7.cellularity-motivicsymm-stable-str}). 
	Therefore the cofibrations in $\motivicAmod$ are effective monomorphisms
	in $\motivicsymmTspectra$.  This takes care of (\ref{def.1.1.2.cell-mod-cats.d})
	since the limits and colimits in $\motivicAmod$ are computed in $\motivicsymmTspectra$.
\end{proof}

\begin{thm}
		\label{thm.2.8.smashT-Quillenequiv}
	Let $T=S^{1}\wedge \gm \in \pointedmotivic$.
	Then the adjunction:
		$$\xymatrix{(-\wedge T, \Omega _{T}, \varphi ):\motivicAmod \ar[r]& \motivicAmod}
		$$
	is a Quillen equivalence.
\end{thm}
\begin{proof}
	Every pointed simplicial presheaf in $\pointedmotivic$ is cofibrant, therefore proposition
	\ref{prop.2.8.motivicA-modules=>M-symmSpt-modelcat}(\ref{prop.2.8.motivicA-modules=>M-symmSpt-modelcat.c})
	implies that $-\wedge T:\motivicAmod \rightarrow \motivicAmod$ is a left Quillen functor.
	
	Let $\eta$, $\epsilon$ denote the unit and counit of the
	adjunction $(-\wedge T, \Omega _{T}, \varphi)$.
	By proposition 1.3.13 in \cite{MR1650134}, it suffices to check
	that the following conditions hold:
		\begin{enumerate}
			\item \label{thm.2.8.smashT-Quillenequiv.a}For every cofibrant
						$A$-module $M$ in $\motivicAmod$, the following composition
							$$\xymatrix{M\ar[r]^-{\eta _{M}}& \Omega _{T}(T\wedge M)
								\ar[rr]^-{\Omega _{T}(R^{T\wedge M})}&& \Omega _{T}R(T\wedge M)}
							$$
						is a weak equivalence in $\motivicAmod$, where $R$ denotes a fibrant replacement functor
						in $\motivicAmod$.
			\item \label{thm.2.8.smashT-Quillenequiv.b}For every fibrant
						$A$-module $M$ in $\motivicAmod$, the following composition
							$$\xymatrix{T\wedge Q(\Omega _{T}M)\ar[rr]^-{id\wedge Q^{\Omega _{T}M}}&& 
													T\wedge (\Omega _{T}M)\ar[r]^-{\epsilon _{M}}& M}
							$$
						is a weak equivalence in $\motivicAmod$, where $Q$ denotes a cofibrant replacement functor
						in $\motivicAmod$.
		\end{enumerate}
		
	(\ref{thm.2.8.smashT-Quillenequiv.a}):  By proposition \ref{prop.2.8.motivicAmod-cofibrations=>motivic-cofibrations}
	we have that $M$ is cofibrant in $\motivicsymmTspectra$.
	Thus the result follows immediately from theorems
	\ref{thm.2.6.Tsusp==QuillenequivonSymmTspectra}(\ref{thm.2.6.Tsusp==QuillenequivonSymmTspectra.a})
	and \ref{thm.2.8.motivicsymmetricA-modules}.
	
	(\ref{thm.2.8.smashT-Quillenequiv.b}):  Follows directly from theorem \ref{thm.2.8.motivicsymmetricA-modules},
	proposition \ref{prop.2.8.motivicAmod-cofibrations=>motivic-cofibrations} and theorem 
	\ref{thm.2.6.Tsusp==QuillenequivonSymmTspectra}(\ref{thm.2.6.Tsusp==QuillenequivonSymmTspectra.b}).
\end{proof}

\begin{prop}
		\label{prop.2.8.Acomm=>motivicAmodules-symmmonoidal}
	Let $A$ be a cofibrant commutative ring spectrum in $\motivicsymmTspectra$.  Then
	$\motivicAmod$ is a symmetric monoidal model category in the sense of Hovey
	(see definition \ref{def.mon-mod-cats}).
\end{prop}
\begin{proof}
	This follows directly from theorem \ref{thm.2.7.cellularity-motivicsymm-stable-str},
	proposition \ref{prop.2.6.symmTspectra-monoidalmodelcategory} and 
	\cite[proposition 2.8(2)]{HoveyPreprint}.
\end{proof}

	If $A$ is not commutative then we get a weaker version of the previous proposition.
	
\begin{prop}
		\label{prop.2.8.A=>motivicAmodules-leftrightaction}
	Let $A$ be a cofibrant ring spectrum in $\motivicsymmTspectra$.  Then
	$-\wedge _{A}-$ defines a Quillen adjunction of two variables
	(see definition \ref{def.Quillen-bifunct}) from the motivic model structure for right and left $A$-modules
	to the motivic symmetric stable model structure:
		$$-\wedge _{A}-:\motivicAmod _{r}\times \motivicAmod _{l} \rightarrow \motivicsymmTspectra$$
\end{prop}
\begin{proof}
	We need to show that given a cofibration $i:M\rightarrow M'$ in $\motivicAmod _{r}$ and
	a cofibration $j:N\rightarrow N'$ in $\motivicAmod _{l}$, the induced map
		$$i \Box _{A} j:M\wedge _{A} N' \coprod _{M\wedge _{A}N} M'\wedge _{A}N \rightarrow M'\wedge _{A}N'
		$$
	is a cofibration in $\motivicsymmTspectra$, which is trivial if either $i$ or $j$ are trivial.
	
	Clearly it is enough to do it for the generating cofibrations and trivial cofibrations
	in $\motivicAmod$ (see lemma 3.5 in \cite{MR1734325}).  Theorem \ref{thm.2.8.motivicAmods=>cellular}
	implies that $A\wedge I^{T}_{\Sigma}$, $A\wedge J^{T}_{\Sigma}$ ($I^{T}_{\Sigma}\wedge A$, $J^{T}_{\Sigma}\wedge A$)
	are the sets of generating cofibrations
	and trivial cofibrations for $\motivicAmod _{l}$ (respectively for $\motivicAmod _{r}$), where
	$I^{T}_{\Sigma}$, $J^{T}_{\Sigma}$ denote the sets of generating cofibrations and trivial cofibrations for
	$\motivicsymmTspectra$.
	
	Now
		$$
			I^{T}_{\Sigma}\wedge {A} \Box _{A}A\wedge I^{T}_{\Sigma} =I^{T}_{\Sigma} \Box A\wedge I^{T}_{\Sigma} \\
				\subseteq I^{T}_{\Sigma} \Box I^{T}_{\Sigma} \subseteq I^{T}_{\Sigma}
		$$
	where the equality follows by definition, the first inclusion follows from the fact that $A$ is cofibrant
	in $\motivicsymmTspectra$ and $\motivicsymmTspectra$ is a symmetric monoidal model category
	(see proposition \ref{prop.2.6.symmTspectra-monoidalmodelcategory}), 
	and the last inclusion follows from the fact that $\motivicsymmTspectra$ is a symmetric
	monoidal model category.  A similar argument shows that
		$$
			J^{T}_{\Sigma}\wedge {A} \Box _{A}A\wedge I^{T}_{\Sigma} =J^{T}_{\Sigma} \Box A\wedge I^{T}_{\Sigma} \\
				\subseteq J^{T}_{\Sigma} \Box I^{T}_{\Sigma} \subseteq J^{T}_{\Sigma}
		$$
	and finally, the remaining case follows from this by symmetry.
\end{proof}
		
	In the rest of this section, we assume that $A$ is a commutative ring spectrum 
	with unit in $\symmTspectra$.  The category of  
	$A$-algebras is a symmetric monoidal category, where the monoidal structure coincides with
	the one exisiting on $A$-$\modules$.  Namely,
		$$\xymatrix@R=.5pt{-\wedge _{A}-:A\text{-}\algebras \times A\text{-}\algebras \ar[r]& A\text{-}\algebras \\
											(C,D) \ar@{|->}[r]& C\wedge _{A}D}
		$$
	However, the category of $A$-algebras is not a closed symmetric monoidal category, i.e. the functor
	$C\wedge _{A}-:A\text{-}\algebras \rightarrow A\text{-}\algebras$ does not have in general a right
	adjoint.
		
\begin{prop}
		\label{prop.2.8.Free-algebra-functor}
	We have the following adjunction between the categories of symmetric $T$-spectra and $A$-algebras:
		$$\xymatrix@R=.5pt{(T, U,\varphi ):\symmTspectra \ar[r] & A\text{-}\algebras}
		$$
	where $U(N)=N$ after forgetting the $A$-algebra structure, and 
	$T(X)=A\wedge \coprod _{n\geq 0}X^{\wedge n}$ has a structure of $A$-algebra
	induced by concatenation together with the ring structure on $A$.
\end{prop}
\begin{proof}
	The unit $\eta$ of the adjunction is  
		$$\xymatrix{X \cong \symmspherespectrum \wedge X 
												\ar[r]^-{u_{A}\wedge id} & A\wedge X \ar@{^{(}->}[r]^-{i_{X\wedge A}}&
												U(T(X))=A\wedge \coprod _{n\geq 0}X^{\wedge n}}
		$$
	where $u_{A}$ is the unit of $A$.  On the other hand,
	the counit $\delta$ of the adjunction is induced by
	iterating the map
	that induces the $A$-algebra structure of $B$ 
		$$\xymatrix{A\wedge \coprod B^{\wedge k} \ar@{^{(}->}[r]& A\wedge \coprod _{n\geq 0}B^{\wedge n}
		            \ar[r]^-{\delta_{B}} & B}$$
\end{proof}

\begin{lem}
		\label{lemma.2.8.smallness-Aalgs}
	Let $I_{\Sigma}^{T}$, $J_{\Sigma}^{T}$ be the sets of generating cofibrations and trivial cofibrations
	for the motivic symmetric stable model structure $\motivicsymmTspectra$ of
	symmetric $T$-spectra (see theorem \ref{thm.2.7.cellularity-motivicsymm-stable-str}).
	Then:
		\begin{enumerate}
			\item \label{lemma.2.8.smallness-Aalgs.a}  The domains of $I_{\Sigma}^{T}$ are small relative to
																									$X\wedge I_{\Sigma}^{T}$-cell for every symmetric $T$-spectrum $X$.
			\item \label{lemma.2.8.smallness-Aalgs.b}  The domains of $J_{\Sigma}^{T}$ are small relative to
																									$X\wedge J_{\Sigma}^{T}$-cell for every symmetric $T$-spectrum $X$.
			\item \label{lemma.2.8.smallness-Aalgs.c}	 The maps of $X\wedge J_{\Sigma}^{T}$-cell are weak equivalences
																									for every symmetric $T$-spectrum $X$.
		\end{enumerate}
\end{lem}
\begin{proof}
	Let $I$, $J$ denote the sets of generating cofibrations and trivial cofibrations
	for the category of symmetric $T$-spectra $\symmTspectra$ equipped with the injective
	model structure (see theorem \ref{thm.2.6.injmodstrsymmTspectra}), 
	where the cofibrations and the weak equivalences are defined
	levelwise.  Hence every symmetric $T$-spectrum is cofibrant in the injective model
	structure.  On the other hand, theorem
	\ref{thm.2.6.injmodstrsymmTspectra} implies that the injective model structure is cofibrantly generated and
	that the codomains of the generating cofibrations $I$ are small relative to $I$. 
	Thus, applying \cite[corollary 11.2.4]{MR1944041} 
	we get that every symmetric $T$-spectrum is small with respect to the
	class of level cofibrations.
	
	(\ref{lemma.2.8.smallness-Aalgs.a}):  It suffices to show that every map in $X\wedge I_{\Sigma}^{T}$
	is a level cofibration.  But this follows directly from \cite[proposition 4.19]{MR1787949}, since
	every symmetric $T$-spectrum $X$ is cofibrant in the injective
	model structure.
	
	(\ref{lemma.2.8.smallness-Aalgs.b}):  It suffices to show that every map in $X\wedge J_{\Sigma}^{T}$
	is a level cofibration.  But this is a consequence of  \cite[proposition 4.19]{MR1787949}, since
	every symmetric $T$-spectrum $X$ is cofibrant in the injective
	model structure.
	
	(\ref{lemma.2.8.smallness-Aalgs.c}):  This follows immediately from \cite[proposition 4.19]{MR1787949},
	since every symmetric $T$-spectrum is cofibrant in the injective model structure.
\end{proof}

\begin{thm}
		\label{thm.2.8.motivic-Aalgebras}
	Let $A$ be a cofibrant commutative ring object 
	with unit in $\motivicsymmTspectra$.  Then
	the adjuntion (see proposition \ref{prop.2.8.Free-algebra-functor}):
		$$\xymatrix{(T, U, \varphi ):\motivicsymmTspectra \ar[r]
									& A\text{-}\algebras}
		$$
	induces a model structure for the category
	$A$-$\algebras$ of $A$-algebras, i.e.
	a map $f$ in $A$-$\algebras$ is a fibration or a weak equivalence
	if and only if $U(f)$ is a fibration or a weak equivalence in $\motivicsymmTspectra$
	(see theorem \ref{thm2.6.stablesymmetricmodelstructure}).
		
	This model structure will be called \emph{motivic}, and the category
	of $A$-algebras equipped with the motivic model structure will be denoted
	by $\motivicAalg$.  Furthermore, $\motivicAalg$ is a cofibrantly generated model
	category with the following sets $I_{A\text{-}\algebras}$, $J_{A\text{-}\algebras}$
	of generating cofibrations and trivial cofibrations respectively 
	(see theorem \ref{thm.2.7.cellularity-motivicsymm-stable-str}):  
		\begin{eqnarray*}
			I_{A\text{-}\algebras} &=& T(I^{T}_{\Sigma})
														\\ &=& \bigcup _{k\geq 0}
														\{ T(i):T(F_{k}^{\Sigma}(Y_{+}))\rightarrow
														 T(F_{k}^{\Sigma}((\Delta _{n}^{U})_{+})) \mid U\in (\smoothS), n\geq 0\} \\		
			J_{A\text{-}\algebras} &=& T(J^{T}_{\Sigma})=\{ T(j):T(X)\rightarrow T(Y) \}
		\end{eqnarray*}
	where $j:X\rightarrow Y$ satisfies the following conditions:
		\begin{enumerate}
			\item $j$ is an inclusion of $I_{\Sigma}^{T}$-complexes.
			\item	$j$ is a stable weak equivalence of symmetric $T$-spectra.
			\item	the size of $Y$ as an $I_{\Sigma}^{T}$-complex is less than $\kappa$,
						where $\kappa$ is the regular cardinal described by Hirschhorn in 
						\cite[definition 4.5.3]{MR1944041}.
		\end{enumerate}
\end{thm}		
\begin{proof}
	Theorem \ref{thm.2.7.cellularity-motivicsymm-stable-str} 
	implies that $\motivicsymmTspectra$ is in particular a cofibrantly
	generated model category, and by proposition \ref{prop.2.6.symmTspectra-monoidalmodelcategory} 
	we have that $\motivicsymmTspectra$
	is a symmetric monoidal model category.  Therefore the result
	follows immediately from lemma \ref{lemma.2.8.smallness-Aalgs}
	and \cite[theorem 3.1]{HoveyPreprint}.
\end{proof}

\begin{prop}
		\label{prop.2.8.Aalg-cofibrations=>motivic-cofibrations}
	Let $A$ be a cofibrant commutative ring object 
	with unit in $\motivicsymmTspectra$, and let $f:B\rightarrow B'$
	be a map of $A$-algebras which is a cofibration in the motivic
	model category $\motivicAalg$ of $A$-algebras.  If
	$B$ is cofibrant in $\motivicAmod$,
	then
	$Uf$ is a cofibration in $\motivicAmod$.
\end{prop}
\begin{proof}
	It follows directly from lemma 6.2 in
	\cite{MR1734325}. 
\end{proof}
		
\end{section}

\end{chapter}
\begin{chapter}{Model Structures for the Slice Filtration}
		\label{chap-slice-filtration}

	This chapter contains our main results.
	In section \ref{section-slice-filtration}, we recall Voevodsky's construction of the slice
	filtration in the context of simplicial presheaves.  In section \ref{section.3.2.modstrslicefilt},
	we apply Hirschhorn's localization techniques to the Morel-Voevodsky
	stable model structure $\motivicTspectra$, in order to construct
	three new families of model structures, namely $\qconnectedTspectra$,
	$\weightqTspectra$ and $\qsliceTspectra$.  These  model structures will
	provide a lifting of Voevodsky's slice filtration to the model category setting.
	Furthermore, we will also get a simple description for the exact functors
	$f_{q}$ (($q-1$)-connective cover) and $s_{q}$ ($q$-slice) defined in section \ref{section-slice-filtration}, 
	in terms of a suitable composition of cofibrant and fibrant replacement functors.
	
	In section \ref{section.3.3.symmmodstrslicefilt}, we promote the model structures introduced in section 
	\ref{section.3.2.modstrslicefilt} to
	the setting of symmetric $T$-spectra.  These new model structures will be denoted by
	$\qconnectedsymmTspectra$, $\weightqsymmTspectra$ and $\qslicesymmTspectra$.  We will prove that
	the Quillen adjunction given by the symmetrization and the forgetful functors descends to a
	Quillen equivalence for these three new model structures.  As a consquence we will see
	that the model categories $\qconnectedsymmTspectra$, $\weightqsymmTspectra$ and $\qslicesymmTspectra$
	provide a lifting for Voevodsky's slice filtration and give an alternative description
	for the functors $f_{q}$ and $s_{q}$.  The great technical advantage of these model structures relies on
	the fact that the underlying category is symmetric monoidal.  Hence, we have a natural framework
	to describe the multiplicative properties of the slice filtration.
	
	In section \ref{section.3.4.multiplicativeproperties-slicefiltration},
	we will show that the slice filtration is compatible with the smash product of
	symmetric $T$-spectra.
	
	In section \ref{section.3.5furthermultiplicativeproperties}, we will promote the model structures
	constructed in section \ref{section.3.3.symmmodstrslicefilt} to the category of $A$-modules,
	where $A$ is a cofibrant ring spectrum with unit in $\motivicsymmTspectra$.  We will denote these 
	new model structures by $\qconnectedAmod$, $\weightqAmod$ and $\qsliceAmod$.  These
	new model structures will give an analogue of the slice filtration for
	the motivic stable homotopy category of $A$-modules.
	We will see that when one imposes some natural additonal conditions on the ring spectrum $A$,
	the free $A$-module functor ($A\wedge -$) induces a strict compatibility
	between the slice filtration in the categories of symmetric $T$-spectra and $A$-modules.
	
	In section \ref{section.3.5.applications} we will use all our previous results
	to show that the smash product of symmetric $T$-spectra induces
	natural pairings (in the motivic stable homotopy category) 
	for the functors $f_{q}$ and $s_{q}$.  We will see
	that for every symmetric $T$-spectrum $X$, and for every $q\in \mathbb Z$:
		\begin{enumerate}
			\item $f_{q}^{\Sigma}X$ is a module (up to homotopy) over the ($-1$)-connective cover of the sphere
						spectrum $f_{0}^{\Sigma}\symmspherespectrum$.
			\item $s_{q}^{\Sigma}X$ is a module (up to homotopy) over the zero slice of the sphere
						spectrum $s_{0}^{\Sigma}\symmspherespectrum$.
		\end{enumerate}
	We will verify that the smash product of symmetric $T$-spectra induces
	natural external pairings in the motivic Atiyah-Hirzebruch spectral sequence
	(see definition \ref{def.3.5.Atiyah-Hirzebruch-spectral-sequence}):
			$$\xymatrix@R=0.5pt{E_{r}^{p,q}(Y;X)\otimes E_{r}^{p',q'}(Y';X') \ar[r]& E_{r}^{p+p',q+q'}(Y\wedge Y';X\wedge X')\\
												(\alpha, \beta) \ar@{|->}[r]& \alpha \smile \beta }
			$$	
	We will also see that for an $A$-module $M$, with $A$ a cofibrant ring spectrum with unit in $\motivicsymmTspectra$,
	which also satisfies the additional hypothesis that are required in section 
	\ref{section.3.5furthermultiplicativeproperties}:
		\begin{enumerate}
			\item $f_{q}^{\Sigma}M$ is again an $A$-module in $\motivicsymmTspectra$ (not just up to homotopy,
						but in a very strict sense).
			\item $s_{q}^{\Sigma}X$ is again an $A$-module in $\motivicsymmTspectra$ (not just up to homotopy,
						but in a very strict sense).
		\end{enumerate}
	Then  we will prove that if the ring spectrum $A$ and its unit map $u: \symmspherespectrum \rightarrow A$
	satisfy the conditions that are required in section \ref{section.3.5furthermultiplicativeproperties},
	the free $A$-module functor $A\wedge -$ induces for every $q\in \mathbb Z$ and for every symmetric
	$T$-spectrum $X$, a natural structure of $A$-module (in $\motivicsymmTspectra$, i.e. not just up to homotopy, 
	but in a very strict sense) on its $q$-slice $s_{q}^{\Sigma}(X)$.
	
	Finally, we will be able to prove a conjecture of M. Levine (see \cite[corollary 11.1.3]{MR2365658}),
	which says that if the base scheme $S$ is a perfect field, then for every $q\in \mathbb Z$ and for every symmetric
	$T$-spectrum $X$, its $q$-slice $s_{q}^{\Sigma}(X)$ is naturally equipped with a module structure
	over the motivic Eilenberg-MacLane spectrum $H\mathbb Z$.
	If we restrict the field even further, considering a field of characteristic zero, then as a consequence we will prove that
	all the slices $s_{q}^{\Sigma}X$ are big motives in the sense of Voevodsky.
		
\begin{section}{The Slice Filtration}
		\label{section-slice-filtration}

	Let $\stablehomotopy$ denote the homotopy category associated
	to $\motivicTspectra$.
	We call $\stablehomotopy$ the \emph{motivic stable homotopy category}.
	We will denote by $[-,-]_{Spt}$ the set of maps between two
	objects in $\stablehomotopy$.
	In \cite{MR1977582} Voevodsky constructs the \emph{slice filtration}
	on motivic stable homotopy theory, 
	using sheaves on the Nisnevich site $\nissite$ instead of simplicial
	presheaves as the underlying category.  In this section we recall
	his construction in the context of simplicial presheaves.

\begin{defi}
		\label{def.3.2.cofibrantstablereplacement}
	Let $Q_{s}$ denote a cofibrant replacement functor
	in $\motivicTspectra$; such that for every $T$-spectrum $X$,
	the natural map:
		$$\xymatrix{Q_{s}X \ar[r]^-{Q_{s}^{X}}& X}
		$$
	is a trivial fibration in $\motivicTspectra$.
\end{defi}
	
\begin{prop}
		\label{prop.3.1.stable-homotopy==triangulated-category}
	The motivic stable homotopy category $\stablehomotopy$
	has a structure of triangulated category defined as follows: 
	\begin{enumerate}
		\item The suspension functor $\Sigma _{T}^{1,0}$ is given by
					$$\xymatrix@R=.5pt{-\wedge S^{1}:\stablehomotopy \ar[r]& \stablehomotopy \\
															X \ar@{|->}[r]& Q_{s}X\wedge S^{1}}
					$$
		\item The distinguished triangles are
					isomorphic to triangles of the form
						$$\xymatrix{A \ar[r]^-{i}& B \ar[r]^-{j}& C \ar[r]^-{k}& \Sigma_{T}^{1,0}A}
						$$
					where $i$ is a cofibration in $\motivicTspectra$, and $C$ is the homotopy cofibre of $i$.
	\end{enumerate}
\end{prop}
\begin{proof}
	Theorem \ref{thm.2.4.stableTspectramodelstr} implies in particular that
	$\motivicTspectra$ is a pointed simplicial model category, and
	theorem \ref{thm.2.4.Tloops-Quillenequiv} implies that the adjunction:
		$$\xymatrix{(-\wedge S^{1},\Omega _{S^{1}},\varphi):\motivicTspectra \ar[r]& \motivicTspectra}
		$$
	is a Quillen equivalence.  The result now follows from the work of
	Quillen in \cite[sections I.2 and I.3]{MR0223432} and the work of 
	Hovey in \cite[chapters VI and VII]{MR1650134}
	(see \cite[proposition 7.1.6]{MR1650134}).
\end{proof}

\begin{note}
		\label{note.3.1.notation-gen-suspensions}
	For $n\in \mathbb Z$, $\Sigma _{T}^{n,0}$ will denote
	the $n^{th}$ iteration of the suspension functor if $n\geq 0$
	($\Sigma _{T}^{0,0}=id$) or the $(-n)^{th}$ iteration of the
	desuspension functor for $n<0$.
\end{note}
	
\begin{lem}
		\label{lem.2.4.compact-respects-coproducts-spectra}
	Let $X\in \pointedmotivic$ be a pointed simplicial presheaf which
	is compact in the sense of Jardine (see definition \ref{def.2.3.compactness}), and let
	$F_{n}(X)$ be the $T$-spectrum constructed in definition \ref{def.2.4.infinitesuspension+shift}.
	Consider an arbitrary collection of $T$-spectra $\{ Z_{i}\}_{i\in I}$ indexed by a set $I$.
	Then 
		$$[F_{n}(X),\coprod _{i\in I}Z_{i}]_{Spt}\cong \coprod _{i\in I}\: [F_{n}(X),Z_{i}]_{Spt}
		$$
\end{lem}
\begin{proof}
	If the indexing set $I$ is finite then the claim holds trivially
	since $\stablehomotopy$ is a triangulated category and therefore finite coproducts and
	finite products are canonically isomorphic.  Thus we can assume that the indexing set $I$
	is infinite.
	
	Choosing a well ordering for the set $I$ there exists a unique ordinal
	$\mu$ which is isomorphic to the ordered set $I$ (see \cite[proposition 10.2.7]{MR1944041}).
	We will prove the lemma by transfinite induction, so assume that 
	for every ordinal $\lambda <\mu$, $F_{n}(X)$ commutes in $\stablehomotopy$
	with coproducts indexed by $\lambda$.
	If $\mu =\lambda +1$, i.e. if $\mu$ is the sucessor of $\lambda$, then
		$$\coprod _{\alpha <\lambda + 1}Z_{\alpha}\cong (\coprod _{\alpha <\lambda} Z_{\alpha})\; \coprod Z_{\lambda} 
		$$
	Therefore 
		$$\xymatrix{[F_{n}(X),\coprod _{\alpha <\lambda + 1}Z_{\alpha}]_{Spt}
							\ar[r]^-{\cong} &
							([F_{n}(X),\coprod _{\alpha <\lambda} Z_{\alpha}]_{Spt})\; \coprod
							\; ([F_{n}(X),Z_{\lambda}]_{Spt})}
		$$
		but by the induction hypothesis
		$$[F_{n}(X),\coprod _{\alpha <\lambda} Z_{\alpha}]_{Spt}
			\cong \coprod _{\alpha <\lambda}[F_{n}(X),Z_{\alpha}]_{Spt}
		$$
	thus
		$$[F_{n}(X),\coprod _{\alpha <\lambda + 1}Z_{\alpha}]_{Spt}\cong
			\coprod _{\alpha <\lambda+1}[F_{n}(X),Z_{\alpha}]_{Spt}
		$$
	as we wanted.
	
	It remains to consider the case when $\mu$ is a limit ordinal.
	In this case proposition 10.2.7 in \cite{MR1944041} implies that
	we can recover the map $\ast \rightarrow \coprod _{\alpha <\mu}Z_{\alpha}$
	as the transfinite composition
	of a $\mu$-sequence:
		$$A_{0}\rightarrow A_{1}\rightarrow \cdots \rightarrow A_{\beta}\rightarrow \cdots
		\; \; (\beta <\mu)
		$$
	where $A_{0}=\ast$, $A_{\beta}=\coprod _{\alpha <\beta}Z_{\alpha}$, and the maps in
	the sequence are the obvious ones.  In particular we have that
	$\coprod _{\alpha <\mu}Z_{\alpha}\cong \varinjlim _{\beta <\mu}A_{\beta}$.
	
	Since $X$ is compact, proposition \ref{prop.2.4.compact-respects-colimits-spectra}
	implies that:
		$$[F_{n}(X),\varinjlim _{\beta <\mu}A_{\beta}]_{Spt}\cong \varinjlim _{\beta <\mu}[F_{n}(X),A_{\beta}]_{Spt}
		$$
	Now using the induction hypothesis we have:
		$$[F_{n}(X),A_{\beta}]_{Spt}\cong \coprod _{\alpha <\beta}\: [F_{n}(X),Z_{\alpha}]_{Spt}
		$$
	and using proposition 10.2.7 in \cite{MR1944041} again, we get:
		$$\varinjlim _{\beta <\mu} \coprod _{\alpha <\beta}\: [F_{n}(X),Z_{\alpha}]_{Spt}
			\cong \coprod _{\alpha <\mu}\: [F_{n}(X),Z_{\alpha}]_{Spt}
		$$
	thus
		$$[F_{n}(X),\coprod _{\alpha <\mu}\ Z_{\alpha}]_{Spt}
			\cong \coprod _{\alpha <\mu}\: [F_{n}(X),Z_{\alpha}]_{Spt}
		$$
	as we wanted.
\end{proof}
	
\begin{prop}
		\label{prop.3.1.stablehomotopy=>compactly-generated}
	The motivic stable homotopy category $\stablehomotopy$ is
	a \emph{compactly generated} triangulated category
	in the sense of Neeman (see \cite[definition 1.7]{MR1308405}).  The set
	of compact generators is given by
	(see definition \ref{def.2.4.infinitesuspension+shift}):				
		$$\compactgenerators
		$$
	i.e. the smallest triangulated subcategory of $\stablehomotopy$ 
	closed under  small coproducts and containing
	all the objects in $C$ coincides with $\stablehomotopy$.
\end{prop}
\begin{proof}
	Since $\stablehomotopy$ is closed under  small coproducts,
	we just need to prove the following two claims:
	\begin{enumerate}
		\item	\label{claim-Neeman1}For every $\generatorNRS \in C$; $\generatorNRS$ commutes with coproducts
					in $\stablehomotopy$, i.e. given a family of $T$-spectra $\{ X_{i}\}_{i\in I}$
					indexed by a set $I$ we have:
						$$[\generatorNRS ,\coprod _{i\in I}X_{i}]_{Spt}
							\cong \coprod _{i\in I}[\generatorNRS ,X_{i}]_{Spt}
						$$
		\item	\label{claim-Neeman2}If a $T$ spectrum $X$ has the following property:
					$[\generatorNRS ,X]_{Spt}=0$
					for every $\generatorNRS \in C$, then $X\cong \ast$ in $\stablehomotopy$.
	\end{enumerate}
	
	(\ref{claim-Neeman1}):  Follows immediately from lemma \ref{lem.2.4.compact-respects-coproducts-spectra}
	since we know by proposition \ref{prop.2.4.T=compact} that
	the pointed simplicial presheaves $S^{r}\wedge \gm ^{s}\wedge U_{+}$ are all
	compact in the sense of Jardine.
		
	(\ref{claim-Neeman2}):  Consider the canonical map
	$X\rightarrow \ast$ in $\motivicTspectra$.  
	Corollary \ref{cor.2.4.detecting-wequivs-repsobjects} together with our hypotheses implies
	that $X\rightarrow \ast$ is a weak equivalence in $\motivicTspectra$, 
	therefore $X\cong \ast$ in $\stablehomotopy$
	as we wanted.
\end{proof}

\begin{cor}
		\label{cor.3.1.detecting-isos-in-SH}
	Let $f:X\rightarrow Y$ be a map in $\stablehomotopy$.
	Then $f$ is an isomorphism if and only if
	$f$ induces an isomorphism of abelian groups:
		$$\xymatrix{[\generatorNRS , X]_{Spt} \ar[r]^-{f_{\ast}}& [\generatorNRS , Y]_{Spt}}
		$$
	for every $\generatorNRS \in C$.
\end{cor}
\begin{proof}
	($\Rightarrow$):  If $f$ is an isomorphism in $\stablehomotopy$ it is clear
	that the induced maps $f_{\ast}$ are isomorphisms of abelian groups for every
	$\generatorNRS \in C$.
	
	($\Leftarrow$):  Complete $f$ to a distinguished triangle in $\stablehomotopy$:
		$$\xymatrix{X \ar[r]^-{f}& Y \ar[r]^-{g}& Z \ar[r]^-{h}& \Tsuspfunctor ^{1,0}X}
		$$
	Then $f$ is an isomorphism if and only if
	$Z\cong \ast$ in $\stablehomotopy$.
	
	Now since the functor $[\generatorNRS ,-]_{Spt}$
	is homological, we get the following long exact sequence of abelian groups:
		$$\xymatrix{\vdots \ar[d] &\\
								[\generatorNRS ,X]_{Spt} \ar[d]^-{f_{\ast}} & \\ 
								[\generatorNRS ,Y]_{Spt} \ar[d]^-{g_{\ast}} & \\ 
								[\generatorNRS ,Z]_{Spt} \ar[d]^-{h_{\ast}} & \\
								[\generatorNRS ,\Sigma _{T}^{1,0}X]_{Spt} 
								\ar[d]^{\Sigma _{T}^{1,0}f_{\ast}} & \ar[l]_-{\cong}^-{\Tsuspfunctor ^{1,0}} 
								[F_{n+1}(S^{r}\wedge \gm ^{s+1}\wedge U_{+}),X]_{Spt} \ar[d]^-{f_{\ast}}\\
								[\generatorNRS ,\Sigma _{T}^{1,0}Y]_{Spt} 
								\ar[d] & [F_{n+1}(S^{r}\wedge \gm ^{s+1}\wedge U_{+}),Y]_{Spt} 
								\ar[l]_-{\cong}^-{\Tsuspfunctor ^{1,0}}\\
								\vdots}
		$$
	But by hypothesis all the maps $f_{\ast}$ are isomorphisms, therefore
	$[\generatorNRS ,Z]_{Spt}=0$
	for every $\generatorNRS \in C$.  Since $\stablehomotopy$ is a compactly generated
	triangulated category (see proposition \ref{prop.3.1.stablehomotopy=>compactly-generated})
	with set of compact generators $C$, we have that
	$Z\cong \ast$.  This implies that $f$ is an isomorphism, as we
	wanted.
\end{proof}

\begin{defi}[Voevodsky, cf. \cite{MR1977582}]
		\label{def.3.1.stable-homotopy-effective}
	We define the \emph{effective motivic stable homotopy category}
	$\stablehomotopyeff \subseteq \stablehomotopy$ as the
	smallest triangulated full subcategory of $\stablehomotopy$
	that is closed under small coproducts and
	contains		
		$$\effcompactgenerators
		$$
\end{defi}

\begin{defi}[Voevodsky, cf. \cite{MR1977582}]
		\label{def.3.1.desusp-stable-homotopy-eff}
	Let $q\in \mathbb Z$.  We define
	$\stablehomotopyeffq \subseteq \stablehomotopy$
	as follows:
	\begin{enumerate}
		\item	If $q=0$, we just take $\stablehomotopyeff$.
		\item	If $q\neq 0$, then $\stablehomotopyeffq$ is the
					smallest triangulated full subcategory of $\stablehomotopy$
					that is closed under  small coproducts and
					contains
						$$\qeffcompactgenerators
						$$
	\end{enumerate}
\end{defi}

\begin{defi}[Voevodsky, cf. \cite{MR1977582}]
	The collection of triangulated subcategories $\stablehomotopyeffq$ 
	for $q\in \mathbb Z$ give a filtration
	on $\stablehomotopy$ which is called
	the \emph{slice filtration}, i.e. we have an inductive system of
	full embeddings
		$$\ldots \subseteq \Sigma _{T}^{q+1}\stablehomotopyeff
			\subseteq \stablehomotopyeffq \subseteq \Sigma _{T}^{q-1}\stablehomotopyeff
			\subseteq \ldots
		$$
	and proposition \ref{prop.3.1.stablehomotopy=>compactly-generated}
	implies that the smallest triangulated subcategory of $\stablehomotopy$
	containing $\stablehomotopyeffq$ for all $q\in \mathbb Z$ and closed
	under  small coproducts coincides with $\stablehomotopy$.
\end{defi}

\begin{prop}
		\label{prop.3.1.stable-q-eff-comp-gen}
	For every $q\in \mathbb Z$, $\stablehomotopyeffq$ is a compactly 
	generated triangulated category in the sense of Neeman,
	where the set of compact generators is
		$$\qeffcompactgenerators
		$$
\end{prop}
\begin{proof}
	By construction $\stablehomotopyeffq$ is closed under small coproducts.
	Therefore we just need to check the following two properties:
	\begin{enumerate}
		\item	\label{claim-Neeman1.effq}For every $\generatorNRS \in C_{eff}^{q}$;
					$\generatorNRS$ commutes with coproducts
					in $\stablehomotopyeffq$, i.e. given a family of $T$-spectra $\{ X_{i}\in \stablehomotopyeffq \}_{i\in I}$
					indexed by a set $I$ we have:
						$$\xymatrix{\Hom _{\stablehomotopyeffq}(\generatorNRS ,
												\coprod _{i\in I}X_{i})\ar[d]^-{\cong} \\
												\coprod _{i\in I}\Hom _{\stablehomotopyeffq}
												(\generatorNRS ,X_{i})}
						$$
		\item	\label{claim-Neeman2.effq}If a $T$-spectrum $X\in \stablehomotopyeffq$ has the following property:
						$$\Hom _{\stablehomotopyeffq}(\generatorNRS ,X)=0$$
					for every $\generatorNRS \in C_{eff}^{q}$, then $X\cong \ast$ in $\stablehomotopyeffq$.
	\end{enumerate}
	
	(\ref{claim-Neeman1.effq}):  Follows immediately from proposition
	\ref{prop.3.1.stablehomotopy=>compactly-generated} since $\stablehomotopyeffq$ is in particular
	a full subcategory of $\stablehomotopy$.
	
	(\ref{claim-Neeman2.effq}):  The natural map
	$X\rightarrow \ast$ is an isomorphism in $\stablehomotopyeffq$ if and only if
	for every $Z\in \stablehomotopyeffq$ we get an induced isomorphism of abelian groups
		$$\xymatrix{\Hom _{\stablehomotopyeffq}(Z,X) \ar[r]^-{\cong}& \Hom _{\stablehomotopyeffq}(Z,\ast)=0}
		$$
	and since $\stablehomotopyeffq$ is a full subcategory of $\stablehomotopy$,
	this last condition is equivalent to:
	For every $Z\in \stablehomotopyeffq$ we have an induced isomorphism of abelian groups
		$$\xymatrix{[Z,X]_{Spt} \ar[r]^-{\cong}& [Z,\ast]_{Spt}=0}
		$$
	Let $\mathcal A _{X}$ be the full subcategory of $\stablehomotopy$
	generated by the $T$-spectra $Y$ satisfying the following
	property
		$$\xymatrix{[\Sigma _{T}^{n,0}Y,X]_{Spt} \ar[r]^-{\cong}& [\Sigma _{T}^{n,0}Y,\ast]_{Spt}}=0
		$$
	for all $n\in \mathbb Z$.
	To finish the proof it is enough to show that
	$\stablehomotopyeffq \subseteq \mathcal A _{X}$, and by construction of $\stablehomotopyeffq$,
	it suffices to prove that $\mathcal A _{X}$ is a triangulated subcategory of $\stablehomotopy$
	which is closed under small coproducts and contains the objects
	$\generatorNRS\in C_{eff}^{q}$.
	The claim that $\mathcal A _{X}$ is triangulated follows immediately from the fact
	that the functor $[-,X]_{Spt}$ is cohomological.
	The claim that $\mathcal A _{X}$ is closed under small coproducts follows from
	the universal property of the coproduct.
	Finally by hypothesis $\mathcal A _{X}$ contains the generators
	$\generatorNRS \in C_{eff}^{q}$.
	This finishes the proof.
\end{proof}

\begin{cor}
		\label{cor.3.1.detecting-isos-stablehomeff-q}
	Let $f:X\rightarrow Y$ be a map in $\stablehomotopyeffq$.
	Then $f$ is an isomorphism if and only if
	one of the following equivalent conditions holds:
		\begin{enumerate}
			\item \label{cor.3.1.detecting-isos-stablehomeff-q.a}For every
						$\generatorNRS \in C_{eff}^{q}$,
						$f$ induces an isomorphism of abelian groups:
							$$\xymatrix{\Hom _{\stablehomotopyeffq}(\generatorNRS ,X) \ar[d]^-{f_{\ast}}\\ 
													\Hom _{\stablehomotopyeffq}(\generatorNRS ,Y)}
							$$
			\item \label{cor.3.1.detecting-isos-stablehomeff-q.b}For every
						$\generatorNRS \in C_{eff}^{q}$,
						$f$ induces an isomorphism of abelian groups:
							$$\xymatrix{[\generatorNRS ,X]_{Spt} \ar[r]^-{i_{q}(f)_{\ast}}&
													[\generatorNRS ,Y]_{Spt}}
							$$
		\end{enumerate}
\end{cor}
\begin{proof}
	Since by construction $\stablehomotopyeffq$ contains $C_{eff}^{q}$ and
	it is a full subcategory of $\stablehomotopy$, we get immediately
	that (\ref{cor.3.1.detecting-isos-stablehomeff-q.a}) and 
	(\ref{cor.3.1.detecting-isos-stablehomeff-q.b}) are equivalent.
	
	We will prove (\ref{cor.3.1.detecting-isos-stablehomeff-q.a}).
	It is clear that if $f$ is an isomorphism then the induced maps
	$f_{\ast}$ considered above are all isomorphisms of abelian groups.
	Conversely, assume that all the induced maps:
		$$\xymatrix{\Hom _{\stablehomotopyeffq}(\generatorNRS ,X) \ar[d]^-{f_{\ast}}\\ 
								\Hom _{\stablehomotopyeffq}(\generatorNRS ,Y)}
		$$
	are isomorphisms for $\generatorNRS \in C_{eff}^{q}$.
	Complete the map $f:X\rightarrow Y$ to a distinguished triangle
	in $\stablehomotopyeffq$:
		$$\xymatrix{X \ar[r]^-{f} & Y \ar[r]^-{g}& Z \ar[r]^-{h}& \Sigma _{T}^{1,0}X}$$
	then $f$ is an isomorphism if and only if $Z\cong \ast$ in $\stablehomotopyeffq$.
	Now since the functor $\Hom _{\stablehomotopyeffq}(\generatorNRS ,-)$
	is homological, we get the following long exact sequence of abelian groups:
		$$\xymatrix@C=-3pc{\vdots \ar[d] &\\
								\Hom _{\stablehomotopyeffq}(\generatorNRS ,X) \ar[d]^-{f_{\ast}} & \\ 
								\Hom _{\stablehomotopyeffq}(\generatorNRS ,Y) \ar[d]^-{g_{\ast}} & \\ 
								\Hom _{\stablehomotopyeffq}(\generatorNRS ,Z) \ar[dd]^-{h_{\ast}} & \\
								& \Hom _{\stablehomotopyeffq}(F_{n+1}(S^{r}\wedge \gm ^{s+1}\wedge U_{+}),X)
								\ar[dd]^-{f_{\ast}} \ar[dl]^-{\cong}_-{\Tsuspfunctor ^{1,0}}\\
								\Hom _{\stablehomotopyeffq}(\generatorNRS ,\Sigma _{T}^{1,0}X) 
								\ar[dd]^{\Sigma _{T}^{1,0}f_{\ast}} & \\ 								  
								& \Hom _{\stablehomotopyeffq}(F_{n+1}(S^{r}\wedge \gm ^{s+1}\wedge U_{+}),Y)
								\ar[ld]^-{\cong}_-{\Tsuspfunctor ^{1,0}}\\
								\Hom _{\stablehomotopyeffq}(\generatorNRS ,\Sigma _{T}^{1,0}Y) 
								\ar[d]  & \\
								\vdots &}
		$$
	But by hypothesis all the maps $f_{\ast}$ are isomorphisms, therefore
		$$\Hom _{\stablehomotopyeffq}(\generatorNRS ,Z)=0$$
	for every $\generatorNRS \in C_{eff}^{q}$.  Since $\stablehomotopyeffq$ is a compactly generated
	triangulated category (see proposition \ref{prop.3.1.stable-q-eff-comp-gen})
	with set of compact generators $C_{eff}^{q}$, we have that
	$Z\cong \ast$.  This implies that $f$ is an isomorphism, as we
	wanted.
\end{proof}

\begin{prop}
		\label{prop.3.1.adjunctions-stableeffq}
	For every $q\in \mathbb Z$ the inclusion 
		$$\xymatrix{i_{q}:\stablehomotopyeffq \ar[r]& \stablehomotopy}
		$$
	has a right adjoint
		$$\xymatrix{r_{q}:\stablehomotopy \ar[r]& \stablehomotopyeffq}
		$$
	which is also an exact functor.
\end{prop}
\begin{proof}
	We have that $\stablehomotopyeffq$ is a compactly generated triangulated category
	(see proposition \ref{prop.3.1.stable-q-eff-comp-gen}), and it is clear that
	the inclusion $i_{q}$ is an exact functor which preserves coproducts.
	Then the existence of the exact right adjoint $r_{q}$ follows from
	theorem 4.1 in \cite{MR1308405}.
\end{proof}

\begin{rmk}
		\label{rmk.3.1.unit=iso}
	\begin{enumerate}	
		\item	Since the inclusion $i_{q}:\stablehomotopyeffq \rightarrow \stablehomotopy$
					is a full embedding,
					we have that the unit of the adjunction $id\stackrel{\tau}{\rightarrow} r_{q}i_{q}$
					is an isomorphism of functors.
		\item	We define $f_{q}=i_{q}r_{q}$.  Then clearly $f_{q+1}f_{q}=f_{q+1}$ and
					there exists a canonical natural transformation $f_{q+1}\rightarrow f_{q}$.
	\end{enumerate}
\end{rmk}

\begin{prop}
		\label{prop.3.1.detecting-fq-isos}
	Fix $q\in \mathbb Z$, and let $g:X\rightarrow Y$ be a map in $\stablehomotopy$.
	Then $f_{q}(g):f_{q}X\rightarrow f_{q}Y$ is
	an isomorphism in $\stablehomotopy$ if and only if
	for every $\generatorNRS \in C_{eff}^{q}$ the induced map:
		$$\xymatrix{[\generatorNRS ,X]_{Spt} \ar[r]^-{g_{\ast}}_{\cong}& [\generatorNRS ,Y]_{Spt}}
		$$
	is an isomorphism of abelian groups.
\end{prop}
\begin{proof}
	We have that $f_{q}=i_{q}r_{q}$, where $i_{q}:\stablehomotopyeffq \rightarrow \stablehomotopy$
	is a full embedding.  Therefore, $f_{q}(g)$ is an isomorphism in $\stablehomotopy$
	if and only if $r_{q}(g)$ is an isomorphism in $\stablehomotopyeffq$.
		
	Hence, corollary \ref{cor.3.1.detecting-isos-stablehomeff-q} implies that $f_{q}(g)$ is an isomorphism
	if and only if for every $\generatorNRS \in C_{eff}^{q}$ the induced map:
		$$\xymatrix{\Hom _{\stablehomotopyeffq}(\generatorNRS ,X) \ar[d]^-{r_{q}(g){\ast}}\\
								\Hom _{\stablehomotopyeffq}(\generatorNRS ,Y)}
		$$
	is an isomorphism.  Fix $\generatorNRS \in C_{eff}^{q}$.
	Finally since $i_{q}, r_{q}$ are adjoint functors
	and $C_{eff}^{q}\subseteq \stablehomotopyeffq$, we have the following commutative diagram,
	where the vertical arrows are all isomorphisms:
		$$\xymatrix@C=-1.5pc{\Hom _{\stablehomotopyeffq}(\generatorNRS ,r_{q}X) \ar[dr]^-{r_{q}(g)_{\ast}}\ar[dd]_-{\cong}&\\
								& \Hom _{\stablehomotopyeffq}(\generatorNRS ,r_{q}Y) \ar[dd]^-{\cong}\\
								[\generatorNRS ,X]_{Spt} \ar[rd]_-{g_{\ast}}& \\
								&[\generatorNRS ,Y]_{Spt}}
		$$
	Therefore, $f_{q}(g)$ is an isomorphism if and only if for every $\generatorNRS \in C_{eff}^{q}$
	the induced map:
		$$\xymatrix{[\generatorNRS ,X]_{Spt} \ar[r]^-{g_{\ast}}_{\cong}& [\generatorNRS ,Y]_{Spt}}
		$$
	is an isomorphism, as we wanted.
\end{proof}

\begin{prop}
		\label{prop.3.1.counit-properties}
	For every $q\in \mathbb Z$ the counit of the adjunction constructed in proposition \ref{prop.3.1.adjunctions-stableeffq},
	$f_{q}=i_{q}r_{q}\stackrel{\theta}{\rightarrow} id$,
	has the following property:
	
	For any $T$-spectrum $X$, and for any compact generator $\generatorNRS \in C_{eff}^{q}$,
	the map $f_{q}X\stackrel{\theta _{X}}{\rightarrow} X$
	in $\stablehomotopy$ induces an isomorphism of
	abelian groups:
		$$\xymatrix{[\generatorNRS , f_{q}X]_{Spt} \ar[r]_-{\cong}^-{\theta _{X\ast}}& [\generatorNRS , X]_{Spt} }
		$$
\end{prop}
\begin{proof}
	Let $\generatorNRS$ be an arbitrary element in $C_{eff}^{q}$.
	Since $\generatorNRS \in \stablehomotopyeffq$ for $n,r,s\geq 0$ with $s-n\geq q$,
	we get the following commutative diagram:
		$$\xymatrix{[\generatorNRS , f_{q}X]_{Spt} \ar[r]^-{\theta _{X\ast}} \ar@{=}[d]& [\generatorNRS , X]_{Spt} \ar@{=}[d]\\ 
								[i_{q}(\generatorNRS), i_{q}r_{q}X]_{Spt} \ar[r]_-{\theta _{X\ast}}& [i_{q}(\generatorNRS), X]_{Spt}}
		$$
	Now using the adjunction between $i_{q}$ and $r_{q}$ we have
	the following commutative diagram:
		$$\xymatrix@C=-3pc{[i_{q}(\generatorNRS), i_{q}r_{q}X]_{Spt} \ar[dr]^-{\theta _{X\ast}} \ar[dd]^-{\cong}& \\
								& [i_{q}(\generatorNRS), X]_{Spt} \ar[dd]_-{\cong}\\
								\Hom _{\stablehomotopyeffq}(\generatorNRS, r_{q}i_{q}r_{q}X) \ar[dr]_-{r_{q}(\theta _{X})_{\ast}}
								&\\
								& \Hom _{\stablehomotopyeffq}(\generatorNRS , r_{q}X)\\
								\Hom _{\stablehomotopyeffq}(\generatorNRS, r_{q}X) \ar[dr]_-{id} 
								\ar[uu]^-{\tau _{r_{q}X}}_-{\cong} &\\
								& \Hom _{\stablehomotopyeffq}(\generatorNRS , r_{q}X) \ar[uu]_-{id}}
		$$
	where $\tau$ is the unit of the adjunction between $i_{q}$ and $r_{q}$.
	This shows that $\theta_{X\ast}$ is an isomorphism, as we wanted.
\end{proof}

\begin{thm}[Voevodsky, cf. \cite{MR1977582}]
		\label{thm.3.1.slicefiltration}
	For every $q\in \mathbb Z$ there exist exact functors
		$$\xymatrix{s_{q}:\stablehomotopy \ar[r]& \stablehomotopy}
		$$
	together with natural transformations
		$$\xymatrix@R=.6pt{\pi _{q}:f_{q} \ar[r]& s_{q} \\
								\sigma _{q}:s_{q} \ar[r]& \Sigma _{T}^{1,0}f_{q+1}}
		$$
	such that the following conditions hold:
	\begin{enumerate}
		\item	\label{thm.3.1.slicefiltration.a}Given any $T$-spectrum $X$, 
					we get the following
					distinguished triangle in $\stablehomotopy$
					\begin{equation}
						\label{equation.3.1.slicefitration}
							\xymatrix{f_{q+1}X \ar[r]& f_{q}X \ar[r]^-{\pi_{q}}& 
							 s_{q}X \ar[r]^-{\sigma_{q}}& \Sigma _{T}^{1,0}f_{q+1}X}
					\end{equation}
		\item	\label{thm.3.1.slicefiltration.b}For any $T$-spectrum $X$, 
					$s_{q}X$ is in $\stablehomotopyeffq$.
		\item	\label{thm.3.1.slicefiltration.c}For any $T$-spectrum $X$, 
					and for any $T$ spectrum $Y$ in 
					$\stablehomotopyeffqplusone$, 
					$[Y,s_{q}X]_{Spt}=0$.
	\end{enumerate}
\end{thm}
\begin{proof}
	Since the triangulated categories $\Tsuspfunctor ^{q+1}\stablehomotopyeff$ and $\stablehomotopyeffq$
	are both compactly generated (see proposition \ref{prop.3.1.stable-q-eff-comp-gen}), 
	the result follows from 
	propositions 9.1.19 and 9.1.8 in \cite{MR1812507}.	
\end{proof}

\begin{defi}[Voevodsky]
		\label{def.3.1.slicefiltration}
	Given an arbitrary $T$-spectrum $X$,
	the sequence of distinguished triangles (\ref{equation.3.1.slicefitration})
	is called the \emph{slice tower} of $X$.  The $T$-spectrum $s_{q}X$ is called
	the \emph{$q$-slice} of $X$, and the $T$-spectrum $f_{q}X$ is called
	the \emph{($q-1$)-connective cover} of $X$.
\end{defi}

\begin{thm}
		\label{thm.3.1.motivictower}
	For every $q\in \mathbb Z$ there exist exact functors
		$$\xymatrix{s_{<q}:\stablehomotopy \ar[r]& \stablehomotopy}
		$$
	together with natural transformations
		$$\xymatrix@R=.6pt{\pi _{<q}:id \ar[r]& s_{<q} \\
								\sigma _{<q}:s_{<q} \ar[r]& \Sigma _{T}^{1,0}f_{q}}
		$$
	such that the following conditions hold:
	\begin{enumerate}
		\item	\label{thm.3.1.motivictower.a}Given any 
					$T$-spectrum $X$, we get the following
					distinguished triangle in $\stablehomotopy$
					\begin{equation}
						\label{equation.3.1.motivictower}
							\xymatrix{f_{q}X \ar[r]& X \ar[r]^-{\pi_{<q}}& 
							 s_{<q}X \ar[r]^-{\sigma_{<q}}& \Sigma _{T}^{1,0}f_{q}X}
					\end{equation}
		\item	\label{thm.3.1.motivictower.b}For any $T$-spectrum $X$, 
					and for any $T$ spectrum $Y$ in 
					$\stablehomotopyeffq$, 
					$[Y,s_{<q}X]_{Spt}=0$.
	\end{enumerate}
\end{thm}
\begin{proof}
	The result follows from 
	propositions 9.1.19 and 9.1.8 in \cite{MR1812507},
	using the fact that
	the triangulated categories $\stablehomotopyeffq$ and $\stablehomotopy$
	are both compactly generated
	(see propositions \ref{prop.3.1.stablehomotopy=>compactly-generated} and 
	\ref{prop.3.1.stable-q-eff-comp-gen}).	
\end{proof}

\begin{prop}
		\label{prop.3.1.octahedral-axiom}
	Let $X$ be an arbitrary $T$-spectrum.  Then
	for every $q\in \mathbb Z$,
	we have the following commutative diagram, where all
	the rows and columns are distinguished triangles in $\stablehomotopy$:
		$$\xymatrix{f_{q+1}X \ar[rr] \ar@{=}[d] && f_{q}X \ar[rr]^-{\pi _{q}} \ar[d]&& s_{q}X \ar[rr]^-{\sigma _{q}} \ar[d]^-{j_{q}}&& 
								\Tsuspfunctor ^{1,0}f_{q+1}X \ar@{=}[d]\\
								f_{q+1}X \ar[rr] \ar[d]&& X \ar[rr]^-{\pi _{<q+1}} \ar[d]_-{\pi _{<q}}&& s_{<q+1}X \ar[rr]^-{\sigma _{<q+1}} 
								\ar[d]^-{t_{q}}&& \Tsuspfunctor ^{1,0}f_{q+1}X \ar[d]\\
								\ast \ar[rr] \ar[d]&& s_{<q}X \ar@{=}[rr] \ar[d]_-{\sigma _{<q}}&& s_{<q}X \ar[rr] \ar[d]&& \ast \ar[d]\\
								\Tsuspfunctor ^{1,0}f_{q+1}X \ar[rr]&& \Tsuspfunctor ^{1,0}f_{q}X \ar[rr]_-{\Tsuspfunctor ^{1,0}\pi _{q}}&&
								\Tsuspfunctor ^{1,0}s_{q}X \ar[rr]_-{\Tsuspfunctor ^{1,0}\sigma _{q}}&& \Tsuspfunctor ^{2,0}f_{q+1}X}
		$$
\end{prop}
\begin{proof}
	It follows from theorems \ref{thm.3.1.slicefiltration} and \ref{thm.3.1.motivictower}, together
	with the octahedral axiom applied to the following commutative diagram:
		$$\xymatrix{f_{q+1}X \ar[rr] \ar[dr]&& f_{q}X \ar[dl]\\
								& X &}
		$$
\end{proof}

\begin{prop}
		\label{prop.3.1.detecting-s<q-isos}
	Fix $q\in \mathbb Z$ and let $f:X\rightarrow Y$ be a map
	in $\stablehomotopy$.
	Then 
		$$s_{<q}f:s_{<q}X\rightarrow s_{<q}Y$$ 
	is an isomorphism in $\stablehomotopy$ if and only if
	$f$ induces the following isomorphisms
	of abelian groups:
		$$\xymatrix{[\generatorNRS , s_{<q}X]_{Spt} \ar[r]^-{(s_{<q}f)_{\ast}}& 
								[\generatorNRS , s_{<q}Y]_{Spt}}
		$$
	for every $\generatorNRS \notin C_{eff}^{q}$.
\end{prop}
\begin{proof}
	($\Rightarrow$):  Assume that $s_{<q}f$ is an isomorphism.
	Then it is clear that $(s_{<q}f)_{\ast}$ is also an isomorphism
	for every $\generatorNRS \notin C_{eff}^{q}$.
	
	($\Leftarrow$):  Corollary  \ref{cor.3.1.detecting-isos-in-SH} implies that 
	$s_{<q}f$ is an isomorphism in $\stablehomotopy$
	if and only if for every $\generatorNRS \in C$, the induced maps:
		$$\xymatrix{[\generatorNRS , s_{<q}X]_{Spt} \ar[r]^-{(s_{<q}f)_{\ast}}& 
								[\generatorNRS , s_{<q}Y]_{Spt}}
		$$
	are isomorphisms of abelian groups.
	
	But theorem \ref{thm.3.1.motivictower}(\ref{thm.3.1.motivictower.b}) implies that
	for every $\generatorNRS \in C_{eff}^{q}$, we have:
		$$\xymatrix{0\cong[\generatorNRS , s_{<q}X]_{Spt} \ar[r]^-{(s_{<q}f)_{\ast}}_-{\cong}& 
								[\generatorNRS , s_{<q}Y]_{Spt}\cong 0}
		$$
	thus $(s_{<q}f)_{\ast}$ is an isomorphism in this case.
		
	Thus in order to show that $s_{<q}f$ is an isomorphism,
	we only need to check that for every $\generatorNRS \notin C_{eff}^{q}$,
	the induced maps:
		$$\xymatrix{[\generatorNRS , s_{<q}X]_{Spt} \ar[r]^-{(s_{<q}f)_{\ast}}& 
								[\generatorNRS , s_{<q}Y]_{Spt}}
		$$
	are all isomorphisms of abelian groups;
	but this holds by hypothesis.  This finishes the proof.
\end{proof}

\begin{prop}
		\label{prop.3.1.detecting-sq-isos}
	Fix $q\in \mathbb Z$ and let $f:X\rightarrow Y$ be a map
	in $\stablehomotopy$.
	Then 
		$$s_{q}f:s_{q}X\rightarrow s_{q}Y$$ 
	is an isomorphism in $\stablehomotopy$ if and only if
	$f$ induces the following isomorphisms
	of abelian groups:
		$$\xymatrix{[\generatorNRS , s_{q}X]_{Spt} \ar[r]^-{(s_{q}f)_{\ast}}& 
								[\generatorNRS , s_{q}Y]_{Spt}}
		$$
	for every $\generatorNRS \in C_{eff}^{q}$ where $s-n=q$.
\end{prop}
\begin{proof}
	($\Rightarrow$):  Assume that $s_{q}f$ is an isomorphism.
	Then it is clear that $(s_{q}f)_{\ast}$ is also an isomorphism
	for every $\generatorNRS \in C_{eff}^{q}$ with $s-n=q$.
	
	($\Leftarrow$):
	Theorem \ref{thm.3.1.slicefiltration}(\ref{thm.3.1.slicefiltration.b}) implies that
	$s_{q}X$ and $s_{q}Y$ are both in $\stablehomotopyeffq$.
	Therefore using corollary \ref{cor.3.1.detecting-isos-stablehomeff-q}
	and the fact that $\stablehomotopyeffq$ is a full subcategory of
	$\stablehomotopy$, we have that $s_{q}f$ is an isomorphism
	if and only if the maps:
		$$\xymatrix{[\generatorNRS , s_{q}X]_{Spt} \ar[r]^-{(s_{q}f)_{\ast}}& 
								[\generatorNRS , s_{q}Y]_{Spt}}
		$$
	are all isomorphisms of abelian groups for every
	$\generatorNRS \in C_{eff}^{q}$.  
	
	But if
	$s-n\geq q+1$, we have that $\generatorNRS$
	is in fact in $\Tsuspfunctor ^{q+1}\stablehomotopyeff$;
	and using theorem \ref{thm.3.1.slicefiltration}(\ref{thm.3.1.slicefiltration.c}) again,
	we have that in this case:
		$$\xymatrix{0\cong[\generatorNRS , s_{q}X]_{Spt} \ar[r]^-{(s_{q}f)_{\ast}}_-{\cong}& 
								[\generatorNRS , s_{q}Y]_{Spt}\cong 0}
			$$
	
	Thus in order to show that $s_{q}f$ is an isomorphism,
	we only need to check that the maps:
		$$\xymatrix{[\generatorNRS , s_{q}X]_{Spt} \ar[r]^-{(s_{q}f)_{\ast}}& 
								[\generatorNRS , s_{q}Y]_{Spt}}
		$$
	are all isomorphisms of abelian groups,
	for every $\generatorNRS \in C_{eff}^{q}$ with $s-n=q$.  This finishes the proof.
\end{proof}

\end{section}

\begin{section}{Model Structures for the Slice Filtration}
		\label{section.3.2.modstrslicefilt}
	
	Our goal in this section is to use the cellularity of $\motivicTspectra$ 
	(see theorem \ref{thm.2.5.cellularity-motivic-stable-str}), to construct 
	using Hirschhorn's localization techniques, several families of model structures on $\Tspectra$
	via left and right Bousfield localization.  This new model structures will
	provide liftings in a suitable sense for the functors 
		$$f_{q},s_{<q},s_{q}:\stablehomotopy \rightarrow \stablehomotopy$$
	described in section \ref{section-slice-filtration}.
	
	The first family of model structures on $\Tspectra$ will be constructed
	via right Bousfield localization.
	These model structures will have the property
	that the cofibrant replacement functor coincides in a suitable sense with
	the functor $f_{q}$ defined in remark \ref{rmk.3.1.unit=iso}.
	This will provide a natural lifting of Voevodsky's slice filtration to the level
	of model categories.
	
\begin{thm}
		\label{thm.3.2.connective-model-structure}
	Fix $q\in \mathbb Z$.  Consider the following set of objects in $\motivicTspectra$
		$$\qeffcompactgenerators
		$$
	Then the right Bousfield localization of $\motivicTspectra$
	with respect to the class of $C_{eff}^{q}$-colocal equivalences exists
	(see definitions \ref{def1.2.cellularization} and \ref{def1.2.1.rightBousloc}).
	This model structure 
	will be called  \emph{($q-1$)-connected motivic stable}, and
	the category of $T$-spectra  equipped with the ($q-1$)-connected motivic
	stable model structure will be denoted by
	$\qconnectedTspectra$.  Furthermore
	$\qconnectedTspectra$ is a right proper and simplicial model
	category.
	The homotopy category associated to $\qconnectedTspectra$
	will be denoted by $\qconnectedstablehomotopy$.
\end{thm}
\begin{proof}
	Theorems \ref{thm.2.4.stableTspectramodelstr} and
	\ref{thm.2.5.cellularity-motivic-stable-str} imply that
	$\motivicTspectra$  is
	cellular, proper and simplicial.
	Therefore we can apply theorem 5.1.1 in
	\cite{MR1944041} to construct
	the right Bousfield localization of $\motivicTspectra$ with respect to the class of $C_{eff}^{q}$-colocal
	equivalences.
	Using theorem 5.1.1 in
	\cite{MR1944041} again, we have
	that this new model structure  is  right proper and
	simplicial.
\end{proof}

\begin{defi}
		\label{def.3.3.Rq-cofibrant-replacement}
	Fix $q\in \mathbb Z$.  Let $C_{q}$ denote a cofibrant replacement functor in 
	$\qconnectedTspectra$; such that for every $T$-spectrum $X$, the natural map
		$$\xymatrix{C_{q}X \ar[r]^-{C_{q}^{X}}& X}
		$$ 
	is a trivial fibration in $\qconnectedTspectra$, and
	$C_{q}X$ is always a $C_{eff}^{q}$-colocal $T$-spectrum.
\end{defi}

\begin{prop}
		\label{prop.3.2.IQTJ-fibrant-replacement-all-Rq}
	Fix $q\in \mathbb Z$.  Then $IQ_{T}J$ is also a fibrant
	replacement functor in $\qconnectedTspectra$
	(see corollary \ref{cor.2.4.stablyfibrant-injective-model}).
\end{prop}
\begin{proof}
	Since $\qconnectedTspectra$ is the right Bousfield localization of
	$\motivicTspectra$ with respect to the $C_{eff}^{q}$-colocal equivalences, by construction
	we have that the fibrations and the trivial cofibrations are indentical in
	$\qconnectedTspectra$ and $\motivicTspectra$ respectively.  This implies that for every
	$T$-spectrum $X$, $IQ_{T}JX$ is fibrant in $\qconnectedTspectra$, and 
	using \cite[proposition 3.1.5]{MR1944041} we have that
	the natural map:
		$$\xymatrix{X\ar[rr]^{IQ_{T}J^{X}}&& IQ_{T}JX}
		$$
	is a weak equivalence in $\qconnectedTspectra$.  Hence $IQ_{T}J$ is also a fibrant replacement
	functor for $\qconnectedTspectra$.
\end{proof}

\begin{prop}
		\label{prop.3.2.classif-Cqcolocalequivs}
	Fix $q\in \mathbb Z$ and let $f:X\rightarrow Y$ be a map in $\motivicTspectra$.
	Then $f$ is a $C_{eff}^{q}$-colocal equivalence
	if and only if for every $\generatorNRS \in C_{eff}^{q}$,
	$f$ induces the
	following isomorphisms of abelian groups:
		$$\xymatrix{[\generatorNRS ,X]_{Spt} \ar[r]^-{f_{\ast}}& [\generatorNRS ,Y]_{Spt}}
		$$
\end{prop}
\begin{proof}
	($\Rightarrow$): Assume that $f$ is a $C_{eff}^{q}$-colocal equivalence.
	Since all the compact generators $\generatorNRS$ are cofibrant
	in $\motivicTspectra$, we have that $f$ is a
	$C_{eff}^{q}$-colocal equivalence if and only if
	the following maps are weak equivalences of simplicial sets:
		$$\xymatrix{Map(\generatorNRS ,IQ_{T}JX) \ar[d]^-{(IQ_{T}Jf)_{\ast}}\\ 
								Map(\generatorNRS ,IQ_{T}JY)}
		$$
	for every $\generatorNRS \in C_{eff}^{q}$.  Since
	$\motivicTspectra$ is a simplicial model category
	and $\generatorNRS$ is cofibrant,
	we have that $Map(\generatorNRS ,IQ_{T}JX)$ and $Map(\generatorNRS ,IQ_{T}JY)$ are both Kan complexes,  
	thus we get the following commutative diagram
	where the top row and the vertical maps are all isomorphisms
	of abelian groups:
		$$\xymatrix@C=-1.5pc{\pi_{0}Map(\generatorNRS ,IQ_{T}JX) \ar[dr]^-{(IQ_{T}Jf)_{\ast}}_-{\cong} \ar[dd]_-{\cong}&\\ 
								& \pi_{0}Map(\generatorNRS ,IQ_{T}JY) \ar[dd]^-{\cong}\\
								[\generatorNRS ,X]_{Spt} \ar[dr]_-{f_{\ast}}& \\ 
								& [\generatorNRS ,Y]_{Spt}}
		$$
	Therefore
		$$\xymatrix{[\generatorNRS ,X]_{Spt} \ar[r]_-{f_{\ast}}& [\generatorNRS ,Y]_{Spt}}
		$$
	is an isomorphism of abelian groups
	for every $\generatorNRS \in C_{eff}^{q}$, as we wanted.
	
	($\Leftarrow$):  Fix $\generatorNRS \in C_{eff}^{q}$.  Let
	$\omega _{0}$, $\eta _{0}$ denote the base points corresponding to
	$Map_{\ast}(\generatorNRS ,IQ_{T}JX)$ and $Map_{\ast}(\generatorNRS ,IQ_{T}JY)$
	respectively.
	We need to show that the map:
		$$\xymatrix{Map(\generatorNRS ,IQ_{T}JX) \ar[d]^-{(IQ_{T}Jf)_{\ast}}\\ 
								Map(\generatorNRS ,IQ_{T}JY)}
		$$
	is a weak equivalence of simplicial sets.  
	
	We know that the map
		$$j:F_{n+1}(S^{r+1}\wedge \gm ^{s+1}\wedge U_{+})\rightarrow \generatorNRS$$ 
	which is adjoint to the identity map
		$$id:S^{r+1}\wedge \gm ^{s+1}\wedge U_{+} \rightarrow Ev_{n+1}(\generatorNRS)=S^{r+1}\wedge \gm ^{s+1}\wedge U_{+}$$
	is a weak equivalence in $\motivicTspectra$.
	Now since $\generatorNRS$ and $F_{n+1}(S^{r+1}\wedge \gm ^{s+1}\wedge U_{+})$ are both cofibrant
	and $\motivicTspectra$ is a simplicial model category, we can apply Ken Brown's lemma (see lemma \ref{lem1.1.KenBrown})
	to conclude that the horizontal maps in the following commutative diagram are weak
	equivalences of simplicial sets:
		$$\xymatrix@C=-2pc{Map(\generatorNRS ,IQ_{T}JX) 
								\ar[dr]^-{j^{\ast}}\ar[dd]_-{(IQ_{T}Jf)_{\ast}} &\\ 
								& Map(F_{n+1}(S^{r+1}\wedge \gm ^{s+1}\wedge U_{+}) ,IQ_{T}JX) \ar[dd]^-{(IQ_{T}Jf)_{\ast}}\\		
								Map(\generatorNRS ,IQ_{T}JY)
								\ar[dr]_-{j^{\ast}}&\\ 
								& Map(F_{n+1}(S^{r+1}\wedge \gm ^{s+1}\wedge U_{+}) ,IQ_{T}JY)}
		$$
	Hence by the two out of three property for weak equivalences, it is enough to show that
	the following induced map
		$$\xymatrix{Map(F_{n+1}(S^{r+1}\wedge \gm ^{s+1}\wedge U_{+}) ,IQ_{T}JX) \ar[d]^-{(IQ_{T}Jf)_{\ast}}\\
								Map(F_{n+1}(S^{r+1}\wedge \gm ^{s+1}\wedge U_{+}) ,IQ_{T}JY)}
		$$
	is a weak equivalence of simplicial sets.
	
	On the other hand, since $\motivicTspectra$ is a pointed simplicial model category,
	we have that lemma 6.1.2 in \cite{MR1650134} together 
	with remark \ref{rmk.2.4.simp-simppointed-structures-Tspectra}(\ref{rmk.2.4.simp-simppointed-structures-Tspectra.b})
	imply that
	the following diagram is commutative for $k\geq 0$:
		$$\xymatrix@C=-3pc{\pi_{k, \omega _{0}}Map(\generatorNRS ,IQ_{T}JX) \ar[dr]^-{(IQ_{T}Jf)_{\ast}}\ar@{=}[dd]&\\  
								& \pi_{k, \eta _{0}}Map(\generatorNRS ,IQ_{T}JY)\ar@{=}[dd]\\
								\pi_{k, \omega _{0}}Map_{\ast}(\generatorNRS ,IQ_{T}JX) \ar[dr]^-{(IQ_{T}Jf)_{\ast}}\ar[dd]_-{\cong}&\\  
								& \pi_{k, \eta _{0}}Map_{\ast}(\generatorNRS ,IQ_{T}JY)\ar[dd]^-{\cong}\\
								[\generatorNRS \wedge S^{k}, IQ_{T}JX]_{Spt} \ar[rd]^-{(IQ_{T}Jf)_{\ast}}&\\ 
								& [\generatorNRS \wedge S^{k}, IQ_{T}JY]_{Spt}\\
								[\generatorNRS \wedge S^{k}, X]_{Spt} \ar[dr]^-{f_{\ast}} \ar[uu]^-{\cong} \ar[dd]_-{\cong}&\\ 
								& [\generatorNRS \wedge S^{k}, Y]_{Spt} \ar[uu]_-{\cong} \ar[dd]^-{\cong}\\
								[F_{n}(S^{k+r}\wedge \gm ^{s}\wedge U_{+}),X]_{Spt} \ar[dr]_-{f_{\ast}}^-{\cong}&\\ 
								& [F_{n}(S^{k+r}\wedge \gm ^{s}\wedge U_{+}),Y]_{Spt}}
		$$
	but by hypothesis we have that the bottom row is an isomorphism of abelian groups.
	Therefore all the maps in the top row are also isomorphisms.  Then
	for every $\generatorNRS \in C_{eff}^{q}$, the induced map
		$$\xymatrix{Map(\generatorNRS ,IQ_{T}JX) \ar[rr]^-{(IQ_{T}Jf)_{\ast}}&& 
								Map(\generatorNRS ,IQ_{T}JY)}
		$$ 
	is a weak equivalence when it is restricted to the path component of $Map(\generatorNRS ,IQ_{T}JX)$
	containing $\omega _{0}$.  But $\generatorNNRSS$ is also in $C_{eff}^{q}$, therefore
	the following induced map
		$$\xymatrix{Map_{\ast}(S^{1},Map_{\ast}(F_{n+1}(S^{r}\wedge \gm ^{s+1}\wedge U_{+}) ,IQ_{T}JX)) 
								\ar[d]_-{(IQ_{T}Jf)_{\ast}}\\
								Map_{\ast}(S^{1},Map_{\ast}(F_{n+1}(S^{r}\wedge \gm ^{s+1}\wedge U_{+}) ,IQ_{T}JY))}
		$$
	is a weak equivalence of simplicial sets, since taking $S^{1}$-loops kills the path components that do not
	contain the base point.  
	
	Finally, since $\motivicTspectra$ is a simplicial model category we have that the rows in the following
	commutative diagram are isomorphisms:
		$$\xymatrix@C=-6pc{Map_{\ast}(S^{1},Map_{\ast}(F_{n+1}(S^{r}\wedge \gm ^{s+1}\wedge U_{+}) ,IQ_{T}JX)) \ar[dr]^-{\cong} 
								\ar[dd]_-{(IQ_{T}Jf)_{\ast}}&\\
								& Map_{\ast}( F_{n+1}(S^{r}\wedge \gm ^{s+1}\wedge U_{+})\wedge S^{1} ,IQ_{T}JX)
								 \ar[dd]^-{(IQ_{T}Jf)_{\ast}}\\
								Map_{\ast}(S^{1},Map_{\ast}(F_{n+1}(S^{r}\wedge \gm ^{s+1}\wedge U_{+}) ,IQ_{T}JY)) \ar[dr]_-{\cong}&\\ 
								& Map_{\ast}( F_{n+1}(S^{r}\wedge \gm ^{s+1}\wedge U_{+})\wedge S^{1} ,IQ_{T}JY)}
		$$
	Hence the two out of three property for weak equivalences implies that the right vertical map
	is a weak equivalence of simplicial sets.  
	But $F_{n+1}(S^{r}\wedge \gm ^{s+1}\wedge U_{+})\wedge S^{1}$
	is clearly isomorphic to $F_{n+1}(S^{r+1}\wedge \gm ^{s+1}\wedge U_{+})$, therefore
	the induced map
		$$\xymatrix{Map(F_{n+1}(S^{r+1}\wedge \gm ^{s+1}\wedge U_{+}),IQ_{T}JX) \ar[d]^-{(IQ_{T}Jf)_{\ast}}\\ 
								Map(F_{n+1}(S^{r+1}\wedge \gm ^{s+1}\wedge U_{+}),IQ_{T}JY)}
		$$
	is a weak equivalence, as we wanted.
\end{proof}

\begin{cor}
		\label{cor.3.2.detecting-Cqeff-colocalequivs2}
	Fix $q\in \mathbb Z$ and let $f:X\rightarrow Y$ be a map in $\motivicTspectra$.
	Then $f$ is a $C_{eff}^{q}$-colocal equivalence
	if and only if the following map
		$$\xymatrix{r_{q}X \ar[r]^{r_{q}(f)}& r_{q}Y}
		$$
	is an isomorphism in $\stablehomotopyeffq$.
\end{cor}
\begin{proof}
	The result follows immediately from
	proposition \ref{prop.3.2.classif-Cqcolocalequivs} and 
	corollary \ref{cor.3.1.detecting-isos-stablehomeff-q}(\ref{cor.3.1.detecting-isos-stablehomeff-q.b}).
\end{proof}

\begin{cor}
		\label{cor.3.2.detecting-trivialspectra-Cqeff}
	Fix $q\in \mathbb Z$ and let $X$ be an arbitrary $T$-spectrum $X$.  Then $X\cong \ast$
	in $\qconnectedstablehomotopy$ if and only if
	the following condition holds:
	
		For every $\generatorNRS \in C_{eff}^{q}$:
						$$[\generatorNRS, X]_{Spt}\cong 0
						$$
\end{cor}
\begin{proof}
	We have that
	$X$ is isomorphic to $\ast$ in $\qconnectedstablehomotopy$ if and only if
	the map $\ast \rightarrow X$ is a $C_{eff}^{q}$-colocal equivalence.
	But corollary \ref{cor.3.2.detecting-Cqeff-colocalequivs2} implies that
	$\ast \rightarrow X$ is a $C_{eff}^{q}$-colocal equivalence if and only if
		$$\xymatrix{\ast \cong r_{q}(\ast) \ar[r]& r_{q}X}
		$$
	becomes an isomorphism in $\stablehomotopyeffq$.
	
	Finally by corollary \ref{cor.3.1.detecting-isos-stablehomeff-q}(\ref{cor.3.1.detecting-isos-stablehomeff-q.b})
	we have that $\ast \rightarrow r_{q}X$ is an isomorphism in $\stablehomotopyeffq$ if
	and only if for every $\generatorNRS \in C_{eff}^{q}$ the following
	induced maps are isomorphisms of abelian groups:
		$$\xymatrix{0\cong [\generatorNRS, \ast]_{Spt} \ar[r]^-{\cong}& [\generatorNRS ,X]_{Spt}}
		$$
	as we wanted.
\end{proof}

\begin{lem}
		\label{lem.3.2.Cq-colocalequivs-stable-S1desuspension}
	Fix $q\in \mathbb Z$ and let $f:X\rightarrow Y$ be a map in $\motivicTspectra$, then $f$ is a
	$C^{q}_{eff}$-colocal equivalence if and only if
	$\Omega _{S^{1}} IQ_{T}J(f)$ is a $C^{q}_{eff}$-colocal equivalence.
\end{lem}
\begin{proof}
	Assume that $f$ is a $C^{q}_{eff}$-colocal equivalence.
	We need to show that $\Omega _{S^{1}}IQ_{T}J(f)$ is a
	$C^{q}_{eff}$-colocal equivalence.  
	Fix $\generatorNRS \in C_{eff}^{q}$.
	Since $\motivicTspectra$
	is a simplicial model category
	and all the compact generators $\generatorNRS$ are cofibrant,
	we have the following commutative diagram:
		$$\xymatrix@C=-1pc{[\generatorNRS ,\Omega _{S^{1}}IQ_{T}JX]_{Spt} \ar[dr]^-{(IQ_{T}Jf)_{\ast}} \ar[dd]_-{\cong}&\\
								& [\generatorNRS ,\Omega _{S^{1}}IQ_{T}JY]_{Spt} \ar[dd]^-{\cong}\\
								[\generatorNRS \wedge S^{1},X]_{Spt} \ar[dr]^-{f_{\ast}} \ar[dd]_-{\cong}&\\ 
								& [\generatorNRS \wedge S^{1}, Y]_{Spt} \ar[dd]^-{\cong}\\
								[F_{n}(S^{r+1}\wedge \gm ^{s}\wedge U_{+}),X]_{Spt} \ar[dr]_-{f_{\ast}}^-{\cong}&\\ 
								& [F_{n}(S^{r+1}\wedge \gm ^{s}\wedge U_{+}),Y]_{Spt}}
		$$
	but using proposition \ref{prop.3.2.classif-Cqcolocalequivs} and the fact that
	$f$ is a $C_{eff}^{q}$-colocal equivalence, we have that the
	bottom row is an isomorphism, therefore the top row is also an isomorphism.  
	Using proposition \ref{prop.3.2.classif-Cqcolocalequivs}
	again, we have that $\Omega _{S^{1}}IQ_{T}J(f)$ is a
	$C_{eff}^{q}$-colocal equivalence, as we wanted. 
	
	Conversely, assume that $\Omega _{S^{1}IQ_{T}J(f)}$ is a
	$C_{eff}^{q}$-colocal equivalence.
	Fix $\generatorNRS \in C_{eff}^{q}$.
	Proposition \ref{prop.3.2.classif-Cqcolocalequivs}
	implies that the top row in the following commutative diagram is
	an isomorphism:
		$$\xymatrix@C=-2pc{[F_{n+1}(S^{r}\wedge \gm ^{s+1}\wedge U_{+}) ,\Omega _{S^{1}}IQ_{T}JX]_{Spt} \ar[dr]^-{IQ_{T}J(f)_{\ast}}_-{\cong} 
								\ar[dd]_-{\cong}&\\ 
								& [F_{n+1}(S^{r}\wedge \gm ^{s+1}\wedge U_{+}) ,\Omega _{S^{1}}IQ_{T}JY]_{Spt} \ar[dd]^-{\cong}\\
								[F_{n+1}(S^{r}\wedge \gm ^{s+1}\wedge U_{+}) \wedge S^{1},X]_{Spt} \ar[dr]^-{f_{\ast}} &\\ 
								& [F_{n+1}(S^{r}\wedge \gm ^{s+1}\wedge U_{+}) \wedge S^{1}, Y]_{Spt} \\
								[\generatorNRS ,X]_{Spt} \ar[dr]_-{f_{\ast}} \ar[uu]^-{\cong}&\\ 
								& [\generatorNRS ,Y]_{Spt} \ar[uu]_-{\cong}}
		$$
	therefore the bottom row is also an isomorphism.
	Finally using Proposition \ref{prop.3.2.classif-Cqcolocalequivs} again,
	we have that $f$ is a $C_{eff}^{q}$-colocal equivalence.
	This finishes the proof.
\end{proof}

\begin{cor}
		\label{cor.3.3.Suspension=>qconn-Quillen-equiv}
	For every $q\in \mathbb Z$, the adjunction
		$$\xymatrix{(-\wedge S^{1},\Omega _{S^{1}},\varphi):\qconnectedTspectra \ar[r]& 
								\qconnectedTspectra}
		$$
	is a Quillen equivalence.
\end{cor}
\begin{proof}
	Using corollary 1.3.16 in \cite{MR1650134} and
	proposition \ref{prop.3.2.IQTJ-fibrant-replacement-all-Rq} we have that it suffices to
	verify the following two conditions:
		\begin{enumerate}
			\item	\label{cor.3.3.Suspension=>qconn-Quillen-equiv.a}For every cofibrant object $X$
						in $\qconnectedTspectra$, the following composition
							$$\xymatrix{X\ar[r]^-{\eta _{X}}& \Omega _{S^{1}}(X\wedge S^{1}) \ar[rr]^-{\Omega _{S^{1}}IQ_{T}J^{X\wedge S^{1}}}&& 
								\Omega _{S^{1}}IQ_{T}J(X\wedge S^{1})}
							$$
						is a $C_{eff}^{q}$-colocal equivalence.
			\item	\label{cor.3.3.Suspension=>qconn-Quillen-equiv.b}$\Omega _{S^{1}}$ reflects $C_{eff}^{q}$-colocal equivalences
						between fibrant objects in $\qconnectedTspectra$.
		\end{enumerate}
		
	(\ref{cor.3.3.Suspension=>qconn-Quillen-equiv.a}):  By construction $\qconnectedTspectra$ is
	a right Bousfield localization of $\motivicTspectra$, therefore the identity
	functor 
		$$\xymatrix{id:\qconnectedTspectra \ar[r]& \motivicTspectra}
		$$ 
	is a left Quillen functor.
	Thus $X$ is also cofibrant in $\motivicTspectra$.
	Since the adjunction $(-\wedge S^{1},\Omega _{S^{1}},\varphi)$
	is a Quillen equivalence on $\motivicTspectra$, 
	\cite[proposition 1.3.13(b)]{MR1650134} implies that the following composition
	is a weak equivalence in $\motivicTspectra$:
		$$\xymatrix{X\ar[r]^-{\eta _{X}}& \Omega _{S^{1}}(X\wedge S^{1}) \ar[rr]^-{\Omega _{S^{1}}IQ_{T}J^{X\wedge S^{1}}}&& 
								\Omega _{S^{1}}IQ_{T}J(X\wedge S^{1})}
		$$
	Hence using \cite[proposition 3.1.5]{MR1944041} it follows that the composition above
	is a $C_{eff}^{q}$-colocal equivalence.
	
	(\ref{cor.3.3.Suspension=>qconn-Quillen-equiv.b}):  This follows immediately
	from proposition \ref{prop.3.2.IQTJ-fibrant-replacement-all-Rq} and 
	lemma \ref{lem.3.2.Cq-colocalequivs-stable-S1desuspension}.
\end{proof}

\begin{rmk}
		\label{rmk.3.3.smashT-not-descending}
	The adjunction $(\Tsuspfunctor ,\Tloops ,\varphi)$ is
	a Quillen equivalence on $\motivicTspectra$.
	However it does not descend even to a Quillen
	adjunction on the ($q-1$)-connected motivic stable
	model category
	$\qconnectedTspectra$.
\end{rmk}

\begin{cor}
		\label{cor.3.3.qconn=>triangcat}
	For every $q\in \mathbb Z$, $\qconnectedstablehomotopy$
	has the structure of a triangulated category.
\end{cor}
\begin{proof}
	Theorem \ref{thm.3.2.connective-model-structure} implies in particular
	that $\qconnectedTspectra$ is a pointed simplicial model category,
	and corollary \ref{cor.3.3.Suspension=>qconn-Quillen-equiv} implies that
	the adjunction 
		$$(-\wedge S^{1},\Omega _{S^{1}},\varphi):\qconnectedTspectra \rightarrow \qconnectedTspectra$$
	is a Quillen equivalence.  Therefore
	the result follows from the work of Quillen in 
	\cite[sections I.2 and I.3]{MR0223432} and the work of 
	Hovey in \cite[chapters VI and VII]{MR1650134}.
\end{proof}

\begin{prop}
		\label{prop.3.3.cofibrant-replacement=>triangulatedfunctor}
	We have the following adjunction
		$$\xymatrix{(C_{q}, IQ_{T}J, \varphi) :\qconnectedstablehomotopy \ar[r]& \stablehomotopy}
		$$
	between exact functors of triangulated categories.
\end{prop}
\begin{proof}
	Since $\qconnectedTspectra$ is the right Bousfield localization of
	$\motivicTspectra$ with respect to the $C^{q}_{eff}$-colocal equivalences, we have that the
	identity functor $id:\qconnectedTspectra \rightarrow \motivicTspectra$
	is a left Quillen functor.  Therefore we get
	the following adjunction at the level of the associated homotopy categories:
		$$\xymatrix{(C_{q}, IQ_{T}J, \varphi ):\qconnectedstablehomotopy \ar[r]& \stablehomotopy}
		$$
	
	Now proposition 6.4.1 in \cite{MR1650134} implies that
	$C_{q}$ maps cofibre sequences in $\qconnectedstablehomotopy$ to cofibre sequences in
	$\stablehomotopy$.
	Therefore 
	using proposition 7.1.12 in \cite{MR1650134} we have
	that $C_{q}$ and $IQ_{T}J$ are both exact functors between triangulated categories.
\end{proof}

\begin{prop}
		\label{prop.3.2.Cq=>full-embedding}
	Fix $q\in \mathbb Z$.  Then the unit of the adjunction
		$$\xymatrix{(C_{q}, IQ_{T}J, \varphi ):\qconnectedstablehomotopy \ar[r]& \stablehomotopy}
		$$
	$\sigma _{X}:X\rightarrow IQ_{T}JC_{q}X$ is an
	isomorphism in $\qconnectedstablehomotopy$ for every $T$-spectrum $X$, and
	the functor:
		$$C_{q}:\qconnectedstablehomotopy \rightarrow \stablehomotopy
		$$
	is a full embedding.
\end{prop}
\begin{proof}
	For any $T$-spectrum $X$, we have the following commutative diagram in $\qconnectedTspectra$:
		$$\xymatrix{C_{q}X \ar[rr]^-{C_{q}^{X}} \ar[d]_-{IQ_{T}J^{C_{q}X}}&& X \ar[d]^-{IQ_{T}J^{X}}\\
								IQ_{T}JC_{q}X \ar[rr]_-{IQ_{T}J(C_{q}^{X})}&& IQ_{T}JX}
		$$
	where $IQ_{T}J^{C_{q}X}$ is in particular a weak equivalence in $\motivicTspectra$. 
	But since $\qconnectedTspectra$ is the right Bousfield localization
	of $\motivicTspectra$ with respect to the $C_{eff}^{q}$-colocal equivalences, 
	proposition 3.1.5 in \cite{MR1944041}
	implies that $IQ_{T}J^{C_{q}X}$ is also a
	$C_{eff}^{q}$-colocal equivalence.
	
	On the other hand, by construction we have that 
	$C_{q}^{X}$ is a $C_{eff}^{q}$-colocal equivalence.  Therefore
	$IQ_{T}J^{C_{q}X}$ and $C_{q}^{X}$ both
	become isomorphisms in
	$\qconnectedstablehomotopy$.  
	
	Finally, since $\sigma _{X}$ is the following
	composition in $\qconnectedstablehomotopy$:
		$$\xymatrix{X \ar[r]^{(C_{q}^{X})^{-1}}_-{\cong} &  
								C_{q}X \ar[rr]^-{IQ_{T}J^{C_{q}X}}_-{\cong} && IQ_{T}JC_{q}X}
		$$
	it follows that $\sigma _{X}$ is an isomorphism in $\qconnectedstablehomotopy$
	as we wanted.  This also implies that the functor 
		$$C_{q}:\qconnectedstablehomotopy \rightarrow \stablehomotopy
		$$
	is a full embedding.
\end{proof}

\begin{prop}
		\label{prop.3.2.detecting-isos-RCqeff}
	Fix $q\in \mathbb Z$, and let $f:X\rightarrow Y$ be a map in $\qconnectedstablehomotopy$.
	Then $f$ is an isomorphism if and only if 
	the following condition holds:
			
	For every $\generatorNRS \in C_{eff}^{q}$, the induced maps
		$$\xymatrix{[\generatorNRS , C_{q}X]_{Spt} \ar[r]^-{(C_{q}f)_{\ast}}& 
													[\generatorNRS , C_{q}Y]_{Spt}}
		$$
	are all isomorphisms of abelian groups.
\end{prop}
\begin{proof}
	Complete the map $f$ to a distinguished triangle in $\qconnectedstablehomotopy$:
		$$\xymatrix{X \ar[r]^-{f}& Y \ar[r]& Z\ar[r] & \Tsuspfunctor ^{1,0}X}
		$$
	
	We have that
		\begin{equation}
				\label{diagram.prop.3.2.detecting-isos-RCqeff}
			\xymatrix{C_{q}X \ar[r]^-{C_{q}f}& C_{q}Y \ar[r]& C_{q}Z\ar[r] & \Tsuspfunctor ^{1,0}C_{q}X}
		\end{equation}
	is also a distinguished triangle in $\qconnectedstablehomotopy$, therefore $C_{q}f$ is
	an isomorphism in $\qconnectedstablehomotopy$ if and only if $C_{q}Z\cong \ast$ in
	$\qconnectedstablehomotopy$; and since $C_{q}$ is a cofibrant replacement functor
	in $\qconnectedTspectra$ we have that $f$ is an isomorphism in $\qconnectedstablehomotopy$
	if and only if $C_{q}f$ is an isomorphism in $\qconnectedstablehomotopy$.
	
	Hence, $f$ is an isomorphism in $\qconnectedstablehomotopy$ if and only if
	$C_{q}Z\cong \ast$ in $\qconnectedstablehomotopy$.
	Now corollary \ref{cor.3.2.detecting-trivialspectra-Cqeff} implies that
	$C_{q}Z\cong \ast$ in $\qconnectedstablehomotopy$ if and only if for every
	$\generatorNRS \in C_{eff}^{q}$:
		$$[\generatorNRS , C_{q}Z]_{Spt}\cong 0
		$$
	
	But proposition \ref{prop.3.3.cofibrant-replacement=>triangulatedfunctor} 
	implies that the diagram (\ref{diagram.prop.3.2.detecting-isos-RCqeff}) is
	a distinguished triangle in $\stablehomotopy$; and since for every $\generatorNRS \in C_{eff}^{q}$,
	the functor $[\generatorNRS ,-]_{Spt}$ is homological, we get the following long exact
	sequence of abelian groups
		$$\xymatrix{\vdots \ar[d] &\\
								[\generatorNRS ,C_{q}X]_{Spt} \ar[d]^-{(C_{q}f)_{\ast}} & \\ 
								[\generatorNRS ,C_{q}Y]_{Spt} \ar[d] & \\ 
								[\generatorNRS ,C_{q}Z]_{Spt} \ar[d] & \\
								[\generatorNRS ,\Sigma _{T}^{1,0}C_{q}X]_{Spt} 
								\ar[d]^{(\Sigma _{T}^{1,0}C_{q}f)_{\ast}} & \ar[l]_-{\cong}^-{\Tsuspfunctor ^{1,0}} 
								[F_{n+1}(S^{r}\wedge \gm ^{s+1}\wedge U_{+}),C_{q}X]_{Spt} \ar[d]^-{(C_{q}f)_{\ast}}\\
								[\generatorNRS ,\Sigma _{T}^{1,0}C_{q}Y]_{Spt} 
								\ar[d] & [F_{n+1}(S^{r}\wedge \gm ^{s+1}\wedge U_{+}),C_{q}Y]_{Spt} 
								\ar[l]_-{\cong}^-{\Tsuspfunctor ^{1,0}}\\
								\vdots}
		$$
	Therefore $[\generatorNRS ,C_{q}Z]_{Spt}\cong 0$ for every $\generatorNRS \in C_{eff}^{q}$ if
	and only if the induced map
		$$\xymatrix{[\generatorNRS ,C_{q}X]_{Spt} \ar[r]^-{(C_{q}f)_{\ast}}& [\generatorNRS ,C_{q}Y]_{Spt}}
		$$
	is an isomorphism of abelian groups for every $\generatorNRS \in C_{eff}^{q}$.
	This finishes the proof.
\end{proof}

\begin{prop}
		\label{prop.3.2.generating-Cqeffcolocalspectra-1}
	Fix $q\in \mathbb Z$, and let $A$ be an arbitrary $T$-spectrum in $\stablehomotopyeffq$.
	Then $(Q_{s}A)\wedge S^{1}$ is a $C_{eff}^{q}$-colocal $T$-spectrum
	in $\motivicTspectra$.
\end{prop}
\begin{proof}
	Let $\omega _{0}$, $\eta _{0}$ denote the base points corresponding to
	$Map_{\ast}(Q_{s}A,IQ_{T}JX)$ and $Map_{\ast}(Q_{s}A,IQ_{T}JY)$ respectively.
	
	It is clear that $Q_{s}A\cong A$ in
	$\stablehomotopy$; then 
	$Q_{s}A$ is in $\stablehomotopyeffq$, since $A$ is in $\stablehomotopyeffq$ and
	$\stablehomotopyeffq$ is a triangulated subcategory of $\stablehomotopy$.
	
	Since $\motivicTspectra$ is a simplicial model category, we have that $Q_{s}A\wedge S^{1}$ is
	cofibrant in $\motivicTspectra$, hence it suffices to check that for every $C_{eff}^{q}$-colocal equivalence
	$f:X\rightarrow Y$, the induced map
		$$\xymatrix{Map(Q_{s}A\wedge S^{1},IQ_{T}JX) \ar[rr]^-{(IQ_{T}Jf)_{\ast}}&& 
								Map(Q_{s}A\wedge S^{1},IQ_{T}JY)}
		$$
	is a weak equivalence of simplicial sets.  
	
	Now
	corollary \ref{cor.3.3.Suspension=>qconn-Quillen-equiv} together with
	proposition \ref{prop.3.2.IQTJ-fibrant-replacement-all-Rq} imply that
	for every $n\geq 0$, $\Omega_{S^{n}}IQ_{T}J(f)$ is also a $C_{eff}^{q}$-colocal equivalence.
	Hence corollary \ref{cor.3.2.detecting-Cqeff-colocalequivs2} implies that
	$r_{q}\Omega _{S^{n}}IQ_{T}J(f)$ is an isomorphism in $\stablehomotopyeffq$.  Since
	$Q_{s}A\in \stablehomotopyeffq$, we have that $i_{q}Q_{s}A=Q_{s}A$,
	then by
	proposition \ref{prop.3.1.adjunctions-stableeffq} we get the following commutative
	diagram where both rows and the left vertical map are isomorphisms of abelian groups:
		$$\xymatrix{\Hom_{\stablehomotopyeffq}(Q_{s}A,r_{q}\Omega_{S^{n}}IQ_{T}JX) \ar[r]^-{\cong} 
								\ar[d]_-{(r_{q}\Omega_{S^{n}}IQ_{T}Jf)_{\ast}}^-{\cong}& 
								[Q_{s}A,\Omega _{S^{n}}IQ_{T}JX]_{Spt} \ar[d]^-{(\Omega_{S^{n}}IQ_{T}Jf)_{\ast}}\\
								\Hom_{\stablehomotopyeffq}(Q_{s}A,r_{q}\Omega_{S^{n}}IQ_{T}JY)\ar[r]^-{\cong}&
								[Q_{s}A,\Omega _{S^{n}}IQ_{T}JY]_{Spt}}
		$$
	Therefore the right vertical map is also an isomorphism of abelian groups.
	
	Now since $\motivicTspectra$ is a pointed simplicial model category,
	we have that lemma 6.1.2 in \cite{MR1650134} together 
	with remark \ref{rmk.2.4.simp-simppointed-structures-Tspectra}(\ref{rmk.2.4.simp-simppointed-structures-Tspectra.b})
	imply that
	the following diagram is commutative for $n\geq 0$, where all the vertical maps together with
	the bottom row are isomorphisms of abelian groups:
		$$\xymatrix{\pi_{n, \omega _{0}}Map(Q_{s}A ,IQ_{T}JX) \ar[rr]^-{(IQ_{T}Jf)_{\ast}}\ar@{=}[d]&&  
								\pi_{n, \eta _{0}}Map(Q_{s}A ,IQ_{T}JY)\ar@{=}[d]\\
								\pi_{n, \omega _{0}}Map_{\ast}(Q_{s}A ,IQ_{T}JX) \ar[rr]^-{(IQ_{T}Jf)_{\ast}}\ar[d]_-{\cong}&&  
								\pi_{n, \eta _{0}}Map_{\ast}(Q_{s}A ,IQ_{T}JY)\ar[d]^-{\cong}\\
								[Q_{s}A , \Omega _{S^{n}}IQ_{T}JX]_{Spt} \ar[rr]_-{(\Omega _{S^{n}}IQ_{T}Jf)_{\ast}}^-{\cong}&& 
								[Q_{s}A , \Omega _{S^{n}}IQ_{T}JY]_{Spt}}
		$$
	Therefore all the maps in the top row are also isomorphisms.  Thus,
	the induced map
		$$\xymatrix{Map(Q_{s}A ,IQ_{T}JX) \ar[rr]^-{(IQ_{T}Jf)_{\ast}}&& 
								Map(Q_{s}A ,IQ_{T}JY)}
		$$ 
	is a weak equivalence when it is restricted to the path component of $Map(Q_{s}A ,IQ_{T}JX)$
	containing $\omega _{0}$.  This implies that the following induced map
		$$\xymatrix{Map_{\ast}(S^{1},Map_{\ast}(Q_{s}A ,IQ_{T}JX)) 
								\ar[d]_-{(IQ_{T}Jf)_{\ast}}\\
								Map_{\ast}(S^{1},Map_{\ast}(Q_{s}A ,IQ_{T}JY))}
		$$
	is a weak equivalence since taking $S^{1}$-loops kills the path components that do not
	contain the base point.  
	
	Finally, since $\motivicTspectra$ is a simplicial model category we have that the rows in the following
	commutative diagram are isomorphisms:
		$$\xymatrix{Map_{\ast}(S^{1},Map_{\ast}(Q_{s}A ,IQ_{T}JX)) \ar[r]^-{\cong} 
								\ar[d]_-{(IQ_{T}Jf)_{\ast}}&
								Map_{\ast}(Q_{s}A \wedge S^{1},IQ_{T}JX)
								 \ar[d]^-{(IQ_{T}Jf)_{\ast}}\\
								Map_{\ast}(S^{1},Map_{\ast}(Q_{s}A ,IQ_{T}JY)) \ar[r]_-{\cong}& 
								Map_{\ast}(Q_{s}A \wedge S^{1},IQ_{T}JY)}
		$$
	Hence the two out of three property for weak equivalences implies that the right vertical map
	is a weak equivalence of simplicial sets, as we wanted.
\end{proof}

\begin{cor}
		\label{cor.3.2.generating-Cqeff-colocalTspectra2}
	Fix $q\in \mathbb Z$ and let $X$ be an arbitrary $T$-spectrum in $\stablehomotopyeffq$.
	Then $Q_{s}X$ is a $C_{eff}^{q}$-colocal $T$-spectrum in $\motivicTspectra$
\end{cor}
\begin{proof}
	Let $R$ denote a fibrant replacement functor in $\motivicTspectra$ such that
	for every $T$-spectrum $Y$, the natural map
		$$\xymatrix{Y\ar[r]^-{R_{Y}}& RY}
		$$
	is a trivial cofibration in $\motivicTspectra$.  Then
	$RQ_{s}X$ is cofibrant in $\motivicTspectra$.  
	Now the map
		$$\xymatrix{Q_{s}X \ar[r]^-{R_{Q_{s}X}}& RQ_{s}X}
		$$
	is in particular a weak equivalence in $\motivicTspectra$, therefore using \cite[lemma 3.2.1(2)]{MR1944041}
	we get that $Q_{s}X$ is $C_{eff}^{q}$-colocal if and only if $RQ_{s}X$
	is $C_{eff}^{q}$-colocal.  We will show that $RQ_{s}X$ is
	$C_{eff}^{q}$-colocal.
	
	By hypothesis $X$ is in $\stablehomotopyeffq$ and
	it is clear that $Q_{s}X\cong X$ in $\stablehomotopy$.  Hence
	$Q_{s}X$ is also in $\stablehomotopyeffq$ since it is a triangulated
	subcategory of $\stablehomotopy$.
	Therefore $\Omega _{S^{1}}RQ_{s}X$ is
	also in $\stablehomotopyeffq$ since
	$\Omega _{S^{1}}RQ_{s}X$ computes the
	desuspension $\Tsuspfunctor ^{-1,0}Q_{s}X$ of $Q_{s}X$.
	
	Using proposition \ref{prop.3.2.generating-Cqeffcolocalspectra-1} 
	we have that $(Q_{s}\Omega _{S^{1}}RQ_{s}X)\wedge S^{1}$
	is $C_{eff}^{q}$-colocal.  But since the adjunction
		$$\xymatrix{(-\wedge S^{1}, \Omega _{S^{1}},\varphi):\motivicTspectra \ar[r]& \motivicTspectra}
		$$
	is a Quillen equivalence, we have the following weak equivalence in
	$\motivicTspectra$:
		$$\xymatrix{\ar[rrrr]^-{\epsilon _{RQ_{s}X}\circ (Q_{s}^{\Omega _{S^{1}}RQ_{s}X}\wedge id)} 
								(Q_{s}\Omega _{S^{1}}RQ_{s}X) \wedge S^{1} &&&&
								RQ_{s}X  }
		$$
	where $\epsilon$ denotes the counit of the adjunction considered above.
	
	Finally using \cite[lemma 3.2.1(2)]{MR1944041} again,
	we get that $RQ_{s}X$ is $C_{eff}^{q}$-colocal.  This finishes the proof.
\end{proof}

\begin{prop}
		\label{prop.3.2.counit=>replaces-fq}
	Fix $q\in \mathbb Z$ and let $\rho$ be the counit of the adjunction:
		$$\xymatrix{(C_{q}, IQ_{T}J, \varphi ):\qconnectedstablehomotopy \ar[r]& \stablehomotopy}
		$$
	Then for every $T$-spectrum $X$, the map
		$$r_{q}(\rho _{X}):r_{q}C_{q}IQ_{T}JX\rightarrow r_{q}X
		$$
	is an isomorphism in $\stablehomotopyeffq$; and this map induces a
	natural isomorphism 
	between the following exact functors 
		$$\xymatrix{\stablehomotopy \ar@<1ex>[rr]^-{r_{q}C_{q}IQ_{T}J} 
								\ar@<-1ex>[rr]_-{r_{q}} && 
								\stablehomotopyeffq}
		$$
\end{prop}
\begin{proof}
	The naturality of the counit $\rho$, implies that 
	$r_{q}(\rho _{-}):r_{q}C_{q}IQ_{T}J\rightarrow r_{q}$ is a natural transformation.
	Hence, it is enough to show that for every $T$-spectrum $X$,
	$r_{q}(\rho _{X})$ is an isomorphism in $\stablehomotopyeffq$.
	
	Consider the following diagram of $T$-spectra:
		$$\xymatrix{C_{q}IQ_{T}JX \ar[rr]^-{C_{q}^{IQ_{T}JX}}&& IQ_{T}JX && X \ar[ll]_-{IQ_{T}J^{X}}}
		$$
	where $IQ_{T}J^{X}$ is a weak equivalence in $\motivicTspectra$ and
	$C_{q}^{IQ_{T}JX}$ is a $C_{eff}^{q}$-colocal equivalence.
	Then it is clear that $r_{q}(IQ_{T}J^{X})$ is an isomorphism in $\stablehomotopyeffq$.
	On the other hand, corollary \ref{cor.3.2.detecting-Cqeff-colocalequivs2}
	implies that $r_{q}(C_{q}^{IQ_{T}JX})$ is also an isomorphism in $\stablehomotopyeffq$.
	
	And this proves the result, since $\rho _{X}$ is just
	the following composition in $\stablehomotopy$:
		$$\xymatrix{C_{q}IQ_{T}JX \ar[rr]^-{C_{q}^{IQ_{T}JX}}&& IQ_{T}JX \ar[rr]^-{(IQ_{T}J^{X})^{-1}}&& X}
		$$
\end{proof}

\begin{prop}
		\label{prop.3.2.counit-stablehomotopyeffq=>iso-qconnected.b}
	Fix $q\in \mathbb Z$ and let $\theta$ be the counit of the adjunction
		$$\xymatrix{(i_{q},r_{q},\varphi):\stablehomotopyeffq \ar[r]& \stablehomotopyeff}
		$$
	Then for any $T$-spectrum $X$ in $\motivicTspectra$,
	the map 
		$$\xymatrix{IQ_{T}J(\theta _{X}):IQ_{T}J(i_{q}r_{q})X \ar[r]&
								IQ_{T}JX}
		$$
	is an isomorphism in $\qconnectedstablehomotopy$;
	and this map induces a
	natural isomorphism 
	between the following exact functors
		$$\xymatrix{\stablehomotopy \ar@<1ex>[rr]^-{IQ_{T}J(i_{q}r_{q})} 
								\ar@<-1ex>[rr]_-{IQ_{T}J} && 
								\qconnectedstablehomotopy}
		$$
\end{prop}
\begin{proof}
	The naturality of the counit $\theta$, implies that 
	$IQ_{T}J(\theta _{-}):IQ_{T}J(i_{q}r_{q})\rightarrow IQ_{T}J$ is a natural transformation.
	Hence, it is enough to show that for every $T$-spectrum $X$,
	$IQ_{T}J(\theta _{X})$ is an isomorphism in $\qconnectedstablehomotopy$.

	By proposition \ref{prop.3.2.detecting-isos-RCqeff} it is enough to show that
	for every $\generatorNRS \in C_{eff}^{q}$ the induced maps
		$$\xymatrix{[\generatorNRS ,C_{q}IQ_{T}J(i_{q}r_{q})X]_{Spt} 
								\ar[d]^-{C_{q}IQ_{T}J(\theta _{X})_{\ast}}\\
								[\generatorNRS ,C_{q}IQ_{T}JX]_{Spt}}
		$$
	are all isomorphisms of abelian groups.
	
	Consider the following commutative diagram in $\stablehomotopy$:
		$$\xymatrix{C_{q}IQ_{T}J(i_{q}r_{q})X \ar[rr]^-{C_{q}IQ_{T}J(\theta _{X})} 
								\ar[d]_-{C_{q}^{IQ_{T}Ji_{q}r_{q}X}}&& C_{q}IQ_{T}JX \ar[d]^-{C_{q}^{IQ_{T}JX}}\\
								IQ_{T}J(i_{q}r_{q})X \ar[rr]^-{IQ_{T}J(\theta _{X})}&& IQ_{T}JX}
		$$
	where $C_{q}^{IQ_{T}Ji_{q}r_{q}X}$ and $C_{q}^{IQ_{T}JX}$ are by construction
	maps of $T$-spectra and $C_{eff}^{q}$-colocal equivalences.  Therefore proposition
	\ref{prop.3.2.classif-Cqcolocalequivs} implies that for every $\generatorNRS \in C_{eff}^{q}$ the induced maps
		$$\xymatrix{[\generatorNRS ,C_{q}IQ_{T}J(i_{q}r_{q})X]_{Spt} \ar[d]_-{(C_{q}^{IQ_{T}Ji_{q}r_{q}X})_{\ast}}&
								[\generatorNRS ,C_{q}IQ_{T}JX]_{Spt} \ar[d]^-{(C_{q}^{IQ_{T}JX})_{\ast}} \\
								[\generatorNRS ,IQ_{T}J(i_{q}r_{q})X]_{Spt} & [\generatorNRS ,IQ_{T}JX]_{Spt}}
		$$
	are both isomorphisms of abelian groups.
	
	On the other hand, 
	proposition \ref{prop.3.1.counit-properties} implies that we have
	an induced isomorphism of abelian groups:
		$$\xymatrix{[\generatorNRS ,IQ_{T}J(i_{q}r_{q})X]_{Spt} 
								\ar[d]^-{IQ_{T}J(\theta _{X})_{\ast}}\\
								[\generatorNRS ,IQ_{T}JX]_{Spt}}
		$$
	for every $\generatorNRS \in C_{eff}^{q}$.
	
	Finally, this implies that for every $\generatorNRS \in C_{eff}^{q}$, we
	get the following induced isomorphisms of abelian groups
		$$\xymatrix{[\generatorNRS ,C_{q}IQ_{T}J(i_{q}r_{q})X]_{Spt} 
								\ar[d]^-{C_{q}IQ_{T}J(\theta _{X})_{\ast}}\\
								[\generatorNRS ,C_{q}IQ_{T}JX]_{Spt}}
		$$
	as we wanted.
\end{proof}

\begin{prop}
		\label{prop.3.2.counit-stablehomotopyeffq=>iso-qconnected.c}
	Fix $q\in \mathbb Z$, and let $\theta$ be the counit of the adjunction
		$$\xymatrix{(i_{q},r_{q},\varphi):\stablehomotopyeffq \ar[r]& \stablehomotopyeff}
		$$
	Then for any $T$-spectrum $X$ ,
	the map 
		$$\xymatrix{C_{q}IQ_{T}J(\theta _{X}):C_{q}IQ_{T}J(i_{q}r_{q})X \ar[r]&
								C_{q}IQ_{T}JX}
		$$
	is an isomorphism in $\stablehomotopy$;
	and this map induces a
	natural isomorphism 
	between the following exact functors
		$$\xymatrix{\stablehomotopy \ar@<1ex>[rr]^-{C_{q}IQ_{T}J(i_{q}r_{q})} 
								\ar@<-1ex>[rr]_-{C_{q}IQ_{T}J} && 
								\stablehomotopy}
		$$
\end{prop}
\begin{proof}
	The naturality of the counit $\theta$, implies that 
	$C_{q}IQ_{T}J(\theta _{-}):C_{q}IQ_{T}J(i_{q}r_{q})\rightarrow C_{q}IQ_{T}J$ is a natural transformation.
	Hence, it is enough to show that for every $T$-spectrum $X$,
	$C_{q}IQ_{T}J(\theta _{X})$ is an isomorphism in $\stablehomotopy$.
	
	But this follows immediately from proposition \ref{prop.3.2.counit-stablehomotopyeffq=>iso-qconnected.b} 
	together with proposition \ref{prop.3.3.cofibrant-replacement=>triangulatedfunctor}.
\end{proof}

\begin{prop}
		\label{prop.3.2.CqIQtJiqrq-->IQtJiqrq==iso-in-stablehomotopy}
	Fix $q\in \mathbb Z$.  Then for every $T$-spectrum $X$, the natural map
		$$\xymatrix{C_{q}IQ_{T}J(i_{q}r_{q})X \ar[rrr]^-{C_{q}^{IQ_{T}J(i_{q}r_{q})X}}&&& IQ_{T}J(i_{q}r_{q})X}
		$$
	is a weak equivalence in $\motivicTspectra$.
	Therefore we have a natural isomorphism 
	between the following exact functors
		$$\xymatrix{\stablehomotopy \ar@<1ex>[rrr]^-{C_{q}IQ_{T}J(i_{q}r_{q})} 
								\ar@<-1ex>[rrr]_-{IQ_{T}J(i_{q}r_{q})}&&& \stablehomotopy}
		$$
\end{prop}
\begin{proof}
	The naturality of the maps
	$C_{q}^{X}:C_{q}X\rightarrow X$ implies that we have an induced natural transformation
	of functors $C_{q}IQ_{T}J(i_{q}r_{q})\rightarrow IQ_{T}J(i_{q}r_{q})$.  Hence,
	it is enough to show that for every $T$-spectrum $X$,
	$C_{q}^{IQ_{T}J(i_{q}r_{q})X}$ is a weak equivalence in $\motivicTspectra$.
	
	Consider the following commutative diagram in $\motivicTspectra$:
		$$\xymatrix{Q_{s}C_{q}IQ_{T}J(i_{q}r_{q})X \ar[rrr]^-{Q_{s}(C_{q}^{IQ_{T}J(i_{q}r_{q})X})} 
								\ar[d]_-{Q_{s}^{C_{q}IQ_{T}J(i_{q}r_{q})X}}&&&
								Q_{s}IQ_{T}J(i_{q}r_{q})X \ar[d]^-{Q_{s}^{IQ_{T}J(i_{q}r_{q})X}}\\
								C_{q}IQ_{T}J(i_{q}r_{q})X \ar[rrr]_-{C_{q}^{IQ_{T}J(i_{q}r_{q})X}}&&& IQ_{T}J(i_{q}r_{q})X}
		$$
	Since $Q_{s}$ is a cofibrant replacement functor in $\motivicTspectra$, it follows
	that the vertical maps are weak equivalences in $\motivicTspectra$.  Hence by the two
	out of three property for weak equivalences it suffices to show that
	$Q_{s}(C_{q}^{IQ_{T}J(i_{q}r_{q})X})$ is a weak equivalence in $\motivicTspectra$.
	
	On the other hand we have that by construction $C_{q}^{IQ_{T}J(i_{q}r_{q})X}$ is
	a $C_{eff}^{q}$-colocal equivalence, and \cite[proposition 3.1.5]{MR1944041} implies that
	the vertical maps in the diagram above are also $C_{eff}^{q}$-colocal equivalences.
	Then by the two out of three property for $C_{eff}^{q}$-colocal equivalences we have that
	$Q_{s}(C_{q}^{IQ_{T}J(i_{q}r_{q})X})$ is a $C_{eff}^{q}$-colocal equivalence.
	
	Now by construction we have that $C_{q}IQ_{T}J(i_{q}r_{q})X$ is a $C_{eff}^{q}$-colocal
	$T$-spectrum, and that $Q_{s}C_{q}IQ_{T}J(i_{q}r_{q})X$ is cofibrant in $\motivicTspectra$.  Since
	$Q_{s}^{C_{q}IQ_{T}J(i_{q}r_{q})X}$ is in particular a weak
	equivalence in $\motivicTspectra$, using \cite[lemma 3.2.1(2)]{MR1944041} we have that
	$Q_{s}C_{q}IQ_{T}J(i_{q}r_{q})X$ is also a $C_{eff}^{q}$-colocal $T$-spectrum.
	
	It is clear that $IQ_{T}J(i_{q}r_{q})X\cong i_{q}r_{q}X$ in $\stablehomotopy$, therefore
	$IQ_{T}J(i_{q}r_{q})X$ is in $\stablehomotopyeffq$ since
	$\stablehomotopyeffq$ is a triangulated subcategory of $\stablehomotopy$ and $i_{q}r_{q}X$
	is in $\stablehomotopyeffq$.  Then using corollary \ref{cor.3.2.generating-Cqeff-colocalTspectra2},
	we have that $Q_{s}IQ_{T}J(i_{q}r_{q})X$ is a $C_{eff}^{q}$-colocal $T$-spectrum.
	
	Finally we have that $Q_{s}(C_{q}^{IQ_{T}J(i_{q}r_{q})X})$ is a $C_{eff}^{q}$-colocal
	equivalence, and that $Q_{s}C_{q}IQ_{T}J(i_{q}r_{q})X$, $Q_{s}IQ_{T}J(i_{q}r_{q})X$
	are both $C_{eff}^{q}$-colocal $T$-spectra.  Then
	\cite[theorem 3.2.13(2)]{MR1944041} implies that $Q_{s}(C_{q}^{IQ_{T}J(i_{q}r_{q})X})$ is
	a weak equivalence in $\motivicTspectra$, as we wanted.
\end{proof}

\begin{thm}
		\label{thm.3.2.Rq-models-f<q}
	Fix $q\in \mathbb Z$.  Then for every $T$-spectrum $X$, we have the following diagram in 
	$\stablehomotopy$:
		\begin{equation}
					\label{diagram.3.2.fq-lifting}
			\begin{array}{c}
				\xymatrix{f_{q}X=i_{q}r_{q}X \ar[rr]^-{IQ_{T}J^{f_{q}X}}_-{\cong}&& IQ_{T}Jf_{q}X  
									&&& \\ &&&&&\\ 
									&& C_{q}IQ_{T}Jf_{q}X \ar[rrr]_-{C_{q}IQ_{T}J(\theta _{X})}^-{\cong}
									\ar[uu]_-{C_{q}^{IQ_{T}Jf_{q}X}}^-{\cong}
									&&& C_{q}IQ_{T}JX }
			\end{array}
		\end{equation}
	where all the maps are isomorphisms in $\stablehomotopy$.
	This diagram induces a
	natural isomorphism  
	between the following exact functors:
		$$\xymatrix{\stablehomotopy \ar@<1ex>[rr]^-{f_{q}}  
								\ar@<-1ex>[rr]_-{C_{q}IQ_{T}J} && \stablehomotopy}
		$$
\end{thm}
\begin{proof}
	Since
	$IQ_{T}J$ is a fibrant replacement functor in $\motivicTspectra$, 
	it is clear that $IQ_{T}J^{f_{q}X}$ becomes
	an isomorphism in the associated homotopy category $\stablehomotopy$.
	
	The fact that $C_{q}^{IQ_{T}Jf_{q}X}$
	is an isomorphism in $\stablehomotopy$ follows from
	proposition \ref{prop.3.2.CqIQtJiqrq-->IQtJiqrq==iso-in-stablehomotopy}.  
	Finally, proposition \ref{prop.3.2.counit-stablehomotopyeffq=>iso-qconnected.c}
  implies that $C_{q}IQ_{T}J(\theta _{X})$ is also an isomorphism in $\stablehomotopy$.
	This shows that all the maps in the diagram (\ref{diagram.3.2.fq-lifting}) are isomorphisms
	in $\stablehomotopy$, therefore for every $T$-spectrum $X$ we can define the
	following composition in $\stablehomotopy$
		\begin{equation}
					\label{diagram.3.2.fq-lifting.b}
			\begin{array}{c}
				\xymatrix{f_{q}X\ar[rr]^-{IQ_{T}J^{f_{q}X}}_-{\cong}&& IQ_{T}Jf_{q}X
									\ar[dd]^-{(C_{q}^{IQ_{T}Jf_{q}X})^{-1}}_-{\cong}
									&&& \\ &&&&&\\ 
									&& C_{q}IQ_{T}Jf_{q}X \ar[rrr]_-{C_{q}IQ_{T}J(\theta _{X})}^-{\cong}
									&&& C_{q}IQ_{T}JX }
			\end{array}
		\end{equation}
	which is an isomorphism.  The fact that
	$IQ_{T}J$ is a functorial fibrant replacement in $\motivicTspectra$,	
	propositions
	\ref{prop.3.2.CqIQtJiqrq-->IQtJiqrq==iso-in-stablehomotopy} and 
	\ref{prop.3.2.counit-stablehomotopyeffq=>iso-qconnected.c},
	imply all together that the isomorphisms defined in diagram (\ref{diagram.3.2.fq-lifting.b})
	induce a natural isomorphism of functors $f_{q}\stackrel{\cong}{\rightarrow}C_{q}IQ_{T}J$.
	This finishes the proof.
\end{proof}

	Theorem \ref{thm.3.2.Rq-models-f<q} gives the desired lifting to 
	the model category level for the functor $f_{q}$.  Now we proceed to show
	that the homotopy categories $\qconnectedstablehomotopy$ are in fact
	equivalent to the categories $\stablehomotopyeffq$ defined in section \ref{section-slice-filtration}.

	Using
	propositions \ref{prop.3.1.adjunctions-stableeffq} and
	\ref{prop.3.3.cofibrant-replacement=>triangulatedfunctor},
	we get the following diagram of adjunctions:
		\begin{equation}
					\label{diagram3.2.adjunctions-inducing-equiv}
			\begin{array}{c}
				\xymatrix{\qconnectedstablehomotopy \ar[rr]^-{(C_{q}, IQ_{T}J, \varphi )}&& \stablehomotopy \\
									&& \stablehomotopyeffq \ar[u]_-{(i_{q},r_{q},\varphi )} \\
									\qconnectedstablehomotopy \ar@<1ex>[rr]^-{C_{q}} &&
									\stablehomotopy \ar@<1ex>[ll]^-{IQ_{T}J} \ar@<-1ex>[d]_-{r_{q}}\\
									&&\stablehomotopyeffq \ar@<-1ex>[u]_-{i_{q}}}
			\end{array}
		\end{equation}
	where all the functors are exact.

\begin{prop}
		\label{prop.3.2.Rq-lifts-qSH}
	For every $q\in \mathbb Z$ the adjunctions of diagram (\ref{diagram3.2.adjunctions-inducing-equiv})
	induce an equivalence of categories:
		\begin{equation}
				\label{diagram.3.2.Rqconn-equiv-to-stabeffq}
			\xymatrix{\qconnectedstablehomotopy \ar@<1ex>[rr]^-{r_{q}C_{q}}&& 
								\stablehomotopyeffq \ar@<1ex>[ll]^-{IQ_{T}Ji_{q}}}
		\end{equation}
	between $\qconnectedstablehomotopy$ and $\stablehomotopyeffq$.
\end{prop}
\begin{proof}
	It is enough to show the existence of the following natural
	isomorphisms between functors:
		\begin{equation}
					\label{natural-adjunctions-slices}
			\begin{array}{c}
				\xymatrix{id \ar[r]_-{\cong}^-{\epsilon}& (IQ_{T}Ji_{q})(r_{q}C_{q})\\
									(r_{q}C_{q})(IQ_{T}Ji_{q}) \ar[r]^-{\cong}_-{\eta}& id}
			\end{array}
		\end{equation}
	We construct first the natural equivalence $\epsilon$.
	Let $f:X\rightarrow Y$ be a map in $\qconnectedstablehomotopy$.
	Applying the functor $i_{q}r_{q}C_{q}$, we get the following commutative diagram
	in $\stablehomotopy$:
		$$\xymatrix{i_{q}r_{q}C_{q}X \ar[d]^-{i_{q}r_{q}C_{q}f} 
								 \ar[rr]^-{\theta _{C_{q}X}} && C_{q}X  
								\ar[d]_-{C_{q}f} \\
								 i_{q}r_{q}C_{q}Y  
								\ar[rr]_-{\theta _{C_{q}Y}}&& C_{q}Y}
		$$
	where $\theta$ denotes the counit of the adjunction between $i_{q}$ and $r_{q}$.
	Now if we apply the functor $IQ_{T}J$, we have the following commutative diagram
	in $\qconnectedstablehomotopy$:
		$$\xymatrix{IQ_{T}Ji_{q}r_{q}C_{q}X \ar[d]_-{IQ_{T}Ji_{q}r_{q}C_{q}f} 
								 \ar[rrr]^-{IQ_{T}J(\theta _{C_{q}X})}_-{\cong} &&& IQ_{T}JC_{q}X  
								\ar[d]_-{IQ_{T}JC_{q}f} & X \ar[l]_-{\sigma _{X}}^-{\cong} \ar[d]^-{f}\\
								IQ_{T}Ji_{q}r_{q}C_{q}Y  
								\ar[rrr]_-{IQ_{T}J(\theta _{C_{q}Y})}^-{\cong} &&& IQ_{T}JC_{q}Y 
								& Y \ar[l]^-{\sigma _{Y}}_-{\cong}}
		$$
	where $\sigma$ denotes the unit of the adjunction between $C_{q}$ and $IQ_{T}J$.
	But propositions \ref{prop.3.2.counit-stablehomotopyeffq=>iso-qconnected.b}
	and \ref{prop.3.2.Cq=>full-embedding} imply that all the horizontal maps are isomorphisms in
	$\qconnectedstablehomotopy$.
	Now if we define 
		$$\epsilon _{X}=(IQ_{T}J(\theta _{C_{q}X}))^{-1}\circ (\sigma _{X})$$
	we get the natural isomorphism of functors
	$\epsilon : id\rightarrow (IQ_{T}Ji_{q})(r_{q}C_{q})$.
	
	To finish the proof, we proceed to construct the natural equivalence $\eta$.
	Let $f:X\rightarrow Y$ be a map in $\stablehomotopyeffq$.
	Applying the functor $C_{q}IQ_{T}Ji_{q}$, we get the following commutative
	diagram in $\stablehomotopy$:
		$$\xymatrix{C_{q}IQ_{T}Ji_{q}X \ar[rr]^-{\rho _{i_{q\: X}}} \ar[d]_-{C_{q}IQ_{T}Ji_{q}f}
								&& i_{q}X \ar[d]^-{i_{q}f}\\
								C_{q}IQ_{T}Ji_{q}Y \ar[rr]_-{\rho _{i_{q\: Y}}}&& i_{q}Y}
		$$
	where $\rho$ denotes the counit of the adjunction between $C_{q}$ and $IQ_{T}J$.
	Now if we apply the functor $r_{q}$, we have the following commutative
	diagram in $\stablehomotopyeffq$:
		$$\xymatrix{r_{q}C_{q}IQ_{T}Ji_{q}X \ar[rr]^-{r_{q}(\rho _{i_{q\: X}})}_-{\cong} \ar[d]_-{r_{q}C_{q}IQ_{T}Ji_{q}f}
								&& r_{q}i_{q}X \ar[d]_-{r_{q}i_{q}f} & X \ar[d]^-{f}\ar[l]_-{\tau _{X}}^-{\cong}\\
								r_{q}C_{q}IQ_{T}Ji_{q}Y \ar[rr]_-{r_{q}(\rho _{i_{q\: Y}})}^-{\cong}&& r_{q}i_{q}Y & 
								Y\ar[l]^-{\tau _{Y}}_-{\cong}}
		$$
	where $\tau$ denotes the unit of the adjunction between $i_{q}$ and $r_{q}$.
	But proposition \ref{prop.3.2.counit=>replaces-fq}
	and remark \ref{rmk.3.1.unit=iso} imply that all the horizontal maps are isomorphisms in
	$\stablehomotopyeffq$.
	Now if we define 
		$$\eta _{X}=(\tau _{X})^{-1}\circ r_{q}(\rho _{i_{q\: X}})$$
	we get the natural isomorphism of functors
	$\eta : (r_{q}C_{q})(IQ_{T}Ji_{q})\rightarrow id$.
	This finishes the proof.
\end{proof}

\begin{prop}
		\label{prop.3.2.functors-between-Rq}
	Fix $q\in \mathbb Z$.
		\begin{enumerate}
			\item	\label{prop.3.2.functors-between-Rq.a}We have the following
						commutative diagram of left Quillen functors:
							\begin{equation}
										\label{diagram3.2.liftslicefiltration}
								\begin{array}{c}
									\xymatrix{R_{C_{eff}^{q+1}}\motivicTspectra \ar[rr]^-{id} \ar[dr]_-{id}&& 
														\qconnectedTspectra \ar[dl]^-{id}\\
														&\motivicTspectra  &}
								\end{array}
							\end{equation}
			\item	\label{prop.3.2.functors-between-Rq.b}For every $T$-spectrum $X$, the natural map:
							$$\xymatrix{C_{q}C_{q+1}X \ar[rr]^-{C_{q}^{C_{q+1}X}} && C_{q+1}X}
							$$
						is a weak equivalence in $\stablehomotopy$, and it induces a natural equivalence
						$C_{q}^{C_{q+1}-}:C_{q}\circ C_{q+1}\rightarrow C_{q+1}$
						between the following functors:
							$$\xymatrix{\qplusoneconnectedstablehomotopy \ar[rr]^-{C_{q+1}} \ar[dr]_-{C_{q+1}}&& 
													\qconnectedstablehomotopy  \ar[dl]^-{C_{q}}\\
													& \stablehomotopy &}
							$$
			\item	\label{prop.3.2.functors-between-Rq.c}The natural transformation $f_{q+1}X\rightarrow f_{q}X$
						(see theorem \ref{thm.3.1.slicefiltration}(\ref{thm.3.1.slicefiltration.a})) gets canonically identified,
						through the equivalence of categories $r_{q}C_{q}$, $IQ_{T}Ji_{q}$  constructed in proposition \ref{prop.3.2.Rq-lifts-qSH};
						with the following composition in $\stablehomotopy$
							$$\xymatrix{ & C_{q}C_{q+1}IQ_{T}JX \ar[dr]^-{C_{q}(C_{q+1}^{IQ_{T}JX})}&\\
													C_{q+1}IQ_{T}JX \ar[ur]^-{(C_{q}^{C_{q+1}IQ_{T}JX})^{-1}}&& C_{q}IQ_{T}JX}
							$$
						which is induced by the following commutative diagram in $\motivicTspectra$
							\begin{equation}
										\label{diagram3.2.Counit-connected}
								\begin{array}{c}
									\xymatrix{C_{q}C_{q+1}IQ_{T}JX \ar[d]_-{C_{q}^{C_{q+1}IQ_{T}JX}} 
														\ar[rr]^-{C_{q}(C_{q+1}^{IQ_{T}JX})} && C_{q}IQ_{T}JX \ar[d]^-{C_{q}^{IQ_{T}JX}}\\
														C_{q+1}IQ_{T}JX \ar[rr]_-{C_{q+1}^{IQ_{T}JX}} && IQ_{T}JX}
								\end{array}
							\end{equation}
		\end{enumerate}
\end{prop}
\begin{proof}
	(\ref{prop.3.2.functors-between-Rq.a}):  Since $R_{C_{eff}^{q+1}}\motivicTspectra$ and $\qconnectedTspectra$ are both
	right Bousfield localizations of $\motivicTspectra$, by construction
	the identity functor 
		$$\xymatrix@R=0.5pt{R_{C_{eff}^{q+1}}\motivicTspectra \ar[r]^-{id} & \motivicTspectra \\
								\qconnectedTspectra \ar[r]_-{id} & \motivicTspectra }$$ 
	is in both cases
	a left Quillen functor.  To finish the proof,
	it suffices to show
	that the identity functor
		$$id:\qconnectedTspectra \rightarrow R_{C_{eff}^{q+1}}\motivicTspectra$$
	is a right Quillen functor.  Using the universal property of right
	Bousfield localizations (see definition \ref{def-1.rightlocmodcats}), it is enough to check that
	if $f:X\rightarrow Y$ is
	a $C_{eff}^{q}$-colocal equivalence in $\motivicTspectra$ then
	$IQ_{T}J(f)$ is  a
	$C_{eff}^{q+1}$-colocal equivalence.
	But since $IQ_{T}JX$ and $IQ_{T}JY$ are already fibrant in $\motivicTspectra$, we have that
	$IQ_{T}J(f)$ is a $C_{eff}^{q+1}$-colocal equivalence if and only if
	for every $\generatorNRS \in C_{eff}^{q+1}$, the induced map:
		$$\xymatrix{Map(\generatorNRS ,IQ_{T}JX) \ar[d]^-{IQ_{T}J(f)_{\ast}}\\ 
								Map(\generatorNRS ,IQ_{T}JY)}
		$$
	is a weak equivalence of simplicial sets.
	But since $C_{eff}^{q+1}\subseteq C_{eff}^{q}$, and by hypothesis
	$f$ is a $C_{eff}^{q}$-colocal equivalence; we have that all
	the induced maps $IQ_{T}J(f)_{\ast}$ are weak
	equivalences of simplicial sets.  Thus $IQ_{T}J(f)$ is a
	$C_{eff}^{q+1}$-colocal equivalence, as we wanted.
	
	Finally (\ref{prop.3.2.functors-between-Rq.b}) and (\ref{prop.3.2.functors-between-Rq.c}) follow directly from
	proposition \ref{prop.3.2.Rq-lifts-qSH}, theorem \ref{thm.3.2.Rq-models-f<q} together with
	the commutative diagram 
	\eqref{diagram3.2.liftslicefiltration} of left Quillen funtors
	constructed above and \cite[theorem 1.3.7]{MR1650134}.
\end{proof}

\begin{thm}
		\label{thm.3.2.liftingslicefilt-modelcats}
	We have the following commutative diagram
	of left Quillen functors:
		\begin{equation}
					\label{diagram.3.2.liftingVoevodskyslicefilt.a}
			\begin{array}{c}
				\xymatrix{\vdots \ar[d]_-{id} & \\ 
									R_{C_{eff}^{q+1}}\motivicTspectra \ar[d]_-{id} \ar[dr]^-{id} & \\ 
									\qconnectedTspectra \ar[d]_-{id} \ar[r]^-{id} & \motivicTspectra \\
									R_{C_{eff}^{q-1}}\motivicTspectra \ar[d]_-{id} \ar[ur]_-{id} & \\ 
									\vdots & }
			\end{array}
		\end{equation}
	and the associated diagram of homotopy categories:	
		\begin{equation}
					\label{diagram.3.2.liftingVoevodskyslicefilt.b}
			\begin{array}{c}
				\xymatrix{\vdots \ar@<-1ex>[d] &&& \\ 
									R_{C_{eff}^{q+1}}\stablehomotopy \ar@<-1ex>[d]_-{C_{q+1}} \ar@<2ex>[drrr]|-{C_{q+1}} \ar@<-1ex>[u] &&& \\ 
									\qconnectedstablehomotopy \ar@<-1ex>[d]_-{C_{q}} 
									\ar@<1ex>[rrr]|-{C_{q}} \ar@<-1ex>[u]_-{IQ_{T}J} &&& \stablehomotopy 
									\ar@<1ex>[lll]|-{IQ_{T}J} \ar[ulll]|-{IQ_{T}J} \ar@<2ex>[dlll]|-{IQ_{T}J}\\
									R_{C_{eff}^{q-1}}\stablehomotopy \ar@<-1ex>[d] \ar[urrr]|-{C_{q-1}} \ar@<-1ex>[u]_-{IQ_{T}J}&&& \\ 
									\vdots \ar@<-1ex>[u]&& }
			\end{array}
		\end{equation}
	gets canonically identified, through
	the equivalences of categories $r_{q}C_{q}$, $IQ_{T}Ji_{q}$ constructed
	in proposition \ref{prop.3.2.Rq-lifts-qSH}; with Voevodsky's slice filtration:
		\begin{equation}
					\label{diagram.3.2.liftingVoevodskyslicefilt.c}
			\begin{array}{c}
				\xymatrix{\vdots \ar@<-1ex>[d] &&& \\ 
									\Tsuspfunctor ^{q+1}\stablehomotopy \ar@<-1ex>[d]_-{i_{q+1}} \ar@<2ex>[drrr]|-{i_{q+1}} \ar@<-1ex>[u] &&& \\ 
									\stablehomotopyeffq \ar@<-1ex>[d]_-{i_{q}} 
									\ar@<1ex>[rrr]|-{i_{q}} \ar@<-1ex>[u]_-{r_{q}} &&& \stablehomotopy 
									\ar@<1ex>[lll]|-{r_{q}} \ar[ulll]|-{r_{q+1}} \ar@<2ex>[dlll]|-{r_{q-1}}\\
									\Tsuspfunctor ^{q-1}\stablehomotopy \ar@<-1ex>[d] \ar[urrr]|-{i_{q-1}} \ar@<-1ex>[u]_-{r_{q-1}}&&& \\ 
									\vdots \ar@<-1ex>[u]&& }
			\end{array}
		\end{equation}
\end{thm}
\begin{proof}
	Follows immediately from propositions
	\ref{prop.3.2.functors-between-Rq} and \ref{prop.3.2.Rq-lifts-qSH}.
\end{proof}

\begin{rmk}
		\label{rmk.3.2.Rq-not-cellular}
	The drawback of the model structures on $\qconnectedTspectra$ is that
	it is not clear if they are cellular again.  Therefore in order to
	recover a lifting for the slice functors $s_{q}$, we are forced to
	take an indirect approach.
\end{rmk}

	The first step in this new approach will be
	to construct another family of model structures on $\Tspectra$, via
	left Bousfield localization; 
	such that the fibrant replacement functor provides an alternative description
	of the functors $s_{<q}$ defined in theorem \ref{thm.3.1.motivictower}.

\begin{defi}
		\label{def.3.2.simplicial-disks}
	For $r\geq 1$, we define $D^{r}$
	using the following
	pushout
	diagram of simplicial sets:
		$$\xymatrix{S^{r-1} \ar[r]^-{j _{0}} \ar[d]& S^{r-1}\times \Delta ^{1} \ar[d]^-{p}\\
								\ast \ar[r]& D^{r}}
		$$
	where $j _{0}$ is the following composition:
		$$\xymatrix{S^{r-1}\cong S^{r-1}\times \Delta ^{0} \ar[rr]^-{id\times d_{1}}&& S^{r-1}\times \Delta ^{1}}
		$$
	and let $\iota _{1}:S^{r-1}\rightarrow D^{r}$ be the following composition:
		$$\xymatrix{S^{r-1}\cong S^{r-1}\times \Delta ^{0} \ar[rr]^-{id\times
								d_{0}}&& S^{r-1}\times \Delta ^{1}\ar[r]^-{p}& D^{r}}
		$$
\end{defi}

\begin{rmk}
		\label{rmk.3.2.disks-contractible}
	It is clear that the canonical map
	$\ast \rightarrow D^{r}$ is a trivial cofibration
	in the category of pointed simplicial sets.
\end{rmk}

\begin{prop}
		\label{prop.3.2.Disks=>compact}
	For every $r\geq 1,s\geq 0$, and for every
	scheme $U\in \smoothS$; the pointed
	simplicial presheaf on the smooth Nisnevich site
	over $S$
		$$D^{r}\wedge \gm ^{s}\wedge U_{+}
		$$
	has the following properties:
	\begin{enumerate}
		\item	\label{prop.3.2.Disks=>compact.a}it is compact in the sense of Jardine
					(see definition \ref{def.2.3.compactness}).
		\item	\label{prop.3.2.Disks=>compact.b}the canonical map
					$\ast \rightarrow \diskNRS$ is a
					trivial cofibration in $\motivicTspectra$.
		\item	\label{prop.3.2.Disks=>compact.c}the canonical map
					$\diskNRS \rightarrow \ast$ is a
					weak equivalence in $\motivicTspectra$.
	\end{enumerate}
\end{prop}
\begin{proof}
	(\ref{prop.3.2.Disks=>compact.a}):  It is clear from the construction that
	$D^{r}$ has only finitely many non-degenerate simplices.  Therefore the result
	follows from \cite[lemma 2.2]{MR1787949}.
	
	(\ref{prop.3.2.Disks=>compact.b}): Proposition \ref{prop.2.3.cell-prop-pointedmotcat}
	implies that $\pointedmotivic$ is a $\simpsets _{\ast}$-model category; and since 
	$\gm ^{s}\wedge U_{+}$ is cofibrant, we have the following Quillen adjunction:
		$$\xymatrix{\simpsets _{\ast} \ar@<1ex>[rrr]^-{-\wedge \gm ^{s}\wedge U_{+}}&&& 
								\pointedmotivic \ar@<1ex>[lll]^-{Map_{\ast}(\gm ^{s}\wedge U_{+},-)}}
		$$
	But $\ast \rightarrow D^{r}$ is a trivial cofibration of pointed simplicial
	sets, therefore the induced map
		$$\xymatrix{\ast \cong \ast \wedge \gm ^{s}\wedge U_{+} \ar[r]& 
								D^{r}\wedge \gm ^{s}\wedge U_{+}}
		$$
	is a trivial cofibration in $\pointedmotivic$.
	Finally using proposition \ref{prop.2.4.Ev-n=Quillen-stable-functor} we have that
		$$\xymatrix{(F_{n},Ev_{n},\varphi ):\pointedmotivic \ar[r]& \motivicTspectra}
		$$
	is a Quillen adjunction.  Hence the canonical map
		$$\xymatrix{\ast \cong F_{n}(\ast) \ar[r]& \diskNRS}
		$$
	is a trivial cofibration in $\motivicTspectra$, as we wanted.
	
	(\ref{prop.3.2.Disks=>compact.c}):  Follows immediately from (\ref{prop.3.2.Disks=>compact.b})
	and the two out of three property for weak equivalences.
\end{proof}

\begin{prop}
		\label{prop.3.2.disks-kill-spheres}
	For every compact generator $\generatorNRS \in C$
	(see proposition \ref{prop.3.1.stablehomotopy=>compactly-generated}),
	there exists a natural cofibration:
		$$\xymatrix{\generatorNRS \ar[rr]^-{\iota^{U}_{n,r,s}} && \diskNRRS}
		$$
	in $\motivicTspectra$.
\end{prop}
\begin{proof}
	We define $\killNRSmap$ as
	$F_{n}(\iota _{1}\wedge \gm ^{s}\wedge U_{+})$,
	where $\iota _{1}:S^{r}\rightarrow D^{r+1}$ is the map
	constructed in definition \ref{def.3.2.simplicial-disks}.
	
	It is clear that $\iota _{1}$ is a cofibration of pointed simplicial sets,
	therefore the result follows from propositions \ref{prop.2.4.Ev-n=Quillen-stable-functor} and 
	\ref{prop.2.3.cell-prop-pointedmotcat} which imply that
	$F_{n}$ and $-\wedge \gm ^{s} \wedge U_{+}$ are both left Quillen
	functors.
\end{proof}

\begin{thm}
		\label{thm.3.2.Lqmodelstructures}
	Fix $q\in \mathbb Z$, and consider the following
	set of maps in $\motivicTspectra$:
		\begin{eqnarray}
				\label{diagram.3.2.mapsleftloc}
			\qleftlocmaps = \{ \killNRS \: | \\
										 \generatorNRS \in C_{eff}^{q}\} \nonumber
		\end{eqnarray}
	Then the left Bousfield localization of $\motivicTspectra$ with respect to
	the $\qleftlocmaps$-local equivalences exist.  This new model structure
	will be called \emph{weight$^{<q}$ motivic stable}.
	$\weightqTspectra$ will denote the category of $T$-spectra equipped with
	the weight$^{<q}$ motivic stable model
	structure, and $\weightqstablehomotopy$ will denote
	its associated homotopy category.  Furthermore the weight$^{<q}$ motivic stable model structure
	is cellular, left proper and simplicial; with the following
	sets of generating cofibrations and trivial cofibrations
	respectively:
	$$\begin{array}{l}
	 I_{\qleftlocmaps}=I^{T}_{M_{*}}  =\bigcup _{n\geq 0}\{F_{n}(Y_{+}\hookrightarrow (\Delta ^{n}_{U})_{+})\} \\
	 \\
	J_{\qleftlocmaps}  =\{j:A\rightarrow B\}
	\end{array}
	$$
	where $j$ satisfies the following conditions:
	\begin{enumerate}
		\item	$j$ is an inclusion of $I^{T}_{M_{*}}$-complexes.
		\item	$j$ is a $\qleftlocmaps$-local equivalence.
		\item	the size of $B$ as an $I^{T}_{M_{*}}$-complex is less than $\kappa$, 
					where $\kappa$ is the regular cardinal defined by Hirschhorn in \cite[definition 4.5.3]{MR1944041}.
	\end{enumerate}
\end{thm}
\begin{proof}
	Theorems \ref{thm.2.5.cellularity-motivic-stable-str} and 
	\ref{thm.2.4.stableTspectramodelstr} imply that $\motivicTspectra$ is a cellular, proper and simplicial
	model category.  Therefore the existence of the left Bousfield localization
	follows from \cite[theorem 4.1.1]{MR1944041}.  Using \cite[theorem 4.1.1]{MR1944041} again, we have that
	$\weightqTspectra$ is cellular, left proper and simplicial; where the sets
	of generating cofibrations and trivial cofibrations are the ones described above.
\end{proof}

\begin{defi}
		\label{def.3.2.stable-weightq-replacementfunctors}
	Fix $q\in \mathbb Z$ and let $W_{q}$ denote a fibrant replacement functor in
	$\weightqTspectra$, such that the for every $T$-spectrum $X$,
	the natural map:
		$$\xymatrix{X \ar[r]^-{W_{q}^{X}}& W_{q}X}
		$$
	is a trivial cofibration in $\weightqTspectra$, and $W_{q}X$
	is $\qleftlocmaps$-local in $\motivicTspectra$.
\end{defi}

\begin{prop}
		\label{prop.3.2.Qs-cofibrant-replacement-all-L<q}
	Fix $q\in \mathbb Z$.  Then $Q_{s}$ is also a cofibrant
	replacement functor in $\weightqTspectra$, and for every
	$T$-spectrum $X$ the natural map
		$$\xymatrix{Q_{s}X\ar[r]^{Q_{s}^{X}}& X}
		$$
	is a trivial fibration in $\weightqTspectra$.
\end{prop}
\begin{proof}
	Since $\weightqTspectra$ is the left Bousfield localization of
	$\motivicTspectra$ with respect to the $\qleftlocmaps$-local equivalences, by construction
	we have that the cofibrations and the trivial fibrations are indentical in
	$\weightqTspectra$ and $\motivicTspectra$ respectively.  This implies that for every
	$T$-spectrum $X$, $Q_{s}X$ is cofibrant in $\weightqTspectra$, and we also have that
	the natural map
		$$\xymatrix{Q_{s}X\ar[r]^{Q_{s}^{X}}& X}
		$$
	is a trivial fibration in $\weightqTspectra$.  Hence $Q_{s}$ is also a cofibrant replacement
	functor for $\weightqTspectra$.
\end{proof}

\begin{prop}
		\label{prop.3.2.Lq-local-objects-classification}
	Fix $q\in \mathbb Z$ and let $Z$ be an arbitrary $T$-spectrum.
	Then $Z$ is $\qleftlocmaps$-local in $\motivicTspectra$ if and only if
	the following conditions hold:
		\begin{enumerate}
			\item \label{prop.3.2.Lq-local-objects-classification.a}$Z$ is fibrant 
						in $\motivicTspectra$.
			\item \label{prop.3.2.Lq-local-objects-classification.b}For every $\generatorNRS 
						\in C_{eff}^{q}$,
						$[\generatorNRS ,Z]_{Spt}\cong 0$
		\end{enumerate}
\end{prop}
\begin{proof}
	($\Rightarrow$):  Assume that $Z$ is $\qleftlocmaps$-local.
	Then by definition we have that $Z$ must be fibrant in
	$\motivicTspectra$.
	Since all the $T$-spectra $\generatorNRS$ and $\diskNRS$
	are cofibrant, and $Z$ is $\qleftlocmaps$-local; for every
	$\generatorNRS \in C_{eff}^{q}$
	we get the following weak equivalence of simplicial sets:
		$$\xymatrix{Map(\diskNRRS, Z) \ar[rr]^-{(\killNRSmap )^{\ast}}&& Map(\generatorNRS , Z)}
		$$
	Now we have that $\motivicTspectra$ is in particular a simplicial model category, therefore
	we get the following commutative diagram:
		$$\xymatrix{\pi_{0}Map(\diskNRRS ,Z) \ar[rr]^-{(\killNRSmap )^{\ast}}_-{\cong} \ar[d]_-{\cong}
								&& \pi_{0}Map(\generatorNRS , Z) \ar[d]^-{\cong}\\
								[\diskNRRS , Z]_{Spt} \ar[rr]_-{(\killNRSmap )^{\ast}}&& [\generatorNRS , Z]_{Spt}}
		$$
	where the vertical arrows and the top row are isomorphisms.
	Therefore we get the following isomorphism: 
		$$\xymatrix{[\diskNRRS , Z]_{Spt} \ar[rr]^-{(\killNRSmap )^{\ast}}_-{\cong}&& [\generatorNRS , Z]_{Spt}}
		$$
	Finally proposition \ref{prop.3.2.Disks=>compact}(\ref{prop.3.2.Disks=>compact.b})
	implies that $[\diskNRRS ,Z]_{Spt}\cong 0$.  Thus,
	for every $\generatorNRS \in C_{eff}^{q}$ we have that
	$[\generatorNRS ,Z]_{Spt}\cong 0$, as we wanted.
	
	($\Leftarrow$):  Assume that $Z$ satisfies (\ref{prop.3.2.Lq-local-objects-classification.a})
	and (\ref{prop.3.2.Lq-local-objects-classification.b}). Let $\omega _{0}$, $\eta _{0}$
	denote the base points corresponding to the pointed simplicial sets
	$Map_{\ast}(\diskNRRS,Z)$ and $Map_{\ast}(\generatorNRS,Z)$ respectively.
	Since
	$\generatorNRS$ and $\diskNRS$ are always cofibrant, it is enough to show that
	the induced map:
		$$\xymatrix{Map(\diskNRRS, Z) \ar[rr]^-{(\killNRSmap )^{\ast}}&& Map(\generatorNRS , Z)}
		$$
	is a weak equivalence of simplicial sets for every map $\killNRSmap \in \qleftlocmaps$.
	
	Fix $\killNRSmap \in  \qleftlocmaps$.
	By proposition \ref{prop.3.2.Disks=>compact}(\ref{prop.3.2.Disks=>compact.c})
	we know that the map $\diskNRRS \rightarrow \ast$ is a weak equivalence in $\motivicTspectra$.  
	Then Ken Brown's lemma (see lemma \ref{lem1.1.KenBrown})
	together with the fact that $\motivicTspectra$ is a simplicial model category, imply
	that the following map is a weak equivalence of simplicial sets:
		$$\xymatrix{\ast \cong Map(\ast,Z) \ar[r]& Map(\diskNRRS ,Z)}
		$$
	In particular $Map(\diskNRRS ,Z)$ has only one
	path connected component.
	
	Since $\motivicTspectra$ is a simplicial model category, we have the following
	isomorphism of abelian groups
		$$\xymatrix{\pi _{0}Map(\generatorNRS,Z) \ar[r]^-{\cong}&
								[\generatorNRS,Z]_{Spt}}
		$$
	but our hypothesis implies that $[\generatorNRS,Z]_{Spt}\cong 0$, hence
	$\pi _{0}Map(\generatorNRS,Z)\cong 0$, i.e. $Map(\generatorNRS,Z)$ has only one path
	connected component.
	
	Now proposition \ref{prop.3.2.Disks=>compact}(\ref{prop.3.2.Disks=>compact.b}) 
	implies that $\ast \rightarrow \diskNRRS$
	is a trivial cofibration in $\motivicTspectra$, 
	and since $- \wedge S^{1}$ is a left Quillen functor, it follows that
		$$\xymatrix{\ast \cong \ast \wedge S^{k} \ar[r]& \diskNRRS \wedge S^{k}}
		$$
	is also a trivial cofibration for $k\geq 0$.  
	Therefore $[\diskNRRS \wedge S^{k},Z]_{Spt}\cong 0$, and this implies that
	the induced map
	$(\killNRSmap )^{\ast}$ is an isomorphism of abelian groups:
		$$\xymatrix{0\cong [\diskNRRS \wedge S^{k}, Z]_{Spt} \ar[d]_-{(\killNRSmap )^{\ast}}\\ 
								[\generatorNRS \wedge S^{k} , Z]_{Spt}\cong
								[F_{n}(S^{k+r}\wedge \gm ^{s}\wedge U_{+}) , Z]_{Spt}}
		$$
	since by hypothesis $[F_{n}(S^{k+r}\wedge \gm ^{s}\wedge U_{+}) , Z]_{Spt}\cong 0$.
		
	On the other hand, since $\motivicTspectra$ is a pointed simplicial model category,
	we have that lemma 6.1.2 in \cite{MR1650134} together with
	remark \ref{rmk.2.4.simp-simppointed-structures-Tspectra}(\ref{rmk.2.4.simp-simppointed-structures-Tspectra.b})
	imply that
	the following diagram is commutative for $k\geq 0$ and all the vertical arrows are isomorphisms:
		$$\xymatrix@C=-1pc{\pi_{k,\omega _{0}}Map(\diskNRRS , Z) \ar[dr]_-{(\killNRSmap )^{\ast}}\ar@{=}[dd]&\\  
								& \pi_{k,\eta _{0}}Map(\generatorNRS , Z)\ar@{=}[dd]\\
								\pi_{k,\omega _{0}}Map_{\ast}(\diskNRRS , Z) \ar[dr]_-{(\killNRSmap )^{\ast}}\ar[dd]_-{\cong}&\\  
								& \pi_{k,\eta _{0}}Map_{\ast}(\generatorNRS , Z)\ar[dd]^-{\cong}\\
								[\diskNRRS \wedge S^{k}, Z]_{Spt} \ar[dr]_-{(\killNRSmap )^{\ast}} \ar[dd]_-{\cong}&\\ 
								& [\generatorNRS \wedge S^{k}, Z]_{Spt} \ar[dd]^-{\cong} \\
								[\diskNRRS \wedge S^{k}, Z]_{Spt} \ar[dr]_-{(\killNRSmap )^{\ast}} &\\ 
								& [F_{n}(S^{k+r}\wedge \gm ^{s}\wedge U_{+}) , Z]_{Spt}}
		$$
	but we know that the bottom row is always an isomorphism of abelian groups,
	hence the top row is also an isomorphism.  This implies that the map
		$$\xymatrix{Map(\diskNRRS , Z) \ar[rr]^-{(\killNRSmap )^{\ast}}&&  
								Map(\generatorNRS , Z)}
		$$
	is a weak equivalence when it is restricted to the path component of $Map(\diskNRRS,Z)$ containing
	$\omega _{0}$.  However we already know that $Map(\diskNRRS,Z)$ and $Map(\generatorNRS,Z)$ have
	only one path connected component.  This implies that the map defined above
	is a weak equivalence of simplicial sets, as we wanted.
\end{proof}

\begin{cor}
		\label{cor.3.2.m>n===>Ln-local=>Lm-local}
	Let $m,n \in \mathbb Z$
	with $m>n$.  If $Z$ is a $L(<n)$-local $T$-spectrum in $\motivicTspectra$ then
	$Z$ is also $L(<m)$-local in $\motivicTspectra$.
\end{cor}
\begin{proof}
	We have that $C_{eff}^{m}\subseteq C_{eff}^{n}$, since $m>n$.
	The result now
	follows immediately from the characterization of $\qleftlocmaps$-local
	objects given in proposition \ref{prop.3.2.Lq-local-objects-classification}.
\end{proof}

\begin{cor}
		\label{cor.3.2.S1loops-preserves-Lqlocal}
	Fix $q\in \mathbb Z$ and let $Z$ be a fibrant $T$-spectrum in $\motivicTspectra$.
	Then $Z$ is $\qleftlocmaps$-local if and only if
	$\Omega _{S^{1}}Z$ is $\qleftlocmaps$-local.
\end{cor}
\begin{proof}
	($\Rightarrow$):  Assume that $Z$ is $\qleftlocmaps$-local.
	We have that $Z$ is fibrant in $\motivicTspectra$; and since 
	$\motivicTspectra$ is a simplicial model category,
	it follows that $\Omega _{S^{1}}Z$ is also fibrant.
	
	Fix $\generatorNRS \in C_{eff}^{q}$.
	Since $\motivicTspectra$ is a simplicial model category, we have the following
	natural isomorphisms:
		\begin{eqnarray*}
			[\generatorNRS, \Omega _{S^{1}}Z]_{Spt} &\cong & [\generatorNRS \wedge S^{1},Z]_{Spt}\\
			&\cong & [F_{n}(S^{r+1}\wedge \gm ^{s}\wedge U_{+}),Z]_{Spt}
		\end{eqnarray*}
	but proposition \ref{prop.3.2.Lq-local-objects-classification} 
	implies that $[F_{n}(S^{r+1}\wedge \gm ^{s}\wedge U_{+}),Z]_{Spt}\cong 0$, hence
	$[\generatorNRS, \Omega _{S^{1}}Z]_{Spt}\cong 0$ for every $\generatorNRS \in C_{eff}^{q}$.
	Finally, using proposition \ref{prop.3.2.Lq-local-objects-classification} again,
	we have that $\Omega _{S^{1}}Z$ is $\qleftlocmaps$-local, as we wanted.
	
	($\Leftarrow$):  Assume that $\Omega _{S^{1}}Z$ is $\qleftlocmaps$-local.
	Since by hypothesis $Z$ is fibrant in $\motivicTspectra$, 
	proposition \ref{prop.3.2.Lq-local-objects-classification}
	implies that it is enough to show that for every $\generatorNRS \in C_{eff}^{q}$:
		$$[\generatorNRS , Z]_{Spt}\cong 0
		$$
	Fix $\generatorNRS \in C_{eff}^{q}$.  Since $\motivicTspectra$ 
	is a simplicial model category, and $Z$
	is fibrant by hypothesis; we have the following natural
	isomorphisms of abelian groups:
		\begin{eqnarray*}
			[F_{n+1}(S^{r}\wedge \gm ^{s+1}\wedge U_{+}) , \Omega _{S^{1}}Z]_{Spt} &\cong &
			[F_{n+1}(S^{r+1}\wedge \gm ^{s+1}\wedge U_{+}),Z]_{Spt}\\
			&\cong & [\generatorNRS ,Z]_{Spt}
		\end{eqnarray*}
	Now using proposition \ref{prop.3.2.Lq-local-objects-classification} and the fact that
	$\Omega _{S^{1}}Z$ is $\qleftlocmaps$-local, it follows that
	$[F_{n+1}(S^{r}\wedge \gm ^{s+1}\wedge U_{+}) , \Omega _{S^{1}}Z]_{Spt}\cong 0$.
	Therefore, $[\generatorNRS ,Z]_{Spt}\cong 0$ for every $\generatorNRS \in C_{eff}^{q}$,
	as we wanted.
\end{proof}

\begin{cor}
		\label{cor.3.2.S1-preserves-Lqlocal}
	Fix $q\in \mathbb Z$, and let $Z$ be a fibrant
	$T$-spectrum in $\motivicTspectra$.  Then $Z$ is $\qleftlocmaps$-local if and only if
	$IQ_{T}J(Q_{s}Z\wedge S^{1})$ is $\qleftlocmaps$-local.
\end{cor}
\begin{proof}
	($\Rightarrow$):  Assume that $Z$ is $\qleftlocmaps$-local.
	Since $IQ_{T}J(Q_{s}Z\wedge S^{1})$ is fibrant, using proposition
	\ref{prop.3.2.Lq-local-objects-classification} we have that it is enough
	to check that for every $\generatorNRS \in C_{eff}^{q}$,
	$[\generatorNRS ,IQ_{T}J(Q_{s}Z\wedge S^{1})]_{Spt}\cong 0$.
	But since $-\wedge S^{1}$ is a Quillen equivalence, we get the
	following diagram:
		$$\xymatrix@C=-4pc{[\generatorNRS ,IQ_{T}J(Q_{s}Z\wedge S^{1})]_{Spt} \ar[dr]^-{\cong}&\\ 
								& [F_{n+1}(S^{r+1}\wedge \gm ^{s+1}\wedge U_{+}),IQ_{T}J(Q_{s}Z\wedge S^{1})]_{Spt}\\
								[F_{n+1}(S^{r}\wedge \gm ^{s+1}\wedge U_{+}), Z]_{Spt} \ar[dr]^-{\Tsuspfunctor ^{1,0}}_-{\cong}&\\ 
								& [F_{n+1}(S^{r+1}\wedge \gm ^{s+1}\wedge U_{+}), Q_{s}Z\wedge S^{1}]_{Spt} \ar[uu]_-{\cong}}
		$$
	where all the maps are isomorphisms of abelian groups.  Since $Z$
	is $\qleftlocmaps$-local, proposition \ref{prop.3.2.Lq-local-objects-classification}
	implies that $[F_{n+1}(S^{r}\wedge \gm ^{s+1}\wedge U_{+}), Z]_{Spt}\cong 0$.
	Therefore 
		$$[\generatorNRS ,IQ_{T}J(Q_{s}Z\wedge S^{1})]_{Spt}\cong 0$$ 
	for every
	$\generatorNRS \in C_{eff}^{q}$, as we wanted.
	
	($\Leftarrow$):  Assume that $IQ_{T}J(Q_{s}Z\wedge S^{1})$ is $\qleftlocmaps$-local.
	By hypothesis, $Z$ is fibrant; therefore proposition \ref{prop.3.2.Lq-local-objects-classification}
	implies that it is enough to show that for every $\generatorNRS \in C_{eff}^{q}$,
	$[\generatorNRS ,Z]_{Spt}\cong 0$.	
	Since $\motivicTspectra$ is a simplicial model category
	and $-\wedge S^{1}$ is a Quillen equivalence;
	we have the following diagram:
		$$\xymatrix{[\generatorNRS , \Omega _{S^{1}}IQ_{T}J(Q_{s}Z\wedge S^{1})]_{Spt} \ar[d]_-{\cong}& \\
								[\generatorNRS \wedge S^{1}, Q_{s}Z\wedge S^{1}]_{Spt}
								& [\generatorNRS ,Z]_{Spt}\ar[l]^-{\Tsuspfunctor ^{1,0}}_-{\cong}}
		$$
	where all the maps are isomorphisms of abelian groups.
	On the other hand, using corollary \ref{cor.3.2.S1loops-preserves-Lqlocal}
	we have that $\Omega _{S^{1}}IQ_{T}J(Q_{s}Z\wedge S^{1})$ is $\qleftlocmaps$-local.
	Therefore using proposition \ref{prop.3.2.Lq-local-objects-classification} again,
	we have that for every $\generatorNRS \in C_{eff}^{q}$:
		$$[\generatorNRS ,Z]_{Spt}\cong [\generatorNRS , \Omega _{S^{1}}IQ_{T}J(Q_{s}Z\wedge S^{1})]_{Spt}\cong 0
		$$
	and this finishes the proof.
\end{proof}

\begin{cor}
		\label{cor.3.2.detecting-Cq-local-equivalences}
	Fix $q\in \mathbb Z$ and let $f:X\rightarrow Y$ be a map in $\motivicTspectra$.
	Then $f$ is a $\qleftlocmaps$-local equivalence if and only if
	for every $\qleftlocmaps$-local $T$-spectrum $Z$, $f$ induces
	the following isomorphism of abelian groups:
		$$\xymatrix{[Y,Z]_{Spt} \ar[r]^-{f^{\ast}}& [X,Z]_{Spt}}
		$$
\end{cor}
\begin{proof}
	Suppose that $f$ is a $\qleftlocmaps$-local equivalence, then by
	definition the induced map:
		$$\xymatrix{Map(Q_{s}Y,Z) \ar[r]^-{(Q_{s}f)^{\ast}}& Map(Q_{s}X,Z)}
		$$
	is a weak equivalence of simplicial sets for every $\qleftlocmaps$-local
	$T$-spectrum $Z$.  Proposition \ref{prop.3.2.Lq-local-objects-classification}(\ref{prop.3.2.Lq-local-objects-classification.a}) 
	implies that $Z$ is fibrant in $\motivicTspectra$, and
	since $\motivicTspectra$ is in particular a simplicial model category; 
	we get the following commutative diagram, where the top row and all the
	vertical maps are isomorphisms of abelian groups:
		$$\xymatrix{\pi_{0}Map(Q_{s}Y,Z) \ar[r]^-{(Q_{s}f)^{\ast}}_-{\cong} \ar[d]_-{\cong}& 
								\pi_{0}Map(Q_{s}X,Z) \ar[d]^-{\cong}\\
								[Y,Z]_{Spt} \ar[r]_-{f^{\ast}}& [X,Z]_{Spt}}
		$$
	hence $f^{\ast}$ is an isomorphism for every $\qleftlocmaps$-local $T$-spectrum $Z$,
	as we wanted.
	
	Conversely, assume that for every $\qleftlocmaps$-local $T$-spectrum $Z$,
	the induced map
		$$\xymatrix{[Y,Z]_{Spt} \ar[r]^-{f^{\ast}}& [X,Z]_{Spt}}
		$$
	is an isomorphism of abelian groups.
	
	Since $\weightqTspectra$ is the left Bousfield localization
	of $\motivicTspectra$ with
	respect to the $\qleftlocmaps$-local equivalences, we have that the 
	identity functor $id:\Tspectra \rightarrow \weightqTspectra$ is
	a left Quillen functor.  Therefore for every $T$-spectrum $Z$,
	we get the following commutative diagram
	where all the vertical arrows are isomorphisms:
		$$\xymatrix{\Hom _{\weightqstablehomotopy}(Q_{s}Y,Z) \ar[r]^-{(Q_{s}f)^{\ast}} \ar[d]_-{\cong} &
								\Hom _{\weightqstablehomotopy}(Q_{s}X,Z) \ar[d]^-{\cong}\\
								[Y,W_{q}Z]_{Spt} \ar[r]_-{f^{\ast}}^-{\cong}& [X,W_{q}Z]_{Spt}}
		$$
	but $W_{q}Z$ is by construction $\qleftlocmaps$-local, then
	by hypothesis the bottom row is an isomorphism of abelian groups.
	Hence it follows that the induced map:
		$$\xymatrix{\Hom _{\weightqstablehomotopy}(Q_{s}Y,Z) \ar[r]^-{(Q_{s}f)^{\ast}}_{\cong}  &
								\Hom _{\weightqstablehomotopy}(Q_{s}X,Z)}
		$$
	is an isomorphism for every $T$-spectrum $Z$.
	This implies that $Q_{s}f$ is a weak equivalence in $\weightqTspectra$, and since
	$Q_{s}$ is also a cofibrant replacement functor in $\weightqTspectra$, it follows
	that $f$ is a weak equivalence in $\weightqTspectra$.
	Therefore we have that $f$ is a $\qleftlocmaps$-local equivalence, as we wanted.
\end{proof}

\begin{lem}
		\label{lem.3.2.Lq-localequivs-stable-S1suspension}
	Fix $q\in \mathbb Z$ and let $f:X\rightarrow Y$ be a map in $\motivicTspectra$,
	then $f$ is a $\qleftlocmaps$-local equivalence
	if and only if  
		$$Q_{s}f\wedge id: Q_{s}X\wedge S^{1} \rightarrow Q_{s}Y\wedge S^{1}$$ 
	is a
	$\qleftlocmaps$-local equivalence.
\end{lem}
\begin{proof}
	Assume that $f$ is a $\qleftlocmaps$-local equivalence, and let $Z$
	be an arbitrary $\qleftlocmaps$-local $T$-spectrum.  Then 
	corollary \ref{cor.3.2.S1loops-preserves-Lqlocal}
	implies that $\Omega _{S^{1}}Z$ is also $\qleftlocmaps$-local.
	Therefore the induced map
		$$\xymatrix{Map(Q_{s}Y,\Omega _{S^{1}}Z) \ar[r]^-{(Q_{s}f)^{\ast}}& Map(Q_{s}X,\Omega _{S^{1}}Z)}
		$$
	is a weak equivalence of simplicial sets.  Now
	since $\motivicTspectra$
	is a simplicial model category, we have the following commutative diagram:
		$$\xymatrix{Map(Q_{s}Y,\Omega _{S^{1}}Z) \ar[rr]^-{(Q_{s}f)^{\ast}}\ar[d]_-{\cong}&& Map(Q_{s}X,\Omega _{S^{1}}Z) \ar[d]^-{\cong}\\
								Map(Q_{s}Y\wedge S^{1},Z) \ar[rr]^-{(Q_{s}f\wedge id)^{\ast}}&& Map(Q_{s}X\wedge S^{1},Z)}
		$$
	and using the two out of three property for weak equivalences of simplicial sets,
	we have that
		$$\xymatrix{Map(Q_{s}Y\wedge S^{1},Z) \ar[rr]^-{(Q_{s}f\wedge id)^{\ast}}&& Map(Q_{s}X\wedge S^{1},Z)}
		$$
	is a weak equivalence.  Since this holds for every $\qleftlocmaps$-local $T$-spectrum
	$Z$, it follows that
		$$Q_{s}f\wedge id: Q_{s}X\wedge S^{1} \rightarrow Q_{s}Y\wedge S^{1}
		$$
	is a $\qleftlocmaps$-local equivalence, as we wanted.
	
	Conversely, suppose that 
		$$Q_{s}f\wedge id: Q_{s}X\wedge S^{1} \rightarrow Q_{s}Y\wedge S^{1}
		$$
	is a $\qleftlocmaps$-local equivalence.  Let $Z$ be an arbitrary
	$\qleftlocmaps$-local $T$-spectrum. Since $\motivicTspectra$
	is a simplicial model category and
	$-\wedge S^{1}$ is a Quillen equivalence, we get the following
	commutative diagram:
		$$\xymatrix{[Q_{s}Y\wedge S^{1}, IQ_{T}J(Q_{s}Z\wedge S^{1})]_{Spt} \ar[rr]^-{(Q_{s}f\wedge id)^{\ast}} \ar[d]_-{\cong}&& 
								[Q_{s}X\wedge S^{1}, IQ_{T}J(Q_{s}Z\wedge S^{1})]_{Spt}\ar[d]^-{\cong}\\
								[Q_{s}Y\wedge S^{1}, Q_{s}Z\wedge S^{1}]_{Spt} \ar[rr]^-{(Q_{s}f\wedge id)^{\ast}}&&
								[Q_{s}X\wedge S^{1}, Q_{s}Z\wedge S^{1}]_{Spt}\\
								[Y, Z]_{Spt} \ar[rr]_-{f^{\ast}}\ar[u]^-{\cong}_-{\Tsuspfunctor ^{1,0}}&&
								[X, Z]_{Spt}\ar[u]_-{\cong}^-{\Tsuspfunctor ^{1,0}} }
		$$
	Now, corollary \ref{cor.3.2.S1-preserves-Lqlocal} implies that $IQ_{T}J(QZ\wedge S^{1})$ is
	also $\qleftlocmaps$-local.  Therefore using corollary \ref{cor.3.2.detecting-Cq-local-equivalences}
	we have that the top row in the diagram above is an isomorphism of abelian groups.
	This implies that the induced map:
		$$\xymatrix{[Y, Z]_{Spt} \ar[r]^-{f^{\ast}}&
								[X, Z]_{Spt}}
		$$
	is an isomorphism of abelian groups for every $\qleftlocmaps$-local spectrum $Z$.
	Finally using corollary \ref{cor.3.2.detecting-Cq-local-equivalences} again,
	we have that $f:X\rightarrow Y$ is a $\qleftlocmaps$-local equivalence, as we
	wanted.
\end{proof}

\begin{cor}
		\label{cor.3.2.susp-Quillen-equiv-onLq}
	For every $q\in \mathbb Z$, the following adjunction:
		$$\xymatrix{(-\wedge S^{1},\Omega _{S^{1}},\varphi):\weightqTspectra \ar[rr]
								&& \weightqTspectra}
		$$
	is a Quillen equivalence.
\end{cor}
\begin{proof}
	Using corollary 1.3.16 in \cite{MR1650134} and
	proposition \ref{prop.3.2.Qs-cofibrant-replacement-all-L<q} we have that it suffices to
	verify the following two conditions:
		\begin{enumerate}
			\item	\label{cor.3.2.susp-Quillen-equiv-onLq.a}For every fibrant object $X$
						in $\weightqTspectra$, the following composition
							$$\xymatrix{ (Q_{s}\Omega _{S^{1}}X)\wedge S^{1} \ar[rr]^-{Q_{s}^{\Omega _{S^{1}}X}\wedge id}&& 
													 (\Omega _{S^{1}}X)\wedge S^{1} \ar[r]^-{\epsilon _{X}}& X}
							$$
						is a $\qleftlocmaps$-local equivalence.
			\item	\label{cor.3.2.susp-Quillen-equiv-onLq.b}$-\wedge S^{1}$ reflects $\qleftlocmaps$-local equivalences
						between cofibrant objects in $\weightqTspectra$.
		\end{enumerate}
		
	(\ref{cor.3.2.susp-Quillen-equiv-onLq.a}):  By construction $\weightqTspectra$ is
	a left Bousfield localization of $\motivicTspectra$, therefore the identity
	functor 
		$$\xymatrix{id:\weightqTspectra \ar[r]& \motivicTspectra}
		$$ 
	is a right Quillen functor.
	Thus $X$ is also fibrant in $\motivicTspectra$.
	Since the adjunction $(-\wedge S^{1},\Omega _{S^{1}},\varphi)$
	is a Quillen equivalence on $\motivicTspectra$, 
	\cite[proposition 1.3.13(b)]{MR1650134} implies that the following composition
	is a weak equivalence in $\motivicTspectra$:
		$$\xymatrix{(Q_{s}\Omega _{S^{1}}X)\wedge S^{1} \ar[rr]^-{Q_{s}^{\Omega _{S^{1}}X}\wedge id}&&
													(\Omega _{S^{1}}X)\wedge S^{1} \ar[r]^-{\epsilon _{X}}& X}
		$$
	Hence using \cite[proposition 3.1.5]{MR1944041} it follows that the composition above
	is a $\qleftlocmaps$-local equivalence.
	
	(\ref{cor.3.2.susp-Quillen-equiv-onLq.b}):  This follows immediately
	from proposition \ref{prop.3.2.Qs-cofibrant-replacement-all-L<q} and 
	lemma \ref{lem.3.2.Lq-localequivs-stable-S1suspension}.
\end{proof}

\begin{rmk}
		\label{rmk.3.2.Tsuspfunctor-doesnotdescend-to-Lq}
	We have a situation similar to the one
	described in remark
	\ref{rmk.3.3.smashT-not-descending} for the model categories
	$\qconnectedTspectra$; i.e. although
	the adjunction $(\Tsuspfunctor ,\Tloops ,\varphi)$ is
	a Quillen equivalence on $\motivicTspectra$,
	it does not descend even to a Quillen
	adjunction on the weight$^{<q}$ motivic stable
	model category
	$\weightqTspectra$.
\end{rmk}

\begin{cor}
		\label{cor.3.2.stablehomotopyLq==triangulated}
	For every $q\in \mathbb Z$,
	the homotopy category $\weightqstablehomotopy$
	associated to $\weightqTspectra$ has the
	structure of a triangulated category.
\end{cor}
\begin{proof}
	Theorem \ref{thm.3.2.Lqmodelstructures} implies in particular
	that $\weightqTspectra$ is a pointed simplicial model category,
	and  corollary \ref{cor.3.2.susp-Quillen-equiv-onLq} implies that
	the adjunction 
		$$(-\wedge S^{1},\Omega _{S^{1}},\varphi):\weightqTspectra \rightarrow \weightqTspectra$$
	is a Quillen equivalence.  Therefore
	the result follows from the work of Quillen in 
	\cite[sections I.2 and I.3]{MR0223432} and the work of 
	Hovey in \cite[chapters VI and VII]{MR1650134}.
\end{proof}

\begin{cor}
		\label{cor.3.2.Lq=>rightproper}
	For every $q\in \mathbb Z$, $\weightqTspectra$ is
	a right proper model category.
\end{cor}
\begin{proof}
	We need to show that the $\qleftlocmaps$-local equivalences are stable
	under pullback along fibrations in $\weightqTspectra$.
	Consider the following pullback diagram:
		$$\xymatrix{Z \ar[r]^-{w^{\ast}} \ar[d]_-{p^{\ast}}& X \ar[d]^-{p}\\
								W \ar[r]_-{w}& Y}
		$$
	where $p$ is a fibration in $\weightqTspectra$, and $w$ is a
	$\qleftlocmaps$-local equivalence.  Let $F$ be the homotopy fibre
	of $p$.  Then we get the following commutative diagram in $\weightqstablehomotopy$:
		$$\xymatrix{\Omega _{S^{1}}Y \ar[r]^-{q} & F \ar[r]^-{i}& X \ar[r]^-{p}& Y\\
								\Omega _{S^{1}}W \ar[r]_-{r} \ar[u]^-{\Omega _{S^{1}}w}& F \ar[r]_-{j} \ar@{=}[u]& 
								Z \ar[r]_-{p^{\ast}} \ar[u]_-{w^{\ast}}& W \ar[u]_-{w}}
		$$
	Since the rows in the diagram above are both fibre sequences in $\weightqTspectra$,
	it follows that both rows are distinguished triangles in $\weightqstablehomotopy$
	(which has the structure of a triangulated category given by 
	corollary \ref{cor.3.2.stablehomotopyLq==triangulated}).  Now
	$w, id_{F}$ are both isomorphisms in $\weightqstablehomotopy$, hence it follows that
	$w^{\ast}$ is also an isomorphism in $\weightqstablehomotopy$.  Therefore
	$w^{\ast}$ is a $\qleftlocmaps$-local equivalence, as we wanted.
\end{proof}

\begin{prop}
		\label{prop.3.2.Lq-exact-adjunctions}
	For every $q\in \mathbb Z$ we have the following
	adjunction
		$$\xymatrix{(Q_{s},W_{q},\varphi):\stablehomotopy \ar[r]& \weightqstablehomotopy}
		$$
	of exact functors between triangulated categories.
\end{prop}
\begin{proof}
	Since $\weightqTspectra$ is the left Bousfield localization of
	$\motivicTspectra$ with respect to the $\qleftlocmaps$-local equivalences, we have that the
	identity functor $id:\motivicTspectra \rightarrow \weightqTspectra$
	is a left Quillen functor.  Therefore we get
	the following adjunction at the level of the associated homotopy categories:
		$$\xymatrix{(Q_{s}, W_{q}, \varphi ):\stablehomotopy \ar[r]& \weightqstablehomotopy}
		$$
	
	Now proposition 6.4.1 in \cite{MR1650134} implies that
	$Q_{s}$ maps cofibre sequences in $\stablehomotopy$ to cofibre sequences in
	$\weightqstablehomotopy$.
	Therefore 
	using proposition 7.1.12 in \cite{MR1650134} we have
	that $Q_{s}$ and $W_{q}$ are both exact functors between triangulated categories.
\end{proof}

\begin{prop}
		\label{prop.3.2.counit-spt-Lq-properties}
	Fix $q\in \mathbb Z$ and let $\eta _{X}:Q_{s}W_{q}X\rightarrow X$ denote the
	counit of the adjunction 
			$$\xymatrix{(Q_{s},W_{q},\varphi ):\stablehomotopy \ar[r]& \weightqstablehomotopy}
			$$
	Then the following conditions hold:
		\begin{enumerate}
			\item	\label{prop.3.2.counit-spt-Lq-properties.a}For every $T$-spectrum $X$, we have that
						$\eta _{X}$ is an isomorphism in $\weightqstablehomotopy$.
			\item \label{prop.3.2.counit-spt-Lq-properties.b}The exact functor
							$$\xymatrix{W_{q}:\weightqstablehomotopy \ar[r]& \stablehomotopy}
							$$
						is a full embedding of triangulated categories.
		\end{enumerate}
\end{prop}
\begin{proof}
	(\ref{prop.3.2.counit-spt-Lq-properties.a}):  We have that
	$\eta _{X}$ is the following composition in $\weightqstablehomotopy$:
		$$\xymatrix{Q_{s}W_{q}X \ar[r]^-{Q_{s}^{W_{q}X}}& 
								W_{q}X \ar[r]^-{(W_{q}^{X})^{-1}}_{\cong}& X}
		$$
	where $Q_{s}^{W_{q}X}$ is a weak equivalence in $\motivicTspectra$.
	Now \cite[proposition 3.1.5]{MR1944041} implies that
	$Q_{s}^{W_{q}X}$ is a $\qleftlocmaps$-local equivalence, i.e.
	a weak equivalence in $\weightqTspectra$.  Therefore $Q_{s}^{W_{q}X}$
	becomes an isomorphism in $\weightqstablehomotopy$, and this implies that
	$\eta _{X}$ is an isomorphism in $\weightqstablehomotopy$, as we wanted.
	
	(\ref{prop.3.2.counit-spt-Lq-properties.b}):  Follows immediately from
	(\ref{prop.3.2.counit-spt-Lq-properties.a}).
\end{proof}

\begin{prop}
		\label{prop.3.2.Fnrs-inCq-vanishinL<q}
	Fix $q\in \mathbb Z$.  Then for every
	$\generatorNRS \in C_{eff}^{q}$, the 
	map $\ast \rightarrow \generatorNRS$ 
	is a $\qleftlocmaps$-local equivalence in $\motivicTspectra$.
\end{prop}
\begin{proof}
	Let $Z$ be an arbitrary $\qleftlocmaps$-local $T$-spectrum.  Then proposition
	\ref{prop.3.2.Lq-local-objects-classification}(\ref{prop.3.2.Lq-local-objects-classification.b}) 
	implies that the following induced map
		$$\xymatrix{0\cong [\generatorNRS,Z]_{Spt} \ar[r]& [\ast ,Z]_{Spt}\cong 0}
		$$
	is an isomorphism of abelian groups.  Therefore using 
	corollary \ref{cor.3.2.detecting-Cq-local-equivalences},
	it follows that $\ast \rightarrow \generatorNRS$ is a $\qleftlocmaps$-local equivalence.
\end{proof}

\begin{prop}
		\label{prop.3.2.detecting-Qs-isos.2}
	Fix $q\in \mathbb Z$ and let $f:X\rightarrow Y$ be a map in $\weightqstablehomotopy$.
	Then $f$ is an isomorphism
	in $\weightqstablehomotopy$ if and only if
	one of the following equivalent conditions holds:
		\begin{enumerate}
			\item	\label{prop.3.2.detecting-Qs-isos.2.a}The following map
							$$\xymatrix{W_{q}X \ar[r]^-{W_{q}(f)}_-{\cong}& W_{q}Y}
							$$
						is an isomorphism in $\stablehomotopy$.
			\item \label{prop.3.2.detecting-Qs-isos.2.b}For every
						$\generatorNRS \notin C_{eff}^{q}$, the induced map
							$$\xymatrix{[\generatorNRS ,W_{q}X]_{Spt} \ar[r]^-{(W_{q}f)_{\ast}}_-{\cong}& 
													[\generatorNRS ,W_{q}Y]_{Spt}}
							$$
						is an isomorphism of abelian groups.
			\item	\label{prop.3.2.detecting-Qs-isos.2.c}For every $\generatorNRS \notin C_{eff}^{q}$,
						the induced map
							$$\xymatrix{\Hom _{\weightqstablehomotopy}(Q_{s}\generatorNRS, X) \ar[d]^-{f_{\ast}}\\ 
													\Hom _{\weightqstablehomotopy}(Q_{s}\generatorNRS, Y)}
							$$
						is an isomorphism of abelian groups.
		\end{enumerate}
\end{prop}
\begin{proof}
	Proposition \ref{prop.3.2.counit-spt-Lq-properties} implies that 
	$f$ is an isomorphism in $\weightqstablehomotopy$ if
	and only if $W_{q}f$ becomes an isomorphism in $\stablehomotopy$.  Thus
	it only remains to show that (\ref{prop.3.2.detecting-Qs-isos.2.a}),
	(\ref{prop.3.2.detecting-Qs-isos.2.b}) and (\ref{prop.3.2.detecting-Qs-isos.2.c})
	are all equivalent.
	
	(\ref{prop.3.2.detecting-Qs-isos.2.a}) $\Leftrightarrow$ (\ref{prop.3.2.detecting-Qs-isos.2.b})
	Corollary \ref{cor.3.1.detecting-isos-in-SH} implies that 
	$W_{q}f$ is an isomorphism in $\stablehomotopy$ if and only if
	for every $\generatorNRS \in C$ the following induced map
		$$\xymatrix{[\generatorNRS ,W_{q}X]_{Spt} \ar[r]^-{(W_{q}f)_{\ast}}_-{\cong}& 
													[\generatorNRS ,W_{q}Y]_{Spt}}
		$$
	is an isomorphism of abelian groups.  But using proposition
	\ref{prop.3.2.Lq-local-objects-classification}(\ref{prop.3.2.Lq-local-objects-classification.b}) we have that
	for every $\generatorNRS \in C_{eff}^{q}$,
		$$0\cong [\generatorNRS ,W_{q}X]_{Spt}\cong [\generatorNRS, W_{q}Y]_{Spt}
		$$
	since by construction $W_{q}X$ and $W_{q}Y$ are both $\qleftlocmaps$-local $T$-spectra.
	Hence $W_{q}f$ is an isomorphism in $\stablehomotopy$ if and only if
	for every $\generatorNRS \notin C_{eff}^{q}$ the following induced map
		$$\xymatrix{[\generatorNRS ,W_{q}X]_{Spt} \ar[r]^-{(W_{q}f)_{\ast}}_-{\cong}& 
													[\generatorNRS ,W_{q}Y]_{Spt}}
		$$
	is an isomorphism of abelian groups.
	
	(\ref{prop.3.2.detecting-Qs-isos.2.b}) $\Leftrightarrow$ (\ref{prop.3.2.detecting-Qs-isos.2.c})
	By proposition \ref{prop.3.2.Lq-exact-adjunctions} 
	we have the following adjunction between exact functors of triangulated categories:
		$$\xymatrix{(Q_{s},W_{q},\varphi):\stablehomotopy \ar[r]& \weightqstablehomotopy}
		$$
	In particular for every $\generatorNRS \notin C_{eff}^{q}$, we get the following commutative
	diagram, where all the vertical arrows are isomorphisms of abelian groups:
		$$\xymatrix@C=-1.5pc{[\generatorNRS ,W_{q}X]_{Spt} \ar[dr]^-{(W_{q}f)_{\ast}} \ar[dd]_-{\cong}&\\ 
								& [\generatorNRS ,W_{q}Y]_{Spt}\ar[dd]^-{\cong}\\
								\Hom _{\weightqstablehomotopy}(Q_{s}\generatorNRS ,X) \ar[dr]^-{f_{\ast}}&\\
								& \Hom _{\weightqstablehomotopy}(Q_{s}\generatorNRS ,Y)}
		$$
	therefore the top row is an isomorphism if and only if the bottom row is an isomorphism
	of abelian groups, as we wanted.
\end{proof}

\begin{lem}
		\label{lemma.3.2.Lq-local=>fq-vanish}
	Fix $q\in \mathbb Z$ and
	let $Z$ be a $\qleftlocmaps$-local $T$-spectrum.
	Then $f_{q}Z\cong \ast$ in $\stablehomotopy$ 
	(see remark \ref{rmk.3.1.unit=iso}).
\end{lem}
\begin{proof}
	Let $j:\ast \rightarrow Z$ denote the canonical map.
	Proposition \ref{prop.3.1.detecting-fq-isos}
	implies that $f_{q}(j):\ast \cong f_{q}(\ast) \rightarrow f_{q}X$ is
	an isomorphism in $\stablehomotopy$ if and only if
	for every $\generatorNRS \in C_{eff}^{q}$ the induced map
		$$\xymatrix{0\cong [\generatorNRS ,\ast]_{Spt} \ar[r]^-{f_{q}(j)_{\ast}}& [\generatorNRS ,Z]_{Spt}}
		$$
	is an isomorphism of abelian groups.
	Therefore it is enough to show that for every $\generatorNRS \in C_{eff}^{q}$, we have
	$[\generatorNRS ,Z]_{Spt}\cong 0$.  But this follows from proposition
	\ref{prop.3.2.Lq-local-objects-classification}(\ref{prop.3.2.Lq-local-objects-classification.b}),
	since $Z$ is $\qleftlocmaps$-local by hypothesis.
\end{proof}

\begin{cor}
		\label{cor.3.2.Qsfq-vanishesin-LqSH}
	For every $q\in \mathbb Z$, and for every $T$-spectrum $X$, $Q_{s}f_{q}X\cong \ast$
	in $\weightqstablehomotopy$.
\end{cor}
\begin{proof}
	We will show that the map $\ast \rightarrow Q_{s}f_{q}X$
	is an isomorphism in $\weightqstablehomotopy$.  By Yoneda's lemma
	it suffices to check that for every $T$-spectrum $Z$,
	the induced map
		$$\xymatrix{\Hom _{\weightqstablehomotopy}(Q_{s}f_{q}X,Z) \ar[r]& 
								\Hom _{\weightqstablehomotopy}(\ast ,Z)\cong 0}
		$$
	is an isomorphism of abelian groups.  Now propositions \ref{prop.3.2.Lq-exact-adjunctions}
	and \ref{prop.3.1.adjunctions-stableeffq}
	imply that we have the following isomorphisms:
		\begin{eqnarray*}
			\Hom _{\weightqstablehomotopy}(Q_{s}f_{q}X,Z) & \cong &
			[f_{q}X,W_{q}Z]_{Spt}=[i_{q}r_{q}X,W_{q}Z]_{Spt} \\
			& \cong & \Hom _{\stablehomotopyeffq}(r_{q}X, r_{q}W_{q}Z)
		\end{eqnarray*}
	Finally since $i_{q}$ is a full embedding, we have
		$$\Hom _{\stablehomotopyeffq}(r_{q}X, r_{q}W_{q}Z)\cong
			[i_{q}r_{q}X,i_{q}r_{q}W_{q}Z]_{Spt}=[f_{q}X,f_{q}W_{q}Z]_{Spt}
		$$
	and lemma \ref{lemma.3.2.Lq-local=>fq-vanish} implies that $f_{q}W_{q}Z\cong \ast$
	in $\stablehomotopy$.  Hence
		$$\Hom _{\weightqstablehomotopy}(Q_{s}f_{q}X,Z)\cong
			[f_{q}X,f_{q}W_{q}Z]_{Spt}\cong [f_{q}X,\ast]_{Spt}\cong 0
		$$
	as we wanted.
\end{proof}

\begin{prop}
		\label{prop.3.2.Qs&QsSq==isoin-LqSH}
	For every $q\in \mathbb Z$ and for every $T$-spectrum $X$,
	the natural map in $\weightqstablehomotopy$
		$$\xymatrix{Q_{s}X \ar[rr]^-{Q_{s}(\pi _{<q}X)}&& Q_{s}s_{<q}X}
		$$
	is an isomorphism,
	where $\pi _{<q}$ is the natural transformation defined in theorem \ref{thm.3.1.motivictower}.  
	Furthermore, these maps
	induce a natural isomorphism between the following exact functors 
		$$\xymatrix{\stablehomotopy \ar@<1ex>[r]^-{Q_{s}} \ar@<-1ex>[r]_-{Q_{s}s_{<q}}& \weightqstablehomotopy}
		$$	
\end{prop}
\begin{proof}
	The naturality of $\pi _{<q}$ and the fact that $Q_{s}$ is a functor imply that
	the maps $Q_{s}(\pi _{<q}X)$ induce a natural transformation
	$Q_{s}\rightarrow Q_{s}s_{<q}$.  Hence it suffices to show that for every
	$T$-spectrum $X$, the map $Q_{s}(\pi _{<q}X)$ is an isomorphism 
	in $\weightqstablehomotopy$.
	
	Theorem \ref{thm.3.1.motivictower} implies that we have the following distinguished
	triangle in $\stablehomotopy$:
		$$\xymatrix{f_{q}X \ar[r]& X \ar[r]^-{\pi _{<q}X}& s_{<q}X \ar[r]^-{\sigma _{<q}X}& 
								\Tsuspfunctor ^{1,0}f_{q}X}
		$$
	and using proposition \ref{prop.3.2.Lq-exact-adjunctions}, we get the following distinguished
	triangle in $\weightqstablehomotopy$:
		$$\xymatrix{Q_{s}f_{q}X \ar[r]& Q_{s}X \ar[rr]^-{Q_{s}(\pi _{<q}X)}&& Q_{s}s_{<q}X \ar[rr]^-{Q_{s}(\sigma _{<q}X)}&& 
								\Tsuspfunctor ^{1,0}Q_{s}f_{q}X}
		$$
	But corollary \ref{cor.3.2.Qsfq-vanishesin-LqSH} 
	implies that $Q_{s}f_{q}X\cong \ast$ in $\weightqstablehomotopy$,
	therefore $Q_{s}(\pi _{<q}X)$ is an isomorphism in $\weightqstablehomotopy$,
	as we wanted.
\end{proof}

\begin{cor}
		\label{cor.3.2.WqQs===WqQss<q}
	For every $q\in \mathbb Z$ and for every $T$-spectrum $X$,
	the natural map in $\stablehomotopy$
		$$\xymatrix{W_{q}Q_{s}X \ar[rrr]^-{W_{q}Q_{s}(\pi _{<q}X)}&&& W_{q}Q_{s}s_{<q}X}
		$$
	is an isomorphism.  
	Furthermore, these maps
	induce a natural isomorphism between the following exact functors 
		$$\xymatrix{\stablehomotopy \ar@<1ex>[rr]^-{W_{q}Q_{s}} \ar@<-1ex>[rr]_-{W_{q}Q_{s}s_{<q}}&& \stablehomotopy}
		$$
\end{cor}
\begin{proof}
	Since $Q_{s}$, $W_{q}$ are both functors and
	$\pi _{<q}:id\rightarrow s_{<q}$  is a natural transformation
	(see theorem \ref{thm.3.1.motivictower});
	we have that 
	the maps $W_{q}Q_{s}(\pi _{<q}X)$ induce a natural transformation
	$W_{q}Q_{s}\rightarrow W_{q}Q_{s}s_{<q}$.  Therefore it suffices to see that for every
	$T$-spectrum $X$, the map $W_{q}Q_{s}(\pi _{<q}X)$ 
	is an isomorphism in $\stablehomotopy$.
	
	But proposition \ref{prop.3.2.Qs&QsSq==isoin-LqSH} implies that
	the map $Q_{s}(\pi _{<q}X)$ is an isomorphism in $\weightqstablehomotopy$.
	Therefore using proposition \ref{prop.3.2.Lq-exact-adjunctions}, we have that
	$W_{q}Q_{s}(\pi _{<q}X)$ is also an isomorphism in $\stablehomotopy$.
\end{proof}

\begin{lem}
		\label{lemma.3.2.IQJs<qX==>Lq-local}
	Fix $q\in \mathbb Z$.  Then for every $T$-spectrum $X$, $IQ_{T}J(Q_{s}s_{<q}X)$
	is $\qleftlocmaps$-local in $\motivicTspectra$.
\end{lem}
\begin{proof}
	Proposition \ref{prop.3.2.Lq-local-objects-classification} implies 
	that it is enough to show that
	$IQ_{T}J(Q_{s}s_{<q}X)$ satisfies the following properties:
		\begin{enumerate}
			\item	\label{lemma.3.2.IQJs<qX==>Lq-local.a}$IQ_{T}J(Q_{s}s_{<q}X)$ is 
						fibrant in $\motivicTspectra$.
			\item	\label{lemma.3.2.IQJs<qX==>Lq-local.b}For every 
						$\generatorNRS \in C_{eff}^{q}$
							$$[\generatorNRS ,IQ_{T}J(Q_{s}s_{<q}X)]_{Spt}\cong 0$$
		\end{enumerate}
	
	The first condition is obvious since $IQ_{T}J$ is a fibrant
	replacement functor in $\motivicTspectra$.
	
	Fix $\generatorNRS \in C_{eff}^{q}$.  Using theorem \ref{thm.3.1.motivictower}(\ref{thm.3.1.motivictower.b})
	and the fact that $C_{eff}^{q}\subseteq \stablehomotopyeffq$, we have
	that 
		$$[\generatorNRS , s_{<q}X]_{Spt}\cong 0$$ 
	Therefore
		$$[\generatorNRS ,IQ_{T}J(Q_{s}s_{<q}X)]_{Spt}\cong [\generatorNRS , s_{<q}X]_{Spt}\cong 0
		$$
	for every $\generatorNRS \in C_{eff}^{q}$.  This takes care of the second condition and
	finishes the proof.
\end{proof}

\begin{prop}
		\label{prop.3.2.Qss<q-isoto--WqQss<q}
	Fix $q\in \mathbb Z$.  Then for every $T$-spectrum $X$ the natural map
		$$\xymatrix{Q_{s}s_{<q}X \ar[rr]^-{W_{q}^{Q_{s}s_{<q}X}}&& W_{q}Q_{s}s_{<q}X}
		$$
	is a weak equivalence in $\motivicTspectra$.
	Therefore, we have a natural isomorphism 
	between the following exact functors
		$$\xymatrix{\stablehomotopy \ar@<1ex>[rr]^-{Q_{s}s_{<q}} \ar@<-1ex>[rr]_-{W_{q}Q_{s}s_{<q}}&& 
								\stablehomotopy}
		$$
\end{prop}
\begin{proof}
	The naturality of the maps $W_{q}^{X}:X\rightarrow W_{q}X$ implies that we have an
	induced natural transformation of functors $Q_{s}s_{<q}\rightarrow W_{q}Q_{s}s_{<q}$.
	Hence, it is enough to show that for every $T$-spectrum $X$, $W_{q}^{Q_{s}s_{<q}X}$ is
	a weak equivalence in $\motivicTspectra$.
	
	Consider the following commutative diagram in $\motivicTspectra$:
		\begin{equation}
					\label{diagram.prop.3.2.Qss<q-isoto--WqQss<q}
			\begin{array}{c}
				\xymatrix{Q_{s}s_{<q}X \ar[r] \ar[d]_-{W_{q}^{Q_{s}s_{<q}X}}& 
									IQ_{T}J(Q_{s}s_{<q}X)\ar[d]^-{IQ_{T}J(W_{q}^{Q_{s}s_{<q}X})}\\
									W_{q}Q_{s}s_{<q}X \ar[r]& IQ_{T}J(W_{q}Q_{s}s_{<q}X)}
			\end{array}
		\end{equation}
	where the horizontal maps are weak equivalences in $\motivicTspectra$.  
	Hence, the two out of three property for weak equivalences implies
	that it is enough to show that $IQ_{T}J(W_{q}^{Q_{s}s_{<q}X})$ is a weak equivalence
	in $\motivicTspectra$.  
	
	By construction the map $W_{q}^{Q_{s}s_{<q}X}$ is a $\qleftlocmaps$-local equivalence,
	and since the horizontal maps in diagram (\ref{diagram.prop.3.2.Qss<q-isoto--WqQss<q}) are weak equivalences
	in $\motivicTspectra$, it follows from
	\cite[proposition 3.1.5]{MR1944041} that these horizontal maps are also $\qleftlocmaps$-local equivalences.
	Therefore, the two out of three property for $\qleftlocmaps$-local equivalences implies that
	$IQ_{T}J(W_{q}^{Q_{s}s_{<q}X})$ is a $\qleftlocmaps$-local equivalence.
	
	Now lemma \ref{lemma.3.2.IQJs<qX==>Lq-local} implies that $IQ_{T}J(Q_{s}s_{<q}X)$ is
	$\qleftlocmaps$-local.  On the other hand, since the map
		$$\xymatrix{W_{q}Q_{s}s_{<q}X \ar[r]& IQ_{T}J(W_{q}Q_{s}s_{<q}X)}
		$$
	is a weak equivalence in $\motivicTspectra$, $W_{q}Q_{s}s_{<q}X$ is by construction
	$\qleftlocmaps$-local, and $IQ_{T}J(W_{q}Q_{s}s_{<q}X)$, $W_{q}Q_{s}s_{<q}X$ are both fibrant in $\motivicTspectra$;
	it follows from \cite[lemma 3.2.1]{MR1944041} that
	$IQ_{T}J(W_{q}Q_{s}s_{<q}X)$ is also $\qleftlocmaps$-local.
	
	Finally we have a $\qleftlocmaps$-local equivalence
		$$\xymatrix{IQ_{T}J(Q_{s}s_{<q}X) \ar[rrr]^-{IQ_{T}J(W_{q}^{Q_{s}s_{<q}X})}&&& IQ_{T}J(W_{q}Q_{s}s_{<q}X)}
		$$
	where the domain and the codomain are both $\qleftlocmaps$-local.
	Then theorem 3.2.13 in \cite{MR1944041} implies that
	$IQ_{T}J(W_{q}^{Q_{s}s_{<q}X})$ is a weak equivalence in $\motivicTspectra$.
	This finishes the proof.
\end{proof}

\begin{thm}
		\label{thm.3.2.Lq-models-s<q}
	Fix $q\in \mathbb Z$.  Then for every $T$-spectrum $X$, we have the following diagram in 
	$\stablehomotopy$:
		\begin{equation}
					\label{diagram.3.2.s<q-lifting}
			\begin{array}{c}
				\xymatrix{s_{<q}X &&& Q_{s}s_{<q}X \ar[lll]_-{Q_{s}^{s_{<q}X}}^-{\cong} 
									\ar[dd]^-{W_{q}^{Q_{s}s_{<q}X}}_-{\cong} &&&\\
									&&&&&&\\ 
									&&& W_{q}Q_{s}s_{<q}X &&& W_{q}Q_{s}X \ar[lll]_-{W_{q}Q_{s}(\pi _{<q})}^-{\cong}}
			\end{array}
		\end{equation}
	where all the maps are isomorphisms in $\stablehomotopy$.
	This diagram induces a
	natural isomorphism  
	between the following exact functors:
		$$\xymatrix{\stablehomotopy \ar@<1ex>[rr]^-{s_{<q}}  
								\ar@<-1ex>[rr]_-{W_{q}Q_{s}} && \stablehomotopy}
		$$
\end{thm}
\begin{proof}
	Since
	$Q_{s}$ is a cofibrant replacement functor in $\motivicTspectra$, it is clear that $Q_{s}^{s_{<q}X}$ becomes
	an isomorphism in the associated homotopy category $\stablehomotopy$.
	
	The fact that $W_{q}^{Q_{s}s_{<q}X}$
	is an isomorphism in $\stablehomotopy$ follows from
	proposition \ref{prop.3.2.Qss<q-isoto--WqQss<q}.  Finally, corollary \ref{cor.3.2.WqQs===WqQss<q}
	implies that $W_{q}Q_{s}(\pi _{<q})$ is also an isomorphism in $\stablehomotopy$.
	This shows that all the maps in the diagram (\ref{diagram.3.2.s<q-lifting}) are isomorphisms
	in $\stablehomotopy$, therefore for every $T$-spectrum $X$ we can define the
	following composition in $\stablehomotopy$
		\begin{equation}
					\label{diagram.3.2.s<q-lifting.b}
			\begin{array}{c}
				\xymatrix{s_{<q}X \ar[rrr]^-{(Q_{s}^{s_{<q}X})^{-1}}_-{\cong}&&& Q_{s}s_{<q}X  
										\ar[dd]^-{W_{q}^{Q_{s}s_{<q}X}}_-{\cong} &&&\\
										&&&&&&\\ 
										&&& W_{q}Q_{s}s_{<q}X \ar[rrr]^-{(W_{q}Q_{s}(\pi _{<q}))^{-1}}_-{\cong}&&& W_{q}Q_{s}X}
			\end{array}
		\end{equation}
	which is an isomorphism.  The fact that
	$Q_{s}$ is a functorial cofibrant replacement in $\motivicTspectra$,	
	proposition
	\ref{prop.3.2.Qss<q-isoto--WqQss<q} and corollary \ref{cor.3.2.WqQs===WqQss<q},
	imply all together that the isomorphisms defined in diagram (\ref{diagram.3.2.s<q-lifting.b})
	induce a natural isomorphism of functors $s_{<q}\stackrel{\cong}{\rightarrow}W_{q}Q_{s}$.
	This finishes the proof.
\end{proof}

\begin{rmk}
		\label{rmk.3.2.lifting-s<q}
	Theorem \ref{thm.3.2.Lq-models-s<q} gives the desired lifting to the
	model category level
	for the functors $s_{<q}$ defined in theorem \ref{thm.3.1.motivictower}.
\end{rmk}

\begin{prop}
		\label{prop.Lq+1-->Lq}
	For every $q\in \mathbb Z$, we have the following
	commutative diagram of left Quillen functors:
		$$\xymatrix{& \motivicTspectra \ar[dl]_-{id} \ar[dr]^-{id}&\\
								L_{<q+1}\motivicTspectra \ar[rr]_-{id}&& \weightqTspectra}
		$$
\end{prop}
\begin{proof}
	Since $\weightqTspectra$ and $L_{<q+1}\motivicTspectra$ are both
	left Bousfield localizations for $\motivicTspectra$, we have that
	the identity functors:
		$$\xymatrix@R=.5pt{id:\motivicTspectra \ar[r]& \weightqTspectra \\
								id:\motivicTspectra \ar[r]& L_{<q+1}\motivicTspectra}
		$$
	are both left Quillen functors.  Hence, it suffices to show that
		$$\xymatrix{id:L_{<q+1}\motivicTspectra \ar[r]& \weightqTspectra}
		$$
	is a left Quillen functor.  Using the universal property for left
	Bousfield localizations (see definition \ref{def-1.leftlocmodcats}), 
	we have that it is enough to check
	that if $f:X\rightarrow Y$ is a $L(<q+1)$-local equivalence then
	$Q_{s}(f):Q_{s}X\rightarrow Q_{s}Y$ is a $\qleftlocmaps$-local equivalence.
	
	But theorem 3.1.6(c) in \cite{MR1944041} implies that this last condition is equivalent
	to the following one: Let $Z$ be an arbitrary $\qleftlocmaps$-local $T$-spectrum,
	then $Z$ is also $L(<q+1)$-local.
	Finally, this last condition follows immediately from corollary
	\ref{cor.3.2.m>n===>Ln-local=>Lm-local}.
\end{proof}

\begin{cor}
		\label{cor.3.2.q+1-->qstableadj}
	For every $q\in \mathbb Z$, we have the following
	adjunction
		$$\xymatrix{(Q_{s},W_{q},\varphi ):\weightqplusonestablehomotopy \ar[r]& 
								\weightqstablehomotopy}
		$$
	of exact functors between triangulated categories.
\end{cor}
\begin{proof}
	Proposition \ref{prop.Lq+1-->Lq} implies that $id:L_{<q+1}\motivicTspectra \rightarrow \weightqTspectra$
	is a left Quillen functor.  Therefore we get the following
	adjunction at the level of the associated homotopy categories
		$$\xymatrix{(Q_{s},W_{q},\varphi):L_{<q+1}\stablehomotopy \ar[r]& \weightqstablehomotopy}
		$$
	Now proposition 6.4.1 in \cite{MR1650134} implies that
	$Q_{s}$ maps cofibre sequences in $L_{<q+1}\stablehomotopy$ to cofibre sequences in
	$\weightqstablehomotopy$.
	Therefore 
	using proposition 7.1.12 in \cite{MR1650134} we have
	that $Q_{s}$ and $W_{q}$ are both exact functors between triangulated categories.
\end{proof}

\begin{thm}
		\label{thm.3.2.motivic-tower}
	We have the following tower of left Quillen functors:
		\begin{equation}
					\label{diagram.3.2.motivictower-modelcatlevel}
			\begin{array}{c}
				\xymatrix{& \vdots \ar[d]^-{id} \\
									& L_{<q+1}\motivicTspectra \ar[d]^-{id} \\
									\motivicTspectra \ar[ur]^-{id} \ar[r]^-{id} \ar[dr]_-{id}& \weightqTspectra \ar[d]^-{id} \\
									& L_{<q-1}\motivicTspectra \ar[d]^-{id} &\\
									& \vdots}
			\end{array}
		\end{equation}
	together with the corresponding tower of associated homotopy categories:
		\begin{equation}
					\label{diagram.3.2.motivictower-homotopylevel}
			\begin{array}{c}
		  	\xymatrix{&&& \vdots \ar@<-1ex>[d]_-{Q_{s}} \\
									&&& L_{<q+1}\stablehomotopy \ar@<-1ex>[d]_-{Q_{s}} \ar@<-1ex>[u]_-{W_{q+1}} 
									\ar[dlll]|-{W_{q+1}}\\
									\stablehomotopy \ar@<2ex>[urrr]|-{Q_{s}} \ar@<1ex>[rrr]|-{Q_{s}} \ar[drrr]|-{Q_{s}}&&&
									\weightqstablehomotopy \ar@<-1ex>[d]_-{Q_{s}} \ar@<-1ex>[u]_-{W_{q}}
									\ar@<1ex>[lll]|-{W_{q}}\\
									&&& L_{<q-1}\stablehomotopy \ar@<-1ex>[d]_-{Q_{s}} \ar@<-1ex>[u]_-{W_{q-1}}
									\ar@<2ex>[ulll]|-{W_{q-1}}\\
									&&& \vdots \ar@<-1ex>[u]_-{W_{q-2}}}
			\end{array}
		\end{equation}
	Furthermore, the tower (\ref{diagram.3.2.motivictower-homotopylevel})
	satisfies the following properties:
		\begin{enumerate}
			\item	All the categories are triangulated.
			\item	All the functors are exact.
			\item	$Q_{s}$ is a left adjoint for all the functors $W_{q}$. 
		\end{enumerate}
\end{thm}
\begin{proof}
	Follows immediately from propositions
	\ref{prop.3.2.Lq-exact-adjunctions}, \ref{prop.Lq+1-->Lq} and 
	corollary \ref{cor.3.2.q+1-->qstableadj}.
\end{proof}

\begin{rmk}
		\label{rmk.3.2.Lq-cellular-Rq-notcellular}
	The great technical advantage of the categories $\weightqTspectra$
	over the categories $\qconnectedTspectra$ is the
	fact that $\weightqTspectra$ are always cellular, whereas it is not clear
	if $\qconnectedTspectra$
	satisfy the cellularity property.  Therefore we still can apply
	Hirschhorn's localization technology to the categories $\weightqTspectra$.
	This will be the final step in our approach to get the desired
	lifting for the functors $s_{q}$ (see theorem \ref{thm.3.1.slicefiltration})
	to the model category level.
\end{rmk}

\begin{defi}
		\label{def.3.2.Sq-colocal-generators}
	For every $q\in \mathbb Z$, we consider the following set
	of $T$-spectra
		$$\qslicegenerators =\{ \generatorNRS \in C | s-n=q \}\subseteq C_{eff}^{q}
		$$
	(see proposition \ref{prop.3.1.stablehomotopy=>compactly-generated} and 
	definition \ref{def.3.1.desusp-stable-homotopy-eff}).
\end{defi}

\begin{thm}
		\label{thm.3.2.Sq-modelcats}
	Fix $q\in \mathbb Z$.  Then the right Bousfield
	localization of the model category $\weightqplusoneTspectra$
	with respect to the $\qslicegenerators$-colocal equivalences
	exists.  This new model structure will be called \emph{$q$-slice
	motivic stable}.
	$\qsliceTspectra$ will denote the category of $T$-spectra
	equipped with the $q$-slice motivic stable model structure, and 
	$\qslicestablehomotopy$ will denote its associated homotopy category.
	Furthermore the $q$-slice motivic stable model structure is  right proper and simplicial.
\end{thm}
\begin{proof}
	Theorem \ref{thm.3.2.Lqmodelstructures} implies that $\weightqplusoneTspectra$ is a cellular
	and simplicial model category.  On the other hand, corollary
	\ref{cor.3.2.Lq=>rightproper} implies that $\weightqplusoneTspectra$ is right proper.
	Therefore we can apply theorem 5.1.1 in \cite{MR1944041} to construct
	the right Bousfield localization
	of $\weightqplusoneTspectra$ with respect to the $\qslicegenerators$-colocal equivalences.
	Using \cite[theorem 5.1.1]{MR1944041} again, we have that $\qsliceTspectra$ is a right proper
	and simplicial model category.
\end{proof}

\begin{defi}
		\label{def.3.2.cofibrantreplacement-SqTspectra}
	Fix $q\in \mathbb Z$.  Let $P_{q}$ denote a
	functorial cofibrant replacement functor in $\qsliceTspectra$
	such that for every $T$-spectrum $X$, the natural map
		$$\xymatrix{P_{q}X \ar[r]^-{P_{q}^{X}} & X}
		$$
	is a trivial fibration in $\qsliceTspectra$, and
	$P_{q}X$ is a $\qslicegenerators$-colocal $T$-spectrum
	in $\weightqplusoneTspectra$.
\end{defi}

\begin{prop}
		\label{prop.3.2.Wq+1-fibrant-replacement-all-Sq}
	Fix $q\in \mathbb Z$.  Then $W_{q+1}$ is also a fibrant
	replacement functor in $\qsliceTspectra$
	(see definition \ref{def.3.2.stable-weightq-replacementfunctors}), and for every $T$-spectrum $X$
	the natural map
		$$\xymatrix{X \ar[rr]^-{W_{q+1}^{X}}&& W_{q+1}X}
		$$
	is a trivial cofibration in $\qsliceTspectra$.
\end{prop}
\begin{proof}
	Since $\qsliceTspectra$ is the right Bousfield localization of
	$\weightqplusoneTspectra$ with respect to the $\qslicegenerators$-colocal equivalences, by construction
	we have that the fibrations and the trivial cofibrations are indentical in
	$\qsliceTspectra$ and $\weightqplusoneTspectra$ respectively.  This implies that for every
	$T$-spectrum $X$, $W_{q+1}X$ is fibrant in $\qsliceTspectra$, and 
	we also have that
	the natural map:
		$$\xymatrix{X\ar[rr]^{W_{q+1}^{X}}&& W_{q+1}X}
		$$
	is a trivial cofibration in $\qsliceTspectra$.  Hence $W_{q+1}$ is also a fibrant replacement
	functor for $\qsliceTspectra$.
\end{proof}

\begin{prop}
		\label{prop.3.2.classifying-Sq-colocal-equivs}
	Fix $q\in \mathbb Z$ and let $f:X\rightarrow Y$ be a map of $T$-spectra.
	Then $f$ is a $\qslicegenerators$-colocal equivalence
	in $\weightqplusoneTspectra$ if and only if
	for every $\generatorNRS \in \qslicegenerators$
	the induced map
		$$\xymatrix{[\generatorNRS ,W_{q+1}X]_{Spt} \ar[rr]^-{(W_{q+1}f)_{\ast}}_-{\cong}&&
			[\generatorNRS ,W_{q+1}Y]_{Spt}}
		$$
	is an isomorphism of abelian groups.
\end{prop}
\begin{proof}
	($\Rightarrow$): Assume that $f$ is a $\qslicegenerators$-colocal equivalence.
	All the compact generators $\generatorNRS$ are cofibrant
	in $\weightqplusoneTspectra$, since they are cofibrant in $\motivicTspectra$, 
	and the cofibrations are exactly
	the same in both model structures. 
	
	Therefore we have that $f$ is a
	$\qslicegenerators$-colocal equivalence if and only if
	for every $\generatorNRS \in S(q)$
	the following maps are weak equivalences of simplicial sets:
		$$\xymatrix{Map(\generatorNRS ,W_{q+1}X) \ar[rr]^-{(W_{q+1}f)_{\ast}}&& 
			Map(\generatorNRS ,W_{q+1}Y)}
		$$
	Since
	$\weightqplusoneTspectra$
	is a simplicial model category,
	we have that $Map(\generatorNRS ,W_{q+1}X)$ and $Map(\generatorNRS ,W_{q+1}Y)$ are both Kan complexes. 
	Now proposition \ref{prop.3.2.Lq-local-objects-classification}(\ref{prop.3.2.Lq-local-objects-classification.a})
	implies that $W_{q+1}X$, $W_{q+1}Y$ are both fibrant in $\motivicTspectra$, therefore
	since $\motivicTspectra$ is a simplicial model category
	we get the following commutative diagram where the top row and all the vertical
	maps are isomorphisms of abelian groups:
		$$\xymatrix@C=-1.5pc{\pi_{0}Map(\generatorNRS ,W_{q+1}X) \ar[dr]^-{(W_{q+1}f)_{\ast}}_-{\cong} \ar[dd]_-{\cong}&\\ 
								& \pi_{0}Map(\generatorNRS ,W_{q+1}Y) \ar[dd]^-{\cong}\\
								[\generatorNRS ,W_{q+1}X]_{Spt} \ar[dr]_-{(W_{q+1}f)_{\ast}}&\\ 
								&[\generatorNRS ,W_{q+1}Y]_{Spt}}
		$$
	Therefore
		$$\xymatrix{[\generatorNRS ,W_{q+1}X]_{Spt} \ar[rr]^-{(W_{q+1}f)_{\ast}}_-{\cong}&& [\generatorNRS ,W_{q+1}Y]_{Spt}}
		$$
	is an isomorphism of abelian groups
	for every $\generatorNRS \in \qslicegenerators$, as we wanted.
	
	($\Leftarrow$):  Fix $\generatorNRS \in \qslicegenerators$.
	Let $\omega _{0}$, $\eta _{0}$ be the base points corresponding to
	$Map_{\ast}(\generatorNNRSS,W_{q+1}X)$ and $Map_{\ast}(\generatorNNRSS,W_{q+1}Y)$
	respectively.
	We need to show that the map:
		$$\xymatrix{Map(\generatorNRS ,W_{q+1}X) \ar[rr]^-{(W_{q+1}f)_{\ast}}&& Map(\generatorNRS ,W_{q+1}Y)}
		$$
	is a weak equivalence of simplicial sets.  Let 
		$$j:F_{n+1}(S^{r+1}\wedge \gm ^{s+1}\wedge U_{+})\rightarrow \generatorNRS$$ 
	be the adjoint to the identity map
		$$id:S^{r+1}\wedge \gm ^{s+1}\wedge U_{+}\rightarrow Ev_{n+1}\generatorNRS 
			=S^{r+1}\wedge \gm ^{s+1}\wedge U_{+}
		$$
	We know that $j$ is a weak equivalence in $\motivicTspectra$, therefore
	\cite[proposition 3.1.5]{MR1944041} 
	implies that $j$ is a $\qplusoneleftlocmaps$-local equivalence, i.e.
	a weak equivalence in $\weightqplusoneTspectra$.
	Now since $\generatorNRS$ and $F_{n+1}(S^{r+1}\wedge \gm ^{s+1}\wedge U_{+})$ are both cofibrant
	in $\weightqplusoneTspectra$,
	and $\weightqplusoneTspectra$ is a simplicial model category, 
	we can apply Ken Brown's lemma (see lemma \ref{lem1.1.KenBrown})
	to conclude that the horizontal maps in the following commutative diagram are weak
	equivalences of simplicial sets:
		$$\xymatrix{Map(\generatorNRS ,W_{q+1}X) 
								\ar[r]^-{j^{\ast}}\ar[d]_-{(W_{q+1}f)_{\ast}} & 
								Map(F_{n+1}(S^{r+1}\wedge \gm ^{s+1}\wedge U_{+}) ,W_{q+1}X) \ar[d]^-{(W_{q+1}f)_{\ast}}\\		
								Map(\generatorNRS ,W_{q+1}Y)
								\ar[r]_-{j^{\ast}}& Map(F_{n+1}(S^{r+1}\wedge \gm ^{s+1}\wedge U_{+}) ,W_{q+1}Y)}
		$$
	Hence by the two out of three property for weak equivalences, it is enough to show that
	the following induced map
		$$\xymatrix{Map(F_{n+1}(S^{r+1}\wedge \gm ^{s+1}\wedge U_{+}) ,W_{q+1}X) \ar[d]^-{(W_{q+1}f)_{\ast}}\\
								Map(F_{n+1}(S^{r+1}\wedge \gm ^{s+1}\wedge U_{+}) ,W_{q+1}Y)}
		$$
	is a weak equivalence of simplicial sets.
	
	On the other hand, since $\motivicTspectra$ is a pointed simplicial model category
	and  
	$W_{q+1}X$, $W_{q+1}Y$ are both fibrant in $\motivicTspectra$ by
	proposition 
	\ref{prop.3.2.Lq-local-objects-classification}(\ref{prop.3.2.Lq-local-objects-classification.a});
	we have that lemma 6.1.2 in \cite{MR1650134} together 
	with remark \ref{rmk.2.4.simp-simppointed-structures-Tspectra}(\ref{rmk.2.4.simp-simppointed-structures-Tspectra.b})
	imply that
	the following diagram is commutative for $k\geq 0$:
		$$\xymatrix@C=-4pc{\pi_{k, \omega _{0}}Map(\generatorNNRSS ,W_{q+1}X) \ar[dr]^-{(W_{q+1}f)_{\ast}}\ar@{=}[dd]&\\  
								& \pi_{k, \eta _{0}}Map(\generatorNNRSS ,W_{q+1}Y)\ar@{=}[dd]\\
								\pi_{k, \omega _{0}}Map_{\ast}(\generatorNNRSS ,W_{q+1}X) \ar[dr]^-{(W_{q+1}f)_{\ast}}\ar[dd]_-{\cong}&\\  
								& \pi_{k, \eta _{0}}Map_{\ast}(\generatorNNRSS ,W_{q+1}Y)\ar[dd]^-{\cong}\\
								[\generatorNNRSS \wedge S^{k}, W_{q+1}X]_{Spt} \ar[dr]^-{(W_{q+1}f)_{\ast}}\ar[dd]_-{\cong}&\\ 
								& [\generatorNNRSS \wedge S^{k}, W_{q+1}Y]_{Spt} \ar[dd]^-{\cong}\\
								[F_{n+1}(S^{k+r}\wedge \gm ^{s+1}\wedge U_{+}) , W_{q+1}X]_{Spt} \ar[dr]^-{(W_{q+1}f)_{\ast}}&\\ 
								& [F_{n+1}(S^{k+r}\wedge \gm ^{s+1}\wedge U_{+}) , W_{q+1}Y]_{Spt}}
		$$
	but by hypothesis we have that the bottom row is an isomorphism of abelian groups,
	since $F_{n+1}(S^{k+r}\wedge \gm ^{s+1}\wedge U_{+})$ is also in $\qslicegenerators$.
	Therefore all the maps in the top row are also isomorphisms.  Hence, the induced map
		$$\xymatrix{Map(\generatorNNRSS ,W_{q+1}X) \ar[d]^-{(W_{q+1}f)_{\ast}}\\ 
								Map(\generatorNNRSS ,W_{q+1}Y)}
		$$ 
	is a weak equivalence when it is restricted to the path component of $Map(\generatorNNRSS ,W_{q+1}X)$
	containing $\omega _{0}$.  This implies that the following induced map
		$$\xymatrix{Map_{\ast}(S^{1},Map_{\ast}(F_{n+1}(S^{r}\wedge \gm ^{s+1}\wedge U_{+}) ,W_{q+1}X)) 
								\ar[d]^-{(W_{q+1}f)_{\ast}}\\
								Map_{\ast}(S^{1},Map_{\ast}(F_{n+1}(S^{r}\wedge \gm ^{s+1}\wedge U_{+}) ,W_{q+1}Y))}
		$$
	is a weak equivalence since taking $S^{1}$-loops kills the path components that do not
	contain the base point.  
	
	Finally, since $\motivicTspectra$ is a simplicial model category we have that the rows in the following
	commutative diagram are isomorphisms:
		$$\xymatrix@C=-5pc{Map_{\ast}(S^{1},Map_{\ast}(F_{n+1}(S^{r}\wedge \gm ^{s+1}\wedge U_{+}) ,W_{q+1}X)) \ar[dr]^-{\cong} 
								\ar[dd]_-{(W_{q+1}f)_{\ast}}&\\
								& Map_{\ast}(F_{n+1}(S^{r}\wedge \gm ^{s+1}\wedge U_{+})\wedge S^{1} ,W_{q+1}X)
								 \ar[dd]^-{(W_{q+1}f)_{\ast}}\\
								Map_{\ast}(S^{1},Map_{\ast}(F_{n+1}(S^{r}\wedge \gm ^{s+1}\wedge U_{+}) ,W_{q+1}Y)) \ar[dr]_-{\cong}&\\
								& Map_{\ast}(F_{n+1}(S^{r}\wedge \gm ^{s+1}\wedge U_{+})\wedge S^{1} ,W_{q+1}Y)}
		$$
	Hence the two out of three property for weak equivalences implies that the right vertical map
	is a weak equivalence of simplicial sets.  
	But $F_{n+1}(S^{r}\wedge \gm ^{s+1}\wedge U_{+})\wedge S^{1}$
	is clearly isomorphic to $F_{n+1}(S^{r+1}\wedge \gm ^{s+1}\wedge U_{+})$, therefore
	the induced map
		$$\xymatrix{Map(F_{n+1}(S^{r+1}\wedge \gm ^{s+1}\wedge U_{+}),W_{q+1}X) \ar[d]^-{(W_{q+1}f)_{\ast}}\\ 
								Map(F_{n+1}(S^{r+1}\wedge \gm ^{s+1}\wedge U_{+}),W_{q+1}Y)}
		$$
	is a weak equivalence, as we wanted.
\end{proof}

\begin{cor}
		\label{cor.3.2.classifying-Sq-colocal-equivs.b}
	Fix $q\in \mathbb Z$ and let $f:X\rightarrow Y$ be a map of $T$-spectra.
	Then $f$ is a $\qslicegenerators$-colocal equivalence
	in $\weightqplusoneTspectra$ if and only if
		$$\xymatrix{W_{q+1}X \ar[rr]^-{W_{q+1}f}&& W_{q+1}Y}
		$$
	is a $C_{eff}^{q}$-colocal equivalence in $\motivicTspectra$.
\end{cor}
\begin{proof}
	($\Rightarrow$):  Assume that $f$ is a $\qslicegenerators$-colocal equivalence, and
	fix $\generatorNRS \in C_{eff}^{q}$.  By proposition \ref{prop.3.2.classif-Cqcolocalequivs} 
	it suffices to show that the induced map
		\begin{equation}
					\label{equation.3.2.classifying-Sq-colocal-equivs.b}
			\begin{array}{c}
				\xymatrix{[\generatorNRS, W_{q+1}X]_{Spt}\ar[d]_-{(W_{q+1}f)_{\ast}}\\ 
									[\generatorNRS, W_{q+1}Y]_{Spt}}
			\end{array}
		\end{equation}
	is an isomorphism of abelian groups.
	
	Since $\generatorNRS \in C_{eff}^{q}$, we have two possibilities:
		\begin{enumerate}
			\item \label{cond.a.cor.3.2.classifying-Sq-colocal-equivs.b} $s-n=q$, i.e. 
						$\generatorNRS \in \qslicegenerators$.
			\item	\label{cond.b.cor.3.2.classifying-Sq-colocal-equivs.b} $s-n\geq q+1$, i.e.
						$\generatorNRS \in C_{eff}^{q+1}$
		\end{enumerate}
	
	In case (\ref{cond.a.cor.3.2.classifying-Sq-colocal-equivs.b}), 
	proposition \ref{prop.3.2.classifying-Sq-colocal-equivs} implies that the induced map
	in diagram (\ref{equation.3.2.classifying-Sq-colocal-equivs.b}) is an isomorphism of abelian
	groups.
	
	On the other hand, in case (\ref{cond.b.cor.3.2.classifying-Sq-colocal-equivs.b}),
	we have by proposition \ref{prop.3.2.Lq-local-objects-classification}(\ref{prop.3.2.Lq-local-objects-classification.b}) 
	that 
		$$[\generatorNRS ,W_{q+1}X]_{Spt}\cong 0\cong [\generatorNRS ,W_{q+1}Y]_{Spt}$$
	since by construction $W_{q+1}X$ and $W_{q+1}Y$ are both $\qplusoneleftlocmaps$-local $T$-spectra.
	Hence the induced map in diagram (\ref{equation.3.2.classifying-Sq-colocal-equivs.b}) is also
	an isomorphism of abelian groups in this case, as we wanted.
	
	($\Leftarrow$):  Assume that $W_{q+1}f$ is a $C_{eff}^{q}$-colocal equivalence in $\motivicTspectra$, and fix
	$\generatorNRS \in \qslicegenerators$.
	
	Since $\qslicegenerators \subseteq C_{eff}^{q}$, it follows from 
	proposition \ref{prop.3.2.classif-Cqcolocalequivs}
	that the induced map
		$$\xymatrix{[\generatorNRS ,W_{q+1}X]_{Spt} \ar[rr]^-{(W_{q+1}f)_{\ast}}&& 
								[\generatorNRS ,W_{q+1}Y]_{Spt}}
		$$
	is an isomorphism of abelian groups.
	Therefore, proposition \ref{prop.3.2.classifying-Sq-colocal-equivs} implies
	that $f$ is a $\qslicegenerators$-colocal equivalence in $\weightqplusoneTspectra$.  This finishes the proof.
\end{proof}

\begin{lem}
		\label{lem.3.2.Sq-colocalequivs-stable-S1desuspension}
	Fix $q\in \mathbb Z$ and let $f:X\rightarrow Y$ be a map of $T$-spectra, then $f$ is a
	$\qslicegenerators$-colocal equivalence in $\weightqplusoneTspectra$ if and only if
		$$\xymatrix{\Omega _{S^{1}} W_{q+1}(f):\Omega _{S^{1}}W_{q+1}X \ar[r]& \Omega _{S^{1}}W_{q+1}Y}
		$$ 
	is a $\qslicegenerators$-colocal equivalence in $\weightqplusoneTspectra$.
\end{lem}
\begin{proof}
	Assume that $f$ is a $\qslicegenerators$-colocal equivalence.
	We need to show that $\Omega _{S^{1}}W_{q+1}(f)$ is a
	$\qslicegenerators$-colocal equivalence in $\weightqplusoneTspectra$.
	
	Fix $\generatorNRS \in \qslicegenerators$.
	Corollary \ref{cor.3.2.S1loops-preserves-Lqlocal} implies that 
	$\Omega _{S^{1}}W_{q+1}X$ and $\Omega _{S^{1}}W_{q+1}Y$
	are both $\qplusoneleftlocmaps$-local; and proposition
	\ref{prop.3.2.Lq-local-objects-classification}(\ref{prop.3.2.Lq-local-objects-classification.a})
	implies that $\Omega _{S^{1}}W_{q+1}X$
	and $\Omega _{S^{1}}W_{q+1}Y$ are both fibrant in $\motivicTspectra$.  Therefore using the fact that
	$\motivicTspectra$
	is a simplicial model category, we get the following commutative diagram:
		$$\xymatrix@C=-2pc{[\generatorNRS ,\Omega _{S^{1}}W_{q+1}X]_{Spt} 
								\ar[dr]^-{(\Omega _{S^{1}}W_{q+1}f)_{\ast}} \ar[dd]_-{\cong}&\\ 
								& [\generatorNRS ,\Omega _{S^{1}}W_{q+1}Y]_{Spt} \ar[dd]^-{\cong} \\
								[\generatorNRS \wedge S^{1},W_{q+1}X]_{Spt} \ar[dr]^-{(W_{q+1}f)_{\ast}} \ar[dd]_-{\cong}&\\ 
								& [\generatorNRS \wedge S^{1},W_{q+1}Y]_{Spt} \ar[dd]^-{\cong}\\
								[F_{n}(S^{r+1}\wedge \gm ^{s}\wedge U_{+}), W_{q+1}X]_{Spt} \ar[dr]_-{(W_{q+1}f)_{\ast}}^-{\cong}&\\
								& [F_{n}(S^{r+1}\wedge \gm ^{s}\wedge U_{+}), W_{q+1}Y]_{Spt}}
		$$
	but using proposition \ref{prop.3.2.classifying-Sq-colocal-equivs} and the fact that
	$f$ is a $\qslicegenerators$-colocal equivalence, we have that the
	bottom row is an isomorphism, therefore the top row is also an isomorphism.
	Hence, the induced map:
		$$\xymatrix{[\generatorNRS ,\Omega _{S^{1}}W_{q+1}X]_{Spt} \ar[d]^-{(\Omega _{S^{1}}W_{q+1}f)_{\ast}}_-{\cong}\\
									[\generatorNRS ,\Omega _{S^{1}}W_{q+1}Y]_{Spt}}
		$$
	is an isomorphism of abelian groups for every $\generatorNRS \in \qslicegenerators$.
	Finally,
	using proposition \ref{prop.3.2.classifying-Sq-colocal-equivs}
	again, together with the fact that $\Omega _{S^{1}}W_{q+1}X$ and $\Omega _{S^{1}}W_{q+1}Y$ are
	both $\qplusoneleftlocmaps$-local $T$-spectra;
	we have that $\Omega _{S^{1}}W_{q+1}(f)$ is a
	$\qslicegenerators$-colocal equivalence in $\weightqplusoneTspectra$, as we wanted. 
	
	Conversely, assume that $\Omega _{S^{1}}W_{q+1}(f)$ is a
	$\qslicegenerators$-colocal equivalence in $\weightqplusoneTspectra$, and
	fix $\generatorNRS \in \qslicegenerators$.
	Corollary \ref{cor.3.2.S1loops-preserves-Lqlocal} implies that 
	$\Omega _{S^{1}}W_{q+1}X$ and $\Omega _{S^{1}}W_{q+1}Y$
	are both $\qplusoneleftlocmaps$-local; and proposition
	\ref{prop.3.2.Lq-local-objects-classification}(\ref{prop.3.2.Lq-local-objects-classification.a})
	implies that $\Omega _{S^{1}}W_{q+1}X$
	and $\Omega _{S^{1}}W_{q+1}Y$ are both fibrant in $\motivicTspectra$.  Therefore using the fact that
	$\motivicTspectra$
	is a simplicial model category, we get the following commutative diagram:
		$$\xymatrix@C=-3pc{[F_{n+1}(S^{r}\wedge \gm ^{s+1}\wedge U_{+}) ,\Omega _{S^{1}}W_{q+1}X]_{Spt} 
								\ar[dr]^-{(\Omega _{S^{1}}W_{q+1}f)_{\ast}}_-{\cong} 
								\ar[dd]_-{\cong}&\\ 
								& [F_{n+1}(S^{r}\wedge \gm ^{s+1}\wedge U_{+}) ,\Omega _{S^{1}}W_{q+1}Y]_{Spt} \ar[dd]^-{\cong}\\
								[F_{n+1}(S^{r}\wedge \gm ^{s+1}\wedge U_{+})\wedge S^{1} ,W_{q+1}X]_{Spt} 
								\ar[dr]^-{(W_{q+1}f)_{\ast}} \ar[dd]_-{\cong}&\\ 
								& [F_{n+1}(S^{r}\wedge \gm ^{s+1}\wedge U_{+})\wedge S^{1}, W_{q+1}Y]_{Spt} \ar[dd]^-{\cong}\\
								[F_{n+1}(S^{r+1}\wedge \gm ^{s+1}\wedge U_{+}),W_{q+1}X]_{Spt} 
								\ar[dr]^-{(W_{q+1}f)_{\ast}} &\\ 
								& [F_{n+1}(S^{r+1}\wedge \gm ^{s+1}\wedge U_{+}),W_{q+1}Y]_{Spt} \\
								[\generatorNRS ,W_{q+1}X]_{Spt} \ar[dr]_-{(W_{q+1}f)_{\ast}} \ar[uu]^-{\cong}&\\ 
								& [\generatorNRS ,W_{q+1}Y]_{Spt} \ar[uu]_-{\cong}}
		$$
	Since $\Omega _{S^{1}}W_{q+1}f$ is a $\qslicegenerators$-colocal equivalence,
	we have that
	proposition \ref{prop.3.2.classifying-Sq-colocal-equivs}
	together with the fact that $\Omega _{S^{1}}W_{q+1}X$ and $\Omega _{S_{1}}W_{q+1}Y$
	are both $\qplusoneleftlocmaps$-local
	imply that the top row in the diagram above is
	an isomorphism;
	therefore the bottom row is also an isomorphism.
	Thus, the induced map:
		$$\xymatrix{[\generatorNRS ,W_{q+1}X]_{Spt} \ar[rr]^-{(W_{q+1}f)_{\ast}}_-{\cong} &&
									[\generatorNRS ,W_{q+1}Y]_{Spt}}
		$$
	is an isomorphism of abelian groups for every $\generatorNRS \in \qslicegenerators$.
	Now using proposition \ref{prop.3.2.classifying-Sq-colocal-equivs} again,
	we have that $f$ is a $\qslicegenerators$-colocal equivalence.
	This finishes the proof.
\end{proof}

\begin{cor}
		\label{cor.3.2.S1-Quillen-equiv-SqTspectra}
	For every $q\in \mathbb Z$, the adjunction
		$$\xymatrix{(-\wedge S^{1},\Omega _{S^{1}},\varphi):\qsliceTspectra \ar[r] & \qsliceTspectra}
		$$
	is a Quillen equivalence.
\end{cor}
\begin{proof}
	Using corollary 1.3.16 in \cite{MR1650134} and
	proposition \ref{prop.3.2.Wq+1-fibrant-replacement-all-Sq} we have that it suffices to
	verify the following two conditions:
		\begin{enumerate}
			\item	\label{cor.3.2.S1-Quillen-equiv-SqTspectra.a}For every cofibrant object $X$
						in $\qsliceTspectra$, the following composition
							$$\xymatrix{X\ar[r]^-{\eta _{X}}& \Omega _{S^{1}}(X\wedge S^{1}) \ar[rr]^-{\Omega _{S^{1}}W_{q+1}^{X\wedge S^{1}}}&& 
								\Omega _{S^{1}}W_{q+1}(X\wedge S^{1})}
							$$
						is a $\qslicegenerators$-colocal equivalence.
			\item	\label{cor.3.2.S1-Quillen-equiv-SqTspectra.b}$\Omega _{S^{1}}$ reflects $\qslicegenerators$-colocal equivalences
						between fibrant objects in $\qsliceTspectra$.
		\end{enumerate}
		
	(\ref{cor.3.2.S1-Quillen-equiv-SqTspectra.a}):  By construction $\qsliceTspectra$ is
	a right Bousfield localization of $\weightqplusoneTspectra$, therefore the identity
	functor 
		$$\xymatrix{id:\qsliceTspectra \ar[r]& \weightqplusoneTspectra}
		$$ 
	is a left Quillen functor.
	Thus $X$ is also cofibrant in $\weightqplusoneTspectra$.
	Since the adjunction $(-\wedge S^{1},\Omega _{S^{1}},\varphi)$
	is a Quillen equivalence on $\weightqplusoneTspectra$, 
	\cite[proposition 1.3.13(b)]{MR1650134} implies that the following composition
	is a weak equivalence in $\weightqplusoneTspectra$:
		$$\xymatrix{X\ar[r]^-{\eta _{X}}& \Omega _{S^{1}}(X\wedge S^{1}) \ar[rr]^-{\Omega _{S^{1}}W_{q+1}^{X\wedge S^{1}}}&& 
								\Omega _{S^{1}}W_{q+1}(X\wedge S^{1})}
		$$
	Hence using \cite[proposition 3.1.5]{MR1944041} it follows that the composition above
	is a $\qslicegenerators$-colocal equivalence.
	
	(\ref{cor.3.2.S1-Quillen-equiv-SqTspectra.b}):  This follows immediately
	from proposition \ref{prop.3.2.Wq+1-fibrant-replacement-all-Sq} and 
	lemma \ref{lem.3.2.Sq-colocalequivs-stable-S1desuspension}.
\end{proof}

\begin{rmk}
		\label{rmk.3.2.Tsuspfunctor-doesnotdescend-to-Sq}
	We have a situation similar to the one
	described in remarks
	\ref{rmk.3.3.smashT-not-descending} and \ref{rmk.3.2.Tsuspfunctor-doesnotdescend-to-Lq}
	for the model categories
	$\qconnectedTspectra$ and $\weightqTspectra$; i.e. although
	the adjunction $(\Tsuspfunctor ,\Tloops ,\varphi)$ is
	a Quillen equivalence on $\motivicTspectra$,
	it does not descend even to a Quillen
	adjunction on the $q$-slice motivic stable
	model category
	$\qsliceTspectra$.
\end{rmk}

\begin{cor}
		\label{cor.3.2.Sqstablehomotopy==triangulatedcat}
	For every $q\in \mathbb Z$, 
	$\qslicestablehomotopy$ has the
	structure of a triangulated category.
\end{cor}
\begin{proof}
	Theorem \ref{thm.3.2.Sq-modelcats} implies in particular
	that $\qsliceTspectra$ is a pointed simplicial model category,
	and corollary \ref{cor.3.2.S1-Quillen-equiv-SqTspectra} implies that
	the adjunction 
		$$(-\wedge S^{1},\Omega _{S^{1}},\varphi):\qsliceTspectra \rightarrow \qsliceTspectra$$
	is a Quillen equivalence.  Therefore
	the result follows from the work of Quillen in 
	\cite[sections I.2 and I.3]{MR0223432} and the work of 
	Hovey in \cite[chapters VI and VII]{MR1650134}.
\end{proof}

\begin{prop}
		\label{prop.3.2.adjunctions--Sq==>Lq}
	For every $q\in \mathbb Z$ we have
	the following adjunction
		$$\xymatrix{(P_{q},W_{q+1},\varphi ):\qslicestablehomotopy \ar[r]& \weightqplusonestablehomotopy}
		$$
	of exact functors between triangulated categories.
\end{prop}
\begin{proof}
	Since $\qsliceTspectra$ is the right Bousfield localization of
	$\weightqplusoneTspectra$ with respect to the $\qslicegenerators$-colocal equivalences, we have that the
	identity functor $id:\qsliceTspectra \rightarrow \weightqplusoneTspectra$
	is a left Quillen functor.  Therefore we get
	the following adjunction at the level of the associated homotopy categories:
		$$\xymatrix{(P_{q}, W_{q+1}, \varphi ):\qslicestablehomotopy \ar[r]& \weightqplusonestablehomotopy}
		$$
	
	Now proposition 6.4.1 in \cite{MR1650134} implies that
	$P_{q}$ maps cofibre sequences in $\qslicestablehomotopy$ to cofibre sequences in
	$\weightqplusonestablehomotopy$.
	Therefore 
	using proposition 7.1.12 in \cite{MR1650134} we have
	that $P_{q}$ and $W_{q+1}$ are both exact functors between triangulated categories.
\end{proof}

\begin{prop}
		\label{prop.3.2.q-connected--->q-slice===leftQuilllenfunctor}
	Fix $q\in \mathbb Z$.  Then the identity functor
		$$\xymatrix{id:\qsliceTspectra \ar[r]& \qconnectedTspectra}
		$$
	is a right Quillen functor.
\end{prop}
\begin{proof}
	Consider the following diagram of right Quillen functors
		$$\xymatrix{\weightqplusoneTspectra \ar[r]^-{id} \ar[d]_-{id}& \motivicTspectra \ar[r]^-{id}& \qconnectedTspectra \\
								\qsliceTspectra \ar@{-->}[urr]_-{id}&&}
		$$
	By the universal property of right Bousfield localizations 
	(see definition \ref{def-1.rightlocmodcats})
	it suffices to check that if $f:X\rightarrow Y$ is a $S(q)$-colocal equivalence in
	$\weightqplusoneTspectra$, then $W_{q+1}f:W_{q+1}X\rightarrow W_{q+1}Y$ is a
	$C_{eff}^{q}$-colocal equivalence in $\motivicTspectra$.
	But this follows immediately from 
	corollary \ref{cor.3.2.classifying-Sq-colocal-equivs.b}.
\end{proof}

\begin{cor}
		\label{cor.3.2.adjunctions--Rq==>Sq}
	For every $q\in \mathbb Z$ we have
	the following adjunction
		$$\xymatrix{(C_{q},W_{q+1},\varphi ):\qconnectedstablehomotopy \ar[r]& \qslicestablehomotopy}
		$$
	of exact functors between triangulated categories.
\end{cor}
\begin{proof}
	By proposition \ref{prop.3.2.q-connected--->q-slice===leftQuilllenfunctor} the identity functor
	$id:\qconnectedTspectra \rightarrow \qsliceTspectra$ is a left Quillen
	functor.  Therefore we get
	the following adjunction at the level of the associated homotopy categories:
		$$\xymatrix{(C_{q}, W_{q+1}, \varphi ):\qconnectedstablehomotopy \ar[r]& \qslicestablehomotopy}
		$$
	
	Now proposition 6.4.1 in \cite{MR1650134} implies that
	$C_{q}$ maps cofibre sequences in $\qconnectedstablehomotopy$ to cofibre sequences in
	$\qslicestablehomotopy$.
	Therefore 
	using proposition 7.1.12 in \cite{MR1650134} we have
	that $C_{q}$ and $W_{q+1}$ are both exact functors between triangulated categories.
\end{proof}

\begin{lem}
		\label{lem.3.2.A-Sqcofibrant====>ALq+1trivially-cofibrant}
	Fix $q\in \mathbb Z$, and let $A$ be a cofibrant $T$-spectrum
	in $\qsliceTspectra$.  Then the map $\ast \rightarrow A$
	is a trivial cofibration in $\weightqTspectra$.
\end{lem}
\begin{proof}
	Let $Z$ be an arbitrary $\qleftlocmaps$-local $T$-spectrum in $\motivicTspectra$.
	We claim that the map $Z\rightarrow \ast$ is a trivial fibration in $\qsliceTspectra$.
	In effect, using corollary \ref{cor.3.2.m>n===>Ln-local=>Lm-local} 
	we have that $Z$ is $\qplusoneleftlocmaps$-local in $\motivicTspectra$, i.e. a fibrant
	object in $\weightqplusoneTspectra$.  By construction $\qsliceTspectra$ is a right Bousfield
	localization of $\weightqplusoneTspectra$, hence $Z$ is also fibrant in $\qsliceTspectra$.
	Then by proposition \ref{prop.3.2.classifying-Sq-colocal-equivs} 
	it suffices to show that for every
	$\generatorNRS \in S(q)$ (i.e. $s-n=q$):
		$$\xymatrix{0\cong [\generatorNRS ,Z]_{Spt}}
		$$
	But this follows immediately from proposition \ref{prop.3.2.Lq-local-objects-classification}, 
	since $Z$ is $\qleftlocmaps$-local.
	
	Now since $\qsliceTspectra$ is a simplicial model category and $A$ is cofibrant in $\qsliceTspectra$,
	we have that the following map is a trivial fibration of simplicial sets:
		$$\xymatrix{Map(A,Z)\ar[r]& Map(A,\ast)=\ast}
		$$  
	The identity functor 
		$$\xymatrix{id:\qsliceTspectra \ar[r]& \weightqplusoneTspectra}
		$$
	is a left Quillen functor, since $\qsliceTspectra$ is a right Bousfield localization
	of $\weightqplusoneTspectra$.  Therefore $A$ is also cofibrant in $\weightqplusoneTspectra$, and
	since $\weightqplusoneTspectra$ is a left Bousfield localization of $\motivicTspectra$; it follows
	that $A$ is also cofibrant in $\motivicTspectra$.  On the other hand, we have that $Z$ is in particular
	fibrant in $\motivicTspectra$.
	Hence $\pi _{0}Map(A,Z)$ computes $[A,Z]_{Spt}$,
	since $\motivicTspectra$ is a simplicial model category.  But $Map(A,Z)\rightarrow \ast$
	is in particular a weak equivalence of simplicial sets, then
		$$[A,Z]_{Spt}\cong 0
		$$
	for every $\qleftlocmaps$-local $T$-spectrum $Z$.
	Finally, corollary \ref{cor.3.2.detecting-Cq-local-equivalences} 
	implies that $\ast \rightarrow A$ is a weak equivalence in $\weightqTspectra$.
	This finishes the proof, since we already know
	that $A$ is cofibrant in $\weightqTspectra$.
\end{proof}

\begin{lem}
		\label{lem.3.2.CqsqX---->sqX---stableweakequiv}
	Fix $q\in \mathbb Z$.  Then the natural map
		$$\xymatrix{C_{q}s_{q}X \ar[rr]^-{C_{q}^{s_{q}X}} && s_{q}X}
		$$
	is a weak equivalence in $\motivicTspectra$.
\end{lem}
\begin{proof}
	Consider the following commutative diagram in $\motivicTspectra$:
		$$\xymatrix{s_{q}X && Q_{s}s_{q}X \ar[ll]_-{Q_{s}^{s_{q}X}} \\
								C_{q}s_{q}X \ar[u]^-{C_{q}^{s_{q}X}} && 
								C_{q}Q_{s}s_{q}X \ar[u]_-{C_{q}^{Q_{s}s_{q}X}} \ar[ll]^-{C_{q}(Q_{s}^{s_{q}X})}}
		$$
	By construction $C_{q}^{s_{q}X}$, $C_{q}^{Q_{s}s_{q}X}$ are both
	weak equivalences in $\qconnectedTspectra$; and 
	\cite[proposition 3.1.5]{MR1944041} implies that
	$Q_{s}^{s_{q}X}$ is a $C_{eff}^{q}$-colocal equivalence in $\motivicTspectra$,
	i.e. a weak equivalence in $\qconnectedTspectra$.  Then the two out of three property
	for weak equivalences implies that $C_{q}(Q_{s}^{s_{q}X})$ is a weak equivalence
	in $\qconnectedTspectra$.
	
	Now \cite[theorem 3.2.13(2)]{MR1944041} 
	implies that $C_{q}(Q_{s}^{s_{q}X})$ is a weak equivalence in $\motivicTspectra$,
	since $C_{q}s_{q}X$ and $C_{q}Q_{s}s_{q}X$ are by construction $C_{eff}^{q}$-colocal $T$-spectra
	in $\motivicTspectra$.  It is clear that $Q_{s}^{s_{q}X}$ is a weak equivalence in $\motivicTspectra$,
	then by the two out of three property for weak equivalences, it suffices to show that
	$C_{q}^{Q_{s}s_{q}X}$ is a weak equivalence in $\motivicTspectra$.
	
	By theorem \ref{thm.3.1.slicefiltration}(\ref{thm.3.1.slicefiltration.b}) 
	we have that $s_{q}X$ is in $\stablehomotopyeffq$, then corollary
	\ref{cor.3.2.generating-Cqeff-colocalTspectra2} 
	implies that $Q_{s}s_{q}X$ is $C_{eff}^{q}$-colocal in $\motivicTspectra$.
	We already know that $C_{q}^{Q_{s}s_{q}X}$ is a $C_{eff}^{q}$-colocal equivalence in $\motivicTspectra$;
	then \cite[theorem 3.2.13(2)]{MR1944041} 
	implies that $C_{q}^{Q_{s}s_{q}X}$ is also a weak equivalence in $\motivicTspectra$,
	since by construction $C_{q}Q_{s}s_{q}X$ is a $C_{eff}^{q}$-colocal $T$-spectrum.
	This finishes the proof.	
\end{proof}

\begin{lem}
		\label{lem.3.2.sq==>Lq+1-local}
	Fix $q\in \mathbb Z$.  Then for every $T$-spectrum $X$,
	we have that
	$IQ_{T}Js_{q}X$ (see theorem \ref{thm.3.1.slicefiltration}) is $\qplusoneleftlocmaps$-local.	
\end{lem}
\begin{proof}
	Proposition \ref{prop.3.2.Lq-local-objects-classification} 
	implies that it suffices to check that $IQ_{T}Js_{q}X$
	satisfies the following conditions:
		\begin{enumerate}
			\item \label{prop.3.2.sq==>Lq+1-local.a}$IQ_{T}Js_{q}X$ is
						fibrant in $\motivicTspectra$.
			\item	\label{prop.3.2.sq==>Lq+1-local.b}For every $\generatorNRS \in C_{eff}^{q+1}$,
							$$[\generatorNRS ,IQ_{T}Js_{q}X]_{Spt}\cong 0$$
		\end{enumerate}
		
	Condition (\ref{prop.3.2.sq==>Lq+1-local.a}) holds trivially, since
	$IQ_{T}J$ is a fibrant replacement functor in $\motivicTspectra$.
	
	Fix $\generatorNRS \in C_{eff}^{q+1}$.  Since $C_{eff}^{q+1}\subseteq \stablehomotopyeffqplusone$,
	it follows from theorem \ref{thm.3.1.slicefiltration}(\ref{thm.3.1.slicefiltration.c})
	that:
		$$[\generatorNRS, IQ_{T}Js_{q}X]_{Spt}\cong [\generatorNRS, s_{q}X]_{Spt}\cong 0
		$$
	and this takes care of condition (\ref{prop.3.2.sq==>Lq+1-local.b}).
\end{proof}

\begin{lem}
		\label{lem.3.2.CqIQJfq+1---vanishes---Lq+1}
	Fix $q\in \mathbb Z$.  Then for every $T$-spectrum $X$,
	$C_{q}IQ_{T}Jf_{q+1}X \cong \ast$ in $\qslicestablehomotopy$.
\end{lem}
\begin{proof}
	Consider the following commutative diagram in $\motivicTspectra$:
		\begin{equation}
					\label{diagram.lem.3.2.CqIQJfq+1---vanishes---Lq+1}
			\begin{array}{c}
				\xymatrix{Q_{s}f_{q+1}X \ar[rr]^-{Q_{s}^{f_{q+1}X}}&& f_{q+1}X \ar[rr]^-{IQ_{T}J^{f_{q+1}X}}&& IQ_{T}Jf_{q+1}X\\
									C_{q}Q_{s}f_{q+1}X \ar[rr]_-{C_{q}(Q_{s}^{f_{q+1}X})} \ar[u]^-{C_{q}^{Q_{s}f_{q+1}X}}&& 
									C_{q}f_{q+1}X \ar[rr]_-{C_{q}(IQ_{T}J^{f_{q+1}X})} \ar[u]_-{C_{q}^{f_{q+1}X}}
									&& C_{q}IQ_{T}Jf_{q+1}X \ar[u]_-{C_{q}^{IQ_{T}Jf_{q+1}X}}}
			\end{array}
		\end{equation}
	We claim that all the maps in the diagram (\ref{diagram.lem.3.2.CqIQJfq+1---vanishes---Lq+1})
	above are weak equivalences in 
	$\motivicTspectra$.  In effect, it is clear that all the maps in the top row
	are weak equivalences in $\motivicTspectra$.  Hence, by the two out of three
	property for weak equivalences it suffices to show that $C_{q}^{Q_{s}f_{q+1}X}$, 
	$C_{q}(Q_{s}^{f_{q+1}X})$ and $C_{q}(IQ_{T}J^{f_{q+1}X})$
	are all weak equivalences in $\motivicTspectra$.
	
	On the other hand, \cite[proposition 3.1.5]{MR1944041} implies that all the maps in the top
	row are weak equivalences in $\qconnectedTspectra$, and it is clear that all the vertical
	maps are also weak equivalences in $\qconnectedTspectra$.  Thus, by the two out of three property for
	weak equivalences we have that all the maps in the diagram (\ref{diagram.lem.3.2.CqIQJfq+1---vanishes---Lq+1})
	above are weak equivalences
	in $\qconnectedTspectra$.
	
	By construction we have that $C_{q}Q_{s}f_{q+1}X$, $C_{q}f_{q+1}X$ and $C_{q}IQ_{T}Jf_{q+1}X$
	are all $C_{eff}^{q}$-colocal $T$-spectra in $\motivicTspectra$.  Then \cite[theorem 3.2.13(2)]{MR1944041}
	implies that $C_{q}(Q_{s}^{f_{q+1}X})$ and $C_{q}(IQ_{T}J^{f_{q+1}X})$ are both weak
	equivalences in $\motivicTspectra$.
	
	Now, by proposition \ref{prop.3.1.adjunctions-stableeffq} 
	we have that $f_{q+1}X\in \stablehomotopyeffqplusone \subseteq \stablehomotopyeffq$.
	Thus, corollary \ref{cor.3.2.generating-Cqeff-colocalTspectra2} 
	implies that $Q_{s}f_{q+1}X$ is a $C_{eff}^{q}$-colocal $T$-spectrum in $\motivicTspectra$.
	Then using \cite[theorem 3.2.13(2)]{MR1944041}  again,  we have that
	$C_{q}^{Q_{s}f_{q+1}X}$ is a weak equivalence in $\motivicTspectra$ since by construction
	$C_{q}Q_{s}f_{q+1}X$ is a $C_{eff}^{q}$-colocal $T$-spectrum and $C_{q}^{Q_{s}f_{q+1}X}$ is a
	$C_{eff}^{q}$-colocal equivalence in $\motivicTspectra$.
	
	This proves the claim, i.e. all the maps in the diagram (\ref{diagram.lem.3.2.CqIQJfq+1---vanishes---Lq+1})
	above are weak equivalences in
	$\motivicTspectra$.  Then using \cite[proposition 3.1.5]{MR1944041} again,
	we have that all the maps in the diagram (\ref{diagram.lem.3.2.CqIQJfq+1---vanishes---Lq+1})
	above are also weak equivalences in $\qsliceTspectra$.
	Therefore, to finish the proof it is enough to check that $\ast \rightarrow Q_{s}f_{q+1}X$ 
	is a weak equivalence in $\qsliceTspectra$.
	
	But corollary \ref{cor.3.2.Qsfq-vanishesin-LqSH} implies that $\ast \rightarrow Q_{s}f_{q+1}X$
	is a weak equivalence in $\weightqplusoneTspectra$.  Therefore, using \cite[proposition 3.1.5]{MR1944041},
	we have that $\ast \rightarrow Q_{s}f_{q+1}X$ is a $S(q)$-colocal equivalence in
	$\weightqplusoneTspectra$, i.e. a weak equivalence in $\qsliceTspectra$.
	This finishes the proof.
\end{proof}

\begin{prop}
		\label{prop.3.2.sq--natequiv--CqIQJsq}
	Fix $q\in \mathbb Z$.  Then for
	every $T$-spectrum $X$,
	the following maps of $T$-spectra:
		\begin{equation}
					\label{diagram.cor.3.2.sq--natequiv--Wq+1Qssq}
			\begin{array}{c}
				\xymatrix{s_{q}X \ar[rr]^-{IQ_{T}J^{s_{q}X}}&& IQ_{T}Js_{q}X  
									&& C_{q}IQ_{T}Js_{q}X \ar[ll]_-{C_{q}^{IQ_{T}Js_{q}X}}}
			\end{array}
		\end{equation}
	are both
	weak equivalences in $\motivicTspectra$.	
	
	Furthermore, these weak equivalences induce natural isomorphisms 
	between the following exact functors
		$$\xymatrix{\stablehomotopy \ar@<1ex>[rr]^-{s_{q}} \ar@<-1ex>[rr]_-{IQ_{T}Js_{q}}&& 
								\stablehomotopy}
		$$
		
		$$\xymatrix{\stablehomotopy \ar@<1ex>[rr]^-{IQ_{T}Js_{q}} \ar@<-1ex>[rr]_-{C_{q}IQ_{T}Js_{q}}&& 
								\stablehomotopy}
		$$
\end{prop}
\begin{proof}
	The naturality of the maps $IQ_{T}J^{X}:X\rightarrow IQ_{T}JX$ and
	$C_{q}^{X}:C_{q}X\rightarrow X$ implies that we have 
	induced natural transformations of functors $s_{q}\rightarrow IQ_{T}Js_{q}$
	and $C_{q}IQ_{T}Js_{q}\rightarrow IQ_{T}Js_{q}$.
	Hence, it is enough to show that for every $T$-spectrum $X$, 
	$IQ_{T}J^{s_{q}X}$ and $C_{q}^{IQ_{T}Js_{q}X}$ are
	weak equivalences in $\motivicTspectra$.
	
	It is clear that
	$IQ_{T}J^{s_{q}X}$ is a weak equivalence in $\motivicTspectra$, since
	$IQ_{T}J$ is a fibrant replacement functor for $\motivicTspectra$.
	
	We now proceed to show that $C_{q}^{IQ_{T}Js_{q}X}$ is
	a weak equivalence in $\motivicTspectra$.
	Consider the following commutative diagram in $\motivicTspectra$:
		$$\xymatrix{s_{q}X \ar[d]_-{IQ_{T}J^{s_{q}X}}&& C_{q}s_{q}X 
								\ar[ll]_-{C_{q}^{s_{q}X}} \ar[d]^-{C_{q}(IQ_{T}J^{s_{q}X})}\\
								IQ_{T}Js_{q}X && C_{q}IQ_{T}Js_{q}X \ar[ll]^-{C_{q}^{IQ_{T}Js_{q}X}}}
		$$
	Lemma \ref{lem.3.2.CqsqX---->sqX---stableweakequiv} 
	implies that $C_{q}^{s_{q}X}$ is a weak equivalence in $\motivicTspectra$.
	Since we know that $IQ_{T}J^{s_{q}X}$ is always a weak equivalence in $\motivicTspectra$,
	the two out of three property for weak equivalences implies that
	it suffices to check that $C_{q}(IQ_{T}J^{s_{q}X})$ is also a weak equivalence
	in $\motivicTspectra$.
	
	Using \cite[proposition 3.1.5]{MR1944041}, we have that
	$IQ_{T}J^{s_{q}X}$ is a $C_{eff}^{q}$-colocal equivalence.  Then the
	two out of three property for $C_{eff}^{q}$-colocal equivalences implies
	that $C_{q}(IQ_{T}J^{s_{q}X})$ is a $C_{eff}^{q}$-colocal equivalence,
	since by construction $C_{q}^{s_{q}X}$ and $C_{q}^{IQ_{T}Js_{q}X}$
	are both $C_{eff}^{q}$-colocal equivalences.
	
	Finally, by construction $C_{q}s_{q}X$ and $C_{q}IQ_{T}Js_{q}X$ are both
	$C_{eff}^{q}$-colocal, therefore \cite[theorem 3.2.13(2)]{MR1944041}
	implies that
	$C_{q}(IQ_{T}J^{s_{q}X})$ is a weak equivalence in $\motivicTspectra$, as we wanted.
\end{proof}

\begin{prop}
		\label{prop.3.2.Wq+1CqIQJsq-q+1stableweakequiv}
	Fix $q\in \mathbb Z$.  Then  for every $T$-spectrum $X$,
	the natural map:
		$$\xymatrix{C_{q}IQ_{T}Js_{q}X \ar[rr]^-{W_{q+1}^{C_{q}IQ_{T}Js_{q}X}}&&
								W_{q+1}C_{q}IQ_{T}Js_{q}X}
		$$
	is a weak equivalence in $\motivicTspectra$.

	Furtheremore, this weak equivalence induces a natural isomorphism 
	between the following exact functors
		$$\xymatrix{\stablehomotopy \ar@<1ex>[rr]^-{C_{q}IQ_{T}Js_{q}} \ar@<-1ex>[rr]_-{W_{q+1}C_{q}IQ_{T}Js_{q}}&& 
								\stablehomotopy}
		$$
\end{prop}
\begin{proof}
	The naturality of the maps $W_{q+1}^{X}:X\rightarrow W_{q+1}X$
	implies that we have an
	induced natural transformation of functors
	$C_{q}IQ_{T}Js_{q}\rightarrow W_{q+1}C_{q}IQ_{T}Js_{q}$.
	Hence, it is enough to show that for every $T$-spectrum $X$, 
	$W_{q+1}^{C_{q}IQ_{T}Js_{q}}$ is a
	weak equivalence in $\motivicTspectra$.
	
	Consider the following commutative diagram in $\motivicTspectra$:
		$$\xymatrix{IQ_{T}Js_{q}X \ar[rr]^-{W_{q+1}^{IQ_{T}Js_{q}X}} && 
								W_{q+1}IQ_{T}Js_{q}X  \\
								C_{q}IQ_{T}Js_{q}X \ar[rr]_-{W_{q+1}^{C_{q}IQ_{T}Js_{q}X}} \ar[u]^-{C_{q}^{IQ_{T}Js_{q}X}}&& 
								W_{q+1}C_{q}IQ_{T}Js_{q}X \ar[u]_-{W_{q+1}(C_{q}^{IQ_{T}Js_{q}X})}}
		$$
	By construction, $W_{q+1}^{IQ_{T}Js_{q}X}$ is a $\qplusoneleftlocmaps$-local equivalence, and
	$W_{q+1}IQ_{T}Js_{q}X$ is $\qplusoneleftlocmaps$-local in $\motivicTspectra$.
	By lemma \ref{lem.3.2.sq==>Lq+1-local} 
	we have that $IQ_{T}Js_{q}X$ is also $\qplusoneleftlocmaps$-local.  Therefore,
	\cite[theorem 3.2.13(1)]{MR1944041} 
	implies that $W_{q+1}^{IQ_{T}Js_{q}X}$ is a weak equivalence in $\motivicTspectra$.
	
	Now, it follows directly from proposition \ref{prop.3.2.sq--natequiv--CqIQJsq} 
	that $C_{q}^{IQ_{T}Js_{q}X}$ is a weak equivalence
	in $\motivicTspectra$.  Hence by the two out of three property for weak equivalences, it
	suffices to show that $W_{q+1}(C_{q}^{IQ_{T}Js_{q}X})$ is a weak equivalence
	in $\motivicTspectra$.
	
	We already know that $C_{q}^{IQ_{T}Js_{q}X}$ is a weak equivalence in
	$\motivicTspectra$, then using \cite[proposition 3.1.5]{MR1944041} we have that
	$C_{q}^{IQ_{T}Js_{q}X}$ is a $\qplusoneleftlocmaps$-local equivalence.  Then the two
	out of three property for $\qplusoneleftlocmaps$-local equivalences implies that
	$W_{q+1}(C_{q}^{IQ_{T}Js_{q}X})$ is also a $\qplusoneleftlocmaps$-local equivalence,
	since by construction $W_{q+1}^{IQ_{T}Js_{q}X}$ and $W_{q+1}^{C_{q}IQ_{T}Js_{q}X}$ 
	are both $\qplusoneleftlocmaps$-local equivalences.
	
	Finally, by construction $W_{q+1}IQ_{T}Js_{q}X$ and $W_{q+1}C_{q}IQ_{T}Js_{q}X$ 
	are $\qplusoneleftlocmaps$-local in $\motivicTspectra$,
	then \cite[theorem 3.2.13(1)]{MR1944041} implies that $W_{q+1}(C_{q}^{IQ_{T}Js_{q}X})$ 
	is a weak equivalence
	in $\motivicTspectra$, as we wanted.
\end{proof}

\begin{prop}
		\label{prop.3.2.CqWq+1CqIQJsq-q+1stableweakequiv}
	Fix $q\in \mathbb Z$.  Then  for every $T$-spectrum $X$,
	the natural map:
		$$\xymatrix{W_{q+1}C_{q}IQ_{T}Js_{q}X &&&
								C_{q}W_{q+1}C_{q}IQ_{T}Js_{q}X \ar[lll]_-{C_{q}^{W_{q+1}C_{q}IQ_{T}Js_{q}X}}}
		$$
	is a weak equivalence in $\motivicTspectra$.

	Furtheremore, this weak equivalence induces a natural isomorphism  
	between the following exact functors	
		$$\xymatrix{\stablehomotopy \ar@<1ex>[rrr]^-{W_{q+1}C_{q}IQ_{T}Js_{q}} \ar@<-1ex>[rrr]_-{C_{q}W_{q+1}C_{q}IQ_{T}Js_{q}}&&& 
								\stablehomotopy}
		$$
\end{prop}
\begin{proof}
	The naturality of the maps 
	$C_{q}^{X}:C_{q}X\rightarrow X$ implies that we have an 
	induced natural transformation of functors
	$C_{q}W_{q+1}C_{q}IQ_{T}Js_{q}\rightarrow W_{q+1}C_{q}IQ_{T}Js_{q}$.
	Hence, it is enough to show that for every $T$-spectrum $X$, 
	$C_{q}^{W_{q+1}C_{q}IQ_{T}Js_{q}X}$ is a
	weak equivalence in $\motivicTspectra$.
	
	Consider the following commutative diagram in $\motivicTspectra$:
		$$\xymatrix{C_{q}IQ_{T}Js_{q}X \ar[d]_-{W_{q+1}^{C_{q}IQ_{T}Js_{q}X}}&&& 
								C_{q}C_{q}IQ_{T}Js_{q}X 
								\ar[lll]_-{C_{q}^{C_{q}IQ_{T}Js_{q}X}} \ar[d]^-{C_{q}(W_{q+1}^{C_{q}IQ_{T}Js_{q}X})}\\
								W_{q+1}C_{q}IQ_{T}Js_{q}X &&& C_{q}W_{q+1}C_{q}IQ_{T}Js_{q}X \ar[lll]^-{C_{q}^{W_{q+1}C_{q}IQ_{T}Js_{q}X}}}
		$$
	By construction $C_{q}^{C_{q}IQ_{T}Js_{q}X}$ is a $C_{eff}^{q}$-colocal equivalence,
	and $C_{q}IQ_{T}Js_{q}X$, 
	$C_{q}C_{q}IQ_{T}Js_{q}X$ are both $C_{eff}^{q}$-colocal in $\motivicTspectra$.
	Therefore, \cite[theorem 3.2.13(2)]{MR1944041} 
	implies that $C_{q}^{C_{q}IQ_{T}Js_{q}X}$ is a weak equivalence in $\motivicTspectra$.
	
	Now, it follows directly from proposition \ref{prop.3.2.Wq+1CqIQJsq-q+1stableweakequiv} 
	that $W_{q+1}^{C_{q}IQ_{T}Js_{q}X}$ is a weak equivalence
	in $\motivicTspectra$.  Hence by the two out of three property for weak equivalences, it
	suffices to show that $C_{q}(W_{q+1}^{C_{q}IQ_{T}Js_{q}X})$ is a weak equivalence
	in $\motivicTspectra$.
	
	We already know that $W_{q+1}^{C_{q}IQ_{T}Js_{q}X}$ is a weak equivalence in
	$\motivicTspectra$, then using \cite[proposition 3.1.5]{MR1944041} we have that
	$W_{q+1}^{C_{q}IQ_{T}Js_{q}X}$ is a $C_{eff}^{q}$-colocal equivalence.  Then the two
	out of three property for $C_{eff}^{q}$-colocal equivalences implies that
	$C_{q}(W_{q+1}^{C_{q}IQ_{T}Js_{q}X})$ is also a $C_{eff}^{q}$-colocal equivalence,
	since by construction $C_{q}^{C_{q}IQ_{T}Js_{q}X}$ and $C_{q}^{W_{q+1}C_{q}IQ_{T}Js_{q}X}$ 
	are both $C_{eff}^{q}$-colocal equivalences.
	
	Finally, by construction $C_{q}C_{q}IQ_{T}Js_{q}X$ and $C_{q}W_{q+1}C_{q}IQ_{T}Js_{q}X$ 
	are $C_{eff}^{q}$-colocal in $\motivicTspectra$,
	then \cite[theorem 3.2.13(2)]{MR1944041} implies that $C_{q}(W_{q+1}^{C_{q}IQ_{T}Js_{q}X})$ 
	is a weak equivalence
	in $\motivicTspectra$, as we wanted.
\end{proof}

\begin{prop}
		\label{prop.3.2.CqIQJX--CqIQJfqX--CqIQJsqX==>q-colocalequiv}
	Fix $q\in \mathbb Z$.  Then for every $T$-spectrum $X$,
	the following natural maps in $\qconnectedstablehomotopy$
	(see  proposition \ref{prop.3.1.counit-properties} and
	theorem \ref{thm.3.1.slicefiltration}):
		$$\xymatrix{IQ_{T}JX && IQ_{T}Jf_{q}X  \ar[rr]^-{IQ_{T}J(\pi _{q}^{X})} \ar[ll]_-{IQ_{T}J(\theta _{X})}&& 
								IQ_{T}Js_{q}X}
		$$
	become isomorphisms in $\qslicestablehomotopy$
	after applying the functor $C_{q}$:
		$$\xymatrix{C_{q}IQ_{T}JX && C_{q}IQ_{T}Jf_{q}X  \ar[rr]^-{C_{q}IQ_{T}J(\pi _{q}^{X})}_-{\cong} 
								\ar[ll]_-{C_{q}IQ_{T}J(\theta _{X})}^-{\cong}&& 
								C_{q}IQ_{T}Js_{q}X}
		$$
\end{prop}
\begin{proof}
	Proposition \ref{prop.3.2.counit-stablehomotopyeffq=>iso-qconnected.b}
	implies that the map
		$$\xymatrix{IQ_{T}Jf_{q}X 
								\ar[rr]^-{IQ_{T}J(\theta _{X})}_-{\cong}&& IQ_{T}JX}
		$$
	is an isomorphism in $\qconnectedstablehomotopy$.  Hence
	using corollary \ref{cor.3.2.adjunctions--Rq==>Sq} we have that
		$$\xymatrix{C_{q}IQ_{T}Jf_{q}X  \ar[rr]^-{C_{q}IQ_{T}J(\theta _{X})}_-{\cong} 
								&& C_{q}IQ_{T}JX}
		$$
	is an isomorphism in $\qslicestablehomotopy$.
	
	On the other hand, theorem \ref{thm.3.1.slicefiltration}(\ref{thm.3.1.slicefiltration.a}) 
	implies that we have the following distinguished
	triangle in $\stablehomotopy$:
		$$\xymatrix{f_{q+1}X \ar[r]& f_{q}X \ar[r]^-{\pi _{q}^{X}}& s_{q}X \ar[r]& \Tsuspfunctor ^{1,0}f_{q+1}X}
		$$
	Proposition \ref{prop.3.3.cofibrant-replacement=>triangulatedfunctor} implies that after applying $IQ_{T}J$
	we get the following distinguished triangle in $\qconnectedstablehomotopy$
		$$\xymatrix{IQ_{T}Jf_{q+1}X \ar[r]& IQ_{T}Jf_{q}X \ar[rr]^-{IQ_{T}J(\pi _{q}^{X})}&& IQ_{T}Js_{q}X 
								\ar[r]& \Tsuspfunctor ^{1,0}IQ_{T}Jf_{q+1}X}
		$$
	Now corollary \ref{cor.3.2.adjunctions--Rq==>Sq} 
	implies that after applying $C_{q}$ we get the following
	distinguished triangle in $\qslicestablehomotopy$
		$$\xymatrix{C_{q}IQ_{T}Jf_{q+1}X \ar[r]& C_{q}IQ_{T}Jf_{q}X \ar[rr]^-{C_{q}IQ_{T}J(\pi _{q}^{X})}&& C_{q}IQ_{T}Js_{q}X 
								\ar[d]\\ 
								&&& \Tsuspfunctor ^{1,0}C_{q}IQ_{T}Jf_{q+1}X}
		$$
	Therefore it is enough to check that $C_{q}IQ_{T}Jf_{q+1}X\cong \ast$ in
	$\qslicestablehomotopy$.  But this follows directly from
	lemma \ref{lem.3.2.CqIQJfq+1---vanishes---Lq+1}.
\end{proof}

\begin{cor}
		\label{cor.3.2.CqWq+1CqIQJX--CqWq+1CqIQJfqX--CqWq+1CqIQJsqX==>stably-weakequivs}
	Fix $q\in \mathbb Z$.  Then for every $T$-spectrum $X$,
	the following natural maps in $\stablehomotopy$
	(see  proposition \ref{prop.3.1.counit-properties} and
	theorem \ref{thm.3.1.slicefiltration}):
		$$\xymatrix{X && f_{q}X  \ar[rr]^-{\pi _{q}^{X}} \ar[ll]_-{\theta _{X}}&& 
								s_{q}X}
		$$
	become isomorphisms in $\stablehomotopy$
	after applying the functor $C_{q}W_{q+1}C_{q}IQ_{T}J$:
		$$\xymatrix{C_{q}W_{q+1}C_{q}IQ_{T}JX &&&& C_{q}W_{q+1}C_{q}IQ_{T}Jf_{q}X  
								\ar[d]^-{C_{q}W_{q+1}C_{q}IQ_{T}J(\pi _{q}^{X})}_-{\cong} 
								\ar[llll]_-{C_{q}W_{q+1}C_{q}IQ_{T}J(\theta _{X})}^-{\cong}\\ 
								&&&& C_{q}W_{q+1}C_{q}IQ_{T}Js_{q}X}
		$$

	Furtheremore, these maps induce natural isomorphisms 
	between the following exact functors
		$$\xymatrix{\stablehomotopy \ar@<1ex>[rrr]^-{C_{q}W_{q+1}C_{q}IQ_{T}Js_{q}} 
								\ar@<-1ex>[rrr]_-{C_{q}W_{q+1}C_{q}IQ_{T}Jf_{q}}&&& 
								\stablehomotopy}
		$$
		
		$$\xymatrix{\stablehomotopy \ar@<1ex>[rrr]^-{C_{q}W_{q+1}C_{q}IQ_{T}Jf_{q}} 
								\ar@<-1ex>[rrr]_-{C_{q}W_{q+1}C_{q}IQ_{T}J} &&& 
								\stablehomotopy}
		$$
\end{cor}
\begin{proof}
	The naturality of the maps $\pi _{q}^{X}:f_{q}X\rightarrow s_{q}X$ and
	$\theta _{X}:f_{q}X\rightarrow X$ implies that we have 
	induced natural transformations of functors $C_{q}W_{q+1}C_{q}IQ_{T}Jf_{q}\rightarrow C_{q}W_{q+1}C_{q}IQ_{T}Js_{q}$
	and $C_{q}W_{q+1}C_{q}IQ_{T}Jf_{q}\rightarrow C_{q}W_{q+1}C_{q}IQ_{T}J$.
	Hence, it is enough to show that for every $T$-spectrum $X$, 
	$C_{q}W_{q+1}C_{q}IQ_{T}J(\pi _{q}^{X})$ and $C_{q}W_{q+1}C_{q}IQ_{T}J(\theta _{X})$ are
	weak equivalences in $\motivicTspectra$.

	Proposition \ref{prop.3.2.CqIQJX--CqIQJfqX--CqIQJsqX==>q-colocalequiv} 
	implies that the following natural maps
		$$\xymatrix{C_{q}IQ_{T}JX && C_{q}IQ_{T}Jf_{q}X  \ar[rr]^-{C_{q}IQ_{T}J(\pi _{q}^{X})}_-{\cong} 
								\ar[ll]_-{C_{q}IQ_{T}J(\theta _{X})}^-{\cong}&& 
								C_{q}IQ_{T}Js_{q}X}
		$$
	are isomorphisms in $\qslicestablehomotopy$.  Then the result follows immediately from
	corollary \ref{cor.3.2.adjunctions--Rq==>Sq} and 
	proposition \ref{prop.3.3.cofibrant-replacement=>triangulatedfunctor}.
\end{proof}

\begin{thm}
		\label{thm.3.2.Sq-models-sq}
	Fix $q\in \mathbb Z$.  Then for every $T$-spectrum $X$, we have the following diagram in 
	$\stablehomotopy$:
		\begin{equation}
					\label{diagram.3.2.sq-lifting}
			\begin{array}{c}
				\xymatrix{& W_{q+1}C_{q}IQ_{T}Js_{q}X  & \\
									C_{q}IQ_{T}Js_{q}X \ar[ur]^-{W_{q+1}^{C_{q}IQ_{T}Js_{q}X}}_-{\cong} 
									\ar[d]_-{C_{q}^{IQ_{T}Js_{q}X}}^-{\cong}&& 
									C_{q}W_{q+1}C_{q}IQ_{T}Js_{q}X 
									\ar[ul]_-{\; \; \; \; \; \; C_{q}^{W_{q+1}C_{q}IQ_{T}Js_{q}X}}^-{\cong}\\
									IQ_{T}Js_{q}X   &&
									C_{q}W_{q+1}C_{q}IQ_{T}Jf_{q}X \ar[u]^-{C_{q}W_{q+1}C_{q}IQ_{T}J(\pi _{q}^{X})}_-{\cong}
									\ar[d]_-{C_{q}W_{q+1}C_{q}IQ_{T}J(\theta _{X})}^-{\cong}\\
									s_{q}X \ar[u]^-{IQ_{T}J^{s_{q}X}}_-{\cong}&& C_{q}W_{q+1}C_{q}IQ_{T}JX}
			\end{array}
		\end{equation}
	where all the maps are isomorphisms in $\stablehomotopy$.
	Furthermore, this diagram induces a
	natural isomorphism  
	between the following exact functors:
		$$\xymatrix{\stablehomotopy \ar@<1ex>[rrr]^-{s_{q}}  
								\ar@<-1ex>[rrr]_-{C_{q}W_{q+1}C_{q}IQ_{T}J} &&& \stablehomotopy}
		$$
\end{thm}
\begin{proof}	
	The fact that $IQ_{T}J^{s_{q}X}$ and $C_{q}^{IQ_{T}Js_{<q}X}$
	are isomorphisms in $\stablehomotopy$ follows from
	proposition \ref{prop.3.2.sq--natequiv--CqIQJsq}.  
	Now proposition \ref{prop.3.2.Wq+1CqIQJsq-q+1stableweakequiv} implies that
	$W_{q+1}^{C_{q}IQ_{T}Js_{q}X}$ is an isomorphism in $\stablehomotopy$, and proposition
	\ref{prop.3.2.CqWq+1CqIQJsq-q+1stableweakequiv} implies that $C_{q}^{W_{q+1}C_{q}IQ_{T}Js_{q}X}$ is also an
	isomorphism in $\stablehomotopy$.
	Finally, corollary \ref{cor.3.2.CqWq+1CqIQJX--CqWq+1CqIQJfqX--CqWq+1CqIQJsqX==>stably-weakequivs} implies that 
	$C_{q}W_{q+1}C_{q}IQ_{T}J(\pi _{q}^{X})$ and 
	$C_{q}W_{q+1}C_{q}IQ_{T}J(\theta _{X})$ are both isomorphisms in $\stablehomotopy$.
	
	This shows that all the maps in the diagram (\ref{diagram.3.2.sq-lifting}) are isomorphisms
	in $\stablehomotopy$, therefore for every $T$-spectrum $X$ we can define the
	following composition in $\stablehomotopy$
		\begin{equation}
					\label{diagram.3.2.sq-lifting.b}
			\begin{array}{c}
				\xymatrix{& W_{q+1}C_{q}IQ_{T}Js_{q}X  
									\ar[dr]^-{\ \ \ \ \ \ \ \ (C_{q}^{W_{q+1}C_{q}IQ_{T}Js_{q}X})^{-1}}_-{\cong}& \\
									C_{q}IQ_{T}Js_{q}X \ar[ur]^-{W_{q+1}^{C_{q}IQ_{T}Js_{q}X}}_-{\cong} 
									&& 
									C_{q}W_{q+1}C_{q}IQ_{T}Js_{q}X \ar[d]_-{(C_{q}W_{q+1}C_{q}IQ_{T}J(\pi _{q}^{X}))^{-1}}^-{\cong}
									\\
									IQ_{T}Js_{q}X   \ar[u]^-{(C_{q}^{IQ_{T}Js_{q}X})^{-1}}_-{\cong}&&
									C_{q}W_{q+1}C_{q}IQ_{T}Jf_{q}X 
									\ar[d]_-{C_{q}W_{q+1}C_{q}IQ_{T}J(\theta _{X})}^-{\cong}\\
									s_{q}X \ar[u]^-{IQ_{T}J^{s_{q}X}}_-{\cong}&& C_{q}W_{q+1}C_{q}IQ_{T}JX}
			\end{array}
		\end{equation}
	which is an isomorphism. 	
	
	On the other hand, propositions \ref{prop.3.2.sq--natequiv--CqIQJsq},
	\ref{prop.3.2.Wq+1CqIQJsq-q+1stableweakequiv} and \ref{prop.3.2.CqWq+1CqIQJsq-q+1stableweakequiv},
	and corollary \ref{cor.3.2.CqWq+1CqIQJX--CqWq+1CqIQJfqX--CqWq+1CqIQJsqX==>stably-weakequivs}
	imply all together that the isomorphisms defined in diagram (\ref{diagram.3.2.sq-lifting.b})
	induce a natural isomorphism of functors $s_{q}\stackrel{\cong}{\rightarrow}C_{q}W_{q+1}C_{q}IQ_{T}J$.
	This finishes the proof.
\end{proof}

\begin{prop}
		\label{prop.3.2.lifting-map-fq-->sq}
	Fix $q\in \mathbb Z$.  Let $\eta$ denote the unit of the adjuntion
	$(C_{q},W_{q+1},\varphi ):\qconnectedstablehomotopy \rightarrow \qslicestablehomotopy$
	constructed in corollary \ref{cor.3.2.adjunctions--Rq==>Sq}.  Then the natural transformation
	$\pi _{q}:f_{q}\rightarrow s_{q}$ (see theorem \ref{thm.3.1.slicefiltration}) gets canonically identified,
	through the equivalence of categories $r_{q}C_{q}$, $IQ_{T}Ji_{q}$ constructed in proposition \ref{prop.3.2.Rq-lifts-qSH}
	with the following map in $\stablehomotopy$:
		$$\xymatrix{C_{q}IQ_{T}JX \ar[rr]^-{C_{q}(\eta _{IQ_{T}JX})}&& C_{q}W_{q+1}C_{q}IQ_{T}JX}
		$$
\end{prop}
\begin{proof}
	It follows directly from theorem \ref{thm.3.1.slicefiltration},
	corollary \ref{cor.3.2.classifying-Sq-colocal-equivs.b} together with \cite[proposition 9.1.8]{MR1812507}.
\end{proof}

\begin{rmk}
		\label{rmk.3.2.lifting-sq}
	Theorem \ref{thm.3.2.Sq-models-sq} gives the desired lifting to the
	model category level
	for the functors $s_{q}$ defined in theorem \ref{thm.3.1.slicefiltration}; and it 
	completes the program that we started at the beginning of this section,
	where the goal was to get a lifting for the slice
	functors $s_{q}$.
\end{rmk}

\end{section}
\begin{section}{The Symmetric Model Structure for the Slice Filtration}
		\label{section.3.3.symmmodstrslicefilt}

	Our goal now is to lift the model structures constructed in section
	\ref{section.3.2.modstrslicefilt} to the category of symmetric $T$-spectra, in order to have a
	natural framework for the study of the multiplicative properties
	of Voevodsky's slice filtration.
	
	Let $\symmstablehomotopy$ denote the homotopy category associated
	to $\motivicsymmTspectra$.
	We call $\symmstablehomotopy$ the \emph{motivic symmetric stable homotopy category}.
	We will denote by $[-,-]_{Spt}^{\Sigma}$ the set of maps between two
	objects in $\symmstablehomotopy$.

\begin{defi}
		\label{def.3.3.cofibrantstable-replacementfunctor}
	Let $Q_{\Sigma}$ denote
	a cofibrant replacement functor in $\motivicsymmTspectra$; such that for every
	symmetric $T$-spectrum $X$, the natural map
		$$\xymatrix{Q_{\Sigma}X\ar[r]^-{Q_{\Sigma}^{X}}& X}
		$$
	is a trivial fibration in $\motivicsymmTspectra$.
\end{defi}
	
\begin{defi}
		\label{def.3.3.fibrantstable-replacementfunctor}
	Let $R_{\Sigma}$ denote
	a fibrant replacement functor in $\motivicsymmTspectra$; such that for every
	symmetric $T$-spectrum $X$, the natural map
		$$\xymatrix{X\ar[r]^-{R_{\Sigma}^{X}}& R_{\Sigma}X}
		$$
	is a trivial cofibration in $\motivicsymmTspectra$.
\end{defi}
	
\begin{prop}
		\label{prop.3.3.symmstable-homotopy==triangulated-category}
	The motivic symmetric stable homotopy category $\symmstablehomotopy$
	has a structure of triangulated category defined as follows: 
	\begin{enumerate}
		\item The suspension $\Sigma _{T}^{1,0}$ functor is given by
					$$\xymatrix@R=.5pt{-\wedge S^{1}:\symmstablehomotopy \ar[r]& \symmstablehomotopy \\
															X \ar@{|->}[r]& Q_{\Sigma}X\wedge S^{1}}
					$$
		\item The distinguished triangles are
					isomorphic to triangles of the form
						$$\xymatrix{A \ar[r]^-{i}& B \ar[r]^-{j}& C \ar[r]^-{k}& \Sigma_{T}^{1,0}A}
						$$
					where $i$ is a cofibration in $\motivicsymmTspectra$, and $C$ is the homotopy cofibre of $i$.
	\end{enumerate}
\end{prop}
\begin{proof}
	Theorem \ref{thm2.6.stablesymmetricmodelstructure} implies in particular that
	$\motivicsymmTspectra$ is a pointed simplicial model category, and
	theorem \ref{thm.2.6.Tsusp==QuillenequivonSymmTspectra} implies that the adjunction:
		$$\xymatrix{(-\wedge S^{1},\Omega _{S^{1}},\varphi):\motivicsymmTspectra \ar[r]& \motivicsymmTspectra}
		$$
	is a Quillen equivalence.  The result now follows from the work of
	Quillen in \cite[sections I.2 and I.3]{MR0223432} and the work of 
	Hovey in \cite[chapters VI and VII]{MR1650134}
	(see \cite[proposition 7.1.6]{MR1650134}).
\end{proof}

\begin{thm}
		\label{thm.3.3.symmetrization-functor-exact-equivalence}
	The adjunction
		$$\xymatrix{(V,U,\varphi):\motivicTspectra \ar[r]& \motivicsymmTspectra}
		$$
	given by the symmetrization and the forgetful functors, induces
	an adjunction
		$$\xymatrix{(VQ_{s},UR_{\Sigma},\varphi):\stablehomotopy \ar[r]& \symmstablehomotopy}
		$$
	of exact funtors between triangulated categories.
	Furthermore, $VQ_{s}$ and $UR_{\Sigma}$ are both equivalences of categories.
\end{thm}
\begin{proof}
	Theorem \ref{thm.2.6.T-spectra===symmTspectra} implies that the adjunction $(V,U,\varphi)$
	is a Quillen equivalence.
	Therefore we get
	the following adjunction at the level of the associated homotopy categories:
		$$\xymatrix{(VQ_{s}, UR_{\Sigma}, \varphi ):\stablehomotopy \ar[r]& \symmstablehomotopy}
		$$
	
	Now \cite[proposition 1.3.13]{MR1650134} implies that $VQ_{s} ,UR_{\Sigma}$ are both equivalences of
	categories.
	Finally, proposition \ref{prop.2.6.enriched-symmetrization-adjunction}
	together with \cite[proposition 6.4.1]{MR1650134} imply that
	$VQ_{s}$ maps cofibre sequences in $\stablehomotopy$ to cofibre sequences in
	$\symmstablehomotopy$.
	Therefore 
	using proposition 7.1.12 in \cite{MR1650134} we have
	that $VQ_{s}$ and $UR_{\Sigma}$ are both exact functors between triangulated categories.
\end{proof}

\begin{cor}
		\label{cor.3.3.symmetric-fq,s<q,sq}
	Fix $q\in \mathbb Z$.
	\begin{enumerate}
		\item	\label{cor.3.3.symmetric-fq,s<q,sq.a}The exact functor (see remark \ref{rmk.3.1.unit=iso}) 
						$$\xymatrix{f_{q}:\stablehomotopy \ar[r]& \stablehomotopy}
						$$
					gets canonically identified with the following exact functor:
						$$\xymatrix@R=.5pt{\tilde{f}_{q}:\symmstablehomotopy \ar[r]& \symmstablehomotopy \\
															X \ar@{|->}[r]& VQ_{s}(f_{q}(UR_{\Sigma}X))}
						$$
					i.e. $\tilde{f}_{q}=VQ_{s}\circ f_{q}\circ UR_{\Sigma}$.
		\item	\label{cor.3.3.symmetric-fq,s<q,sq.b}The exact functor (see theorem \ref{thm.3.1.motivictower}) 
						$$\xymatrix{s_{<q}:\stablehomotopy \ar[r]& \stablehomotopy}
						$$
					gets canonically identified with the following exact functor:
						$$\xymatrix@R=.5pt{\tilde{s}_{<q}:\symmstablehomotopy \ar[r]& \symmstablehomotopy \\
															X \ar@{|->}[r]& VQ_{s}(s_{<q}(UR_{\Sigma}X))}
						$$
					i.e. $\tilde{s}_{<q}=VQ_{s}\circ s_{<q}\circ UR_{\Sigma}$.
		\item	\label{cor.3.3.symmetric-fq,s<q,sq.c}The exact functor (see theorem \ref{thm.3.1.slicefiltration}) 
						$$\xymatrix{s_{q}:\stablehomotopy \ar[r]& \stablehomotopy}
						$$
					gets canonically identified with the following exact functor:
						$$\xymatrix@R=.5pt{\tilde{s}_{q}:\symmstablehomotopy \ar[r]& \symmstablehomotopy \\
															X \ar@{|->}[r]& VQ_{s}(s_{q}(UR_{\Sigma}X))}
						$$
					i.e. $\tilde{s}_{q}=VQ_{s}\circ s_{q}\circ UR_{\Sigma}$.
	\end{enumerate}
\end{cor}
\begin{proof}
	Follows immediately from theorem \ref{thm.3.3.symmetrization-functor-exact-equivalence}.
\end{proof}

\begin{lem}
		\label{lem.3.3.compact-respects-coproducts-symmspectra}
	Let $X\in \pointedmotivic$ be a pointed simplicial presheaf which
	is compact in the sense of Jardine (see definition \ref{def.2.3.compactness}), and let
	$F_{n}^{\Sigma}(X)$ be the symmetric $T$-spectrum constructed in definition \ref{def.2.6.Fnsigma-Evn-adjunction}.
	Consider an arbitrary collection of symmetric $T$-spectra $\{ Z_{i}\}_{i\in I}$ indexed by a set $I$.
	Then 
		$$[F_{n}^{\Sigma}(X),\coprod _{i\in I}Z_{i}]_{Spt}^{\Sigma}\cong \coprod _{i\in I}\: [F_{n}^{\Sigma}(X),Z_{i}]_{Spt}^{\Sigma}
		$$
\end{lem}
\begin{proof}
	Since every pointed simplicial presheaf in $\pointedmotivic$ is cofibrant and
	$F_{n}^{\Sigma}=V\circ F_{n}$ (see proposition \ref{prop.2.6.symmetrizationfunctorV})
	is a left Quillen functor,
	using theorem
	\ref{thm.2.6.T-spectra===symmTspectra} we have:
		\begin{eqnarray*}
			[F_{n}^{\Sigma}(X),\coprod _{i\in I}Z_{i}]_{Spt}^{\Sigma} &=& [V(F_{n}(X)),\coprod _{i\in I}Z_{i}]_{Spt}^{\Sigma}\\
			&\cong & [F_{n}(X),UR_{\Sigma}(\coprod _{i\in I}Z_{i})]_{Spt}\\
			&\cong & [F_{n}(X),(\coprod _{i\in I}UR_{\Sigma}Z_{i})]_{Spt}
		\end{eqnarray*}
	where the last isomorphism follows from theorem \ref{thm.3.3.symmetrization-functor-exact-equivalence}, 
	which implies in particular that $UR_{\Sigma}:\symmstablehomotopy \rightarrow \stablehomotopy$
	is a left adjoint, since it is an equivalence of categories.
	Now since $X\in \pointedmotivic$ is compact in the sense of Jardine,
	lemma \ref{lem.2.4.compact-respects-coproducts-spectra} implies that:
		\begin{eqnarray*}
			[F_{n}(X),\coprod _{i\in I}UR_{\Sigma}Z_{i}]_{Spt} \cong
		  \coprod _{i\in I}\: [F_{n}(X),UR_{\Sigma}Z_{i}]_{Spt}
		\end{eqnarray*}
	Finally using proposition \ref{prop.2.6.symmetrizationfunctorV}
	and theorem \ref{thm.2.6.T-spectra===symmTspectra} again, we get:
		\begin{eqnarray*}
			[F_{n}^{\Sigma}(X),\coprod _{i\in I}Z_{i}]_{Spt}^{\Sigma} &\cong &
			  \coprod _{i\in I}\: [F_{n}(X),UR_{\Sigma}Z_{i}]_{Spt}\\ &\cong & \coprod _{i\in I}\: [V(F_{n}(X)),Z_{i}]_{Spt}^{\Sigma}\\
				& =& \coprod _{i\in I}\: [F_{n}^{\Sigma}(X),Z_{i}]_{Spt}^{\Sigma}
		\end{eqnarray*}
	as we wanted.
\end{proof}

\begin{prop}
		\label{prop.3.3.symmstablehomotopy=>compactly-generated}
	The motivic symmetric stable homotopy category $\symmstablehomotopy$ is
	a \emph{compactly generated} triangulated category
	in the sense of Neeman (see \cite[definition 1.7]{MR1308405}).  The set
	of compact generators is given by
	(see definition \ref{def.2.6.Fnsigma-Evn-adjunction}):				
		$$\symmcompactgenerators
		$$
	i.e. the smallest triangulated subcategory of $\symmstablehomotopy$ 
	closed under small  coproducts and containing
	all the objects in $C^{\Sigma}$ coincides with $\symmstablehomotopy$.
\end{prop}
\begin{proof}
	Since $\symmstablehomotopy$ is closed under small coproducts,
	we just need to prove the following two claims:
	\begin{enumerate}
		\item	\label{claim-Neeman1-symm}For every $\symmgeneratorNRS \in C^{\Sigma}$; $\symmgeneratorNRS$ commutes with coproducts
					in $\symmstablehomotopy$, i.e. given a family of symmetric $T$-spectra $\{ X_{i}\}_{i\in I}$
					indexed by a set $I$ we have:
						$$[\symmgeneratorNRS ,\coprod _{i\in I}X_{i}]_{Spt}^{\Sigma}
							\cong \coprod _{i\in I}[\symmgeneratorNRS ,X_{i}]_{Spt}^{\Sigma}
						$$
		\item	\label{claim-Neeman2-symm}If a symmetric $T$-spectrum $X$ has the following property:
					$[\symmgeneratorNRS ,X]_{Spt}^{\Sigma}=0$
					for every $\symmgeneratorNRS \in C^{\Sigma}$, then $X\cong \ast$ in $\symmstablehomotopy$.
	\end{enumerate}
	
	(\ref{claim-Neeman1}):  Follows immediately from lemma \ref{lem.3.3.compact-respects-coproducts-symmspectra}
	since we know by proposition \ref{prop.2.4.T=compact} that
	the pointed simplicial presheaves $S^{r}\wedge \gm ^{s}\wedge U_{+}$ are all
	compact in the sense of Jardine.
		
	(\ref{claim-Neeman2}):	Fix $\generatorNRS \in C\subseteq \motivicTspectra$. Using
	theorem \ref{thm.2.6.T-spectra===symmTspectra} we have that:
		$$[\generatorNRS , UX]_{Spt} \cong [\symmgeneratorNRS , X]_{Spt}^{\Sigma} =0 
		$$
	Therefore, proposition \ref{prop.3.1.stablehomotopy=>compactly-generated} 
	implies that the map $UX\rightarrow U(\ast)=\ast$ is a weak
	equivalence in $\motivicTspectra$.
	Hence, \cite[proposition 4.8]{MR1787949} implies that $X\rightarrow \ast$ is also a weak equivalence in
	$\motivicsymmTspectra$, i.e.
	$X\cong \ast$ in $\symmstablehomotopy$.  This finishes the proof.
\end{proof}

\begin{cor}
		\label{cor.3.3.detecting-isos-in-symmSH}
	Let $f:X\rightarrow Y$ be a map in $\symmstablehomotopy$.
	Then $f$ is an isomorphism if and only if
	$f$ induces an isomorphism of abelian groups:
		$$\xymatrix{[\symmgeneratorNRS , X]_{Spt}^{\Sigma} \ar[r]^-{f_{\ast}}& [\symmgeneratorNRS , Y]_{Spt}^{\Sigma}}
		$$
	for every $\symmgeneratorNRS \in C^{\Sigma}$.
\end{cor}
\begin{proof}
	($\Rightarrow$):  If $f$ is an isomorphism in $\symmstablehomotopy$ it is clear
	that the induced maps $f_{\ast}$ are isomorphisms of abelian groups for every
	$\symmgeneratorNRS \in C^{\Sigma}$.
	
	($\Leftarrow$):  Complete $f$ to a distinguished triangle in $\symmstablehomotopy$:
		$$\xymatrix{X \ar[r]^-{f}& Y \ar[r]^-{g}& Z \ar[r]^-{h}& \Tsuspfunctor ^{1,0}X}
		$$
	Then $f$ is an isomorphism if and only if
	$Z\cong \ast$ in $\symmstablehomotopy$.
	
	Now since the functor $[\symmgeneratorNRS ,-]_{Spt}^{\Sigma}$
	is homological, we get the following long exact sequence of abelian groups:
		$$\xymatrix{\vdots \ar[d] &\\
								[\symmgeneratorNRS ,X]_{Spt}^{\Sigma} \ar[d]^-{f_{\ast}} & \\ 
								[\symmgeneratorNRS ,Y]_{Spt}^{\Sigma} \ar[d]^-{g_{\ast}} & \\ 
								[\symmgeneratorNRS ,Z]_{Spt}^{\Sigma} \ar[d]^-{h_{\ast}} & \\
								[\symmgeneratorNRS ,\Sigma _{T}^{1,0}X]_{Spt}^{\Sigma} 
								\ar[d]^{\Sigma _{T}^{1,0}f_{\ast}} & \ar[l]_-{\cong}^-{\Tsuspfunctor ^{1,0}} 
								[F_{n+1}^{\Sigma}(S^{r}\wedge \gm ^{s+1}\wedge U_{+}),X]_{Spt}^{\Sigma} \ar[d]^-{f_{\ast}}\\
								[\symmgeneratorNRS ,\Sigma _{T}^{1,0}Y]_{Spt}^{\Sigma} 
								\ar[d] & [F_{n+1}^{\Sigma}(S^{r}\wedge \gm ^{s+1}\wedge U_{+}),Y]_{Spt}^{\Sigma} 
								\ar[l]_-{\cong}^-{\Tsuspfunctor ^{1,0}}\\
								\vdots}
		$$
	But by hypothesis all the maps $f_{\ast}$ are isomorphisms, therefore
	$[\symmgeneratorNRS ,Z]_{Spt}^{\Sigma}=0$
	for every $\symmgeneratorNRS \in C^{\Sigma}$.  Since $\symmstablehomotopy$ is a compactly generated
	triangulated category (see proposition \ref{prop.3.3.symmstablehomotopy=>compactly-generated})
	with set of compact generators $C^{\Sigma}$, we have that
	$Z\cong \ast$.  This implies that $f$ is an isomorphism, as we
	wanted.
\end{proof}

\begin{thm}
		\label{thm.3.3.symmetricconnective-model-structure}
	Fix $q\in \mathbb Z$.  Consider the following set of objects in $\motivicsymmTspectra$
	(see theorem \ref{thm.3.2.connective-model-structure}):
		$$\qeffsymmcompactgenerators
		$$
	Then the right Bousfield localization of $\motivicsymmTspectra$
	with respect to the class of $C_{eff}^{q,\Sigma}$-colocal equivalences exists
	(see definitions \ref{def1.2.cellularization} and \ref{def1.2.1.rightBousloc}).
	This model structure 
	will be called  \emph{($q-1$)-connected motivic symmetric stable}, and
	the category of symmetric $T$-spectra equipped with the ($q-1$)-connected motivic symmetric
	stable model structure will be denoted by
	$\qconnectedsymmTspectra$.  Furthermore
	$\qconnectedsymmTspectra$ is a right proper and simplicial model
	category.
	The homotopy category associated to $\qconnectedsymmTspectra$
	will be denoted by $\qconnectedsymmstablehomotopy$.
\end{thm}
\begin{proof}
	Theorems \ref{thm2.6.stablesymmetricmodelstructure} and
	\ref{thm.2.7.cellularity-motivicsymm-stable-str} imply that
	$\motivicsymmTspectra$ is a
	cellular, proper and simplicial model category.
	Therefore we can apply theorem 5.1.1 in
	\cite{MR1944041} to construct
	the right Bousfield localization of $\motivicsymmTspectra$ with respect to the class of $C_{eff}^{q,\Sigma}$-colocal
	equivalences.
	Using theorem 5.1.1 in
	\cite{MR1944041} again, we have
	that this new model structure  is  right proper and
	simplicial.
\end{proof}

\begin{defi}
		\label{def.3.3.Cqsigma-cofibrant-replacement}
	Fix $q\in \mathbb Z$.  Let $C_{q}^{\Sigma}$ denote a cofibrant replacement functor in 
	$\qconnectedsymmTspectra$; such that for every symmetric $T$-spectrum $X$, the natural map
		$$\xymatrix{C_{q}^{\Sigma}X \ar[r]^-{C_{q}^{\Sigma ,X}}& X}
		$$ 
	is a trivial fibration in $\qconnectedsymmTspectra$, and
	$C_{q}^{\Sigma}X$ is always $C_{eff}^{q,\Sigma}$-colocal in $\motivicsymmTspectra$.
\end{defi}

\begin{prop}
		\label{prop.3.3.Rsigma-fibrant-replacement-all-Rq}
	Fix $q\in \mathbb Z$.  Then $R_{\Sigma}$ is also a fibrant
	replacement functor in $\qconnectedsymmTspectra$
	(see definition \ref{def.3.3.fibrantstable-replacementfunctor}),
	and for every symmetric $T$-spectrum $X$ the natural map
		$$\xymatrix{X\ar[r]^-{R_{\Sigma}^{X}}& R_{\Sigma}X}
		$$
	is a trivial cofibration in $\qconnectedsymmTspectra$.
\end{prop}
\begin{proof}
	Since $\qconnectedsymmTspectra$ is the right Bousfield localization of
	$\motivicsymmTspectra$ with respect to the $C_{eff}^{q,\Sigma}$-colocal equivalences, by construction
	we have that the fibrations and the trivial cofibrations are indentical in
	$\qconnectedsymmTspectra$ and $\motivicsymmTspectra$ respectively.  This implies that for every symmetric
	$T$-spectrum $X$, $R_{\Sigma}X$ is fibrant in $\qconnectedsymmTspectra$,
	and we also have that
	the natural map
		$$\xymatrix{X\ar[r]^{R_{\Sigma}^{X}}& R_{\Sigma}X}
		$$
	is a trivial cofibration in $\qconnectedsymmTspectra$.  Hence $R_{\Sigma}$ is also a fibrant replacement
	functor for $\qconnectedsymmTspectra$.
\end{proof}

\begin{prop}
		\label{prop.3.3.f-Cqeff-colocal=iff=Uf-Cqeff-colocal}
	Fix $q\in \mathbb Z$.  Then
	a map of symmetric $T$-spectra $f:X\rightarrow Y$ is a $C^{q,\Sigma}_{eff}$-colocal equivalence
	in $\motivicsymmTspectra$
	if and only if the underlying map
	$UR_{\Sigma}(f):UR_{\Sigma}X\rightarrow UR_{\Sigma}Y$ is a $C^{q}_{eff}$-colocal equivalence
	in $\motivicTspectra$.
\end{prop}
\begin{proof}
	Consider $\symmgeneratorNRS \in C_{eff}^{q,\Sigma}$.
	Using the enriched adjunctions of proposition
	\ref{prop.2.6.enriched-symmetrization-adjunction}, 
	we get the following commutative diagram where all
	the vertical arrows are isomorphisms:
		$$\xymatrix{Map\: _{\Sigma}(\symmgeneratorNRS ,R_{\Sigma}X) \ar[r]^-{R_{\Sigma}f_{\ast}} \ar@{=}[d]
								& Map\: _{\Sigma}(\symmgeneratorNRS ,R_{\Sigma}Y) \ar@{=}[d]\\
								Map\: _{\Sigma}(V(\generatorNRS) ,R_{\Sigma}X) \ar[r]^-{R_{\Sigma}f_{\ast}} \ar[d]_-{\cong}
								& Map\: _{\Sigma}(V(\generatorNRS) ,R_{\Sigma}Y) \ar[d]^-{\cong}\\
								Map (\generatorNRS ,UR_{\Sigma}X) \ar[r]_-{UR_{\Sigma}f_{\ast}}& Map (\generatorNRS ,UR_{\Sigma}Y)}
		$$
	Since $UR_{\Sigma}X$ and $UR_{\Sigma}Y$ are both fibrant in $\motivicTspectra$, we have that $UR_{\Sigma}(f)$
	is a $C_{eff}^{q}$-colocal equivalence in $\motivicTspectra$ if and only if
	the bottom row in the diagram above is a weak equivalence of simplicial sets
	for every $\generatorNRS \in C_{eff}^{q}$.  By the two out of three property for
	weak equivalences we have that this happens if and only if the top
	row in the diagram above is a weak equivalence for every $\symmgeneratorNRS \in C_{eff}^{q,\Sigma}$.
	But this last condition holds if and only if $f$ is a $C_{eff}^{q,\Sigma}$-colocal equivalence
	in $\motivicsymmTspectra$.  This finishes the proof.
\end{proof}

\begin{prop}
		\label{prop.3.3.detecting-symmetric-Cqeff-colocal-equivalences}
	Fix $q\in \mathbb Z$, and let $f:X\rightarrow Y$
	be a map of symmetric $T$-spectra.  Then $f$
	is a $C_{eff}^{q,\Sigma}$-colocal equivalence in $\motivicsymmTspectra$ if and only if
	for every $\symmgeneratorNRS \in C_{eff}^{q,\Sigma}$, the induced map:
		$$\xymatrix{[\symmgeneratorNRS , X]_{Spt}^{\Sigma} \ar[r]^-{f_{\ast}}& [\symmgeneratorNRS , Y]_{Spt}^{\Sigma}}
		$$
	is an isomorphism of abelian groups.
\end{prop}
\begin{proof}
	By proposition \ref{prop.3.3.f-Cqeff-colocal=iff=Uf-Cqeff-colocal}, $f$ is a $C_{eff}^{q,\Sigma}$-colocal
	equivalence in $\motivicsymmTspectra$ if and only if $UR_{\Sigma}(f)$ is a $C_{eff}^{q}$-colocal equivalence
	in $\motivicTspectra$.  Using proposition \ref{prop.3.2.classif-Cqcolocalequivs} 
	we have that $UR_{\Sigma}(f)$ is a $C_{eff}^{q}$-colocal equivalence if and
	only if for every $\generatorNRS \in C_{eff}^{q}$, the induced map
		$$\xymatrix{[\generatorNRS , UR_{\Sigma}X]_{Spt} \ar[rr]^-{UR_{\Sigma}(f)_{\ast}}&& [\generatorNRS , UR_{\Sigma}Y]_{Spt}}
		$$
	is an isomorphism of abelian groups.
	
	Now theorem \ref{thm.2.6.T-spectra===symmTspectra} implies that we have 
	the following commutative diagram, where all the vertical arrows
	are isomorphisms:
		$$\xymatrix{[\generatorNRS , UR_{\Sigma}X]_{Spt} \ar[rr]^-{UR_{\Sigma}(f)_{\ast}} \ar[d]_-{\cong}&&
								[\generatorNRS , UR_{\Sigma}Y]_{Spt} \ar[d]^-{\cong}\\
								[V(\generatorNRS), X]_{Spt}^{\Sigma} \ar[rr]^-{f_{\ast}} \ar@{=}[d]&& 
								[V(\generatorNRS), Y]_{Spt}^{\Sigma} \ar@{=}[d]\\
								[\symmgeneratorNRS, X]_{Spt}^{\Sigma} \ar[rr]^-{f_{\ast}}&& [\symmgeneratorNRS, Y]_{Spt}^{\Sigma}}
		$$
	Therefore $f$ is a $C_{eff}^{q,\Sigma}$-colocal equivalence if and only if for every 
	$\symmgeneratorNRS \in C_{eff}^{q,\Sigma}$, the bottom row is an isomorphism of abelian groups.
	This finishes the proof.
\end{proof}

\begin{lem}
		\label{lem.3.3.Cqsymmetric-colocalequivs-stable-S1desuspension}
	Fix $q\in \mathbb Z$, and let $f:X\rightarrow Y$ be a map of symmetric $T$-spectra.
	Then $f$ is a
	$C^{q,\Sigma}_{eff}$-colocal equivalence in $\motivicsymmTspectra$ if and only if
	$\Omega _{S^{1}} R_{\Sigma}f$ is a $C^{q,\Sigma}_{eff}$-colocal equivalence in $\motivicsymmTspectra$.
\end{lem}
\begin{proof}
	It follows from proposition \ref{prop.3.3.f-Cqeff-colocal=iff=Uf-Cqeff-colocal} that $f$ is a $C_{eff}^{q,\Sigma}$-colocal
	equivalence in $\motivicsymmTspectra$ if and only if $UR_{\Sigma}f$ is a $C_{eff}^{q}$-colocal equivalence
	in $\motivicTspectra$.  Since $UR_{\Sigma}X , UR_{\Sigma}Y$ are both fibrant in $\motivicTspectra$, using
	lemma \ref{lem.3.2.Cq-colocalequivs-stable-S1desuspension} 
	we have that $UR_{\Sigma}f$ is a $C_{eff}^{q}$-colocal equivalence if and
	only if $\Omega _{S^{1}}UR_{\Sigma}f=U(\Omega_{S^{1}}R_{\Sigma}f)$ is a
	$C_{eff}^{q}$-colocal equivalence in $\motivicTspectra$.
	
	Finally, since $\Omega _{S^{1}}R_{\Sigma}X , \Omega _{S^{1}}R_{\Sigma}Y$ are both
	fibrant in $\motivicsymmTspectra$,
	we have by proposition \ref{prop.3.3.f-Cqeff-colocal=iff=Uf-Cqeff-colocal} that
	$U(\Omega _{S^{1}}R_{\Sigma}f)$ is a $C_{eff}^{q}$-colocal equivalence if and only if
	$\Omega _{S^{1}}R_{\Sigma}f$ is a $C_{eff}^{q,\Sigma}$-colocal equivalence.
	This finishes the proof.
\end{proof}

\begin{cor}
		\label{cor.3.3.Suspension=>qconnsymm-Quillen-equiv}
	Fix $q\in \mathbb Z$.  Then the adjunction
		$$\xymatrix{(-\wedge S^{1},\Omega _{S^{1}},\varphi):\qconnectedsymmTspectra \ar[r]& 
								\qconnectedsymmTspectra}
		$$
	is a Quillen equivalence.
\end{cor}
\begin{proof}
	Using corollary 1.3.16 in \cite{MR1650134} and
	proposition \ref{prop.3.3.Rsigma-fibrant-replacement-all-Rq} we have that it suffices to
	verify the following two conditions:
		\begin{enumerate}
			\item	\label{cor.3.3.Suspension=>qconnsymm-Quillen-equiv.a}For every cofibrant object $X$
						in $\qconnectedsymmTspectra$, the following composition
							$$\xymatrix{X\ar[r]^-{\eta _{X}}& \Omega _{S^{1}}(X\wedge S^{1}) \ar[rr]^-{\Omega _{S^{1}}R_{\Sigma}^{X \wedge S^{1}}}&& 
								\Omega _{S^{1}}R_{\Sigma}(X\wedge S^{1})}
							$$
						is a $C_{eff}^{q,\Sigma}$-colocal equivalence.
			\item	\label{cor.3.3.Suspension=>qconnsymm-Quillen-equiv.b}$\Omega _{S^{1}}$ reflects $C_{eff}^{q,\Sigma}$-colocal equivalences
						between fibrant objects in $\qconnectedsymmTspectra$.
		\end{enumerate}
		
	(\ref{cor.3.3.Suspension=>qconnsymm-Quillen-equiv.a}):  By construction $\qconnectedsymmTspectra$ is
	a right Bousfield localization of $\motivicsymmTspectra$, therefore the identity
	functor 
		$$\xymatrix{id:\qconnectedsymmTspectra \ar[r]& \motivicsymmTspectra}
		$$ 
	is a left Quillen functor.
	Thus $X$ is also cofibrant in $\motivicsymmTspectra$.
	Since the adjunction $(-\wedge S^{1},\Omega _{S^{1}},\varphi)$
	is a Quillen equivalence on $\motivicsymmTspectra$, 
	\cite[proposition 1.3.13(b)]{MR1650134} implies that the following composition
	is a weak equivalence in $\motivicsymmTspectra$:
		$$\xymatrix{X\ar[r]^-{\eta _{X}}& \Omega _{S^{1}}(X\wedge S^{1}) \ar[rr]^-{\Omega _{S^{1}}R_{\Sigma}^{X\wedge S^{1}}}&& 
								\Omega _{S^{1}}R_{\Sigma}(X\wedge S^{1})}
		$$
	Hence using \cite[proposition 3.1.5]{MR1944041} it follows that the composition above
	is a $C_{eff}^{q,\Sigma}$-colocal equivalence.
	
	(\ref{cor.3.3.Suspension=>qconnsymm-Quillen-equiv.b}):  This follows immediately
	from proposition \ref{prop.3.3.Rsigma-fibrant-replacement-all-Rq} and 
	lemma \ref{lem.3.3.Cqsymmetric-colocalequivs-stable-S1desuspension}.
\end{proof}

\begin{rmk}
		\label{rmk.3.3.smashT-not-symmdescending}
	The adjunction $(\Tsuspfunctor ,\Tloops ,\varphi)$ is
	a Quillen equivalence on $\motivicsymmTspectra$.
	However it does not descend even to a Quillen
	adjunction on the ($q-1$)-connected motivic symmetric stable
	model category
	$\qconnectedsymmTspectra$.
\end{rmk}

\begin{cor}
		\label{cor.3.3.symmqconn=>triangcat}
	Fix $q\in \mathbb Z$.  Then $\qconnectedsymmstablehomotopy$
	has the structure of a triangulated category.
\end{cor}
\begin{proof}
	Theorem \ref{thm.3.3.symmetricconnective-model-structure} implies in particular
	that $\qconnectedsymmTspectra$ is a pointed simplicial model category,
	and corollary \ref{cor.3.3.Suspension=>qconnsymm-Quillen-equiv} implies that
	the adjunction 
		$$(-\wedge S^{1},\Omega _{S^{1}},\varphi):\qconnectedsymmTspectra \rightarrow \qconnectedsymmTspectra$$
	is a Quillen equivalence.  Therefore
	the result follows from the work of Quillen in 
	\cite[sections I.2 and I.3]{MR0223432} and the work of 
	Hovey in \cite[chapters VI and VII]{MR1650134}.
\end{proof}

\begin{prop}
		\label{prop.3.3.symmcofibrant-replacement=>triangulatedfunctor}
	For every $q\in \mathbb Z$, we have the following adjunction
		$$\xymatrix{(C_{q}^{\Sigma}, R_{\Sigma}, \varphi) :\qconnectedsymmstablehomotopy \ar[r]& \symmstablehomotopy}
		$$
	between exact functors of triangulated categories.
\end{prop}
\begin{proof}
	Since $\qconnectedsymmTspectra$ is the right Bousfield localization of
	$\motivicsymmTspectra$ with respect to the $C^{q,\Sigma}_{eff}$-colocal equivalences, we have that the
	identity functor $id:\qconnectedsymmTspectra \rightarrow \motivicsymmTspectra$
	is a left Quillen functor.  Therefore we get
	the following adjunction at the level of the associated homotopy categories:
		$$\xymatrix{(C_{q}^{\Sigma}, R_{\Sigma}, \varphi ):\qconnectedsymmstablehomotopy \ar[r]& \symmstablehomotopy}
		$$
	
	Now proposition 6.4.1 in \cite{MR1650134} implies that
	$C_{q}^{\Sigma}$ maps cofibre sequences in $\qconnectedsymmstablehomotopy$ to cofibre sequences in
	$\symmstablehomotopy$.
	Therefore 
	using proposition 7.1.12 in \cite{MR1650134} we have
	that $C_{q}^{\Sigma}$ and $R_{\Sigma}$ are both exact functors between triangulated categories.
\end{proof}

\begin{thm}
		\label{thm.3.3.symmetrization-qconnected-Quillen-equiv}
	Fix $q \in \mathbb Z$.  Then the adjunction
		$$\xymatrix{(V,U,\varphi):\qconnectedTspectra \ar[r]& \qconnectedsymmTspectra}
		$$
	given by the symmetrization and the forgetful functors is a
	Quillen equivalence.
\end{thm}
\begin{proof}
	Proposition \ref{prop.3.3.f-Cqeff-colocal=iff=Uf-Cqeff-colocal} 
	together with the universal property for right Bousfield localizations
	(see definition \ref{def-1.rightlocmodcats})
	imply that 
		$$\xymatrix{U:\qconnectedsymmTspectra \ar[r]& \qconnectedTspectra}$$ 
	is a right
	Quillen functor.
	Using corollary 1.3.16 in \cite{MR1650134} and
	proposition \ref{prop.3.3.Rsigma-fibrant-replacement-all-Rq} we have that it suffices to
	verify the following two conditions:
		\begin{enumerate}
			\item	\label{thm.3.3.symmetrization-qconnected-Quillen-equiv.a}For every cofibrant object $X$
						in $\qconnectedTspectra$, the following composition
							$$\xymatrix{X\ar[r]^-{\eta _{X}}& UV(X) \ar[rr]^-{UR_{\Sigma}^{VX}}&& 
								UR_{\Sigma}V(X)}
							$$
						is a weak equivalence in $\qconnectedTspectra$.
			\item	\label{thm.3.3.symmetrization-qconnected-Quillen-equiv.b}$U$ reflects weak equivalences
						between fibrant objects in $\qconnectedsymmTspectra$.
		\end{enumerate}
		
	(\ref{thm.3.3.symmetrization-qconnected-Quillen-equiv.a}):  
	By construction $\qconnectedTspectra$ is
	a right Bousfield localization of $\motivicTspectra$, therefore the identity
	functor 
		$$\xymatrix{id:\qconnectedTspectra \ar[r]& \motivicTspectra}
		$$ 
	is a left Quillen functor.
	Thus $X$ is also cofibrant in $\motivicTspectra$.
	Since the adjunction $(V,U,\varphi)$
	is a Quillen equivalence between $\motivicTspectra$  and $\motivicsymmTspectra$, 
	\cite[proposition 1.3.13(b)]{MR1650134} implies that the following composition
	is a weak equivalence in $\motivicTspectra$:
		$$\xymatrix{X\ar[r]^-{\eta _{X}}& UV (X) \ar[rr]^-{UR_{\Sigma}^{V X}}&& 
								UR_{\Sigma}V(X)}
		$$
	Hence using \cite[proposition 3.1.5]{MR1944041} it follows that the composition above
	is a $C_{eff}^{q}$-colocal equivalence in $\motivicTspectra$, i.e. a weak equivalence
	in $\qconnectedTspectra$.
	
	(\ref{thm.3.3.symmetrization-qconnected-Quillen-equiv.b}):  This follows immediately
	from propositions \ref{prop.3.3.Rsigma-fibrant-replacement-all-Rq} and 
	\ref{prop.3.3.f-Cqeff-colocal=iff=Uf-Cqeff-colocal}.
\end{proof}

\begin{cor}
		\label{cor.3.3.symmetrization-functor-exact-equivalenceRq}
	Fix $q\in \mathbb Z$.  Then the adjunction
		$$\xymatrix{(V,U,\varphi):\qconnectedTspectra \ar[r]& \qconnectedsymmTspectra}
		$$
	given by the symmetrization and the forgetful functors, induces
	an adjunction
		$$\xymatrix{(VC_{q},UR_{\Sigma},\varphi):\qconnectedstablehomotopy \ar[r]& \qconnectedsymmstablehomotopy}
		$$
	of exact funtors between triangulated categories.
	Furthermore, $VC_{q}$ and $UR_{\Sigma}$ are both equivalences of categories.
\end{cor}
\begin{proof}
	Theorem \ref{thm.3.3.symmetrization-qconnected-Quillen-equiv} implies that the adjunction $(V,U,\varphi)$
	is a Quillen equivalence.
	Therefore we get
	the following adjunction at the level of the associated homotopy categories:
		$$\xymatrix{(VC_{q}, UR_{\Sigma}, \varphi ):\qconnectedstablehomotopy \ar[r]& \qconnectedsymmstablehomotopy}
		$$
	
	Now \cite[proposition 1.3.13]{MR1650134} implies that $VC_{q} ,UR_{\Sigma}$ are both equivalences of
	categories.
	Finally, proposition \ref{prop.2.6.enriched-symmetrization-adjunction}
	together with \cite[proposition 6.4.1]{MR1650134} imply that
	$VC_{q}$ maps cofibre sequences in $\qconnectedstablehomotopy$ to cofibre sequences in
	$\qconnectedsymmstablehomotopy$.
	Therefore 
	using proposition 7.1.12 in \cite{MR1650134} we have
	that $VC_{q}$ and $UR_{\Sigma}$ are both exact functors between triangulated categories.
\end{proof}

	Now it is very easy to find the desired lifting 
	for the functor $\tilde{f}_{q}:\symmstablehomotopy \rightarrow \symmstablehomotopy$
	(see corollary \ref{cor.3.3.symmetric-fq,s<q,sq}(\ref{cor.3.3.symmetric-fq,s<q,sq.a})) to the
	model category level.
	
\begin{lem}
		\label{lem.3.3.homotopycoherence==>liftings}
	Fix $q\in \mathbb Z$.
	\begin{enumerate}
		\item \label{lem.3.3.homotopycoherence==>liftings.a}Let $X$ be an arbitrary 
					$T$-spectrum in $\qconnectedTspectra$.
					Then the following maps in $\motivicsymmTspectra$
						$$\xymatrix{VQ_{s}(C_{q}X)\ar[rr]^-{V(Q_{s}^{C_{q}X})}&& 
												VC_{q}X && C_{q}^{\Sigma}(VC_{q}X) \ar[ll]_-{C_{q}^{\Sigma ,VC_{q}X}}}
						$$
					induce natural isomorphisms between the functors:
						$$C_{q}^{\Sigma}\circ VC_{q}, VC_{q}, VQ_{s}\circ C_{q}:\qconnectedstablehomotopy
							\rightarrow \symmstablehomotopy$$
						
						$$\xymatrix{& \qconnectedsymmstablehomotopy \ar[dr]^-{C_{q}^{\Sigma}}&\\
												 \qconnectedstablehomotopy
												\ar[rr]^{VC_{q}} \ar[ur]^-{VC_{q}} \ar[dr]_-{C_{q}}&& \symmstablehomotopy \\
												& \stablehomotopy \ar[ur]_-{VQ_{s}}&}
						$$
					Given a $T$-spectrum $X$
					$$\xymatrix{\alpha _{X}:VQ_{s}(C_{q}X)\ar[r]^-{\cong}& C_{q}^{\Sigma}(VC_{q}X)}$$ 
					will denote
					the isomorphism in $\symmstablehomotopy$ corresponding to the natural isomorphism
					between $VQ_{s}\circ C_{q}$
					and $C_{q}^{\Sigma}\circ VC_{q}$.
		\item \label{lem.3.3.homotopycoherence==>liftings.b}Let $X$ be an arbitrary 
					symmetric $T$-spectrum.
					Then the following maps in $\qconnectedTspectra$
						$$\xymatrix{IQ_{T}J(UR_{\Sigma}X)&& 
												UR_{\Sigma}X \ar[ll]_-{IQ_{T}J^{UR_{\Sigma}X}} \ar[rr]^-{U(R_{\Sigma}^{R_{\Sigma}X})}&& 
												UR_{\Sigma}(R_{\Sigma}X) }
						$$
					induce natural isomorphisms between the functors:
						$$IQ_{T}J\circ UR_{\Sigma}, UR_{\Sigma}, UR_{\Sigma}\circ R_{\Sigma}:\symmstablehomotopy
							\rightarrow \qconnectedstablehomotopy$$
						
						$$\xymatrix{& \stablehomotopy \ar[dr]^-{IQ_{T}J}&\\
												 \symmstablehomotopy
												\ar[rr]^{UR_{\Sigma}} \ar[ur]^-{UR_{\Sigma}} \ar[dr]_-{R_{\Sigma}}&& \qconnectedstablehomotopy \\
												& \qconnectedsymmstablehomotopy \ar[ur]_-{UR_{\Sigma}}&}
						$$
					Given a symmetric $T$-spectrum $X$
					$$\xymatrix{\beta _{X}:IQ_{T}J(UR_{\Sigma}X)\ar[r]^-{\cong}& UR_{\Sigma}(R_{\Sigma}X)}$$ 
					will denote
					the isomorphism in $\qconnectedstablehomotopy$ corresponding to the natural isomorphism
					between $IQ_{T}J\circ UR_{\Sigma}$
					and $UR_{\Sigma}\circ R_{\Sigma}$.
	\end{enumerate}
\end{lem}	
\begin{proof}
	(\ref{lem.3.3.homotopycoherence==>liftings.a}):  Follows immediately from 
	theorem 1.3.7 in \cite{MR1650134} and the following
	commutative diagram of left Quillen functors:
		$$\xymatrix{\qconnectedTspectra \ar[r]^-{V} \ar[d]_-{id}& \qconnectedsymmTspectra \ar[d]^-{id}\\
								\motivicTspectra \ar[r]_-{V}& \motivicsymmTspectra}
		$$
		
	(\ref{lem.3.3.homotopycoherence==>liftings.b}):  Follows immediately from the dual of
	theorem 1.3.7 in \cite{MR1650134} and the following
	commutative diagram of right Quillen functors:
		$$\xymatrix{\qconnectedTspectra & \qconnectedsymmTspectra \ar[l]_-{U}\\
								\motivicTspectra \ar[u]^-{id}& \motivicsymmTspectra \ar[l]^-{U} \ar[u]_-{id}}
		$$
\end{proof}
	
\begin{thm}
		\label{thm.3.3.symmRq-models-symmf<q}
	Fix $q\in \mathbb Z$, and let $X$ be an arbitrary symmetric $T$-spectrum.
	\begin{enumerate}
		\item \label{thm.3.3.symmRq-models-symmf<q.a}The
					diagram (\ref{diagram.3.2.fq-lifting})
					in theorem \ref{thm.3.2.Rq-models-f<q} induces
					the following diagram in 
					$\symmstablehomotopy$:
						\begin{equation}
										\label{diagram.3.3.symmfq-lifting.a}
							\begin{array}{c}
									\xymatrix{VQ_{s}(IQ_{T}Jf_{q}(UR_{\Sigma}X)) &&& VQ_{s}(C_{q}IQ_{T}Jf_{q}(UR_{\Sigma}X)) 
														\ar[lll]_-{VQ_{s}(C_{q}^{IQ_{T}Jf_{q}(UR_{\Sigma}X)})}^-{\cong}
														\ar[dd]_-{VQ_{s}(C_{q}IQ_{T}J(\theta _{UR_{\Sigma}X}))}^-{\cong}\\ 
														&&&\\ 
														\tilde{f}_{q}X=VQ_{s}(f_{q}(UR_{\Sigma}X)) \ar[uu]_-{VQ_{s}(IQ_{T}J^{f_{q}UR_{\Sigma}X})}^-{\cong}
														&&& VQ_{s}(C_{q}IQ_{T}J(UR_{\Sigma}X)) }
							\end{array}
						\end{equation}
					where all the maps are isomorphisms in $\symmstablehomotopy$.
					Furthermore, this diagram induces a
					natural isomorphism  
					between the following exact functors:
						$$\xymatrix{\symmstablehomotopy \ar@<1ex>[rrr]^-{\tilde{f}_{q}}  
												\ar@<-1ex>[rrr]_-{VQ_{s}\circ C_{q}IQ_{T}J \circ UR_{\Sigma}} &&& \symmstablehomotopy}
						$$
		\item	\label{thm.3.3.symmRq-models-symmf<q.b}Let $\epsilon$ be the
					counit of the adjunction (see corollary \ref{cor.3.3.symmetrization-functor-exact-equivalenceRq}):
						$$\xymatrix{(VC_{q},UR_{\Sigma},\varphi):\qconnectedstablehomotopy \ar[r]&
												\qconnectedsymmstablehomotopy}
						$$						
					Then we have the following diagram in $\symmstablehomotopy$
					(see lemma \ref{lem.3.3.homotopycoherence==>liftings}):
						\begin{equation}
										\label{diagram.3.3.symmfq-lifting.b}
							\begin{array}{c}
									\xymatrix{C_{q}^{\Sigma}(VC_{q}(IQ_{T}J(UR_{\Sigma}X))) \ar[rr]^-{C_{q}^{\Sigma}VC_{q}(\beta _{X})}_-{\cong}&&
														C_{q}^{\Sigma}(VC_{q}(UR_{\Sigma}(R_{\Sigma}X))) \ar[dd]^-{C_{q}^{\Sigma}(\epsilon _{R_{\Sigma}X})}_-{\cong}\\
														&&\\
														VQ_{s}(C_{q}IQ_{T}J(UR_{\Sigma}X))\ar[uu]^-{\alpha _{IQ_{T}J(UR_{\Sigma}X)}}_-{\cong}
														&& C_{q}^{\Sigma}R_{\Sigma}X=f_{q}^{\Sigma}X}
							\end{array}
						\end{equation}
					where all the maps are isomorphisms in $\symmstablehomotopy$.  This diagram induces a natural isomorphism
					between the following exact functors:
						$$\xymatrix{\symmstablehomotopy \ar@<1ex>[rrr]^-{VQ_{s}\circ C_{q}IQ_{T}J\circ UR_{\Sigma}}  
												\ar@<-1ex>[rrr]_-{C_{q}^{\Sigma}R_{\Sigma}=f_{q}^{\Sigma}} &&& \symmstablehomotopy}
						$$
		\item	\label{thm.3.3.symmRq-models-symmf<q.c}Combining the diagrams (\ref{diagram.3.3.symmfq-lifting.a}) and 
					(\ref{diagram.3.3.symmfq-lifting.b}) above
					we get a natural isomorphism 
					between the following exact functors:
						$$\xymatrix{\symmstablehomotopy \ar@<1ex>[r]^-{\tilde{f}_{q}}  
												\ar@<-1ex>[r]_-{f_{q}^{\Sigma}} & \symmstablehomotopy}
						$$
	\end{enumerate}
\end{thm}
\begin{proof}
	It is clear that it suffices to prove only the first two claims.

	(\ref{thm.3.3.symmRq-models-symmf<q.a}):  Follows immediately from theorems
	\ref{thm.3.2.Rq-models-f<q} and \ref{thm.3.3.symmetrization-functor-exact-equivalence}.
	
	(\ref{thm.3.3.symmRq-models-symmf<q.b}):  Follows immediately from lemma \ref{lem.3.3.homotopycoherence==>liftings}
	and corollary \ref{cor.3.3.symmetrization-functor-exact-equivalenceRq}.
\end{proof}

\begin{prop}
		\label{prop.3.3.lifting-map-fq-->idsymm}
	Fix $q\in \mathbb Z$.  Let $\epsilon$ denote the counit of the adjuntion
	$(C_{q}^{\Sigma},R_{\Sigma},\varphi ):\qconnectedsymmstablehomotopy \rightarrow \symmstablehomotopy$
	constructed in proposition \ref{prop.3.3.symmcofibrant-replacement=>triangulatedfunctor}.  Then the natural transformation
	$\theta _{q}:f_{q}\rightarrow id$ (see proposition \ref{prop.3.1.counit-properties}) gets canonically identified,
	through the equivalence of categories $r_{q}C_{q}$, $IQ_{T}Ji_{q}$, $VC_{q}$, $UR_{\Sigma}$,
	$VQ_{s}$ and $UR_{\Sigma}$
	constructed in proposition \ref{prop.3.2.Rq-lifts-qSH}, corollary 
	\ref{cor.3.3.symmetrization-functor-exact-equivalenceRq} and
	theorem \ref{thm.3.3.symmetrization-functor-exact-equivalence};
	with $\theta _{q}^{\Sigma}=\epsilon$.
\end{prop}
\begin{proof}
	By construction $\theta _{q}$ is the counit of the adjunction
		$$(i_{q},r_{q},\varphi ):\stablehomotopyeffq \rightarrow \stablehomotopy$$
	(see proposition \ref{prop.3.1.counit-properties}).
	The result follows immediately from proposition \ref{prop.3.2.Rq-lifts-qSH}, corollary 
	\ref{cor.3.3.symmetrization-functor-exact-equivalenceRq} and theorem
	\ref{thm.3.3.symmetrization-functor-exact-equivalence}.
\end{proof}

	The functor $f_{q}^{\Sigma}$ gives the desired lifting for
	the functor $\tilde{f}_{q}$
	to the model category level, and it will be used in the
	study of the multiplicative properties of Voevodsky's slice filtration.

\begin{prop}
		\label{prop.3.2.functors-between-Rqsymm}
	Fix $q\in \mathbb Z$.
		\begin{enumerate}
			\item	\label{prop.3.2.functors-between-Rqsymm.a}We have the following
						commutative diagram of left Quillen functors:
							\begin{equation}
										\label{diagram3.3.liftslicefiltrationsymm}
								\begin{array}{c}
									\xymatrix{\qplusoneconnectedsymmTspectra \ar[rr]^-{id} \ar[dr]_-{id}&& 
														\qconnectedsymmTspectra \ar[dl]^-{id}\\
														&\motivicsymmTspectra  &}
								\end{array}
							\end{equation}
			\item	\label{prop.3.2.functors-between-Rqsymm.b}For every symmetric $T$-spectrum $X$, the natural map:
							$$\xymatrix{C_{q}^{\Sigma}C_{q+1}^{\Sigma}X \ar[rr]^-{C_{q}^{\Sigma ,C_{q+1}^{\Sigma}X}} && C_{q+1}^{\Sigma}X}
							$$
						is a weak equivalence in $\symmstablehomotopy$, and it induces a natural equivalence
						$C_{q}^{\Sigma ,C_{q+1}^{\Sigma}-}:C_{q}^{\Sigma}\circ C_{q+1}^{\Sigma}\rightarrow C_{q+1}^{\Sigma}$
						between the following functors:
							$$\xymatrix{\qplusoneconnectedsymmstablehomotopy \ar[rr]^-{C_{q+1}^{\Sigma}} \ar[dr]_-{C_{q+1}^{\Sigma}}&& 
													\qconnectedsymmstablehomotopy  \ar[dl]^-{C_{q}^{\Sigma}}\\
													& \symmstablehomotopy &}
							$$
			\item	\label{prop.3.2.functors-between-Rqsymm.c}The natural transformation $f_{q+1}X\rightarrow f_{q}X$
						(see theorem \ref{thm.3.1.slicefiltration}(\ref{thm.3.1.slicefiltration.a})) gets canonically identified,
						through the equivalence of categories $r_{q}C_{q}$, $IQ_{T}Ji_{q}$, $VC_{q}$ and $UR_{\Sigma}$
						constructed in proposition \ref{prop.3.2.Rq-lifts-qSH} and corollary 
						\ref{cor.3.3.symmetrization-functor-exact-equivalenceRq};
						with the following composition $\rho _{q}^{X}:f_{q+1}^{\Sigma}X\rightarrow f_{q}^{\Sigma}X$
						in $\symmstablehomotopy$
							$$\xymatrix{ & C_{q}^{\Sigma}C_{q+1}^{\Sigma}R_{\Sigma}X \ar[dr]^-{C_{q}^{\Sigma}(C_{q+1}^{\Sigma ,R_{\Sigma}X})}&\\
													C_{q+1}^{\Sigma}R_{\Sigma}X \ar[ur]^-{(C_{q}^{\Sigma ,C_{q+1}^{\Sigma}R_{\Sigma}X})^{-1}}&& C_{q}^{\Sigma}R_{\Sigma}X}
							$$
						which is induced by the following commutative diagram in $\motivicsymmTspectra$
							\begin{equation}
										\label{diagram3.3.Counit-connected}
								\begin{array}{c}
									\xymatrix{C_{q}^{\Sigma}C_{q+1}^{\Sigma}R_{\Sigma}X \ar[d]_-{C_{q}^{\Sigma ,C_{q+1}^{\Sigma}R_{\Sigma}X}} 
														\ar[rr]^-{C_{q}^{\Sigma}(C_{q+1}^{\Sigma ,R_{\Sigma}X})} && C_{q}^{\Sigma}R_{\Sigma}X 
														\ar[d]^-{C_{q}^{\Sigma ,R_{\Sigma}X}}\\
														C_{q+1}^{\Sigma}R_{\Sigma}X \ar[rr]_-{C_{q+1}^{\Sigma ,R_{\Sigma}X}} && R_{\Sigma}X}
								\end{array}
							\end{equation}
		\end{enumerate}
\end{prop}

\begin{proof}
	(\ref{prop.3.2.functors-between-Rqsymm.a}):  Since $\qplusoneconnectedsymmTspectra$ and $\qconnectedsymmTspectra$ are both
	right Bousfield localizations of $\motivicsymmTspectra$, by construction
	the identity functor 
		$$id: \qplusoneconnectedsymmTspectra \rightarrow \motivicsymmTspectra$$
		$$id:\qconnectedsymmTspectra \rightarrow \motivicsymmTspectra$$ 
	is in both cases
	a left Quillen functor.  To finish the proof,
	it suffices to show
	that the identity functor
		$$id:\qconnectedsymmTspectra \rightarrow \qplusoneconnectedsymmTspectra$$
	is a right Quillen functor.  Using the universal property of right
	Bousfield localizations (see definition \ref{def-1.rightlocmodcats}), it is enough to check that
	if $f:X\rightarrow Y$ is
	a $C_{eff}^{q,\Sigma}$-colocal equivalence in $\motivicsymmTspectra$ then
	$R_{\Sigma}f$ is  a
	$C_{eff}^{q+1,\Sigma}$-colocal equivalence.
	But since $R_{\Sigma}X$ and $R_{\Sigma}Y$ are already fibrant in $\motivicsymmTspectra$, we have that
	$R_{\Sigma}(f)$ is a $C_{eff}^{q+1,\Sigma}$-colocal equivalence if and only if
	for every $\symmgeneratorNRS \in C_{eff}^{q+1,\Sigma}$, the induced map:
		$$\xymatrix{Map\: _{\Sigma}(\symmgeneratorNRS ,R_{\Sigma}X) \ar[d]^-{(R_{\Sigma}f)_{\ast}}\\ 
								Map\: _{\Sigma}(\symmgeneratorNRS ,R_{\Sigma}Y)}
		$$
	is a weak equivalence of simplicial sets.
	But since $C_{eff}^{q+1,\Sigma}\subseteq C_{eff}^{q,\Sigma}$, and by hypothesis
	$f$ is a $C_{eff}^{q,\Sigma}$-colocal equivalence; we have that all
	the induced maps $(R_{\Sigma}f)_{\ast}$ are weak
	equivalences of simplicial sets.  Thus $R_{\Sigma}f$ is a
	$C_{eff}^{q+1,\Sigma}$-colocal equivalence, as we wanted.
	
	Finally (\ref{prop.3.2.functors-between-Rqsymm.b}) and (\ref{prop.3.2.functors-between-Rqsymm.c})
	follow directly from
	proposition \ref{prop.3.2.Rq-lifts-qSH}, corollary \ref{cor.3.3.symmetrization-functor-exact-equivalenceRq},
	theorems \ref{thm.3.2.Rq-models-f<q}, \ref{thm.3.3.symmRq-models-symmf<q} together with
	the commutative diagram 
	\eqref{diagram3.3.liftslicefiltrationsymm} of left Quillen funtors
	constructed above and \cite[theorem 1.3.7]{MR1650134}.
\end{proof}

\begin{thm}
		\label{thm.3.3.liftingslicefilt-modelcatssymm}
	We have the following commutative diagram
	of left Quillen functors:
		\begin{equation}
					\label{diagram.3.3.liftingVoevodskyslicefiltsymm.a}
			\begin{array}{c}
				\xymatrix{\vdots \ar[d]_-{id} & \\ 
									\qplusoneconnectedsymmTspectra \ar[d]_-{id} \ar[dr]^-{id} & \\ 
									\qconnectedsymmTspectra \ar[d]_-{id} \ar[r]^-{id} & \motivicsymmTspectra \\
									\qminusoneconnectedsymmTspectra \ar[d]_-{id} \ar[ur]_-{id} & \\ 
									\vdots & }
			\end{array}
		\end{equation}
	and the associated diagram of homotopy categories:	
		\begin{equation}
					\label{diagram.3.3.liftingVoevodskyslicefiltsymm.b}
			\begin{array}{c}
				\xymatrix{\vdots \ar@<-1ex>[d] &&& \\ 
									R_{C_{eff}^{q+1}}\symmstablehomotopy \ar@<-1ex>[d]_-{C_{q+1}^{\Sigma}} 
									\ar@<2ex>[drrr]|-{C_{q+1}^{\Sigma}} \ar@<-1ex>[u] &&& \\ 
									\qconnectedsymmstablehomotopy \ar@<-1ex>[d]_-{C_{q}^{\Sigma}} 
									\ar@<1ex>[rrr]|-{C_{q}^{\Sigma}} \ar@<-1ex>[u]_-{R_{\Sigma}} &&& \symmstablehomotopy 
									\ar@<1ex>[lll]|-{R_{\Sigma}} \ar[ulll]|-{R_{\Sigma}} \ar@<2ex>[dlll]|-{R_{\Sigma}}\\
									R_{C_{eff}^{q-1}}\symmstablehomotopy \ar@<-1ex>[d] \ar[urrr]|-{C_{q-1}^{\Sigma}} \ar@<-1ex>[u]_-{R_{\Sigma}}&&& \\ 
									\vdots \ar@<-1ex>[u]&& }
			\end{array}
		\end{equation}
	gets canonically identified, through
	the equivalences of categories $r_{q}C_{q}$, $IQ_{T}Ji_{q}$, 
	$VC_{q}$ and $UR_{\Sigma}$ constructed
	in proposition \ref{prop.3.2.Rq-lifts-qSH} and corollary \ref{cor.3.3.symmetrization-functor-exact-equivalenceRq}; 
	with Voevodsky's slice filtration:
		\begin{equation}
					\label{diagram.3.3.liftingVoevodskyslicefiltsymm.c}
			\begin{array}{c}
				\xymatrix{\vdots \ar@<-1ex>[d] &&& \\ 
									\Tsuspfunctor ^{q+1}\stablehomotopy \ar@<-1ex>[d]_-{i_{q+1}} \ar@<2ex>[drrr]|-{i_{q+1}} \ar@<-1ex>[u] &&& \\ 
									\stablehomotopyeffq \ar@<-1ex>[d]_-{i_{q}} 
									\ar@<1ex>[rrr]|-{i_{q}} \ar@<-1ex>[u]_-{r_{q}} &&& \stablehomotopy 
									\ar@<1ex>[lll]|-{r_{q}} \ar[ulll]|-{r_{q+1}} \ar@<2ex>[dlll]|-{r_{q-1}}\\
									\Tsuspfunctor ^{q-1}\stablehomotopy \ar@<-1ex>[d] \ar[urrr]|-{i_{q-1}} \ar@<-1ex>[u]_-{r_{q-1}}&&& \\ 
									\vdots \ar@<-1ex>[u]&& }
			\end{array}
		\end{equation}
\end{thm}
\begin{proof}
	Follows immediately from proposition
	\ref{prop.3.2.functors-between-Rqsymm}, corollary \ref{cor.3.3.symmetrization-functor-exact-equivalenceRq}
	and theorem \ref{thm.3.2.liftingslicefilt-modelcats}.
\end{proof}

\begin{thm}
		\label{thm.3.3.Lqsymmetricmodelstructures}
	Fix $q\in \mathbb Z$.  Consider the following
	set of maps in $\motivicsymmTspectra$ (see theorem \ref{thm.3.2.Lqmodelstructures}):
		\begin{eqnarray}
				\label{diagram.3.3.symmmapsleftloc}
			\qleftsymmlocmaps = \{ \symmkillNRS \: | \\
										 \symmgeneratorNRS \in C_{eff}^{q,\Sigma}\} \nonumber
		\end{eqnarray}
	Then the left Bousfield localization of $\motivicsymmTspectra$ with respect to
	the $\qleftsymmlocmaps$-local equivalences exists.  This new model structure
	will be called \emph{weight$^{<q}$ motivic symmetric stable}.
	$\weightqsymmTspectra$ will denote the category of symmetric $T$-spectra equipped with
	the weight$^{<q}$ motivic symmetric stable model
	structure, and $\weightqsymmstablehomotopy$ will denote
	its associated homotopy category.  Furthermore the weight$^{<q}$ motivic symmetric stable model structure
	is cellular, left proper and simplicial; with the following
	sets of generating cofibrations and trivial cofibrations
	respectively:
	$$\begin{array}{l}
	 I_{\qleftsymmlocmaps}=I^{T}_{\Sigma}  =\bigcup _{n\geq 0}\{F_{n}^{\Sigma}(Y_{+}\hookrightarrow (\Delta ^{n}_{U})_{+})\} \\
	 \\
	J_{\qleftsymmlocmaps}  =\{j:A\rightarrow B\}
	\end{array}
	$$
	where $j$ satisfies the following conditions:
	\begin{enumerate}
		\item	$j$ is an inclusion of $I^{T}_{\Sigma}$-complexes.
		\item	$j$ is a $\qleftsymmlocmaps$-local equivalence.
		\item	the size of $B$ as an $I^{T}_{\Sigma}$-complex is less than $\kappa$, 
					where $\kappa$ is the regular cardinal defined by Hirschhorn in \cite[definition 4.5.3]{MR1944041}.
	\end{enumerate}
\end{thm}
\begin{proof}
	Theorems \ref{thm.2.7.cellularity-motivicsymm-stable-str} and 
	\ref{thm2.6.stablesymmetricmodelstructure} imply that $\motivicsymmTspectra$ is a cellular, proper and simplicial
	model category.  Therefore the existence of the left Bousfield localization
	follows from \cite[theorem 4.1.1]{MR1944041}.  Using \cite[theorem 4.1.1]{MR1944041} again, we have that
	$\weightqsymmTspectra$ is cellular, left proper and simplicial; where the sets
	of generating cofibrations and trivial cofibrations are the ones described above.
\end{proof}

\begin{defi}
		\label{def.3.3.stable-weightqsymm-replacementfunctors}
	Fix $q\in \mathbb Z$.  Let $W_{q}^{\Sigma}$ denote a fibrant replacement functor in
	$\weightqsymmTspectra$, such that the for every symmetric $T$-spectrum $X$,
	the natural map:
		$$\xymatrix{X \ar[r]^-{W_{q}^{\Sigma ,X}}& W_{q}^{\Sigma}X}
		$$
	is a trivial cofibration in $\weightqsymmTspectra$, and $W_{q}^{\Sigma}X$
	is $\qleftsymmlocmaps$-local in $\motivicsymmTspectra$.
\end{defi}

\begin{prop}
		\label{prop.3.3.Q-cofibrant-replacement-all-symmetricL<q}
	Fix $q\in \mathbb Z$.  Then $Q_{\Sigma}$ is also a cofibrant
	replacement functor in $\weightqsymmTspectra$, and for every symmetric
	$T$-spectrum $X$ the natural map
		$$\xymatrix{Q_{\Sigma}X\ar[r]^{Q_{\Sigma}^{X}}& X}
		$$
	is a trivial fibration in $\weightqsymmTspectra$.
\end{prop}
\begin{proof}
	Since $\weightqsymmTspectra$ is the left Bousfield localization of
	$\motivicsymmTspectra$ with respect to the $\qleftsymmlocmaps$-local equivalences, by construction
	we have that the cofibrations and the trivial fibrations are indentical in
	$\weightqsymmTspectra$ and $\motivicsymmTspectra$ respectively.  This implies that for every symmetric
	$T$-spectrum $X$, $Q_{\Sigma}X$ is cofibrant in $\weightqsymmTspectra$, and we also have that
	the natural map
		$$\xymatrix{Q_{\Sigma}X\ar[r]^{Q_{\Sigma}^{X}}& X}
		$$
	is a trivial fibration in $\weightqsymmTspectra$.  Hence $Q_{\Sigma}$ is also a cofibrant replacement
	functor for $\weightqsymmTspectra$.
\end{proof}

\begin{prop}
		\label{prop.3.3.Z-symmLq-local-iff-UZ-Lq-local}
	Fix $q\in \mathbb Z$.
	Then a symmetric $T$-spectrum $Z$ is $\qleftsymmlocmaps$-local in $\motivicsymmTspectra$
	if and only if $UZ$ is $\qleftlocmaps$-local in $\motivicTspectra$.
\end{prop}
\begin{proof}
	We have that $Z$ is $\qleftsymmlocmaps$-local if and only if
	$Z$ is fibrant in $\motivicsymmTspectra$ and for every 
		$$\killNRS \in \qleftlocmaps$$
	the induced map
		$$\xymatrix{Map\: _{\Sigma}(V(\diskNRRS),Z) \ar[d]^-{V(\killNRSmap)^{\ast}}\\
								Map\: _{\Sigma}(V(\generatorNRS),Z)}
		$$
	is a weak equivalence of simplicial sets.
	
	On the other hand, we have that $UZ$ is $\qleftlocmaps$-local in $\motivicTspectra$ if and only
	if $UZ$ is fibrant in $\motivicTspectra$ and for every $\killNRS \in \qleftlocmaps$,
	the induced map
		$$\xymatrix{Map(\diskNRRS,UZ) \ar[rr]^-{(\killNRSmap)^{\ast}}&& Map(\generatorNRS,UZ)}
		$$
	is a weak equivalence of simplicial sets.
	
	Then the result follows from the following facts:
		\begin{enumerate}
			\item	By definition, $Z$ is fibrant in $\motivicsymmTspectra$ if $UZ$ is
						fibrant in $\motivicTspectra$.
			\item	Proposition \ref{prop.2.6.enriched-symmetrization-adjunction}, which implies that the adjunction
							$$\xymatrix{(V,U,\varphi):\motivicTspectra \ar[r]& \motivicsymmTspectra}
							$$
						is enriched in the category of simplicial sets.
		\end{enumerate}
\end{proof}

\begin{prop}
		\label{prop.3.3.Lqsymm-local-objects-classification}
	Fix $q\in \mathbb Z$, and let $Z$ be a symmetric $T$-spectrum.
	Then $Z$ is $\qleftsymmlocmaps$-local in $\motivicsymmTspectra$ if and only if
	the following conditions hold:
		\begin{enumerate}
			\item \label{prop.3.3.Lqsymm-local-objects-classification.a}$Z$ is fibrant 
						in $\motivicsymmTspectra$.
			\item \label{prop.3.3.Lqsymm-local-objects-classification.b}For every $\symmgeneratorNRS 
						\in C_{eff}^{q,\Sigma}$,
						$[\symmgeneratorNRS ,Z]_{Spt}^{\Sigma}\cong 0$
		\end{enumerate}
\end{prop}
\begin{proof}
	Follows directly from propositions \ref{prop.3.3.Z-symmLq-local-iff-UZ-Lq-local} 
	and \ref{prop.3.2.Lq-local-objects-classification},
	together with the fact that $(V,U,\varphi):\motivicTspectra \rightarrow \motivicsymmTspectra$
	is a Quillen adjunction.
\end{proof}

\begin{cor}
		\label{cor.3.3.S1loops-preserves-symmLqlocal}
	Fix $q\in \mathbb Z$, and let $Z$ be a fibrant symmetric $T$-spectrum in
	$\motivicsymmTspectra$.
	Then $Z$ is $\qleftsymmlocmaps$-local in $\motivicsymmTspectra$ if and only if
	$\Omega _{S^{1}}Z$ is $\qleftsymmlocmaps$-local in $\motivicsymmTspectra$.
\end{cor}
\begin{proof}
	By proposition \ref{prop.3.3.Z-symmLq-local-iff-UZ-Lq-local} we have that
	$Z$ is $\qleftsymmlocmaps$-local if and only if $UZ$ is $\qleftlocmaps$-local
	in $\motivicTspectra$.  Now corollary \ref{cor.3.2.S1loops-preserves-Lqlocal} implies
	that $UZ$ is $\qleftlocmaps$-local if and only if $\Omega _{S^{1}}UZ=U(\Omega _{S^{1}}Z)$ is
	$\qleftlocmaps$-local.
	
	Therefore using proposition \ref{prop.3.3.Z-symmLq-local-iff-UZ-Lq-local} again, we get
	that $Z$ is $\qleftsymmlocmaps$-local if and only if $\Omega _{S^{1}}Z$ is
	$\qleftsymmlocmaps$-local.
\end{proof}

\begin{cor}
		\label{cor.3.3.S1-preserves-symmLqlocal}
	Fix $q\in \mathbb Z$, and let $Z$ be a fibrant
	symmetric $T$-spectrum in $\motivicsymmTspectra$.  Then $Z$ is $\qleftsymmlocmaps$-local 
	in $\motivicsymmTspectra$ if and only if
	$R_{\Sigma}(Q_{\Sigma}Z \wedge S^{1})$ is $\qleftsymmlocmaps$-local in $\motivicsymmTspectra$.
\end{cor}
\begin{proof}
	($\Rightarrow$):  Assume that $Z$ is $\qleftsymmlocmaps$-local.
	Since $R_{\Sigma}(Q_{\Sigma}Z \wedge S^{1})$ is fibrant, using proposition
	\ref{prop.3.3.Lqsymm-local-objects-classification} we have that it is enough
	to check that for every $\symmgeneratorNRS \in C_{eff}^{q,\Sigma}$,
	$[\symmgeneratorNRS ,R_{\Sigma}(Q_{\Sigma}Z \wedge S^{1})]_{Spt}^{\Sigma}\cong 0$.
	But since $- \wedge S^{1}$ is a Quillen equivalence, we get the
	following diagram:
		$$\xymatrix@C=-2pc{[\symmgeneratorNRS ,R_{\Sigma}(Q_{\Sigma}Z \wedge S^{1})]_{Spt}^{\Sigma} \ar[dr]^-{\cong}&\\ 
								& [F_{n+1}^{\Sigma}(S^{r+1}\wedge \gm ^{s+1}\wedge U_{+}),R_{\Sigma}(Q_{\Sigma}Z \wedge S^{1})]_{Spt}^{\Sigma}\\
								[F_{n+1}^{\Sigma}(S^{r}\wedge \gm ^{s+1}\wedge U_{+}), Z]_{Spt}^{\Sigma} \ar[dr]^-{\Tsuspfunctor ^{1,0}}_-{\cong}&\\ 
								& [F_{n+1}^{\Sigma}(S^{r+1}\wedge \gm ^{s+1}\wedge U_{+}), Q_{\Sigma}Z \wedge S^{1}]_{Spt}^{\Sigma} \ar[uu]_-{\cong}}
		$$
	where all the maps are isomorphisms of abelian groups.  Since $Z$
	is $\qleftsymmlocmaps$-local, proposition \ref{prop.3.3.Lqsymm-local-objects-classification}
	implies that $[F_{n+1}^{\Sigma}(S^{r}\wedge \gm ^{s+1}\wedge U_{+}), Z]_{Spt}^{\Sigma}\cong 0$.
	Therefore 
		$$[\symmgeneratorNRS ,R_{\Sigma}(Q_{\Sigma}Z \wedge S^{1})]_{Spt}^{\Sigma}\cong 0$$ 
	for every
	$\symmgeneratorNRS \in C_{eff}^{q,\Sigma}$, as we wanted.
	
	($\Leftarrow$):  Assume that $R_{\Sigma}(Q_{\Sigma}Z \wedge S^{1})$ is $\qleftsymmlocmaps$-local.
	By hypothesis, $Z$ is fibrant; therefore proposition \ref{prop.3.3.Lqsymm-local-objects-classification}
	implies that it is enough to show that for every $\symmgeneratorNRS \in C_{eff}^{q,\Sigma}$,
	$[\symmgeneratorNRS ,Z]_{Spt}^{\Sigma}\cong 0$.	
	Since $\motivicsymmTspectra$ is a simplicial model category
	and $- \wedge S^{1}$ is a Quillen equivalence;
	we have the following diagram:
		$$\xymatrix{[\symmgeneratorNRS , \Omega _{S^{1}}R_{\Sigma}(Q_{\Sigma}Z \wedge S^{1})]_{Spt}^{\Sigma} \ar[d]^-{\cong}&\\
								[\symmgeneratorNRS \wedge S^{1}, Q_{\Sigma}Z \wedge S^{1}]_{Spt}^{\Sigma}
								& [\symmgeneratorNRS ,Z]_{Spt}^{\Sigma}\ar[l]_-{\Tsuspfunctor ^{1,0}}^-{\cong}}
		$$
	where all the maps are isomorphisms of abelian groups.
	On the other hand, using corollary \ref{cor.3.3.S1loops-preserves-symmLqlocal}
	we have that $\Omega _{S^{1}}R_{\Sigma}(Q_{\Sigma}Z \wedge S^{1})$ is $\qleftsymmlocmaps$-local.
	Therefore using proposition \ref{prop.3.3.Lqsymm-local-objects-classification} again,
	we have that for every $\symmgeneratorNRS \in C_{eff}^{q,\Sigma}$:
		$$[\symmgeneratorNRS ,Z]_{Spt}^{\Sigma}\cong [\symmgeneratorNRS , 
			\Omega _{S^{1}}R_{\Sigma}(Q_{\Sigma}Z \wedge S^{1})]_{Spt}^{\Sigma}\cong 0
		$$
	and this finishes the proof.
\end{proof}

\begin{cor}
		\label{cor.3.3.detecting-Cq-symmlocal-equivalences}
	Fix $q\in \mathbb Z$, and let $f:X\rightarrow Y$ be a map of symmetric $T$-spectra.
	Then $f$ is a $\qleftsymmlocmaps$-local equivalence in $\motivicTspectra$ if and only if
	for every $\qleftsymmlocmaps$-local symmetric $T$-spectrum $Z$, $f$ induces
	the following isomorphism of abelian groups:
		$$\xymatrix{[Y,Z]_{Spt}^{\Sigma} \ar[r]^-{f^{\ast}}& [X,Z]_{Spt}^{\Sigma}}
		$$
\end{cor}
\begin{proof}
	Suppose that $f$ is a $\qleftsymmlocmaps$-local equivalence, then by
	definition the induced map:
		$$\xymatrix{Map\: _{\Sigma}(Q_{\Sigma}Y,Z) \ar[r]^-{(Q_{\Sigma}f)^{\ast}}& Map\: _{\Sigma}(Q_{\Sigma}X,Z)}
		$$
	is a weak equivalence of simplicial sets for every $\qleftsymmlocmaps$-local symmetric
	$T$-spectrum $Z$.  Proposition \ref{prop.3.3.Lqsymm-local-objects-classification}(\ref{prop.3.3.Lqsymm-local-objects-classification.a}) 
	implies that $Z$ is fibrant in $\motivicsymmTspectra$, and
	since $\motivicsymmTspectra$ is in particular a simplicial model category; 
	we get the following commutative diagram, where the top row and all the
	vertical maps are isomorphisms of abelian groups:
		$$\xymatrix{\pi_{0}Map\: _{\Sigma}(Q_{\Sigma}Y,Z) \ar[r]^-{(Q_{\Sigma}f)^{\ast}}_-{\cong} \ar[d]_-{\cong}& 
								\pi_{0}Map\: _{\Sigma}(Q_{\Sigma}X,Z) \ar[d]^-{\cong}\\
								[Y,Z]_{Spt}^{\Sigma} \ar[r]_-{f^{\ast}}& [X,Z]_{Spt}^{\Sigma}}
		$$
	hence $f^{\ast}$ is an isomorphism for every $\qleftsymmlocmaps$-local symmetric $T$-spectrum $Z$,
	as we wanted.
	
	Conversely, assume that for every $\qleftsymmlocmaps$-local symmetric $T$-spectrum $Z$,
	the induced map
		$$\xymatrix{[Y,Z]_{Spt}^{\Sigma} \ar[r]^-{f^{\ast}}& [X,Z]_{Spt}^{\Sigma}}
		$$
	is an isomorphism of abelian groups.
	
	Since $\weightqsymmTspectra$ is the left Bousfield localization
	of $\motivicsymmTspectra$ with
	respect to the $\qleftsymmlocmaps$-local equivalences, we have that the 
	identity functor $id:\motivicsymmTspectra \rightarrow \weightqsymmTspectra$ is
	a left Quillen functor.  Therefore for every symmetric $T$-spectrum $Z$,
	we get the following commutative diagram
	where all the vertical arrows are isomorphisms:
		$$\xymatrix{\Hom _{\weightqsymmstablehomotopy}(Q_{\Sigma}Y,Z) \ar[r]^-{(Q_{\Sigma}f)^{\ast}} \ar[d]_-{\cong} &
								\Hom _{\weightqsymmstablehomotopy}(Q_{\Sigma}X,Z) \ar[d]^-{\cong}\\
								[Y,W_{q}^{\Sigma}Z]_{Spt}^{\Sigma} \ar[r]_-{f^{\ast}}^-{\cong}& [X,W_{q}^{\Sigma}Z]_{Spt}^{\Sigma}}
		$$
	but $W_{q}^{\Sigma}Z$ is by construction $\qleftsymmlocmaps$-local, then
	by hypothesis the bottom row is an isomorphism of abelian groups.
	Hence it follows that the induced map:
		$$\xymatrix{\Hom _{\weightqsymmstablehomotopy}(Q_{\Sigma}Y,Z) \ar[r]^-{(Q_{\Sigma}f)^{\ast}}_{\cong}  &
								\Hom _{\weightqsymmstablehomotopy}(Q_{\Sigma}X,Z)}
		$$
	is an isomorphism for every symmetric $T$-spectrum $Z$.
	This implies that $Q_{\Sigma}f$ is a weak equivalence in $\weightqsymmTspectra$, and since
	$Q_{\Sigma}$ is also a cofibrant replacement functor in $\weightqsymmTspectra$, it follows
	that $f$ is a weak equivalence in $\weightqsymmTspectra$.
	Therefore we have that $f$ is a $\qleftsymmlocmaps$-local equivalence, as we wanted.
\end{proof}

\begin{lem}
		\label{lem.3.3.symmLq-localequivs-stable-S1suspension}
	Fix $q\in \mathbb Z$, and let $f:X\rightarrow Y$ be a map of symmetric $T$-spectra.
	Then $f$ is a $\qleftsymmlocmaps$-local equivalence in $\motivicsymmTspectra$
	if and only if  
		$$Q_{\Sigma}f \wedge id: Q_{\Sigma}X \wedge S^{1} \rightarrow Q_{\Sigma}Y \wedge S^{1}$$ 
	is a
	$\qleftsymmlocmaps$-local equivalence in $\motivicsymmTspectra$.
\end{lem}
\begin{proof}
	Assume that $f$ is a $\qleftsymmlocmaps$-local equivalence, and let $Z$
	be an arbitrary $\qleftsymmlocmaps$-local symmetric $T$-spectrum.  Then 
	corollary \ref{cor.3.3.S1loops-preserves-symmLqlocal}
	implies that $\Omega _{S^{1}}Z$ is also $\qleftsymmlocmaps$-local.
	Therefore the induced map
		$$\xymatrix{Map\: _{\Sigma}(Q_{\Sigma}Y,\Omega _{S^{1}}Z) \ar[r]^-{(Q_{\Sigma}f)^{\ast}}& 
								Map\: _{\Sigma}(Q_{\Sigma}X,\Omega _{S^{1}}Z)}
		$$
	is a weak equivalence of simplicial sets.  Now
	since $\motivicsymmTspectra$
	is a simplicial model category, we have the following commutative diagram:
		$$\xymatrix{Map\: _{\Sigma}(Q_{\Sigma}Y,\Omega _{S^{1}}Z) \ar[rr]^-{(Q_{\Sigma}f)^{\ast}}\ar[d]_-{\cong}&& 
								Map\: _{\Sigma}(Q_{\Sigma}X,\Omega _{S^{1}}Z) \ar[d]^-{\cong}\\
								Map\: _{\Sigma}(Q_{\Sigma}Y \wedge S^{1},Z) \ar[rr]^-{(Q_{\Sigma}f \wedge id)^{\ast}}&& 
								Map\: _{\Sigma}(Q_{\Sigma}X \wedge S^{1},Z)}
		$$
	and using the two out of three property for weak equivalences of simplicial sets,
	we have that
		$$\xymatrix{Map\: _{\Sigma}(Q_{\Sigma}Y \wedge S^{1},Z) \ar[rr]^-{(Q_{\Sigma}f \wedge id)^{\ast}}&& 
								Map\: _{\Sigma}(Q_{\Sigma}X \wedge S^{1},Z)}
		$$
	is a weak equivalence.  Since this holds for every $\qleftsymmlocmaps$-local symmetric $T$-spectrum
	$Z$, it follows that
		$$Q_{\Sigma}f \wedge id: Q_{\Sigma}X \wedge S^{1} \rightarrow Q_{\Sigma}Y \wedge S^{1}
		$$
	is a $\qleftsymmlocmaps$-local equivalence, as we wanted.
	
	Conversely, suppose that 
		$$Q_{\Sigma}f \wedge id: Q_{\Sigma}X \wedge S^{1} \rightarrow Q_{\Sigma}Y \wedge S^{1}
		$$
	is a $\qleftsymmlocmaps$-local equivalence.  Let $Z$ be an arbitrary
	$\qleftsymmlocmaps$-local symmetric $T$-spectrum. Since $\motivicsymmTspectra$
	is a simplicial model category and
	$-\wedge S^{1}$ is a Quillen equivalence, we get the following
	commutative diagram:
		$$\xymatrix{[Q_{\Sigma}Y \wedge S^{1}, R_{\Sigma}(Q_{\Sigma}Z \wedge S^{1})]_{Spt}^{\Sigma} 
								\ar[rr]^-{(Q_{\Sigma}f \wedge id)^{\ast}} \ar[d]_-{\cong}&& 
								[Q_{\Sigma}X \wedge S^{1}, R_{\Sigma}(Q_{\Sigma}Z \wedge S^{1})]_{Spt}^{\Sigma}\ar[d]^-{\cong}\\
								[Q_{\Sigma}Y \wedge S^{1}, Q_{\Sigma}Z \wedge S^{1}]_{Spt}^{\Sigma} \ar[rr]^-{(Q_{\Sigma}f \wedge id)^{\ast}}&&
								[Q_{\Sigma}X \wedge S^{1}, Q_{\Sigma}Z \wedge S^{1}]_{Spt}^{\Sigma}\\
								[Y, Z]_{Spt}^{\Sigma} \ar[rr]_-{f^{\ast}}\ar[u]^-{\cong}_-{\Tsuspfunctor ^{1,0}}&&
								[X, Z]_{Spt}^{\Sigma}\ar[u]_-{\cong}^-{\Tsuspfunctor ^{1,0}} }
		$$
	Now, corollary \ref{cor.3.3.S1-preserves-symmLqlocal} implies that $R_{\Sigma}(Q_{\Sigma}Z \wedge S^{1})$ is
	also $\qleftsymmlocmaps$-local.  Therefore using corollary \ref{cor.3.3.detecting-Cq-symmlocal-equivalences}
	we have that the top row in the diagram above is an isomorphism of abelian groups.
	This implies that the induced map:
		$$\xymatrix{[Y, Z]_{Spt}^{\Sigma} \ar[r]^-{f^{\ast}}&
								[X, Z]_{Spt}^{\Sigma}}
		$$
	is an isomorphism of abelian groups for every $\qleftsymmlocmaps$-local symmetric spectrum $Z$.
	Finally using corollary \ref{cor.3.3.detecting-Cq-symmlocal-equivalences} again,
	we have that $f:X\rightarrow Y$ is a $\qleftsymmlocmaps$-local equivalence, as we
	wanted.
\end{proof}

\begin{cor}
		\label{cor.3.3.susp-Quillen-equiv-onsymmLq}
	For every $q\in \mathbb Z$, the following adjunction:
		$$\xymatrix{(-\wedge S^{1},\Omega _{S^{1}},\varphi):\weightqsymmTspectra \ar[rr]
								&& \weightqsymmTspectra}
		$$
	is a Quillen equivalence.
\end{cor}
\begin{proof}
	Using corollary 1.3.16 in \cite{MR1650134} and
	proposition \ref{prop.3.3.Q-cofibrant-replacement-all-symmetricL<q} we have that it suffices to
	verify the following two conditions:
		\begin{enumerate}
			\item	\label{cor.3.3.susp-Quillen-equiv-onsymmLq.a}For every fibrant object $X$
						in $\weightqsymmTspectra$, the following composition
							$$\xymatrix{(Q_{\Sigma}\Omega _{S^{1}}X)\wedge S^{1} \ar[rr]^-{Q_{\Sigma}^{\Omega _{S^{1}}X}\wedge id}&& 
													(\Omega _{S^{1}}X)\wedge S^{1} \ar[r]^-{\epsilon _{X}}& X}
							$$
						is a $\qleftsymmlocmaps$-local equivalence.
			\item	\label{cor.3.3.susp-Quillen-equiv-onsymmLq.b}$-\wedge S^{1}$ reflects $\qleftsymmlocmaps$-local equivalences
						between cofibrant objects in $\weightqsymmTspectra$.
		\end{enumerate}
		
	(\ref{cor.3.3.susp-Quillen-equiv-onsymmLq.a}):  By construction $\weightqsymmTspectra$ is
	a left Bousfield localization of $\motivicsymmTspectra$, therefore the identity
	functor 
		$$\xymatrix{id:\weightqsymmTspectra \ar[r]& \motivicsymmTspectra}
		$$ 
	is a right Quillen functor.
	Thus $X$ is also fibrant in $\motivicsymmTspectra$.
	Since the adjunction $(-\wedge S^{1},\Omega _{S^{1}},\varphi)$
	is a Quillen equivalence on $\motivicsymmTspectra$, 
	\cite[proposition 1.3.13(b)]{MR1650134} implies that the following composition
	is a weak equivalence in $\motivicsymmTspectra$:
		$$\xymatrix{(Q_{\Sigma}\Omega _{S^{1}}X)\wedge S^{1} \ar[rr]^-{Q_{\Sigma}^{\Omega _{S^{1}}X}\wedge id}&&
													(\Omega _{S^{1}}X) \wedge S^{1} \ar[r]^-{\epsilon _{X}}& X}
		$$
	Hence using \cite[proposition 3.1.5]{MR1944041} it follows that the composition above
	is a $\qleftsymmlocmaps$-local equivalence.
	
	(\ref{cor.3.3.susp-Quillen-equiv-onsymmLq.b}):  This follows immediately
	from proposition \ref{prop.3.3.Q-cofibrant-replacement-all-symmetricL<q} and 
	lemma \ref{lem.3.3.symmLq-localequivs-stable-S1suspension}.
\end{proof}

\begin{rmk}
		\label{rmk.3.3.Tsuspfunctor-doesnotdescend-to-symmLq}
	We have a situation similar to the one
	described in remark
	\ref{rmk.3.3.smashT-not-symmdescending} for the model categories
	$\qconnectedsymmTspectra$; i.e. although
	the adjunction $(\Tsuspfunctor ,\Tloops ,\varphi)$ is
	a Quillen equivalence on $\motivicsymmTspectra$,
	it does not descend even to a Quillen
	adjunction on the weight$^{<q}$ motivic symmetric stable
	model category
	$\weightqsymmTspectra$.
\end{rmk}

\begin{cor}
		\label{cor.3.3.stablehomotopysymmLq==triangulated}
	For every $q\in \mathbb Z$,
	the homotopy category $\weightqsymmstablehomotopy$
	associated to $\weightqsymmTspectra$ has the
	structure of a triangulated category.
\end{cor}
\begin{proof}
	Theorem \ref{thm.3.3.Lqsymmetricmodelstructures} implies in particular
	that $\weightqsymmTspectra$ is a pointed simplicial model category,
	and  corollary \ref{cor.3.3.susp-Quillen-equiv-onsymmLq} implies that
	the adjunction 
		$$(-\wedge S^{1},\Omega _{S^{1}},\varphi):\weightqsymmTspectra \rightarrow \weightqsymmTspectra$$
	is a Quillen equivalence.  Therefore
	the result follows from the work of Quillen in 
	\cite[sections I.2 and I.3]{MR0223432} and the work of 
	Hovey in \cite[chapters VI and VII]{MR1650134}.
\end{proof}

\begin{cor}
		\label{cor.3.3.symmLq=>rightproper}
	For every $q\in \mathbb Z$, $\weightqsymmTspectra$ is
	a right proper model category.
\end{cor}
\begin{proof}
	We need to show that the $\qleftsymmlocmaps$-local equivalences are stable
	under pullback along fibrations in $\weightqsymmTspectra$.
	Consider the following pullback diagram:
		$$\xymatrix{Z \ar[r]^-{w^{\ast}} \ar[d]_-{p^{\ast}}& X \ar[d]^-{p}\\
								W \ar[r]_-{w}& Y}
		$$
	where $p$ is a fibration in $\weightqsymmTspectra$, and $w$ is a
	$\qleftsymmlocmaps$-local equivalence.  Let $F$ be the homotopy fibre
	of $p$.  Then we get the following commutative diagram in $\weightqsymmstablehomotopy$:
		$$\xymatrix{\Omega _{S^{1}}Y \ar[r]^-{q} & F \ar[r]^-{i}& X \ar[r]^-{p}& Y\\
								\Omega _{S^{1}}W \ar[r]_-{r} \ar[u]^-{\Omega _{S^{1}}w}& F \ar[r]_-{j} \ar@{=}[u]& 
								Z \ar[r]_-{p^{\ast}} \ar[u]_-{w^{\ast}}& W \ar[u]_-{w}}
		$$
	Since the rows in the diagram above are both fibre sequences in $\weightqsymmTspectra$,
	it follows that both rows are distinguished triangles in $\weightqsymmstablehomotopy$
	(which has the structure of a triangulated category given by 
	corollary \ref{cor.3.3.stablehomotopysymmLq==triangulated}).  Now
	$w, id_{F}$ are both isomorphisms in $\weightqsymmstablehomotopy$, hence it follows that
	$w^{\ast}$ is also an isomorphism in $\weightqsymmstablehomotopy$.  Therefore
	$w^{\ast}$ is a $\qleftsymmlocmaps$-local equivalence, as we wanted.
\end{proof}

\begin{prop}
		\label{prop.3.3.symmLq-exact-adjunctions}
	For every $q\in \mathbb Z$ we have the following
	adjunction
		$$\xymatrix{(Q_{\Sigma},W_{q}^{\Sigma},\varphi):\symmstablehomotopy \ar[r]& \weightqsymmstablehomotopy}
		$$
	of exact functors between triangulated categories.
\end{prop}
\begin{proof}
	Since $\weightqsymmTspectra$ is the left Bousfield localization of
	$\motivicsymmTspectra$ with respect to the $\qleftsymmlocmaps$-local equivalences, we have that the
	identity functor $id:\motivicsymmTspectra \rightarrow \weightqsymmTspectra$
	is a left Quillen functor.  Therefore we get
	the following adjunction at the level of the associated homotopy categories:
		$$\xymatrix{(Q_{\Sigma}, W_{q}^{\Sigma}, \varphi ):\symmstablehomotopy \ar[r]& \weightqsymmstablehomotopy}
		$$
	
	Now proposition 6.4.1 in \cite{MR1650134} implies that
	$Q_{\Sigma}$ maps cofibre sequences in $\symmstablehomotopy$ to cofibre sequences in
	$\weightqsymmstablehomotopy$.
	Therefore 
	using proposition 7.1.12 in \cite{MR1650134} we have
	that $Q_{\Sigma}$ and $W_{q}^{\Sigma}$ are both exact functors between triangulated categories.
\end{proof}

\begin{lem}
		\label{lem.2.3.L<q-local==stable-under-symmetrization}
	Fix $q\in \mathbb Z$, and let $X$ be
	a $\qleftlocmaps$-local spectrum in $\motivicTspectra$.  Then $Q_{s}X$ and
	$UR_{\Sigma}VQ_{s}X$ are also $\qleftlocmaps$-local
	in $\motivicTspectra$.
\end{lem}
\begin{proof}
	Since $X$ is $\qleftlocmaps$-local, it follows that
	$X$ is fibrant in $\motivicTspectra$.
	By definition we have that the natural map
		$$\xymatrix{Q_{s}X \ar[r]^-{Q_{s}^{X}}& X}
		$$
	is a trivial fibration in $\motivicTspectra$, therefore
	$Q_{s}X$ is also fibrant in $\motivicTspectra$.  Hence
	\cite[lemma 3.2.1(a)]{MR1944041} implies that $Q_{s}X$ is
	$\qleftlocmaps$-local.
	
	Since the adjunction $(V,U,\varphi)$ is a Quillen equivalence
	between $\motivicTspectra$ and $\motivicsymmTspectra$, 
	we have that $UR_{\Sigma}VQ_{s}X$ is fibrant
	in $\motivicTspectra$, and \cite[proposition 1.3.13(b)]{MR1650134} implies that the composition
		$$\xymatrix{Q_{s}X \ar[r]^-{\eta _{Q_{s}X}}& UV(Q_{s}X) \ar[rr]^-{U(R_{\Sigma}^{VQ_{s}X})}&& UR_{\Sigma}VQ_{s}X}
		$$
	is a weak equivalence in $\motivicTspectra$.  Since we already know that $Q_{s}X$ is
	$\qleftlocmaps$-local, using \cite[lemma 3.2.1(a)]{MR1944041} again we get
	that $UR_{\Sigma}VQ_{s}X$ is also $\qleftlocmaps$-local in $\motivicTspectra$.
	This finishes the proof.
\end{proof}

\begin{prop}
		\label{prop.3.3.V-reflects-Lq-localequivalences}
	Fix $q\in \mathbb Z$, and let
	$f:X\rightarrow Y$ be a map in $\motivicTspectra$.  Then $f$ is a $\qleftlocmaps$-local
	equivalence in $\motivicTspectra$ if and only if $VQ_{s}f$ is a $\qleftsymmlocmaps$-local
	equivalence in $\motivicsymmTspectra$.
\end{prop}
\begin{proof}
	($\Rightarrow$):  Assume that $f$ is a $\qleftlocmaps$-local equivalence, and let
	$Z$ be an arbitrary $\qleftsymmlocmaps$-local symmetric $T$-spectrum.  
	Then $Z$ is fibrant in $\motivicsymmTspectra$, and using theorem
	\ref{thm.2.6.T-spectra===symmTspectra} we get the following commutative diagram where all the vertical arrows are
	isomorphisms:
		$$\xymatrix{[VQ_{s}Y,Z]_{Spt}^{\Sigma} \ar[rr]^-{(VQ_{s}f)^{\ast}}\ar[d]_-{\cong}&& [VQ_{s}X,Z]_{Spt}^{\Sigma} \ar[d]^-{\cong}\\
								[Y, UZ]_{Spt} \ar[rr]_-{f^{\ast}}&& [X, UZ]_{Spt}}
		$$
	By proposition \ref{prop.3.3.Z-symmLq-local-iff-UZ-Lq-local} 
	we have that $UZ$ is $\qleftlocmaps$-local in $\motivicTspectra$, hence corollary
	\ref{cor.3.2.detecting-Cq-local-equivalences} implies that the bottom row in the diagram above is always an isomorphism.
	Therefore the top row in the diagram above is an isomorphism 
	for every $\qleftsymmlocmaps$-local symmetric $T$-spectrum $Z$, then by 
	corollary \ref{cor.3.3.detecting-Cq-symmlocal-equivalences}
	it follows that $VQ_{s}f$ is a $\qleftsymmlocmaps$-local equivalence in $\motivicsymmTspectra$.
	
	($\Leftarrow$):  Assume that $VQ_{s}f$ is a $\qleftsymmlocmaps$-local equivalence
	in $\motivicsymmTspectra$, and let $Z$ be an arbitrary $\qleftlocmaps$-local $T$-spectrum
	in $\motivicTspectra$. We need to show that the induced map:
		$$\xymatrix{Map(Q_{s}Y,Z) \ar[r]^-{(Q_{s}f)^{\ast}}& Map(Q_{s}X,Z)}
		$$
	is a weak equivalence of simplicial sets.
	
	But theorem \ref{thm.2.6.T-spectra===symmTspectra} implies that
	the adjunction $(V,U,\varphi)$ is a Quillen equivalence between $\motivicTspectra$ and
	$\motivicsymmTspectra$, therefore using \cite[proposition 1.3.13(b)]{MR1650134} we have that all the maps
	in the following diagram are weak equivalences in $\motivicTspectra$:
		$$\xymatrix{Z & Q_{s}Z \ar[l]_-{Q_{s}^{Z}} \ar[rrrr]^-{U(R_{\Sigma}^{VQ_{s}Z})\circ \eta _{Q_{s}Z}}
								&&&& UR_{\Sigma}VQ_{s}Z}
		$$
		
	Lemma \ref{lem.2.3.L<q-local==stable-under-symmetrization} implies in particular that
	$Z, Q_{s}Z, UR_{\Sigma}VQ_{s}Z$ are all fibrant in $\motivicTspectra$.  Now using the fact
	that $\motivicTspectra$ is a simplicial model category together with Ken Brown's lemma
	(see lemma \ref{lemma1.1.KenBrown2}) and the two out of three property for weak equivalences, we have that it
	suffices to prove that the induced map:
		$$\xymatrix{Map(Q_{s}Y,UR_{\Sigma}VQ_{s}Z) \ar[r]^-{(Q_{s}f)^{\ast}}& Map(Q_{s}X,UR_{\Sigma}VQ_{s}Z)}
		$$
	is a weak equivalence of simplicial sets.  Using the enriched adjunctions
	of proposition \ref{prop.2.6.enriched-symmetrization-adjunction}, 
	we get the following commutative diagram where all the vertical arrows
	are isomorphisms:
		$$\xymatrix{Map(Q_{s}Y,UR_{\Sigma}VQ_{s}Z) \ar[rr]^-{(Q_{s}f)^{\ast}} \ar[d]_-{\cong}&&
								Map(Q_{s}X,UR_{\Sigma}VQ_{s}Z) \ar[d]^-{\cong}\\
								Map\: _{\Sigma}(VQ_{s}Y, R_{\Sigma}VQ_{s}Z)\ar[rr]_-{(VQ_{s}f)^{\ast}}&&
								Map\: _{\Sigma}(VQ_{s}X, R_{\Sigma}VQ_{s}Z)}
		$$
	Finally, lemma \ref{lem.2.3.L<q-local==stable-under-symmetrization} implies that
	$UR_{\Sigma}VQ_{s}Z$ is $\qleftlocmaps$-local in $\motivicTspectra$, therefore by 
	proposition \ref{prop.3.3.Z-symmLq-local-iff-UZ-Lq-local}
	we have that $R_{\Sigma}VQ_{s}Z$ is $\qleftsymmlocmaps$-local in $\motivicsymmTspectra$.
	Since $VQ_{s}f$ is a $\qleftsymmlocmaps$-local equivalence and $VQ_{s}X, VQ_{s}Y$ are
	both cofibrant in $\motivicsymmTspectra$, it follows that the bottom row in the diagram above
	is a weak equivalence of simplicial sets.  This implies that the top row is also
	a weak equivalence of simplicial sets, as we wanted.
\end{proof}

\begin{thm}
		\label{thm.3.3.symmetrization-Quillen-equivalence-weight<q}
	For every $q\in \mathbb Z$, the adjunction
		$$\xymatrix{(V,U,\varphi):\weightqTspectra \ar[r]& \weightqsymmTspectra}
		$$
	given by the symmetrization and the forgetful functor is a Quillen equivalence.
\end{thm}
\begin{proof}
	Proposition \ref{prop.3.3.V-reflects-Lq-localequivalences} 
	together with the universal property for left Bousfield localizations
	(see definition \ref{def-1.leftlocmodcats})
	imply that 
		$$\xymatrix{V:\weightqTspectra \ar[r]& \weightqsymmTspectra}$$ 
	is a left
	Quillen functor.
	Using corollary 1.3.16 in \cite{MR1650134} and
	proposition \ref{prop.3.2.Qs-cofibrant-replacement-all-L<q} we have that it suffices to
	verify the following two conditions:
		\begin{enumerate}
			\item	\label{thm.3.3.symmetrization-Quillen-equivalence-weight<q.a}For every fibrant object $X$
						in $\weightqsymmTspectra$, the following composition
							$$\xymatrix{VQ_{s}U(X) \ar[rr]^-{V(Q_{s}^{UX})}&& 
													VU(X)\ar[r]^-{\epsilon _{X}}& X}
							$$
						is a weak equivalence in $\weightqsymmTspectra$.
			\item	\label{thm.3.3.symmetrization-Quillen-equivalence-weight<q.b}$V$ reflects weak equivalences
						between cofibrant objects in $\weightqTspectra$.
		\end{enumerate}
		
	(\ref{thm.3.3.symmetrization-Quillen-equivalence-weight<q.a}):
	By construction $\weightqsymmTspectra$ is
	a left Bousfield localization of $\motivicsymmTspectra$, therefore the identity
	functor 
		$$\xymatrix{id:\weightqsymmTspectra \ar[r]& \motivicsymmTspectra}
		$$ 
	is a right Quillen functor.
	Thus $X$ is also fibrant in $\motivicsymmTspectra$.
	Since the adjunction $(V,U,\varphi)$
	is a Quillen equivalence between $\motivicTspectra$ and $\motivicsymmTspectra$, 
	\cite[proposition 1.3.13(b)]{MR1650134} implies that the following composition
	is a weak equivalence in $\motivicsymmTspectra$:
		$$\xymatrix{VQ_{s}U(X) \ar[rr]^-{V(Q_{s}^{UX})}&&
													VU(X)\ar[r]^-{\epsilon _{X}}& X}
		$$
	Hence using \cite[proposition 3.1.5]{MR1944041} it follows that the composition above
	is a $\qleftsymmlocmaps$-local equivalence.
	
	(\ref{thm.3.3.symmetrization-Quillen-equivalence-weight<q.b}):  This follows immediately
	from propositions \ref{prop.3.2.Qs-cofibrant-replacement-all-L<q} and 
	\ref{prop.3.3.V-reflects-Lq-localequivalences}.
\end{proof}

\begin{cor}
		\label{cor.3.3.symmetrization-functor-exact-equivalenceL<q}
	Fix $q\in \mathbb Z$.  Then the adjunction
		$$\xymatrix{(V,U,\varphi):\weightqTspectra \ar[r]& \weightqsymmTspectra}
		$$
	given by the symmetrization and the forgetful functors, induces
	an adjunction
		$$\xymatrix{(VQ_{s},UW_{q}^{\Sigma},\varphi):\weightqstablehomotopy \ar[r]& \weightqsymmstablehomotopy}
		$$
	of exact funtors between triangulated categories.
	Furthermore, $VQ_{s}$ and $UW_{q}^{\Sigma}$ are both equivalences of categories.
\end{cor}
\begin{proof}
	Theorem \ref{thm.3.3.symmetrization-Quillen-equivalence-weight<q} implies that the adjunction $(V,U,\varphi)$
	is a Quillen equivalence.
	Therefore we get
	the following adjunction at the level of the associated homotopy categories:
		$$\xymatrix{(VQ_{s}, UW_{q}^{\Sigma}, \varphi ):\weightqstablehomotopy \ar[r]& \weightqsymmstablehomotopy}
		$$
	
	Now \cite[proposition 1.3.13]{MR1650134} implies that $VQ_{s} ,UW_{q}^{\Sigma}$ are both equivalences of
	categories.
	Finally, proposition \ref{prop.2.6.enriched-symmetrization-adjunction}
	together with \cite[proposition 6.4.1]{MR1650134} imply that
	$VQ_{s}$ maps cofibre sequences in $\weightqstablehomotopy$ to cofibre sequences in
	$\weightqsymmstablehomotopy$.
	Therefore 
	using proposition 7.1.12 in \cite{MR1650134} we have
	that $VQ_{s}$ and $UW_{q}^{\Sigma}$ are both exact functors between triangulated categories.
\end{proof}

	Now it is very easy to find the desired lifting 
	for the functor $\tilde{s}_{<q}:\symmstablehomotopy \rightarrow \symmstablehomotopy$
	(see corollary \ref{cor.3.3.symmetric-fq,s<q,sq}(\ref{cor.3.3.symmetric-fq,s<q,sq.b})) to the
	model category level.
	
\begin{lem}
		\label{lem.3.3.homotopycoherence==>liftingssymmweight<q}
	Fix $q\in \mathbb Z$, and let $X$ be an arbitrary symmetric $T$-spectrum.
	\begin{enumerate}
		\item \label{lem.3.3.homotopycoherence==>liftingssymmweight<q.a}The following maps in $\weightqsymmTspectra$
						$$\xymatrix{Q_{\Sigma}(VQ_{s}X)\ar[rr]^-{Q_{\Sigma}^{VQ_{s}X}}&& 
												VQ_{s}X &&  VQ_{s}(Q_{s}X) \ar[ll]_-{V(Q_{s}^{Q_{s}X})}}
						$$
					induce natural isomorphisms between the functors:
						$$Q_{\Sigma}\circ VQ_{s}, VQ_{s}, VQ_{s}\circ Q_{s}:\stablehomotopy
							\rightarrow \weightqsymmstablehomotopy$$
						
						$$\xymatrix{& \symmstablehomotopy \ar[dr]^-{Q_{\Sigma}}&\\
												 \stablehomotopy
												\ar[rr]^{VQ_{s}} \ar[ur]^-{VQ_{s}} \ar[dr]_-{Q_{s}}&& \weightqsymmstablehomotopy \\
												& \weightqstablehomotopy \ar[ur]_-{VQ_{s}}&}
						$$
					Given a $T$-spectrum $X$
					$$\xymatrix{\kappa _{X}:Q_{\Sigma}(VQ_{s}X)\ar[r]^-{\cong}& VQ_{s}(Q_{s}X)}$$ 
					will denote
					the isomorphism in $\weightqsymmstablehomotopy$ corresponding to the natural isomorphism
					between $Q_{\Sigma}\circ VQ_{s}$
					and $VQ_{s}\circ Q_{s}$.
		\item \label{lem.3.3.homotopycoherence==>liftingssymmweight<q.b}The following maps in $\motivicTspectra$
						$$\xymatrix{UR_{\Sigma}(W_{q}^{\Sigma}X)&& 
												UW_{q}^{\Sigma}X \ar[ll]_-{U(R_{\Sigma}^{W_{q}^{\Sigma}X})} \ar[rr]^-{W_{q}^{UW_{q}^{\Sigma}X}}&& 
												W_{q}(UW_{q}^{\Sigma}X) }
						$$
					induce natural isomorphisms between the functors:
						$$UR_{\Sigma}\circ W_{q}^{\Sigma}, UW_{q}^{\Sigma}, W_{q}\circ UW_{q}^{\Sigma}:\weightqsymmstablehomotopy
							\rightarrow \stablehomotopy$$
						
						$$\xymatrix{& \symmstablehomotopy \ar[dr]^-{UR_{\Sigma}}&\\
												 \weightqsymmstablehomotopy
												\ar[rr]^{UW_{q}^{\Sigma}} \ar[ur]^-{W_{q}^{\Sigma}} \ar[dr]_-{UW_{q}^{\Sigma}}&& \stablehomotopy \\
												& \weightqstablehomotopy \ar[ur]_-{W_{q}}&}
						$$
					Given a symmetric $T$-spectrum $X$
					$$\xymatrix{\mu _{X}:UR_{\Sigma}(W_{q}^{\Sigma}X)\ar[r]^-{\cong}& W_{q}(UW_{q}^{\Sigma}X)}$$ 
					will denote
					the isomorphism in $\stablehomotopy$ corresponding to the natural isomorphism
					between $UR_{\Sigma}\circ W_{q}^{\Sigma}$
					and $W_{q}\circ UW_{q}^{\Sigma}$.
	\end{enumerate}
\end{lem}	
\begin{proof}
	(\ref{lem.3.3.homotopycoherence==>liftingssymmweight<q.a}):  Follows immediately from 
	theorem 1.3.7 in \cite{MR1650134} and the following
	commutative diagram of left Quillen functors:
		$$\xymatrix{\motivicTspectra \ar[r]^-{V} \ar[d]_-{id}& \motivicsymmTspectra \ar[d]^-{id}\\
								\weightqTspectra \ar[r]_-{V}& \weightqsymmTspectra}
		$$
		
	(\ref{lem.3.3.homotopycoherence==>liftingssymmweight<q.b}):  Follows immediately from the dual of
	theorem 1.3.7 in \cite{MR1650134} and the following
	commutative diagram of right Quillen functors:
		$$\xymatrix{\motivicTspectra & \motivicsymmTspectra \ar[l]_-{U}\\
								\weightqTspectra \ar[u]^-{id}& \weightqsymmTspectra \ar[l]^-{U} \ar[u]_-{id}}
		$$
\end{proof}
	
\begin{thm}
		\label{thm.3.3.symmLq-models-symms<q}
	Fix $q\in \mathbb Z$, and let $X$ be an arbitrary symmetric $T$-spectrum.
	\begin{enumerate}
		\item \label{thm.3.3.symmLq-models-symms<q.a}The diagram (\ref{diagram.3.2.s<q-lifting})
					in theorem \ref{thm.3.2.Lq-models-s<q} induces
					the following diagram in 
					$\symmstablehomotopy$:
						\begin{equation}
										\label{diagram.3.3.symms<q-lifting.a}
							\begin{array}{c}
									\xymatrix{VQ_{s}(Q_{s}s_{<q}(UR_{\Sigma}X)) \ar[dd]_-{VQ_{s}(Q_{s}^{s_{<q}UR_{\Sigma}X})}^-{\cong}
														\ar[dr]^-{\ \ \ \ \ \ VQ_{s}(W_{q}^{Q_{s}s_{<q}(UR_{\Sigma}X)})}_-{\cong}&\\ 
														& VQ_{s}(W_{q}Q_{s}s_{<q}(UR_{\Sigma}X)) \\  
														\tilde{s}_{<q}X=VQ_{s}(s_{<q}(UR_{\Sigma}X))&\\
														& VQ_{s}(W_{q}Q_{s}(UR_{\Sigma}X)) 
														\ar[uu]_-{VQ_{s}(W_{q}Q_{s}(\pi _{<q} ^{UR_{\Sigma}X}))}^-{\cong}}
							\end{array}
						\end{equation}
					where all the maps are isomorphisms in $\symmstablehomotopy$.
					Furthermore, this diagram induces a
					natural isomorphism  
					between the following exact functors:
						$$\xymatrix{\symmstablehomotopy \ar@<1ex>[rrr]^-{\tilde{s}_{<q}}  
												\ar@<-1ex>[rrr]_-{VQ_{s}\circ W_{q}Q_{s} \circ UR_{\Sigma}} &&& \symmstablehomotopy}
						$$
		\item	\label{thm.3.3.symmLq-models-symms<q.b}Let $\eta$ be the
					unit of the adjunction (see corollary \ref{cor.3.3.symmetrization-functor-exact-equivalenceL<q}):
						$$\xymatrix{(VQ_{s},UW_{q}^{\Sigma},\varphi):\weightqstablehomotopy \ar[r]&
												\weightqsymmstablehomotopy}
						$$						
					Then we have the following diagram in $\stablehomotopy$
					(see lemma \ref{lem.3.3.homotopycoherence==>liftingssymmweight<q}):
						\begin{equation}
										\label{diagram.3.3.symms<q-lifting.b}
							\begin{array}{c}
									\xymatrix{W_{q}(UW_{q}^{\Sigma}(Q_{\Sigma}(VQ_{s}X))) \ar[rr]^-{W_{q}UW_{q}^{\Sigma}(\kappa _{X})}_-{\cong}&&
														W_{q}(UW_{q}^{\Sigma}(VQ_{s}(Q_{s}X))) \ar[dd]^-{W_{q}(\eta _{Q_{s}X})^{-1}}_-{\cong}\\
														&&\\
														UR_{\Sigma}(W_{q}^{\Sigma}Q_{\Sigma}(VQ_{s}X))\ar[uu]^-{\mu _{Q_{\Sigma}(VQ_{s}X)}}_-{\cong}
														&& W_{q}Q_{s}X}
							\end{array}
						\end{equation}
					where all the maps are isomorphisms in $\symmstablehomotopy$.  This diagram induces a natural isomorphism
					between the following exact functors:
						$$\xymatrix{\stablehomotopy \ar@<1ex>[rrr]^-{UR_{\Sigma}\circ W_{q}^{\Sigma}Q_{\Sigma}\circ VQ_{s}}  
												\ar@<-1ex>[rrr]_-{W_{q}Q_{s}} &&& \stablehomotopy}
						$$
		\item	\label{thm.3.3.symmLq-models-symms<q.c}Let $\epsilon$
					denote the counit of the adjunction (see theorem \ref{thm.3.3.symmetrization-functor-exact-equivalence}):
						$$\xymatrix{(VQ_{s},UR_{\Sigma},\varphi):\stablehomotopy \ar[r]&
												\symmstablehomotopy}
						$$
					and let $\gamma$ denote the natural isomorphism constructed above in (\ref{thm.3.3.symmLq-models-symms<q.b}).
					Then we have the following diagram in $\symmstablehomotopy$:
						\begin{equation}
										\label{diagram.3.3.symms<q-lifting.c}
							\begin{array}{c}
									\xymatrix{VQ_{s}(UR_{\Sigma}W_{q}^{\Sigma}Q_{\Sigma}VQ_{s}(UR_{\Sigma}X)) 
														\ar[dr]^-{\ \ \epsilon _{W_{q}^{\Sigma}Q_{\Sigma}VQ_{s}(UR_{\Sigma}X)}}_-{\cong}&\\
														& W_{q}^{\Sigma}Q_{\Sigma}VQ_{s}(UR_{\Sigma}X) \ar[dd]^-{W_{q}^{\Sigma}Q_{\Sigma}(\epsilon _{X})}_-{\cong}\\
														VQ_{s}(W_{q}Q_{s}(UR_{\Sigma}X))\ar[uu]^-{VQ_{s}(\gamma _{UR_{\Sigma}X})^{-1}}_-{\cong}&\\
														&W_{q}^{\Sigma}Q_{\Sigma}X=s_{<q}^{\Sigma}X}
							\end{array}
						\end{equation}
					where all the maps are isomorphisms in $\symmstablehomotopy$.  This diagram induces a natural isomorphism
					between the following exact functors:
						$$\xymatrix{\symmstablehomotopy \ar@<1ex>[rrr]^-{VQ_{s}\circ W_{q}Q_{s}\circ UR_{\Sigma}}  
												\ar@<-1ex>[rrr]_-{W_{q}^{\Sigma}Q_{\Sigma}=s_{<q}^{\Sigma}} &&& \symmstablehomotopy}
						$$
		\item	\label{thm.3.3.symmLq-models-symms<q.d}Combining the diagrams (\ref{diagram.3.3.symms<q-lifting.a}) and 
					(\ref{diagram.3.3.symms<q-lifting.c}) above
					we get a natural isomorphism 
					between the following exact functors:
						$$\xymatrix{\symmstablehomotopy \ar@<1ex>[rr]^-{\tilde{s}_{<q}}  
												\ar@<-1ex>[rr]_-{s_{<q}^{\Sigma}} && \symmstablehomotopy}
						$$
	\end{enumerate}
\end{thm}
\begin{proof}
	It is clear that it suffices to prove only the first three claims.
	
	(\ref{thm.3.3.symmLq-models-symms<q.a}):  Follows immediately from theorems
	\ref{thm.3.2.Lq-models-s<q} and \ref{thm.3.3.symmetrization-functor-exact-equivalence}.
	
	(\ref{thm.3.3.symmLq-models-symms<q.b}):  Follows immediately from lemma \ref{lem.3.3.homotopycoherence==>liftingssymmweight<q}
	and corollary \ref{cor.3.3.symmetrization-functor-exact-equivalenceL<q}.
	
	(\ref{thm.3.3.symmLq-models-symms<q.c}):  Follows immediately from (\ref{thm.3.3.symmLq-models-symms<q.b}) above,
	and theorem \ref{thm.3.3.symmetrization-functor-exact-equivalence}.
\end{proof}

	The functor $s_{<q}^{\Sigma}$ gives the desired lifting for
	the functor $\tilde{s}_{<q}$
	to the model category level.

\begin{prop}
		\label{prop.Lq+1-->Lqsymm}
	For every $q\in \mathbb Z$, we have the following
	commutative diagram of left Quillen functors:
		$$\xymatrix{& \motivicsymmTspectra \ar[dl]_-{id} \ar[dr]^-{id}&\\
								L_{<q+1}\motivicsymmTspectra \ar[rr]_-{id}&& \weightqsymmTspectra}
		$$
\end{prop}
\begin{proof}
	Since $\weightqsymmTspectra$ and $L_{<q+1}\motivicsymmTspectra$ are both
	left Bousfield localizations for $\motivicsymmTspectra$, we have that
	the identity functors:
		$$\xymatrix@R=.5pt{id:\motivicsymmTspectra \ar[r]& \weightqsymmTspectra \\
								id:\motivicsymmTspectra \ar[r]& L_{<q+1}\motivicsymmTspectra}
		$$
	are both left Quillen functors.  Hence, it suffices to show that
		$$\xymatrix{id:L_{<q+1}\motivicsymmTspectra \ar[r]& \weightqsymmTspectra}
		$$
	is a left Quillen functor.  Using the universal property for left
	Bousfield localizations (see definition \ref{def-1.leftlocmodcats}), 
	we have that it is enough to check
	that if $f:X\rightarrow Y$ is a $L^{\Sigma}(<q+1)$-local equivalence then
	$Q_{\Sigma}f:Q_{\Sigma}X\rightarrow Q_{\Sigma}Y$ is a $\qleftsymmlocmaps$-local equivalence.
	
	But theorem 3.1.6(c) in \cite{MR1944041} implies that this last condition is equivalent
	to the following one: Let $Z$ be an arbitrary $\qleftsymmlocmaps$-local symmetric $T$-spectrum,
	then $Z$ is also $L^{\Sigma}(<q+1)$-local.
	Finally, this last condition follows immediately from proposition 
	\ref{prop.3.3.Z-symmLq-local-iff-UZ-Lq-local} and corollary
	\ref{cor.3.2.m>n===>Ln-local=>Lm-local}.
\end{proof}

\begin{cor}
		\label{cor.3.3.q+1-->qstableadjsymm}
	For every $q\in \mathbb Z$, we have the following
	adjunction
		$$\xymatrix{(Q_{\Sigma},W_{q}^{\Sigma},\varphi ):\weightqplusonesymmstablehomotopy \ar[r]& 
								\weightqsymmstablehomotopy}
		$$
	of exact functors between triangulated categories.
\end{cor}
\begin{proof}
	Proposition \ref{prop.Lq+1-->Lqsymm} implies that $id:L_{<q+1}\motivicsymmTspectra \rightarrow \weightqsymmTspectra$
	is a left Quillen functor.  Therefore we get the following
	adjunction at the level of the associated homotopy categories
		$$\xymatrix{(Q_{\Sigma},W_{q}^{\Sigma},\varphi):L_{<q+1}\symmstablehomotopy \ar[r]& \weightqsymmstablehomotopy}
		$$
	Now proposition 6.4.1 in \cite{MR1650134} implies that
	$Q_{\Sigma}$ maps cofibre sequences in $L_{<q+1}\symmstablehomotopy$ to cofibre sequences in
	$\weightqsymmstablehomotopy$.
	Therefore 
	using proposition 7.1.12 in \cite{MR1650134} we have
	that $Q_{\Sigma}$ and $W_{q}^{\Sigma}$ are both exact functors between triangulated categories.
\end{proof}

\begin{thm}
		\label{thm.3.3.motivic-towersymm}
	We have the following tower of left Quillen functors:
		\begin{equation}
					\label{diagram.3.3.motivictower-modelcatlevelsymm}
			\begin{array}{c}
				\xymatrix{& \vdots \ar[d]^-{id} \\
									& L_{<q+1}\motivicsymmTspectra \ar[d]^-{id} \\
									\motivicsymmTspectra \ar[ur]^-{id} \ar[r]^-{id} \ar[dr]_-{id}& \weightqsymmTspectra \ar[d]^-{id} \\
									& L_{<q-1}\motivicsymmTspectra \ar[d]^-{id} &\\
									& \vdots}
			\end{array}
		\end{equation}
	together with the corresponding tower of associated homotopy categories:
		\begin{equation}
					\label{diagram.3.3.motivictower-homotopylevelsymm}
			\begin{array}{c}
		  	\xymatrix{&&& \vdots \ar@<-1ex>[d]_-{Q_{\Sigma}} \\
									&&& L_{<q+1}\symmstablehomotopy \ar@<-1ex>[d]_-{Q_{\Sigma}} \ar@<-1ex>[u]_-{W_{q+1}^{\Sigma}} 
									\ar[dlll]|-{W_{q+1}^{\Sigma}}\\
									\symmstablehomotopy \ar@<2ex>[urrr]|-{Q_{\Sigma}} \ar@<1ex>[rrr]|-{Q_{\Sigma}} \ar[drrr]|-{Q_{\Sigma}}&&&
									\weightqsymmstablehomotopy \ar@<-1ex>[d]_-{Q_{\Sigma}} \ar@<-1ex>[u]_-{W_{q}^{\Sigma}}
									\ar@<1ex>[lll]|-{W_{q}^{\Sigma}}\\
									&&& L_{<q-1}\symmstablehomotopy \ar@<-1ex>[d]_-{Q_{\Sigma}} \ar@<-1ex>[u]_-{W_{q-1}^{\Sigma}}
									\ar@<2ex>[ulll]|-{W_{q-1}^{\Sigma}}\\
									&&& \vdots \ar@<-1ex>[u]_-{W_{q-2}^{\Sigma}}}
			\end{array}
		\end{equation}
	The tower (\ref{diagram.3.3.motivictower-homotopylevelsymm})
	gets canonically identified, through the equivalences of categories
	$VQ_{s}$, $UR_{\Sigma}$ and $UW_{q}^{\Sigma}$ constructed in 
	theorem \ref{thm.3.3.symmetrization-functor-exact-equivalence} and 
	corollary \ref{cor.3.3.symmetrization-functor-exact-equivalenceL<q};
	with the tower (\ref{diagram.3.2.motivictower-homotopylevel}) defined in theorem
	\ref{thm.3.2.motivic-tower}.  Moreover, this tower also
	satisfies the following properties:
		\begin{enumerate}
			\item	All the categories are triangulated.
			\item	All the functors are exact.
			\item	$Q_{\Sigma}$ is a left adjoint for all the functors $W_{q}^{\Sigma}$. 
		\end{enumerate}
\end{thm}
\begin{proof}
	Follows immediately from propositions
	\ref{prop.3.3.symmLq-exact-adjunctions}, \ref{prop.Lq+1-->Lqsymm}, 
	corollary \ref{cor.3.3.q+1-->qstableadjsymm} together with
	theorem \ref{thm.3.3.symmetrization-functor-exact-equivalence} and 
	corollary \ref{cor.3.3.symmetrization-functor-exact-equivalenceL<q}.
\end{proof}

\begin{defi}
		\label{def.3.3.symmetricSq-colocal-generators}
	For every $q\in \mathbb Z$, we consider the following set
	of symmetric $T$-spectra
		$$\symmqslicegenerators =\{ \symmgeneratorNRS \in C^{\Sigma} | s-n=q \}\subseteq C_{eff}^{q,\Sigma}
		$$
	(see proposition \ref{prop.3.1.stablehomotopy=>compactly-generated} and 
	definition \ref{def.3.1.desusp-stable-homotopy-eff}).
\end{defi}

\begin{thm}
		\label{thm.3.3.symmetricSq-modelcats}
	Fix $q\in \mathbb Z$.  Then the right Bousfield
	localization of the model category $\weightqplusonesymmTspectra$
	with respect to the $\symmqslicegenerators$-colocal equivalences
	exists.  This new model structure will be called \emph{$q$-slice
	motivic symmetric stable}.
	$\qslicesymmTspectra$ will denote the category of symmetric $T$-spectra
	equipped with the $q$-slice motivic symmetric stable model structure, and 
	$\qslicesymmstablehomotopy$ will denote its associated homotopy category.
	Furthermore, the $q$-slice motivic symmetric stable model structure is right proper and simplicial.
\end{thm}
\begin{proof}
	Theorem \ref{thm.3.3.Lqsymmetricmodelstructures} implies that $\weightqplusonesymmTspectra$ is a cellular
	and simplicial model category.  On the other hand, corollary
	\ref{cor.3.3.symmLq=>rightproper} implies that $\weightqplusonesymmTspectra$ is right proper.
	Therefore we can apply theorem 5.1.1 in \cite{MR1944041} to construct
	the right Bousfield localization
	of $\weightqplusonesymmTspectra$ with respect to the $\symmqslicegenerators$-colocal equivalences.
	Using \cite[theorem 5.1.1]{MR1944041} again, we have that $\qslicesymmTspectra$ is a right proper
	and simplicial model category.
\end{proof}

\begin{defi}
		\label{def.3.3.Pqsigma-cofibrant-replacement}
	Fix $q\in \mathbb Z$.  Let $P_{q}^{\Sigma}$ denote a cofibrant replacement functor in 
	$\qslicesymmTspectra$; such that for every symmetric $T$-spectrum $X$, the natural map
		$$\xymatrix{P_{q}^{\Sigma}X \ar[r]^-{P_{q}^{\Sigma ,X}}& X}
		$$ 
	is a trivial fibration in $\qslicesymmTspectra$, and
	$P_{q}^{\Sigma}X$ is always a $\symmqslicegenerators$-colocal symmetric $T$-spectrum
	in $\weightqplusonesymmTspectra$.
\end{defi}

\begin{prop}
		\label{prop.3.3.Wsigmaqplusone-fibrant-replacement-all-Sq}
	Fix $q\in \mathbb Z$.  Then $W^{\Sigma}_{q+1}$ is also a fibrant
	replacement functor in $\qslicesymmTspectra$
	(see definition \ref{def.3.3.stable-weightqsymm-replacementfunctors}),
	and for every symmetric $T$-spectrum $X$ the natural map
		$$\xymatrix{X\ar[rr]^-{W^{\Sigma ,X}_{q+1}}&& W^{\Sigma}_{q+1}X}
		$$
	is a trivial cofibration in $\qslicesymmTspectra$.
\end{prop}
\begin{proof}
	Since $\qslicesymmTspectra$ is the right Bousfield localization of
	$\weightqplusonesymmTspectra$ with respect to the $\symmqslicegenerators$-colocal equivalences, by construction
	we have that the fibrations and the trivial cofibrations are indentical in
	$\qslicesymmTspectra$ and $\weightqplusonesymmTspectra$ respectively.  This implies that for every symmetric
	$T$-spectrum $X$, $W^{\Sigma}_{q+1}X$ is fibrant in $\qslicesymmTspectra$,
	and we also have that
	the natural map
		$$\xymatrix{X\ar[rr]^{W^{\Sigma ,X}_{q+1}}&& W^{\Sigma}_{q+1}X}
		$$
	is a trivial cofibration in $\qslicesymmTspectra$.  Hence $W^{\Sigma}_{q+1}$ is also a fibrant replacement
	functor for $\qslicesymmTspectra$.
\end{proof}

\begin{prop}
		\label{prop.3.3.f-Sqeff-colocal=iff=Uf-Sqeff-colocal}
	Fix $q\in \mathbb Z$, and let $f:X\rightarrow Y$ be a map in $\weightqplusonesymmTspectra$.
	Then $f$ is a $\symmqslicegenerators$-colocal equivalence
	in $\weightqplusonesymmTspectra$
	if and only if the underlying map
	$UW_{q+1}^{\Sigma}(f):UW_{q+1}^{\Sigma}X\rightarrow UW_{q+1}^{\Sigma}Y$ is a $\qslicegenerators$-colocal equivalence
	in $\weightqplusoneTspectra$.
\end{prop}
\begin{proof}
	Consider $\symmgeneratorNRS \in \symmqslicegenerators$.
	Using the enriched adjunctions of proposition
	\ref{prop.2.6.enriched-symmetrization-adjunction}, 
	we get the following commutative diagram where
	the vertical arrows are all isomorphisms:
		$$\xymatrix@C=-1pc{Map\: _{\Sigma}(\symmgeneratorNRS ,W_{q+1}^{\Sigma}X) \ar[dr]^-{W_{q+1}^{\Sigma}f_{\ast}} \ar@{=}[dd]&\\
								& Map\: _{\Sigma}(\symmgeneratorNRS ,W_{q+1}^{\Sigma}Y) \ar@{=}[dd]\\
								Map\: _{\Sigma}(V(\generatorNRS) ,W_{q+1}^{\Sigma}X) \ar[dr]^-{W_{q+1}^{\Sigma}f_{\ast}} \ar[dd]_-{\cong}&\\
								& Map\: _{\Sigma}(V(\generatorNRS) ,W_{q+1}^{\Sigma}Y) \ar[dd]^-{\cong}\\
								Map (\generatorNRS ,UW_{q+1}^{\Sigma}X) \ar[dr]_-{UW_{q+1}^{\Sigma}f_{\ast}}&\\ 
								& Map (\generatorNRS ,UW_{q+1}^{\Sigma}Y)}
		$$
	Since $UW_{q+1}^{\Sigma}X$ and $UW_{q+1}^{\Sigma}Y$ are both fibrant in $\weightqplusoneTspectra$, 
	we have that $UW_{q+1}^{\Sigma}(f)$
	is a $\qslicegenerators$-colocal equivalence in $\weightqplusoneTspectra$ if and only if
	the bottom row in the diagram above is a weak equivalence of simplicial sets
	for every $\generatorNRS \in \qslicegenerators$.  By the two out of three property for
	weak equivalences we have that this happens if and only if the top
	row in the diagram above is a weak equivalence for every $\symmgeneratorNRS \in \symmqslicegenerators$.
	But this last condition holds if and only if $f$ is a $\symmqslicegenerators$-colocal equivalence
	in $\weightqplusonesymmTspectra$.  This finishes the proof.
\end{proof}

\begin{prop}
		\label{prop.3.3.detecting-symmetric-Sqeff-colocal-equivalences}
	Fix $q\in \mathbb Z$, and let $f:X\rightarrow Y$ be a map in $\weightqplusonesymmTspectra$.
	Then $f$ is a $\symmqslicegenerators$-colocal equivalence in 
	$\weightqplusonesymmTspectra$ if and only if
	for every $\symmgeneratorNRS \in \symmqslicegenerators$, the induced map:
		$$\xymatrix{[\symmgeneratorNRS , W_{q+1}^{\Sigma}X]_{Spt}^{\Sigma} 
								\ar[rr]^-{(W_{q+1}^{\Sigma}f)_{\ast}}&& [\symmgeneratorNRS , W_{q+1}^{\Sigma}Y]_{Spt}^{\Sigma}}
		$$
	is an isomorphism of abelian groups.
\end{prop}
\begin{proof}
	By proposition \ref{prop.3.3.f-Sqeff-colocal=iff=Uf-Sqeff-colocal}, $f$ is a $\symmqslicegenerators$-colocal
	equivalence in $\weightqplusonesymmTspectra$ if and only if $UW_{q+1}^{\Sigma}(f)$ is a $\qslicegenerators$-colocal equivalence
	in $\weightqplusoneTspectra$.  Since $UW_{q+1}^{\Sigma}X ,UW_{q+1}^{\Sigma}Y$ are both
	fibrant in $\weightqplusoneTspectra$,
	using proposition \ref{prop.3.2.classifying-Sq-colocal-equivs} 
	we have that $UW_{q+1}^{\Sigma}(f)$ is a $\qslicegenerators$-colocal equivalence if and
	only if for every $\generatorNRS \in \qslicegenerators$, the induced map
		$$\xymatrix{[\generatorNRS , UW_{q+1}^{\Sigma}X]_{Spt} \ar[rr]^-{UW_{q+1}^{\Sigma}(f)_{\ast}}&& 
								[\generatorNRS , UW_{q+1}^{\Sigma}Y]_{Spt}}
		$$
	is an isomorphism of abelian groups.
	
	Now since $W_{q+1}^{\Sigma}X ,W_{q+1}^{\Sigma}Y$ are also fibrant in $\motivicsymmTspectra$, 
	theorem \ref{thm.2.6.T-spectra===symmTspectra} implies that we have 
	the following commutative diagram, where all the vertical arrows
	are isomorphisms:
		$$\xymatrix{[\generatorNRS , UW_{q+1}^{\Sigma}X]_{Spt} \ar[dr]_-{(UW_{q+1}^{\Sigma}f)_{\ast}} \ar[dd]_-{\cong}&\\
								& [\generatorNRS , UW_{q+1}^{\Sigma}Y]_{Spt} \ar[dd]^-{\cong}\\
								[V(\generatorNRS), W_{q+1}^{\Sigma}X]_{Spt}^{\Sigma} \ar[dr]_-{(W_{q+1}^{\Sigma}f)_{\ast}} \ar@{=}[dd]&\\ 
								& [V(\generatorNRS), W_{q+1}^{\Sigma}Y]_{Spt}^{\Sigma} \ar@{=}[dd]\\
								[\symmgeneratorNRS, W_{q+1}^{\Sigma}X]_{Spt}^{\Sigma} \ar[dr]_-{(W_{q+1}^{\Sigma}f)_{\ast}}&\\ 
								& [\symmgeneratorNRS, W_{q+1}^{\Sigma}Y]_{Spt}^{\Sigma}}
		$$
	Therefore $f$ is a $\symmqslicegenerators$-colocal equivalence if and only if for every 
	$\symmgeneratorNRS \in \symmqslicegenerators$, the bottom row is an isomorphism of abelian groups.
	This finishes the proof.
\end{proof}

\begin{cor}
		\label{cor.3.3.classifying-Sq-symmetriccolocal-equivs.b}
	Fix $q\in \mathbb Z$ and let $f:X\rightarrow Y$ be a map of symmetric $T$-spectra.
	Then $f$ is a $\symmqslicegenerators$-colocal equivalence
	in $\weightqplusonesymmTspectra$ if and only if
		$$\xymatrix{W_{q+1}^{\Sigma}X \ar[rr]^-{W_{q+1}^{\Sigma}f}&& W_{q+1}^{\Sigma}Y}
		$$
	is a $C_{eff}^{q,\Sigma}$-colocal equivalence in $\motivicsymmTspectra$.
\end{cor}
\begin{proof}
	($\Rightarrow$):  Assume that $f$ is a $\symmqslicegenerators$-colocal equivalence, and
	fix $\symmgeneratorNRS \in C_{eff}^{q,\Sigma}$.  By proposition \ref{prop.3.3.detecting-symmetric-Cqeff-colocal-equivalences} 
	it suffices to show that the induced map
		\begin{equation}
					\label{equation.3.3.classifying-Sq-symmetriccolocal-equivs.b}
			\begin{array}{c}
				\xymatrix{[\symmgeneratorNRS, W_{q+1}^{\Sigma}X]_{Spt}^{\Sigma}\ar[d]_-{(W_{q+1}^{\Sigma}f)_{\ast}}\\ 
									[\symmgeneratorNRS, W_{q+1}^{\Sigma}Y]_{Spt}^{\Sigma}}
			\end{array}
		\end{equation}
	is an isomorphism of abelian groups.
	
	Since $\symmgeneratorNRS \in C_{eff}^{q,\Sigma}$, we have two possibilities:
		\begin{enumerate}
			\item \label{cond.a.cor.3.3.classifying-Sq-symmetriccolocal-equivs.b} $s-n=q$, i.e. 
						$\symmgeneratorNRS \in \symmqslicegenerators$.
			\item	\label{cond.b.cor.3.3.classifying-Sq-symmetriccolocal-equivs.b} $s-n\geq q+1$, i.e.
						$\symmgeneratorNRS \in C_{eff}^{q+1,\Sigma}$
		\end{enumerate}
	
	In case (\ref{cond.a.cor.3.3.classifying-Sq-symmetriccolocal-equivs.b}), 
	proposition \ref{prop.3.3.detecting-symmetric-Sqeff-colocal-equivalences} implies that the induced map
	in diagram (\ref{equation.3.3.classifying-Sq-symmetriccolocal-equivs.b}) is an isomorphism of abelian
	groups.
	
	On the other hand, in case (\ref{cond.b.cor.3.3.classifying-Sq-symmetriccolocal-equivs.b}),
	we have by proposition \ref{prop.3.3.Lqsymm-local-objects-classification}(\ref{prop.3.3.Lqsymm-local-objects-classification.b}) 
	that 
		$$[\symmgeneratorNRS ,W_{q+1}^{\Sigma}X]_{Spt}^{\Sigma}\cong 0\cong [\symmgeneratorNRS ,W_{q+1}^{\Sigma}Y]_{Spt}^{\Sigma}$$
	since by construction $W_{q+1}^{\Sigma}X$ and $W_{q+1}^{\Sigma}Y$ are both 
	$\qplusoneleftsymmlocmaps$-local symmetric $T$-spectra.
	Hence the induced map in diagram (\ref{equation.3.3.classifying-Sq-symmetriccolocal-equivs.b}) is also
	an isomorphism of abelian groups in this case, as we wanted.
	
	($\Leftarrow$):  Assume that $W_{q+1}^{\Sigma}f$ is a $C_{eff}^{q,\Sigma}$-colocal equivalence in $\motivicsymmTspectra$, and fix
	$\symmgeneratorNRS \in \symmqslicegenerators$.
	
	Since $\symmqslicegenerators \subseteq C_{eff}^{q,\Sigma}$, it follows from 
	proposition \ref{prop.3.3.detecting-symmetric-Cqeff-colocal-equivalences}
	that the induced map
		$$\xymatrix{[\symmgeneratorNRS ,W_{q+1}^{\Sigma}X]_{Spt}^{\Sigma} \ar[rr]^-{(W_{q+1}^{\Sigma}f)_{\ast}}&& 
								[\symmgeneratorNRS ,W_{q+1}^{\Sigma}Y]_{Spt}^{\Sigma}}
		$$
	is an isomorphism of abelian groups.
	Therefore, proposition \ref{prop.3.3.detecting-symmetric-Sqeff-colocal-equivalences} implies
	that $f$ is a $\symmqslicegenerators$-colocal equivalence in $\weightqplusonesymmTspectra$.  This finishes the proof.
\end{proof}

\begin{lem}
		\label{lem.3.3.Sqsymmetric-colocalequivs-stable-S1desuspension}
	Fix $q\in \mathbb Z$, and let $f:X\rightarrow Y$ be a map in $\weightqplusonesymmTspectra$.  Then $f$ is a
	$\symmqslicegenerators$-colocal equivalence in $\weightqplusonesymmTspectra$ if and only if
	$\Omega _{S^{1}} W_{q+1}^{\Sigma}f$ is a $\symmqslicegenerators$-colocal equivalence
	in $\weightqplusonesymmTspectra$.
\end{lem}
\begin{proof}
	It follows from proposition \ref{prop.3.3.f-Sqeff-colocal=iff=Uf-Sqeff-colocal} that 
	$f$ is a $\symmqslicegenerators$-colocal
	equivalence in $\weightqplusonesymmTspectra$ if and only if $UW_{q+1}^{\Sigma}f$ is a 
	$\qslicegenerators$-colocal equivalence
	in $\weightqplusoneTspectra$.  Since $UW_{q+1}^{\Sigma}X , UW_{q+1}^{\Sigma}Y$ 
	are both fibrant in $\weightqplusoneTspectra$, using
	lemma \ref{lem.3.2.Sq-colocalequivs-stable-S1desuspension} 
	we have that $UW_{q+1}^{\Sigma}f$ is a $\qslicegenerators$-colocal equivalence if and
	only if $\Omega _{S^{1}}UW_{q+1}^{\Sigma}f=U(\Omega_{S^{1}}W_{q+1}^{\Sigma}f)$ is a
	$\qslicegenerators$-colocal equivalence.
	
	Finally, since $\Omega _{S^{1}}W_{q+1}^{\Sigma}X , \Omega _{S^{1}}W_{q+1}^{\Sigma}Y$ are both
	fibrant in $\weightqplusonesymmTspectra$,
	we have by proposition \ref{prop.3.3.f-Sqeff-colocal=iff=Uf-Sqeff-colocal} that
	$U(\Omega _{S^{1}}W_{q+1}^{\Sigma}f)$ is a $\symmqslicegenerators$-colocal equivalence if and only if
	$\Omega _{S^{1}}W_{q+1}^{\Sigma}f$ is a $\symmqslicegenerators$-colocal equivalence.
	This finishes the proof.
\end{proof}

\begin{cor}
		\label{cor.3.3.Suspension=>qslicesymm-Quillen-equiv}
	Fix $q\in \mathbb Z$.  Then the adjunction
		$$\xymatrix{(-\wedge S^{1},\Omega _{S^{1}},\varphi):\qslicesymmTspectra \ar[r]& 
								\qslicesymmTspectra}
		$$
	is a Quillen equivalence.
\end{cor}
\begin{proof}
	Using corollary 1.3.16 in \cite{MR1650134} and
	proposition \ref{prop.3.3.Wsigmaqplusone-fibrant-replacement-all-Sq} we have that it suffices to
	verify the following two conditions:
		\begin{enumerate}
			\item	\label{cor.3.3.Suspension=>qslicesymm-Quillen-equiv.a}For every cofibrant object $X$
						in $\qslicesymmTspectra$, the following composition
							$$\xymatrix{X\ar[r]^-{\eta _{X}}& \Omega _{S^{1}}(X\wedge S^{1}) \ar[rr]^-{\Omega _{S^{1}}W_{q+1}^{\Sigma ,X\wedge S^{1}}}&& 
								\Omega _{S^{1}}W_{q+1}^{\Sigma}(X\wedge S^{1})}
							$$
						is a $\symmqslicegenerators$-colocal equivalence.
			\item	\label{cor.3.3.Suspension=>qslicesymm-Quillen-equiv.b}$\Omega _{S^{1}}$ reflects $\symmqslicegenerators$-colocal equivalences
						between fibrant objects in $\qslicesymmTspectra$.
		\end{enumerate}
		
	(\ref{cor.3.3.Suspension=>qslicesymm-Quillen-equiv.a}):  By construction $\qslicesymmTspectra$ is
	a right Bousfield localization of $\weightqplusonesymmTspectra$, therefore the identity
	functor 
		$$\xymatrix{id:\qslicesymmTspectra \ar[r]& \weightqplusonesymmTspectra}
		$$ 
	is a left Quillen functor.
	Thus $X$ is also cofibrant in $\weightqplusonesymmTspectra$.
	Since the adjunction $(-\wedge S^{1},\Omega _{S^{1}},\varphi)$
	is a Quillen equivalence on $\weightqplusonesymmTspectra$, 
	\cite[proposition 1.3.13(b)]{MR1650134} implies that the following composition
	is a weak equivalence in $\weightqplusonesymmTspectra$:
		$$\xymatrix{X\ar[r]^-{\eta _{X}}& \Omega _{S^{1}}(X\wedge S^{1}) \ar[rr]^-{\Omega _{S^{1}}W_{q+1}^{\Sigma ,X\wedge S^{1}}}&& 
								\Omega _{S^{1}}W_{q+1}^{\Sigma}(X\wedge S^{1})}
		$$
	Hence using \cite[proposition 3.1.5]{MR1944041} it follows that the composition above
	is a $\symmqslicegenerators$-colocal equivalence.
	
	(\ref{cor.3.3.Suspension=>qslicesymm-Quillen-equiv.b}):  This follows immediately
	from proposition \ref{prop.3.3.Wsigmaqplusone-fibrant-replacement-all-Sq} and 
	lemma \ref{lem.3.3.Sqsymmetric-colocalequivs-stable-S1desuspension}.
\end{proof}

\begin{rmk}
		\label{rmk.3.3.smashT-not-symmqslicedescending}
	The adjunction $(\Tsuspfunctor ,\Tloops ,\varphi)$ is
	a Quillen equivalence on $\motivicsymmTspectra$.
	However it does not descend even to a Quillen
	adjunction on the $q$-slice motivic symmetric stable
	model category
	$\qslicesymmTspectra$.
\end{rmk}

\begin{cor}
		\label{cor.3.3.symmqslice=>triangcat}
	For every $q\in \mathbb Z$,
	$\qslicesymmstablehomotopy$
	has the structure of a triangulated category.
\end{cor}
\begin{proof}
	Theorem \ref{thm.3.3.symmetricSq-modelcats} implies in particular
	that $\qslicesymmTspectra$ is a pointed simplicial model category,
	and corollary \ref{cor.3.3.Suspension=>qslicesymm-Quillen-equiv} implies that
	the adjunction 
		$$(-\wedge S^{1},\Omega _{S^{1}},\varphi):\qslicesymmTspectra \rightarrow \qslicesymmTspectra$$
	is a Quillen equivalence.  Therefore
	the result follows from the work of Quillen in 
	\cite[sections I.2 and I.3]{MR0223432} and the work of 
	Hovey in \cite[chapters VI and VII]{MR1650134}.
\end{proof}

\begin{prop}
		\label{prop.3.3.symmPqcofibrant-replacement=>triangulatedfunctor}
	Fix $q\in \mathbb Z$.  Then we have the following adjunction
		$$\xymatrix{(P_{q}^{\Sigma}, W_{q+1}^{\Sigma}, \varphi) :\qslicesymmstablehomotopy \ar[r]& \weightqplusonesymmstablehomotopy}
		$$
	between exact functors of triangulated categories.
\end{prop}
\begin{proof}
	Since $\qslicesymmTspectra$ is the right Bousfield localization of
	$\weightqplusonesymmTspectra$ with respect to the $\symmqslicegenerators$-colocal equivalences, we have that the
	identity functor $id:\qslicesymmTspectra \rightarrow \weightqplusonesymmTspectra$
	is a left Quillen functor.  Therefore we get
	the following adjunction at the level of the associated homotopy categories:
		$$\xymatrix{(P_{q}^{\Sigma}, W_{q+1}^{\Sigma}, \varphi ):\qslicesymmstablehomotopy \ar[r]& \weightqplusonesymmstablehomotopy}
		$$
	
	Now proposition 6.4.1 in \cite{MR1650134} implies that
	$P_{q}^{\Sigma}$ maps cofibre sequences in $\qslicesymmstablehomotopy$ to cofibre sequences in
	$\weightqplusonesymmstablehomotopy$.
	Therefore 
	using proposition 7.1.12 in \cite{MR1650134} we have
	that $P_{q}^{\Sigma}$ and $W_{q+1}^{\Sigma}$ are both exact functors between triangulated categories.
\end{proof}

\begin{prop}
		\label{prop.3.3.symmetricq-connected--->q-slice===leftQuilllenfunctor}
	Fix $q\in \mathbb Z$.  Then the identity functor
		$$\xymatrix{id:\qslicesymmTspectra \ar[r]& \qconnectedsymmTspectra}
		$$
	is a right Quillen functor.
\end{prop}
\begin{proof}
	Consider the following diagram of right Quillen functors
		$$\xymatrix{\weightqplusonesymmTspectra \ar[r]^-{id} \ar[d]_-{id}& \motivicsymmTspectra \ar[r]^-{id}& \qconnectedsymmTspectra \\
								\qslicesymmTspectra \ar@{-->}[urr]_-{id}&&}
		$$
	By the universal property of right Bousfield localizations 
	(see definition \ref{def-1.rightlocmodcats})
	it suffices to check that if $f:X\rightarrow Y$ is a $\symmqslicegenerators$-colocal equivalence in
	$\weightqplusonesymmTspectra$, then $W_{q+1}^{\Sigma}f:W_{q+1}^{\Sigma}X\rightarrow W_{q+1}^{\Sigma}Y$ is a
	$C_{eff}^{q,\Sigma}$-colocal equivalence in $\motivicsymmTspectra$.
	But this follows immediately from 
	corollary \ref{cor.3.3.classifying-Sq-symmetriccolocal-equivs.b}.
\end{proof}

\begin{cor}
		\label{cor.3.3.adjunctions--symmetricRq==>Sq}
	For every $q\in \mathbb Z$ we have
	the following adjunction
		$$\xymatrix{(C_{q}^{\Sigma},W_{q+1}^{\Sigma},\varphi ):\qconnectedsymmstablehomotopy \ar[r]& \qslicesymmstablehomotopy}
		$$
	of exact functors between triangulated categories.
\end{cor}
\begin{proof}
	By proposition \ref{prop.3.3.symmetricq-connected--->q-slice===leftQuilllenfunctor} the identity functor
	$id:\qconnectedsymmTspectra \rightarrow \qslicesymmTspectra$ is a left Quillen
	functor.  Therefore we get
	the following adjunction at the level of the associated homotopy categories:
		$$\xymatrix{(C_{q}^{\Sigma}, W_{q+1}^{\Sigma}, \varphi ):\qconnectedsymmstablehomotopy \ar[r]& \qslicesymmstablehomotopy}
		$$
	
	Now proposition 6.4.1 in \cite{MR1650134} implies that
	$C_{q}^{\Sigma}$ maps cofibre sequences in $\qconnectedsymmstablehomotopy$ to cofibre sequences in
	$\qslicesymmstablehomotopy$.
	Therefore 
	using proposition 7.1.12 in \cite{MR1650134} we have
	that $C_{q}^{\Sigma}$ and $W_{q+1}^{\Sigma}$ are both exact functors between triangulated categories.
\end{proof}

\begin{lem}
		\label{lem.3.3.A-symmSqcofibrant====>AsymmLq+1trivially-cofibrant}
	Fix $q\in \mathbb Z$, and let $A$ be a cofibrant symmetric $T$-spectrum
	in $\qslicesymmTspectra$.  Then the map $\ast \rightarrow A$
	is a trivial cofibration in $\weightqsymmTspectra$.
\end{lem}
\begin{proof}
	Let $Z$ be an arbitrary $\qleftsymmlocmaps$-local symmetric $T$-spectrum in $\motivicsymmTspectra$.
	We claim that the map $Z\rightarrow \ast$ is a trivial fibration in $\qslicesymmTspectra$.
	In effect, using proposition \ref{prop.3.3.Z-symmLq-local-iff-UZ-Lq-local} and
	corollary \ref{cor.3.2.m>n===>Ln-local=>Lm-local} 
	we have that $Z$ is $\qplusoneleftsymmlocmaps$-local in $\motivicsymmTspectra$, i.e. a fibrant
	object in $\weightqplusonesymmTspectra$.  By construction $\qslicesymmTspectra$ is a right Bousfield
	localization of $\weightqplusonesymmTspectra$, hence $Z$ is also fibrant in $\qslicesymmTspectra$.
	Then by proposition \ref{prop.3.3.detecting-symmetric-Sqeff-colocal-equivalences} 
	it suffices to show that for every
	$\symmgeneratorNRS \in S^{\Sigma}(q)$ (i.e. $s-n=q$):
		$$\xymatrix{0\cong [\generatorNRS ,Z]_{Spt}^{\Sigma}}
		$$
	But this follows immediately from proposition \ref{prop.3.3.Lqsymm-local-objects-classification}, 
	since $Z$ is $\qleftsymmlocmaps$-local.
	
	Now since $\qslicesymmTspectra$ is a simplicial model category and $A$ is cofibrant in $\qslicesymmTspectra$,
	we have that the following map is a trivial fibration of simplicial sets:
		$$\xymatrix{Map\: _{\Sigma}(A,Z)\ar[r]& Map\: _{\Sigma}(A,\ast)=\ast}
		$$  
	The identity functor 
		$$\xymatrix{id:\qslicesymmTspectra \ar[r]& \weightqplusonesymmTspectra}
		$$
	is a left Quillen functor, since $\qslicesymmTspectra$ is a right Bousfield localization
	of $\weightqplusonesymmTspectra$.  Therefore $A$ is also cofibrant in $\weightqplusonesymmTspectra$, and
	since $\weightqplusonesymmTspectra$ is a left Bousfield localization of $\motivicsymmTspectra$; it follows
	that $A$ is also cofibrant in $\motivicsymmTspectra$.  On the other hand, we have that $Z$ is in particular
	fibrant in $\motivicsymmTspectra$.
	Hence $\pi _{0}Map\: _{\Sigma}(A,Z)$ computes $[A,Z]_{Spt}^{\Sigma}$,
	since $\motivicsymmTspectra$ is a simplicial model category.  But $Map\: _{\Sigma}(A,Z)\rightarrow \ast$
	is in particular a weak equivalence of simplicial sets, then
		$$[A,Z]_{Spt}^{\Sigma}\cong 0
		$$
	for every $\qleftsymmlocmaps$-local symmetric $T$-spectrum $Z$.
	Finally, corollary \ref{cor.3.3.detecting-Cq-symmlocal-equivalences} 
	implies that $\ast \rightarrow A$ is a weak equivalence in $\weightqsymmTspectra$.
	This finishes the proof, since we already know
	that $A$ is cofibrant in $\weightqsymmTspectra$.
\end{proof}

\begin{thm}
		\label{thm.3.3.symmetrization-qslice-Quillen-equiv}
	Fix $q\in \mathbb Z$.  Then the adjunction
		$$\xymatrix{(V,U,\varphi):\qsliceTspectra \ar[r]& \qslicesymmTspectra}
		$$
	given by the symmetrization and the forgetful functors is a
	Quillen equivalence.
\end{thm}
\begin{proof}
	Proposition \ref{prop.3.3.f-Sqeff-colocal=iff=Uf-Sqeff-colocal} 
	together with the universal property for right Bousfield localizations
	(see definition \ref{def-1.rightlocmodcats})
	imply that 
		$$\xymatrix{U:\qslicesymmTspectra \ar[r]& \qsliceTspectra}$$ 
	is a right
	Quillen functor.	
	Using corollary 1.3.16 in \cite{MR1650134} and
	proposition \ref{prop.3.3.Wsigmaqplusone-fibrant-replacement-all-Sq} we have that it suffices to
	verify the following two conditions:
		\begin{enumerate}
			\item	\label{thm.3.3.symmetrization-qslice-Quillen-equiv.a}For every cofibrant object $X$
						in $\qsliceTspectra$, the following composition
							$$\xymatrix{X\ar[r]^-{\eta _{X}}& UV(X) \ar[rr]^-{UW_{q+1}^{\Sigma ,VX}}&& 
								UW_{q+1}^{\Sigma}V(X)}
							$$
						is a weak equivalence in $\qsliceTspectra$.
			\item	\label{thm.3.3.symmetrization-qslice-Quillen-equiv.b}$U$ reflects weak equivalences
						between fibrant objects in $\qslicesymmTspectra$.
		\end{enumerate}
		
	(\ref{thm.3.3.symmetrization-qslice-Quillen-equiv.a}):    
	By construction $\qsliceTspectra$ is
	a right Bousfield localization of $\weightqplusoneTspectra$, therefore the identity
	functor 
		$$\xymatrix{id:\qsliceTspectra \ar[r]& \weightqplusoneTspectra}
		$$ 
	is a left Quillen functor.
	Thus $X$ is also cofibrant in $\weightqplusoneTspectra$.
	Since the adjunction $(V,U,\varphi)$
	is a Quillen equivalence between $\weightqplusoneTspectra$  and $\weightqplusonesymmTspectra$, 
	\cite[proposition 1.3.13(b)]{MR1650134} implies that the following composition
	is a weak equivalence in $\weightqplusoneTspectra$:
		$$\xymatrix{X\ar[r]^-{\eta _{X}}& UV (X) \ar[rr]^-{UW_{q+1}^{\Sigma ,VX}}&& 
								UW_{q+1}^{\Sigma}V(X)}
		$$
	Hence using \cite[proposition 3.1.5]{MR1944041} it follows that the composition above
	is a $\qslicegenerators$-colocal equivalence in $\weightqplusoneTspectra$, i.e. a weak equivalence
	in $\qsliceTspectra$.
	
	(\ref{thm.3.3.symmetrization-qslice-Quillen-equiv.b}):  This follows immediately
	from propositions \ref{prop.3.3.Wsigmaqplusone-fibrant-replacement-all-Sq} and 
	\ref{prop.3.3.f-Sqeff-colocal=iff=Uf-Sqeff-colocal}.
\end{proof}

\begin{cor}
		\label{cor.3.3.symmetrization-functor-exact-equivalenceSq}
	Fix $q\in \mathbb Z$.  Then the adjunction
		$$\xymatrix{(V,U,\varphi):\qsliceTspectra \ar[r]& \qslicesymmTspectra}
		$$
	given by the symmetrization and the forgetful functors, induces
	an adjunction
		$$\xymatrix{(VP_{q},UW_{q+1}^{\Sigma},\varphi):\qslicestablehomotopy \ar[r]& \qslicesymmstablehomotopy}
		$$
	of exact funtors between triangulated categories.
	Furthermore, $VP_{q}$ and $UW_{q+1}^{\Sigma}$ are both equivalences of categories.
\end{cor}
\begin{proof}
	Theorem \ref{thm.3.3.symmetrization-qslice-Quillen-equiv} implies that the adjunction $(V,U,\varphi)$
	is a Quillen equivalence.
	Therefore we get
	the following adjunction at the level of the associated homotopy categories:
		$$\xymatrix{(VP_{q}, UW_{q+1}^{\Sigma}, \varphi ):\qslicestablehomotopy \ar[r]& \qslicesymmstablehomotopy}
		$$
	
	Now \cite[proposition 1.3.13]{MR1650134} implies that $VP_{q} ,UW_{q+1}^{\Sigma}$ are both equivalences of
	categories.
	Finally, proposition \ref{prop.2.6.enriched-symmetrization-adjunction}
	together with \cite[proposition 6.4.1]{MR1650134} imply that
	$VP_{q}$ maps cofibre sequences in $\qslicestablehomotopy$ to cofibre sequences in
	$\qslicesymmstablehomotopy$.
	Therefore 
	using proposition 7.1.12 in \cite{MR1650134} we have
	that $VP_{q}$ and $UW_{q+1}^{\Sigma}$ are both exact functors between triangulated categories.
\end{proof}

	Now it is very easy to find the desired lifting 
	for the functor $s_{q}^{\Sigma}:\symmstablehomotopy \rightarrow \symmstablehomotopy$
	(see corollary \ref{cor.3.3.symmetric-fq,s<q,sq}(\ref{cor.3.3.symmetric-fq,s<q,sq.c})) to the
	model category level.
	
\begin{lem}
		\label{lem.3.3.homotopycoherence==>liftings-qslice}
	Fix $q\in \mathbb Z$.
	\begin{enumerate}
		\item \label{lem.3.3.homotopycoherence==>liftings-qslice.a}Let $X$ be an arbitrary 
					$T$-spectrum in $\qconnectedTspectra$.
					Then the following maps in $\qslicesymmTspectra$
						$$\xymatrix{VP_{q}(C_{q}X)\ar[rr]^-{V(P_{q}^{C_{q}X})}&& 
												VC_{q}X && C_{q}^{\Sigma}(VC_{q}X) \ar[ll]_-{C_{q}^{\Sigma ,VC_{q}X}}}
						$$
					induce natural isomorphisms between the functors:
						$$C_{q}^{\Sigma}\circ VC_{q}, VC_{q}, VP_{q}\circ C_{q}:\qconnectedstablehomotopy
							\rightarrow \qslicesymmstablehomotopy$$
						
						$$\xymatrix{& \qconnectedsymmstablehomotopy \ar[dr]^-{C_{q}^{\Sigma}}&\\
												 \qconnectedstablehomotopy
												\ar[rr]^{VC_{q}} \ar[ur]^-{VC_{q}} \ar[dr]_-{C_{q}}&& \qslicesymmstablehomotopy \\
												& \qslicestablehomotopy \ar[ur]_-{VP_{q}}&}
						$$
					Given a $T$-spectrum $X$
						$$\xymatrix{\sigma _{X}:VP_{q}(C_{q}X)\ar[r]^-{\cong}& C_{q}^{\Sigma}(VC_{q}X)}$$ 
					will denote
					the isomorphism in $\qslicesymmstablehomotopy$ corresponding to the natural isomorphism
					between $VP_{q}\circ C_{q}$
					and $C_{q}^{\Sigma}\circ VC_{q}$.
		\item \label{lem.3.3.homotopycoherence==>liftings-qslice.b}Let $X$ be an arbitrary 
					symmetric $T$-spectrum in $\qslicesymmTspectra$.
					Then the following maps in $\qconnectedTspectra$
						$$\xymatrix{W_{q+1}(UW_{q+1}^{\Sigma}X)&& 
												UW_{q+1}^{\Sigma}X \ar[ll]_-{W_{q+1}^{UW_{q+1}^{\Sigma}X}} \ar[rr]^-{U(R_{\Sigma}^{W_{q+1}^{\Sigma}X})}&& 
												UR_{\Sigma}(W_{q+1}^{\Sigma}X) }
						$$
					induce natural isomorphisms between the functors:
						$$W_{q+1}\circ UW_{q+1}^{\Sigma}, UW_{q+1}^{\Sigma}, UR_{\Sigma}\circ W_{q+1}^{\Sigma}:\qslicesymmstablehomotopy
							\rightarrow \qconnectedstablehomotopy$$
						
						$$\xymatrix{& \qslicestablehomotopy \ar[dr]^-{W_{q+1}}&\\
												 \qslicesymmstablehomotopy
												\ar[rr]^{UW_{q+1}^{\Sigma}} \ar[ur]^-{UW_{q+1}^{\Sigma}} \ar[dr]_-{W_{q+1}^{\Sigma}}&& \qconnectedstablehomotopy \\
												& \qconnectedsymmstablehomotopy \ar[ur]_-{UR_{\Sigma}}&}
						$$
					Given a symmetric $T$-spectrum $X$
					$$\xymatrix{\tau _{X}:W_{q+1}(UW_{q+1}^{\Sigma}X)\ar[r]^-{\cong}& UR_{\Sigma}(W_{q+1}^{\Sigma}X)}$$ 
					will denote
					the isomorphism in $\qconnectedstablehomotopy$ corresponding to the natural isomorphism
					between $W_{q+1}\circ UW_{q+1}^{\Sigma}$
					and $UR_{\Sigma}\circ W_{q+1}^{\Sigma}$.
	\end{enumerate}
\end{lem}	
\begin{proof}
	(\ref{lem.3.3.homotopycoherence==>liftings-qslice.a}):  Follows immediately from 
	theorem 1.3.7 in \cite{MR1650134} and the following
	commutative diagram of left Quillen functors:
		$$\xymatrix{\qconnectedTspectra \ar[r]^-{V} \ar[d]_-{id}& \qconnectedsymmTspectra \ar[d]^-{id}\\
								\qsliceTspectra \ar[r]_-{V}& \qslicesymmTspectra}
		$$
		
	(\ref{lem.3.3.homotopycoherence==>liftings-qslice.b}):  Follows immediately from the dual of
	theorem 1.3.7 in \cite{MR1650134} and the following
	commutative diagram of right Quillen functors:
		$$\xymatrix{\qconnectedTspectra & \qconnectedsymmTspectra \ar[l]_-{U}\\
								\qsliceTspectra \ar[u]^-{id}& \qslicesymmTspectra \ar[l]^-{U} \ar[u]_-{id}}
		$$
\end{proof}

\begin{lem}
		\label{lem.3.3.homotopycoherence.b==>liftings-qslice}
		Fix $q \in \mathbb Z$.
		Let $X$ be an arbitrary
		$T$-spectrum, and let $\eta$ be the unit of the adjunction
		(see corollary \ref{cor.3.3.symmetrization-functor-exact-equivalenceSq}):
			$$\xymatrix{(VP_{q}, UW_{q+1}^{\Sigma}, \varphi):\qslicestablehomotopy \ar[r]& \qslicesymmstablehomotopy}
			$$						
		Then we have the following diagram in $\qconnectedstablehomotopy$
		(see lemma \ref{lem.3.3.homotopycoherence==>liftings-qslice}):
			$$\xymatrix{W_{q+1}UW_{q+1}^{\Sigma}VP_{q}C_{q}X 
									\ar[rrr]^-{W_{q+1}UW_{q+1}^{\Sigma}(\sigma _{X})}_-{\cong}&&&
									W_{q+1}UW_{q+1}^{\Sigma}C_{q}^{\Sigma}VC_{q}X
									\ar[dd]^-{\tau _{C_{q}^{\Sigma}VC_{q}X}}_-{\cong}\\
									&&&\\
									W_{q+1}C_{q}X \ar[uu]^-{W_{q+1}(\eta _{C_{q}X})}_-{\cong}
									&&& UR_{\Sigma}W_{q+1}^{\Sigma}C_{q}^{\Sigma}VC_{q}X}
			$$
		where all the maps are isomorphisms in $\qconnectedstablehomotopy$.  This diagram
		induces a natural isomorphism 
		between the following exact functors:
			$$\xymatrix{\qconnectedstablehomotopy \ar@<1ex>[rrr]^-{W_{q+1}C_{q}}  
									\ar@<-1ex>[rrr]_-{UR_{\Sigma}W_{q+1}^{\Sigma}C_{q}^{\Sigma}VC_{q}} &&& 
									\qconnectedstablehomotopy}
			$$
\end{lem}
\begin{proof}
	Follows immediately from lemma
	\ref{lem.3.3.homotopycoherence==>liftings-qslice} and 
	corollary \ref{cor.3.3.symmetrization-functor-exact-equivalenceSq}.
\end{proof}
	
\begin{thm}
		\label{thm.3.3.symmSq-models-symm-sq}
	Fix $q\in \mathbb Z$, and let $X$ be an arbitrary symmetric $T$-spectrum.
	\begin{enumerate}
		\item \label{thm.3.3.symmSq-models-symm-sq.a}The diagram (\ref{diagram.3.2.sq-lifting})
					in theorem \ref{thm.3.2.Sq-models-sq} induces
					the following diagram in 
					$\symmstablehomotopy$:
						\begin{equation}
										\label{diagram.thm.3.3.symmSq-models-symm-sq.a}
							\begin{array}{c}
									\xymatrix{\tilde{s}_{q}X=VQ_{s}(s_{q}UR_{\Sigma}X) \ar[d]_-{VQ_{s}(IQ_{T}J^{s_{q}UR_{\Sigma}X})}^-{\cong}\\
														VQ_{s}(IQ_{T}Js_{q}UR_{\Sigma}X)\\
														VQ_{s}(C_{q}IQ_{T}Js_{q}UR_{\Sigma}X) \ar[u]^-{VQ_{s}(C_{q}^{IQ_{T}Js_{q}UR_{\Sigma}X})}_-{\cong}
														\ar[d]_-{VQ_{s}(W_{q+1}^{C_{q}IQ_{T}Js_{q}UR_{\Sigma}X})}^-{\cong}\\ 
														VQ_{s}(W_{q+1}C_{q}IQ_{T}Js_{q}UR_{\Sigma}X) \\ 
														VQ_{s}(C_{q}W_{q+1}C_{q}IQ_{T}Js_{q}UR_{\Sigma}X) 					 
														\ar[u]^-{VQ_{s}(C_{q}^{W_{q+1}C_{q}IQ_{T}Js_{q}UR_{\Sigma}X})}_-{\cong}\\
														VQ_{s}(C_{q}W_{q+1}C_{q}IQ_{T}Jf_{q}UR_{\Sigma}X) 
														\ar[u]^-{VQ_{s}(C_{q}W_{q+1}C_{q}IQ_{T}J(\pi _{q}^{UR_{\Sigma}X}))}_-{\cong}
														\ar[d]_-{VQ_{s}(C_{q}W_{q+1}C_{q}IQ_{T}J(\theta _{UR_{\Sigma}X}))}^-{\cong}\\
														VQ_{s}(C_{q}W_{q+1}C_{q}IQ_{T}JUR_{\Sigma}X)}
							\end{array}
						\end{equation}
					where all the maps are isomorphisms in $\symmstablehomotopy$.
					This diagram induces a
					natural isomorphism  
					between the following exact functors:
						$$\xymatrix{\symmstablehomotopy \ar@<1ex>[rrrr]^-{\tilde{s}_{q}}  
												\ar@<-1ex>[rrrr]_-{VQ_{s}\circ C_{q}W_{q+1}C_{q}IQ_{T}J\circ UR_{\Sigma}} &&&& \symmstablehomotopy}
						$$
		\item	\label{thm.3.3.symmSq-models-symm-sq.b}Let $\epsilon$
					denote the counit of the adjunction (see corollary \ref{cor.3.3.symmetrization-functor-exact-equivalenceRq}):
						$$\xymatrix{(VC_{q},UR_{\Sigma},\varphi):\qconnectedstablehomotopy \ar[r]&
												\qconnectedsymmstablehomotopy}
						$$
					and let $\delta$ denote the natural isomorphism constructed in lemma \ref{lem.3.3.homotopycoherence.b==>liftings-qslice}.					
					Then we have the following diagram in $\symmstablehomotopy$
					(see lemmas \ref{lem.3.3.homotopycoherence.b==>liftings-qslice} and
					\ref{lem.3.3.homotopycoherence==>liftings}):
						\begin{equation}
										\label{diagram.thm.3.3.symmSq-models-symm-sq.b}
							\begin{array}{c}
									\xymatrix{C_{q}^{\Sigma}W_{q+1}^{\Sigma}C_{q}^{\Sigma}R_{\Sigma}X=s_{q}^{\Sigma}X\\
														C_{q}^{\Sigma}(VC_{q}UR_{\Sigma}W_{q+1}^{\Sigma}C_{q}^{\Sigma}R_{\Sigma}X)
														\ar[u]^-{C_{q}^{\Sigma}(\epsilon _{W_{q+1}^{\Sigma}C_{q}^{\Sigma}R_{\Sigma}X})}_-{\cong}\\
														C_{q}^{\Sigma}VC_{q}UR_{\Sigma}W_{q+1}^{\Sigma}C_{q}^{\Sigma}(VC_{q}UR_{\Sigma}R_{\Sigma}X)
														\ar[u]^-{C_{q}^{\Sigma}VC_{q}UR_{\Sigma}W_{q+1}^{\Sigma}C_{q}^{\Sigma}(\epsilon _{R_{\Sigma}X})}_-{\cong}\\
														C_{q}^{\Sigma}VC_{q}(W_{q+1}C_{q})UR_{\Sigma}R_{\Sigma}X 
														\ar[u]^-{C_{q}^{\Sigma}VC_{q}(\delta _{UR_{\Sigma}R_{\Sigma}X})}_-{\cong}\\
														C_{q}^{\Sigma}VC_{q}W_{q+1}C_{q}(IQ_{T}J UR_{\Sigma}X)
														\ar[u]^-{C_{q}^{\Sigma}VC_{q}W_{q+1}C_{q}(\beta _{X})}_-{\cong}\\
														VQ_{s}C_{q}(W_{q+1}C_{q}IQ_{T}J UR_{\Sigma}X)
														\ar[u]^-{\alpha _{W_{q+1}C_{q}IQ_{T}JUR_{\Sigma}X}}_-{\cong}}
							\end{array}
						\end{equation}
					where all the maps are isomorphisms in $\symmstablehomotopy$.  This diagram induces a natural isomorphism
					between the following exact functors:
						$$\xymatrix{\symmstablehomotopy \ar@<1ex>[rrrr]^-{VQ_{s}\circ C_{q}W_{q+1}C_{q}IQ_{T}J\circ UR_{\Sigma}}  
												\ar@<-1ex>[rrrr]_-{C_{q}^{\Sigma}W_{q+1}^{\Sigma}C_{q}^{\Sigma}R_{\Sigma}=s_{q}^{\Sigma}} &&&& \symmstablehomotopy}
						$$
		\item	\label{thm.3.3.symmSq-models-symm-sq.c}Combining the diagrams (\ref{diagram.thm.3.3.symmSq-models-symm-sq.a}) and 
					(\ref{diagram.thm.3.3.symmSq-models-symm-sq.b}) above
					we get a natural isomorphism 
					between the following exact functors:
						$$\xymatrix{\symmstablehomotopy \ar@<1ex>[r]^-{\tilde{s}_{q}}  
												\ar@<-1ex>[r]_-{s_{q}^{\Sigma}} & \symmstablehomotopy}
						$$
	\end{enumerate}
\end{thm}
\begin{proof}
	It is clear that it suffices to prove only the first two claims.
	
	(\ref{thm.3.3.symmSq-models-symm-sq.a}):  Follows immediately from theorems
	\ref{thm.3.2.Sq-models-sq} and \ref{thm.3.3.symmetrization-functor-exact-equivalence}.
	
	(\ref{thm.3.3.symmSq-models-symm-sq.b}):  Follows immediately from lemmas
	\ref{lem.3.3.homotopycoherence==>liftings} and \ref{lem.3.3.homotopycoherence.b==>liftings-qslice}
	together with corollary \ref{cor.3.3.symmetrization-functor-exact-equivalenceRq}.
\end{proof}

\begin{prop}
		\label{prop.3.3.lifting-map-fq-->sqsymm}
	Fix $q\in \mathbb Z$.  Let $\eta$ denote the unit of the adjuntion
	$(C_{q}^{\Sigma},W_{q+1}^{\Sigma},\varphi ):\qconnectedsymmstablehomotopy \rightarrow \qslicesymmstablehomotopy$
	constructed in corollary \ref{cor.3.3.adjunctions--symmetricRq==>Sq}.  Then the natural transformation
	$\pi _{q}:f_{q}\rightarrow s_{q}$ (see theorem \ref{thm.3.1.slicefiltration}) gets canonically identified,
	through the equivalence of categories $r_{q}C_{q}$, $IQ_{T}Ji_{q}$, $VC_{q}$ and $UR_{\Sigma}$
	constructed in proposition \ref{prop.3.2.Rq-lifts-qSH} and corollary 
	\ref{cor.3.3.symmetrization-functor-exact-equivalenceRq};
	with the following map $\pi _{q}^{\Sigma}:f_{q}^{\Sigma}\rightarrow s_{q}^{\Sigma}$ in $\symmstablehomotopy$:
		$$\xymatrix{C_{q}^{\Sigma}R_{\Sigma}X \ar[rr]^-{C_{q}^{\Sigma}(\eta _{R_{\Sigma}X})}&& 
								C_{q}^{\Sigma}W_{q+1}^{\Sigma}C_{q}^{\Sigma}R_{\Sigma}X}
		$$
\end{prop}
\begin{proof}
	The result follows easily from proposition \ref{prop.3.2.lifting-map-fq-->sq},
	corollaries \ref{cor.3.3.symmetrization-functor-exact-equivalenceRq},
	\ref{cor.3.3.symmetrization-functor-exact-equivalenceSq} and
	theorem \ref{thm.3.3.symmSq-models-symm-sq}.
\end{proof}

	The functor $s_{q}^{\Sigma}$ gives the desired lifting for
	the functor $\tilde{s}_{q}$
	to the model category level, and it will be the main ingredient for the
	study of the multiplicative properties of Voevodsky's slice filtration.
	This completes the program that we started at the beginning of this section.
	
\end{section}
\begin{section}{Multiplicative Properties of the Slice Filtration}
		\label{section.3.4.multiplicativeproperties-slicefiltration}
		
	Our goal in this section is to show that the smash product of spectra is compatible in
	a suitable sense with the slice filtration.  To establish this compatibility in a
	formal way, we will use 
	the model structures constructed in section \ref{section.3.3.symmmodstrslicefilt}.

\begin{lem}
		\label{lem.3.5.spherespectrum-0-connected-0-slice-cofibrant}
	The sphere spectrum $\symmspherespectrum$ is a cofibrant object
	in $\zeroconnectedsymmTspectra$, $\zeroslicesymmTspectra$
	and $\motivicsymmTspectra$.
\end{lem}
\begin{proof}
	By proposition
	\ref{prop.3.3.symmetricq-connected--->q-slice===leftQuilllenfunctor} 
	and theorem \ref{thm.3.3.symmetricconnective-model-structure} we have that
	it is enough to show that
	$\symmspherespectrum$ is  cofibrant in $\zeroconnectedsymmTspectra$.
	
	Now, corollary \ref{cor.3.2.generating-Cqeff-colocalTspectra2} 
	implies that
	$\spherespectrum$ is a $C_{eff}^{0}$-colocal $T$-spectrum in $\motivicTspectra$,
	since $\spherespectrum \in \stablehomotopyeff$.  Then using 
	\cite[theorem 5.1.1(2)]{MR1944041} we have that
	$\spherespectrum$ is a cofibrant object in $\zeroconnectedTspectra$, and this implies
	that $\symmspherespectrum=V(\spherespectrum )$ is also cofibrant in $\zeroconnectedsymmTspectra$,
	since the symmetrization functor
		$$\xymatrix{V:\zeroconnectedTspectra \ar[r]& \zeroconnectedsymmTspectra}
		$$
	is a left Quillen functor.
\end{proof}
	
\begin{lem}
		\label{lem.3.4.cofibrant-preserves-fibrations}
	Fix $p, q \in \mathbb Z$, and let $A$ be a
	symmetric $T$-spectrum.
	\begin{enumerate}
		\item \label{lem.3.4.cofibrant-preserves-fibrations.a}If $A$ is cofibrant in $\pconnectedsymmTspectra$, 
					then the functor $\inthomsymmTspectra (A,-)$ maps fibrations in
					$\pplusqconnectedsymmTspectra$ to fibrations in $\qconnectedsymmTspectra$.
		\item \label{lem.3.4.cofibrant-preserves-fibrations.b}If $A$ is cofibrant in $\motivicsymmTspectra$, 
					then the functor $\inthomsymmTspectra (A,-)$ maps fibrations in
					$\pplusqconnectedsymmTspectra$ to fibrations in $\qconnectedsymmTspectra$.
	\end{enumerate}
\end{lem}
\begin{proof}
	Since $\pconnectedsymmTspectra$, $\qconnectedsymmTspectra$ and $\pplusqconnectedsymmTspectra$
	are all right Bousfield localizations of $\motivicsymmTspectra$, we have that the fibrations in
	all these model structures coincide and also the identity functor			
		$$\xymatrix{id: \pconnectedsymmTspectra \ar[r]& \motivicsymmTspectra}
		$$
	is a left Quillen functor.  Therefore if $A$ is cofibrant in $\pconnectedsymmTspectra$,
	then $A$ is also cofibrant in $\motivicsymmTspectra$.  Hence it suffices to prove
	(\ref{lem.3.4.cofibrant-preserves-fibrations.b}).
	
	So assume that $A$ is cofibrant in $\motivicsymmTspectra$, and let
	$f:X\rightarrow Y$ be an arbitrary fibration in $\pplusqconnectedsymmTspectra$.
	Then using corollary \ref{cor.2.6.smash-cofibrant-leftQuillenfunctor} 
	together with the fact that $A$ is cofibrant
	in $\motivicsymmTspectra$, we get that 
		$$\xymatrix{\inthomsymmTspectra (A,X) \ar[r]^-{f_{\ast}}& \inthomsymmTspectra (A,Y)}$$
	is a fibration in $\motivicsymmTspectra$, or equivalently a
	fibration in $\qconnectedsymmTspectra$.
\end{proof}
	
\begin{lem}
		\label{lem.3.4.generator-kills-trivial-fibres}
	Fix $p, q \in \mathbb Z$, and let $A=\symmgeneratorNRS$ be an arbitrary element
	in $S^{\Sigma}(p)$, i.e. $s-n=p$.  Assume that $F$ is a symmetric
	$T$-spectrum  such that the map
	$F\rightarrow \ast$ is a trivial fibration in $\pplusqconnectedsymmTspectra$.
	Then $\pi :\inthomsymmTspectra (A,F)\rightarrow \ast$ is a trivial fibration in $\qconnectedsymmTspectra$.
\end{lem}
\begin{proof}
	Since $A$ is cofibrant in $\motivicsymmTspectra$,
	it follows directly from lemma \ref{lem.3.4.cofibrant-preserves-fibrations} that
	$\pi$ is a fibration in $\qconnectedsymmTspectra$.  Thus,
	it only remains to show that $\pi$ is a weak equivalence in $\qconnectedsymmTspectra$.
	
	Fix $\symmgeneratorKLM \in C_{eff}^{q,\Sigma}$.
	By construction $\pplusqconnectedsymmTspectra$ is a right Bousfield localization of
	$\motivicsymmTspectra$, therefore $F$ is also fibrant in $\motivicsymmTspectra$.
	Since $A$ is cofibrant and $F$ is fibrant in $\motivicsymmTspectra$, corollary
	\ref{cor.2.6.smash-cofibrant-leftQuillenfunctor}  implies that we have the following natural
	isomorphism of abelian groups:
		$$[\symmgeneratorKLM ,\inthomsymmTspectra (A,F)]_{Spt}^{\Sigma} \cong  
			[\symmgeneratorKLM \wedge A, F]_{Spt}^{\Sigma}$$
	and proposition \ref{prop.2.6.free-symmetric-functors--compatiblewithsmash} implies that:
		\begin{eqnarray*}
			\symmgeneratorKLM \wedge A &= &
			\symmgeneratorKLM \wedge \symmgeneratorNRS\\
			&\cong & \symmgeneratorKLMNRS
		\end{eqnarray*}
	But clearly $\symmgeneratorKLMNRS \in C_{eff}^{p+q,\Sigma}$, and since 
	$F\rightarrow \ast$ is a weak equivalence in $\pplusqconnectedsymmTspectra$,
	we have by proposition \ref{prop.3.3.detecting-symmetric-Cqeff-colocal-equivalences}:
		\begin{eqnarray*}
			0 &\cong & [\symmgeneratorKLMNRS ,F]_{Spt}^{\Sigma}\\
			&\cong & [\symmgeneratorKLM ,\inthomsymmTspectra (A,F)]_{Spt}^{\Sigma}
		\end{eqnarray*}
	Finally, using proposition \ref{prop.3.3.detecting-symmetric-Cqeff-colocal-equivalences} again,
	we get that $\pi$ is a weak equivalence in 
	$\qconnectedsymmTspectra$, as we wanted.
\end{proof}

\begin{lem}
		\label{lem.3.4.cofibrant-kills-trivial-fibres}
	Fix $p, q \in \mathbb Z$, and let $A$ be cofibrant symmetric $T$-spectrum
	in $\pconnectedsymmTspectra$.  Assume that $F$ is a symmetric
	$T$-spectrum  such that the map
	$F\rightarrow \ast$ is a trivial fibration in $\pplusqconnectedsymmTspectra$.
	Then $\pi :\inthomsymmTspectra (A,F)\rightarrow \ast$ is a trivial fibration in $\qconnectedsymmTspectra$.
\end{lem}
\begin{proof}
	Since $A$ is cofibrant in $\pconnectedsymmTspectra$, it follows from
	lemma \ref{lem.3.4.cofibrant-preserves-fibrations}(\ref{lem.3.4.cofibrant-preserves-fibrations.a})
	that $\pi$ is a fibration in $\qconnectedsymmTspectra$.  Thus,
	it only remains to show that $\pi$ is a weak equivalence in $\qconnectedsymmTspectra$.
	
	Fix $\symmgeneratorNRS \in C_{eff}^{q,\Sigma}$.
	Then lemma \ref{lem.3.4.generator-kills-trivial-fibres} together with
	proposition \ref{prop.3.2.functors-between-Rqsymm}(\ref{prop.3.2.functors-between-Rqsymm.a}) imply that
		$$\xymatrix{\inthomsymmTspectra (\symmgeneratorNRS ,F)\ar[r]& \ast}
		$$
	is a trivial fibration in $\pconnectedsymmTspectra$.  Now, since $A$ is cofibrant
	in $\pconnectedsymmTspectra$ which is in particular a simplicial model category, we have
	that the induced map:
		$$\xymatrix{Map\: _{\Sigma}(A,\inthomsymmTspectra (\symmgeneratorNRS ,F))\ar[r]& Map\: _{\Sigma}(A,\ast)=\ast}
		$$
	is a trivial fibration of simplicial sets.
	Finally using the enriched adjunctions of proposition \ref{prop.2.6.enriched-adjunc-symmspectraloops}, 
	this last map gets canonically identified with
		$$\xymatrix{Map\: _{\Sigma}(\symmgeneratorNRS ,\inthomsymmTspectra (A,F))\ar[d]\\ 
								Map\: _{\Sigma}(\symmgeneratorNRS ,\ast)=\ast}
		$$
	Since $\inthomsymmTspectra (A,F)$ is in particular fibrant in $\motivicsymmTspectra$, by definition we have that
	$\pi$ is a $C_{eff}^{q,\Sigma}$-colocal equivalence in $\motivicsymmTspectra$, i.e. a weak equivalence
	in $\qconnectedsymmTspectra$.  This finishes the proof.
\end{proof}

\begin{thm}
		\label{thm.3.4.smash-Quillenbifunctor-RpxRq----->Rqplusp}
	Fix $p,q \in \mathbb Z$.  Then the smash product of symmetric $T$-spectra
		$$\xymatrix{-\wedge -:
									\pconnectedsymmTspectra \times \qconnectedsymmTspectra \ar[r]& \pplusqconnectedsymmTspectra}
		$$
	is a Quillen bifunctor in the sense of Hovey (see definition \ref{def.Quillen-bifunct}).
\end{thm}
\begin{proof}
	By lemma \ref{lem.cond-Quillen-bifunc}, it is enough to prove the following claim:
	
		Given a cofibration $i:A\rightarrow B$ in $\pconnectedsymmTspectra$ and a fibration
		$f:X\rightarrow Y$ in $\pplusqconnectedsymmTspectra$, the induced map
			$$\xymatrix{(i^{\ast},f_{\ast}):\inthomsymmTspectra (B,X)\ar[r]&
									\inthomsymmTspectra (A,X)\times _{\inthomsymmTspectra (A,Y)}\inthomsymmTspectra (B,Y)}
			$$
		is a fibration in $\qconnectedsymmTspectra$ which is trivial if either $i$ or $f$ is a weak
		equivalence.
		
		Since $\pconnectedsymmTspectra$, $\qconnectedsymmTspectra$ and $\pplusqconnectedsymmTspectra$
		are all right Bousfield localizations of $\motivicsymmTspectra$, we have that the fibrations in
		all these model structures coincide and also the identity functor
			\begin{equation}
					\label{diagram.thm.3.4.smash-Quillenbifunctor-RpxRq----->Rqplusp}
				\xymatrix{id: \pconnectedsymmTspectra \ar[r]& \motivicsymmTspectra}
			\end{equation}
		is a left Quillen functor.  Hence it follows that $i$ is cofibration in $\motivicsymmTspectra$
		and $f$ is fibration in $\motivicsymmTspectra$.  Then 
		proposition \ref{prop.2.6.symmTspectra-monoidalmodelcategory} implies that
		$(i^{\ast},f_{\ast})$ is a fibration in $\motivicsymmTspectra$, or equivalently a fibration
		in $\qconnectedsymmTspectra$.
		
		Now assume that $i$ is a trivial cofibration in $\pconnectedsymmTspectra$.  Since the identity functor
		considered in (\ref{diagram.thm.3.4.smash-Quillenbifunctor-RpxRq----->Rqplusp}) 
		above is a left Quillen functor, we have that $i$ is also a trivial
		cofibration in $\motivicsymmTspectra$.  Hence using proposition 
		\ref{prop.2.6.symmTspectra-monoidalmodelcategory} again, we have that
		$(i^{\ast},f_{\ast})$ is in particular a weak equivalence in $\motivicsymmTspectra$.
		Then \cite[proposition 3.1.5]{MR1944041} implies that $(i^{\ast},f_{\ast})$ is also a weak equivalence
		in $\qconnectedsymmTspectra$.
		
		Finally, assume that $f$ is a trivial fibration in $\pplusqconnectedsymmTspectra$.
		Consider the following commutative diagrams
			$$\xymatrix{F \ar[r]^-{\kappa} \ar[d]& \ast \ar[d]& A \ar[r]^-{i} \ar[d]& B \ar[d]\\
									X \ar[r]_-{f}& Y & \ast \ar[r]_-{\iota}& B/A}
			$$
		where the diagram on the left is a pullback in $\pplusqconnectedsymmTspectra$ and
		the diagram on the right is a pushout in $\pconnectedsymmTspectra$.  We already know that
		the map $(i^{\ast},f_{\ast})$ is a fibration in $\qconnectedsymmTspectra$, therefore
		it is clear that $\inthomsymmTspectra (B/A,F)$ is the homotopy fibre of
		$(i^{\ast},f_{\ast})$ in $\qconnectedsymmTspectra$.  On the other hand,
		it is clear that $\kappa$ is a trivial fibration in $\pplusqconnectedsymmTspectra$ and
		$\iota$ is a cofibration in $\pconnectedsymmTspectra$.
		
		By corollary \ref{cor.3.3.symmqconn=>triangcat} 
		we have that the homotopy category associated to $\qconnectedsymmTspectra$ is
		triangulated, hence to check that $(i^{\ast},f_{\ast})$ is a weak equivalence in
		$\qconnectedsymmTspectra$ it is enough to show that the map 
			$$\xymatrix{\pi :\inthomsymmTspectra (B/A,F)\ar[r]& \ast}$$
		is a weak equivalence in $\qconnectedsymmTspectra$.  But this follows immediately from lemma
		\ref{lem.3.4.cofibrant-kills-trivial-fibres}.
\end{proof}

\begin{lem}
		\label{lem.3.4.cofibrant-preserves-cofibrationsLq}
	Fix $p, q, r \in \mathbb Z$.  Let $i:A\rightarrow B$ be a cofibration in 
	$\weightpsymmTspectra$, and let $j:C\rightarrow D$ be a cofibration in
	$\weightqsymmTspectra$.  Then
		$$\xymatrix{B\wedge C \coprod _{A\wedge C}A\wedge D \ar[rr]^-{i\Box j}&& B\wedge D}
		$$
	is a cofibration in $\weightrsymmTspectra$.
\end{lem}
\begin{proof}
	Since $\weightpsymmTspectra$, $\weightqsymmTspectra$ and $\weightrsymmTspectra$
	are all left Bousfield localizations of $\motivicsymmTspectra$, we have that the cofibrations in
	these four model categories coincide.
	
	Then the result follows
	immediately from proposition \ref{prop.2.6.symmTspectra-monoidalmodelcategory}.
\end{proof}

\begin{lem}
		\label{lem.3.4.loops-generator-Lqplusplocal--->Lplocal}
	Fix $p,q \in \mathbb Z$.  Let $A=\symmgeneratorNRS$ be an arbitrary element
	in $S^{\Sigma}(p)$,
	i.e. $s-n=p$, and let $Z$ be an arbitrary $\pplusqleftsymmlocmaps$-local symmetric $T$-spectrum
	in $\motivicsymmTspectra$.
	Then $\inthomsymmTspectra (A,Z)$ is a $\qleftsymmlocmaps$-local symmetric $T$-spectrum
	in $\motivicsymmTspectra$. 
\end{lem}
\begin{proof}
	By proposition \ref{prop.3.3.Lqsymm-local-objects-classification} 
	it is enough to check that the following two conditions hold:
		\begin{enumerate}
			\item	\label{lem.3.4.loops-generator-Lqplusplocal--->Lplocal.a}$\inthomsymmTspectra (A,Z)$ 
						is fibrant in $\motivicsymmTspectra$.
			\item	\label{lem.3.4.loops-generator-Lqplusplocal--->Lplocal.b}For every 
						$\symmgeneratorKLM \in C_{eff}^{q,\Sigma}$
							$$[\symmgeneratorKLM ,\inthomsymmTspectra (A,Z)]_{Spt}^{\Sigma}\cong 0$$
		\end{enumerate}
		
	Since $Z$ is $\pplusqleftsymmlocmaps$-local in $\motivicsymmTspectra$, we have that $Z$ is in particular
	fibrant in $\motivicsymmTspectra$.  Now corollary \ref{cor.2.6.smash-cofibrant-leftQuillenfunctor} 
	together with the fact that $A$ is cofibrant
	in $\motivicsymmTspectra$ imply that $\inthomsymmTspectra (A,Z)$ is fibrant in $\motivicsymmTspectra$.
	This takes care of the first condition.
	
	Fix $\symmgeneratorKLM \in C_{eff}^{q,\Sigma}$.  Since $A$ is cofibrant and $Z$ is fibrant
	in $\motivicsymmTspectra$, it follows from corollary \ref{cor.2.6.smash-cofibrant-leftQuillenfunctor}
	that we have the following natural isomorphism of abelian groups:
	$$[\symmgeneratorKLM ,\inthomsymmTspectra (A,Z)]_{Spt}^{\Sigma} \cong
		[\symmgeneratorKLM \wedge A ,Z]_{Spt}^{\Sigma}
	$$
	Using proposition \ref{prop.2.6.free-symmetric-functors--compatiblewithsmash} 
	we have the following isomorphisms of symmetric $T$-spectra:
		\begin{eqnarray*}
			\symmgeneratorKLM \wedge A &= &
			\symmgeneratorKLM \wedge \symmgeneratorNRS\\
			&\cong & \symmgeneratorKLMNRS
		\end{eqnarray*}
	But clearly $\symmgeneratorKLMNRS \in C_{eff}^{p+q,\Sigma}$.  Since
	$Z$ is a $\pplusqleftsymmlocmaps$-local in $\motivicsymmTspectra$, 
	proposition \ref{prop.3.3.Lqsymm-local-objects-classification}
	implies:
		\begin{eqnarray*}
			0 & \cong & [\symmgeneratorKLMNRS ,Z]_{Spt}^{\Sigma}\\
			& \cong & [\symmgeneratorKLM ,\inthomsymmTspectra (A,Z)]_{Spt}^{\Sigma}
		\end{eqnarray*}
	This finishes the proof.
\end{proof}

\begin{lem}
		\label{lem.3.4.loops-trivcofibrant-Lqplusplocal--->Lplocal}
	Fix $p,q \in \mathbb Z$.  
	Let $A$ be a symmetric
	$T$-spectrum  such that the map
	$\ast \rightarrow A$ is a trivial cofibration in $\weightpsymmTspectra$, 
	and let $Z$ be an arbitrary $\pplusqleftsymmlocmaps$-local symmetric $T$-spectrum
	in $\motivicsymmTspectra$.
	Then $\inthomsymmTspectra (A,Z)$ is a $\qleftsymmlocmaps$-local symmetric $T$-spectrum
	in $\motivicsymmTspectra$.
\end{lem}
\begin{proof}
	By proposition \ref{prop.3.3.Lqsymm-local-objects-classification} 
	it is enough to check that the following two conditions hold:
		\begin{enumerate}
			\item	\label{lem.3.4.loops-trivcofibrant-Lqplusplocal--->Lplocal.a}$\inthomsymmTspectra (A,Z)$ 
						is fibrant in $\motivicsymmTspectra$.
			\item	\label{lem.3.4.loops-trivcofibrant-Lqplusplocal--->Lplocal.b}For every 
						$\symmgeneratorNRS \in C_{eff}^{q,\Sigma}$
							$$[\symmgeneratorNRS ,\inthomsymmTspectra (A,Z)]_{Spt}^{\Sigma}\cong 0$$
		\end{enumerate}
	
	Since $Z$ is $\pplusqleftsymmlocmaps$-local in $\motivicsymmTspectra$, we have that $Z$ is in particular
	fibrant in $\motivicsymmTspectra$.  By construction $\weightpsymmTspectra$ is a left Bousfield localization
	of $\motivicsymmTspectra$, then it follows that $A$ is cofibrant in $\motivicsymmTspectra$.
	Therefore, corollary \ref{cor.2.6.smash-cofibrant-leftQuillenfunctor} 
	implies that $\inthomsymmTspectra (A,Z)$ is fibrant in $\motivicsymmTspectra$.
	This takes care of the first condition.
	
	Fix $\symmgeneratorNRS \in C_{eff}^{q,\Sigma}$.  Using
	lemma \ref{lem.3.4.loops-generator-Lqplusplocal--->Lplocal}
	together with proposition \ref{prop.Lq+1-->Lqsymm},
	we get that the induced map
		$$\xymatrix{\inthomsymmTspectra (\symmgeneratorNRS ,Z)\ar[r]& \ast}
		$$
	is a fibration in $\weightpsymmTspectra$.
	Since $\weightpsymmTspectra$ is a simplicial model category and
	$\ast \rightarrow A$ is a trivial cofibration in $\weightpsymmTspectra$,
	it follows that the following map is a trivial fibration of simplicial sets:
		$$\xymatrix{Map\: _{\Sigma}(A, \inthomsymmTspectra (\symmgeneratorNRS ,Z))\ar[r]&
								Map\: _{\Sigma}(A,\ast)=\ast}
		$$
	Finally using the enriched adjunctions of 
	proposition \ref{prop.2.6.enriched-adjunc-symmspectraloops},
	the map above becomes:
		$$\xymatrix{Map\: _{\Sigma}(\symmgeneratorNRS ,\inthomsymmTspectra (A,Z))\ar[d]\\
								Map\: _{\Sigma}(\symmgeneratorNRS ,\ast)=\ast}
		$$
	which is in particular a weak equivalence of simplicial sets.
	We already know that $\inthomsymmTspectra (A,Z)$ is fibrant in $\motivicsymmTspectra$,
	and we have that $\symmgeneratorNRS$ is cofibrant in $\motivicsymmTspectra$.  Since
	$\motivicsymmTspectra$ is a simplicial model category, we have that
		\begin{eqnarray*}
			0 & \cong & \pi _{0}Map\: _{\Sigma}(\symmgeneratorNRS ,\inthomsymmTspectra (A,Z)) \\
			& \cong & [\symmgeneratorNRS ,\inthomsymmTspectra (A,Z)]_{Spt}^{\Sigma} 
		\end{eqnarray*}
	for every $\symmgeneratorNRS \in C_{eff}^{q,\Sigma}$.
	This finishes the proof.
\end{proof}

\begin{lem}
		\label{lem.3.4.generator-kills-trivial-cofibres}
	Fix $p, q \in \mathbb Z$, and let $A=\symmgeneratorNRS$ be an arbitrary element
	in $S^{\Sigma}(p)$, i.e. $s-n=p$.  Assume that $C$ is a symmetric
	$T$-spectrum  such that the map
	$\ast \rightarrow C$ is a trivial cofibration in $\weightqsymmTspectra$.
	Then $\iota :\ast \rightarrow C\wedge A$ is a trivial cofibration in $\weightpplusqsymmTspectra$.
\end{lem}
\begin{proof}
	Since $A$ is cofibrant in $\motivicsymmTspectra$
	and $\weightpsymmTspectra$ is a left Bousfield localization of $\motivicsymmTspectra$,
	we have that $A$ is also cofibrant in $\weightpsymmTspectra$.  Then
	it follows directly from lemma \ref{lem.3.4.cofibrant-preserves-cofibrationsLq} that
	$\iota$ is a cofibration in $\weightpplusqsymmTspectra$.  Thus,
	it only remains to show that $\iota$ is a weak equivalence in $\weightpplusqsymmTspectra$.
	
	Let $Z$ be an arbitrary $\pplusqleftsymmlocmaps$-local $T$-spectrum in $\motivicsymmTspectra$.
	Then by lemma \ref{lem.3.4.loops-generator-Lqplusplocal--->Lplocal}, we have that
	$\inthomsymmTspectra (A,Z)$ is $\qleftsymmlocmaps$-local in $\motivicsymmTspectra$.
	Now corollary \ref{cor.3.3.detecting-Cq-symmlocal-equivalences} implies that
		$$[C,\inthomsymmTspectra (A,Z)]_{Spt}^{\Sigma} \cong 0
		$$
	
	But $A, C$ are cofibrant in $\motivicsymmTspectra$ and $Z$ is in particular fibrant
	in $\motivicsymmTspectra$, then using corollary \ref{cor.2.6.smash-cofibrant-leftQuillenfunctor} 
	we get the following isomorphism:
		$$[C\wedge A,Z]_{Spt}^{\Sigma}\cong [C,\inthomsymmTspectra (A,Z)]_{Spt}^{\Sigma}\cong 0
		$$
	Hence the induced map
		$$\xymatrix{0\cong [C\wedge A,Z]_{Spt}^{\Sigma} \ar[r]^-{\iota ^{\ast}}& [\ast ,Z]_{Spt}^{\Sigma}\cong 0}
		$$
	is an isomorphism for every $\pplusqleftsymmlocmaps$-local $T$-spectrum $Z$.  Thus,
	using corollary \ref{cor.3.3.detecting-Cq-symmlocal-equivalences} again, 
	we have that $\iota$ is a $\pplusqleftsymmlocmaps$-local equivalence.
	This finishes the proof.
\end{proof}

\begin{lem}
		\label{lem.3.4.trivialcofibrant-kills-trivial-cofibres}
	Fix $p, q \in \mathbb Z$.  Assume that $A, C$ are symmetric
	$T$-spectra  such that 
	$\ast \rightarrow A$ is a trivial cofibration in $\weightpsymmTspectra$, and
	$\ast \rightarrow C$ is a trivial cofibration in $\weightqsymmTspectra$.
	Then $\iota :\ast \rightarrow C\wedge A$ is a trivial cofibration in $\weightpplusqsymmTspectra$.
\end{lem}
\begin{proof}
	Since $A$ is in particular cofibrant in $\weightpsymmTspectra$,
	it follows directly from lemma \ref{lem.3.4.cofibrant-preserves-cofibrationsLq} that
	$\iota$ is a cofibration in $\weightpplusqsymmTspectra$.  Thus,
	it only remains to show that $\iota$ is a weak equivalence in $\weightpplusqsymmTspectra$.
	
	Let $Z$ be an arbitrary $\pplusqleftsymmlocmaps$-local $T$-spectrum in $\motivicsymmTspectra$.
	Then by lemma \ref{lem.3.4.loops-trivcofibrant-Lqplusplocal--->Lplocal}, we have that
	$\inthomsymmTspectra (A,Z)$ is $\qleftsymmlocmaps$-local in $\motivicsymmTspectra$.
	Now corollary \ref{cor.3.3.detecting-Cq-symmlocal-equivalences} implies that
		$$[C,\inthomsymmTspectra (A,Z)]_{Spt}^{\Sigma} \cong 0
		$$
	
	But $A, C$ are in particular cofibrant in $\motivicsymmTspectra$, and $Z$ is in particular fibrant
	in $\motivicsymmTspectra$, then using corollary \ref{cor.2.6.smash-cofibrant-leftQuillenfunctor} 
	we get the following isomorphism:
		$$[C\wedge A,Z]_{Spt}^{\Sigma}\cong [C,\inthomsymmTspectra (A,Z)]_{Spt}^{\Sigma}\cong 0
		$$
	Hence the induced map
		$$\xymatrix{0\cong [C\wedge A,Z]_{Spt}^{\Sigma} \ar[r]^-{\iota ^{\ast}}& [\ast ,Z]_{Spt}^{\Sigma}\cong 0}
		$$
	is an isomorphism for every $\pplusqleftsymmlocmaps$-local $T$-spectrum $Z$.  Thus,
	using corollary \ref{cor.3.3.detecting-Cq-symmlocal-equivalences} again, 
	we have that $\iota$ is a $\pplusqleftsymmlocmaps$-local equivalence.
	This finishes the proof.
\end{proof}

\begin{lem}
	\label{lem.3.4.generator-kills-trivial-fibresSq}
	Fix $p, q \in \mathbb Z$, and let $A=\symmgeneratorNRS$ be an arbitrary element
	in $S^{\Sigma}(p)$, i.e. $s-n=p$.  Assume that $F$ is a symmetric
	$T$-spectrum  such that the map
	$F \rightarrow \ast$ is a trivial fibration in $\pplusqslicesymmTspectra$.
	Then $\pi :\inthomsymmTspectra (A,F) \rightarrow \ast$ is a trivial fibration 
	in $\qslicesymmTspectra$.
\end{lem}
\begin{proof}
	$F$ is fibrant in $\weightpplusqplusonesymmTspectra$, since by construction
	$\pplusqslicesymmTspectra$ is a right Bousfield localization of $\weightpplusqplusonesymmTspectra$.
	Applying lemma \ref{lem.3.4.loops-generator-Lqplusplocal--->Lplocal}, we get that 
	$\inthomsymmTspectra (A,F)$ is fibrant in $\weightqplusonesymmTspectra$; and since
	$\qslicesymmTspectra$ is a right Bousfield localization of $\weightqplusonesymmTspectra$,
	it follows that $\inthomsymmTspectra (A,F)$ is fibrant in $\qslicesymmTspectra$.
	
	By proposition \ref{prop.3.3.detecting-symmetric-Sqeff-colocal-equivalences} 
	it only remains to check that for every $\symmgeneratorKLM \in S^{\Sigma}(q)$,
	i.e. $l-j=q$,
		$$[\symmgeneratorKLM ,\inthomsymmTspectra (A,F)]_{Spt}^{\Sigma}\cong 0
		$$
	Since $A$ is cofibrant in $\motivicsymmTspectra$
	and $F$ is in particular fibrant in $\motivicsymmTspectra$, corollary 
	\ref{cor.2.6.smash-cofibrant-leftQuillenfunctor} implies that we have the following
	natural isomorphism of abelian groups:
		$$[\symmgeneratorKLM ,\inthomsymmTspectra (A,F)]_{Spt}^{\Sigma}\cong
			[\symmgeneratorKLM \wedge A ,F]_{Spt}^{\Sigma}
		$$
	But using proposition \ref{prop.2.6.free-symmetric-functors--compatiblewithsmash} we get:
		\begin{eqnarray*}
			\symmgeneratorKLM \wedge A &= & \symmgeneratorKLM \wedge \symmgeneratorNRS \\
														 &\cong & \symmgeneratorKLMNRS
		\end{eqnarray*}
	and it is clear that $\symmgeneratorKLMNRS \in S^{\Sigma}(p+q)$.
	
	Finally, since $F\rightarrow \ast$ is a trivial
	fibration in $\pplusqslicesymmTspectra$, using proposition \ref{prop.3.3.detecting-symmetric-Sqeff-colocal-equivalences} 
	we get that for every $\symmgeneratorKLM \in S^{\Sigma}(q)$:
		\begin{eqnarray*}
			0 &\cong & [\symmgeneratorKLMNRS ,F]_{Spt}^{\Sigma} \\ 
			&\cong & [\symmgeneratorKLM ,\inthomsymmTspectra (A,F)]_{Spt}^{\Sigma}			
		\end{eqnarray*}
	as we wanted.
\end{proof}

\begin{lem}
	\label{lem.3.4.trivialcofibrant-kills-trivial-fibresSq}
	Fix $p, q \in \mathbb Z$. 
	Assume that $A$ is a cofibrant symmetric $T$-spectrum
	in $\pslicesymmTspectra$, and $F$ is a symmetric
	$T$-spectrum  such that the map
	$F \rightarrow \ast$ is a trivial fibration in $\pplusqslicesymmTspectra$.
	Then $\pi :\inthomsymmTspectra (A,F) \rightarrow \ast$ is a trivial fibration 
	in $\qslicesymmTspectra$.
\end{lem}
\begin{proof}
	$F$ is fibrant in $\weightpplusqplusonesymmTspectra$,
	since by construction $\pplusqslicesymmTspectra$ 
	is a right Bousfield localization of
	$\weightpplusqplusonesymmTspectra$.  Now,
	lemma \ref{lem.3.3.A-symmSqcofibrant====>AsymmLq+1trivially-cofibrant} implies that
	$\ast \rightarrow A$ is a trivial cofibration in $\weightpsymmTspectra$.
	Applying lemma \ref{lem.3.4.loops-trivcofibrant-Lqplusplocal--->Lplocal}, we get that 
	$\inthomsymmTspectra (A,F)$ is fibrant in $\weightqplusonesymmTspectra$; and since
	$\qslicesymmTspectra$ is a right Bousfield localization of $\weightqplusonesymmTspectra$,
	it follows that $\inthomsymmTspectra (A,F)$ is fibrant in $\qslicesymmTspectra$.
	
	Fix $\symmgeneratorNRS \in S^{\Sigma}(q)$, i.e $s-n=q$.  Applying 
	lemma \ref{lem.3.4.generator-kills-trivial-fibresSq} we have that
		$$\xymatrix{\inthomsymmTspectra (\symmgeneratorNRS ,F)\ar[r]& \ast}
		$$
	is a trivial fibration in $\pslicesymmTspectra$ which is in particular a simplicial
	model category.  Therefore the induced map
		$$\xymatrix{Map\: _{\Sigma}(A,\inthomsymmTspectra (\symmgeneratorNRS ,F))\ar[r]& Map\: _{\Sigma}(A,\ast )=\ast}
		$$
	is a trivial fibration of simplicial sets.  Finally using the enriched adjunctions of 
	proposition \ref{prop.2.6.enriched-adjunc-symmspectraloops}, the map above becomes:
		$$\xymatrix{Map\: _{\Sigma}(\symmgeneratorNRS,\inthomsymmTspectra (A ,F))\ar[d]\\ 
								Map\: _{\Sigma}(\symmgeneratorNRS ,\ast)=\ast}
		$$
	which is in particular a weak equivalence of simplicial sets.  We already know
	that $\inthomsymmTspectra (A,F)$ is fibrant in $\weightqplusonesymmTspectra$,
	then by definition it follows that $\pi$ is a $S^{\Sigma}(q)$-colocal equivalence in
	$\weightqplusonesymmTspectra$,
	i.e. a weak equivalence in $\qslicesymmTspectra$.  This finishes the proof.
\end{proof}

\begin{thm}
		\label{thm.3.4.smash-Quillenbifunctor-SpxSq----->Sqplusp}
	Fix $p,q \in \mathbb Z$.  Then the smash product of symmetric $T$-spectra
		$$\xymatrix{-\wedge -:
									\pslicesymmTspectra \times \qslicesymmTspectra \ar[r]& \pplusqslicesymmTspectra}
		$$
	is a Quillen bifunctor in the sense of Hovey (see definition \ref{def.Quillen-bifunct}).
\end{thm}
\begin{proof}
	Since
		$$\xymatrix{-\wedge -:
									\pslicesymmTspectra \times \qslicesymmTspectra \ar[r]& \pplusqslicesymmTspectra}
		$$
	is an adjunction of two variables (see lemma \ref{lem.cond-Quillen-bifunc}),
	it follows that it is enough to prove the following two claims:
		\begin{enumerate}
			\item \label{thm.3.4.smash-Quillenbifunctor-SpxSq----->Sqplusp.a}Let $i:A\rightarrow B$ be a cofibration
						in $\pslicesymmTspectra$, and let $j:C\rightarrow D$ be a cofibration in $\qslicesymmTspectra$.
						Assume that either $i$ or $j$ is trivial.  Then
							$$\xymatrix{B\wedge C \coprod _{A\wedge C} A\wedge D \ar[rr]^-{i\Box j}&& B\wedge D}
							$$
						is a trivial cofibration in $\pplusqslicesymmTspectra$.
			\item \label{thm.3.4.smash-Quillenbifunctor-SpxSq----->Sqplusp.b}Let $i:A\rightarrow B$
						be a cofibration in $\pslicesymmTspectra$, and let $p:X\rightarrow Y$ be a trivial fibration
						in $\pplusqslicesymmTspectra$.  Then
							$$\xymatrix{\inthomsymmTspectra (B,X)\ar[rr]^-{(i^{\ast},p_{\ast})}&&
													\inthomsymmTspectra (B,Y) \times _{\inthomsymmTspectra (A,Y)} \inthomsymmTspectra (A,X)}
							$$
						is a trivial fibration in $\qslicesymmTspectra$.
		\end{enumerate}
		
	(\ref{thm.3.4.smash-Quillenbifunctor-SpxSq----->Sqplusp.a}):  By symmetry, it is enough to consider the
	case where $i$ is a cofibration in $\pslicesymmTspectra$, and $j$ is a trivial cofibration
	in $\qslicesymmTspectra$.  Since $\qslicesymmTspectra$ and $\pslicesymmTspectra$ are 
	right Bousfield localizations of $\weightqplusonesymmTspectra$ and $\weightpplusonesymmTspectra$
	respectively, we have that the identity functor
		$$\xymatrix@R=.5pt{id:\qslicesymmTspectra \ar[r]& \weightqplusonesymmTspectra \\
											 id:\pslicesymmTspectra \ar[r]& \weightpplusonesymmTspectra}
		$$
	is in both cases a left Quillen functor.  This implies in particular that
	$i$ is a cofibration in $\weightpplusonesymmTspectra$ and $j$ is a cofibration
	in $\weightqplusonesymmTspectra$.  Then by 
	lemma \ref{lem.3.4.cofibrant-preserves-cofibrationsLq} we have that
	$i\Box j$ is a cofibration in $\weightpplusqplusonesymmTspectra$.
	
	By construction $\pplusqslicesymmTspectra$ is a right Bousfield localization of
	$\weightpplusqplusonesymmTspectra$, hence the trivial cofibrations in both model
	structures are exactly the same.  Thus, it only remains to show that
	$i\Box j$ is a weak equivalence in $\weightpplusqplusonesymmTspectra$.
	
	Consider the following pushout diagrams in $\motivicTspectra$:
		$$\xymatrix{A \ar[r]^-{i} \ar[d]& B \ar[d] & C \ar[r]^-{j} \ar[d] & D \ar[d]\\
								\ast \ar[r]_-{\iota}& B/A & \ast \ar[r]_-{\kappa} & D/C}
		$$
	By construction
	$\qslicesymmTspectra$ is a right Bousfield localization of $\weightqplusonesymmTspectra$;
	therefore the trivial cofibrations coincide in both model structures.  This implies that
	$j$ and $\kappa$ are both trivial cofibrations in $\weightqplusonesymmTspectra$.
	On the other hand, it is clear that $\iota$ is a cofibration in $\pslicesymmTspectra$.
	Then lemma \ref{lem.3.3.A-symmSqcofibrant====>AsymmLq+1trivially-cofibrant} 
	implies that $\iota$ is a trivial cofibration in $\weightpsymmTspectra$.
	Using lemma \ref{lem.3.4.trivialcofibrant-kills-trivial-cofibres}, 
	we get that the map $\ast \rightarrow (B/A)\wedge (D/C)$ is a trivial
	cofibration in $\weightpplusqplusonesymmTspectra$.
	
	Finally, since $i\Box j$ is a cofibration in $\weightpplusqplusonesymmTspectra$, it follows that
	$(B/A)\wedge (D/C)$ is the homotopy cofibre of $i\Box j$ in $\weightpplusqplusonesymmTspectra$.
	But corollary \ref{cor.3.3.stablehomotopysymmLq==triangulated} 
	implies that the homotopy category associated to $\weightpplusqplusonesymmTspectra$
	is triangulated.  Therefore $i\Box j$ is a trivial cofibration in $\weightpplusqplusonesymmTspectra$,
	since its homotopy cofibre is contractible.
	
	(\ref{thm.3.4.smash-Quillenbifunctor-SpxSq----->Sqplusp.b}):  Using
	(\ref{thm.3.4.smash-Quillenbifunctor-SpxSq----->Sqplusp.a}) above together with the fact that
		$$\xymatrix{-\wedge -:
									\pslicesymmTspectra \times \qslicesymmTspectra \ar[r]& \pplusqslicesymmTspectra}
		$$
	is an adjunction of two variables, we have that $(i^{\ast},p_{\ast})$ is a fibration in $\qslicesymmTspectra$.
	Thus, it only remains to show that $(i^{\ast},p_{\ast})$ is a weak equivalence
	in $\qslicesymmTspectra$.
	
	Consider the following diagrams in $\motivicTspectra$:
		$$\xymatrix{A \ar[r]^-{i} \ar[d]& B \ar[d]& F \ar[r]^-{\kappa} \ar[d]& \ast \ar[d]\\
								\ast \ar[r]_-{\iota}& B/A & X \ar[r]_-{p}& Y}
		$$
	where the diagram on the left is a pushout square and the diagram on the right is a
	pullback square.  It is clear that $\iota$ is a cofibration in $\pslicesymmTspectra$ and that
	$\kappa$ is a trivial fibration in $\pplusqslicesymmTspectra$.
	Then lemma \ref{lem.3.4.trivialcofibrant-kills-trivial-fibresSq} implies that 
	$\inthomsymmTspectra (B/A,F)\rightarrow \ast$ is a trivial fibration in 
	$\qslicesymmTspectra$.
		
	We already know that $(i^{\ast},p_{\ast})$ is a fibration in $\qslicesymmTspectra$,
	therefore $\inthomsymmTspectra (B/A,F)$ is the homotopy fibre of $(i^{\ast},p_{\ast})$
	in $\qslicesymmTspectra$.
	Finally, by corollary \ref{cor.3.3.symmqslice=>triangcat} 
	we have that the homotopy category associated to $\qslicesymmTspectra$
	is triangulated.  Therefore it follows that $(i^{\ast},p_{\ast})$ is a trivial
	fibration in $\qslicesymmTspectra$, since its homotopy fibre is contractible.
\end{proof}

\end{section}
\begin{section}{Further Multiplicative Properties of the Slice Filtration}
		\label{section.3.5furthermultiplicativeproperties}

	In this section $A$ will always denote a cofibrant ring spectrum with unit in $\motivicsymmTspectra$.
	The goal in this section is to use the motivic model structure $\motivicAmod$ for the category of
	$A$-modules (see section \ref{section.2.8Modules-Algebras}) together with the
	model structures for the category of symmetric $T$-spectra
	constructed in section \ref{section.3.3.symmmodstrslicefilt}
	(which provide a lifting of the slice filtration to the model category level),
	in order to get an analogue of the slice filtration for the category of $A$-modules.
	The main results of this section guarantee that under suitable conditions, the ($q-1$)-connective
	cover $f_{q}^{\Sigma}(M)$, $s_{<q}^{\Sigma}(M)$ and the $q$-slice	$s_{q}^{\Sigma}(M)$
	of an arbitrary $A$-module $M$, inherit a natural structure of $A$-module; and that
	the unit map $u:\symmspherespectrum \rightarrow A$ satisfying some natural additional conditions,
	induces for every symmetric $T$-spectrum $X$
	a natural structure of $A$-module on its $q$-slice $s_{q}^{\Sigma}(X)$.
	
	Let $\Amodstablehomotopy$ denote the homotopy category associated
	to $\motivicAmod$.
	We call $\Amodstablehomotopy$ the \emph{motivic stable homotopy category of $A$-modules}.
	We will denote by $[-,-]_{m}$ the set of maps between two
	objects in $\Amodstablehomotopy$.

\begin{defi}
		\label{def.3.5.cofibrantstable-replacementfunctor}
	Let $Q_{m}$ denote
	a cofibrant replacement functor in $\motivicAmod$; such that for every
	$A$-module $M$, the natural map
		$$\xymatrix{Q_{m}M\ar[r]^-{Q_{m}^{M}}& M}
		$$
	is a trivial fibration in $\motivicAmod$.
\end{defi}
	
\begin{defi}
		\label{def.3.5.fibrantstable-replacementfunctor}
	Let $R_{m}$ denote
	a fibrant replacement functor in $\motivicAmod$; such that for every
	$A$-module $M$, the natural map
		$$\xymatrix{M\ar[r]^-{R_{m}^{M}}& R_{m}M}
		$$
	is a trivial cofibration in $\motivicAmod$.
\end{defi}
	
\begin{prop}
		\label{prop.3.5.symmstable-homotopy==triangulated-category}
	The motivic stable homotopy category of $A$-modules $\Amodstablehomotopy$
	has a structure of triangulated category defined as follows: 
	\begin{enumerate}
		\item The suspension $\Sigma _{T}^{1,0}$ functor is given by
					$$\xymatrix@R=.5pt{-\wedge S^{1}:\Amodstablehomotopy \ar[r]& \Amodstablehomotopy \\
															M \ar@{|->}[r]& Q_{m}M\wedge S^{1}}
					$$
		\item The distinguished triangles are
					isomorphic to triangles of the form
						$$\xymatrix{M \ar[r]^-{i}& N \ar[r]^-{j}& O \ar[r]^-{k}& \Sigma_{T}^{1,0}M}
						$$
					where $i$ is a cofibration in $\motivicAmod$, and $O$ is the homotopy cofibre of $i$.
	\end{enumerate}
\end{prop}
\begin{proof}
	By proposition \ref{prop.2.8.motivicA-modules=>M-symmSpt-modelcat}(\ref{prop.2.8.motivicA-modules=>M-symmSpt-modelcat.b})
	we have that
	$\motivicAmod$ is a pointed simplicial model category, and
	theorem \ref{thm.2.8.smashT-Quillenequiv} implies that the adjunction:
		$$\xymatrix{(-\wedge S^{1},\Omega _{S^{1}},\varphi):\motivicAmod \ar[r]& \motivicAmod}
		$$
	is a Quillen equivalence.  The result now follows from the work of
	Quillen in \cite[sections I.2 and I.3]{MR0223432} and the work of 
	Hovey in \cite[chapters VI and VII]{MR1650134}
	(see \cite[proposition 7.1.6]{MR1650134}).
\end{proof}

\begin{thm}
		\label{thm.3.5.freemodulefunctor-leftQuillen}
	The adjunction
		$$\xymatrix{(A\wedge -,U,\varphi):\motivicsymmTspectra \ar[r]& \motivicAmod}
		$$
	defined in proposition \ref{prop.2.8.Free-module-functor}, induces
	an adjunction
		$$\xymatrix{(A\wedge Q_{\Sigma}-,UR_{m},\varphi):\symmstablehomotopy \ar[r]& \Amodstablehomotopy}
		$$
	of exact funtors between triangulated categories.
\end{thm}
\begin{proof}
	The proof is exactly the same as in
	theorem \ref{thm.3.3.symmetrization-functor-exact-equivalence}.  We leave
	the details to the reader.
\end{proof}

\begin{lem}
		\label{lem.3.5.compact-respects-coproducts-symmspectra}
	Let $X\in \pointedmotivic$ be a pointed simplicial presheaf which
	is compact in the sense of Jardine (see definition \ref{def.2.3.compactness}), and let
	$F_{n}^{\Sigma}(X)$ be the symmetric $T$-spectrum constructed in definition \ref{def.2.6.Fnsigma-Evn-adjunction}.
	Consider an arbitrary collection of $A$-modules $\{ M_{i}\}_{i\in I}$ indexed by a set $I$.
	Then 
		$$[A\wedge F_{n}^{\Sigma}(X),\coprod _{i\in I}M_{i}]_{m}\cong \coprod _{i\in I}\: 
			[A\wedge F_{n}^{\Sigma}(X),M_{i}]_{m}
		$$
\end{lem}
\begin{proof}
	The proof is exactly the same as in
	lemma \ref{lem.3.3.compact-respects-coproducts-symmspectra}.  We leave
	the details to the reader.
\end{proof}

\begin{prop}
		\label{prop.3.5.symmstablehomotopy=>compactly-generated}
	The motivic stable homotopy category of $A$-modules $\Amodstablehomotopy$ is
	a \emph{compactly generated} triangulated category
	in the sense of Neeman (see \cite[definition 1.7]{MR1308405}).  The set
	of compact generators is given by
	(see definition \ref{def.2.6.Fnsigma-Evn-adjunction}):				
		$$\Amodcompactgenerators
		$$
	i.e. the smallest triangulated subcategory of $\Amodstablehomotopy$ 
	closed under small  coproducts and containing
	all the objects in $C^{m}$ coincides with $\Amodstablehomotopy$.
\end{prop}
\begin{proof}
	The proof is exactly the same as in 
	proposition \ref{prop.3.3.symmstablehomotopy=>compactly-generated}.
	We leave the details to the reader.
\end{proof}

\begin{cor}
		\label{cor.3.5.detecting-isos-in-symmSH}
	Let $f:M\rightarrow N$ be a map in $\Amodstablehomotopy$.
	Then $f$ is an isomorphism if and only if
	$f$ induces an isomorphism of abelian groups:
		$$\xymatrix{[A\wedge \symmgeneratorNRS , M]_{m} \ar[r]^-{f_{\ast}}& [A\wedge \symmgeneratorNRS , N]_{m}}
		$$
	for every $A\wedge \symmgeneratorNRS \in C^{m}$.
\end{cor}
\begin{proof}
	The proof is exactly the same as in 
	corollary \ref{cor.3.3.detecting-isos-in-symmSH}.
	We leave the details to the reader.
\end{proof}

	In the rest of this section some results will be just stated without proof.
	In every case, the proof is exactly the same as the one given in section \ref{section.3.3.symmmodstrslicefilt},
	taking into consideration 
	all that has been proved so far in this section together with proposition
	\ref{prop.2.8.enriched-freemodule-adjunction}, the cellularity for the motivic model category of
	$A$-modules (see theorem \ref{thm.2.8.motivicAmods=>cellular}), and the fact
	that the generators $A\wedge \symmgeneratorNRS \in C^{m}$ are all cofibrant in $\motivicAmod$
	(this follows immediately from theorem \ref{thm.2.8.motivicsymmetricA-modules}).

\begin{thm}
		\label{thm.3.5.symmetricconnective-model-structure}
	Fix $q\in \mathbb Z$.  Consider the following set of objects in $\motivicAmod$
	(see theorem \ref{thm.3.2.connective-model-structure}):
		$$C_{eff}^{q,m}=\Amodqeffsymmcompactgenerators
		$$
	Then the right Bousfield localization of $\motivicAmod$
	with respect to the class of $C_{eff}^{q,m}$-colocal equivalences exists
	(see definitions \ref{def1.2.cellularization} and \ref{def1.2.1.rightBousloc}).
	This model structure 
	will be called  \emph{($q-1$)-connected motivic stable}, and
	the category of $A$-modules equipped with the ($q-1$)-connected motivic
	stable model structure will be denoted by
	$\qconnectedAmod$.  Furthermore
	$\qconnectedAmod$ is a right proper and simplicial model
	category.
	The homotopy category associated to $\qconnectedAmod$
	will be denoted by $\qconnectedAmodstablehomotopy$.
\end{thm}

\begin{rmk}
		\label{rmk.3.5.usualargumentsnotgood}
	Notice that we can not use the adjuntion 
		$$(A\wedge -,U,\varphi ):\qconnectedsymmTspectra \rightarrow \qconnectedAmod$$
	for the construction of
	$\qconnectedAmod$, since we do not know if the model structure for $\qconnectedsymmTspectra$
	is cofibrantly generated.
\end{rmk}

\begin{defi}
		\label{def.3.5.Cqsigma-cofibrant-replacement}
	Fix $q\in \mathbb Z$.  Let $C_{q}^{m}$ denote a cofibrant replacement functor in 
	$\qconnectedAmod$; such that for every $A$-module $M$, the natural map
		$$\xymatrix{C_{q}^{m}M \ar[r]^-{C_{q}^{m,M}}& M}
		$$ 
	is a trivial fibration in $\qconnectedAmod$, and
	$C_{q}^{m}M$ is always $C_{eff}^{q,m}$-colocal in $\motivicAmod$.
\end{defi}

\begin{prop}
		\label{prop.3.5.Rsigma-fibrant-replacement-all-Rq}
	Fix $q\in \mathbb Z$.  Then $R_{m}$ is also a fibrant
	replacement functor in $\qconnectedAmod$
	(see definition \ref{def.3.5.fibrantstable-replacementfunctor}),
	and for every $A$-module $M$ the natural map
		$$\xymatrix{M\ar[r]^-{R_{m}^{M}}& R_{m}M}
		$$
	is a trivial cofibration in $\qconnectedAmod$.
\end{prop}

\begin{prop}
		\label{prop.3.5.Suspension=>qconnsymm-Quillen-equiv}
	Fix $q\in \mathbb Z$.  Then the adjunction
		$$\xymatrix{(-\wedge S^{1},\Omega _{S^{1}},\varphi):\qconnectedAmod \ar[r]& 
								\qconnectedAmod}
		$$
	is a Quillen equivalence, and $\qconnectedAmodstablehomotopy$ has the structure
	of a triangulated category.
\end{prop}

\begin{prop}
		\label{prop.3.5.symmcofibrant-replacement=>triangulatedfunctor}
	For every $q\in \mathbb Z$, we have the following Quillen adjunction
		$$\xymatrix{(id,id,\varphi ):\qconnectedAmod \ar[r] & \motivicAmod}
		$$
	which induces an adjunction
		$$\xymatrix{(C_{q}^{m}, R_{m}, \varphi) :\qconnectedAmodstablehomotopy \ar[r]& \Amodstablehomotopy}
		$$
	between exact functors of triangulated categories.
\end{prop}

\begin{thm}
		\label{thm.3.5.symmetrization-qconnected-Quillen-equiv}
	Fix $q \in \mathbb Z$.  Then the adjunction (see theorem \ref{thm.3.2.connective-model-structure})
		$$\xymatrix{(A\wedge -,U,\varphi):\qconnectedsymmTspectra \ar[r]& \qconnectedAmod}
		$$
	given by the free $A$-module and the forgetful functors is a
	Quillen adjunction, and it induces
	an adjunction
		$$\xymatrix{(A\wedge C_{q}^{\Sigma}-,UR_{m},\varphi):\qconnectedsymmstablehomotopy \ar[r]& \qconnectedAmod}
		$$
	of exact funtors between triangulated categories.
\end{thm}
	
\begin{prop}
		\label{prop.3.5.functors-between-Rqsymm}
	Fix $q\in \mathbb Z$.  We have the following
	commutative diagram of left Quillen functors:
		\begin{equation}
				\label{diagram3.5.liftslicefiltrationsymm.a}
					\begin{array}{c}
						\xymatrix{\qplusoneconnectedsymmTspectra \ar[rr]^-{id} \ar[dr]_-{id} \ar[dd]_-{A\wedge -}&& 
														\qconnectedsymmTspectra \ar[dl]^-{id} \ar[dd]^-{A\wedge -} \\
														&\motivicsymmTspectra  \ar[dd]^(.3){A\wedge -}&\\
														\qplusoneconnectedAmod \ar[rr]^(.3){id} \ar[dr]_-{id}&& 
														\qconnectedAmod \ar[dl]^-{id}\\
														&\motivicAmod  &}
					\end{array}
		\end{equation}
	and the following associated commutative diagrams of homotopy categories:
		\begin{equation}
	 		\label{diagram3.5.liftslicefiltrationsymm.b}
	 			\begin{array}{c}
	 				\xymatrix{& \qconnectedAmodstablehomotopy \ar[dr]^-{C_{q}^{m}}&\\
												 \qconnectedsymmstablehomotopy
												\ar[rr]^{A\wedge C_{q}^{\Sigma}-} \ar[ur]^-{A\wedge C_{q}^{\Sigma}-} \ar[dr]_-{C_{q}^{\Sigma}}&& \Amodstablehomotopy \\
											& \symmstablehomotopy \ar[ur]_-{A\wedge Q_{\Sigma}-}&}
				\end{array}
		\end{equation}
		
		\begin{equation}
	 		\label{diagram3.5.liftslicefiltrationsymm.c}
	 			\begin{array}{c}
	 				\xymatrix{& \symmstablehomotopy \ar[dr]^-{R_{\Sigma}}&\\
												 \Amodstablehomotopy
												\ar[rr]^{UR_{m}} \ar[ur]^-{UR_{m}} \ar[dr]_-{R_{m}}&& \qconnectedsymmstablehomotopy \\
												& \qconnectedAmodstablehomotopy \ar[ur]_-{UR_{m}}&}
				\end{array}
		\end{equation}
\end{prop}

\begin{thm}
		\label{thm.3.5.liftingslicefilt-modelcatssymm}
	We have the following commutative diagram
	of left Quillen functors:
		\begin{equation}
					\label{diagram.3.5.liftingVoevodskyslicefiltsymm.a}
			\begin{array}{c}
				\xymatrix@C=0.3pc{& \vdots \ar[d]_-{id} & \vdots \ar[d]_-{id} & \\ 
									& \qplusoneconnectedsymmTspectra \ar[d]_-{id} \ar[r]_-{A\wedge -} \ar[dl]_-{id}& 
									\qplusoneconnectedAmod \ar[d]_-{id} \ar[dr]^-{id} & \\ 
									\motivicsymmTspectra & \qconnectedsymmTspectra \ar[d]_-{id} \ar[r]_-{A\wedge -} \ar[l]_-{id}& 
									\qconnectedAmod \ar[d]_-{id} \ar[r]^-{id} & \motivicAmod \\
									& \qminusoneconnectedsymmTspectra \ar[d]_-{id} \ar[r]_-{A\wedge -} \ar[ul]^-{id}& 
									\qminusoneconnectedAmod \ar[d]_-{id} \ar[ur]_-{id} & \\ 
									& \vdots & \vdots & }
			\end{array}
		\end{equation}
	and the following associated commutative diagrams of homotopy categories:	
		\begin{equation}
					\label{diagram.3.5.liftingVoevodskyslicefiltsymm.b}
			\begin{array}{c}
				\xymatrix@C=.3pc{&&& \vdots \ar[d] && \vdots \ar[d] &&& \\
									&&&
				 					R_{C_{eff}^{q+1}}\symmstablehomotopy \ar[d]_-{C_{q+1}^{\Sigma}} 
				 					\ar[dlll]|-{C_{q+1}^{\Sigma}}  \ar[rr]_-{A\wedge Q_{\Sigma}-} &&
									R_{C_{eff}^{q+1}}\Amodstablehomotopy \ar[d]_-{C_{q+1}^{m}} 
									\ar[drrr]|-{C_{q+1}^{m}}  &&& \\ 
									\symmstablehomotopy &&&
									\qconnectedsymmstablehomotopy \ar[d]_-{C_{q}^{\Sigma}} 
									\ar[lll]|-{C_{q}^{\Sigma}}  \ar[rr]_-{A\wedge Q_{\Sigma}-} &&  
									\qconnectedAmodstablehomotopy \ar[d]_-{C_{q}^{m}} 
									\ar[rrr]|-{C_{q}^{m}}  &&& \Amodstablehomotopy \\
									&&& R_{C_{eff}^{q-1}}\symmstablehomotopy \ar[d] \ar[ulll]|-{C_{q-1}^{\Sigma}} 
									\ar[rr]_-{A\wedge Q_{\Sigma}-} &&
									R_{C_{eff}^{q-1}}\Amodstablehomotopy \ar[d] \ar[urrr]|-{C_{q-1}^{m}} &&& \\ 
									&&& \vdots  &&
									\vdots &&& }
			\end{array}
		\end{equation}
		
	\begin{equation}
					\label{diagram.3.5.liftingVoevodskyslicefiltsymm.c}
			\begin{array}{c}
				\xymatrix@C=.3pc{&&& \vdots  && \vdots  &&& \\
									&&&
				 					R_{C_{eff}^{q+1}}\symmstablehomotopy  \ar[u]
				 					&& \ar[ll]^-{UR_{m}} 
									R_{C_{eff}^{q+1}}\Amodstablehomotopy  \ar[u]
									&&&  \\ 
									\symmstablehomotopy  \ar[urrr]|-{R_{\Sigma}} \ar[rrr]|-{R_{\Sigma}} \ar[drrr]|-{R_{\Sigma}}&&&
									\qconnectedsymmstablehomotopy \ar[u]^-{R_{\Sigma}}  
									&& \ar[ll]^-{UR_{m}}  
									\qconnectedAmodstablehomotopy  \ar[u]^-{R_{m}}   
								  &&& \Amodstablehomotopy  \ar[ulll]|-{R_{m}} \ar[lll]|-{R_{m}} \ar[dlll]|-{R_{m}} \\
									&&& R_{C_{eff}^{q-1}}\symmstablehomotopy  \ar[u]^-{R_{\Sigma}}  
									&&
									R_{C_{eff}^{q-1}}\Amodstablehomotopy  \ar[u]^-{R_{m}}  \ar[ll]^-{UR_{m}} &&& \\ 
									&&& \vdots  \ar[u] &&
									\vdots \ar[u] &&& }
			\end{array}
		\end{equation}
\end{thm}

\begin{thm}
		\label{thm.3.5.smash-Quillenbifunctor-RpxRq----->Rqplusp-coeffs}
	Fix $p,q \in \mathbb Z$.  Let $A$ be a commutative cofibrant ring spectrum in $\motivicsymmTspectra$.
	Then the symmetric monoidal structure for the category of $A$-modules, induces
	the following Quillen bifunctor in the sense of Hovey (see definition \ref{def.Quillen-bifunct}).
		$$\xymatrix{-\wedge _{A}-:
									\pconnectedAmod \times \qconnectedAmod \ar[r]& \pplusqconnectedAmod}
		$$
\end{thm}
\begin{proof}
	The proof is similar to the one given for theorem \ref{thm.3.4.smash-Quillenbifunctor-RpxRq----->Rqplusp}.
	We leave the details to the reader.
\end{proof}

	If the ring $A$ is not commutative, then it needs to satisfy some
	additional conditions
	in order to get a weaker version of
	the previous result (see 
	theorem \ref{thm.3.5.notcommutativesmash-Quillenbifunctor-RpxRq----->Rqplusp-coeffs}).
	
\begin{lem}
		\label{lemma.3.5.changecoefficients-Quillenfunctor}
	Fix $q\in \mathbb Z$.  Let $f:A\rightarrow A'$ be a map 
	between cofibrant ring spectra in $\motivicsymmTspectra$,
	which is compatible with the ring structures.
		\begin{enumerate}
			\item \label{lemma.3.5.changecoefficients-Quillenfunctor.a}  Then the adjunction:
								$$(A'\wedge _{A}-, U,\varphi ):\qconnectedAmod \rightarrow \qconnectedAAmod
								$$
							is a Quillen adjunction.  
			\item \label{lemma.3.5.changecoefficients-Quillenfunctor.b}  Furthermore, a map 
							$w:M\rightarrow M'$ in $\qconnectedAAmod$ is a weak equivalence if
							and only if $Uw$ is a weak equivalence in $\qconnectedAmod$.
		\end{enumerate}
\end{lem}
\begin{proof}
	(\ref{lemma.3.5.changecoefficients-Quillenfunctor.a}):  Lemma \ref{lemma.2.5.changecoefficients-Quillenfunctor} 
	implies that $U:\motivicAAmod \rightarrow \motivicAmod$ is a right Quillen functor.
	Consider the following commutative diagram of right Quillen functors:
		$$\xymatrix{\motivicAAmod \ar[r]^-{U} \ar[d]_-{id} & \motivicAmod \ar[d]^-{id}\\
								\qconnectedAAmod \ar@{-->}[r]_-{U}& \qconnectedAmod}
		$$
	then the universal property of right Bousfield localizations together with proposition
	\ref{prop.2.8.enriched-freemodule-adjunction2} imply that the dotted arrow $U$ is a right Quillen functor.
	
	(\ref{lemma.3.5.changecoefficients-Quillenfunctor.b}):  Let $R_{m}$, $R_{m'}$ denote fibrant replacement functors
	in $\motivicAmod$ and $\motivicAAmod$ respectively, and let $N$ be an arbitrary $A'$-module.  We have the following
	commutative diagram in $\motivicAmod$:
		$$\xymatrix{N  \ar[d]_-{R_{m'}^{N}} \ar[r]^-{R_{m}^{N}}& 
								R_{m}N \ar[d]^-{R_{m}(R_{m'}^{N})}\\
								R_{m'}N \ar[r]_-{R_{m}^{R_{m'}N}} & R_{m}R_{m'}N }
		$$
	lemma \ref{lemma.2.5.changecoefficients-Quillenfunctor} 
	implies that all the maps in the diagram above are weak equivalences in $\motivicAmod$.  
	  	
	Now fix $\generatorNRS \in C_{eff}^{q,\Sigma}$.  Using
	the naturality of the diagram above together with 
	proposition \ref{prop.2.8.enriched-freemodule-adjunction2}, we get
	the following commutative diagram of simplicial sets:	
		$$\xymatrix@C=-6pc{Map_{A'\text{-}\modules}(A'\wedge \generatorNRS ,R_{m'}M) \ar[dr]^-{(R_{m'}w)_{\ast}} \ar[dd]_-{\cong}&\\  
								& Map_{A'\text{-}\modules}(A'\wedge \generatorNRS ,R_{m'}M')\ar[dd]^-{\cong}\\
								Map_{A\text{-}\modules}(A\wedge \generatorNRS ,UR_{m'}M) \ar[dr]^-{(UR_{m'}w)_{\ast}}
								\ar[dd]_-{(R_{m}^{UR_{m'}M})_{\ast}}&\\  
								& Map_{A\text{-}\modules}(A\wedge \generatorNRS ,UR_{m'}M')\ar[dd]^-{(R_{m}^{UR_{m'}M'})_{\ast}}\\
								Map_{A\text{-}\modules}(A\wedge \generatorNRS ,R_{m}UR_{m'}M) \ar[dr]^-{(R_{m}UR_{m'}w)_{\ast}}
								\ar@{<-}[dd]_-{(UR_{m}(R_{m'}^{M}))_{\ast}}&\\  
								& Map_{A\text{-}\modules}(A\wedge \generatorNRS ,R_{m}UR_{m'}M')
								\ar@{<-}[dd]^-{(UR_{m}(R_{m'}^{M'}))_{\ast}}\\
								Map_{A\text{-}\modules}(A\wedge \generatorNRS ,UR_{m}M) \ar[dr]^-{(UR_{m}w)_{\ast}}&\\  
								& Map_{A\text{-}\modules}(A\wedge \generatorNRS ,UR_{m}M')}
		$$
	where the top vertical arrows are isomorphisms of simplicial sets.
	But $\motivicAmod$, $\motivicAAmod$ are simplicial
	model categories (see proposition 
	\ref{prop.2.8.motivicA-modules=>M-symmSpt-modelcat}(\ref{prop.2.8.motivicA-modules=>M-symmSpt-modelcat.b})) 
	and the natural maps $R_{m}^{UR_{m'}M}$, $UR_{m}(R_{m'}^{M})$, $R_{m}^{UR_{m'}M'}$ and $UR_{m}(R_{m'}^{M'})$ 
	are all weak equivalences between fibrant objects, thus by Ken Brown's lemma (see lemma \ref{lemma1.1.KenBrown2})
	all the vertical arrows are weak equivalences
	of simplicial sets. 
	
	Therefore, the top row is a weak equivalence of simplicial sets
	if and only if the bottom row is a weak equivalence of simplicial sets.
	This proves the claim.	
\end{proof}

\begin{prop}
		\label{prop.3.5.changecoefficients-Quillenequiv}  
	Fix $q\in \mathbb Z$.  Let $f:A\rightarrow A'$ be
	a map between cofibrant ring spectra in $\motivicsymmTspectra$, which is
	compatible with the ring structures.
	Assume that $f$ is a weak equivalence in $\motivicsymmTspectra$.  Then the adjunction
		$$(A'\wedge _{A}-, U,\varphi ):\qconnectedAmod \rightarrow \qconnectedAAmod
		$$
	is a Quillen equivalence.
\end{prop}
\begin{proof}
	We have shown in lemma 
	\ref{lemma.3.5.changecoefficients-Quillenfunctor}(\ref{lemma.3.5.changecoefficients-Quillenfunctor.a}) that
		$$(A'\wedge _{A} -, U,\varphi ):\qconnectedAmod \rightarrow \qconnectedAAmod$$
	is a Quillen adjunction.
		
	Now 
	let $\eta$, $\epsilon$ denote the unit and counit of the
	adjunction $(A' \wedge _{A}-, U, \varphi)$.
	By corollary 1.3.16(c) in \cite{MR1650134}, it suffices to check
	that the following conditions hold:
		\begin{enumerate}
			\item \label{prop.3.5.changecoefficients-Quillenequiv.a}  For every cofibrant
						$A$-module $M$ in $\qconnectedAmod$, the following composition
							$$\xymatrix{M\cong A\wedge _{A}M \ar[rr]^-{\eta _{M}=f\wedge _{A}id} && A'\wedge _{A}M
								\ar[rr]^-{R_{m'}^{A'\wedge _{A}M}}&& R_{m'}(A'\wedge _{A}M)}
							$$
						is a weak equivalence in $\qconnectedAmod$, where $R_{m'}$ denotes a fibrant replacement functor
						in $\motivicAAmod$ (see proposition \ref{prop.3.5.Rsigma-fibrant-replacement-all-Rq}).
			\item \label{prop.3.5.changecoefficients-Quillenequiv.b}  $U$ reflects
						weak equivalences between fibrant objects in $\qconnectedAAmod$.
		\end{enumerate}
		
	(\ref{prop.3.5.changecoefficients-Quillenequiv.a}):  Since $id:\qconnectedAmod \rightarrow \motivicAmod$
	is a left Quillen functor, we have that $M$ is also cofibrant in $\motivicAmod$.  Hence, proposition
	\ref{prop.2.8.invarianceofcoefficients} implies that $R_{m'}^{A'\wedge _{A}M}\circ \eta _{M}$ 
	is a weak equivalence in $\motivicAmod$.  Finally, by \cite[proposition 3.1.5]{MR1944041}
	we have that $R_{m'}^{A'\wedge _{A}M}\circ \eta _{M}$ is a weak equivalence in $\qconnectedAmod$,
	as we wanted.

	(\ref{prop.3.5.changecoefficients-Quillenequiv.b}): This follows immediately
	from lemma \ref{lemma.3.5.changecoefficients-Quillenfunctor}(\ref{lemma.3.5.changecoefficients-Quillenfunctor.b}).
\end{proof}

\begin{prop}
		\label{prop.3.5.Acolocal=>colocalcofibrations=>cofibrations}
	Let $A$ be a cofibrant ring spectrum in $\motivicsymmTspectra$, which is also
	$C_{eff}^{0,\Sigma}$-colocal (equivalently cofibrant in $\zeroconnectedsymmTspectra$).  
	Then for every $q\in \mathbb Z$, and for every cofibration $f:M\rightarrow N$
	in $\qconnectedAmod$ we have
	that $f$ is also a cofibration in $\qconnectedsymmTspectra$.
\end{prop}
\begin{proof}
	Let (see theorem \ref{thm.2.8.motivicAmods=>cellular})
		\begin{eqnarray*}
			\overline{\Lambda (K)} = J_{A\text{-}\modules}\; \cup \; \{ A\wedge \symmgeneratorNRS \otimes \partial \Delta ^{k}
								\rightarrow  \\ 
								A\wedge \symmgeneratorNRS \otimes  \Delta ^{k} \mid s-n\geq q,  k\geq 0\}
		\end{eqnarray*}

	Since $\motivicAmod$ is in particular a simplicial model category
	(see proposition \ref{prop.2.8.motivicA-modules=>M-symmSpt-modelcat}(\ref{prop.2.8.motivicA-modules=>M-symmSpt-modelcat.b})), 
	using definitions 5.2.1, 16.3.1 and propositions 5.3.6, 16.1.3 
	in \cite{MR1944041}, we have that $f$ is a retract of a cofibration $g:M\rightarrow O$ in $\motivicAmod$
	for which there is a weak equivalence $h:O \rightarrow P$ in $\motivicAmod$
	such that the composition $h\circ g$ is a relative $\overline{\Lambda (K)}$-cell complex.
	
	It is clear that it is enough to check that $g$ is a cofibration in $\qconnectedsymmTspectra$.
	Now, using lemma 5.3.4  
	in \cite{MR1944041}, we have that this follows from:
		\begin{enumerate}
			\item	\label{prop.3.5.Acolocal=>colocalcofibrations=>cofibrations.a}  $g$ is a cofibration in $\motivicsymmTspectra$.
			\item	\label{prop.3.5.Acolocal=>colocalcofibrations=>cofibrations.b}	$h$ is a weak equivalence in $\motivicsymmTspectra$.
			\item	\label{prop.3.5.Acolocal=>colocalcofibrations=>cofibrations.c}  $h\circ g$ is a cofibration in $\qconnectedsymmTspectra$.
		\end{enumerate}
	
	(\ref{prop.3.5.Acolocal=>colocalcofibrations=>cofibrations.a}):  This follows directly from
	proposition \ref{prop.2.8.motivicAmod-cofibrations=>motivic-cofibrations}.
	
	(\ref{prop.3.5.Acolocal=>colocalcofibrations=>cofibrations.b}):  This follows directly from
	theorem \ref{thm.2.8.motivicsymmetricA-modules}.
	
	(\ref{prop.3.5.Acolocal=>colocalcofibrations=>cofibrations.c}): Let $\mathcal C$ denote the class of cofibrations
	in $\qconnectedsymmTspectra$.  Theorem
	\ref{thm.2.8.motivicAmods=>cellular} implies that $J_{A\text{-}\modules}$
	is a set of generating trivial cofibrations for $\motivicAmod$, therefore
	proposition \ref{prop.2.8.motivicAmod-cofibrations=>motivic-cofibrations} 
	and theorem \ref{thm.2.8.motivicsymmetricA-modules} 
	imply that all the maps in $J_{A\text{-}\modules}$
	are trivial cofibrations in $\motivicsymmTspectra$.
	But $\qconnectedsymmTspectra$ is a right Bousfield localization with respect to $\motivicsymmTspectra$,
	hence all the maps in $J_{A\text{-}\modules}$ are also trivial cofibrations in $\qconnectedsymmTspectra$.
	We have that in particular
		$$J_{A\text{-}\modules}=A\wedge J^{T}_{\Sigma}=\{ id\wedge j:A\wedge X\rightarrow A\wedge Y\}$$
	is contained in $\mathcal C$.  On the other hand, by hypothesis the map $\ast \rightarrow A$ is 
	a cofibration in $R_{C^{0}_{eff}}\motivicsymmTspectra$ and by construction 
	$\ast \rightarrow \generatorNRS$ are cofibrations in $\qconnectedsymmTspectra$ for $s-n\geq q$.
	Then theorem \ref{thm.3.4.smash-Quillenbifunctor-RpxRq----->Rqplusp} together with
	the fact that $\qconnectedsymmTspectra$ is a simplicial model category
	(see theorem \ref{thm.3.3.symmetricconnective-model-structure}) imply that
		\begin{eqnarray*}
			\{ A\wedge \symmgeneratorNRS \otimes \partial \Delta ^{k}
								\rightarrow 
								A\wedge \symmgeneratorNRS \otimes  \Delta ^{k} \mid \\ 
								s-n\geq q, k\geq 0\}
		\end{eqnarray*}
	is also contained in $\mathcal C$.
	Therefore, we have that all the maps in 
	$\overline{\Lambda(K)}$ are contained in $\mathcal C$.
	
	Finally since
	limits and colimits in $A\text{-}\modules$
	are computed in $\symmTspectra$, we have that $h\circ g$ is a relative $\mathcal C$-cell complex
	in $\symmTspectra$, 
	and since $\mathcal C$ is clearly closed under
	coproducts, pushouts and filtered colimits, we have that
	$h\circ g$ is a cofibration in $\qconnectedsymmTspectra$.
\end{proof}

\begin{thm}
		\label{thm.3.5.notcommutativesmash-Quillenbifunctor-RpxRq----->Rqplusp-coeffs}
	Fix $p,q \in \mathbb Z$.  Let $A$ be a cofibrant ring spectrum in $\motivicsymmTspectra$,
	which is also cofibrant in $\zeroconnectedsymmTspectra$.  Then
	$-\wedge _{A}-$ defines a Quillen adjunction of two variables
	(see definition \ref{def.Quillen-bifunct}) from the ($p-1$)-connected motivic model structure for right $A$-modules
	and the ($q-1$)-connected motivic model structure for left $A$-modules
	to the ($p+q-1$)-connected motivic symmetric stable model structure:
		$$\xymatrix{-\wedge _{A}-:
									\pconnectedAmod _{r}\times \qconnectedAmod _{l}\ar[r]& \pplusqconnectedsymmTspectra}
		$$
\end{thm}
\begin{proof}
	By lemma \ref{lem.cond-Quillen-bifunc}, it is enough to prove the following claim:
	
	Given a cofibration $i:N\rightarrow N'$ in $\qconnectedAmod _{l}$ and a fibration
		$f:X\rightarrow Y$ in $\pplusqconnectedsymmTspectra$, the induced map
			$$\xymatrix{\inthomsymmTspectra (N',X)\ar[d]^-{(i^{\ast},f_{\ast})}\\
									\inthomsymmTspectra (N,X)\times _{\inthomsymmTspectra (N,Y)}\inthomsymmTspectra (N',Y)}
			$$
		is a fibration in $\pconnectedAmod _{r}$ which is trivial if either $i$ or $f$ is a weak
		equivalence.
		
	However, proposition \ref{prop.3.5.Acolocal=>colocalcofibrations=>cofibrations}
	and lemma \ref{lemma.3.5.changecoefficients-Quillenfunctor}(\ref{lemma.3.5.changecoefficients-Quillenfunctor.b})
	imply that $i$ is also a cofibration in $\qconnectedsymmTspectra$,
	which is trivial if $i$ is a weak equivalence in $\qconnectedAmod _{l}$.  Now,
	it follows from theorem \ref{thm.3.4.smash-Quillenbifunctor-RpxRq----->Rqplusp} 
	that $(i^{\ast},f_{\ast})$ is a fibration in $\pconnectedsymmTspectra$
	which is trivial if either $i$ or $f$ is a weak equivalence.  By
	lemma \ref{lemma.3.5.changecoefficients-Quillenfunctor}(\ref{lemma.3.5.changecoefficients-Quillenfunctor.b})
	we have that it suffices to check that $(i^{\ast},f_{\ast})$ is a fibration in $\pconnectedAmod _{r}$.
	
	Proposition \ref{prop.3.2.functors-between-Rqsymm} implies that $(i^{\ast},f_{\ast})$
	is a fibration in $\motivicsymmTspectra$, hence it follows from 
	theorem \ref{thm.2.8.motivicsymmetricA-modules}
	that $(i^{\ast},f_{\ast})$ is also a fibration in $\motivicAmod _{r}$.  However,
	by construction $\pconnectedAmod _{r}$ is a right Bousfield localization with respect to
	$\motivicAmod _{r}$, therefore the classes of fibrations in both model structures are identical.
	Thus $(i^{\ast},f_{\ast})$ is a fibration in $\pconnectedAmod _{r}$, as we wanted.
\end{proof}

\begin{thm}
		\label{thm.3.5.inheritingmodulestructures1}
	Fix $q\in \mathbb Z$.  Let $A$ be a cofibrant ring spectrum in $\motivicsymmTspectra$,
	which is also $C_{eff}^{0,\Sigma}$-colocal in $\motivicsymmTspectra$
	(equivalently cofibrant in $\zeroconnectedsymmTspectra$), and let $M$ be an arbitrary $A$-module.
	Then the solid arrows in the following commutative diagram:
		\begin{equation}
					\label{thm.3.5.inheritingmodulestructures1.a}
			\begin{array}{c}
				\xymatrix@C=1.4pc{C_{q}^{\Sigma}R_{m}M \ar@{-->}[dd]^-{C_{q}^{\Sigma ,R_{m}M}} \ar[rr]^-{C_{q}^{\Sigma}(R_{m}^{M})}&& 
									C_{q}^{\Sigma}M \ar@{-->}[dd]_-{C_{q}^{\Sigma ,M}} && C_{q}^{\Sigma}C_{q}^{m}M \ar[ll]_-{C_{q}^{\Sigma}(C_{q}^{m,M})} 
									\ar[rr]^-{C_{q}^{\Sigma}(R_{m}^{C_{q}^{m}M})} \ar@{-->}[dd]_-{C_{q}^{\Sigma ,C_{q}^{m}M}}&& 
									C_{q}^{\Sigma}R_{m}C_{q}^{m}M \ar[dd]_-{C_{q}^{\Sigma ,R_{m}C_{q}^{m}M}}\\
									&&&&&& \\
								  R_{m}M \ar@{-->}[rr]_-{R_{m}^{M}} && M && 
								  C_{q}^{m}M \ar@{-->}[ll]^-{C_{q}^{m,M}} \ar@{-->}[rr]_-{R_{m}^{C_{q}^{m}M}} && R_{m}C_{q}^{m}M}
			\end{array}
		\end{equation}
	induce a natural equivalence between the functors:
		\begin{equation}
					\label{thm.3.5.inheritingmodulestructures1.b}
			\begin{array}{c}
				\xymatrix{& \Amodstablehomotopy \ar[dr]^-{UR_{m}}&\\
												 \qconnectedAmodstablehomotopy
												 \ar[ur]^-{C^{m}_{q}} \ar[dr]_-{UR_{m}}&& \symmstablehomotopy \\
												& \qconnectedsymmstablehomotopy \ar[ur]_-{C^{\Sigma}_{q}}&}
			\end{array}
		\end{equation}
\end{thm}
\begin{proof}
	Clearly it suffices to show that $C_{q}^{\Sigma}(R_{m}^{M})$, $C_{q}^{\Sigma}(C_{q}^{m,M})$, 
	$C_{q}^{\Sigma}(R_{m}^{C_{q}^{m}M})$ and 
	$C_{q}^{\Sigma ,R_{m}C_{q}^{m}M}$ are all weak equivalences in $\motivicsymmTspectra$.
	
	Proposition \ref{prop.3.5.Rsigma-fibrant-replacement-all-Rq} 
	implies that $R_{m}^{M}$ is a weak equivalence in $\qconnectedAmod$, then applying lemma
	\ref{lemma.3.5.changecoefficients-Quillenfunctor}(\ref{lemma.3.5.changecoefficients-Quillenfunctor.b})
	to the unit map $\symmspherespectrum \rightarrow A$
	we have that $R_{m}^{M}$ is a weak equivalence in $\qconnectedsymmTspectra$.  By construction
	$C_{q}^{\Sigma ,R_{m}M}$ and $C_{q}^{\Sigma ,M}$ are both weak equivalences in $\qconnectedsymmTspectra$.  Hence, the two out of three
	property for weak equivalences implies that $C_{q}^{\Sigma}(R_{m}^{M})$ is a weak equivalence
	in $\qconnectedsymmTspectra$.  However, $C_{q}^{\Sigma}R_{m}M$ and $C_{q}^{\Sigma}M$ are both $C_{eff}^{q,\Sigma}$-colocal;
	therefore \cite[theorem 3.2.13(2)]{MR1944041}
	implies that $C_{q}^{\Sigma}(R_{m}^{M})$ is a weak equivalence in 
	$\motivicsymmTspectra$.
	
	Using lemma \ref{lemma.3.5.changecoefficients-Quillenfunctor}(\ref{lemma.3.5.changecoefficients-Quillenfunctor.b})
	again, we have that $C_{q}^{m,M}$ is a weak equivalence in $\qconnectedsymmTspectra$.  But
	$C_{q}^{\Sigma ,M}$ and $C_{q}^{\Sigma ,C_{q}^{m}M}$ are both weak equivalences in $\qconnectedsymmTspectra$,
	thus the two out of three property for weak equivalences implies that $C_{q}^{\Sigma}(C_{q}^{m,M})$
	is also a weak equivalence in $\qconnectedsymmTspectra$.  However, by construction $C_{q}^{\Sigma}M$ and $C_{q}^{\Sigma}C_{q}^{m}M$
	are both $C_{eff}^{q, \Sigma}$-colocal; hence by \cite[theorem 3.2.13(2)]{MR1944041}
	we have that $C_{q}^{\Sigma}(C_{q}^{m,M})$ is a weak equivalence in $\motivicsymmTspectra$.
	
	By proposition \ref{prop.3.5.Rsigma-fibrant-replacement-all-Rq} 
	we have that $R_{m}^{C_{q}^{m}M}$ is a weak equivalence in $\qconnectedAmod$, then lemma
	\ref{lemma.3.5.changecoefficients-Quillenfunctor}(\ref{lemma.3.5.changecoefficients-Quillenfunctor.b})
	implies that $R_{m}^{C_{q}^{m}M}$ is a weak equivalence in $\qconnectedsymmTspectra$.  Now,
	$C_{q}^{\Sigma ,C_{q}^{m}M}$ and $C_{q}^{\Sigma ,R_{m}C_{q}^{m}M}$ 
	are both weak equivalences in $\qconnectedsymmTspectra$.  Thus, the two out of three
	property for weak equivalences implies that $C_{q}^{\Sigma}(R_{m}^{C_{q}^{m}M})$ is a weak equivalence
	in $\qconnectedsymmTspectra$.  However,
	by construction, $C_{q}^{\Sigma}C_{q}^{m}M$ and $C_{q}^{\Sigma}R_{m}C_{q}^{m}M$ are both $C_{eff}^{q,\Sigma}$-colocal;
	therefore \cite[theorem 3.2.13(2)]{MR1944041}
	implies that $C_{q}^{\Sigma}(R_{m}^{C_{q}^{m}M})$ is a weak equivalence in 
	$\motivicsymmTspectra$.
	
	We already know that $C_{q}^{\Sigma}(R_{m}^{C_{q}^{m}M})$ is a weak equivalence in $\motivicsymmTspectra$,
	and definition \ref{def.3.5.fibrantstable-replacementfunctor} together with 
	theorem \ref{thm.2.8.motivicsymmetricA-modules} imply
	that $R_{m}^{C_{q}^{m}M}$ is also a weak equivalence in $\motivicsymmTspectra$.
	Therefore, to show that $C_{q}^{\Sigma ,R_{m}C_{q}^{m}M}$ is a weak equivalence in $\motivicsymmTspectra$,
	it suffices to check that $C_{q}^{\Sigma ,C_{q}^{m}M}$ is a weak equivalence in $\motivicsymmTspectra$.
	Now, by construction we have that $C_{q}^{\Sigma ,C_{q}^{m}M}$ is a $C_{eff}^{q,\Sigma}$-colocal equivalence
	in $\motivicsymmTspectra$
	and that $C_{q}^{\Sigma}C_{q}^{m}M$ is a $C_{eff}^{q,\Sigma}$-colocal spectrum, thus by \cite[theorem 3.2.13(2)]{MR1944041}
	it only remains to show that $C_{q}^{m}M$ 
	is $C_{eff}^{q, \Sigma}$-colocal.  But this follows from
	our hypothesis which says that $A$ is $C_{eff}^{0,\Sigma}$-colocal together with
	proposition \ref{prop.3.5.Acolocal=>colocalcofibrations=>cofibrations}.
	This finishes the proof. 
\end{proof}

\begin{thm}
		\label{thm.3.5.connectedinvarianceofcoefficients}
	Fix $q\in \mathbb Z$.  Let $f:A\rightarrow A'$ be a map between 
	cofibrant ring spectra in $\motivicsymmTspectra$,
	which is compatible with the ring structures.
	Assume that there exists $p\in \mathbb Z$ such that $A$, $A'$ are
	both $C_{eff}^{p,\Sigma}$-colocal in $\motivicsymmTspectra$ and
	$f$ is a weak equivalence in $\pconnectedsymmTspectra$ (equivalently in $\pconnectedAmod$).
	Then $f$ induces a Quillen equivalence between the ($q-1$)-connected motivic 
	stable model structures of $A$ and $A'$ modules:
		$$\xymatrix{(A'\wedge _{A} -, U,\varphi ):\qconnectedAmod \ar[r] & \qconnectedAAmod}
		$$
\end{thm}
\begin{proof}
	Since $A$ and $A'$ are $C_{eff}^{p,\Sigma}$-colocal, \cite[theorem 3.2.13(2)]{MR1944041}
	implies that $f$ is a weak equivalence in $\motivicsymmTspectra$.  Therefore,
	the result follows directly from proposition
	\ref{prop.3.5.changecoefficients-Quillenequiv}.
\end{proof}

\begin{thm}
		\label{thm.3.5.Lqsymmetricmodelstructures}
	Fix $q\in \mathbb Z$.  Consider the following
	set of maps in $\motivicAmod$ (see theorem \ref{thm.3.3.Lqsymmetricmodelstructures}):
		\begin{eqnarray}
				\label{diagram.3.5.symmmapsleftloc}
			\qleftAmodlocmaps = \{ \AmodkillNRS \: | 
										 \symmgeneratorNRS \in C_{eff}^{q,\Sigma}\} \nonumber
		\end{eqnarray}
	Then the left Bousfield localization of $\motivicAmod$ with respect to
	the $\qleftAmodlocmaps$-local equivalences exists.  This new model structure
	will be called \emph{weight$^{<q}$ motivic stable}.
	$\weightqAmod$ will denote the category of $A$-modules equipped with
	the weight$^{<q}$ motivic stable model
	structure, and $\weightqAmodstablehomotopy$ will denote
	its associated homotopy category.  Furthermore the weight$^{<q}$ motivic stable model structure
	is cellular, left proper and simplicial; with the following
	sets of generating cofibrations and trivial cofibrations
	respectively:
	$$\begin{array}{l}
	 I_{\qleftAmodlocmaps}=I_{A\text{-}\modules}  =\bigcup _{n\geq 0}\{A\wedge F_{n}^{\Sigma}(Y_{+}\hookrightarrow (\Delta ^{n}_{U})_{+})\} \\
	 \\
	J_{\qleftAmodlocmaps}  =\{j:A\rightarrow B\}
	\end{array}
	$$
	where $j$ satisfies the following conditions:
	\begin{enumerate}
		\item	$j$ is an inclusion of $I^{T}_{A\text{-}\modules}$-complexes.
		\item	$j$ is a $\qleftAmodlocmaps$-local equivalence.
		\item	the size of $B$ as an $I^{T}_{A\text{-}\modules}$-complex is less than $\kappa$, 
					where $\kappa$ is the regular cardinal defined by Hirschhorn in \cite[definition 4.5.3]{MR1944041}.
	\end{enumerate}
\end{thm}

\begin{rmk}
		\label{rmk.3.5.usualargumentspartiallynotgood}
	Notice that the model category $\weightqsymmTspectra$ is not
	a symmetric monoidal model category, i.e. the smash product and
	the model structure are not compatible, therefore in general
	it is not possible to use the adjuntion 
		$$(A\wedge -,U,\varphi ):\symmTspectra \rightarrow A\text{-}\modules $$
	for the construction of a model structure on the category of $A$-modules.
	However, if $A$ satisfies additional conditions 
	(see proposition \ref{prop.3.5.smash-induces-weightq-modelstructures}) then
	the adjunction above induces a model structure on the category of $A$-modules
	which coincides with
	$\weightqAmod$ (see proposition \ref{prop.3.5.smash-induces-weightq-modelstructures} 
	and theorem \ref{thm.3.5.weightq-induced==leftBousfieldlocalization}).
\end{rmk}

\begin{defi}
		\label{def.3.5.stable-weightqsymm-replacementfunctors}
	Fix $q\in \mathbb Z$.  Let $W_{q}^{m}$ denote a fibrant replacement functor in
	$\weightqAmod$, such that the for every $A$-module $M$,
	the natural map:
		$$\xymatrix{M \ar[r]^-{W_{q}^{m, M}}& W_{q}^{m}M}
		$$
	is a trivial cofibration in $\weightqAmod$, and $W_{q}^{m}M$
	is $\qleftAmodlocmaps$-local in $\motivicAmod$.
\end{defi}

\begin{prop}
		\label{prop.3.5.Q-cofibrant-replacement-all-symmetricL<q}
	Fix $q\in \mathbb Z$.  Then $Q_{m}$ is also a cofibrant 
	(see definition \ref{def.3.5.cofibrantstable-replacementfunctor})
	replacement functor in $\weightqAmod$, and for every $A$-module
	$M$ the natural map
		$$\xymatrix{Q_{m}M\ar[r]^{Q_{m}^{M}}& M}
		$$
	is a trivial fibration in $\weightqAmod$.
\end{prop}

\begin{prop}
		\label{prop.3.5.Z-AmodLq-local-iff-UZ-symmLq-local}
	Fix $q\in \mathbb Z$.
	Then an $A$-module $M$ is $\qleftAmodlocmaps$-local in $\motivicAmod$
	if and only if $UM$ is $\qleftsymmlocmaps$-local in $\motivicsymmTspectra$.
\end{prop}

\begin{prop}
		\label{prop.3.5.susp-Quillen-equiv-onsymmLq}
	For every $q\in \mathbb Z$, the following adjunction:
		$$\xymatrix{(-\wedge S^{1},\Omega _{S^{1}},\varphi):\weightqAmod \ar[rr]
								&& \weightqAmod}
		$$
	is a Quillen equivalence, and the homotopy category $\weightqAmodstablehomotopy$
	associated to $\weightqAmod$ has the
	structure of a triangulated category.
\end{prop}

\begin{cor}
		\label{cor.3.5.symmLq=>rightproper}
	For every $q\in \mathbb Z$, $\weightqAmod$ is
	a right proper model category.
\end{cor}

\begin{prop}
		\label{prop.3.5.symmLq-exact-adjunctions}
	For every $q\in \mathbb Z$ we have the following Quillen
	adjunction
		$$\xymatrix{(id,id,\varphi ):\motivicAmod \ar[r]& \weightqAmod}
		$$
	which induces and adjunction
		$$\xymatrix{(Q_{m},W_{q}^{m},\varphi):\Amodstablehomotopy \ar[r]& \weightqAmodstablehomotopy}
		$$
	of exact functors between triangulated categories.
\end{prop}

\begin{thm}
		\label{thm.3.5.symmetrization-Quillen-equivalence-weight<q}
	For every $q\in \mathbb Z$, the adjunction
		$$\xymatrix{(A\wedge -,U,\varphi):\weightqsymmTspectra \ar[r]& \weightqAmod}
		$$
	given by the free $A$-module and the forgetful functor is a Quillen adjunction,
	and it induces
	an adjunction
		$$\xymatrix{(A\wedge Q_{s}-,UW_{q}^{m},\varphi):\weightqsymmstablehomotopy \ar[r]& \weightqAmodstablehomotopy}
		$$
	of exact funtors between triangulated categories.
\end{thm}
		
\begin{prop}
		\label{prop.3.5.Lq+1-->Lqsymm}
	Fix  $q\in \mathbb Z$.  We have the following
	commutative diagram of left Quillen functors:
		\begin{equation}
				\label{diagram.prop.3.5.Lq+1-->Lqsymm.a}
			\begin{array}{c}
				\xymatrix{& \motivicsymmTspectra \ar[dl]_-{id} \ar[dr]^-{id} \ar[dd]^(0.3){A\wedge -}&\\
									L_{<q+1}\motivicsymmTspectra \ar[rr]_(0.3){id} \ar[dd]^-{A\wedge -}
									&& \weightqsymmTspectra \ar[dd]^-{A\wedge -} \\
									& \motivicAmod \ar[dl]_-{id} \ar[dr]^-{id}&\\
									L_{<q+1}\motivicAmod \ar[rr]_-{id}&& \weightqAmod}
			\end{array}
		\end{equation}
	and the following associated commutative diagrams of homotopy categories:
		\begin{equation}
				\label{diagram.prop.3.5.Lq+1-->Lqsymm.b}
			\begin{array}{c}
				\xymatrix{& \Amodstablehomotopy \ar[dr]^-{Q_{m}}&\\
											 \symmstablehomotopy
												\ar[rr]^{A\wedge Q_{\Sigma}-} \ar[ur]^-{A\wedge Q_{\Sigma}-} \ar[dr]_-{Q_{\Sigma}}&& \weightqAmodstablehomotopy \\
												& \weightqsymmstablehomotopy \ar[ur]_-{A\wedge Q_{\Sigma}-}&}
			\end{array}
		\end{equation}
		
		\begin{equation}
				\label{diagram.prop.3.5.Lq+1-->Lqsymm.c}
			\begin{array}{c}
				\xymatrix{& \Amodstablehomotopy \ar[dr]^-{UR_{m}}&\\
												 \weightqAmodstablehomotopy
												\ar[rr]^{UW_{q}^{m}} \ar[ur]^-{W_{q}^{m}} \ar[dr]_-{UW_{q}^{m}}&& \symmstablehomotopy \\
												& \weightqsymmstablehomotopy \ar[ur]_-{W_{q}^{\Sigma}}&}
			\end{array}
		\end{equation}
\end{prop}

\begin{thm}
		\label{thm.3.5.motivic-towersymm}
	We have the following commutative diagram
	of left Quillen functors:
		\begin{equation}
					\label{diagram.3.5.thm.3.5.motivic-towersymm.a}
			\begin{array}{c}
				\xymatrix@C=0.6pc{& \vdots \ar[d]_-{id} & \vdots \ar[d]_-{id} & \\ 
									& L_{<q+1}\motivicsymmTspectra \ar[d]_-{id} \ar[r]_-{A\wedge -} \ar@{<-}[dl]_-{id}& 
									L_{<q+1}\motivicAmod \ar[d]_-{id} \ar@{<-}[dr]^-{id} & \\ 
									\motivicsymmTspectra & \weightqsymmTspectra \ar[d]_-{id} \ar[r]_-{A\wedge -} \ar@{<-}[l]_-{id}& 
									\weightqAmod \ar[d]_-{id} \ar@{<-}[r]^-{id} & \motivicAmod \\
									& L_{<q-1}\motivicsymmTspectra \ar[d]_-{id} \ar[r]_-{A\wedge -} \ar@{<-}[ul]^-{id}& 
									L_{<q-1}\motivicAmod \ar[d]_-{id} \ar@{<-}[ur]_-{id} & \\ 
									& \vdots & \vdots & }
			\end{array}
		\end{equation}
	and the following associated commutative diagrams of homotopy categories:	
		\begin{equation}
					\label{diagram.3.5.thm.3.5.motivic-towersymm.b}
			\begin{array}{c}
				\xymatrix@C=0.3pc{&&& \vdots \ar[d] && \vdots \ar[d] &&& \\
									&&&
				 					L_{<q+1}\symmstablehomotopy \ar[d]_-{Q_{\Sigma}} 
				 					\ar@{<-}[dlll]|-{Q_{\Sigma}}  \ar[rr]_-{A\wedge Q_{\Sigma}-} &&
									L_{<q+1}\Amodstablehomotopy \ar[d]_-{Q_{m}} 
									\ar@{<-}[drrr]|-{Q_{m}}  &&& \\ 
									\symmstablehomotopy &&&
									\weightqsymmstablehomotopy \ar[d]_-{Q_{\Sigma}} 
									\ar@{<-}[lll]|-{Q_{\Sigma}}  \ar[rr]_-{A\wedge Q_{\Sigma}-} &&  
									\weightqAmodstablehomotopy \ar[d]_-{Q_{m}} 
									\ar@{<-}[rrr]|-{Q_{m}}  &&& \Amodstablehomotopy \\
									&&& L_{<q-1}\symmstablehomotopy \ar[d] \ar@{<-}[ulll]|-{Q_{\Sigma}} 
									\ar[rr]_-{A\wedge Q_{\Sigma}-} &&
									L_{<q-1}\Amodstablehomotopy \ar[d] \ar@{<-}[urrr]|-{Q_{m}} &&& \\ 
									&&& \vdots  &&
									\vdots &&& }
			\end{array}
		\end{equation}
		
	\begin{equation}
					\label{diagram.3.5.thm.3.5.motivic-towersymm.c}
			\begin{array}{c}
				\xymatrix@C=0.3pc{&&& \vdots  && \vdots  &&& \\
									&&&
				 					L_{<q+1}\symmstablehomotopy  \ar[u]
				 					&& \ar[ll]^-{UW^{m}_{q+1}} 
									L_{<q+1}\Amodstablehomotopy  \ar[u]
									&&&  \\ 
									\symmstablehomotopy  \ar@{<-}[urrr]|-{W^{\Sigma}_{q+1}} \ar@{<-}[rrr]|-{W^{\Sigma}_{q}} \ar@{<-}[drrr]|-{W^{\Sigma}_{q-1}}&&&
									\weightqsymmstablehomotopy \ar[u]^-{W^{\Sigma}_{q}}  
									&& \ar[ll]^-{UW^{m}_{q}}  
									\weightqAmodstablehomotopy  \ar[u]^-{W^{m}_{q}}   
								  &&& \Amodstablehomotopy  \ar@{<-}[ulll]|-{W^{m}_{q+1}} \ar@{<-}[lll]|-{W^{m}_{q}} \ar@{<-}[dlll]|-{W^{m}_{q-1}} \\
									&&& L_{<q-1}\symmstablehomotopy  \ar[u]^-{W^{\Sigma}_{q-1}}  
									&&
									L_{<q-1}\Amodstablehomotopy  \ar[u]^-{W^{m}_{q-1}}  \ar[ll]^-{UW^{m}_{q-1}} &&& \\ 
									&&& \vdots  \ar[u] &&
									\vdots \ar[u] &&& }
			\end{array}
		\end{equation}
\end{thm}

\begin{lem}
		\label{lemma.3.5.changecoefficients-Quillenfunctor2}
	Fix $q \in \mathbb Z$, and let $f:A\rightarrow  A'$ 
	be a map between cofibrant ring spectra $\motivicsymmTspectra$, 
	which is compatible with the ring structures. If $g: M\rightarrow N$
	is a $\qleftAmodlocmaps$-local equivalence in 
	$\motivicAmod$ then 
	$id\wedge _{A}Q_{m}g:A'\wedge _{A}Q_{m}M \rightarrow A'\wedge_{A}Q_{m}N$ 
	is a $\qleftAAmodlocmaps$-local equivalence in 
	$\motivicAAmod$, where $Q_{m}$ denotes a
	cofibrant replacement functor in $\motivicAmod$.
\end{lem}
\begin{proof}
	Let $Z$ be an arbitrary 
	$\qleftAAmodlocmaps$-local $A'$-module in $\motivicAAmod$.
	Lemma \ref{lemma.2.5.changecoefficients-Quillenfunctor} implies that $A'\wedge _{A}Q_{m}M$, 
	$A'\wedge _{A}Q_{m}N$ are both cofibrant in $\motivicAAmod$.
	Therefore it suffices to show that the induced map
		$$\xymatrix{Map_{A'\text{-}\modules}(A'\wedge _{A}Q_{m}N,Z) 
								\ar[rr]^-{(id\wedge _{A}Q_{m}g)^{\ast}} && 
								Map_{A'\text{-}\modules}(A'\wedge _{A}Q_{m}M,Z)}$$
	is a weak equivalence of simplical sets. However, using proposition \ref{prop.2.8.enriched-freemodule-adjunction2} 
	we get the following commutative diagram, 
	where the vertical maps are isomorphisms of
	simplicial sets
		$$\xymatrix{Map_{A'\text{-}\modules}(A'\wedge _{A}Q_{m}N,Z) 
								\ar[rr]^-{(id\wedge _{A}Q_{m}g)^{\ast}} \ar[d]_-{\cong} &&
								Map_{A'\text{-}\modules}(A'\wedge _{A}Q_{m}M,Z) \ar[d]^-{\cong}\\
								Map_{A\text{-}\modules}(Q_{m}N,UZ) \ar[rr]_-{(Q_{m}g)^{\ast}}
								&& Map_{A\text{-}\modules}(Q_{m}M,UZ)}$$
	Finally, proposition \ref{prop.3.5.Z-AmodLq-local-iff-UZ-symmLq-local} implies that $UZ$ is 
	$\qleftAmodlocmaps$-local in $\motivicAmod$, therefore
	the bottom row is a weak equivalence of simplicial sets, since by hypothesis $g$
	is a $\qleftAmodlocmaps$-local equivalence in 
	$\motivicAmod$. Hence, the two out of three property
	for weak equivalences implies that the top row is also a weak equivalence of
	simplicial sets, as we wanted.
\end{proof}

\begin{prop}
		\label{prop.3.5.changecoefficients-Quillenfunctor2}
	Fix $q\in \mathbb Z$.  Let $f:A\rightarrow A'$ be a map 
	between cofibrant ring spectra in $\motivicsymmTspectra$,
	which is compatible with the ring structures.
	Then the adjunction:
		$$(A'\wedge _{A}-, U,\varphi ):\weightqAmod \rightarrow \weightqAAmod
		$$
	is a Quillen adjunction.
\end{prop}
\begin{proof}
	Lemma \ref{lemma.2.5.changecoefficients-Quillenfunctor} 
	implies that $A'\wedge _{A} -:\motivicAmod \rightarrow \motivicAAmod$ is
	a left Quillen functor. Consider the following commutative diagram of left Quillen
	functors:
		$$\xymatrix{\motivicAmod \ar[r]^-{A'\wedge _{A}-} \ar[d]_-{id}&
								\motivicAAmod \ar[d]^-{id}\\
								\weightqAmod \ar@{-->}[r]_-{A'\wedge _{A}-} &
								\weightqAAmod}$$
	then the universal property of left Bousfield localizations together with lemma
	\ref{lemma.3.5.changecoefficients-Quillenfunctor2} 
	imply that the dotted arrow $A'\wedge _{A}-$ is a left Quillen functor.
\end{proof}

\begin{lem}
		\label{lemma.3.5.changecoeffs-preserves-localobjects}
	Fix $q \in \mathbb Z$, and let $f:A\rightarrow A'$ 
	be a map between cofibrant
	ring spectra in $\motivicsymmTspectra$, 
	which is compatible with the ring structures. If f is a weak
	equivalence in $\motivicsymmTspectra$
	(equivalently in $\motivicAmod$), then for every 
	$\qleftAmodlocmaps$-local
	$A$-module $M$ in $\motivicAmod$, we have that $Q_{m}X$ 
	and $UR_{m'}(A'\wedge _{A}Q_{m}X)$ are also
	$\qleftAmodlocmaps$-local in $\motivicAmod$, 
	where $Q_{m}$ denotes a cofibrant replacement functor
	in $\motivicAmod$ and $R_{m'}$ denotes a fibrant replacement functor 
	in $\motivicAAmod$.
\end{lem}
\begin{proof}
	Since $M$ is $\qleftAmodlocmaps$-local, it follows that 
	$M$ is fibrant in $\motivicAmod$.
	By definition we have that the natural map
		$$\xymatrix{Q_{m}M \ar[r]^-{Q_{m}^{M}} & M}$$
	is a trivial fibration in $\motivicAmod$, 
	therefore $Q_{m}M$ is also fibrant in $\motivicAmod$.
	Hence \cite[lemma 3.2.1(a)]{MR1944041} implies that $Q_{m}M$ 
	is $\qleftAmodlocmaps$-local.
	Proposition \ref{prop.2.8.invarianceofcoefficients} implies that the adjunction 
	$(A'\wedge _{A}-,U, \varphi)$ is a Quillen equivalence
	between $\motivicAmod$ and $\motivicAAmod$, therefore we have that 
	$UR_{m'}(A'\wedge _{A}Q_{m}M)$ is fibrant in $\motivicAmod$, 
	and \cite[proposition 1.3.13(b)]{MR1650134} implies that
	the composition
		$$\xymatrix{Q_{m}M \ar[r]^-{\eta _{Q_{m}M}}& U(A'\wedge _{A}Q_{m}M)
 								\ar[rr]^-{U(R_{m'}^{A'\wedge _{A}Q_{m}M})} && UR_{m'}(A'\wedge _{A}Q_{m}M)}$$
	is a weak equivalence in $\motivicAmod$. Since we already know 
	that $Q_{m}M$ is $\qleftAmodlocmaps$-local, 
	using \cite[lemma 3.2.1(a)]{MR1944041} again we get that 
	$UR_{m'}(A'\wedge _{A}Q_{m}M)$ is also
	$\qleftAmodlocmaps$-local in $\motivicAmod$. This finishes the proof.
\end{proof}

\begin{lem}
		\label{lemma.3.5.changecoeffs-preserves-Quillenequivs}
	Fix $q \in \mathbb Z$, and let $f:A\rightarrow A'$ 
	be a map between cofibrant
	ring spectra $\motivicsymmTspectra$, 
	which is compatible with the ring structures. If $f$ is a weak
	equivalence in $\motivicsymmTspectra$ 
	(equivalently in $\motivicAmod$), then $g:M\rightarrow N$ is a 
	$\qleftAmodlocmaps$-local equivalence in $\motivicAmod$ 
	if and only if 
	$id\wedge _{A}Q_{m}g:A'\wedge _{A}Q_{m}M \rightarrow A'\wedge _{A}Q_{m}N$ 
	is a $\qleftAAmodlocmaps$-local equivalence in $\motivicAAmod$, 
	where $Q_{m}$ denotes a
	cofibrant replacement functor in $\motivicAmod$.
\end{lem}
\begin{proof}
	($\Rightarrow$): It follows directly from lemma \ref{lemma.3.5.changecoefficients-Quillenfunctor2}.
	
	($\Leftarrow$): Assume that $id\wedge _{A}Q_{m}g$ is a 
	$\qleftAAmodlocmaps$-local equivalence in $\motivicAAmod$,
	and let $Z$ be an arbitrary $\qleftAmodlocmaps$-local $A$-module in 
	$\motivicAmod$. We need to
	show that the induced map:
		$$\xymatrix{Map_{A\text{-}\modules}(Q_{m}N,Z) \ar[r]^-{(Q_{m}g)^{\ast}} &
								Map_{A\text{-}\modules}(Q_{m}M,Z)}$$
	is a weak equivalence of simplicial sets.
	
	But proposition \ref{prop.2.8.invarianceofcoefficients} implies that the adjunction $(A'\wedge _{A}-,U, \varphi)$ 
	is a Quillen equivalence between 
	$\motivicAmod$ and $\motivicAAmod$, therefore using 
	\cite[proposition 1.3.13(b)]{MR1650134} 
	we have that all the maps in the following diagram are weak equivalences
	in $\motivicAmod$:
		$$\xymatrix{Z & \ar[l]_-{Q_{m}^{Z}} Q_{m}Z
								\ar[rrrr]^-{U(R_{m'}^{A'\wedge _{A}Q_{m}Z})\circ \eta _{Q_{m}Z}} &&&&
								UR_{m'}(A'\wedge _{A}Q_{m}Z)}$$
	where $R_{m'}$ denotes a fibrant replacement functor in $\motivicAAmod$. 
	Lemma \ref{lemma.3.5.changecoeffs-preserves-localobjects}
	implies in particular that $Z$, $Q_{m}Z$, $UR_{m'}(A'\wedge _{A}Q_{m}Z)$ 
	are all fibrant in $\motivicAmod$.
	Now using the fact that $\motivicAmod$ is a simplicial model category together with
	Ken Brown's lemma (see lemma \ref{lemma1.1.KenBrown2}) and the two out of three property for weak
	equivalences, we have that it suffices to prove that the induced map:
		$$\xymatrix{Map_{A\text{-}\modules}(Q_{m}N,UR_{m'}(A'\wedge _{A}Q_{m}Z)) 
								\ar[d]^-{(Q_{m}g)^{\ast}} \\
								Map_{A\text{-}\modules}(Q_{m}M,UR_{m'}(A'\wedge _{A}Q_{m}Z))}$$
	is a weak equivalence of simplicial sets. Using the enriched adjunctions of proposition
	\ref{prop.2.8.enriched-freemodule-adjunction2}, we get the following commutative diagram where all the vertical arrows
	are isomorphisms:
		$$\xymatrix@C=-4pc{Map_{A\text{-}\modules}(Q_{m}N,UR_{m'}(A'\wedge _{A}Q_{m}Z))
								\ar[dr]^-{(Q_{m}g)^{\ast}} \ar[dd]_{\cong} & \\
								& Map_{A\text{-}\modules}(Q_{m}M,UR_{m'}(A'\wedge _{A}Q_{m}Z)) \ar[dd]^-{\cong}\\
								Map_{A'\text{-}\modules}(A'\wedge _{A}Q_{m}N,R_{m'}(A'\wedge _{A}Q_{m}Z))
								\ar[dr]^-{(id\wedge _{A}Q_{m}g)^{\ast}} & \\
								& Map_{A'\text{-}\modules}(A'\wedge _{A}Q_{m}M,R_{m'}(A'\wedge _{A}Q_{m}Z))}$$
	Finally, lemma \ref{lemma.3.5.changecoeffs-preserves-localobjects} implies that $UR_{m'}(A'\wedge _{A}Q_{m}Z)$ 
	is $\qleftAmodlocmaps$-local in $\motivicAmod$,
	therefore by proposition \ref{prop.3.5.Z-AmodLq-local-iff-UZ-symmLq-local} we have that $R_{m'}(A'\wedge _{A}Q_{m}Z)$ 
	is $\qleftAAmodlocmaps$-local in
	$\motivicAAmod$. Since $id\wedge _{A}Q_{m}g$ is a $\qleftAAmodlocmaps$-local 
	equivalence and $A'\wedge _{A}Q_{m}M$, $A'\wedge _{A}Q_{m}N$
	are both cofibrant in $\motivicAAmod$, it follows that the bottom row in the diagram
	above is a weak equivalence of simplicial sets. This implies that the top row
	is also a weak equivalence of simplicial sets, as we wanted.
\end{proof}

\begin{prop}
		\label{prop.3.5.changecoefficients-Quillenfunctor2a}
	Fix $q\in \mathbb Z$.  Let $f:A\rightarrow A'$ be a map 
	between cofibrant ring spectra in $\motivicsymmTspectra$,
	which is compatible with the ring structures.
	If $f$ is
	a weak equivalence in $\motivicsymmTspectra$ 
	then the adjunction 
		$$(A'\wedge _{A}-, U,\varphi ):\weightqAmod \rightarrow \weightqAAmod
		$$
	is a Quillen equivalence.  
\end{prop}
\begin{proof}
	Proposition \ref{prop.3.5.changecoefficients-Quillenfunctor2} 
	implies that the adjuntion $(A'\wedge _{A} -,U,\varphi )$ is a
	Quillen adjunction.
	Using corollary 1.3.16 in \cite{MR1650134} and proposition \ref{prop.3.5.Q-cofibrant-replacement-all-symmetricL<q} 
	we have that it suffices
	to verify the following two conditions:
		\begin{enumerate}
			\item \label{prop.3.5.changecoefficients-Quillenfunctor2a.a} For every fibrant $A'$-module $M$ in 
										 $\weightqAAmod$, the following composition
												$$\xymatrix{A'\wedge _{A} Q_{m}UM \ar[rr]^-{id\wedge _{A}(Q_{m}^{UM})} &&
																		A'\wedge _{A}UM \ar[r]^-{\epsilon _{M}} & M}$$
											is a weak equivalence in $\weightqAAmod$, where $Q_{m}$ denotes
											a cofibrant replacement functor in $\motivicAmod$ (see proposition
											\ref{prop.3.5.Q-cofibrant-replacement-all-symmetricL<q}).
			\item \label{prop.3.5.changecoefficients-Quillenfunctor2a.b} $A'\wedge _{A}-$ reflects weak equivalences 
										 between cofibrant $A$-modules in $\weightqAmod$.
		\end{enumerate}
		
	(\ref{prop.3.5.changecoefficients-Quillenfunctor2a.a}): By construction $\weightqAAmod$ is a left Bousfield localization 
	of $\motivicAAmod$,
	therefore the identity functor
		$$id:\weightqAAmod \rightarrow \motivicAAmod$$
	is a right Quillen functor. Thus $M$ is also fibrant in $\motivicAAmod$. Proposition
	\ref{prop.2.8.invarianceofcoefficients} implies that the adjunction $(A'\wedge _{A}-,U,\varphi )$ 
	is a Quillen equivalence between
	$\motivicAmod$ and $\motivicAAmod$, hence using 
	\cite[proposition 1.3.13(b)]{MR1650134} we have
	that the following composition is a weak equivalence in $\motivicAAmod$:
		$$\xymatrix{A'\wedge _{A}Q_{m}UM \ar[rr]^-{id\wedge _{A}(Q_{m}^{UM})} &&
								A'\wedge _{A}UM \ar[r]^-{\epsilon _{M}} & M}$$
	Therefore \cite[proposition 3.1.5]{MR1944041} 
	implies that the composition above is a $\qleftAAmodlocmaps$-local equivalence.
	
	(\ref{prop.3.5.changecoefficients-Quillenfunctor2a.b}): This follows immediately from 
	proposition \ref{prop.3.5.Q-cofibrant-replacement-all-symmetricL<q} and lemma \ref{lemma.3.5.changecoeffs-preserves-Quillenequivs}.
\end{proof}

\begin{prop}
		\label{prop.3.5.smash-induces-weightq-modelstructures}
	Fix $q\in \mathbb Z$.  Let $A$ be a cofibrant ring spectrum in $\motivicsymmTspectra$,
	which is also cofibrant in $\zeroslicesymmTspectra$.  Then
	the adjunction 
		$$(A\wedge -,U,\varphi ):\symmTspectra \rightarrow A\text{-}\modules
		$$
	between symmetric $T$-spectra and $A$-modules, together with the model structure
	$\weightqsymmTspectra$ (see theorem \ref{thm.3.3.Lqsymmetricmodelstructures}), 
	induces a model structure on $A\text{-}\modules$, which we
	will denote by $\weightqAmodinduced$; i.e. a map $f:M\rightarrow N$ of $A$-modules is a fibration
	or a weak equivalence in $\weightqAmodinduced$ if and only if $Uf$ is a fibration or a weak
	equivalence in $\weightqsymmTspectra$.  Furthermore, the model category $\weightqAmodinduced$
	is cofibrantly generated, with the following sets of
	generating cofibrations and trivial cofibrations respectively:
		\begin{eqnarray*}
	 		I_{\qleftAmodlocmaps} &= & I_{A\text{-}\modules}  = A\wedge I^{T}_{\Sigma} \\
	 											&= & \bigcup _{k\geq 0}
														\{ id\wedge i:A\wedge F_{k}^{\Sigma}(Y_{+})\rightarrow
														A\wedge F_{k}^{\Sigma}((\Delta _{n}^{U})_{+}) \mid U\in (\smoothS), n\geq 0\}\\
			\widetilde{J_{\qleftAmodlocmaps}} &= & \{id\wedge j:A\wedge X\rightarrow A\wedge Y\}
		\end{eqnarray*}
	where $j:X\rightarrow Y$ satisfies the following conditions:
	\begin{enumerate}
		\item	$j$ is an inclusion of $I^{T}_{\Sigma}$-complexes in $\weightqsymmTspectra$.
		\item	$j$ is a $\qleftsymmlocmaps$-local equivalence in $\motivicsymmTspectra$.
		\item	the size of $Y$ as an $I^{T}_{\Sigma}$-complex is less than $\kappa$, 
					where $\kappa$ is the regular cardinal defined by Hirschhorn in \cite[definition 4.5.3]{MR1944041}.
	\end{enumerate}
\end{prop}
\begin{proof}
	Using a result of D. Kan (see theorem 11.3.2 in \cite{MR1944041}), we have
	that it is enough to prove that the following conditions hold:
		\begin{enumerate}
			\item \label{prop.3.5.smash-induces-weightq-modelstructures.a}  The domains 
						of $I_{\qleftAmodlocmaps}$ (respectively $J_{\qleftAmodlocmaps}$)
						are small relative to the $I_{\qleftAmodlocmaps}$-cells 
						(respectively $J_{\qleftAmodlocmaps}$-cells) in the category
						of $A$-modules.
			\item \label{prop.3.5.smash-induces-weightq-modelstructures.b}  $U$ maps relative 
						$J_{\qleftAmodlocmaps}$-cell complexes to weak equivalences in $\weightqsymmTspectra$.
		\end{enumerate}
		
	(\ref{prop.3.5.smash-induces-weightq-modelstructures.a}):  By adjointness it suffices to check that the domains of
	$I_{\qleftsymmlocmaps}$ (respectively $J_{\qleftsymmlocmaps}$)
	are small relative to the $I_{\qleftAmodlocmaps}$-cells (respectively $J_{\qleftAmodlocmaps}$-cells)
	in $\symmTspectra$.  Theorem \ref{thm.3.3.Lqsymmetricmodelstructures}
	implies that $\weightqsymmTspectra$ is in particular a cofibrantly generated model category
	with the sets $I_{\qleftsymmlocmaps}$ and $J_{\qleftsymmlocmaps}$ as generating cofibrations
	and trivial cofibrations, therefore by \cite[proposition 2.1.16]{MR1650134} it only remains
	to show that all the maps in $I_{\qleftAmodlocmaps}$-cells (respectively $J_{\qleftAmodlocmaps}$-cells)
	are cofibrations (respectively trivial cofibrations) in $\weightqsymmTspectra$.
	
	Since $A$ is in particular cofibrant in $\motivicsymmTspectra$ and the cofibrations
	in $\motivicsymmTspectra$ and $\weightqsymmTspectra$ are identical,
	proposition \ref{prop.2.6.symmTspectra-monoidalmodelcategory}
	implies that all the maps in $I_{\qleftAmodlocmaps}$ are cofibrations  in $\weightqsymmTspectra$.
	However, the class of cofibrations is closed
	under coproducts and filtered colimits, and the limits and colimits in the category of $A$-modules
	are computed in $\symmTspectra$, hence
	all the maps in $I_{\qleftAmodlocmaps}$-cells
	are cofibrations in $\weightqsymmTspectra$.

	By hypothesis $A$ 
	is cofibrant in $\zeroslicesymmTspectra$, and
	every map $j$ in $J_{\qleftsymmlocmaps}$ is clearly a trivial cofibration
	in $\weightqsymmTspectra$.  Since $\qminusoneslicesymmTspectra$ is
	a right Bousfield localization with respect to $\weightqsymmTspectra$,
	we have that every map $j$ in $J_{\qleftsymmlocmaps}$
	is also a trivial cofibration in $\qminusoneslicesymmTspectra$.
	Therefore, theorem \ref{thm.3.4.smash-Quillenbifunctor-SpxSq----->Sqplusp} 
	implies that all the maps in $J_{\qleftAmodlocmaps}$ are trivial cofibrations
	in $\qminusoneslicesymmTspectra$, and since $\qminusoneslicesymmTspectra$ is a right Bousfield
	localization with respect to $\weightqsymmTspectra$; we get that all the maps in $J_{\qleftAmodlocmaps}$
	are also trivial cofibrations in $\weightqsymmTspectra$.   
	Finally, since the class of trivial cofibrations is closed
	under coproducts and filtered colimits, and the limits and colimits in the category of $A$-modules
	are computed in $\symmTspectra$,
	we have that all the maps in $J_{\qleftAmodlocmaps}$-cells
	are also trivial cofibrations in $\weightqsymmTspectra$.
	
	(\ref{prop.3.5.smash-induces-weightq-modelstructures.b}):  We have shown that every map
	in $J_{\qleftAmodlocmaps}$-cells is a trivial cofibration in $\weightqsymmTspectra$.
	In particular,
	every relative $J_{\qleftAmodlocmaps}$-cell complex is a weak equivalence in $\weightqsymmTspectra$,
	as we wanted.
\end{proof}

\begin{rmk}
		\label{rmk.3.5.usualargumentnotgood}
	Notice that we can not use the same argument as in 
	theorem \ref{thm.2.8.motivicsymmetricA-modules}
	to construct the model structure $\weightqAmodinduced$, since
	the model category $\weightqsymmTspectra$ is not a symmetric monoidal
	model category, i.e. the monoidal structure on symmetric $T$-spectra
	is not compatible with the model structure on $\weightqsymmTspectra$.
	Therefore, the hypothesis of $A$ being cofibrant in $\zeroslicesymmTspectra$
	is really necessary.
\end{rmk}

\begin{lem}
		\label{lemma.3.5.inducedweightq-simplicial}
	Fix $q\in \mathbb Z$.  Let $A$ be a cofibrant ring spectrum in $\motivicsymmTspectra$,
	which is also cofibrant in $\zeroslicesymmTspectra$.  Then
	the model category $\weightqAmodinduced$ described in
	proposition \ref{prop.3.5.smash-induces-weightq-modelstructures} is simplicial.
\end{lem}
\begin{proof}
	Since the cotensor objects $N^{K}$ for the simplicial structure
	are identical in $\weightqAmodinduced$ and $\weightqsymmTspectra$, the results follows from 
	proposition
	\ref{prop.3.5.smash-induces-weightq-modelstructures} and theorem
	\ref{thm.3.3.Lqsymmetricmodelstructures} which implies in particular that $\weightqsymmTspectra$
	is a simplicial model category.
\end{proof}

\begin{thm}
		\label{thm.3.5.weightq-induced==leftBousfieldlocalization}
	Fix $q\in \mathbb Z$.  Let $A$ be a cofibrant ring spectrum in $\motivicsymmTspectra$,
	which is also cofibrant in $\zeroslicesymmTspectra$.  Then
	the model structures $\weightqAmod$ (see theorem \ref{thm.3.5.Lqsymmetricmodelstructures})
	and $\weightqAmodinduced$ (see proposition \ref{prop.3.5.smash-induces-weightq-modelstructures})
	on the category of $A$-modules
	are identical.
\end{thm}
\begin{proof}
	Theorem \ref{thm.3.5.Lqsymmetricmodelstructures} and proposition \ref{prop.3.5.smash-induces-weightq-modelstructures}
	imply that both $\weightqAmod$ and $\weightqAmodinduced$ have
		$$\bigcup _{k\geq 0}
														\{ id\wedge i:A\wedge F_{k}^{\Sigma}(Y_{+})\rightarrow
														A\wedge F_{k}^{\Sigma}((\Delta _{n}^{U})_{+}) \mid U\in (\smoothS), n\geq 0\}
		$$
	as set of generating cofibrations.  Hence the cofibrations in $\weightqAmod$ and $\weightqAmodinduced$
	are exactly the same.  It suffices to check that the weak equivalences in both model structures are identical.
	
	However, theorem \ref{thm.3.3.Lqsymmetricmodelstructures} and
	lemma \ref{prop.3.5.smash-induces-weightq-modelstructures} imply that $\weightqAmod$ and $\weightqAmodinduced$
	are both simplicial model categories.  Therefore, corollary
	\ref{cor.1.1.4.detect-weak-equiv.2}(\ref{cor.1.1.4.detect-weak-equiv.2.b}) implies that
	it is enough to show  that the fibrant objects in $\weightqAmod$ and $\weightqAmodinduced$ coincide.
	But this follows directly from propositions \ref{prop.3.5.Z-AmodLq-local-iff-UZ-symmLq-local} and
	\ref{prop.3.5.smash-induces-weightq-modelstructures}.
\end{proof}

\begin{thm}
		\label{thm.3.5.inheritmodelstructures2}
	Fix $q\in \mathbb Z$.  Let $A$ be a cofibrant ring spectrum in $\motivicsymmTspectra$,
	which is also cofibrant in $\zeroslicesymmTspectra$, and let $M$ be an arbitrary $A$-module.
	Then the solid arrows in the following commutative diagram:
		\begin{equation}
					\label{thm.3.5.inheritingmodulestructures2.a}
			\begin{array}{c}
				\xymatrix@C=1.4pc{Q_{\Sigma}R_{m}M \ar[dd]_-{Q_{\Sigma}^{R_{m}M}}&& 
									\ar@{-->}[ll]^-{Q_{\Sigma}(R_{m}^{M})} Q_{\Sigma}M \ar@{-->}[dd]_-{Q_{\Sigma}^{M}}&& 
									\ar@{-->}[ll]^-{Q_{\Sigma}(Q_{m}^{M})} Q_{\Sigma}Q_{m}M \ar@{-->}[dd]_-{Q_{\Sigma}^{Q_{m}M}}
									\ar@{-->}[rr]_-{Q_{\Sigma}(W_{q}^{m,Q_{m}M})} && Q_{\Sigma}W_{q}^{m}Q_{m}M 
									\ar@{-->}[dd]_-{Q_{\Sigma}^{W_{q}^{m}Q_{m}M}}\\
									&&&&&& \\
									R_{m}M && \ar[ll]^-{R_{m}^{M}} M && \ar[ll]^-{Q_{m}^{M}} Q_{m}M 
									\ar[rr]_-{W_{q}^{m,Q_{m}M}} && W_{q}^{m}Q_{m}M}
			\end{array}
		\end{equation}
	induce a natural equivalence between the functors:
		\begin{equation}
					\label{thm.3.5.inheritingmodulestructures2.b}
			\begin{array}{c}
				\xymatrix{& \weightqAmodstablehomotopy \ar[dr]^-{UW_{q}^{m}}&\\
												 \Amodstablehomotopy
												 \ar[ur]^-{Q_{m}} \ar[dr]_-{UR_{m}}&& \weightqsymmstablehomotopy \\
												& \symmstablehomotopy \ar[ur]_-{Q_{\Sigma}}&}
			\end{array}
		\end{equation}
\end{thm}
\begin{proof}
	It suffices to show that all the maps $W_{q}^{m,Q_{m}M}$, $Q_{m}^{M}$, $R_{m}^{M}$ and
	$Q_{\Sigma}^{R_{m}M}$ are weak equivalences in $\weightqsymmTspectra$.  Proposition  
	\ref{prop.3.3.Q-cofibrant-replacement-all-symmetricL<q} implies that 
	$Q_{\Sigma}^{R_{m}M}$ is a weak equivalence in $\weightqsymmTspectra$.
	
	Since $A$ is cofibrant in $\zeroslicesymmTspectra$, theorem \ref{thm.3.5.weightq-induced==leftBousfieldlocalization} 
	and proposition \ref{prop.3.5.smash-induces-weightq-modelstructures}
	imply that it is enough to show that $W_{q}^{m,Q_{m}M}$, $Q_{m}^{M}$ and $R_{m}^{M}$
	are weak equivalences in $\weightqAmod$.
	By construction (see definition \ref{def.3.5.stable-weightqsymm-replacementfunctors})
	$W_{q}^{m,Q_{m}M}$ is a weak equivalence in $\weightqAmod$,
	and proposition \ref{prop.3.5.Q-cofibrant-replacement-all-symmetricL<q} implies that 
	$Q_{m}^{M}$ is a weak equivalence in $\weightqAmod$.
	Finally, by construction (see definition \ref{def.3.5.fibrantstable-replacementfunctor})
	$R_{m}^{M}$ is a weak equivalence in $\motivicAmod$,
	and \cite[proposition 3.1.5]{MR1944041} implies that $R_{m}^{M}$ is also a weak equivalence in $\weightqAmod$.
	This finishes the proof.
\end{proof}

\begin{thm}
		\label{thm.3.5.weightqinvarianceofcoefficients}
	Fix $q\in \mathbb Z$.  Let $f:A\rightarrow A'$ be a map 
	between cofibrant ring spectra in $\motivicsymmTspectra$,
	which is compatible with the ring structures.
	Assume that one of the following conditions holds:
		\begin{enumerate}
			\item	\label{thm.3.5.weightqinvarianceofcoefficients.a}  $f$ is a weak equivalence in $\motivicsymmTspectra$
								(equivalently in $\motivicAmod$).
			\item \label{thm.3.5.weightqinvarianceofcoefficients.b} There exists $p\in \mathbb Z$
								such that $A$, $A'$ are both $L^{m}(<p)$-local 
								in $\motivicsymmTspectra$ and $f$ is
								a weak equivalence in $\weightpsymmTspectra$.
			\item	\label{thm.3.5.weightqinvarianceofcoefficients.c} There exists $p\in \mathbb Z$
								such that $A$, $A'$ are both $C_{eff}^{p,\Sigma}$-colocal
								in $\motivicsymmTspectra$ and $f$ is a weak equivalence in
								$\pconnectedsymmTspectra$ (equivalently in $\pconnectedAmod$).
		\end{enumerate}
	Then $f$ induces a Quillen equivalence between the weight$^{<q}$ motivic stable model structures of $A$ and $A'$ modules:
		$$\xymatrix{(A'\wedge _{A} -, U,\varphi ):\weightqAmod \ar[r] & \weightqAAmod}
		$$
\end{thm}
\begin{proof}
	(\ref{thm.3.5.weightqinvarianceofcoefficients.a}):  This is just proposition
	\ref{prop.3.5.changecoefficients-Quillenfunctor2a}.
	
	(\ref{thm.3.5.weightqinvarianceofcoefficients.b}):  Since $A$ and $A'$ are $L^{m}(<p)$-local in $\motivicsymmTspectra$,
	\cite[theorem 3.2.13(1)]{MR1944041} implies that $f$ is a weak equivalence in $\motivicsymmTspectra$.
	Therefore the result follows from proposition \ref{prop.3.5.changecoefficients-Quillenfunctor2a}.
	
	(\ref{thm.3.5.weightqinvarianceofcoefficients.c}):  Since $A$ and $A'$ are  $C_{eff}^{p,\Sigma}$-colocal
	in $\motivicsymmTspectra$, using \cite[theorem 3.2.13(2)]{MR1944041} we have that $f$ is a weak equivalence in $\motivicsymmTspectra$.
	Thus, the result follows from proposition \ref{prop.3.5.changecoefficients-Quillenfunctor2a}.
\end{proof}

\begin{defi}
		\label{def.3.5.symmetricSq-colocal-generators}
	For every $q\in \mathbb Z$, we consider the following set
	of $A$-modules
		$$\Amodqslicegenerators =\{ A\wedge \symmgeneratorNRS \in C^{m} | s-n=q \}\subseteq C_{eff}^{q,m}
		$$
	(see  
	definition \ref{def.3.3.symmetricSq-colocal-generators}).
\end{defi}

\begin{thm}
		\label{thm.3.5.symmetricSq-modelcats}
	Fix $q\in \mathbb Z$.  Then the right Bousfield
	localization of the model category $\weightqplusoneAmod$
	with respect to the $\Amodqslicegenerators$-colocal equivalences
	exists.  This new model structure will be called \emph{$q$-slice
	motivic stable}.
	$\qsliceAmod$ will denote the category of $A$-modules
	equipped with the $q$-slice motivic stable model structure, and 
	$\qsliceAmodstablehomotopy$ will denote its associated homotopy category.
	Furthermore, the $q$-slice motivic stable model structure is right proper and simplicial.
\end{thm}

\begin{rmk}
		\label{rmk.3.5.usualargumentsnotgood2}
	Notice that
	we can not use the adjuntion $(A\wedge -,U,\varphi ):\qslicesymmTspectra \rightarrow \qsliceAmod$ 
	for the construction of
	$\qsliceAmod$, since we do not know if the model structure for $\qslicesymmTspectra$
	is cofibrantly generated.
\end{rmk}

\begin{defi}
		\label{def.3.5.Pqsigma-cofibrant-replacement}
	Fix $q\in \mathbb Z$.  Let $P_{q}^{m}$ denote a cofibrant replacement functor in 
	$\qsliceAmod$; such that for every $A$-module $M$, the natural map
		$$\xymatrix{P_{q}^{m}M \ar[r]^-{P_{q}^{m,M}}& M}
		$$ 
	is a trivial fibration in $\qsliceAmod$, and
	$P_{q}^{m}M$ is always a $\Amodqslicegenerators$-colocal $A$-module
	in $\weightqplusoneAmod$.
\end{defi}

\begin{prop}
		\label{prop.3.5.Wsigmaqplusone-fibrant-replacement-all-Sq}
	Fix $q\in \mathbb Z$.  Then $W^{m}_{q+1}$ is also a fibrant
	replacement functor in $\qsliceAmod$
	(see definition \ref{def.3.5.stable-weightqsymm-replacementfunctors}),
	and for every $A$-module $M$ the natural map
		$$\xymatrix{M\ar[rr]^-{W^{m,M}_{q+1}}&& W^{m}_{q+1}M}
		$$
	is a trivial cofibration in $\qsliceAmod$.
\end{prop}

\begin{cor}
		\label{cor.3.5.classifying-Sq-symmetriccolocal-equivs.b}
	Fix $q\in \mathbb Z$ and let $f:M\rightarrow N$ be a map of $A$-modules.
	Then $f$ is a $\Amodqslicegenerators$-colocal equivalence
	in $\weightqplusoneAmod$ if and only if
		$$\xymatrix{W_{q+1}^{m}M \ar[rr]^-{W_{q+1}^{m}f}&& W_{q+1}^{m}N}
		$$
	is a $C_{eff}^{q,m}$-colocal equivalence in $\motivicAmod$.
\end{cor}

\begin{prop}
		\label{prop.3.5.Suspension=>qslicesymm-Quillen-equiv}
	Fix $q\in \mathbb Z$.  Then the adjunction
		$$\xymatrix{(-\wedge S^{1},\Omega _{S^{1}},\varphi):\qsliceAmod \ar[r]& 
								\qsliceAmod}
		$$
	is a Quillen equivalence, and $\qsliceAmodstablehomotopy$
	has the structure of a triangulated category.
\end{prop}

\begin{prop}
		\label{prop.3.5.symmPqcofibrant-replacement=>triangulatedfunctor}
	Fix $q\in \mathbb Z$.  Then we have the following adjunction
		$$\xymatrix{(P_{q}^{m}, W_{q+1}^{m}, \varphi) :\qsliceAmodstablehomotopy \ar[r]& \weightqplusoneAmodstablehomotopy}
		$$
	between exact functors of triangulated categories.
\end{prop}

\begin{prop}
		\label{prop.3.5.symmetricq-connected--->q-slice===leftQuilllenfunctor}
	Fix $q\in \mathbb Z$.  Then the identity functor
		$$\xymatrix{id:\qsliceAmod \ar[r]& \qconnectedAmod}
		$$
	is a right Quillen functor, and it induces the following adjunction
		$$\xymatrix{(C_{q}^{m},W_{q+1}^{m},\varphi ):\qconnectedAmodstablehomotopy \ar[r]& \qsliceAmodstablehomotopy}
		$$
	of exact functors between triangulated categories.
\end{prop}

\begin{lem}
		\label{lem.3.5.A-symmSqcofibrant====>AsymmLq+1trivially-cofibrant}
	Fix $q\in \mathbb Z$, and let $M$ be a cofibrant $A$-module
	in $\qsliceAmod$.  Then the map $\ast \rightarrow M$
	is a trivial cofibration in $\weightqAmod$.
\end{lem}

\begin{thm}
		\label{thm.3.5.symmetrization-qslice-Quillen-equiv}
	Fix $q\in \mathbb Z$.  Then the adjunction
		$$\xymatrix{(A\wedge -,U,\varphi):\qslicesymmTspectra \ar[r]& \qsliceAmod}
		$$
	given by the free $A$-module and the forgetful functors is a
	Quillen adjunction, and it induces
	an adjunction
		$$\xymatrix{(A\wedge P_{q}^{\Sigma}-,UW_{q+1}^{m},\varphi):\qslicesymmstablehomotopy \ar[r]& \qsliceAmodstablehomotopy}
		$$
	of exact funtors between triangulated categories.
\end{thm}

\begin{prop}
		\label{prop.3.5.homotopycoherence==>liftings-qslice}
	Fix $q\in \mathbb Z$.  We have the following
	commutative diagram of left Quillen functors:
		\begin{equation}
				\label{diagram3.5.homotopycoherence==>liftings-qslice.a}
					\begin{array}{c}
						\xymatrix{\qconnectedsymmTspectra \ar[r]^-{A\wedge -} \ar[d]_-{id}& \qconnectedAmod \ar[d]^-{id}\\
								\qslicesymmTspectra \ar[r]_-{A\wedge -}& \qsliceAmod}
					\end{array}
		\end{equation}
	and the following associated commutative diagrams of homotopy categories:
		\begin{equation}
	 		\label{diagram3.5.homotopycoherence==>liftings-qslice.b}
	 			\begin{array}{c}
	 				\xymatrix{& \qconnectedAmodstablehomotopy \ar[dr]^-{C_{q}^{m}}&\\
												 \qconnectedsymmstablehomotopy
												\ar[rr]^{A\wedge C_{q}^{\Sigma}-} \ar[ur]^-{A\wedge C_{q}^{\Sigma}-} 
												\ar[dr]_-{C_{q}^{\Sigma}}&& \qsliceAmodstablehomotopy \\
												& \qslicesymmstablehomotopy \ar[ur]_-{A\wedge P_{q}^{\Sigma}-}&}
				\end{array}
		\end{equation}
		
		\begin{equation}
	 		\label{diagram3.5.homotopycoherence==>liftings-qslice.c}
	 			\begin{array}{c}
	 				\xymatrix{& \qslicesymmstablehomotopy \ar[dr]^-{W_{q+1}^{\Sigma}}&\\
												 \qsliceAmodstablehomotopy
												\ar[rr]^{UW_{q+1}^{m}} \ar[ur]^-{UW_{q+1}^{m}} \ar[dr]_-{W_{q+1}^{m}}&& \qconnectedsymmstablehomotopy \\
												& \qconnectedAmodstablehomotopy \ar[ur]_-{UR_{m}}&}
				\end{array}
		\end{equation}
\end{prop}

\begin{thm}
		\label{thm.3.5.smash-Quillenbifunctor-SpxSq----->Sqplusp-coefss}
	Fix $p,q \in \mathbb Z$.  Let $A$ be a commutative cofibrant ring spectrum in $\motivicsymmTspectra$.
	Then the symmetric monoidal structure for the category of $A$-modules, induces
	the following Quillen bifunctor in the sense of Hovey (see definition \ref{def.Quillen-bifunct}).
		$$\xymatrix{-\wedge _{A}-:
									\psliceAmod \times \qsliceAmod \ar[r]& \pplusqsliceAmod}
		$$
\end{thm}
\begin{proof}
	The proof is similar to the one given for theorem \ref{thm.3.4.smash-Quillenbifunctor-SpxSq----->Sqplusp}.
	We leave the details to the reader.
\end{proof}

	If the ring $A$ is not commutative, then it needs to satisfy some
	additional conditions
	in order to get a weaker version of
	the previous result (see theorem \ref{thm.3.5.notcommutativesmash-Quillenbifunctor-SpxSq----->Sqplusp-coefss}).

\begin{lem}
		\label{lemma.3.5.changecoefficients-slice-Quillenfunctor}
	Fix $q\in \mathbb Z$.  Let $f:A\rightarrow A'$ be a map 
	between cofibrant ring spectra in $\motivicsymmTspectra$,
	which is compatible with the ring structures.
	Then the adjunction:
		$$(A'\wedge _{A}-, U,\varphi ):\qsliceAmod \rightarrow \qsliceAAmod
		$$
	is a Quillen adjunction.  
\end{lem}
\begin{proof}
	Proposition \ref{prop.3.5.changecoefficients-Quillenfunctor2} 
	implies that $U:\weightqplusoneAAmod \rightarrow \weightqplusoneAmod$ is a right Quillen functor.
	Consider the following commutative diagram of right Quillen functors:
		$$\xymatrix{\weightqplusoneAAmod \ar[r]^-{U} \ar[d]_-{id} & \weightqplusoneAmod \ar[d]^-{id}\\
								\qsliceAAmod \ar@{-->}[r]_-{U}& \qsliceAmod}
		$$
	then the universal property of right Bousfield localizations together with proposition
	\ref{prop.2.8.enriched-freemodule-adjunction2} imply that the dotted arrow $U$ is a right Quillen functor.
\end{proof}

\begin{prop}
		\label{prop.3.5.changecoefficients-slice-Quillenequiv}  Fix $q\in \mathbb Z$.  Let $f:A\rightarrow A'$ be
	a map between cofibrant ring spectra in $\motivicsymmTspectra$, which is
	compatible with the ring structures.
	Assume that $f$ is a weak equivalence in $\motivicsymmTspectra$.  Then the adjunction
		$$(A'\wedge _{A}-, U,\varphi ):\qsliceAmod \rightarrow \qsliceAAmod
		$$
	is a Quillen equivalence.
\end{prop}
\begin{proof}
	We have shown in lemma 
	\ref{lemma.3.5.changecoefficients-slice-Quillenfunctor} that
		$$(A'\wedge _{A} -, U,\varphi ):\qsliceAmod \rightarrow \qsliceAAmod$$
	is a Quillen adjunction.
		
	Now 
	let $\eta$, $\epsilon$ denote the unit and counit of the
	adjunction $(A' \wedge _{A}-, U, \varphi)$.
	By corollary 1.3.16(c) in \cite{MR1650134}, it suffices to check
	that the following conditions hold:
		\begin{enumerate}
			\item \label{prop.3.5.changecoefficients-slice-Quillenequiv.a}  For every cofibrant
						$A$-module $M$ in $\qsliceAmod$, the following composition
							$$\xymatrix{M\cong A\wedge _{A}M \ar[rr]^-{\eta _{M}=f\wedge _{A}id} && A'\wedge _{A}M
								\ar[rr]^-{W^{m',A'\wedge _{A}M}_{q+1}}&& W_{q+1}^{m'}(A'\wedge _{A}M)}
							$$
						is a weak equivalence in $\qsliceAmod$, where $W^{m'}_{q+1}$ denotes a fibrant replacement functor
						in $\qsliceAAmod$ (see proposition \ref{prop.3.5.Wsigmaqplusone-fibrant-replacement-all-Sq}).
			\item \label{prop.3.5.changecoefficients-slice-Quillenequiv.b}  $U$ reflects
						weak equivalences between fibrant objects in $\qsliceAAmod$.
		\end{enumerate}
		
	(\ref{prop.3.5.changecoefficients-slice-Quillenequiv.a}):  Since $id:\qsliceAmod \rightarrow \weightqplusoneAmod$
	is a left Quillen functor, we have that $M$ is also cofibrant in $\weightqplusoneAmod$.  Hence, theorem
	\ref{thm.3.5.weightqinvarianceofcoefficients}(\ref{thm.3.5.weightqinvarianceofcoefficients.a}) 
	implies that $W^{m',A'\wedge _{A}M}_{q+1}\circ \eta _{M}$ 
	is a weak equivalence in $\weightqplusoneAmod$.  Finally, by \cite[proposition 3.1.5]{MR1944041}
	we have that $W^{m',A'\wedge _{A}M}_{q+1}\circ \eta _{M}$ is a weak equivalence in $\qsliceAmod$,
	as we wanted.

	(\ref{prop.3.5.changecoefficients-slice-Quillenequiv.b}): Let $g:M\rightarrow N$ be a map
	between fibrant $A'$-modules in $\qsliceAAmod$, such that $Ug$ is a weak equivalence in $\qsliceAmod$.
	
	Fix $\symmgeneratorNRS \in S^{\Sigma}(q)$ (see definition \ref{def.3.3.symmetricSq-colocal-generators}).
	Using the enriched adjunctions of proposition \ref{prop.2.8.enriched-freemodule-adjunction2},
	we get the following commutative diagram of simplicial sets where the vertical arrows are 
	isomorphisms
		$$\xymatrix@C=-4pc{Map_{A'\text{-}\modules}(A'\wedge \symmgeneratorNRS ,M) \ar[dd]_-{\cong} \ar[dr]^-{g_{\ast}}&\\ 
											& Map_{A'\text{-}\modules}(A'\wedge \symmgeneratorNRS ,N) \ar[dd]^-{\cong}\\
											Map_{A\text{-}\modules}(A\wedge \symmgeneratorNRS , UM) \ar[dr]_-{(Ug)_{\ast}}&\\ 
											& Map_{A\text{-}\modules}(A\wedge \symmgeneratorNRS , UN)}
		$$
	Now $M$ and $N$ are both fibrant in $\weightqplusoneAAmod$
	(this follows from proposition \ref{prop.3.5.Wsigmaqplusone-fibrant-replacement-all-Sq}), 
	hence proposition \ref{prop.3.5.changecoefficients-Quillenfunctor2}
	implies that $UM$ and $UN$ are also fibrant in $\weightqplusoneAmod$.  Therefore, the bottom
	row in the diagram above is a weak equivalence of simplicial sets, since by hypothesis $Ug$
	is a weak equivalence in $\qsliceAmod$.  Finally, by the two out of three property for weak
	equivalences we get that the top row is also a weak equivalence of simplicial sets, and this
	implies that $g$ is a weak equivalence in $\qsliceAAmod$, since $M$ and $N$ are both fibrant
	in $\weightqplusoneAAmod$.	
\end{proof}

\begin{lem}
		\label{lemma.3.5.forgetful-detects-reflects-sliceequivs}
	Fix $q\in \mathbb Z$.  Let $f:A\rightarrow A'$ be a map 
	between cofibrant ring spectra in $\motivicsymmTspectra$,
	which is compatible with the ring structures.
	Assume that $A$ and $A'$ are cofibrant in $\zeroslicesymmTspectra$.
	Then $w:M\rightarrow M'$ is a weak equivalence in $\qsliceAAmod$  if
	and only if $Uw$ is a weak equivalence in $\qsliceAmod$.
\end{lem}
\begin{proof}
	Let $W^{m}_{q+1}$, $W^{m'}_{q+1}$ denote fibrant replacement functors
	in $\weightqplusoneAmod$ and $\weightqplusoneAAmod$ respectively, and let $N$ be an arbitrary $A'$-module.  We have the following
	commutative diagram in $\weightqplusoneAmod$:
		$$\xymatrix{N  \ar[d]_-{W^{m',N}_{q+1}} \ar[rr]^-{W^{m,N}_{q+1}}&& 
								W^{m}_{q+1}N \ar[d]^-{W^{m}_{q+1}(W^{m',N}_{q+1})}\\
								W^{m'}_{q+1}N \ar[rr]_-{W^{m,W^{m'}_{q+1}N}_{q+1}} && W^{m}_{q+1}W^{m'}_{q+1}N }
		$$
	since $A$, $A'$ are both cofibrant in $\zeroslicesymmTspectra$,
	theorem \ref{thm.3.5.weightq-induced==leftBousfieldlocalization} and
	proposition \ref{prop.3.5.smash-induces-weightq-modelstructures} 
	imply that all the maps in the diagram above are weak equivalences in $\weightqplusoneAmod$.  
	  	
	Now fix $\generatorNRS \in S^{m}(q)$ (see definition \ref{def.3.5.symmetricSq-colocal-generators}).  Using
	the naturality of the diagram above together with 
	proposition \ref{prop.2.8.enriched-freemodule-adjunction2}, we get
	the following commutative diagram of simplicial sets:	
		$$\xymatrix@C=-8pc{Map_{A'\text{-}\modules}(A'\wedge \generatorNRS ,W^{m'}_{q+1}M) 
								\ar[dr]^-{(W^{m'}_{q+1}w)_{\ast}} \ar[dd]_-{\cong}&\\  
								& Map_{A'\text{-}\modules}(A'\wedge \generatorNRS ,W^{m'}_{q+1}M')\ar[dd]^-{\cong}\\
								Map_{A\text{-}\modules}(A\wedge \generatorNRS ,UW^{m'}_{q+1}M) \ar[dr]^-{(UW^{m'}_{q+1}w)_{\ast}}
								\ar[dd]_-{(W^{m,UW^{m'}_{q+1}M}_{q+1})_{\ast}}&\\  
								& Map_{A\text{-}\modules}(A\wedge \generatorNRS ,UW^{m'}_{q+1}M')\ar[dd]^-{(W^{m,UW^{m'}_{q+1}M'}_{q+1})_{\ast}}\\
								Map_{A\text{-}\modules}(A\wedge \generatorNRS ,W^{m}_{q+1}UW^{m'}_{q+1}M) \ar[dr]^-{(W^{m}_{q+1}(UW^{m'}_{q+1}w))_{\ast}}
								\ar@{<-}[dd]_-{(UW^{m}_{q+1}(W^{m',M}_{q+1}))_{\ast}}&\\  
								& Map_{A\text{-}\modules}(A\wedge \generatorNRS ,W^{m}_{q+1}UW^{m'}_{q+1}M')
								\ar@{<-}[dd]^-{(UW^{m}_{q+1}(W^{m',M'}_{q+1}))_{\ast}}\\
								Map_{A\text{-}\modules}(A\wedge \generatorNRS ,UW^{m}_{q+1}M) \ar[dr]^-{(UW^{m}_{q+1}w)_{\ast}}&\\  
								& Map_{A\text{-}\modules}(A\wedge \generatorNRS ,UW^{m}_{q+1}M')}
		$$
	where the top vertical arrows are isomorphisms of simplicial sets.
	But $\weightqplusoneAmod$, $\weightqplusoneAAmod$ are simplicial
	model categories (see theorem 
	\ref{thm.3.5.Lqsymmetricmodelstructures}) 
	and the natural maps $W^{m,UW^{m'}_{q+1}M}_{q+1}$, $UW^{m}_{q+1}(W^{m',M}_{q+1})$, 
	$W^{m,UW^{m',M'}_{q+1}}_{q+1}$ and $UW^{m}_{q+1}(W^{m',M'}_{q+1})$ 
	are all weak equivalences between fibrant objects, 
	thus by Ken Brown's lemma (see lemma \ref{lemma1.1.KenBrown2})
	all the vertical arrows are weak equivalences
	of simplicial sets. 
	
	Therefore, the top row is a weak equivalence of simplicial sets
	if and only if the bottom row is a weak equivalence of simplicial sets.
	This proves the claim.
\end{proof}

\begin{prop}
		\label{prop.3.5.Acolocal=>qslicelocalcofibrations=>cofibrations}
	Let $A$ be a cofibrant ring spectrum in $\motivicsymmTspectra$, which is also
	cofibrant in $\zeroslicesymmTspectra$.  Then for every $q\in \mathbb Z$, and for every cofibration $f:M\rightarrow N$
	in $\qsliceAmod$ we have
	that $f$ is also a cofibration in $\qslicesymmTspectra$.
\end{prop}
\begin{proof}
	Let (see theorem \ref{thm.3.5.Lqsymmetricmodelstructures})
		\begin{eqnarray*}
			\overline{\Lambda (K)} = J_{\qplusoneleftAmodlocmaps}\; \cup \; \{ A\wedge \symmgeneratorNRS \otimes \partial \Delta ^{k}
								\rightarrow  \\ 
								A\wedge \symmgeneratorNRS \otimes  \Delta ^{k} \mid s-n=q,  k\geq 0\}
		\end{eqnarray*}

	Since $\weightqplusoneAmod$ is in particular a simplicial model category
	(see theorem \ref{thm.3.5.Lqsymmetricmodelstructures}), 
	using definitions 5.2.1, 16.3.1 and propositions 5.3.6, 16.1.3 
	in \cite{MR1944041}, we have that $f$ is a retract of a cofibration $g:M\rightarrow O$ in $\weightqplusoneAmod$
	for which there is a weak equivalence $h:O \rightarrow P$ in $\weightqplusoneAmod$
	such that the composition $h\circ g$ is a relative $\overline{\Lambda (K)}$-cell complex.
	
	It is clear that it is enough to check that $g$ is a cofibration in $\qslicesymmTspectra$.
	Now, using lemma 5.3.4  
	in \cite{MR1944041}, we have that this follows from:
		\begin{enumerate}
			\item	\label{prop.3.5.Acolocal=>qslicelocalcofibrations=>cofibrations.a}  $g$ is a cofibration in $\weightqplusonesymmTspectra$.
			\item	\label{prop.3.5.Acolocal=>qslicelocalcofibrations=>cofibrations.b}	$h$ is a weak equivalence in $\weightqplusonesymmTspectra$.
			\item	\label{prop.3.5.Acolocal=>qslicelocalcofibrations=>cofibrations.c}  $h\circ g$ is a cofibration in $\qslicesymmTspectra$.
		\end{enumerate}
	
	(\ref{prop.3.5.Acolocal=>colocalcofibrations=>cofibrations.a}):  Since $\weightqplusoneAmod$ is a left
	Bousfield localization with respect to $\motivicAmod$, we have that the cofibrations are exactly the same
	in both model structures.  Hence $g$ is a cofibration in $\motivicAmod$, and
	proposition \ref{prop.2.8.motivicAmod-cofibrations=>motivic-cofibrations} implies that
	$g$ is also a cofibration in $\motivicsymmTspectra$.  But  $\weightqplusonesymmTspectra$ is
	a left Bousfield localization with respect to $\motivicsymmTspectra$, therefore
	$g$ is a cofibration in $\weightqplusonesymmTspectra$.
	
	(\ref{prop.3.5.Acolocal=>colocalcofibrations=>cofibrations.b}):  Since $A$ is cofibrant in $\zeroslicesymmTspectra$,
	theorem \ref{thm.3.5.weightq-induced==leftBousfieldlocalization} 
	and proposition \ref{prop.3.5.smash-induces-weightq-modelstructures} 
	imply that $h$ is a weak equivalence in $\weightqplusonesymmTspectra$.
	
	(\ref{prop.3.5.Acolocal=>colocalcofibrations=>cofibrations.c}): Let $\mathcal C$ denote the class of cofibrations
	in $\qslicesymmTspectra$.  Theorem
	\ref{thm.3.5.Lqsymmetricmodelstructures} implies that $J_{\qplusoneleftAmodlocmaps}$
	is a set of generating trivial cofibrations for $\weightqplusoneAmod$, and
	since $A$ is cofibrant in $\zeroslicesymmTspectra$, theorem
	\ref{thm.3.5.weightq-induced==leftBousfieldlocalization} together 
	with \ref{prop.3.5.smash-induces-weightq-modelstructures} imply that
	all the maps in $J_{\qplusoneleftAmodlocmaps}$ are weak equivalences
	in $\weightqplusonesymmTspectra$.  
	
	Now, $\weightqplusoneAmod$ is a
	left Bousfield localization with respect to $\motivicAmod$, thus all the maps in
	$J_{\qplusoneleftAmodlocmaps}$ are cofibrations in $\motivicAmod$, and
	proposition \ref{prop.2.8.motivicAmod-cofibrations=>motivic-cofibrations} implies that 
	the maps in $J_{\qplusoneleftAmodlocmaps}$ are also
	cofibrations in $\motivicsymmTspectra$.  However, $\weightqplusonesymmTspectra$
	is a left Bousfield localization with respect to $\motivicsymmTspectra$,
	hence all the maps in $J_{\qplusoneleftAmodlocmaps}$ are cofibrations
	in $\weightqplusonesymmTspectra$.
	
	Therefore, all the maps in $J_{\qplusoneleftAmodlocmaps}$
	are trivial cofibrations in $\weightqplusonesymmTspectra$.
	But $\qslicesymmTspectra$ is a right Bousfield localization with respect to $\weightqplusonesymmTspectra$,
	hence all the maps in $J_{\qplusoneleftAmodlocmaps}$ are also trivial cofibrations in $\qslicesymmTspectra$.
	We have that in particular $J_{\qplusoneleftAmodlocmaps}$
	is contained in $\mathcal C$.
	On the other hand, by construction 
	$\ast \rightarrow \generatorNRS$ are cofibrations in $\qconnectedsymmTspectra$ for $s-n= q$,
	thus, proposition \ref{prop.3.3.symmetricq-connected--->q-slice===leftQuilllenfunctor} 
	implies that $\ast \rightarrow \generatorNRS$ are also cofibrations in $\qslicesymmTspectra$
	for $s-n=q$.
	By hypothesis the map $\ast \rightarrow A$ is 
	a cofibration in $\zeroslicesymmTspectra$,
	Then theorem \ref{thm.3.4.smash-Quillenbifunctor-SpxSq----->Sqplusp} together with
	the fact that $\qslicesymmTspectra$ is a simplicial model category
	(see theorem \ref{thm.3.3.symmetricSq-modelcats}) imply that
		\begin{eqnarray*}
			\{ A\wedge \symmgeneratorNRS \otimes \partial \Delta ^{k}
								\rightarrow 
								A\wedge \symmgeneratorNRS \otimes  \Delta ^{k} \mid \\ 
								s-n= q, k\geq 0\}
		\end{eqnarray*}
	is also contained in $\mathcal C$.
	Therefore, we have that all the maps in 
	$\overline{\Lambda(K)}$ are contained in $\mathcal C$.
	
	Finally since
	limits and colimits in $A\text{-}\modules$
	are computed in $\symmTspectra$, we have that $h\circ g$ is a relative $\mathcal C$-cell complex
	in $\symmTspectra$, 
	and since $\mathcal C$ is clearly closed under
	coproducts, pushouts and filtered colimits, we have that
	$h\circ g$ is a cofibration in $\qslicesymmTspectra$.
\end{proof}

\begin{thm}
		\label{thm.3.5.notcommutativesmash-Quillenbifunctor-SpxSq----->Sqplusp-coefss}
	Fix $p,q \in \mathbb Z$.  Let $A$ be a cofibrant ring spectrum in $\motivicsymmTspectra$,
	which is also cofibrant in $\zeroslicesymmTspectra$.  Then
	$-\wedge _{A}-$ defines a Quillen adjunction of two variables
	(see definition \ref{def.Quillen-bifunct}) from the $p$-slice motivic model structure for right $A$-modules
	and the $q$-slice motivic model structure for left $A$-modules
	to the ($p+q$)-slice motivic symmetric stable model structure:
		$$\xymatrix{-\wedge _{A}-:
									\psliceAmod _{r}\times \qsliceAmod _{l}\ar[r]& \pplusqslicesymmTspectra}
		$$
\end{thm}
\begin{proof}
	By lemma \ref{lem.cond-Quillen-bifunc}, it is enough to prove the following claim:
	
	Given a cofibration $i:N\rightarrow N'$ in $\qsliceAmod _{l}$ and a fibration
		$f:X\rightarrow Y$ in $\pplusqslicesymmTspectra$, the induced map
			$$\xymatrix{\inthomsymmTspectra (N',X)\ar[d]^-{(i^{\ast},f_{\ast})}\\
									\inthomsymmTspectra (N,X)\times _{\inthomsymmTspectra (N,Y)}\inthomsymmTspectra (N',Y)}
			$$
		is a fibration in $\psliceAmod _{r}$ which is trivial if either $i$ or $f$ is a weak
		equivalence.
		
	However, proposition \ref{prop.3.5.Acolocal=>qslicelocalcofibrations=>cofibrations}
	and lemma \ref{lemma.3.5.forgetful-detects-reflects-sliceequivs} ($\symmspherespectrum$ is cofibrant
	in $\zeroslicesymmTspectra$ by lemma \ref{lem.3.5.spherespectrum-0-connected-0-slice-cofibrant})
	imply that $i$ is also a cofibration in $\qslicesymmTspectra$,
	which is trivial if $i$ is a weak equivalence in $\qsliceAmod _{l}$.  Now,
	it follows from theorem \ref{thm.3.4.smash-Quillenbifunctor-SpxSq----->Sqplusp} 
	that $(i^{\ast},f_{\ast})$ is a fibration in $\pslicesymmTspectra$
	which is trivial if either $i$ or $f$ is a weak equivalence.  By
	lemma \ref{lemma.3.5.forgetful-detects-reflects-sliceequivs}
	we have that it suffices to check that $(i^{\ast},f_{\ast})$ is a fibration in $\psliceAmod _{r}$.
	
	By definition $\pslicesymmTspectra$ is a right Bousfield localization
	with respect to $\weightpplusonesymmTspectra$, hence the fibrations in
	both model structures coincide.   This implies that $(i^{\ast},f_{\ast})$
	is a fibration in $\weightpplusonesymmTspectra$.  Now, proposition
	\ref{prop.3.5.smash-induces-weightq-modelstructures} and theorem \ref{thm.3.5.weightq-induced==leftBousfieldlocalization}
	imply that $(i^{\ast},f_{\ast})$ is also a fibration in $\weightpplusoneAmod _{r}$
	since we are assuming that $A$ is cofibrant in $\zeroslicesymmTspectra$.
	However, by construction $\psliceAmod _{r}$ is a right Bousfield localization with respect to
	$\weightpplusoneAmod _{r}$, therefore the classes of fibrations in both model structures are identical.
	Thus $(i^{\ast},f_{\ast})$ is a fibration in $\psliceAmod _{r}$, as we wanted.
\end{proof}

\begin{lem}
		\label{lemma.3.5.changecoeffs-preserves-cofibrantequivs}
	Fix $q\in \mathbb Z$.  Let $f:A\rightarrow A'$ be a map 
	between cofibrant ring spectra in $\motivicsymmTspectra$,
	which is compatible with the ring structures.
	Assume that $A$ and $A'$ are cofibrant in $\zeroslicesymmTspectra$.
	Furthermore, assume that $A'$ is also cofibrant in $\motivicAmod$.
	If
	$f$ is a weak equivalence in $\zeroslicesymmTspectra$, then
	for every cofibrant $A$-module $M$ in $\qsliceAmod$, the induced map
		$$\xymatrix{M\cong A\wedge _{A}M \ar[rr]^-{f\wedge _{A}id} && A'\wedge _{A}M}
		$$
	is a weak equivalence in $\qsliceAmod$.
\end{lem}
\begin{proof}
	Lemma \ref{lem.3.5.spherespectrum-0-connected-0-slice-cofibrant}
	implies that $\symmspherespectrum$ is cofibrant in $\zeroslicesymmTspectra$
	and $A$ is by hypothesis cofibrant in $\zeroslicesymmTspectra$, thus
	by lemma \ref{lemma.3.5.forgetful-detects-reflects-sliceequivs}
	it suffices to check that $f\wedge _{A}id$ is a weak equivalence in $\qslicesymmTspectra$.
	
	Using lemma \ref{lemma.3.5.forgetful-detects-reflects-sliceequivs} again, we get that
	$f$ is a weak equivalence in $\zerosliceAmod$.
	Now, $M$ is cofibrant in $\qsliceAmod$ and $f:A\rightarrow A'$ may be considered
	as a map of right $A$-modules; therefore theorem 
	\ref{thm.3.5.notcommutativesmash-Quillenbifunctor-SpxSq----->Sqplusp-coefss} together with Ken Brown's lemma
	(see lemma \ref{lem1.1.KenBrown}) 
	imply that it suffices to show that $A$ and $A'$ are both cofibrant
	in $\zerosliceAmod$.
	
	We have that $\symmspherespectrum$ is cofibrant in $\zeroslicesymmTspectra$
	by lemma \ref{lem.3.5.spherespectrum-0-connected-0-slice-cofibrant}, therefore
	theorem \ref{thm.3.5.symmetrization-qslice-Quillen-equiv} 
	implies that $A$ is cofibrant in $\zerosliceAmod$.
	
	Finally, since $\zerosliceAmod$ is a right Bousfield localization with respect
	to $\weightoneAmod$, \cite[proposition 3.2.2(2)]{MR1944041} 
	implies that to show that $A'$ is cofibrant in $\zerosliceAmod$
	it suffices to check that $A'$ is cofibrant in $\weightoneAmod$ and that $f$
	is a weak equivalence in $\weightoneAmod$.  
	
	On the other hand, $\weightoneAmod$
	is a left Bousfield localization with respect to $\motivicAmod$, hence
	$A'$ is cofibrant in $\weightoneAmod$ since by hypothesis
	$A'$ is cofibrant in $\motivicAmod$. Now, we are assuming that $f$ is
	a weak equivalence in $\zeroslicesymmTspectra$ and that $A$, $A'$ are both
	cofibrant in $\zeroslicesymmTspectra$; therefore, \cite[theorem 3.2.13(2)]{MR1944041}
	implies that $f$ is also a weak equivalence in $\weightonesymmTspectra$.
	But since $A$ is cofibrant in $\zeroslicesymmTspectra$ 
	and lemma \ref{lem.3.5.spherespectrum-0-connected-0-slice-cofibrant}
	implies that $\symmspherespectrum$ is also cofibrant in $\zeroslicesymmTspectra$, we can apply 
	proposition \ref{prop.3.5.smash-induces-weightq-modelstructures} and
	theorem \ref{thm.3.5.weightq-induced==leftBousfieldlocalization}
	to conclude that $f$ is a weak equivalence in $\weightoneAmod$,
	as we wanted.	
\end{proof}

\begin{prop}
		\label{prop.3.5.sliceinvarienceofcoeffs-slicecofibrant}
	Fix $q\in \mathbb Z$.  Let $f:A\rightarrow A'$ be a map 
	between cofibrant ring spectra in $\motivicsymmTspectra$,
	which is compatible with the ring structures.
	Assume that $A$ and $A'$ are cofibrant in $\zeroslicesymmTspectra$.
	Furthermore, assume that $A'$ is also cofibrant in $\motivicAmod$.
	If $f$ is a weak equivalence in $\zeroslicesymmTspectra$, then
	it induces a Quillen equivalence between the $q$-slice motivic stable model structures of $A$ and $A'$ modules:
		$$\xymatrix{(A'\wedge _{A} -, U,\varphi ):\qsliceAmod \ar[r] & \qsliceAAmod}
		$$
\end{prop}
\begin{proof}
	We have shown in lemma 
	\ref{lemma.3.5.changecoefficients-slice-Quillenfunctor} that
		$$(A'\wedge _{A} -, U,\varphi ):\qsliceAmod \rightarrow \qsliceAAmod$$
	is a Quillen adjunction.
		
	Now 
	let $\eta$, $\epsilon$ denote the unit and counit of the
	adjunction $(A' \wedge _{A}-, U, \varphi)$.
	By corollary 1.3.16(c) in \cite{MR1650134}, it suffices to check
	that the following conditions hold:
		\begin{enumerate}
			\item \label{prop.3.5.sliceinvarienceofcoeffs-slicecofibrant.a}  For every cofibrant
						$A$-module $M$ in $\qsliceAmod$, the following composition
							$$\xymatrix{M\cong A\wedge _{A}M \ar[rr]^-{\eta _{M}=f\wedge _{A}id} && A'\wedge _{A}M
								\ar[rr]^-{W^{m',A'\wedge _{A}M}_{q+1}}&& W_{q+1}^{m'}(A'\wedge _{A}M)}
							$$
						is a weak equivalence in $\qsliceAmod$, where $W^{m'}_{q+1}$ denotes a fibrant replacement functor
						in $\qsliceAAmod$ (see proposition \ref{prop.3.5.Wsigmaqplusone-fibrant-replacement-all-Sq}).
			\item \label{prop.3.5.sliceinvarienceofcoeffs-slicecofibrant.b}  $U$ reflects
						weak equivalences between fibrant objects in $\qsliceAAmod$.
		\end{enumerate}
		
	(\ref{prop.3.5.sliceinvarienceofcoeffs-slicecofibrant.a}):  Lemma \ref{lemma.3.5.changecoeffs-preserves-cofibrantequivs} 
	implies that $f\wedge _{A}id$
	is a weak equivalence in $\qsliceAmod$, and lemma \ref{lemma.3.5.forgetful-detects-reflects-sliceequivs} 
	implies that $W^{m',A'\wedge _{A}M}_{q+1}$ is also
	a weak equivalence in $\qsliceAmod$.  Therefore, the result follows from the two out
	of three property for weak equivalences.
	
	(\ref{prop.3.5.sliceinvarienceofcoeffs-slicecofibrant.b}): This follows
	immediately from lemma \ref{lemma.3.5.forgetful-detects-reflects-sliceequivs}.
\end{proof}

\begin{thm}
		\label{thm.3.5.inheritmodelstructures-slices}
	Fix $q\in \mathbb Z$.  Let $A$ be a cofibrant ring spectrum in $\motivicsymmTspectra$,
	which is also cofibrant in $\zeroconnectedsymmTspectra$, and let $M$ be an arbitrary $A$-module.
	Then the solid arrows in the following commutative diagram:
		\begin{equation}
					\label{thm.3.5.inheritmodelstructures-slices.a}
			\begin{array}{c}
				\xymatrix@C=1.4pc{C_{q}^{\Sigma}W_{q+1}^{m}C_{q}^{m}M \ar@{-->}[dd]_-{C_{q}^{\Sigma ,W_{q+1}^{m}C_{q}^{m}M}} && 
									\ar@{-->}[ll]^-{C_{q}^{\Sigma}(W_{q+1}^{m,C_{q}^{m}M})} 
									\ar[rr]^-{C_{q}^{\Sigma}(C_{q}^{m,M})} C_{q}^{\Sigma}C_{q}^{m}M  
									\ar[dd]_-{C_{q}^{\Sigma ,C_{q}^{m}M}}&& C_{q}^{\Sigma}M 
									\ar[rr]^-{C_{q}^{\Sigma}(R_{m}^{M})} \ar@{-->}[dd]_-{C_{q}^{\Sigma ,M}} && C_{q}^{\Sigma}R_{m}M 
									\ar@{-->}[dd]_-{C_{q}^{\Sigma ,R_{m}M}}\\
									&&&&&& \\
									W_{q+1}^{m}C_{q}^{m}M && \ar[ll]^-{W_{q+1}^{m,C_{q}^{m}M}} C_{q}^{m}M \ar@{-->}[rr]_-{C_{q}^{m,M}}&&  M 
									\ar@{-->}[rr]_-{R_{m}^{M}} && R_{m}M}
			\end{array}
		\end{equation}
	induce a natural equivalence between the functors:
		\begin{equation}
					\label{thm.3.5.inheritmodelstructures-slices.b}
			\begin{array}{c}
				\xymatrix{& \qconnectedsymmstablehomotopy \ar[dr]^-{C_{q}^{\Sigma}}&\\
												 \qconnectedAmodstablehomotopy
												 \ar[ur]^-{UR_{m}} \ar[dr]_-{C^{m}_{q}}&& \qslicesymmstablehomotopy \\
												& \qsliceAmodstablehomotopy \ar[ur]_-{UW_{q+1}^{m}}&}
			\end{array}
		\end{equation}
\end{thm}
\begin{proof}
	Clearly, it is enough to prove that the maps $W_{q+1}^{m,C_{q}^{m}M}$, $C_{q}^{\Sigma ,C_{q}^{m}M}$, 
	$C_{q}^{\Sigma}(C_{q}^{m,M})$ and $C_{q}^{\Sigma}(R_{m}^{M})$ are all weak equivalences in $\qslicesymmTspectra$.
	
	Lemma \ref{lem.3.5.spherespectrum-0-connected-0-slice-cofibrant} implies that 
	$\symmspherespectrum$ is cofibrant in $\zeroslicesymmTspectra$,
	and proposition \ref{prop.3.3.symmetricq-connected--->q-slice===leftQuilllenfunctor} 
	implies that $A$ is also cofibrant in $\zeroslicesymmTspectra$.
	
	Now, proposition \ref{prop.3.5.Wsigmaqplusone-fibrant-replacement-all-Sq} implies that $W_{q+1}^{m}$ is 
	a fibrant replacement functor in $\qsliceAmod$, then 
	using lemma \ref{lemma.3.5.forgetful-detects-reflects-sliceequivs} we get that
	$W_{q+1}^{m,C_{q}^{m}M}$ is a weak equivalence in $\qslicesymmTspectra$.
	
	By construction $\qslicesymmTspectra$ is a right Bousfield localization with
	respect to $\weightqplusonesymmTspectra$, and on the other hand,
	$\weightqplusonesymmTspectra$ is a left Bousfield localization with respect
	to $\motivicsymmTspectra$.  Hence, \cite[proposition 3.1.5]{MR1944041} implies that it suffices to show
	that the remaining maps $C_{q}^{\Sigma ,C_{q}^{m}M}$, $C_{q}^{\Sigma}(C_{q}^{m,M})$ 
	and $C_{q}^{\Sigma}(R_{m}^{M})$ are weak equivalences in $\motivicsymmTspectra$.
	We will show that this is the case.
	
	Since $A$ is cofibrant in $\zeroconnectedsymmTspectra$, 
	proposition \ref{prop.3.5.Acolocal=>colocalcofibrations=>cofibrations}
	implies that $C_{q}^{m}M$ is cofibrant in $\qconnectedsymmTspectra$, and
	$C_{q}^{\Sigma ,C_{q}^{m}M}$ is by definition a weak equivalence in $\qconnectedsymmTspectra$;
	therefore \cite[theorem 3.2.13(2)]{MR1944041} implies that $C_{q}^{\Sigma ,C_{q}^{m}M}$
	is a weak equivalence in $\motivicsymmTspectra$, since $C_{q}^{\Sigma}C_{q}^{m}M$ is
	also cofibrant in $\qconnectedsymmTspectra$.
	
	Since $C_{q}^{\Sigma}C_{q}^{m}M$ and $C_{q}^{\Sigma}M$ are both cofibrant in $\qconnectedsymmTspectra$ by construction, 
	using theorem 3.2.13(2) in \cite{MR1944041} we get that if $C_{q}^{\Sigma}(C_{q}^{m,M})$
	is a weak equivalence in $\qconnectedsymmTspectra$, then it is also a weak equivalence
	in $\motivicsymmTspectra$.  But it is clear that $C_{q}^{\Sigma ,C_{q}^{m}M}$ and $C_{q}^{\Sigma ,M}$
	are both weak equivalences in $\qconnectedsymmTspectra$, then by the two out of three
	property of weak equivalences, it is enough to check that the map $C_{q}^{m,M}$
	is a weak equivalence in $\qconnectedsymmTspectra$.  Applying lemma 
	\ref{lemma.3.5.changecoefficients-Quillenfunctor} we get that $C_{q}^{m,M}$ is a weak equivalence in $\qconnectedsymmTspectra$,
	since $C_{q}^{m,M}$ is by construction a weak equivalence in $\qconnectedAmod$.
	
	Since $C_{q}^{\Sigma}M$ and $C_{q}^{\Sigma}R_{m}M$ are both cofibrant in $\qconnectedsymmTspectra$ by construction, 
	using theorem 3.2.13(2) in \cite{MR1944041} again, we get that if $C_{q}^{\Sigma}(R_{m}^{M})$
	is a weak equivalence in $\qconnectedsymmTspectra$ then it is also a weak equivalence
	in $\motivicsymmTspectra$.  But it is clear that $C_{q}^{\Sigma ,M}$ and $C_{q}^{\Sigma , R_{m}M}$
	are both weak equivalences in $\qconnectedsymmTspectra$, then by the two out of three
	property of weak equivalences, it is enough to check that the map $R_{m}^{M}$
	is a weak equivalence in $\qconnectedsymmTspectra$.  However, theorem \ref{thm.2.8.motivicsymmetricA-modules} 
	and definition \ref{def.3.5.fibrantstable-replacementfunctor} imply that 
	$R_{m}^{M}$ is a weak equivalence in $\motivicsymmTspectra$, and by
	\cite[proposition 3.1.5]{MR1944041} we have that $R_{m}^{M}$ is a weak equivalence in $\qconnectedsymmTspectra$.
	This finishes the proof.
\end{proof}

\begin{thm}
		\label{thm.3.5.sliceinvarianceofcoefficients}
	Fix $q\in \mathbb Z$.  Let $f:A\rightarrow A'$ be a map 
	between cofibrant ring spectra in $\motivicsymmTspectra$,
	which is compatible with the ring structures.
	Assume that one of the following conditions holds:
		\begin{enumerate}
			\item	\label{thm.3.5.sliceinvarianceofcoefficients.a}  $f$ is a weak equivalence in $\motivicsymmTspectra$.
			\item \label{thm.3.5.sliceinvarianceofcoefficients.b} There exists $p\in \mathbb Z$
								such that $A$, $A'$ are both $L^{m}(<p)$-local 
								in $\motivicsymmTspectra$ and $f$ is
								a weak equivalence in $\weightpsymmTspectra$.
			\item	\label{thm.3.5.sliceinvarianceofcoefficients.c} There exists $p\in \mathbb Z$
								such that $A$, $A'$ are both $C_{eff}^{p,\Sigma}$-colocal
								in $\motivicsymmTspectra$ and $f$ is a weak equivalence in
								$\pconnectedsymmTspectra$ (equivalently in $\pconnectedAmod$).
			\item	\label{thm.3.5.sliceinvarianceofcoefficients.d} $A$, $A'$ are both cofibrant
								in $\zeroconnectedsymmTspectra$, $A'$ is also cofibrant in $\motivicAmod$
								and $f$ is a weak equivalence in
								$\zeroslicesymmTspectra$.
			\item	\label{thm.3.5.sliceinvarianceofcoefficients.e} $A$, $A'$ are both cofibrant
								in $\zeroslicesymmTspectra$, $A'$ is also cofibrant in $\motivicAmod$
								and $f$ is a weak equivalence in
								$\zeroslicesymmTspectra$.
		\end{enumerate}
	Then $f$ induces a Quillen equivalence between the $q$-slice motivic stable model structures of $A$ and $A'$ modules:
		$$\xymatrix{(A'\wedge _{A} -, U,\varphi ):\qsliceAmod \ar[r] & \qsliceAAmod}
		$$
\end{thm}
\begin{proof}
	(\ref{thm.3.5.sliceinvarianceofcoefficients.a}):  This is just 
	proposition \ref{prop.3.5.changecoefficients-slice-Quillenequiv}.
	
	(\ref{thm.3.5.sliceinvarianceofcoefficients.b}):  Since $A$ and $A'$ are $L^{m}(<p)$-local in $\motivicsymmTspectra$,
	\cite[theorem 3.2.13(1)]{MR1944041} implies that $f$ is a weak equivalence in $\motivicsymmTspectra$.
	Therefore the result follows from proposition \ref{prop.3.5.changecoefficients-slice-Quillenequiv}.
	
	(\ref{thm.3.5.sliceinvarianceofcoefficients.c}):  Since $A$ and $A'$ are  $C_{eff}^{p,\Sigma}$-colocal
	in $\motivicsymmTspectra$, using \cite[theorem 3.2.13(2)]{MR1944041} we have that $f$ is a weak equivalence in $\motivicsymmTspectra$.
	Thus, the result follows from proposition \ref{prop.3.5.changecoefficients-slice-Quillenequiv}.
	
	(\ref{thm.3.5.sliceinvarianceofcoefficients.d}):  Proposition
	\ref{prop.3.3.symmetricq-connected--->q-slice===leftQuilllenfunctor} implies that $A$ and $A'$ are both
	cofibrant in $\zeroslicesymmTspectra$, therefore the result follows from 
	proposition \ref{prop.3.5.sliceinvarienceofcoeffs-slicecofibrant}.
	
	(\ref{thm.3.5.sliceinvarianceofcoefficients.e}):  This is just proposition 
	\ref{prop.3.5.sliceinvarienceofcoeffs-slicecofibrant}.
\end{proof}

\end{section}
\begin{section}{Applications}
		\label{section.3.5.applications}
		
	In this section we will describe some of the consequences that follow from
	the compatibility of the slice filtration with the smash product of symmetric
	$T$-spectra in the sense of theorems \ref{thm.3.4.smash-Quillenbifunctor-RpxRq----->Rqplusp} and 
	\ref{thm.3.4.smash-Quillenbifunctor-SpxSq----->Sqplusp}, as well as those
	that follow from the compatibility between the slice filtration
	on the categories of symmetric $T$-spectra and $A$-modules in the
	sense of propositions \ref{prop.3.5.functors-between-Rqsymm}, \ref{prop.3.5.Lq+1-->Lqsymm}, 
	\ref{prop.3.5.homotopycoherence==>liftings-qslice} and
	theorems \ref{thm.3.5.inheritingmodulestructures1}, \ref{thm.3.5.inheritmodelstructures2},
	\ref{thm.3.5.inheritmodelstructures-slices}, \ref{thm.3.5.sliceinvarianceofcoefficients}.

\begin{prop}
		\label{prop.3.5.R0--S0---symmetricmonoidalmodelcategories}
	The model categories $\zeroconnectedsymmTspectra$ and $\zeroslicesymmTspectra$
	are both symmetric monoidal (with respect to the smash product of symmetric $T$-spectra)
	model categories in the sense of Hovey (see definition \ref{def.mon-mod-cats}).
\end{prop}
\begin{proof}
	Follows directly from lemma \ref{lem.3.5.spherespectrum-0-connected-0-slice-cofibrant},
	together with theorems
	\ref{thm.3.4.smash-Quillenbifunctor-RpxRq----->Rqplusp} and
	\ref{thm.3.4.smash-Quillenbifunctor-SpxSq----->Sqplusp} 
\end{proof}

\begin{thm}
		\label{thm.3.5.homotopycategories-zero-connected-slice-symmetric-monoidal}
	The triangulated categories $\symmstablehomotopy$,
	$\zeroconnectedsymmstablehomotopy$ and
	$\zeroslicesymmstablehomotopy$ inherit a natural symmetric monoidal structure
	from the smash product of symmetric $T$-spectra.  The symmetric monoidal structure is
	defined as follows:
		\begin{enumerate}
			\item $$\xymatrix@R=.5pt{-\wedge ^{\mathbf L}-:\symmstablehomotopy \times \symmstablehomotopy
												\ar[r]& \symmstablehomotopy \\
												(X,Y) \ar@{|->}[r]& Q_{\Sigma}X\wedge Q_{\Sigma}Y}
						$$
			\item	$$\xymatrix@R=.5pt{-\wedge ^{\mathbf L}-:\zeroconnectedsymmstablehomotopy \times \zeroconnectedsymmstablehomotopy
												\ar[r]& \zeroconnectedsymmstablehomotopy \\
												(X,Y) \ar@{|->}[r]& C_{0}^{\Sigma}X\wedge C_{0}^{\Sigma}Y}
						$$
			\item	$$\xymatrix@R=.5pt{-\wedge ^{\mathbf L}-:\zeroslicesymmstablehomotopy \times \zeroslicesymmstablehomotopy
												\ar[r]& \zeroslicesymmstablehomotopy \\
												(X,Y) \ar@{|->}[r]& P_{0}^{\Sigma}X\wedge P_{0}^{\Sigma}Y}
						$$
		\end{enumerate}
\end{thm}
\begin{proof}
	Follows directly from propositions \ref{prop.2.6.symmTspectra-monoidalmodelcategory}
	and \ref{prop.3.5.R0--S0---symmetricmonoidalmodelcategories},
	together with theorem \ref{thm.symmmon-descends-homot}.
\end{proof}

\begin{prop}
		\label{prop.3.5.cofibrantreplacementfunctors--strongsymmetricmonoidal}
	The following exact functors between triangulated categories are both
	strong symmetric monoidal:
		$$\xymatrix@R=.5pt{C_{0}^{\Sigma}:\zeroconnectedsymmstablehomotopy \ar[r]& \zeroslicesymmstablehomotopy \\
											 C_{0}^{\Sigma}:\zeroconnectedsymmstablehomotopy \ar[r]& \symmstablehomotopy}
		$$
\end{prop}
\begin{proof}
	Propositions \ref{prop.2.6.symmTspectra-monoidalmodelcategory} and 
	\ref{prop.3.5.R0--S0---symmetricmonoidalmodelcategories} imply that
	$\symmTspectra$, $\zeroconnectedsymmTspectra$ and $\zeroslicesymmTspectra$ are
	all symmetric monoidal model categories in the sense of Hovey.  Now,
	using proposition \ref{prop.3.3.symmetricq-connected--->q-slice===leftQuilllenfunctor} 
	and theorem \ref{thm.3.3.symmetricconnective-model-structure} we have that
	the following adjunctions
		$$\xymatrix@R=.5pt{(id, id,\varphi):\zeroconnectedsymmTspectra \ar[r]& \zeroslicesymmTspectra \\
											 (id, id,\varphi):\zeroconnectedsymmTspectra \ar[r]& \symmTspectra}
		$$
	are both symmetric monoidal Quillen adunctions
	(see definition \ref{def.mon-Quillen-func}).  The result then
	follows immediately from theorem 4.3.3 in \cite{MR1650134}.
\end{proof}

\begin{cor}
		\label{cor.3.5.fibrantreplacementfunctors--laxsymmetricmonoidal}
	The following exact functors between triangulated categories are both
	lax symmetric monoidal:
		$$\xymatrix@R=.5pt{R_{\Sigma}:\symmstablehomotopy \ar[r]& \zeroconnectedsymmstablehomotopy \\
											 W_{1}^{\Sigma}:\zeroslicesymmstablehomotopy \ar[r]& \zeroconnectedsymmstablehomotopy}
		$$
\end{cor}
\begin{proof}
	By proposition \ref{prop.3.3.symmcofibrant-replacement=>triangulatedfunctor} and 
	corollary \ref{cor.3.3.adjunctions--symmetricRq==>Sq} we have the following adjunctions
		$$\xymatrix@R=.5pt{(C_{0}^{\Sigma}, R_{\Sigma}, \varphi):\zeroconnectedsymmstablehomotopy \ar[r]& \symmstablehomotopy \\
												(C_{0}^{\Sigma}, W_{1}^{\Sigma}, \varphi):\zeroconnectedsymmstablehomotopy \ar[r]& \zeroslicesymmstablehomotopy}
		$$
	Using proposition \ref{prop.3.5.cofibrantreplacementfunctors--strongsymmetricmonoidal}
	we have that the left adjoints for $R_{\Sigma}$ and $W_{1}^{\Sigma}$ are both strong symmetric
	monoidal.  Finally by standard results in category theory we get that the right adjoints
	$R_{\Sigma}$ and $W_{1}^{\Sigma}$ are both lax symmetric monoidal
	(see \cite[theorem 1.5]{MR0360749}).
\end{proof}

\begin{prop}
		\label{prop.3.5.modules---over--zero-slice-connected}
	Fix $q\in \mathbb Z$.  Then the smash product of symmetric $T$-spectra induces the
	following Quillen adjunctions of two variables:
		\begin{enumerate}
			\item	\label{prop.3.5.modules---over--zero-slice-connected.a}  $\qconnectedsymmTspectra$
						is a $\zeroconnectedsymmTspectra$-model category in the sense of Hovey
						(see definition \ref{def.module-modcats}).
			\item	\label{prop.3.5.modules---over--zero-slice-connected.b} $\qslicesymmTspectra$ is
						a $\zeroslicesymmTspectra$-model category in the sense of Hovey.
			\item	\label{prop.3.5.modules---over--zero-slice-connected.c} $\motivicsymmTspectra$ is
						a $\zeroconnectedsymmTspectra$-model category in the sense of Hovey.
			\item	\label{prop.3.5.modules---over--zero-slice-connected.d} $\qslicesymmTspectra$ is
						a $\zeroconnectedsymmTspectra$-model category in the sense of Hovey.
		\end{enumerate}
\end{prop}
\begin{proof}
	(\ref{prop.3.5.modules---over--zero-slice-connected.a}):  Follows immediately from
	lemma \ref{lem.3.5.spherespectrum-0-connected-0-slice-cofibrant} and 
	theorem \ref{thm.3.4.smash-Quillenbifunctor-RpxRq----->Rqplusp}.
	
	(\ref{prop.3.5.modules---over--zero-slice-connected.b}):  Follows immediately from
	lemma \ref{lem.3.5.spherespectrum-0-connected-0-slice-cofibrant} and 
	theorem \ref{thm.3.4.smash-Quillenbifunctor-SpxSq----->Sqplusp}.
	
	(\ref{prop.3.5.modules---over--zero-slice-connected.c}): Follows from
	proposition
	\ref{prop.2.6.symmTspectra-monoidalmodelcategory} and
	theorem \ref{thm.3.3.symmetricconnective-model-structure} which imply that the following composition is a
	Quillen adjunction of two variables:
		$$\xymatrix{\zeroconnectedsymmTspectra \times \motivicsymmTspectra \ar[rr]^-{(id,id)}&&
								\motivicsymmTspectra \times \motivicsymmTspectra \ar[d]^-{-\wedge -}\\ 
								&& \motivicsymmTspectra}
		$$
	
	(\ref{prop.3.5.modules---over--zero-slice-connected.d}):  Follows from
	proposition
	\ref{prop.3.3.symmetricq-connected--->q-slice===leftQuilllenfunctor} and
	theorem \ref{thm.3.4.smash-Quillenbifunctor-SpxSq----->Sqplusp} which imply that the following composition is a
	Quillen adjunction of two variables:
		$$\xymatrix{\zeroconnectedsymmTspectra \times \qslicesymmTspectra \ar[rr]^-{(id,id)}&&
								\zeroslicesymmTspectra \times \qslicesymmTspectra \ar[d]^-{-\wedge -}\\ 
								&& \qslicesymmTspectra}
		$$
\end{proof}

\begin{thm}
		\label{thm.3.5.homotopymodules-over-zero-slice-connected}
	Fix $q\in \mathbb Z$.  Then the smash product of symmetric $T$-spectra
	induces the following natural module structures
	(see definition \ref{def.leftmods/moncats}):
		\begin{enumerate}
			\item \label{thm.3.5.homotopymodules-over-zero-slice-connected.a} The triangulated category
						$\qconnectedsymmstablehomotopy$
						has a natural structure of $\zeroconnectedsymmstablehomotopy$-module,
						defined as follows:
							$$\xymatrix@R=.5pt{-\wedge ^{\mathbf L}-:\zeroconnectedsymmstablehomotopy \times \qconnectedsymmstablehomotopy
															\ar[r]& \qconnectedsymmstablehomotopy \\
															(X,Y) \ar@{|->}[r]& C_{0}^{\Sigma}X\wedge C_{q}^{\Sigma}Y}
							$$
			\item \label{thm.3.5.homotopymodules-over-zero-slice-connected.b} The triangulated category
						$\qslicesymmstablehomotopy$ has a natural structure of
						$\zeroslicesymmstablehomotopy$-module, defined as follows:
							$$\xymatrix@R=.5pt{-\wedge ^{\mathbf L}-:\zeroslicesymmstablehomotopy \times \qslicesymmstablehomotopy
															\ar[r]& \qslicesymmstablehomotopy \\
															(X,Y) \ar@{|->}[r]& P_{0}^{\Sigma}X\wedge P_{q}^{\Sigma}Y}
							$$
			\item \label{thm.3.5.homotopymodules-over-zero-slice-connected.c} The triangulated category
						$\symmstablehomotopy$ has a natural structure of
						$\zeroconnectedsymmstablehomotopy$-module, defined as follows:
							$$\xymatrix@R=.5pt{-\wedge ^{\mathbf L}-:\zeroconnectedsymmstablehomotopy \times \symmstablehomotopy
															\ar[r]& \symmstablehomotopy \\
															(X,Y) \ar@{|->}[r]& C_{0}^{\Sigma}X\wedge Q_{\Sigma}Y}
							$$
			\item \label{thm.3.5.homotopymodules-over-zero-slice-connected.d} The triangulated category
						$\qslicesymmstablehomotopy$ has a natural structure of
						$\zeroconnectedsymmstablehomotopy$-module, defined as follows:
							$$\xymatrix@R=.5pt{-\wedge ^{\mathbf L}-:\zeroconnectedsymmstablehomotopy \times \qslicesymmstablehomotopy
															\ar[r]& \qslicesymmstablehomotopy \\
															(X,Y) \ar@{|->}[r]& C_{0}^{\Sigma}X\wedge P_{q}^{\Sigma}Y}
							$$
		\end{enumerate}
\end{thm}
\begin{proof}
	Follows directly from lemma \ref{lem.3.5.spherespectrum-0-connected-0-slice-cofibrant}, 
	proposition \ref{prop.3.5.modules---over--zero-slice-connected} and
	\cite[theorem 4.3.4]{MR1650134}.
\end{proof}

\begin{thm}
		\label{thm.3.5.bilinear-pairings}
	Fix $p, q\in \mathbb Z$.  Then the smash product of symmetric $T$-spectra
	induces the following adjunctions of two variables
	(see definition \ref{def.adj-two-vars}):
		\begin{enumerate}
			\item \label{thm.3.5.bilinear-pairings.a} We have the following
							adjunction of two variables, which is also a bilinear pairing:
							$$\xymatrix@R=.5pt{-\wedge ^{\mathbf L}-:\pconnectedsymmstablehomotopy \times \qconnectedsymmstablehomotopy
															\ar[r]& \pplusqconnectedsymmstablehomotopy \\
															(X,Y) \ar@{|->}[r]& C_{p}^{\Sigma}X\wedge C_{q}^{\Sigma}Y}
							$$
			\item \label{thm.3.5.bilinear-pairings.b} We have the following
							adjunction of two variables, which is also a bilinear pairing:
							$$\xymatrix@R=.5pt{-\wedge ^{\mathbf L}-:\pslicesymmstablehomotopy \times \qslicesymmstablehomotopy
															\ar[r]& \pplusqslicesymmstablehomotopy \\
															(X,Y) \ar@{|->}[r]& P_{p}^{\Sigma}X\wedge P_{q}^{\Sigma}Y}
							$$
		\end{enumerate}
\end{thm}
\begin{proof}
	(\ref{thm.3.5.bilinear-pairings.a}):  By theorem \ref{thm.3.4.smash-Quillenbifunctor-RpxRq----->Rqplusp}
	we have that
		$$\xymatrix{-\wedge -:
									\pconnectedsymmTspectra \times \qconnectedsymmTspectra \ar[r]& \pplusqconnectedsymmTspectra}
		$$
	is a Quillen bifunctor.  Then proposition \ref{prop.Quillenbifunc-descends-homot} implies that
		$$\xymatrix@R=.5pt{-\wedge ^{\mathbf L}-:\pconnectedsymmstablehomotopy \times \qconnectedsymmstablehomotopy
											\ar[r]& \pplusqconnectedsymmstablehomotopy \\
											(X,Y) \ar@{|->}[r]& C_{p}^{\Sigma}X\wedge C_{q}^{\Sigma}Y}
		$$
	is an adjunction of two variables.  Finally, since the coproduct of two cofibrant objects
	is always cofibrant, and $X\wedge (Y\coprod Z)$ is canonically isomorphic in $\symmTspectra$ to
	$(X\wedge Y)\coprod (X\wedge Z)$, we get that the pairing $-\wedge ^{\mathbf L}-$ is bilinear.
	
	(\ref{thm.3.5.bilinear-pairings.b}):  By theorem \ref{thm.3.4.smash-Quillenbifunctor-SpxSq----->Sqplusp}
	we have that
		$$\xymatrix{-\wedge -:
									\pslicesymmTspectra \times \qslicesymmTspectra \ar[r]& \pplusqslicesymmTspectra}
		$$
	is a Quillen bifunctor.  Then proposition \ref{prop.Quillenbifunc-descends-homot} implies that
		$$\xymatrix@R=.5pt{-\wedge ^{\mathbf L}-:\pslicesymmstablehomotopy \times \qslicesymmstablehomotopy
											\ar[r]& \pplusqslicesymmstablehomotopy \\
											(X,Y) \ar@{|->}[r]& P_{p}^{\Sigma}X\wedge P_{q}^{\Sigma}Y}
		$$
	is an adjunction of two variables.  Finally, since the coproduct of two cofibrant objects
	is always cofibrant and $X\wedge (Y\coprod Z)$ is canonically isomorphic in $\symmTspectra$ to
	$(X\wedge Y)\coprod (X\wedge Z)$, we get that the pairing $-\wedge ^{\mathbf L}-$ is bilinear.				
\end{proof}

\begin{prop}
		\label{prop.3.5.change-of-base-induced-pairings}
	Fix $p, q\in \mathbb Z$, and let $X, Y$ be two arbitrary symmetric
	$T$-spectra.
		\begin{enumerate}
			\item	\label{prop.3.5.change-of-base-induced-pairings.a}There exists a natural
						bilinear isomorphism in $\symmstablehomotopy$:
							$$\xymatrix{C_{p}^{\Sigma}X\wedge C_{q}^{\Sigma}Y \ar[r]^-{m_{1}^{X,Y}}&
													C_{p+q}^{\Sigma}(C_{p}^{\Sigma}X\wedge C_{q}^{\Sigma}Y)}
							$$
			
			\item \label{prop.3.5.change-of-base-induced-pairings.b}There exists a natural bilinear map
						in $\pplusqconnectedsymmstablehomotopy$:
							$$\xymatrix{C_{p}^{\Sigma}R_{\Sigma}X\wedge C_{q}^{\Sigma}R_{\Sigma}Y \ar[r]^-{m_{2}^{X,Y}}&
													R_{\Sigma}(X\wedge Y)}
							$$
			
			\item	\label{prop.3.5.change-of-base-induced-pairings.c}There exists a natural bilinear
						isomorphism in $\pplusqslicesymmstablehomotopy$:
							$$\xymatrix{C_{p}^{\Sigma}X\wedge C_{q}^{\Sigma}Y \ar[r]^-{m_{3}^{X,Y}}&
													C_{p+q}^{\Sigma}(C_{p}^{\Sigma}X\wedge C_{q}^{\Sigma} Y)}
							$$
			
			\item	\label{prop.3.5.change-of-base-induced-pairings.d}There exists a natural bilinear map
						in $\pplusqconnectedsymmstablehomotopy$:
							$$\xymatrix{C_{p}^{\Sigma}W_{p+1}^{\Sigma}X\wedge 
													C_{q}^{\Sigma}W_{q+1}^{\Sigma}Y \ar[r]^-{m_{4}^{X,Y}}&
													W_{p+q+1}^{\Sigma}(X\wedge Y)}
							$$
		\end{enumerate}
\end{prop}
\begin{proof}
	(\ref{prop.3.5.change-of-base-induced-pairings.a}): Theorems \ref{thm.3.4.smash-Quillenbifunctor-RpxRq----->Rqplusp} and
	\ref{thm.3.3.symmetricconnective-model-structure}
	imply that we have the following commutative diagram
	of Quillen bifunctors:
		$$\xymatrix{\pconnectedsymmTspectra \times \qconnectedsymmTspectra  \ar[d]_-{-\wedge -} \ar[dr]^-{-\wedge -}&\\
								\pplusqconnectedsymmTspectra \ar[r]_-{id}& \motivicsymmTspectra}
		$$
	Using \cite[theorem 1.3.7]{MR1650134} we get the natural isomorphism $m_{1}$, which is bilinear since
	the functors $C_{p}^{\Sigma}$, $C_{q}^{\Sigma}$, $C_{p+q}^{\Sigma}$ are all exact and the smash product
	is bilinear.
	
	(\ref{prop.3.5.change-of-base-induced-pairings.b}):  By proposition
	\ref{prop.3.3.symmcofibrant-replacement=>triangulatedfunctor} we have the following adjunctions:
		$$\xymatrix@R=.5pt{(C_{p}^{\Sigma},R_{\Sigma},\varphi):\pconnectedsymmstablehomotopy \ar[r]& \symmstablehomotopy \\
											 (C_{q}^{\Sigma},R_{\Sigma},\varphi):\qconnectedsymmstablehomotopy \ar[r]& \symmstablehomotopy}
		$$
	Let $\epsilon _{p}$, $\epsilon _{q}$ denote the respective counits, and let $\check{m}_{2}^{X,Y}$ be the
	following composition in $\symmstablehomotopy$:
		$$\xymatrix{C_{p+q}^{\Sigma}(C_{p}^{\Sigma}R_{\Sigma}X \wedge C_{q}^{\Sigma}R_{\Sigma}Y)
								\ar[rrr]^-{(m_{1}^{R_{\Sigma}X,R_{\Sigma}Y})^{-1}}
								\ar@{-->}[d]_-{\check{m}_{2}^{X,Y}}&&&
								C_{p}^{\Sigma}R_{\Sigma}X\wedge C_{q}^{\Sigma}R_{\Sigma}Y \ar[dlll]^-{\epsilon _{p}^{X}\wedge \epsilon _{q}^{Y}}\\ 
								X\wedge Y &&& }
		$$
	Then using the adjunction between $C_{p+q}^{\Sigma}$ and $R_{\Sigma}$ considered above,
	we define $m_{2}^{X,Y}$ as the adjoint of $\check{m}_{2}^{X,Y}$.  The naturality of $m_{2}$ follows
	from:
		\begin{enumerate} 
			\item the naturality of $m_{1}$
			\item the naturality of the fibrant replacement functor, and 
			\item the naturality of the counits $\epsilon _{p}$ and $\epsilon _{q}$.
		\end{enumerate}
	Finally we have that $m_{2}$ is bilinear since:
		\begin{enumerate}
			\item	$m_{1}$ is bilinear
			\item the functors $C_{p}^{\Sigma}$, $C_{q}^{\Sigma}$, $C_{p+q}^{\Sigma}$ and $R_{\Sigma}$
						are all exact, and
			\item	the smash product is bilinear.
		\end{enumerate}
		
	(\ref{prop.3.5.change-of-base-induced-pairings.c}):Theorem \ref{thm.3.4.smash-Quillenbifunctor-RpxRq----->Rqplusp} and
	proposition
	\ref{prop.3.3.symmetricq-connected--->q-slice===leftQuilllenfunctor} 
	imply that we have the following commutative diagram
	of Quillen bifunctors:
		$$\xymatrix{\pconnectedsymmTspectra \times \qconnectedsymmTspectra  \ar[d]_-{-\wedge -} \ar[rd]^-{-\wedge -}&\\
								\pplusqconnectedsymmTspectra \ar[r]_-{id}& \pplusqslicesymmTspectra}
		$$
	Using \cite[theorem 1.3.7]{MR1650134} we get the natural isomorphism $m_{3}$, which is bilinear since
	the functors $C_{p}^{\Sigma}$, $C_{q}^{\Sigma}$, $C_{p+q}^{\Sigma}$ are all exact and the smash product
	is bilinear.
	
	(\ref{prop.3.5.change-of-base-induced-pairings.d}):By corollary
	\ref{cor.3.3.adjunctions--symmetricRq==>Sq} we have the following adjunctions:
		$$\xymatrix@R=.5pt{(C_{p}^{\Sigma},W_{p+1}^{\Sigma},\varphi):\pconnectedsymmstablehomotopy \ar[r]& \qslicesymmstablehomotopy \\
											 (C_{q}^{\Sigma},W_{q+1}^{\Sigma},\varphi):\qconnectedsymmstablehomotopy \ar[r]& \qslicesymmstablehomotopy}
		$$
	Let $\epsilon _{p}$, $\epsilon _{q}$ denote the respective counits, and let $\check{m}_{4}^{X,Y}$ be the
	following composition in $\qslicesymmstablehomotopy$:
		$$\xymatrix{C_{p+q}^{\Sigma}(C_{p}^{\Sigma}W_{p+1}^{\Sigma}X \wedge C_{q}^{\Sigma}W_{q+1}^{\Sigma}Y)
								\ar[rrr]^-{(m_{3}^{W_{p+1}^{\Sigma}X,W_{q+1}^{\Sigma}Y})^{-1}}
								\ar@{-->}[d]_-{\check{m}_{4}^{X,Y}}&&&
								C_{p}^{\Sigma}W_{p+1}^{\Sigma}X\wedge C_{q}^{\Sigma}W_{q+1}^{\Sigma}Y 
								\ar[dlll]^-{\epsilon _{p}^{X}\wedge \epsilon _{q}^{Y}}\\
								X\wedge Y &&&}
		$$
	Then using the adjunction between $C_{p+q}^{\Sigma}$ and $W_{p+q+1}^{\Sigma}$ considered above,
	we define $m_{4}^{X,Y}$ as the adjoint of $\check{m}_{4}^{X,Y}$.  The naturality of $m_{4}$ follows
	from:
		\begin{enumerate} 
			\item the naturality of $m_{3}$
			\item the naturality of the fibrant replacement functors, and 
			\item the naturality of the counits $\epsilon _{p}$ and $\epsilon _{q}$.
		\end{enumerate}
	Finally we have that $m_{4}$ is bilinear since:
		\begin{enumerate}
			\item	$m_{3}$ is bilinear
			\item the functors $C_{p}^{\Sigma}$, $C_{q}^{\Sigma}$, $C_{p+q}^{\Sigma}$, $W_{p+1}^{\Sigma}$,
						$W_{q+1}^{\Sigma}$ and $W_{p+q+1}^{\Sigma}$
						are all exact, and
			\item	the smash product is bilinear.
		\end{enumerate}
\end{proof}

\begin{thm}
		\label{thm.3.5.external-pairings}
	Fix $p,q\in \mathbb Z$.  Then the smash product of symmetric $T$-spectra induces the following natural pairings
	(external products):
		\begin{enumerate}
			\item	\label{thm.3.5.external-pairings.a}For every couple of symmetric
						$T$-spectra $X$, $Y$ we have the following natural map in $\symmstablehomotopy$:
							$$\xymatrix{f_{p}^{\Sigma}X\wedge f_{q}^{\Sigma}Y \ar[rr]^-{\cup ^{c}_{p,q}} \ar@{=}[d]&& 
													f_{p+q}^{\Sigma}(X\wedge Y) \ar@{=}[d]\\
													C_{p}^{\Sigma}R_{\Sigma}X\wedge C_{q}^{\Sigma}R_{\Sigma}Y 
													\ar[dr]^-{m_{1}^{R_{\Sigma}X,R_{\Sigma}Y}}_-{\cong}&& 
													C_{p+q}^{\Sigma}R_{\Sigma}(X\wedge Y) \\
													&C_{p+q}^{\Sigma}(C_{p}^{\Sigma}R_{\Sigma}X\wedge C_{q}^{\Sigma}R_{\Sigma}Y) 
													\ar[ur]^-{C_{p+q}^{\Sigma}(m_{2}^{X,Y})}&}
							$$
						(see proposition \ref{prop.3.5.change-of-base-induced-pairings})
						which induces a bilinear natural transformation between the functors:
							$$\xymatrix@R=.5pt{\symmstablehomotopy \times \symmstablehomotopy \ar[r]& \symmstablehomotopy \\
																	(X,Y)\ar@{|->}[r]& f_{p}^{\Sigma}X\wedge f_{q}^{\Sigma}Y }
							$$
							
							$$\xymatrix@R=.5pt{\symmstablehomotopy \times \symmstablehomotopy \ar[r]& \symmstablehomotopy \\
																	(X,Y)\ar@{|->}[r]& f_{p+q}^{\Sigma}(X\wedge Y)}
							$$
			\item	\label{thm.3.5.external-pairings.b}For every couple of symmetric
						$T$-spectra $X$, $Y$ we have the following natural map in $\symmstablehomotopy$:
							$$\xymatrix@C=-3pc{s_{p}^{\Sigma}X\wedge s_{q}^{\Sigma}Y\ar[dr]^-{\cup ^{s}_{p,q}} \ar@{=}[d]& \\
																C_{p}^{\Sigma}W_{p+1}^{\Sigma}C_{p}^{\Sigma}R_{\Sigma}X\wedge 
																C_{q}^{\Sigma}W_{q+1}^{\Sigma}C_{q}^{\Sigma}R_{\Sigma}Y 
																\ar[dd]_-{m_{1}^{W_{p+1}^{\Sigma}C_{p}^{\Sigma}R_{\Sigma}X,W_{q+1}^{\Sigma}C_{q}^{\Sigma}R_{\Sigma}Y}}^-{\cong}& 
																s_{p+q}^{\Sigma}(X\wedge Y) \ar@{=}[d]\\
																& C_{p+q}^{\Sigma}W_{p+q+1}^{\Sigma}C_{p+q}^{\Sigma}R_{\Sigma}(X\wedge Y)\\
																C_{p+q}^{\Sigma}(C_{p}^{\Sigma}W_{p+1}^{\Sigma}C_{p}^{\Sigma}R_{\Sigma}X\wedge
																C_{q}^{\Sigma}W_{q+1}^{\Sigma}C_{q}^{\Sigma}R_{\Sigma}Y) 
																\ar[dd]_-{C_{p+q}^{\Sigma}(m_{4}^{C_{p}^{\Sigma}R_{\Sigma}X,C_{q}^{\Sigma}R_{\Sigma}Y})}&& \\
																& \\
																C_{p+q}^{\Sigma}W_{p+q+1}^{\Sigma}(C_{p}^{\Sigma}R_{\Sigma}X\wedge C_{q}^{\Sigma}R_{\Sigma}Y)
																\ar[dr]_<(0.15){C_{p+q}^{\Sigma}W_{p+q+1}^{\Sigma}(m_{3}^{R_{\Sigma}X,R_{\Sigma}Y})
																\ \ \ \ \ \ \ \ \ \ }^-{\cong} & \\
																& C_{p+q}^{\Sigma}W_{p+q+1}^{\Sigma}C_{p+q}^{\Sigma}(C_{p}^{\Sigma}R_{\Sigma}X\wedge C_{q}^{\Sigma}R_{\Sigma}Y)
																\ar[uuuu]_-{C_{p+q}^{\Sigma}W_{p+q+1}^{\Sigma}C_{p+q}^{\Sigma}(m_{2}^{X,Y})}}
							$$
						(see proposition \ref{prop.3.5.change-of-base-induced-pairings})
						which induces a bilinear natural transformation between the following functors:
							$$\xymatrix@R=.5pt{\symmstablehomotopy \times \symmstablehomotopy \ar[r]& \symmstablehomotopy \\
																	(X,Y)\ar@{|->}[r]& s_{p}^{\Sigma}X\wedge s_{q}^{\Sigma}Y }
							$$
							
							$$\xymatrix@R=.5pt{\symmstablehomotopy \times \symmstablehomotopy \ar[r]& \symmstablehomotopy \\
																	(X,Y)\ar@{|->}[r]& s_{p+q}^{\Sigma}(X\wedge Y)}
							$$
		\end{enumerate}  
\end{thm}
\begin{proof}
	(\ref{thm.3.5.external-pairings.a}):  Follows immediately from (\ref{prop.3.5.change-of-base-induced-pairings.a}) 
	and (\ref{prop.3.5.change-of-base-induced-pairings.b})
	in proposition \ref{prop.3.5.change-of-base-induced-pairings}.
	
	(\ref{thm.3.5.external-pairings.b}):  Follows immediately from (\ref{prop.3.5.change-of-base-induced-pairings.a}),
	(\ref{prop.3.5.change-of-base-induced-pairings.d}), (\ref{prop.3.5.change-of-base-induced-pairings.c}) and
	(\ref{prop.3.5.change-of-base-induced-pairings.b}) in proposition \ref{prop.3.5.change-of-base-induced-pairings}.
\end{proof}

\begin{thm}
		\label{thm.3.5.external-pairings-compatible}
	Fix $p,q\in \mathbb Z$.  Then the pairings $\cup ^{c}_{p,q}$ and $\cup ^{s}_{p,q}$ constructed in theorem
	\ref{thm.3.5.external-pairings} are compatible with the natural transformations $\rho $ and
	$\pi ^{\Sigma}$ (see propositions 
	\ref{prop.3.2.functors-between-Rqsymm}(\ref{prop.3.2.functors-between-Rqsymm.c}) and
	\ref{prop.3.3.lifting-map-fq-->sqsymm}) 
	in the following sense:
		\begin{enumerate}
			\item	\label{thm.3.5.external-pairings-compatible.a}For every couple of symmetric $T$-spectra 
						$X$, $Y$; the following diagram is commutative
						in $\symmstablehomotopy$:
							$$\xymatrix{f_{p+1}^{\Sigma}X\wedge f_{q}^{\Sigma}Y \ar[rr]^-{\rho _{p}^{X}\wedge id} 
													\ar[d]_-{\cup ^{c}_{p+1,q}}&& f_{p}^{\Sigma}X\wedge f_{q}^{\Sigma}Y \ar[d]^{\cup ^{c}_{p,q}}\\
													f_{p+q+1}^{\Sigma}(X\wedge Y) \ar[rr]_-{\rho _{p+q+1}^{X\wedge Y}}&& f_{p+q}^{\Sigma}(X\wedge Y)}
							$$
			\item	\label{thm.3.5.external-pairings-compatible.b}For every couple of symmetric $T$-spectra 
						$X$, $Y$; the following diagram is commutative
						in $\symmstablehomotopy$:
							$$\xymatrix{f_{p}^{\Sigma}X\wedge f_{q+1}^{\Sigma}Y \ar[rr]^-{id \wedge \rho _{q}^{Y}} 
													\ar[d]_-{\cup ^{c}_{p,q+1}}&& f_{p}^{\Sigma}X\wedge f_{q}^{\Sigma}Y \ar[d]^{\cup ^{c}_{p,q}}\\
													f_{p+q+1}^{\Sigma}(X\wedge Y) \ar[rr]_-{\rho _{p+q+1}^{X\wedge Y}}&& f_{p+q}^{\Sigma}(X\wedge Y)}
							$$
			\item	\label{thm.3.5.external-pairings-compatible.c}For every couple of symmetric $T$-spectra 
						$X$, $Y$; the following diagram is commutative
						in $\symmstablehomotopy$:
							$$\xymatrix{f_{p}^{\Sigma}X\wedge f_{q}^{\Sigma}Y \ar[rr]^-{\pi _{p}^{\Sigma ,X}\wedge \pi _{q}^{\Sigma ,Y}} 
													\ar[d]_-{\cup ^{c}_{p,q}}&& s_{p}^{\Sigma}X\wedge s_{q}^{\Sigma}Y \ar[d]^-{\cup ^{s}_{p,q}}\\
													f_{p+q}^{\Sigma}(X\wedge Y) \ar[rr]_-{\pi _{p+q}^{\Sigma ,X\wedge Y}}&& s_{p+q}^{\Sigma}(X\wedge Y)}
							$$
		\end{enumerate}
\end{thm}
\begin{proof}
	(\ref{thm.3.5.external-pairings-compatible.a}):  This follows from the following commutative diagram of left
	Quillen (bi)functors, together with the construction of the external pairing $\cup ^{c}$ given in theorem
	\ref{thm.3.5.external-pairings}(\ref{thm.3.5.external-pairings.a}) and the construction of the natural
	transformation $\rho$ given in proposition \ref{prop.3.2.functors-between-Rqsymm}(\ref{prop.3.2.functors-between-Rqsymm.c}):
		$$\xymatrix{\pplusoneconnectedsymmTspectra \times \qconnectedsymmTspectra \ar[rr]^-{id\times id} \ar[dr]^-{-\wedge -}
								\ar[dd]_-{-\wedge -}&& \pconnectedsymmTspectra \times
								\qconnectedsymmTspectra \ar[dl]_-{-\wedge -} \ar[dd]^-{-\wedge -}& \\
								& \motivicsymmTspectra &\\
								\pplusoneplusqconnectedsymmTspectra \ar[rr]_-{id} \ar[ur]_-{id}&& 
								\pplusqconnectedsymmTspectra \ar[ul]_-{id}}
		$$
		
	(\ref{thm.3.5.external-pairings-compatible.b}):  This follows from the following commutative diagram of left
	Quillen (bi)functors, together with the construction of the external pairing $\cup ^{c}$ given in theorem
	\ref{thm.3.5.external-pairings}(\ref{thm.3.5.external-pairings.a}) and the construction of the natural
	transformation $\rho$ given in proposition \ref{prop.3.2.functors-between-Rqsymm}(\ref{prop.3.2.functors-between-Rqsymm.c}):
		$$\xymatrix{\pconnectedsymmTspectra \times \qplusoneconnectedsymmTspectra \ar[rr]^-{id\times id} \ar[dr]^-{-\wedge -}
								\ar[dd]_-{-\wedge -}&& \pconnectedsymmTspectra \times
								\qconnectedsymmTspectra \ar[dl]_-{-\wedge -} \ar[dd]^-{-\wedge -}& \\
								& \motivicsymmTspectra &\\
								\pplusoneplusqconnectedsymmTspectra \ar[rr]_-{id} \ar[ur]_-{id}&& 
								\pplusqconnectedsymmTspectra \ar[ul]_-{id}}
		$$
		
	(\ref{thm.3.5.external-pairings-compatible.c}):  This follows from the following commutative diagram of left
	Quillen (bi)functors, together with the construction of the external pairings $\cup ^{c}$, $\cup ^{s}$ given in theorem
	\ref{thm.3.5.external-pairings}(\ref{thm.3.5.external-pairings.a})-({\ref{thm.3.5.external-pairings.b}})
	and the construction of the natural
	transformation $\pi ^{\Sigma}$ given in proposition \ref{prop.3.3.lifting-map-fq-->sqsymm}:
		$$\xymatrix{\motivicsymmTspectra \times \motivicsymmTspectra \ar[r]^-{-\wedge -}& \motivicsymmTspectra \\
								\pconnectedsymmTspectra \times \qconnectedsymmTspectra 
								\ar[u]^-{id\times id} \ar[d]_-{id\times id} \ar[r]^-{-\wedge -}
								& \pplusqconnectedsymmTspectra \ar[u]_-{id} \ar[d]^-{id}\\
								\pslicesymmTspectra \times \qslicesymmTspectra \ar[r]_-{-\wedge -}& \pplusqslicesymmTspectra}
		$$	
\end{proof}

\begin{defi}
		\label{def.3.5.total-slicefunctor}
	Consider the following functors:
		$$\xymatrix@R=.5pt{f^{\Sigma}:\symmstablehomotopy \ar[r]& \symmstablehomotopy \\
											X \ar@{|->}[r]&\bigoplus _{q\in \mathbb Z}f_{q}^{\Sigma}X}
		$$
		
		$$\xymatrix@R=.5pt{s^{\Sigma}:\symmstablehomotopy \ar[r]& \symmstablehomotopy \\
											X \ar@{|->}[r]&\bigoplus _{q\in \mathbb Z}s_{q}^{\Sigma}X}
		$$
\end{defi}

\begin{prop}
		\label{prop.3.5.total-slice====exact-functor}
	\begin{enumerate}
		\item	\label{prop.3.5.total-slice====exact-functor.a}The functor
					$f^{\Sigma}:\symmstablehomotopy \rightarrow \symmstablehomotopy$
					is an exact functor.
		\item	\label{prop.3.5.total-slice====exact-functor.b}The functor
					$s^{\Sigma}:\symmstablehomotopy \rightarrow \symmstablehomotopy$
					is an exact functor.
	\end{enumerate}
\end{prop}
\begin{proof}
	(\ref{prop.3.5.total-slice====exact-functor.a}):Theorem \ref{thm.3.3.symmRq-models-symmf<q}(\ref{thm.3.3.symmRq-models-symmf<q.c}) 
	implies  that all the functors
	$f_{q}^{\Sigma}$ are exact.  Therefore $f^{\Sigma}=\oplus _{q\in \mathbb Z}f_{q}^{\Sigma}$ is
	also an exact functor, since the coproduct of a collection of distinguished triangles is a
	distinguished triangle.
	
	(\ref{prop.3.5.total-slice====exact-functor.b}):Theorem \ref{thm.3.3.symmSq-models-symm-sq}(\ref{thm.3.3.symmSq-models-symm-sq.c}) 
	implies  that all the functors
	$s_{q}^{\Sigma}$ are exact.  Therefore $s^{\Sigma}=\oplus _{q\in \mathbb Z}s_{q}^{\Sigma}$ is
	also an exact functor, since the coproduct of a collection of distinguished triangles is a
	distinguished triangle.
\end{proof}

\begin{thm}
		\label{thm.3.5.qslices-modules-over-zeroslice}
	Fix $q\in \mathbb Z$.  Let $X$ be a ring spectrum in $\symmstablehomotopy$ and let
	$M$ be an $X$-module.
		\begin{enumerate}
			\item \label{thm.3.5.qslices-modules-over-zeroslice.a}The ($-1$)-connective cover of $X$,
						$f_{0}^{\Sigma}X$ 
						(see theorem \ref{thm.3.3.symmRq-models-symmf<q}(\ref{thm.3.3.symmRq-models-symmf<q.c})) 
						also has the structure of a ring
						spectrum in $\symmstablehomotopy$. 
			
			\item \label{thm.3.5.qslices-modules-over-zeroslice.b}The ($q-1$)-connective cover of $M$,
						$f_{q}^{\Sigma}M$ is a module in $\symmstablehomotopy$ over
						the ($-1$)-connective cover of $X$, $f_{0}^{\Sigma}X$.
						
			\item \label{thm.3.5.qslices-modules-over-zeroslice.c}The coproduct of all the connective covers
						of $X$, $f^{\Sigma}X$ has the structure of a graded ring spectrum in $\symmstablehomotopy$.
			
			\item \label{thm.3.5.qslices-modules-over-zeroslice.d}The coproduct of all the connective covers
						of $M$, $f^{\Sigma}M$ is a graded module in $\symmstablehomotopy$ over the
						graded ring $f^{\Sigma}X$.
			
			\item \label{thm.3.5.qslices-modules-over-zeroslice.e}The zero slice of $X$,
						$s_{0}^{\Sigma}X$ 
						(see theorem \ref{thm.3.3.symmSq-models-symm-sq}(\ref{thm.3.3.symmSq-models-symm-sq.c})) 
						also has the structure of a ring spectrum in $\symmstablehomotopy$.
			
			\item \label{thm.3.5.qslices-modules-over-zeroslice.f}The $q$-slice of $M$,
						$s_{q}^{\Sigma}M$ is a module in $\symmstablehomotopy$
						over the zero slice of $X$, 
						$s_{0}^{\Sigma}X$.
						
			\item \label{thm.3.5.qslices-modules-over-zeroslice.g}The coproduct of all the slices
						of $X$, $s^{\Sigma}X$ has the structure of a graded ring spectrum in $\symmstablehomotopy$.
						
			\item \label{thm.3.5.qslices-modules-over-zeroslice.h}The coproduct of all the slices
						of $M$, $s^{\Sigma}M$ is a graded module in $\symmstablehomotopy$ over the
						graded ring $s^{\Sigma}X$.
		\end{enumerate}
\end{thm}
\begin{proof}
	We have that (\ref{thm.3.5.qslices-modules-over-zeroslice.a}) and (\ref{thm.3.5.qslices-modules-over-zeroslice.e})
	follow immediately from proposition
	\ref{prop.3.5.cofibrantreplacementfunctors--strongsymmetricmonoidal} and corollary
	\ref{cor.3.5.fibrantreplacementfunctors--laxsymmetricmonoidal}.  On the other hand,	
	(\ref{thm.3.5.qslices-modules-over-zeroslice.b}), (\ref{thm.3.5.qslices-modules-over-zeroslice.c})
	and (\ref{thm.3.5.qslices-modules-over-zeroslice.d}) follow directly from
	theorem \ref{thm.3.5.external-pairings}(\ref{thm.3.5.external-pairings.a}).
	Finally,
	(\ref{thm.3.5.qslices-modules-over-zeroslice.f}), (\ref{thm.3.5.qslices-modules-over-zeroslice.g})
	and (\ref{thm.3.5.qslices-modules-over-zeroslice.h}) follow
	directly from
	theorem \ref{thm.3.5.external-pairings}(\ref{thm.3.5.external-pairings.b}).
\end{proof}

\begin{thm}
		\label{thm.3.5.qslices-modules-over-zeroslice--spherespectrum}
	Fix $q\in \mathbb Z$, and let $X$ be an arbitrary symmetric $T$-spectrum.
		\begin{enumerate}
			\item \label{thm.3.5.qslices-modules-over-zeroslice--spherespectrum.a}The ($-1$)-connective cover of the sphere spectrum,
						$f_{0}^{\Sigma}\symmspherespectrum$ 
						has the structure of a ring
						spectrum in $\symmstablehomotopy$.
				
			\item \label{thm.3.5.qslices-modules-over-zeroslice--spherespectrum.b}The ($q-1$)-connective cover of $X$,
						$f_{q}^{\Sigma}X$ is a module in $\symmstablehomotopy$ over
						the ($-1$)-connective cover of the sphere spectrum, $f_{0}^{\Sigma}\symmspherespectrum$. 
			
			\item \label{thm.3.5.qslices-modules-over-zeroslice--spherespectrum.c}The coproduct
						of all the connective covers of the sphere spectrum,
						$f^{\Sigma}\symmspherespectrum$ 
						has the structure of a graded ring
						spectrum in $\symmstablehomotopy$.
						
			\item \label{thm.3.5.qslices-modules-over-zeroslice--spherespectrum.d}The coproduct of all the connective covers of $X$,
						$f^{\Sigma}X$ is a graded module in $\symmstablehomotopy$ over
						the graded ring $f^{\Sigma}\symmspherespectrum$.
			
			\item \label{thm.3.5.qslices-modules-over-zeroslice--spherespectrum.e}The zero slice of the sphere spectrum,
						$s_{0}^{\Sigma}\symmspherespectrum$ 
						has the structure of a ring spectrum in $\symmstablehomotopy$.
						
			\item \label{thm.3.5.qslices-modules-over-zeroslice--spherespectrum.f}The $q$-slice of $X$,
						$s_{q}^{\Sigma}X$ is a module in $\symmstablehomotopy$
						over the zero slice of the sphere spectrum, 
						$s_{0}^{\Sigma}\symmspherespectrum$.
						
			\item \label{thm.3.5.qslices-modules-over-zeroslice--spherespectrum.g}The coproduct of all the slices of the sphere spectrum,
						$s^{\Sigma}\symmspherespectrum$ 
						has the structure of a graded ring spectrum in $\symmstablehomotopy$.
						
			\item \label{thm.3.5.qslices-modules-over-zeroslice--spherespectrum.h}The coproduct of all the slices of $X$,
						$s^{\Sigma}X$ is a graded module in $\symmstablehomotopy$
						over the graded ring 
						$s^{\Sigma}\symmspherespectrum$.
		\end{enumerate}
\end{thm}
\begin{proof}
	It is clear that the sphere spectrum $\symmspherespectrum$ is a ring spectrum in
	$\symmstablehomotopy$, and  by construction we have that every symmetric $T$-spectrum $X$
	is a module in $\symmstablehomotopy$ over the sphere spectrum.
	
	The result then follows immediately from theorem \ref{thm.3.5.qslices-modules-over-zeroslice}. 
\end{proof}

	Using the slice filtration, it is possible to construct a spectral sequence which is an analogue of the classical
	Atiyah-Hirzebruch spectral sequence in algebraic topology.
	
\begin{defi}[Motivic Atiyah-Hirzebruch Spectral Sequence]
		\label{def.3.5.Atiyah-Hirzebruch-spectral-sequence}
	Let $X$, $Y$ be a pair of symmetric $T$-spectra.  Then the collection of distinguished triangles
	in $\symmstablehomotopy$ (see theorem \ref{thm.3.1.slicefiltration} and
	propositions \ref{prop.3.2.functors-between-Rqsymm}(\ref{prop.3.2.functors-between-Rqsymm.c}),
	\ref{prop.3.3.lifting-map-fq-->sqsymm}):
		$$\xymatrix{f_{q+1}^{\Sigma}X \ar[r]^-{\rho _{q}^{X}}& f_{q}^{\Sigma}X \ar[r]^-{\pi _{q}^{\Sigma ,X}}& 
								s_{q}^{\Sigma}X \ar[r]^-{\sigma _{q}^{\Sigma ,X}}& \Sigma _{T}^{1,0}f_{q+1}^{\Sigma}X}
		$$
	generates an exact couple $(D_{1}^{p,q}(Y;X),E_{1}^{p,q}(Y;X))$, where:
		\begin{enumerate}
			\item	$D_{1}^{p,q}=[Y,\Sigma _{T}^{p+q,0}f_{p}^{\Sigma}X]^{\Sigma}_{Spt}$, and
			\item	$E_{1}^{p,q}(Y;X)=[Y,\Sigma _{T}^{p+q,0}s_{p}^{\Sigma}X]^{\Sigma}_{Spt}$.
		\end{enumerate}
\end{defi}

	The compatibility of the slice filtration with the smash product of symmetric $T$-spectra implies that
	the smash product of symmetric $T$-spectra induces a pairing of spectral sequences:
	
\begin{thm}
		\label{thm.3.5.pairings-Atiyah-Hirzebruch-ss}
	Let $X$, $X'$, $Y$, $Y'$ be symmetrict $T$-spectra.  Then the smash product of symmetric $T$-spectra induces
	the following natural external pairings in the motivic Atiyah-Hirzebruch spectral sequence:
		$$\xymatrix@R=0.5pt{E_{r}^{p,q}(Y;X)\otimes E_{r}^{p',q'}(Y';X') \ar[r]& E_{r}^{p+p',q+q'}(Y\wedge Y';X\wedge X')\\
												(\alpha, \beta) \ar@{|->}[r]& \alpha \smile \beta }
		$$
	where $\alpha :Y\rightarrow \Sigma _{T}^{p+q,0}s^{\Sigma}_{p}X$, $\beta :Y'\rightarrow \Sigma _{T}^{p'+q',0}s^{\Sigma}_{p'}X'$
	and $\alpha \smile \beta$ is the following composition
	(see theorem \ref{thm.3.5.external-pairings}(\ref{thm.3.5.external-pairings.b})):
		$$\xymatrix{Y\wedge Y' \ar[r]^-{\alpha \wedge \beta}& \Sigma _{T}^{p+q,0}s^{\Sigma}_{p}X \wedge 
												\Sigma _{T}^{p'+q',0}s^{\Sigma}_{p'}X'
												\ar[d]^-{\cong} &&\\
												&\Sigma _{T}^{p+p'+q+q',0}s^{\Sigma}_{p}X\wedge s^{\Sigma}_{p'}X 
												\ar[drr]^-{\Sigma _{T}^{p+p'+q+q',0}\circ \cup ^{s}_{p,p'}} &&\\
												&&& \Sigma _{T}^{p+p'+q+q',0}s^{\Sigma}_{p+p'}X\wedge X'}
		$$
\end{thm}
\begin{proof}
	Using the naturality of the external pairings $\cup ^{c}_{p,q}$, $\cup ^{s}_{p,q}$
	(see theorem \ref{thm.3.5.external-pairings}) and theorem \ref{thm.3.5.external-pairings-compatible}, 
	the result follows immediately from the work
	of Massey \cite{MR0060829} together with \cite[proposition 14.3]{MR1949356}.
\end{proof}

\begin{defi}
		\label{def.3.5.slice-cover-with-coefficients}
	Fix $q\in \mathbb Z$.  Let $A$ be a cofibrant ring spectrum
	with unit in $\motivicsymmTspectra$.  
		\begin{enumerate}
			\item \label{def.3.5.slice-cover-with-coefficients.a}  Let $f_{q}^{m}$ denote the following composition of exact
						functors between triangulated categories (see 
						proposition \ref{prop.3.5.symmcofibrant-replacement=>triangulatedfunctor})
							$$\xymatrix{\Amodstablehomotopy \ar[dr]_-{R_{m}} \ar[rr]^-{f_{q}^{m}} && \Amodstablehomotopy \\
													 &\qconnectedAmodstablehomotopy \ar[ur]_-{C_{q}^{m}} &}
							$$
			\item \label{def.3.5.slice-cover-with-coefficients.b}  Let $s_{<q}^{m}$ denote the following composition of exact
						functors between triangulated categories (see
						proposition \ref{prop.3.5.symmLq-exact-adjunctions})
							$$\xymatrix{\Amodstablehomotopy \ar[dr]_-{Q_{m}} \ar[rr]^-{s_{<q}^{m}} && \Amodstablehomotopy \\
													 &\weightqAmodstablehomotopy \ar[ur]_-{W_{q}^{m}} &}
							$$
			\item \label{def.3.5.slice-cover-with-coefficients.c}  Let $s_{q}^{m}$ denote the following composition of exact
						functors between triangulated categories (see propositions \ref{prop.3.5.symmcofibrant-replacement=>triangulatedfunctor}
						and \ref{prop.3.5.symmetricq-connected--->q-slice===leftQuilllenfunctor})
							$$\xymatrix{\Amodstablehomotopy \ar[d]_-{R_{m}} \ar[rr]^-{s_{q}^{m}}&& \Amodstablehomotopy\\
													\qconnectedAmodstablehomotopy \ar[r]_-{C_{q}^{m}} & \qsliceAmodstablehomotopy \ar[r]_-{W_{q+1}^{m}} 
													& \qconnectedAmodstablehomotopy \ar[u]_-{C_{q}^{m}}}
							$$
		\end{enumerate}
\end{defi}

\begin{rmk}
		\label{rmk.3.5.colocal===>covers-slices-inherit-module-structures}
	Notice that the following two theorems (\ref{thm.3.5.colocal===>covers-slices-inherit-module-structures}
	and \ref{thm.3.5.change-of-coefficients})
	are much stronger than theorem \ref{thm.3.5.qslices-modules-over-zeroslice};
	since the module structures in the latter are defined just up to homotopy 
	(i.e. they make sense in $\symmstablehomotopy$), whereas the module structures in the first case
	are strict (i.e. they are defined in the
	model category $\motivicsymmTspectra$).
\end{rmk}

\begin{thm}
		\label{thm.3.5.colocal===>covers-slices-inherit-module-structures}
	Fix $q\in \mathbb Z$.  Let $A$ be a cofibrant ring spectrum 
	with unit in $\motivicsymmTspectra$.
		\begin{enumerate}
			\item \label{thm.3.5.colocal===>covers-slices-inherit-module-structures.a}  If $A$ is cofibrant 
						in $\zeroconnectedsymmTspectra$, then
						the functor 
						$f_{q}^{\Sigma}\circ UR_{m}$ (see theorems \ref{thm.3.3.symmRq-models-symmf<q} and \ref{thm.3.5.freemodulefunctor-leftQuillen})
							$$\xymatrix{\Amodstablehomotopy \ar[r]^-{UR_{m}}& \symmstablehomotopy \ar[r]^-{f_{q}^{\Sigma}} &
													\symmstablehomotopy}$$
						factors through $\Amodstablehomotopy$ (see definition 
						\ref{def.3.5.slice-cover-with-coefficients}(\ref{def.3.5.slice-cover-with-coefficients.a}))
							$$\xymatrix{\Amodstablehomotopy \ar[r]^-{UR_{m}} \ar@{-->}[d]_-{f_{q}^{m}} & \symmstablehomotopy \ar[d]^-{f_{q}^{\Sigma}}\\
													\Amodstablehomotopy \ar[r]_-{UR_{m}} & \symmstablehomotopy}$$
						i.e. for every $A$-module
						$M$, its ($q-1$)-connective cover $f_{q}^{\Sigma}(M)$ inherits a natural strict
						structure of $A$-module in $\motivicsymmTspectra$.
			\item	\label{thm.3.5.colocal===>covers-slices-inherit-module-structures.b}  If $A$ is cofibrant 
						in $\zeroslicesymmTspectra$, then
						the functor $s_{<q}^{\Sigma}\circ UR_{m}$ (see theorems \ref{thm.3.3.symmLq-models-symms<q}
						and \ref{thm.3.5.freemodulefunctor-leftQuillen})
							$$\xymatrix{\Amodstablehomotopy \ar[r]^-{UR_{m}} & \symmstablehomotopy \ar[r]^-{s_{<q}^{\Sigma}} & 
													\symmstablehomotopy}$$
						factors through $\Amodstablehomotopy$ (see definition 
						\ref{def.3.5.slice-cover-with-coefficients}(\ref{def.3.5.slice-cover-with-coefficients.b}))
							$$\xymatrix{\Amodstablehomotopy \ar[r]^-{UR_{m}} \ar@{-->}[d]_-{s_{<q}^{m}} & \symmstablehomotopy \ar[d]^-{s_{<q}^{\Sigma}}\\
													\Amodstablehomotopy \ar[r]_-{UR_{m}}& \symmstablehomotopy}$$
						i.e. for every $A$-module
						$M$, $s_{<q}^{\Sigma}(M)$ inherits a natural strict
						structure of $A$-module in $\motivicsymmTspectra$.
			\item	\label{thm.3.5.colocal===>covers-slices-inherit-module-structures.c}  If $A$ is cofibrant 
						in $\zeroconnectedsymmTspectra$, then
						the functor $s_{q}^{\Sigma}\circ UR_{m}$ (see theorems \ref{thm.3.3.symmSq-models-symm-sq}
						and \ref{thm.3.5.freemodulefunctor-leftQuillen})
						  $$\xymatrix{\Amodstablehomotopy \ar[r]^-{UR_{m}} & \symmstablehomotopy \ar[r]^-{s_{q}^{\Sigma}} &
						  						\symmstablehomotopy}$$
						factors through $\Amodstablehomotopy$ (see definition 
						\ref{def.3.5.slice-cover-with-coefficients}(\ref{def.3.5.slice-cover-with-coefficients.c}))
							$$\xymatrix{\Amodstablehomotopy \ar[r]^-{UR_{m}} \ar@{-->}[d]_-{s_{q}^{m}} & \symmstablehomotopy \ar[d]^-{s_{q}^{\Sigma}}\\
													 \Amodstablehomotopy \ar[r]_-{UR_{m}} & \symmstablehomotopy}$$
						i.e. for every $A$-module
						$M$, its $q$-slice $s_{q}^{\Sigma}(M)$ inherits a natural strict
						structure of $A$-module in $\motivicsymmTspectra$.
		\end{enumerate}
\end{thm}
\begin{proof}
	(\ref{thm.3.5.colocal===>covers-slices-inherit-module-structures.a}):  By construction
	(see theorem \ref{thm.3.3.symmRq-models-symmf<q}) the functor $f_{q}^{\Sigma}$
	is defined as the following composition
		$$\xymatrix{\symmstablehomotopy \ar[dr]_-{R_{\Sigma}} \ar[rr]^-{f_{q}^{\Sigma}} && \symmstablehomotopy \\
													 &\qconnectedsymmstablehomotopy \ar[ur]_-{C_{q}^{\Sigma}} &}
							$$
	Since we are assuming that $A$ is cofibrant in $\zeroconnectedsymmTspectra$
	(equivalently $C_{eff}^{0,\Sigma}$-colocal in $\motivicsymmTspectra$),
	the result follows directly from diagram (\ref{diagram3.5.liftslicefiltrationsymm.c})
	in proposition \ref{prop.3.5.functors-between-Rqsymm} and
	theorem \ref{thm.3.5.inheritingmodulestructures1}.
	
	(\ref{thm.3.5.colocal===>covers-slices-inherit-module-structures.b}):  By construction
	(see theorem \ref{thm.3.3.symmLq-models-symms<q}) the functor $s_{<q}^{\Sigma}$
	is defined as the following composition
		$$\xymatrix{\symmstablehomotopy \ar[dr]_-{Q_{\Sigma}} \ar[rr]^-{s_{<q}^{\Sigma}} && \symmstablehomotopy \\
													 &\weightqsymmstablehomotopy \ar[ur]_-{W_{q}^{\Sigma}} &}
							$$
	Since we are assuming that $A$ is cofibrant in $\zeroslicesymmTspectra$,
	the result follows directly from theorem \ref{thm.3.5.inheritmodelstructures2}
	and diagram (\ref{diagram.prop.3.5.Lq+1-->Lqsymm.c})
	in proposition \ref{prop.3.5.Lq+1-->Lqsymm}.
	
	(\ref{thm.3.5.colocal===>covers-slices-inherit-module-structures.c}):  By construction
	(see theorem \ref{thm.3.3.symmSq-models-symm-sq}) the functor $s_{q}^{\Sigma}$
	is defined as the following composition
		$$\xymatrix{\symmstablehomotopy \ar[d]_-{R_{\Sigma}} \ar[rr]^-{s_{q}^{\Sigma}}&& \symmstablehomotopy\\
													\qconnectedsymmstablehomotopy \ar[r]_-{C_{q}^{\Sigma}} & \qslicesymmstablehomotopy \ar[r]_-{W_{q+1}^{\Sigma}} 
													& \qconnectedsymmstablehomotopy \ar[u]_-{C_{q}^{\Sigma}}}
							$$
	Since we are assuming that $A$ is cofibrant in $\zeroconnectedsymmTspectra$
	(equivalently $C_{eff}^{0,\Sigma}$-colocal in $\motivicsymmTspectra$),
	the result is a consequence of diagram (\ref{diagram3.5.liftslicefiltrationsymm.c})
	in proposition \ref{prop.3.5.functors-between-Rqsymm},
	theorem \ref{thm.3.5.inheritmodelstructures-slices},
	diagram (\ref{diagram3.5.homotopycoherence==>liftings-qslice.c}) 
	in proposition \ref{prop.3.5.homotopycoherence==>liftings-qslice}
	and theorem \ref{thm.3.5.inheritingmodulestructures1}.
\end{proof}

\begin{thm}
		\label{thm.3.5.change-of-coefficients}
	Fix $q\in \mathbb Z$.  Let $A$ be a cofibrant ring spectrum 
	with unit in $\motivicsymmTspectra$.
	Consider the following composition of exact functors between triangulated categories 
	(see proposition \ref{prop.3.3.symmcofibrant-replacement=>triangulatedfunctor},
	corollary \ref{cor.3.3.adjunctions--symmetricRq==>Sq},
	theorem \ref{thm.3.5.symmetrization-qslice-Quillen-equiv} and propositions
	\ref{prop.3.5.symmcofibrant-replacement=>triangulatedfunctor},
	\ref{prop.3.5.symmetricq-connected--->q-slice===leftQuilllenfunctor})
		\begin{equation}
					\label{diagram.thm.3.5.change-of-coefficients}
			\begin{array}{c}
				\xymatrix{\symmstablehomotopy \ar[d]_-{R_{\Sigma}} && \symmstablehomotopy\\
									\qconnectedsymmstablehomotopy \ar[d]_-{C_{q}^{\Sigma}}  && \Amodstablehomotopy \ar[u]_-{UR_{m}}\\
									\qslicesymmstablehomotopy \ar[r]_-{A\wedge P_{q}^{\Sigma}-} & \qsliceAmodstablehomotopy \ar[r]_-{W_{q+1}^{m}} 
									& \qconnectedAmodstablehomotopy \ar[u]_-{C_{q}^{m}}}
			\end{array}
		\end{equation}
	Furthermore, assume that $A$ is cofibrant in $\zeroconnectedsymmTspectra$ and the unit map 
	$u:\symmspherespectrum \rightarrow A$
	is a weak equivalence in $\zeroslicesymmTspectra$.
	Then $u$ induces a natural equivalence between $s_{q}^{\Sigma}$
	(see theorem \ref{thm.3.3.symmSq-models-symm-sq}) and the functor defined above
	in diagram (\ref{diagram.thm.3.5.change-of-coefficients}),
	i.e. for every symmetric $T$-spectrum $X$, its $q$-slice 
	$s_{q}^{\Sigma}(X)$ is equipped with a natural strict structure of $A$-module
	in $\motivicsymmTspectra$.
\end{thm}
\begin{proof}
	The functor $s_{q}^{\Sigma}$
	(see theorem \ref{thm.3.3.symmSq-models-symm-sq}) 
	is defined as the following composition
		$$\xymatrix{\symmstablehomotopy \ar[d]_-{R_{\Sigma}} \ar[rr]^-{s_{q}^{\Sigma}}&& \symmstablehomotopy\\
								\qconnectedsymmstablehomotopy \ar[r]_-{C_{q}^{\Sigma}} & \qslicesymmstablehomotopy \ar[r]_-{W_{q+1}^{\Sigma}} 
								& \qconnectedsymmstablehomotopy \ar[u]_-{C_{q}^{\Sigma}}}
							$$
	By hypothesis $A$ is cofibrant in $\zeroconnectedsymmTspectra$, and
	lemma \ref{lem.3.5.spherespectrum-0-connected-0-slice-cofibrant} implies that 
	$\symmspherespectrum$ is also cofibrant in $\zeroconnectedsymmTspectra$.
	Since the unit map $u:\symmspherespectrum \rightarrow A$ is assumed to be a weak equivalence
	in $\zeroslicesymmTspectra$, it follows from theorem 
	\ref{thm.3.5.sliceinvarianceofcoefficients}(\ref{thm.3.5.sliceinvarianceofcoefficients.d}) that
	the adjunction 
		$$(A\wedge-, U,\varphi):\qslicesymmTspectra \rightarrow \qsliceAmod$$
	is a Quillen equivalence.  Therefore the functor $s_{q}^{\Sigma}$ is naturally isomorphic to
	the following composition
		$$\xymatrix{\symmstablehomotopy \ar[d]_-{R_{\Sigma}} && \symmstablehomotopy\\
								\qconnectedsymmstablehomotopy \ar[d]_-{C_{q}^{\Sigma}}  && \qconnectedsymmstablehomotopy \ar[u]_-{C_{q}^{\Sigma}}\\
								\qslicesymmstablehomotopy \ar[r]_-{A\wedge P_{q}^{\Sigma}-} & \qsliceAmodstablehomotopy \ar[r]_-{UW_{q+1}^{m}} 
								& \qslicesymmstablehomotopy \ar[u]_-{W_{q+1}^{\Sigma}}}
							$$
	Now, proposition \ref{prop.3.3.symmetricq-connected--->q-slice===leftQuilllenfunctor} 
	implies that $A$ is cofibrant in $\zeroslicesymmTspectra$, therefore
	using diagram (\ref{diagram3.5.homotopycoherence==>liftings-qslice.c})
	in proposition \ref{prop.3.5.homotopycoherence==>liftings-qslice}, 
	we get that the functor $s_{q}^{\Sigma}$ becomes naturally isomorphic
	to the following composition
		$$\xymatrix{\symmstablehomotopy \ar[d]_-{R_{\Sigma}} && \symmstablehomotopy\\
								\qconnectedsymmstablehomotopy \ar[d]_-{C_{q}^{\Sigma}}  && \qconnectedsymmstablehomotopy \ar[u]_-{C_{q}^{\Sigma}}\\
								\qslicesymmstablehomotopy \ar[r]_-{A\wedge P_{q}^{\Sigma}-} & \qsliceAmodstablehomotopy \ar[r]_-{W_{q+1}^{m}} 
								& \qconnectedAmodstablehomotopy \ar[u]_-{UR_{m}}}
							$$
	Finally, since $A$ is cofibrant in $\zeroconnectedsymmTspectra$ we can apply
	theorem \ref{thm.3.5.inheritingmodulestructures1} and we get that $s_{q}^{\Sigma}$
	is naturally isomorphic to the following composition
		$$\xymatrix{\symmstablehomotopy \ar[d]_-{R_{\Sigma}} && \symmstablehomotopy\\
								\qconnectedsymmstablehomotopy \ar[d]_-{C_{q}^{\Sigma}}  && \Amodstablehomotopy \ar[u]_-{UR_{m}}\\
								\qslicesymmstablehomotopy \ar[r]_-{A\wedge P_{q}^{\Sigma}-} & \qsliceAmodstablehomotopy \ar[r]_-{W_{q+1}^{m}} 
								& \qconnectedAmodstablehomotopy \ar[u]_-{C_{q}^{m}}}
							$$
	This finishes the proof.
\end{proof}

\begin{lem}
		\label{lemma.3.6.connecting-functors-with-model-structures}
	Fix $q\in \mathbb Z$.  Let $g:X\rightarrow Y$ be a map
	between cofibrant spectra in $\motivicsymmTspectra$.
		\begin{enumerate}
			\item	\label{lemma.3.6.connecting-functors-with-model-structures.a}  We have that $X$ is cofibrant
										in $\qconnectedsymmTspectra$ if and only if the natural map
										$\theta _{q}^{\Sigma ,X}:f_{q}^{\Sigma}(X)\rightarrow X$ (see 
										proposition \ref{prop.3.3.lifting-map-fq-->idsymm}) is an isomorphism
										in $\symmstablehomotopy$.
			\item	\label{lemma.3.6.connecting-functors-with-model-structures.b}  The map $g$
										is a weak equivalence in $\qconnectedsymmTspectra$ if and only if
										the induced map $f_{q}^{\Sigma}(g):f_{q}^{\Sigma}(X) \rightarrow f_{q}^{\Sigma}(Y)$
										is a weak equivalence in $\motivicsymmTspectra$.
			\item	\label{lemma.3.6.connecting-functors-with-model-structures.c}  If $X\cong s_{q}^{\Sigma}(X)$
										in $\symmstablehomotopy$, then $X$ is cofibrant in $\qconnectedsymmTspectra$
										and $X\cong f_{q}^{\Sigma}(X)$ in $\symmstablehomotopy$.
			\item	\label{lemma.3.6.connecting-functors-with-model-structures.d}  Furthermore, assume that
										$X$, $Y$ are both cofibrant in $\qconnectedsymmTspectra$.  Then $g$ is a weak equivalence
										in $\qslicesymmTspectra$ if and only if the induced map 
										$s_{q}^{\Sigma}(g):s_{q}^{\Sigma}(X) \rightarrow s_{q}^{\Sigma}(Y)$
										is a weak equivalence in $\motivicsymmTspectra$.
		\end{enumerate}
\end{lem}
\begin{proof}
	(\ref{lemma.3.6.connecting-functors-with-model-structures.a}):  Proposition \ref{prop.3.3.lifting-map-fq-->idsymm}
	implies that the natural map $\theta _{q}^{\Sigma}$ is just the counit
	of the adjunction (see proposition \ref{prop.3.3.symmcofibrant-replacement=>triangulatedfunctor})
		$$(C_{q}^{\Sigma},R_{\Sigma},\varphi ):\qconnectedsymmstablehomotopy \rightarrow \symmstablehomotopy
		$$ 
	Consider the following diagram in $\motivicsymmTspectra$
		$$\xymatrix{X \ar[r]^-{R_{\Sigma}^{X}}& R_{\Sigma}X & C_{q}^{\Sigma}R_{\Sigma}X \ar[l]_-{C_{q}^{\Sigma, R_{\Sigma}X}}}
		$$
	By construction (see definition \ref{def.3.3.fibrantstable-replacementfunctor}) 
	$R_{\Sigma}^{X}$ is always a weak equivalence in $\motivicsymmTspectra$;
	therefore, $\theta _{q}^{\Sigma ,X}$ is an isomorphism in $\symmstablehomotopy$
	if and only if $C_{q}^{\Sigma ,R_{\Sigma}X}$ is a weak equivalence in $\motivicsymmTspectra$.
	
	On the other hand, $X$ is cofibrant in $\motivicsymmTspectra$ and by construction $R_{\Sigma}^{X}$
	is a trivial cofibration in $\motivicsymmTspectra$ (see definition \ref{def.3.3.fibrantstable-replacementfunctor}), 
	hence we get that $R_{\Sigma}^{X}$ is cofibrant in $\motivicsymmTspectra$. Therefore
	\cite[proposition 3.2.2]{MR1944041} implies that $X$ is cofibrant in $\qconnectedsymmTspectra$
	if and only if $R_{\Sigma}X$ is cofibrant in $\qconnectedsymmTspectra$
	
	Finally, $C_{q}^{\Sigma}$ is a cofibrant replacement functor in $\qconnectedsymmTspectra$
	(see definition \ref{def.3.3.Cqsigma-cofibrant-replacement}).
	Hence, \cite[theorem 3.2.13(2)]{MR1944041} and \cite[proposition 3.2.2]{MR1944041} imply that 
	$C_{q}^{\Sigma ,R_{\Sigma}X}$ is a weak equivalence in $\motivicsymmTspectra$
	if and only if $R_{\Sigma}X$ is cofibrant in $\qconnectedsymmTspectra$.
	
	(\ref{lemma.3.6.connecting-functors-with-model-structures.b}):  By construction, we have that
	$f_{q}^{\Sigma}=C_{q}^{\Sigma}\circ R_{\Sigma}$ (see theorem \ref{thm.3.3.symmRq-models-symmf<q}).
	Consider the following commutative
	diagram in $\qconnectedsymmTspectra$
		$$\xymatrix{X \ar[r]^-{g} \ar[d]_-{R_{\Sigma}^{X}}& Y \ar[d]^-{R_{\Sigma}^{Y}}\\
								R_{\Sigma}X \ar[r]^-{R_{\Sigma}(g)}& R_{\Sigma}Y\\
								C_{q}^{\Sigma}R_{\Sigma}X \ar[u]^-{C_{q}^{\Sigma ,R_{\Sigma}X}} \ar[r]_-{C_{q}^{\Sigma}R_{\Sigma}(g)}& 
								C_{q}^{\Sigma}R_{\Sigma}Y \ar[u]_-{C_{q}^{\Sigma ,R_{\Sigma}Y}}}
		$$
	Proposition \ref{prop.3.3.Rsigma-fibrant-replacement-all-Rq} and 
	definition \ref{def.3.3.Cqsigma-cofibrant-replacement} 
	imply that all the vertical arrows are weak equivalences in $\qconnectedsymmTspectra$.
	Hence, using the two out of three property for weak equivalences we get that
	the top row is a weak equivalence in $\qconnectedsymmTspectra$
	if and only if the bottom row is a weak equivalence in $\qconnectedsymmTspectra$.
	
	On the other hand, by construction $C_{q}^{\Sigma}R_{\Sigma}X$, $C_{q}^{\Sigma}R_{\Sigma}Y$
	are both cofibrant in $\qconnectedsymmTspectra$ (see definition \ref{def.3.3.Cqsigma-cofibrant-replacement});
	thus,
	\cite[theorem 3.2.13(2)]{MR1944041} and \cite[proposition 3.1.5]{MR1944041}
	imply that $C_{q}^{\Sigma}R_{\Sigma}(g)$ is a weak equivalence in $\qconnectedsymmTspectra$
	if and only if $C_{q}^{\Sigma}R_{\Sigma}(g)$ is a weak equivalence in $\motivicsymmTspectra$.
	
	(\ref{lemma.3.6.connecting-functors-with-model-structures.c}):  By (\ref{lemma.3.6.connecting-functors-with-model-structures.a})
	above, it suffices to show that $X$ is cofibrant in $\qconnectedsymmTspectra$.
	Since we are assuming that $X$ is cofibrant in $\motivicsymmTspectra$
	and $X\cong s_{q}^{\Sigma}(X)$ in $\symmstablehomotopy$, \cite[proposition 3.2.2]{MR1944041}
	implies that it is enough to check that $s_{q}^{\Sigma}(X)$ is cofibrant in $\qconnectedsymmTspectra$.
	
	However, by definition $s_{q}^{\Sigma}=C_{q}^{\Sigma}\circ W_{q+1}^{\Sigma} \circ C_{q}^{\Sigma} \circ R_{\Sigma}$
	(see theorem \ref{thm.3.3.symmSq-models-symm-sq}), 
	and by construction $C_{q}^{\Sigma}$ is a cofibrant replacement functor 
	in $\qconnectedsymmTspectra$ (see definition \ref{def.3.3.Cqsigma-cofibrant-replacement}).
	Therefore, $s_{q}^{\Sigma}(X)$ is always cofibrant in $\qconnectedsymmTspectra$, as we wanted.
	
	(\ref{lemma.3.6.connecting-functors-with-model-structures.d}):  By construction, we have that
	$s_{q}^{\Sigma}=C_{q}^{\Sigma}\circ W_{q+1}^{\Sigma} \circ C_{q}^{\Sigma} \circ R_{\Sigma}$
	(see theorem \ref{thm.3.3.symmSq-models-symm-sq}).
	Consider the following commutative
	diagram in $\motivicsymmTspectra$
		$$\xymatrix{X \ar[rr]^-{g} \ar[d]_-{R_{\Sigma}^{X}}&& Y \ar[d]^-{R_{\Sigma}^{Y}}\\
								R_{\Sigma}X \ar[rr]^-{R_{\Sigma}(g)}&& R_{\Sigma}Y\\
								C_{q}^{\Sigma}R_{\Sigma}X \ar[u]^-{C_{q}^{\Sigma ,R_{\Sigma}X}} 
								\ar[d]_-{W_{q+1}^{\Sigma, C_{q}^{\Sigma}R_{\Sigma}X}} \ar[rr]_-{C_{q}^{\Sigma}R_{\Sigma}(g)}&& 
								C_{q}^{\Sigma}R_{\Sigma}Y \ar[u]_-{C_{q}^{\Sigma ,R_{\Sigma}Y}} 
								\ar[d]^-{W_{q+1}^{\Sigma, C_{q}^{\Sigma}R_{\Sigma}Y}}\\
								W_{q+1}^{\Sigma}C_{q}^{\Sigma}R_{\Sigma}X  
								\ar[rr]_-{W_{q+1}^{\Sigma}C_{q}^{\Sigma}R_{\Sigma}(g)}&& 
								W_{q+1}^{\Sigma}C_{q}^{\Sigma}R_{\Sigma}Y \\
								C_{q}^{\Sigma}W_{q+1}^{\Sigma}C_{q}^{\Sigma}R_{\Sigma}X \ar[u]^-{C_{q}^{\Sigma ,W_{q+1}^{\Sigma}C_{q}^{\Sigma}R_{\Sigma}X}} 
								\ar[rr]_-{s_{q}^{\Sigma}(g)}&& 
								C_{q}^{\Sigma}W_{q+1}^{\Sigma}C_{q}^{\Sigma}R_{\Sigma}Y \ar[u]_-{C_{q}^{\Sigma ,W_{q+1}^{\Sigma}C_{q}^{\Sigma}R_{\Sigma}Y}}}
		$$
	We claim that $C_{q}^{\Sigma}R_{\Sigma}(g)$ is a weak equivalence
	in $\qslicesymmTspectra$ if and only if $s_{q}^{\Sigma}(g)$ is a weak equivalence in $\motivicsymmTspectra$.
	In effect, corollary \ref{cor.3.3.classifying-Sq-symmetriccolocal-equivs.b} 
	implies that $C_{q}^{\Sigma}R_{\Sigma}(g)$ is a weak equivalence in
	$\qslicesymmTspectra$ if and only if $W_{q+1}^{\Sigma}C_{q}^{\Sigma}R_{\Sigma}(g)$
	is a weak equivalence in $\qconnectedsymmTspectra$.  But $C_{q}^{\Sigma}$
	is by construction a cofibrant replacement functor in $\qconnectedsymmTspectra$
	(see definition \ref{def.3.3.Cqsigma-cofibrant-replacement}); thus, $W_{q+1}^{\Sigma}C_{q}^{\Sigma}R_{\Sigma}(g)$
	is a weak equivalence in $\qconnectedsymmTspectra$ if and only if $s_{q}^{\Sigma}(g)$ is a
	weak equivalence in $\qconnectedsymmTspectra$.  Finally, since $s_{q}^{\Sigma}(X)$, $s_{q}^{\Sigma}(Y)$
	are always cofibrant in $\qconnectedsymmTspectra$, we have that
	\cite[theorem 3.2.13(2)]{MR1944041} and \cite[proposition 3.1.5]{MR1944041} imply that
	$s_{q}^{\Sigma}(g)$ is a weak equivalence in $\qconnectedsymmTspectra$
	if and only if $s_{q}^{\Sigma}(g)$ is a weak equivalence in $\motivicsymmTspectra$.
	
	Now, the two out of three property for weak equivalences implies that
	it is enough to show that the maps $R_{\Sigma}^{X}$, $R_{\Sigma}^{Y}$,
	$C_{q}^{\Sigma ,R_{\Sigma}X}$ and $C_{q}^{\Sigma ,R_{\Sigma}Y}$ are
	all weak equivalences in $\qslicesymmTspectra$.  But $\qslicesymmTspectra$
	is a right Bousfield localization with respect to $\weightqplusonesymmTspectra$,
	and similarly $\weightqplusonesymmTspectra$ is a left Bousfield localization with respect to
	$\motivicsymmTspectra$; thus, \cite[proposition 3.1.5]{MR1944041} implies that
	it is enough to check that $R_{\Sigma}^{X}$, $R_{\Sigma}^{Y}$,
	$C_{q}^{\Sigma ,R_{\Sigma}X}$ and $C_{q}^{\Sigma ,R_{\Sigma}Y}$ are weak equivalences in 
	$\motivicsymmTspectra$.
	
	By construction the maps $R_{\Sigma}^{X}$, $R_{\Sigma}^{Y}$ are trivial cofibrations in
	$\motivicsymmTspectra$ (see definition \ref{def.3.3.fibrantstable-replacementfunctor}); 
	hence, they are in particular weak equivalences in 
	$\motivicsymmTspectra$, and $R_{\Sigma}X$, $R_{\Sigma}Y$ are both cofibrant in
	$\motivicsymmTspectra$ since we are assuming that $X$ and $Y$ are cofibrant in $\motivicsymmTspectra$.
	Now \cite[proposition 3.2.2]{MR1944041} implies that $R_{\Sigma}X$ and $R_{\Sigma}Y$
	are also cofibrant in $\qconnectedsymmTspectra$, since by hypothesis $X$, $Y$ are both
	cofibrant in $\qconnectedsymmTspectra$.
	Therefore,  \cite[theorem 3.2.13(2)]{MR1944041} and \cite[proposition 3.1.5]{MR1944041} imply that
	$C_{q}^{\Sigma ,R_{\Sigma}X}$ and $C_{q}^{\Sigma ,R_{\Sigma}Y}$ are weak equivalences in 
	$\motivicsymmTspectra$, if and only if they are weak equivalences in $\qconnectedsymmTspectra$,
	but this is clear since $C_{q}^{\Sigma}$ is a cofibrant replacement functor in $\qconnectedsymmTspectra$
	(see definition \ref{def.3.3.Cqsigma-cofibrant-replacement}).
\end{proof}

	The next theorem proves a conjecture of M. Levine \cite[corollary 11.1.3]{MR2365658}.
	
\begin{thm}
		\label{thm.3.5.q-slices=====HZ-modules}
	Fix $q\in \mathbb Z$.  Let $H\mathbb Z$ denote the motivic Eilenberg-MacLane spectrum
	in $\motivicsymmTspectra$ (see \cite[example 8.2.2(2)]{MR2365658}), and
	assume that the base scheme is a perfect field $k$. 
	Then for every
	symmetric $T$-spectrum $X$ in $\ksymmstablehomotopy$:
		\begin{itemize}
			\item	The $q$-slice of $X$, $s_{q}^{\Sigma}X$ has a natural
						structure of $H\mathbb Z$-module in $\motivicsymmTspectra$, i.e. $s_{q}^{\Sigma}X$
						is in a natural way an object in the motivic stable homotopy
						category of $H\mathbb Z$-modules
						$\stableHZmodules$.
		\end{itemize}
\end{thm}
\begin{proof}
	The work of M. Levine (see \cite[theorem 10.5.1]{MR2365658}) implies
	that $s_{0}^{\Sigma}(u)$ is a weak equivalence in $\motivicsymmTspectra$, where
	$u$ denotes the unit map $u:\symmspherespectrum \rightarrow H\mathbb Z$
	for the commutative ring spectrum $H\mathbb Z$ in $\motivicsymmTspectra$.
	
	On the other hand, theorem \ref{thm.2.8.motivic-Aalgebras}, proposition \ref{prop.2.8.Aalg-cofibrations=>motivic-cofibrations}
	and lemma \ref{lem.3.5.spherespectrum-0-connected-0-slice-cofibrant},
	imply that we can assume that $H\mathbb Z$ is
	cofibrant in $\motivicsymmTspectra$.
	
	Furthermore,  lemma 10.4.1 in \cite{MR2365658} implies that 
	$s_{0}^{\Sigma}(H\mathbb Z)\cong H\mathbb Z$ in $\symmstablehomotopy$; hence by 
	lemma \ref{lemma.3.6.connecting-functors-with-model-structures}(\ref{lemma.3.6.connecting-functors-with-model-structures.c})
	we get that $H\mathbb Z$ is cofibrant in $\zeroconnectedsymmTspectra$.
	On the other hand, lemma \ref{lem.3.5.spherespectrum-0-connected-0-slice-cofibrant} implies that 
	$\symmspherespectrum$ is also cofibrant in $\zeroconnectedsymmTspectra$.
	
	Therefore, lemma 
	\ref{lemma.3.6.connecting-functors-with-model-structures}(\ref{lemma.3.6.connecting-functors-with-model-structures.d}) implies that
	$u:\symmspherespectrum \rightarrow H\mathbb Z$ is a weak equivalence in
	$\zeroslicesymmTspectra$.
	Thus, the result follows immediately from theorem \ref{thm.3.5.change-of-coefficients}.
\end{proof}

	The motivic stable model category of $H\mathbb Z$-modules has been studied
	in detail by R{\"o}ndigs and {\O}stv{\ae}r (see \cite{MR2435654}), as
	a consequence of their work we get that the slices may be interpreted as motives
	in the sense of Voevodsky.
	
\begin{thm}
		\label{thm.3.5.q-slices=====big-motives}
	Let $k$ be a field of characteristic zero.  Then for every $q\in \mathbb Z$ and for every
	symmetric $T$-spectrum $X$ in $\ksymmstablehomotopy$:
		\begin{itemize}
			\item	The $q$-slice of $X$, $s_{q}^{\Sigma}X$
						is a big motive (see \cite{MR1764202}, \cite[section 2.3]{MR2435654}) 
						in the sense of Voevodsky .
		\end{itemize}
\end{thm}
\begin{proof}
	The work of R{\"o}ndigs and {\O}stv{\ae}r in \cite{MR2435654}, shows
	in particular that over a field of characteristic zero, the motivic stable homotopy category $\stableHZmodules$ of modules over
	the motivic Eilenberg-MacLane spectrum $H\mathbb Z$ is equivalent to Voevodsky's big
	category of motives $DM_{k}$, where the equivalence preserves the monoidal and triangulated
	structures (see \cite[theorem 1.1]{MR2435654}). 
	
	Therefore, the result is an immediate consequence of theorem \ref{thm.3.5.q-slices=====HZ-modules}
\end{proof}

\end{section}

\end{chapter}

\backmatter
\bibliography{biblio_pablo}

\begin{thebibliography}{10}

\bibitem{MR0354652}
{\em Th\'eorie des topos et cohomologie \'etale des sch\'emas. {T}ome 1:
  {T}h\'eorie des topos}.
\newblock Springer-Verlag, Berlin, 1972.
\newblock S\'eminaire de G\'eom\'etrie Alg\'ebrique du Bois-Marie 1963--1964
  (SGA 4), Dirig\'e par M. Artin, A. Grothendieck, et J. L. Verdier. Avec la
  collaboration de N. Bourbaki, P. Deligne et B. Saint-Donat, Lecture Notes in
  Mathematics, Vol. 269.

\bibitem{MR2102294}
W.~G. Dwyer, P.~S. Hirschhorn, D.~M. Kan, and J.~H. Smith.
\newblock {\em Homotopy limit functors on model categories and homotopical
  categories}, volume 113 of {\em Mathematical Surveys and Monographs}.
\newblock American Mathematical Society, Providence, RI, 2004.

\bibitem{MR1949356}
E.~M. Friedlander and A.~Suslin.
\newblock The spectral sequence relating algebraic {$K$}-theory to motivic
  cohomology.
\newblock {\em Ann. Sci. \'Ecole Norm. Sup. (4)}, 35(6):773--875, 2002.

\bibitem{MR1650938}
P.~G. Goerss and J.~F. Jardine.
\newblock Localization theories for simplicial presheaves.
\newblock {\em Canad. J. Math.}, 50(5):1048--1089, 1998.

\bibitem{MR1711612}
P.~G. Goerss and J.~F. Jardine.
\newblock {\em Simplicial homotopy theory}, volume 174 of {\em Progress in
  Mathematics}.
\newblock Birkh\"auser Verlag, Basel, 1999.

\bibitem{Over-Und}
P.~Hirschhorn.
\newblock Overcategories and undercategories of model categories.
\newblock {\em Preprint 2005}.

\bibitem{MR1944041}
P.~S. Hirschhorn.
\newblock {\em Model categories and their localizations}, volume~99 of {\em
  Mathematical Surveys and Monographs}.
\newblock American Mathematical Society, Providence, RI, 2003.

\bibitem{MR2197578}
J.~Hornbostel.
\newblock Localizations in motivic homotopy theory.
\newblock {\em Math. Proc. Cambridge Philos. Soc.}, 140(1):95--114, 2006.

\bibitem{HoveyPreprint}
M.~Hovey.
\newblock Monoidal model categories.
\newblock {\em Preprint 1998}.

\bibitem{MR1650134}
M.~Hovey.
\newblock {\em Model categories}, volume~63 of {\em Mathematical Surveys and
  Monographs}.
\newblock American Mathematical Society, Providence, RI, 1999.

\bibitem{MR1860878}
M.~Hovey.
\newblock Spectra and symmetric spectra in general model categories.
\newblock {\em J. Pure Appl. Algebra}, 165(1):63--127, 2001.

\bibitem{MR1695653}
M.~Hovey, B.~Shipley, and J.~Smith.
\newblock Symmetric spectra.
\newblock {\em J. Amer. Math. Soc.}, 13(1):149--208, 2000.

\bibitem{MR906403}
J.~F. Jardine.
\newblock Simplicial presheaves.
\newblock {\em J. Pure Appl. Algebra}, 47(1):35--87, 1987.

\bibitem{MR1787949}
J.~F. Jardine.
\newblock Motivic symmetric spectra.
\newblock {\em Doc. Math.}, 5:445--553 (electronic), 2000.

\bibitem{MR0360749}
G.~M. Kelly.
\newblock Doctrinal adjunction.
\newblock In {\em Category Seminar (Proc. Sem., Sydney, 1972/1973)}, pages
  257--280. Lecture Notes in Math., Vol. 420. Springer, Berlin, 1974.

\bibitem{MR2365658}
M.~Levine.
\newblock The homotopy coniveau tower.
\newblock {\em J. Topol.}, 1(1):217--267, 2008.

\bibitem{MR0060829}
W.~S. Massey.
\newblock Products in exact couples.
\newblock {\em Ann. of Math. (2)}, 59:558--569, 1954.

\bibitem{MR1813224}
F.~Morel and V.~Voevodsky.
\newblock {${\bf A}\sp 1$}-homotopy theory of schemes.
\newblock {\em Inst. Hautes \'Etudes Sci. Publ. Math.}, (90):45--143 (2001),
  1999.

\bibitem{MR1308405}
A.~Neeman.
\newblock The {G}rothendieck duality theorem via {B}ousfield's techniques and
  {B}rown representability.
\newblock {\em J. Amer. Math. Soc.}, 9(1):205--236, 1996.

\bibitem{MR1812507}
A.~Neeman.
\newblock {\em Triangulated categories}, volume 148 of {\em Annals of
  Mathematics Studies}.
\newblock Princeton University Press, Princeton, NJ, 2001.

\bibitem{MR0223432}
D.~G. Quillen.
\newblock {\em Homotopical algebra}.
\newblock Lecture Notes in Mathematics, No. 43. Springer-Verlag, Berlin, 1967.

\bibitem{MR2435654}
O.~R{\"o}ndigs and P.~A. {\O}stv{\ae}r.
\newblock Modules over motivic cohomology.
\newblock {\em Adv. Math.}, 219(2):689--727, 2008.

\bibitem{MR1734325}
S.~Schwede and B.~E. Shipley.
\newblock Algebras and modules in monoidal model categories.
\newblock {\em Proc. London Math. Soc. (3)}, 80(2):491--511, 2000.

\bibitem{MR1764202}
V.~Voevodsky.
\newblock Triangulated categories of motives over a field.
\newblock In {\em Cycles, transfers, and motivic homology theories}, volume 143
  of {\em Ann. of Math. Stud.}, pages 188--238. Princeton Univ. Press,
  Princeton, NJ, 2000.

\bibitem{MR1977582}
V.~Voevodsky.
\newblock Open problems in the motivic stable homotopy theory. {I}.
\newblock In {\em Motives, polylogarithms and Hodge theory, Part I (Irvine, CA,
  1998)}, volume~3 of {\em Int. Press Lect. Ser.}, pages 3--34. Int. Press,
  Somerville, MA, 2002.

\end{thebibliography}
\bibliographystyle{abbrv}

\end{document}